\documentclass[twoside,12pt]{article}
\DeclareMathSizes{20}{30}{16}{12}
\usepackage{float}
\usepackage{amsmath}
\usepackage{mathtools}
\usepackage{mathrsfs}
\usepackage{stmaryrd}
\usepackage{bbm}
\usepackage{yfonts}
\usepackage{amsfonts}
\usepackage{tcolorbox}
\usepackage{tikz-cd}
\usepackage[utf8]{inputenc}

\usepackage{caption}
\usepackage{indentfirst}
\usepackage[Symbol]{upgreek}
\usepackage{enumerate}   
\usepackage{upref}
\usepackage{hyphenat}
\usepackage[absolute,overlay]{textpos}
\usepackage{pdfpages}

\usepackage{graphicx}
\usepackage[margin=0.8 in]{geometry}
\usepackage{fancyhdr}
\usepackage{hyperref}
\hypersetup{
    colorlinks,
    citecolor=black,
    filecolor=black,
    linkcolor=black,
    urlcolor=black
}
\newcommand{\N}{\mathbb{N}}
\newcommand{\Z}{\mathbb{Z}}
\newcommand{\R}{\mathbb{R}}
\newcommand{\ra}{\rightarrow}

\newcommand{\I}{\mathfrak{I}}
\newcommand{\calA}{\mathcal{A}}
\newcommand{\calI}{\mathcal{I}}
\newcommand{\calJ}{\mathcal{J}}
\newcommand{\calE}{\mathcal{E}}
\newcommand{\calK}{\mathcal{K}}
\newcommand{\calC}{\mathcal{C}}
\newcommand{\calT}{\mathcal{T}}
\newcommand{\calD}{\mathcal{D}}

\newcommand{\calR}{\mathcal{R}}
\newcommand{\calB}{\mathcal{B}}

\newcommand{\calM}{\mathcal{M}}
\newcommand{\calN}{\mathcal{N}}
\newcommand{\calP}{\mathcal{P}}
\newcommand{\calQ}{\mathcal{Q}}
\newcommand{\calV}{\mathcal{V}}

\newcommand{\J}{\mathfrak{J}}
\newcommand{\op}{\operatorname}
\newcommand{\w}{\widehat}

\newcommand{\Ox}{\mathcal{O}}
\newcommand{\F}{\mathcal{F}}
\newcommand{\calG}{\mathcal{G}}

\newcommand{\ov}{\overline}

\newcommand{\frakA}{\mathfrak{A}}
\newcommand{\frakB}{\mathfrak{B}}
\newcommand{\frakX}{\mathfrak{X}}
\newcommand{\frakE}{\mathfrak{E}}
\newcommand{\frakF}{\mathfrak{F}}
\newcommand{\frakU}{\mathfrak{U}}
\newcommand{\frakV}{\mathfrak{V}}

\newcommand{\frakS}{\mathcal{S}}
\newcommand{\frakY}{\mathfrak{Y}}
\newcommand{\frakT}{\mathfrak{T}}
\newcommand{\frakZ}{\mathfrak{Z}}

\newcommand{\frakG}{\mathfrak{G}}

\newcommand{\Sp}{\operatorname{Spec}}

\long\def\/*#1*/{}

\usepackage{unicode-math}
\usepackage{fontspec}

\usepackage{sectsty}
\sectionfont{\centering}

\usepackage[toc,page]{appendix}
\usepackage{amsthm}

\newtheorem*{proof*}{Proof}

\newtheorem{theorem}[subsection]{Theorem}

\newtheorem{theorem*}{Theorem}
\newtheorem{proposition}[subsection]{Proposition}
\newtheorem{corollaire}[subsection]{Corollary}
\newtheorem{lemma}[subsection]{Lemma}

\theoremstyle{definition}

\newtheorem{definition}[subsection]{Definition}
\newtheorem{remark}[subsection]{Remarks}

\newtheorem{exmp}[subsection]{Examples}

\theoremstyle{definition}
\newtheorem{parag}[subsection]{}

\numberwithin{equation}{subsection}

\title{Logarithmic Cartier Transform}
\author{Sami Fersi}
\date{}

\AtEndDocument{\bigskip{\footnotesize%
  \textsc{Sami Fersi, Laboratoire Alexander Grothendieck, Institut des Hautes Études Scientifiques, 35 Route de Chartres, 91440 Bures-sur-Yvette, France} \par
  \textit{E-mail address}: \texttt{fersi@ihes.fr} \par
}}
\bibliography{./references/ref}
\begin{document}

\maketitle

\begin{abstract}
We generalize the Cartier transform of Ogus and Vologodsky to log smooth schemes. More precisely, we generalize a local version of this transform, due to Shiho, and a topos-theoretic version, due to Oyama. Let $k$ be a perfect field of positive characteristic $p$ and equip $\op{Spec}k$ with the trivial log structure. For a log smooth morphism of logarithmic schemes $X \ra S,$ where $S$ is log flat and locally of finite type over $\op{Spec}k,$ we obtain, under the assumption that the exact relative Frobenius lifts over the Witt vectors of $k,$ a fully faithful functor from the category of quasi-coherent modules on the base change $X'=X\times_{S,F_S}S$ of $X$ by the Frobenius $F_S$ of $S,$ equipped with a quasi-nilpotent Higgs field, to the category of quasi-coherent modules on $X$ equipped with a quasi-nilpotent integrable connection.
In another direction, we construct crystalline-like topoi and subcategories of crystals $\calC'$ and $\underline{\calC},$ equivalent respectively to modules with Higgs fields and integrable connections, and a fully faithful functor $\calC' \ra \underline{\calC}.$ Since the Frobenius morphism is not, in general, flat in the log smooth setting, it is not clear that these functors are essentially surjective. To address this issue, we refine the topoi and crystals mentioned above by endowing them with an indexed structure, inspired by Lorenzon’s extension of Cartier descent to smooth logarithmic schemes. Using the Azumaya property of the ring of logarithmic differential operators, we then obtain an equivalence between the corresponding categories of indexed crystals, thereby generalizing the Cartier transform.
\end{abstract}

\tableofcontents

\section{Introduction}

\begin{parag}
In this manuscript, our goal is to generalize the \emph{Cartier transform} of Ogus and Vologodsky to the logarithmic case. In the classical smooth case, this transform can be seen as a generalization of two results for schemes in positive characteristic: \emph{Cartier descent} and \emph{Deligne-Illusie decomposition of the de Rham complex}.
For a smooth morphism of schemes $X\ra S$ of positive characteristic $p>0,$ Cartier descent is an equivalence between the category of quasi-coherent $\Ox_X$-modules equipped with an integrable connection with vanishing $p$-curvature, and the category of quasi-coherent $\Ox_{X'}$-modules, where $X'=X\times_{S,F_S}S$ is the base change of $X$ by the absolute Frobenius $F_S$ of $S$ (\cite{Katz} 5.1). Under this equivalence, an $\Ox_{X'}$-module $\calE'$ corresponds to the module $F_{X/S}^*\calE',$ where $F_{X/S}:X\ra X'$ is the relative Frobenius, equipped with the canonical connection
$$
\nabla^{\mathrm{can}}:
\begin{array}[t]{clc}
F_{X/S}^*\calE' & \ra & \left (F_{X/S}^*\calE'\right ) \otimes_{\Ox_X}\Omega^1_{X/S}\\
aF_{X/S}^* x & \mapsto & \left (F_{X/S}^*x\right ) \otimes da,
\end{array}
$$
for local sections $a$ and $x$ of $\Ox_X$ and $\calE'$ respectively.
On the other hand, Deligne and Illusie proved the following: if the absolute Frobenius $F_S$ lifts over $W_2(k)$ and if the relative Frobenius $F_{X/S}$ lifts to a morphism between smooth schemes over $W_2(k),$ then there exists a quasi-isomorphism
$$
\bigoplus_{i\in \N} \Omega^i_{X'/S}[-i] \xrightarrow{\sim} F_{X/S*}\Omega^{\bullet}_{X/S},
$$
inducing the Cartier isomorphism on cohomologies. If only $X'$ lifts to a smooth scheme over $W_2(k),$ then there exists an isomorphism in $D(\Ox_{X'})$
$$
\bigoplus_{0\le i<p} \Omega_{X'/S}^i[-i] \xrightarrow{\sim} F_{X/S*}\tau_{<p} \Omega^{\bullet}_{X/S},
$$
where $\tau_{<p}$ is the truncation functor up to $p-1.$ The \emph{Cartier transform} of Ogus and Vologodsky comes as a generalization of these two results. More precisly, the Cartier transform is an equivalence between the categories of $\Ox_{X'}$-modules equipped with Higgs fields and $\Ox_X$-modules equipped with integrable connections, both satisfying certain nilpotence conditions. In addition, this transform is compatible with the natural cohomologies.

Shiho \cite{Shiho} provided a local construction of the Cartier transform under the assumption that the relative Frobenius $F_{X/S}$ lifts to a morphism of smooth formal schemes over the Witt vectors $W(k)$ of a perfect field $k$ of characteristic $p.$ More precisely, he interpreted modules with Higgs fields and integrable connections in terms of stratifications relative to certain groupoids. Using faithfully flat descent, he then proved the equivalence for modules with stratifications and deduced the one for Higgs fields and connections. The groupoids are defined as PD-envelopes of dilatations of the ideal of the diagonal embedding of the formal scheme. Oyama used similarly constructed groupoids to provide a crystalline-like interpretation of the Cartier transform \cite{Oyama}. As opposed to Shiho's, Oyama's construction is global and doesn't depend on a lifting of $F_{X/S}.$

While the Deligne-Illusie decomposition of the de Rham complex generalizes to the log smooth case without problems, as was proved by Kato in (\cite{Kat89} 4.12), one encounters a major problem coming from residues when attempting to generalize Cartier descent, as was noted by Ogus in (\cite{Ogus94} 1.3). Lorenzon provided a solution to this by introducing, for a morphism of logarithmic schemes $f:X\ra S,$ \emph{indexed algebras} $\calA_X$ and $\calB_{X/S}$ \cite{Lor2000}. These indexed algebras are topos-theoretic versions of graded algebras in the sense that each local section possesses an \emph{index}, adding two local sections is only possible if they have the same index and multiplying them adds the corresponding indices. The indexed algebra $\calA_X$ is naturally equipped with a connection $d_{\calA_X}$ and $\calB_{X/S}$ is taken to be its sheaf of horizontal sections. Lorenzon proves that, if $f$ is log smooth, then locally free $\calB_{X/S}$-modules of finite type are equivalent to locally free $\calA_X$-modules of finite type, equipped with an integrable connection, which has vanishing $p$-curvature and is compatible with $d_{\calA_X}.$
 
In this thesis, we start by generalizing Shiho functor to log smooth logarithmic schemes. When attempting to define groupoids similar to those introduced by Shiho, one major problem we encounter is the \emph{exactification} of the diagonal immersion: decomposing the diagonal immersion into an exact immersion followed by a log étale morphism. This is possible locally or, globally if the logarithmic schemes are equipped with \emph{frames}, a notion introduced by Kato and Saito in \cite{Saito04}. We hence work systematically with framed logarithmic schemes and we can thus construct the two groupoids corresponding to Higgs fields and integrable connections and the functor between them. We show that it is fully faithful using a log flat descent theorem for morphisms, that we prove based on results by Kato. We don't have an equivalence at this level though. We also generalize Oyama's work: we construct two crystalline-like topoi, a morphism between them and subcategories of crystals $\calC$ and $\calC'.$ We prove that the pullback functor preserves crystals and so it induces a functor $\calC' \ra \calC$ that gives the Cartier transform. We also prove that it is fully faithful.
In the smooth case, faithfully flat descent is crucial in proving that these two functors are essentially surjective. Using the same argument is not possible in the log smooth case. That is why we introduce the indexed algebras and modules of Lorenzon and shift to an \emph{Azumaya algebra approach}. Ohkawa and Schepler introduce an \emph{indexed sheaf of differential operators} $\widetilde{\calD}_{X/S}$ and prove that it is an Azumaya algebra over its center (\cite{Schepler} and \cite{Ohkawa}). We review their work in more detail and then construct indexed versions of Oyama crystals that represent $\calB_{X/S}$-modules with Higgs fields and $\calA_X$-modules with integrable connections. We finally construct an equivalence of categories between them.
\end{parag}

\begin{parag}
We now describe our work in more details. We fix a prime number $p.$ For a morphism $X\ra S$ of fs logarithmic schemes of characteristic $p,$ let $X''$ be the base change, in the category of fine logarithmic schemes, of $X$ by the absolute Frobenius morphism $F_S:S\ra S$ of $S.$ Let $F_{X/S}$ be the relative Frobenius morphism of $X$ with respect to $S$ i.e. the unique morphism $X\ra X''$ making the following diagram commutative
\begin{equation}
\begin{tikzcd}
X\ar{dr}{F_{X/S}}\ar[bend right=-30]{drr}{F_X}\ar[swap,bend right=30]{ddr} & & \\
 & X''\ar{r}\ar{d} & X\ar{d} \\
 & S\ar{r}{F_S} & S
\end{tikzcd}
\end{equation}
Recall that $F_{X/S}$ is \emph{weakly inseparable} (\cite{Ogus2018} III 2.4), so $F_{X/S}$ factors uniquely as a composition of a log étale morphism $G:X'\ra X''$ and an inseparable morphism $F:X\ra X'$ (\cite{Ogus2018} IV 3.3.8): 
\begin{equation}
\begin{tikzcd}
X\ar{r}{F}\ar{dr}\ar[bend right=30]{rdd} & X'\ar{d}{G}\ar{dr}{\pi} & \\
 & X''\ar{d}\ar{r} & X\ar{d} \\
 & S \ar{r}{F_S} & S
\end{tikzcd}
\end{equation}
The morphism $F:X\ra X'$ is called the \emph{exact relative Frobenius of $X$ with respect to $S$} (\cite{Ogus2018} IV 3.3.9).
In section 3, we prove that the exact relative Frobenius of a log smooth morphism of fs logarithmic schemes is log flat \eqref{thmlogflat}. We also define, for an fs logarithmic scheme $T$ of characteristic $p,$ the logarithmic scheme theoretic image of the Frobenius $F_T$ \eqref{era3logimagedef}, denoted by $\underline{T}.$ We then prove, in \ref{era4prop15}, that the absolute Frobenius $F_T$ factors through $\underline{T}:$
\begin{equation}\label{IranG}
F_T:T\xrightarrow{G_T} \underline{T} \hookrightarrow T.
\end{equation}
We also prove, in \ref{propfT/S} and under the assumption that $S=\Sp k$ equipped with the trivial logarithmic structure, that the relative Frobenius $F_{T/S}:T \ra T'$ factors through $\underline{T}':$
\begin{equation}\label{Iranf}
T\xrightarrow{f_{T/S}} \underline{T}' \hookrightarrow T'.
\end{equation}
\end{parag}

\begin{parag}
We fix a perfect field $k$ of characteristic $p$ and we denote its ring of Witt vectors by $W(k).$ If a gothic letter $\frakX$ denotes a logarithmic $p$-adic formal scheme over $W(k),$ then the corresponding roman letter $X=\frakX \times_{\op{Spf}W(k)} \op{Spec}k$ will denote its special fiber.
We consider a log smooth morphism $f:\frakX \ra \frakS$ of fs $p$-adic logarithmic formal schemes, flat over $\op{Spf}W(k).$ Denote by $f_1:X \ra S$ its special fiber.
Following Shiho's construction, we suppose that the exact relative Frobenius  $F_1:X\ra X'$ lifts to an $\frakS$-morphism of logarithmic formal schemes $F:\frakX\ra \frakX',$ such that $\frakX'$ is flat over $\op{Spf}W(k).$ The flatness condition implies that there exists a unique morphism $\frac{dF}{p}$ fitting into the commutative diagram
$$
\begin{tikzcd}
F^*\omega^1_{\frakX'/\frakS} \ar{r}{dF} \ar[swap]{dr}{ \frac{dF}{p}} & p\omega^1_{\frakX / \frakS} \\
& \omega^1_{\frakX / \frakS} \ar{u}{\times p}.
\end{tikzcd}
$$
Let $n\ge 0$ be an integer. Let $\calE'$ be an $\Ox_{\frakX'}$-module equipped with a $p^{n+1}$-connection
$$
\nabla':\calE' \ra \calE' \otimes_{\Ox_{\frakX'}} \omega^1_{\frakX' / \frakS},
$$
i.e. an additive morphism satisfying, for local sections $a$ and $x$ of $\Ox_{\frakX'}$ and $\calE'$ respectively, the modified Leibniz rule
$$
\nabla'(ax)=a\nabla'(x)+p^{n+1}x\otimes da.
$$
We obtain an $\Ox_{\frakX}$-module $\calE=F^*\calE'$ equipped with a $p^n$-connection
$$
\nabla:\begin{array}[t]{clc}
\calE & \ra & \calE \otimes_{\Ox_{\frakX}} \omega^1_{\frakX / \frakS} \\
a\otimes x' & \mapsto &  \zeta(a\otimes x')+F^*x'\otimes da, 
\end{array}
$$
where $a$ and $x'$ are local sections of $\Ox_{\frakX}$ and $\calE'$ respectively and $\zeta$ is the composition
$$
\zeta:\calE \xrightarrow{F^*\nabla'} \calE \otimes_{\Ox_{\frakX}}F^*\omega^1_{\frakX' / \frakS} \xrightarrow{\op{Id}_{\calE} \otimes \frac{dF}{p}} \calE \otimes_{\Ox_{\frakX}}\omega^1_{\frakX / \frakS}.
$$
If $\nabla'$ is integrable then so is $\nabla.$
This construction thus yields a functor
$$
p^{n+1}\text{-}\op{MIC}(\frakX'/\frakS) \ra p^n\text{-}\op{MIC}(\frakX/\frakS),
$$
where $p^{n}\text{-}\op{MIC}(\frakX/\frakS)$ is the category of $\Ox_{\frakX}$-modules equipped with integrable $p^{n}$-connections.
In particular, if $n=0$ and if we work modulo $p^l$ for some positive integer $l,$ we have a functor
\begin{equation}\label{intro1}
\Phi_l:p\text{-}\op{MIC}(\frakX_l'/\frakS_l) \ra \op{MIC}(\frakX_l/\frakS_l),
\end{equation}
where $\frakX_l$ (resp $\frakS_l$ )denotes the logarithmic scheme obtained from $\frakX$ (resp $\frakS$) by reduction modulo $p^l.$
To study this functor, we interpret modules equipped with quasi-nilpotent integrable $p^n$-connections as modules with stratifications with respect to certain groupoids. To introduce these groupoids, we need to exactify the diagonal immersion. This is possible using \emph{frames}, a notion introduced by Saito and Kato \cite{Saito04}.
\end{parag}

\begin{parag}
Let $\boldsymbol{\op{L}}$ be the category of fine saturated logarithmic schemes. For an fs monoid $M,$ we denote by $[M]$ the presheaf
$$
[M]:\begin{array}[t]{clc}
\boldsymbol{\op{L}} & \ra & \boldsymbol{\op{Sets}} \\
T & \mapsto & \op{Hom} \left (M,\Gamma(T,\ov{\calM}_T) \right ).
\end{array}
$$
A frame on a logarithmic scheme or a logarithmic formal scheme $T$ is a morphism of presheaves $T \ra [M],$ or equivalently a morphism of monoids $M \ra \Gamma\left (T,\ov{\calM}_T \right ),$ 
that lifts, étale locally on $T,$ to a chart $T \ra \Gamma \left (T,\calM_T \right ).$ An important result we use is a direct generalization for formal schemes of a result by Kato and Saito (\cite{Saito04} 4.2.8): suppose that $\frakX$ and $\frakS$ are fs logarithmic $p$-adic formal schemes equipped with frames $Q\ra \Gamma(\frakX,\ov{\calM}_{\frakX})$ and $P\ra \Gamma(\frakS,\ov{\calM}_{\frakS})$ respectively and that we are given a morphism of monoids $\theta:P \ra Q$ such that $(f,\theta):(\frakX,Q) \ra (\frakS,P)$ is a morphism of framed logarithmic formal schemes, i.e. the diagram
$$
\begin{tikzcd}
P \ar{r}{\theta} \ar{d} & Q\ar{d} \\
\Gamma \left (\frakS,\ov{\calM}_{\frakS} \right ) \ar[swap]{r}{f^{\flat}} & \Gamma \left (\frakX,\ov{\calM}_{\frakX} \right )
\end{tikzcd}
$$
is commutative. The presheaf on $\boldsymbol{\op{L}}$
$$
\frakX \times_{\frakS,[Q]}^{\op{log}} \frakX: T \mapsto \frakX(T) \times_{\frakS(T)\times [Q](T)}^{\op{log}}\frakX(T),
$$
is then representable by a logarithmic $p$-adic formal scheme which is affine and log étale over $\frakX \times_{\frakS}^{\op{log}}\frakX.$
Furthermore, the diagonal immersion $\frakX \ra \frakX \times_{\frakS}^{\op{log}}\frakX$ factors as
$$
\begin{tikzcd}
 & \frakX \times_{\frakS,[Q]}^{\op{log}} \frakX \ar{d} \\
\frakX \ar{r} \ar{ur} & \frakX \times_{\frakS}^{\op{log}}\frakX,
\end{tikzcd}
$$
where $\frakX \ra \frakX \times_{\frakS,[Q]}^{\op{log}} \frakX$ is strict and the conormal sheaf of the immersion
$$\frakX_l \ra \frakX_l \times_{\frakS_l,[Q]}^{\op{log}} \frakX_l$$
is canonically isomorphic to the module of logarithmic differentials $\omega^1_{\frakX_l/\frakS_l},$ for all positive integers $l.$
\end{parag}

\begin{parag}
Let $X \ra S$ be a morphism of fs logarithmic schemes of characteristic $p$ and denote by $\calT_{X/S}=\mathscr{Hom}_{\Ox_X} \left ( \omega^1_{X/S}, \Ox_X \right )$ the dual of the module of logarithmic differentials $\omega^1_{X/S}.$ As shown in (\cite{Ogus94} 1.2.1), the Lie algebra of derivations $\calT_{X/S}$ can be equipped with a restricted Lie algebra structure via an operation
$$
\partial \mapsto \partial^{(p)}.
$$
Seeing derivations as differential operators of $\calD_{X/S},$ we obtain a morphism
$$
F_X^*\calT_{X/S} \ra \calD_{X/S},\ \partial \mapsto \partial^p-\partial^{(p)}.
$$
By adjunction, this induces a morphism
\begin{equation*}
\psi:\calT_{X'/S} \ra F_{*}\calD_{X/S}, 
\end{equation*}
where $F$ is the exact relative Frobenius. For any derivation $\partial,$ $\psi(\partial)$ is in the center of $F_*\calD_{X/S}$ and so we obtain a morphism of $\Ox_X$-algebras
\begin{equation}\label{intropsi}
\psi:S^{\bullet}\calT_{X'/S} \ra F_{*}\calD_{X/S}, 
\end{equation}
map $\psi$ is called the \emph{$p$-curvature map}.
\end{parag}

\begin{parag}\label{intro16}
We equip $\op{Spf}W(k)$ with the trivial logarithmic structure. We say that a morphism of logarithmic $p$-adic formal schemes over $\op{Spf}W(k)$ is log flat if it is so modulo $p^l$ for all $l\ge 1.$ We consider an fs logarithmic $p$-adic formal scheme $\frakS$ log flat and locally of finite type over $\op{Spf}W(k).$ Note that this implies that $\frakS$ is flat over $\op{Spf}W(k)$ (\ref{dxuflat} and \ref{Wflat}). We also consider a log smooth morphism $f:(\frakX,Q) \ra (\frakS,P)$ of framed fs logarithmic $p$-adic formal schemes. Let $\calI$ be the ideal of the exact diagonal immersion $\frakX \ra \frakY=\frakX \times_{\frakS,[Q]}^{\op{log}}\frakX.$
For a positive integer $n,$ let $R_{\frakX,n}$ (resp. $Q_{\frakX}$) be the dilatation of $\calI+(p)$ with respect to $p^n$ (resp. the dilatation of the ideal locally generated by $\{p,a^p,a\in \calI \}$ with respect to $p$). In other words, $R_{\frakX,n}$ is the largest open formal subscheme of the admissible blow-up of $\calI+(p^n)$ such that $\left ( \calI+(p^n) \right ) \Ox_{R_{\frakX,n}}=(p^n).$ We then consider, for every positive integer $l,$ the PD-envelope $\left (P_{\frakX/\frakS,n} \right )_l$ of the immersion $\frakX_l \ra \left (R_{\frakX,n} \right )_l$ and equip it with the logarithmic structure pull back of that of $\frakY_l.$ We consider the inductive limit
$$
P_{\frakX/\frakS,n} =\lim\limits_{\substack{\longrightarrow \\l\ge 1}}\left (P_{\frakX/\frakS,n}\right )_l.
$$
The logarithmic schemes $\left (P_{\frakX/\frakS,n} \right )_l$ provide stratified interpretations for $p^n$-connections. More precisely, we prove that the logarithmic schemes $R_{\frakX,n},$ $Q_{\frakX}$ and $P_{\frakX/\frakS,n}$ have natural structures of formal groupoids and that the category $n\text{-}\op{MHS}(\frakX_l/\frakS_l)$ of $\Ox_{\frakX_l}$-modules equipped with stratifications relative to the formal groupoid $\left (P_{\frakX/\frakS,n} \right )_l,$ is equivalent to the category $p^n\text{-}\op{MIC}^{\text{qn}}\left (\frakX_l/\frakS_l\right )$ of $\Ox_{\frakX_l}$-modules equipped with a quasi-nilpotent integrable $p^n$-connection. We also prove the following proposition
\begin{theorem}[\ref{thm1237}]
Let $f:(\frakX,Q) \ra (\frakS,P)$ be a log smooth morphism of framed fs logarithmic $p$-adic formal schemes such that $\frakS$ is log flat and locally of finite type over $\op{Spf}W(k).$ Denote by $X\ra S$ its special fiber. Let $\calD_{X/S}$ be the sheaf of logarithmic differential operators, $\w{\Gamma}^{\bullet}\calT_{X'/S}$ the completed PD-algebra of the tangent sheaf $\calT_{X'/S}.$ We consider the $p$-curvature map $S^{\bullet}\calT_{X'/S} \ra \calD_{X/S}$ \eqref{intropsi} and 
$$\calD^{\gamma}_{X/S}=\calD_{X/S}\otimes_{S^{\bullet} \calT_{X'/S}} \widehat{\Gamma}^{\bullet}\calT_{X'/S}.$$
The following tensor categories are canonically equivalent:
\begin{enumerate}
\item The category of $\Ox_X$-modules equipped with a stratification relative to $R_{\frakX,1}$ (resp. $Q_{\frakX}$).
\item The category of locally PD-nilpotent $\widehat{\Gamma}^{\bullet}\calT_{X/S}$-modules (resp. locally PD-nilpotent $\calD^{\gamma}_{X/S}$-modules).
\end{enumerate}
\end{theorem}
Furthermore, if we suppose that the exact relative Frobenius $F_1:X\ra X'$ lifts to a morphism of framed logarithmic formal $\frakS$-schemes $F:\frakX \ra \frakX',$ where $\frakX'$ is log smooth over $\frakS,$ then $F$ induces two morphisms
\begin{equation}\label{Muzannu}
\varphi:P_{\frakX/\frakS,0} \ra P_{\frakX'/\frakS,1},\ \nu:Q_{\frakX}  \ra R_{\frakX',1}.
\end{equation}
In addition, in \ref{lem96} by proving that the ideal of the diagonal $\frakX \ra \frakX\times_{\frakX'}^{\op{log}}\frakX$ has a unique PD-structure, we prove the existence of a morphism
$$
\Psi:\frakX\times_{\frakX'}^{\op{log}}\frakX \ra P_{\frakX/\frakS,0}.
$$
These morphisms fit into commutative diagrams
$$
\begin{tikzcd}
 & & \frakX'\ar{d}{\iota'} \\
 & P_{\frakX/\frakS,0} \ar{r}{\varphi} \ar{d} & P_{\frakX'/\frakS,1}\ar{d} \\
\frakX\times_{\frakX'}^{\op{log}}\frakX \ar{ur}{\psi} \ar[bend right=-30]{uurr} \ar{r} & \frakX\times^{\op{log}}_{\frakS,[Q]}\frakX \ar{r}{F^2} & \frakX'\times_{\frakS,[Q']}^{\op{log}}\frakX'
\end{tikzcd}
$$
and
$$
\begin{tikzcd}
Q_{\frakX} \ar{r}{\nu} \ar{d} & R_{\frakX',1} \ar{d} \\
\frakX\times^{\op{log}}_{\frakS,[Q]}\frakX \ar{r}{F^2} & \frakX'\times_{\frakS,[Q']}^{\op{log}}\frakX'.
\end{tikzcd}
$$
The morphism $\varphi$ allows us to define a functor
\begin{equation}\label{Psiintro}
\Psi_l:\begin{array}[t]{clc}
1\text{-}\op{MHS}\left (\frakX_l'/\frakS_l\right ) & \ra & 0\text{-}\op{MHS}\left (\frakX_l/\frakS_l\right ) \\
(\calE',\epsilon') & \mapsto & (F_l^*\calE',\varphi_l^*\epsilon').
\end{array}
\end{equation}
\begin{theorem}[\ref{THMX1}, \ref{THMX2}]
Let $f:(\frakX,Q) \ra (\frakS,P)$ be a log smooth morphism of framed fs logarithmic $p$-adic formal schemes such that $\frakS$ is log flat and locally of finite type over $\op{Spf}W(k).$ Denote by $X\ra S$ its special fiber. We suppose that the exact relative Frobenius $X\ra X'$ lifts to a morphism of framed fs $p$-adic logarithmic formal schemes $F:(\frakX,Q) \ra (\frakX',Q')$ over $(\frakS,P),$ such that $\frakX'$ is log smooth over $\frakS.$ Let $n$ be a positive integer. The functors $\Phi_l$ \eqref{intro1} and $\Psi_l$ \eqref{Psiintro} induce fully faithful functors between quasi-coherent and quasi-nilpotent objects, fitting into a commutative diagram
$$
\begin{tikzcd}
1\text{-}\op{MHS}^{qcoh}(\frakX_l'/\frakS_l) \ar{rr}{\Psi_l} \ar[swap,sloped]{d}{\sim} & & 0\text{-}\op{MHS}^{qcoh}(\frakX_l/\frakS_l) \ar[sloped]{d}{\sim} \\
p\text{-}\op{MIC}^{qcoh,qn}(\frakX_l'/\frakS_l) \ar{rr}{\Phi_l} & & \op{MIC}^{qcoh,qn}(\frakX_l/\frakS_l).
\end{tikzcd}
$$
\end{theorem}
To show the full faithfulness, we prove a kind of logarithmic flat descent theorem for quasi-coherent modules, slightly generalizing a result of Kato for structural rings \cite{Kat19}:
\begin{theorem}[\ref{eraflogflatdescent}]
Let $T$ be an fs logarithmic scheme, $(T_i \xrightarrow{f_i} T)_{i\in I}$ a log flat covering and $\calE$ a quasi-coherent $\Ox_T$-module. For any $i,j\in I,$ let $f_{ij}:T_i\times_T^{\op{log}}T_j \ra T$ be the canonical morphism. Then, the sequence of $\Ox_T$-modules
$$0 \ra \calE \ra \prod_{i\in I}f_{i*}f_i^*\calE \ra \prod_{i,j\in I}f_{ij*}f_{ij}^*\calE,$$
where the last arrow is the difference between the morphisms induced by the projections $T_i\times_T^{\op{log}}T_j \ra T_i$ and $T_i\times_T^{\op{log}}T_j \ra T_j,$ is exact.
\end{theorem}
\end{parag}

\begin{parag}\label{110}
In the classical smooth case, Oyama provided a crystalline-like interpretation of the Cartier transform. We follow his work and provide a similar interpretation in the log smooth case. Let $W(k)$ be the ring of Witt vectors of a perfect field $k$ of positive characteristic $p.$ We equip $\op{Spf}W(k)$ with the trivial logarithmic structure. Let $\frakS$ be a log flat and locally of finite type fs logarithmic $p$-adic formal scheme over $\op{Spf}W(k).$ Denote by $S$ its special fiber and consider a log smooth morphism $X\ra S$ of fs logarithmic schemes.
We define categories $\calE(X/\frakS)$ and $\underline{\calE}(X/\frakS)$ \eqref{defforintro} as follows: an object of $\calE(X/\frakS)$ (resp. $\underline{\calE}(X/\frakS)$) is a triple $(U,\frakT,u)$ consisting of an étale strict morphism of logarithmic schemes $U \ra X,$ a log flat fs logarithmic $p$-adic formal $\frakS$-scheme $\frakT$ and an $S$-morphism of logarithmic schemes $u:T \ra U$ (resp. $u:\underline{T}\ra U$) which is affine as a morphism of schemes, where $\underline{T}$ is the logarithmic scheme-theoretic image of the absolute Frobenius $F_T$ of $T.$
A morphism $(U_1,\frakT_1,u_1) \ra (U_2,\frakT_2,u_2)$ in $\calE(X/\frakS)$ (resp. $\underline{\calE}(X/\frakS)$) is a pair $(f,g)$ consisting of an $\frakS$-morphism $f:\frakT_1 \ra \frakT_2$ and an $X$-morphism $g:U_1 \ra U_2$ such that $u_2\circ f_1=g\circ u_1$ (resp. $u_2 \circ \underline{f_1}=g\circ u_2$), where $f_1:T_1 \ra T_2$ is the morphism induced by $f$ by reduction modulo $p$ and $\underline{f_1}:\underline{T_1} \ra \underline{T_2}$ is induced by $f_1.$
We define, in \ref{parettop} and \ref{era4logflattop}, two topologies on each one of these two categories: the étale topology and the log flat topology and we denote by $\widetilde{\calE}(X/\frakS),$ $\widetilde{\underline{\calE}}(X/\frakS),$ $\widetilde{\calE}_{lf}(X/\frakS)$ and $\widetilde{\underline{\calE}}_{lf}(X/\frakS)$ the corresponding topoi. We equip them with the rings
$$
\Ox_{\calE(X/\frakS)}:(U,\frakT,u) \mapsto \Gamma\left (T,\Ox_T \right ),\ \Ox_{\underline{\calE}(X/\frakS)}:(U,\frakT,u) \mapsto \Gamma\left (T,\Ox_T \right ).
$$
A morphism $(f,g):(U_1,\frakT_1,u_1) \ra (U_2,\frakT_2,u_2)$ of $\widetilde{\calE}(X/\frakS)$ (resp. $\widetilde{\underline{\calE}}(X/\frakS)$) induces a morphism of ringed topoi
$$
\widetilde{f}:\left (U_{1,\text{ét}},u_{1*}\Ox_{T_1} \right ) \ra \left (U_{2,\text{ét}},u_{2*}\Ox_{T_2} \right ).
$$
We prove that the data of a module $\F$ of $\widetilde{\calE}(X/\frakS)$ (resp. $\widetilde{\underline{\calE}}(X/\frakS)$) is equivalent to the data, for every object $(U,\frakT,u),$ of a module $\F_{(U,\frakT,u)}$ of $\left (U_{\text{ét}},u_*\Ox_T \right ),$ and, for every morphism $(f,g):(U_1,\frakT_1,u_1) \ra (U_2,\frakT_2,u_2),$ a morphism
$$
c_{\F,(f,g)}:\widetilde{f}^*\F_{(U_2,\frakT_2,u_2)} \ra \F_{(U_1,\frakT_1,u_1)},
$$
satisfying certain conditions \eqref{lindescentdata}. A module $\F$ is said to be a crystal if $c_{\F,(f,g)}$ is an isomorphism for every morphism $(f,g).$ The category of crystals of $\widetilde{\calE}(X/\frakS)$ (resp. $\widetilde{\underline{\calE}}(X/\frakS)$) will be denoted $\calC(X/\frakS)$ (resp. $\underline{\calC}(X/\frakS)$). Consider the following assumption:
\begin{itemize}
\item[$(H)$] The logarithmic formal scheme $\frakS$ is equipped with a frame $\frakS \ra [P]$ and the morphism $X\ra S$ lifts to a log smooth morphism $f:(\frakX,Q) \ra (\frakS,P)$ of framed fs logarithmic $p$-adic formal schemes.
\end{itemize}
Under this assumption, we prove, in \ref{equivRQ}, that we have equivalences of categories
$$
\calC(X/\frakS) \xrightarrow{\sim} \begin{Bmatrix} \Ox_{X}\text{-modules\ with\ an}\\ \calR_1\text{-stratification} \end{Bmatrix},
$$
$$
\underline{\calC}(X/\frakS) \xrightarrow{\sim} \begin{Bmatrix} \Ox_{X}\text{-modules\ with\ a}\\ \calQ_1\text{-stratification}. \end{Bmatrix},
$$
where $\calQ_1$ and $\calR_1$ are the hopf algebras corresponding to the special fibers of $Q_{\frakX}$ and $R_{\frakX,1}$ \eqref{intro16} respectively.

For an object $(U,\frakT,u)$ of $\underline{\calE}(X/\frakS),$ we prove in \ref{Shihoparag1} that there exists a unique morphism $T \ra U''$ fitting into the commutative diagram
$$
\begin{tikzcd}
T\ar{r}{G_T} \ar{dr}\ar[bend right=30]{rdd} & \underline{T} \ar{dr}{u} \ar[hook]{r} & T\ar[bend right=-30]{dd} \\
 & U''\ar{d}\ar{r} & U\ar{d} \\
 & S \ar{r}{F_S} & S.
\end{tikzcd}
$$
where $U''$ is the base change of $U\ra S$ by $F_S$ in the category of fine logarithmic schemes and $G_T$ is given in \eqref{IranG}. We then construct, in \ref{propShiho185}, a morphism $\Delta_u:T \ra U'$ fitting into the commutative diagram
$$
\begin{tikzcd}
U \ar{rr} & & U' \ar{d} \\
\underline{T} \ar[swap]{u}{u} \ar[hook]{r} & T \ar{ur}{\Delta_{u}} \ar{r} & U'',
\end{tikzcd}
$$
We use it to define a functor
$$
\rho:\begin{array}[t]{clc}
\underline{\calE}(X/\frakS) & \ra & \calE(X'/\frakS) \\
(U,\frakT,u) & \mapsto & \left ( U',\frakT,\Delta_u \right ).
\end{array}
$$
We prove that if $X\ra S$ is of Cartier type, then $\rho$ is continuous and cocontinuous for both the étale and log flat topologies:
\begin{proposition}[\ref{thmproof2}]
If $X\ra S$ is of Cartier type, then $\rho$ is continuous and cocontinuous for both the étale and log flat topologies and hence induces morphisms of topoi
$$
C_{X/\frakS}:\widetilde{\underline{\calE}}(X/\frakS) \ra \widetilde{\calE}(X'/\frakS),
$$
$$
C_{X/\frakS,lf}:\widetilde{\underline{\calE}}_{lf}(X/\frakS) \ra \widetilde{\calE}_{lf}(X'/\frakS).
$$
such that $C_{X/\frakS}^{-1}$ and $C_{X/\frakS,lf}^{-1}$ are composition by $\rho.$
\end{proposition}
\end{parag}

\begin{parag}
Keep the assumptions of \ref{110} and suppose furthermore that $\frakS=\op{Spf}W(k)$ equipped with the trivial logarithmic structure. For an object $(U,\frakT,u)$ of $\underline{\calE}(X/\frakS),$ we have a commutative diagram
$$
\begin{tikzcd}
U \ar[swap]{d}{F_{U/S}} & & \underline{T}\ar[swap]{ll}{u} \ar[swap]{d}{F_{\underline{T}/S}} \ar[hook]{rr} & & T\ar{d}{F_{T/S}} \ar{dll}{f_{T/S}} \\
U' & & \underline{T}' \ar{ll}{u'} \ar[hook]{rr} & & T',
\end{tikzcd}
$$
where $f_{T/S}$ is defined in \eqref{Iranf}.
The functor $\rho$ is then given by
$$
\rho:\begin{array}[t]{clc}
\underline{\calE}(X/\frakS) & \ra & \calE(X'/\frakS) \\
(U,\frakT,u) & \mapsto & \left ( U',\frakT,u'\circ f_{T/S} \right ).
\end{array}
$$
We prove in this case that $C_{X/\frakS,lf}$ is an equivalence of topoi and that the inverse image functor $C_{X/\frakS}^{-1}$ induces a fully faithful functor between crystals:
\begin{theorem}[\ref{lemdirectimage}]
If $\frakS=\op{Spf}W(k)$ equipped with the trivial logarithmic structure, $S$ is its special fiber and $X$ is a log smooth fs logarithmic $S$-scheme, then the functor $C_{X/\frakS}^*$ induces a fully faithful functor
$$\calC^{\text{qcoh}}(X'/\frakS) \ra \underline{\calC}^{\text{qcoh}}(X/\frakS).$$
\end{theorem}
To get the essential surjectivity, we introduce the indexed algebras $\calA_X$ and $\calB_{X/S}$ of Lorenzon, putting an emphasis on the topos-theoretical nature of these objects.
\end{parag}

\begin{parag}
Let $(X,\Ox)$ be a ringed topos, $\calI$ a sheaf of monoids of $X,$ $p_1,p_2:\calI^2 \ra \calI$ the canonical projections and $\sigma:\calI^2 \ra \calI$ the addition map. We have a localization morphism
$$
j_{\calI}:X_{/\calI} \ra X.
$$
The morphisms $p_1,p_2$ and $\sigma$ also induce morphisms of ringed topoi, abusively denoted by
$$
p_1,p_2,\sigma: \left (X_{/\calI^2},j_{\calI^2}^*\Ox \right ) \ra \left (X_{/\calI},j_{\calI}^*\Ox \right ).
$$
An $\calI$-indexed algebra is then an $j_{\calI}^*\Ox$-module equipped with a multiplication
$$
\calA \boxtimes \calA :=p_1^*\calA \otimes_{j_{\calI^2}^*\Ox}p_2^*\calA \ra  \sigma^*\calA,
$$
satisfying certain associativity, commutativity and unity properties \eqref{inddef}. Similarly, given an $\calI$-indexed algebra $\calA,$ an $\calI$-indexed $\calA$-module is an $j_{\calI}^*\Ox$-module $\calE$ equipped with a morphism
$$
\calA \boxtimes \calE \ra \sigma^*\calE,
$$
satisfying certain properties. We then introduce indexed versions of the endomorphisms sheaf \eqref{inddef2} and the tensor product \eqref{indtensprod}. Given an $\calI$-indexed algebra $\calA$ and $\calI$-indexed $\calA$-modules $\calE$ and $\F,$ we define an $\calI$-indexed algebra
$$
\mathscr{End}_{\calA}\left (\calE,\F\right)
$$
and an $\calI$-indexed $\calA$-module
$$
\calE \circledast_{\calA} \F.
$$
\end{parag}

\begin{parag}
A logarithmic scheme $X$ gives rise to an exact sequence of abelian groups
$$
0 \ra \Ox_X^* \xrightarrow{\lambda} \calM_X^{gp} \ra \upmu:=\ov{\calM}_X^{gp} \ra 0,
$$
which, in turn, gives an invertible $\Ox_{X,\upmu}$-module
$$
\calA_X=\left (\calM^{gp} \times \Ox_{X,\upmu} \right )_{/\sim},
$$
where $(m,ab) \sim (\lambda(a)+m,b)$ for any local sections $a,$ $b$ and $m$ of $\Ox_{X,\upmu}^*,$ $\Ox_{X,\upmu}$ and $\calM_X^{gp}$ respectively. We equip this module $\calA_X$ with a structure of $\upmu$-indexed algebra. Given a morphism of logarithmic schemes $X\ra S,$ Lorenzon defines a canonical connection
$$
d_{\calA_X}:\calA_X \ra \calA_X\otimes_{\Ox_{X,\upmu}} \omega^1_{X/S,\upmu}
$$
and denotes by $\calB_{X/S}$ its kernel, which inherits a structure of $\upmu$-indexed algebra from $\calA_X$ \eqref{Bgroupe2}. He also proves that, if $X\ra S$ is log smooth, then $\calA_X$ is a locally free $\upmu$-indexed $\calB_{X/S}$-module of rank $p^d,$ where $d$ is the rank of $\omega^1_{X/S}$ \eqref{thmlor}.
\end{parag}

\begin{parag}
In \ref{krazcor158}, we recall a result by Schepler and Ohkawa about the Azumaya property of the indexed algebra $\widetilde{\calD}_{X/S}:=\calA_X \otimes_{\Ox_{X,\upmu}} \calD_{X/S,\upmu},$ where $\upmu=\ov{\calM}_X^{gp}.$ More precisely, we consider a log smooth morphism $f:(\frakX,Q) \ra (\frakS,P)$ of framed fs $p$-adic logarithmic formal schemes such that $\frakS$ is log flat and locally of finite type over $\op{Spf}W(k).$ We suppose that the exact relative Frobenius $F_1:X\ra X'$ lifts to an $(\frakS,P)$-morphism of framed fs $p$-adic logarithmic formal schemes $F:(\frakX,Q) \ra (\frakX',Q')$ such that $\frakX'$ is log smooth over $\frakS.$ If we consider $\calD_{X/S}$ as an algebra over $S^{\bullet}\calT_{X'/S}$ via the $p$-curvature map, then $\widetilde{\calD}_{X/S}$ is a $\upmu$-indexed Azumaya algebra over its center $\calB_{X/S}\otimes_{\Ox_{X',\upmu}}S^{\bullet}\calT_{X/S,\upmu}.$
Ohkawa then constructs a splitting of $\widetilde{\calD}_{X/S}$ over $\w{S}^{\bullet}\calT_{X'/S,\upmu},$ where $\w{S}^{\bullet}\calT_{X'/S}$ is the completion with respect to the ideal $\bigoplus_{n\ge 1}S^{n}\calT_{X'/S}.$ We deduce the following proposition:
\begin{proposition}[\ref{propKoko1717}]
Let $f:(\frakX,Q) \ra (\frakS,P)$ be a log smooth morphism of framed fs $p$-adic logarithmic formal schemes such that $\frakS$ is log flat and locally of finite type over $\op{Spf}W(k).$ We suppose that the exact relative Frobenius $F_1:X \ra X'$ lifts to an $(\frakS,P)$-morphism of framed fs $p$-adic logarithmic formal schemes $F:(\frakX,Q) \ra (\frakX',Q')$ such that $\frakX'$ is log smooth over $\frakS.$ Denote by $X\ra S$ its special fiber.
Let $\widehat{\Gamma}^{\bullet}\calT_{X'/S}$ be the completed PD-algebra of $\calT_{X'/S}$ and $\widehat{\Gamma}=\left (\widehat{\Gamma}^{\bullet}\calT_{X'/S} \right )_{\upmu}.$
Then we have a Morita equivalence:
$$
\Phi_{\upmu}:\begin{Bmatrix}\mathrm{Locally\ PD}\text{-}\mathrm{nilpotent}\\
\upmu\text{-}\mathrm{indexed\ }\left (\calB_{X/S}\otimes_{\Ox_{X',\upmu}}\widehat{\Gamma}\right )\text{-}\\ \mathrm{modules} \end{Bmatrix} \xrightarrow{\sim} \begin{Bmatrix} \mathrm{Locally\ PD}\text{-}\mathrm{nilpotent\ }\upmu\text{-}\mathrm{indexed}\\ \left ((\calB_{X/S}\otimes_{\Ox_{X',\upmu}}\widehat{\Gamma})\circledast_{\calB_{X/S}\otimes_{\Ox_{X',\upmu}}S}\widetilde{\calD}_{X/S} \right )\text{-} \\ \mathrm{modules}\end{Bmatrix}.
$$
\end{proposition}
\end{parag}

\begin{parag}
We now define indexed versions of Oyama topoi. In what follows, we consider the ring of Witt vectors $W(k)$ of a perfect field $k$ of characteristic $p.$ We equip $\op{Spf}W(k)$ with the trivial logarithmic structure and we consider a log flat and locally of finite type fs logarithmic $p$-adic formal scheme $\frakS$ over $\op{Spf}W(k).$ We denote by $S$ its special fiber and consider a log smooth morphism $X\ra S$ of fs logarithmic schemes which is of Cartier type. The functors
$$
u':\begin{array}[t]{clc}
\calE(X'/\frakS) & \ra & \text{ét}_{/X'} \\
(U,\frakT,u) & \mapsto & U
\end{array},\ 
\underline{u}:\begin{array}[t]{clc}
\underline{\calE}(X/\frakS) & \ra & \text{ét}_{/X} \\
(U,\frakT,u) & \mapsto & U
\end{array}
$$
induce morphisms of topoi
\begin{equation}
u':\widetilde{\calE}(X'/\frakS) \ra X'_{\text{ét}},\ \underline{u}:\widetilde{\underline{\calE}}(X/\frakS) \ra X_{\text{ét}}.
\end{equation}
We set
$$
\upmu_X=\ov{\calM}_X^{gp},\ \upmu'=u'^{-1}\upmu_X,\ \underline{\upmu}=\underline{u}^{-1}\upmu_X.
$$
The localized topoi
$$
\calE(X'/\frakS)_{/\upmu'},\ \underline{\calE}(X/\frakS)_{/\underline{\upmu}}
$$
will play the role of Oyama topoi in the indexed framework. Just as in the nonindexed case, the data of a module $\F$ of $\widetilde{\calE}(X'/\frakS)_{/\upmu'}$ (resp. $\widetilde{\underline{\calE}}(X/\frakS)_{/\underline{\upmu}}$) is equivalent to the data, for every object $(U,\frakT,u),$ of a module $\F_{(U,\frakT,u)}$ of $\left (U_{\text{ét}/\ov{\calM}_U^{gp}},\left (u_*\Ox_T \right )_{\ov{\calM}_U^{gp}} \right ),$ and, for every morphism $(f,g):(U_1,\frakT_1,u_1) \ra (U_2,\frakT_2,u_2),$ a morphism
$$
\widetilde{f}_{\upmu}^*\F_{(U_2,\frakT_2,u_2)} \ra \F_{(U_1,\frakT_1,u_1)},
$$
satisfying certain conditions \eqref{indlindescentdata2}.

A local section $s\in \Gamma(U,\upmu_X)$ naturally induces morphisms of topoi
\begin{equation}
\alpha_s':\widetilde{\calE}(U'/\frakS)  \ra  \widetilde{\calE}(X'/\frakS)_{/\upmu'},\ 
\underline{\alpha}_s:\widetilde{\underline{\calE}}(U/\frakS) \ra  \widetilde{\underline{\calE}}(X/\frakS)_{/\underline{\upmu}}.
\end{equation}
We then define a crystal of $\widetilde{\calE}(X'/\frakS)_{/\upmu'}$ (resp. $\widetilde{\underline{\calE}}(X/\frakS)_{/\underline{\upmu}}$) to be a module $\F$ of $\widetilde{\calE}(X'/\frakS)_{/\upmu'}$ (resp. $\widetilde{\underline{\calE}}(X/\frakS)_{/\underline{\upmu}}$) such that 
$$
\alpha_s'^{-1}\F\ (\text{resp.}\ \underline{\alpha}_s^{-1}\F)
$$
is a crystal.
\end{parag}

\begin{parag}
Similarly, we can define subsites $\calE_{\mathrm{strict}}(X'/\frakS)_{/\upmu'}$ and $\underline{\calE}_{\mathrm{strict}}(X/\frakS)_{/\underline{\upmu}}$ of $\calE(X'/\frakS)_{/\upmu'}$ and $\underline{\calE}(X/\frakS)_{/\underline{\upmu}}$ respectively, consisting of objects $(U,\frakT,u,s)$ such that $u$ is strict. We then prove, in \ref{equivShihoind}, that the crystals on these subsites extend naturally to crystals on $\calE_{\mathrm{strict}}(X'/\frakS)_{/\upmu'}$ and $\underline{\calE}_{\mathrm{strict}}(X/\frakS)_{/\underline{\upmu}}$ respectively:
\begin{proposition}[\ref{equivShihoind}]
Let $\frakS$ be an fs logarithmic $p$-adic formal scheme, log flat and locally of finite type over $\op{Spf}W(k).$ Let $X$ be a log smooth fs logarithmic scheme over $S.$
The restriction functor
$$
\widetilde{\calE}(X'/\frakS)_{/\upmu'} \ra \widetilde{\calE}_{\mathrm{strict}}(X'/\frakS)_{/\upmu'}
$$
$$
\left ( \text{resp.}\ \widetilde{\underline{\calE}}(X/\frakS)_{/\underline{\upmu}} \ra \widetilde{\underline{\calE}}_{\mathrm{strict}}(X/\frakS) \right )_{/\underline{\upmu}}
$$
induces an equivalence between the subcategories of crystals
$$
\begin{Bmatrix}\text{Crystals\ of\ }\widetilde{\calE}(X'/\frakS)_{/\upmu'} \end{Bmatrix} \xrightarrow{\sim} \begin{Bmatrix}\text{Crystals\ of\ }\widetilde{\calE}_{\mathrm{strict}}(X'/\frakS)_{/\upmu'} \end{Bmatrix}
$$
$$
\left (\text{resp.}\ \begin{Bmatrix}\text{Crystals\ of\ }\widetilde{\underline{\calE}}(X/\frakS)_{/\underline{\upmu}} \end{Bmatrix} \xrightarrow{\sim} \begin{Bmatrix}\text{Crystals\ of\ }\widetilde{\underline{\calE}}_{\mathrm{strict}}(X/\frakS)_{/\underline{\upmu}} \end{Bmatrix}
\right )
$$
\end{proposition}
We can define a crystal $\frakA_X$ of $\widetilde{\underline{\calE}}_{\mathrm{strict}}(X/\frakS)_{/\underline{\upmu}}$ satisfying, for any object $(U,\frakT,u)$ of $\underline{\calE}_{\mathrm{strict}}(X/\frakS),$
$$
\frakA_{X,(U,\frakT,u)}=\calA_T.
$$
This crystal then extends to a crystal $\frakA_X$ of $\widetilde{\underline{\calE}}(X/\frakS)_{/\underline{\upmu}}.$
Consider the assumption
\begin{itemize}
\item[$(H)$] The morphism $X\ra S$ lifts to a log smooth morphism $f:(\frakX,Q) \ra (\frakS,P)$ of framed fs logarithmic $p$-adic formal schemes and $X'$ lifts to an fs framed logarithmic $p$-adic formal scheme $(\frakX',Q')$ over $(\frakS,P),$ such that $\frakX' \ra \frakS$ is log smooth.
\end{itemize}
We then prove, in \ref{equivRQind} and under the hypothesis $(H),$ that there exists a canonical equivalence of categories between the category of crystals of $\widetilde{\underline{\calE}}(X/\frakS)_{/\underline{\upmu}}$ and the category of $\Ox_{X,\upmu_X}$-modules equipped with a $\calQ_{1,\upmu_X}$-stratification. We also prove that, under $(H),$ there exists a canonical equivalence of categories between the category of crystals of $\widetilde{\calE}(X'/\frakS)_{/\upmu'}$ and the category of $\Ox_{X',\upmu_X}$-modules equipped with an $\calR_{1,\upmu_X}'$-stratification. In addition, these equivalences preserve structures of $\upmu_X$-indexed algebras and $\upmu_X$-indexed modules.
Let $q_1,q_2:Q_{\frakX} \ra \frakX$ and $q_1',q_2':R_{\frakX',1} \ra \frakX'$ be the canonical projections. The projections $q_1$ and $q_2$ are strict and so they induce, by (\cite{Lor2000} I 3.1 page 273), isomorphisms of $\upmu_X$-indexed algebras
$$c_i:q_i^*\calA_X \xrightarrow{\sim} \calA_{Q_1}.$$
We get a $\calQ_{1,\upmu_X}$-stratification
\begin{equation}\label{intro2}
\epsilon_{\calA_X,Q}:q_2^*\calA_X \xrightarrow{c_2} \calA_{Q_1} \xrightarrow{c_1^{-1}} q_1^*\calA_X.
\end{equation}
Under the hypothesis $(H),$ the crystal $\frakA_X$ previously defined, corresponds to $\left (\calA_X,\epsilon_{\calA_X,Q} \right ).$
Now consider the $\upmu_X$-indexed $\Ox_{X',\upmu_X}$-module $\calB_{X/S}.$ The isomorphism that exchanges factors $q_2'^*\calB_{X/S} \xrightarrow{\sim} q_1'^*\calB_{X/S}$ is an $\calR_{1,\upmu_X}'$-stratification on $\calB_{X/S}$ and so it corresponds to a crystal $\frakB_X$ of $\widetilde{\calE}(X'/\frakS)_{/\upmu'}.$ We check in \ref{propBglue} that the definition of this crystal $\frakB_X$ doesn't depend on the local assumption $(H)$ and can hence be defined globally as a crystal of $\widetilde{\calE}(X'/\frakS)_{/\upmu'}$ without the need of the assumption $(H).$ The crystals $\frakA_X$ and $\frakB_X$ are naturally equipped with structures of indexed algebras \eqref{equivRQind}.
We will be interested in crystals of $\upmu'$-indexed $\frakB_X$-modules and $\underline{\upmu}$-indexed $\frakA_X$-modules. Under the hypothesis $(H),$ these crystals correspond to modules with \emph{admissible} stratifications:
\begin{enumerate}
\item Let $\calE$ be a $\upmu_X$-indexed $\calA_{X}$-module equipped with a $\calQ_{1,\upmu_X}$-stratification
$$\epsilon:\calQ_{1,\upmu_X}\otimes_{\Ox_{X,\upmu_X}}\calE \ra \calE\otimes_{\Ox_{X,\upmu_X}}\calQ_{1,\upmu_X}.$$
The stratification $\epsilon$ is said to be \emph{admissible} if it is compatible with $\epsilon_{\calA_X,Q}$ \eqref{intro2}.
\item A $\upmu_X$-indexed $\calB_{X/S}$-module $\calE'$ equipped with an $\calR_{1,\upmu_X}'$-stratification
$$\epsilon':\calR_{1,\upmu_X}'\otimes_{\Ox_{X',\upmu_X}}\calE' \ra \calE'\otimes_{\Ox_{X',\upmu_X}}\calR_{1,\upmu_X}'$$
is said to be \emph{admissible} if $\epsilon'$ is $\calB_{X/S}$-linear.
\end{enumerate}
We prove, in \ref{equivRQAB}, that the canonical equivalence of categories between the category of crystals of $\widetilde{\underline{\calE}}(X/\frakS)_{/\underline{\upmu}}$ (resp. $\widetilde{\calE}(X'/\frakS)_{/\upmu'}$) and the category of $\Ox_{X,\upmu_X}$-modules equipped with a $\calQ_{1,\upmu_X}$-stratification (resp. $\Ox_{X',\upmu_X}$-modules equipped with an $\calR_{1,\upmu_X}'$-stratification), induces a canonical equivalence of categories between the full subcategory of crystals of $\upmu$-indexed $\frakA_X$-modules of $\widetilde{\underline{\calE}}(X/\frakS)_{/\underline{\upmu}}$ (resp. crystals of $\upmu'$-indexed $\frakB_X$-modules of $\widetilde{\calE}(X'/\frakS)_{/\upmu'}$) and the full subcategory of $\upmu_X$-indexed $\calA_X$-modules equipped with an admissible $\calQ_{1,\upmu_X}$-stratification (resp. $\upmu_X$-indexed $\calB_{X/S}$-modules equipped with an admissible $\calR_{1,\upmu_X}'$-stratification).
\end{parag}

\begin{parag}
In addition to the hypothesis $(H),$ suppose that the exact relative Frobenius $F_1:X\ra X'$ lifts to an $(\frakS,P)$-morphism of framed logarithmic formal schemes $F:(\frakX,Q) \ra (\frakX',Q'),$ such that $\frakX'\ra \frakS$ is log smooth. Consider the morphism of Hopf algebras $\calV:\calR_1' \ra \calQ_1$ induced by $\nu:Q_{\frakX}\ra R_{\frakX',1}$ \eqref{Muzannu}. In \ref{propkraz1934}, we construct a commutative diagram
\begin{equation}
\begin{tikzcd}
\begin{Bmatrix}\text{Crystals\ of\ }\widetilde{\calE}(X'/\frakS)_{/\upmu'} \end{Bmatrix} \ar{r}{C^{-1}_{X/\frakS,\upmu}} \ar[swap, sloped]{d}{\sim} & \begin{Bmatrix}\text{Crystals\ of\ }\widetilde{\underline{\calE}}(X/\frakS)_{/\underline{\upmu}} \end{Bmatrix}  \ar[sloped]{d}{\sim} \\
\begin{Bmatrix}\Ox_{X',\upmu_X}\text{-modules\ with\ an}\\ \calR_{1,\upmu_X}'\text{-stratification} \end{Bmatrix} \ar{r}{\Psi_{\upmu}} & \begin{Bmatrix}\Ox_{X,\upmu_X}\text{-modules\ with\ a}\\ \calQ_{1,\upmu_X}\text{-stratification} \end{Bmatrix},
\end{tikzcd}
\end{equation}
where the vertical equivalences were given above, $C_{X/\frakS,\upmu}$ is induced by the morphism of topoi $C_{X/\frakS}$ and $\Psi_{\upmu}(\calE',\epsilon') = \left (F_1^*\calE',\calV^*\epsilon'\right )$ for an $\Ox_{X',\upmu_X}$-module $\calE'$ equipped with an $\calR_{1,\upmu_X}'$-stratification $\epsilon'.$
Finally, we prove, always under the assumption of the lifting of the Frobenius, the following theorem:
\begin{theorem}[\ref{THMXX}]
Let $\frakS$ be an fs logarithmic $p$-adic formal scheme which is log flat and locally of finite type over $\op{Spf}W(k)$ and denote its special fiber by $S.$ Let $X$ be an fs logarithmic scheme which is log smooth and of Cartier type over $S.$ Suppose that $\frakS$ is equipped with a frame $\frakS \ra [P]$ and that the exact relative Frobenius $F_1:X\ra X'$ lifts to an $(\frakS,P)$-morphism of framed logarithmic formal schemes $F:(\frakX,Q) \ra (\frakX',Q'),$ such that $\frakX\ra \frakS$ and $\frakX'\ra \frakS$ are log smooth. Then, there exists a commutative diagram
\begin{equation}
\begin{tikzcd}
\begin{Bmatrix}\mathrm{Crystals\ of\ }\upmu'\text{-}\mathrm{indexed}\\ \frakB_X\text{-}\mathrm{modules\ of\ }\widetilde{\calE}(X'/\frakS)_{/\upmu'} \end{Bmatrix} \ar{r}{C_{X/\frakS,\upmu}^{-1}} \ar[swap,sloped]{d}{\sim} & \begin{Bmatrix}\mathrm{Crystals\ of\ }\underline{\upmu}\text{-}\mathrm{indexed}\\ \frakA_X\text{-}\mathrm{modules\ of\ }\widetilde{\underline{\calE}}(X/\frakS)_{/\underline{\upmu}} \end{Bmatrix} \ar[sloped]{d}{\sim}\\
\begin{Bmatrix}\upmu_X\text{-}\mathrm{indexed\ }\calB_{X/S}\text{-}\mathrm{modules\ with}\\ \mathrm{an\ admissible}\\ \calR_{1,\upmu_X}'\text{-}\mathrm{stratification} \end{Bmatrix} \ar[swap,sloped]{d}{\sim} \ar{r} & \begin{Bmatrix}\upmu_X\text{-}\mathrm{indexed\ }\calA_X\text{-}\mathrm{modules\ with}\\ \mathrm{an\ admissible}\\Q_{1,\upmu_X}\text{-}\mathrm{stratification} \end{Bmatrix} \ar[sloped]{d}{\sim} \\
\begin{Bmatrix}\mathrm{Locally\ PD}\text{-}\mathrm{nilpotent\ } \upmu_X\text{-}\mathrm{indexed}\\\left (\left (\w{\Gamma}^{\bullet}\calT_{X'/S}\right )_{\upmu_X}\otimes_{\Ox_{X',\upmu_X}}\calB_{X/S} \right )\text{-}\mathrm{modules} \end{Bmatrix} \ar{r}{\Phi_{\upmu}} & \begin{Bmatrix}\mathrm{Locally\ PD}\text{-}\mathrm{nilpotent\ }\upmu_X\text{-}\mathrm{indexed}\\\widetilde{\calD}_{X/S}^{\gamma}\text{-}\mathrm{modules} \end{Bmatrix}.
\end{tikzcd}
\end{equation}
Since $\Phi_{\upmu}$ is an equivalence of categories, the functor $C_{X/\frakS,\upmu}^{-1}$ is also an equivalence of categories.
\end{theorem}
Finally, we deduce an equivalence of categories between crystals of $\frakB_X$-modules and crystals of $\frakA_X$-modules:
\begin{theorem}[\ref{THMXXYY}]
Let $\frakS$ be an fs logarithmic $p$-adic formal scheme which is log flat and locally of finite type over $\op{Spf}W(k)$ and denote its special fiber by $S.$ Let $X$ be an fs logarithmic scheme which is log smooth and of Cartier type over $S.$ Then, we have an equivalence of categories
$$
\begin{Bmatrix}\mathrm{Crystals\ of\ }\upmu'\text{-}\mathrm{indexed}\\ \frakB_X\text{-}\mathrm{modules\ of\ }\widetilde{\calE}(X'/\frakS)_{/\upmu'} \end{Bmatrix}  \ra \begin{Bmatrix}\mathrm{Crystals\ of\ }\underline{\upmu}\text{-}\mathrm{indexed}\\ \frakA_X\text{-}\mathrm{modules\ of\ }\widetilde{\underline{\calE}}(X/\frakS)_{/\underline{\upmu}} \end{Bmatrix}.
$$
\end{theorem}
Note that, in this last theorem, we do not require any lifting assumption neither on the exact relative Frobenius nor on $X.$
\end{parag}

\textbf{Acknowledgements.} This manuscript is the result of my PhD thesis, prepared at Université Paris-Saclay and IHÉS. I would like to express my gratitude to my PhD advisor Ahmed Abbes for introducing me to this topic and for his guidance and patience during this journey. 
I also would like to thank Atsushi Shiho and Daxin Xu for their helpful comments and suggestions.

\newpage

\section{Notations and conventions}

\begin{parag}
In this manuscript, we fix a universe $\mathbb{U}$ and a prime number $p.$
\end{parag}

\begin{parag}\label{Not9}
For multi-indices $I=(I_1,\hdots,I_d),J=(J_1,\hdots,J_d)\in \N^d,$ we denote by $I!$ the product
$$I!=\prod_{i=1}I_i!$$
and by $\begin{pmatrix}I \\ J \end{pmatrix}$ the binomial coefficient
$$
\begin{pmatrix}I \\ J \end{pmatrix}=\prod_{i=1}^d\begin{pmatrix}I_i \\ J_i \end{pmatrix}.
$$
\end{parag}

\begin{parag}\label{Not8}
If $X$ is a scheme or a formal scheme, we denote by $|X|$ its underlying topological space.
\end{parag}

\begin{parag}
If $\F$ is a sheaf of sets on a site $X,$ the notation $x\in \F$ means that there exists an object $U$ of $X$ such that $x\in \F(U).$
\end{parag}

\begin{parag}
For morphisms of fs logarithmic schemes $X\ra S$ and $Y\ra S,$ we denote by $X\times_S^{\op{log}}Y$ their fiber product in the category of fs logarithmic schemes.
\end{parag}

\begin{parag}\label{Not2}
If $A$ is a ring, $a\in A$ and $M$ is an $A$-module, we denote by $M_{/a\text{-tor}}$ the quotient of $M$ by the submodule of $a$-torsion of $M$ i.e. the submodule consisting of elements $x\in M$ such that $a^kx=0$ for some positive integer $k.$
\end{parag}

\begin{parag}\label{Not1}
The logarithmic structures considered in this manuscript are all defined on the small étale site.  
For a morphism of logarithmic schemes $X\ra S,$ we denote by $\omega^1_{X/S}$ the sheaf of logarithmic differentials of $X$ over $S$ of order $1$ (\cite{Ogus2018} IV 1.2.4). For any positive integer $i,$ we set $\omega^i_{X/S}=\Lambda^i\omega^1_{X/S}.$
\end{parag}

\begin{parag}\label{Not3}
For a logarithmic scheme $T$ of positive characteristic $p,$ we denote by $F_T$ the absolute Frobenius morphism of $T,$ i.e. the morphism given by the identity on the underlying topological spaces, $\Ox_T\ra \Ox_T,\ x\mapsto x^p$ and $\calM_T\ra \calM_T,\ m\mapsto pm.$
\end{parag}

\begin{parag}\label{Not4}
Let $R$ be a ring and $P$ a monoid, we denote by $R[P]$ the free $R$-algebra on $P$ and by $A_R[P]$ the scheme $\Sp R[P]$ equipped with the logarithmic structure associated with the canonical chart $P\ra R[P].$ Let $p$ be a prime number and $n\ge 1$ an integer. We denote by $A_n[P]$ (resp. $A[P]$) the logarithmic scheme $A_{\Z/p^n\Z}[P]$ (resp. $A_{\Z}[P]$).
\end{parag}

\begin{parag}\label{Not5}
If $M$ is a monoid, we denote by $M^{\times}$ the subgroup of $M$ consisting of invertible elements, by $\ov{M}$ the quotient $M/M^{\times}$ and by $M^{gp}$ the group associated to $M.$ We also denote by $M^{int}$ the image of $M$ by the canonical morphism $M\ra M^{gp}.$ Similarly, if $\calM$ is a sheaf of monoids, we denote by $\calM^{\times}$ the subsheaf of invertible sections, by $\ov{\calM}$ the quotient sheaf $\calM/\calM^{\times}$ i.e. the sheaf associated to the presheaf $U\mapsto \ov{\calM(U)}=\calM(U)/\calM(U)^{\times},$ by $\calM^{gp}$ the sheaf associated to the presheaf $U\mapsto \calM(U)^{gp}$ and by $\calM^{int}$ the sheaf associated to the presheaf $U\mapsto \calM(U)^{int}.$
\end{parag}

\begin{parag}
A morphism of monoids (resp. sheaves of monoids) is said to be surjective if it is so as a map of sets (resp. morphism of sheaves of sets).
\end{parag}

\begin{parag}\label{Not7}
We say that a morphism $u:M\ra N$ of monoids is \emph{strict} if the induced morphism $\ov{u}:\ov{M} \ra \ov{N}$ is an isomorphism.
\end{parag}

\begin{parag}
For a scheme $X,$ we denote by $\text{ét}_{/X}$ the small étale site on $X$ and by $X_{\text{ét}}$ the small étale topos on $X.$
\end{parag}

\begin{parag}\label{Not6}
Following \href{https://stacks.math.columbia.edu/tag/03DL}{Definition 03DL}, for a ringed site $(C,\Ox),$ we say that an $\Ox$-module $\calE$ is quasi-coherent if, for every object $U$ of $C,$ there exists a covering $(U_i\ra U)$ such that, for every $i,$ the restricted module $\calE_{|U_i}$ admits an exaact sequence of $\Ox_{|U_i}$-modules of the form:
$$\bigoplus_{j\in J}\Ox_{|U_i} \ra \bigoplus_{k\in K} \Ox_{|U_i} \ra \calE_{|U_i} \ra 0.$$
For any scheme $X,$ we denote by $\op{QCoh}(X_{\text{zar}})$ (resp. $\op{QCoh}(X_{\text{ét}})$) the category of quasi-coherent $\Ox_X$-modules on the small Zariski site of $X$ (resp. the category of quasi-coherent $\Ox_X$-modules on the small étale site of $X$).
Recall, by \href{https://stacks.math.columbia.edu/tag/03DX}{Proposition 03DX}, that $\op{QCoh}(X_{\text{zar}})$ and $\op{QCoh}(X_{\text{ét}})$ are canonically equivalent.
\end{parag}

\newpage

\section{Frobenuiserie}

\begin{parag}\label{PFrob}
Let $f:X\ra S$ be a morphism of fine logarithmic schemes of characteristic $p.$ Let $X''$ be the base change, in the category of fine logarithmic schemes, of $X$ by the absolute Frobenius morphism $F_S:S\ra S$ of $S$ (\ref{Not3}). Let $F_{X/S}$ be the relative Frobenius morphism of $X$ with respect to $S$ i.e. the unique morphism $X\ra X''$ making the following diagram commutative
\begin{equation}\label{diag481}
\begin{tikzcd}
X\ar{dr}{F_{X/S}}\ar[bend right=-30]{drr}{F_X}\ar[swap,bend right=30]{ddr}{f} & & \\
 & X''\ar{r}\ar{d} & X\ar{d}{f} \\
 & S\ar{r}{F_S} & S
\end{tikzcd}
\end{equation}
Recall that $F_{X/S}$ is \emph{weakly inseparable} (\cite{Ogus2018} III 2.4), so $F_{X/S}$ factors uniquely as a composition of a log étale morphism $G:X'\ra X''$ and an inseparable morphism $F:X\ra X'$ (\cite{Ogus2018} IV 3.3.8): 
\begin{equation}\label{diag51}
\begin{tikzcd}
X\ar{r}{F}\ar{dr}\ar[bend right=30]{rdd} & X'\ar{d}{G}\ar{dr}{\pi} & \\
 & X''\ar{d}\ar{r} & X\ar{d}{f} \\
 & S \ar{r}{F_S} & S
\end{tikzcd}
\end{equation}
The morphism $F:X\ra X'$ is called the \emph{exact relative Frobenius of $X$ with respect to $S$} (\cite{Ogus2018} IV 3.3.9).

In the rest of this manuscript, we introduce the following notations: if $P$ is a monoid, we denote by
$$F_P:P \ra P,\ x\mapsto px$$
the Frobenius morphism of $P.$ If $\theta:P\ra Q$ is a morphism of fine monoids, we set $Q''=(Q\oplus_{P,F_P}P)^{int}$ and denote by $Q'$ the inverse image of $Q$ by
\begin{equation}\label{eqfrobmon1}
Q''^{gp} \ra Q^{gp},\ (x,y) \mapsto px+\theta^{gp}(y).
\end{equation}
Denote by $v:Q' \ra Q$ the morphism induced by \eqref{eqfrobmon1}. 
If $\theta:P\ra Q$ is a chart for $f,$ then, by the proof of (\cite{Kat89} 4.10), the canonical morphism $X''\ra A_1[Q'']$ is a chart, $X'=X''\times_{A_1[Q'']}A_1[Q']$ and $v:Q' \ra Q$ is a chart of the exact relative Frobenius $F:X\ra X'.$
Note that if $Q$ is saturated, then $Q'$ is clearly saturated. It follows that if $X$ and $S$ are fs, then so is $X'.$
\end{parag}

\begin{proposition}\label{piunivhomeo}
The morphism $\pi:X' \ra X$ \eqref{diag51} is a universal homeomorphism.
\end{proposition}

\begin{proof}
The Frobenius $F_S$ is weakly inseparable (\cite{Ogus2018} III 2.4) hence so is the canonical projection $X'' \ra X$ (\cite{Ogus2018} 2.4.8). The relative Frobenius $X\ra X''$ and the exact relative Frobenius $X\ra X'$ are also weakly inseparable so they are universal homeomorphisms. Hence so is the canonical morphism $X' \ra X''$ by (\cite{SP} \href{https://stacks.math.columbia.edu/tag/0H2M}{Lemma 0H2M}). It follows that $\pi:X' \ra X$ is a universal homeomorphism.
\end{proof}

\begin{definition}\label{defkummer}~ 
\begin{enumerate}
\item A morphism $u:M\ra N$ of fs monoids is said to be \emph{Kummer} if it is injective and for any $y\in N,$ there exists a positive integer $n$ such that $ny\in u(M).$
\item A morphism $f:X\ra Y$ of fs logarithmic schemes is said to be \emph{of Kummer type} if, for any geometric point $\ov{x}$ of $X$ and $\ov{y}=f(\ov{x}),$ the morphism of monoids
$$f^{\flat}_{\ov{x}}:\ov{\calM}_{Y,\ov{y}} \ra \ov{\calM}_{X,\ov{x}}$$
is Kummer.
\end{enumerate}
\end{definition}

\begin{proposition}\label{propfrob12}
Let $\theta:P\ra Q$ a morphism of fine monoids such that $\theta^{gp}$ is injective and the torsion subgroup of $\op{coker}\theta^{gp}$ is finite of order coprime with $p.$ Let $Q'$ and $v:Q' \ra Q$ be as defined in \ref{PFrob}. Then $v$ is injective.
\end{proposition}

\begin{proof}
Let $(x,y),(x',y')\in Q'$ such that $$v(x,y)=v(x',y').$$
Then
$$px+\theta^{gp}(y)=px'+\theta^{gp}(y')\in Q$$
and
$$p(x-x')=\theta^{gp}(y'-y).$$
It follows that $p(x-x')=0$ in $\op{coker}\theta^{gp}.$ And since the torsion subgroup of $\op{coker}\theta^{gp}$ has a finite order coprime with $p,$ we deduce that $x=x'+\theta^{gp}(t)$ for a certain element $t\in P^{gp}.$ Then
$$\theta^{gp}(y+pt)=\theta^{gp}(y)+p(x-x')=\theta^{gp}(y').$$
Since $\theta^{gp}$ is injective, we get $y+pt=y'$ and so, in $Q',$
$$(x,y)=(x'+\theta^{gp}(t),y)=(x',y+pt)=(x',y').$$ 
\end{proof}

\begin{corollaire}\label{propfrob13}
Let $\theta:P\ra Q$ a morphism of fs monoids such that $\theta^{gp}$ is injective and the torsion subgroup of $\op{coker}\theta^{gp}$ is finite of order coprime with $p.$ Let $Q'$ and $v:Q' \ra Q$ be as defined in \ref{PFrob}. Then $v$ is Kummer \eqref{defkummer}.
\end{corollaire}

\begin{proof}
The injectivity of $v$ follows from \ref{propfrob12} and by definition, for every $x\in Q,$ $(x,0)\in Q'$ and $v(x,0)=px.$
\end{proof}

\begin{proposition}\label{FKummer}
Let $f:X\ra S$ be a morphism of fs logarithmic schemes of characteristic $p$ and $F:X\ra X'$ the exact relative Frobenius. Then $F$ is of Kummer type \eqref{defkummer}.
\end{proposition}

\begin{proof}
Let $\ov{x}\ra X$ be a geometric point and $\ov{x}'=F(\ov{x}).$ The morphism $\ov{F}_{\ov{x}}^{\flat}:\ov{\calM}_{X',\ov{x}'} \ra \ov{\calM}_{X,\ov{x}}$ is by definition exact. By (\cite{Ogus2018} I 4.3.10), to prove that it is Kummer, it is sufficient to prove that it is small i.e. for all $z\in \ov{\calM}_{X,\ov{x}}^{gp}$ there exists a positive integer $n$ and $y\in \ov{\calM}_{X',\ov{x}'}^{gp}$ such that $nz=\left (\ov{F}_{\ov{x}}^{\flat}\right )^{gp} \left (y \right )$ (\cite{Ogus2018} I 4.3.1). This is immediate since, for every local section $m$ of $\calM_{X},$ we have $pm=F^{\flat}\left (\pi^{\flat}(m) \right ),$ where $\pi:X'\ra X$ is given in \eqref{diag51}.
\end{proof}

\begin{parag}
In the non logarithmic case, if $X\ra S$ is a smooth morphism of schemes of characteristic $p,$ then the relative Frobenius morphism of $X$ with respect to $S$ is flat. We end this section by proving that a similar result holds in the logarithmic case: we recall the definition of \emph{log flatness} and prove that if $X\ra S$ is a log smooth morphism of fine logarithmic schemes of characteristic $p,$ then the exact relative Frobenius morphism of $X$ with respect to $S$ is log flat. 
\end{parag}

\begin{definition}[\cite{Ogus2018} IV 4.1.1]\label{loggfflat}
We say that a morphism of fine logarithmic schemes $f:X\ra S$ is \emph{log flat} if, fppf locally on $X$ and $S,$ there exists a chart $\theta:P\ra Q$ of $f$ such that:
\begin{enumerate}
\item $\theta$ is an injective morphism of fine monoids.
\item The morphism $X \ra S\times_{A[P]}A[Q],$ induced by $f$ and $X\ra A[Q],$ is flat as a morphism of schemes.
\end{enumerate} 
\end{definition}

\begin{theorem}\label{thmlogflat}
Let $f:X\ra S$ be a log smooth morphism of fine logarithmic schemes of characteristic $p.$ Then the exact relative Frobenius $F:X\ra X'$ \eqref{diag51} is log flat.
\end{theorem}

\begin{proof}
By (\cite{Kat89} 3.5), there exists, étale locally on $X$ and $S,$ a chart $\theta:P \ra Q$ of $f$ such that $\theta^{gp}$ is injective, the torsion subgroup of $\op{coker}\theta^{gp}$ has a finite order coprime with $p$ and the morphism $g:X \ra S\times_{A_1[P]}A_1[Q]$ induced by $f$ and $A_1[\theta]$ is smooth as a morphism of schemes. Set $T=S\times_{A_1[P]}A_1[Q]$ and $F_T:T\ra T$ the absolute Frobenius morphism of $T.$ The morphism $g:X\ra T$ is smooth so the relative Frobenius
$$F_{X/T}:X\ra X\times_{T,F_T}T$$
is flat. Let $X''$ be the fiber product of $f:X\ra S$ and $F_S:S\ra S$ in the category of fine logarithmic schemes,
$$F_P:P \ra P,\ x\mapsto px,$$
$Q''=(Q\oplus_{P,F_P}P)^{int}$ and $Q'$ the inverse image of $Q$ by
$$v:Q''^{gp} \ra Q^{gp},\ (x,y) \mapsto px+\theta^{gp}(y).$$
By the proof of (\cite{Kat89} 4.10), the canonical morphism $X''\ra A_1[Q'']$ is a chart, $X'=X''\times_{A_1[Q'']}A_1[Q']$ and $v:Q' \ra Q$ is a chart of the exact relative Frobenius $F:X\ra X'.$ By the lemma \ref{lemX333} below, there exists a canonical isomorphism $X'\times_{A_1[Q']}A_1[Q] \xrightarrow{\sim} X\times_{T,F_T}T$ and the morphism $G:X\ra X'\times_{A_1[Q']}A_1[Q],$ induced by $F:X\ra X'$ and $X\ra A_1[Q],$ identifies with $F_{X/T}.$ So $G$ is flat. It remains to prove that the chart $v:Q'\ra Q$ is injective, which we did in \ref{propfrob12}.
\end{proof}

\begin{corollaire}\label{corlogflat}
Let $f:X\ra S$ be a log smooth morphism of fine logarithmic locally Noetherian schemes of characteristic $p$ and $F:X\ra X'$ the exact relative Frobenius \eqref{diag51}. If $F^{\flat}:F^{-1}\calM_X \ra \calM_{X'}$ is integral, then the underlying morphism of schemes of $F$ is faithfully flat.
\end{corollaire}

\begin{proof}
This is an immediate consequence of \ref{thmlogflat} and \cite{Ogus2018} IV 4.3.5.
\end{proof}

\begin{lemma}\label{lemX333}
Keep the same hypothesis and notation of \ref{thmlogflat} and its proof. Then there exists a canonical isomorphism
$$X'\times_{A_1[Q']}A_1[Q]\xrightarrow{\sim}X\times_{T,F_T}T.$$
\end{lemma}

\begin{proof}
Clearly,
$$X'\times_{A_1[Q']}A_1[Q]=X''\times_{A_1[Q'']}A_1[Q].$$
Consider the cartesian diagram
$$
\begin{tikzcd}
T\ar{r}{q_1}\ar[swap]{d}{q_2} & S\ar{d}{\beta} \\
A_1[Q] \ar{r}{A_1[\theta]} & A_1[P]
\end{tikzcd}
$$
and the following commutative diagram
$$
\begin{tikzcd}
X\times_{T,F_T}T \ar{r}\ar[swap]{d}  & X''\times_{A_1[Q'']}A_1[Q] \ar{d} \ar{r} & X'' \ar{d} \ar{r} & X\ar{d} \\
T\ar{r} \ar[swap,bend right=20]{rrr}{F_T} & T''\times_{A_1[Q'']}A_1[Q] \ar{r} & T''\ar{r} & T 
\end{tikzcd}
$$
The big and two right squares are cartesian, so the left square is also cartesian and it is thus sufficient to prove that the morphism
$$\varphi:T\ra T''\times_{A_1[Q'']}A_1[Q],$$
defined by the relative Frobenius $F_{T/S}:T\ra T''$ and $q_2:T\ra A_1[Q],$ is an isomorphism. For that, we exhibit an inverse morphism.
Consider the two cartesian diagrams
$$
\begin{tikzcd}
T'' \ar{r}{p_1}\ar[swap]{d}{p_2} & T\ar{d}{q_1} & & T''\times_{A_1[Q'']}A_1[Q] \ar{r}{\rho_1} \ar[swap]{d}{\rho_2} & T''\ar{d}\\
S\ar{r}{F_S} & S & & A_1[Q] \ar{r} & A_1[Q'']
\end{tikzcd}
$$
The morphism $\varphi$ fits into the commutative diagram
$$
\begin{tikzcd}
T \ar[bend right=-20]{rrd}{F_{T/S}} \ar[swap,bend right=30]{ddr}{q_2} \ar{dr}{\varphi} &  & \\
 & T''\times_{A_1[Q'']}A_1[Q]\ar{r}{\rho_1}\ar{d}{\rho_2} & T''\ar{d} \\
 & A_1[Q]\ar{r} & A_1[Q'']
\end{tikzcd}
$$
Consider the morphism
$$\psi:T''\times_{A_1[Q'']}A_1[Q] \ra T$$
defined in the diagram
$$
\begin{tikzcd}
T''\times_{A_1[Q'']}A_1[Q] \ar{r}{\rho_1} \ar[swap,bend right=30]{ddr}{\rho_2} \ar{dr}{\psi} & T''\ar{dr}{p_2} & \\
 & T\ar{r}{q_1}\ar{d}{q_2} & S\ar{d}{\beta} \\
 & A_1[Q]\ar{r}{A_1[\theta]} & A_1[P]
\end{tikzcd}
$$
The fact that $\psi\circ \varphi=\op{Id}_T$ follows immediately from the definitions of $\varphi$ and $\psi.$ Let us check that $\varphi\circ \psi=\op{Id}_{T''\times_{A_1[Q'']}A_1[Q]}.$ By definition, we have
$$\rho_2\circ \varphi\circ \psi=q_2\circ \psi=\rho_2.$$
Now, we just need to check that
$$\rho_1\circ \varphi \circ \psi=\rho_1.$$
First, we have
$$p_2\circ \rho_1 \circ \varphi \circ \psi=p_2\circ F_{T/S}\circ \psi=q_1\circ \psi=p_2\circ \rho_1.$$
It is thus sufficient to check that
$$p_1\circ \rho_1\circ \varphi\circ \psi=p_1\circ \rho_1,$$
which is equivalent to the two following equalities
\begin{alignat}{2}
q_1\circ p_1\circ \rho_1\circ \varphi\circ \psi=q_1\circ p_1\circ \rho_1 \label{Ja1} \\
q_2\circ p_1\circ \rho_1\circ \varphi\circ \psi=q_2\circ p_1\circ \rho_1. \label{Ja2}
\end{alignat}
The equality (\ref{Ja1}) follows from
\begin{alignat*}{2}
q_1\circ p_1\circ \rho_1\circ \varphi\circ \psi &= F_S\circ p_2\circ \rho_1\circ \varphi\circ \psi \\
&= F_S\circ p_2\circ \rho_1 \\
&= q_1 \circ p_1 \circ \rho_1.
\end{alignat*}
Now, consider the morphism
$$F_Q:Q\ra Q,\ x\mapsto px.$$
It fits into the following commutative diagram
$$
\begin{tikzcd}
T''\times_{A_1[Q'']}A_1[Q] \ar{r}{\rho_1} \ar[swap]{d}{\rho_2} & T''\ar{d}{p_2}\ar{r}{p_1} & T\ar{d}{q_2} \\
A_1[Q] \ar{r} \ar[swap,bend right=20]{rr}{A_1[F_Q]} & A_1[Q''] \ar{r} & A_1[Q] 
\end{tikzcd}
$$
the equality (\ref{Ja2}) follows then from
\begin{alignat*}{2}
q_2 \circ p_1 \circ \rho_1 \circ \varphi \circ \psi &= q_2 \circ p_1 \circ F_{T/S} \circ \psi \\
&= q_2 \circ F_T \circ \psi \\
&= A_1[F_Q] \circ q_2 \circ \psi \\
&= A_1[F_Q] \circ \rho_2 \\
&= q_2 \circ p_1 \circ \rho_1.
\end{alignat*}
This finishes the proof.
\end{proof}

\begin{parag}\label{dxuflat}
Let $k$ be a perfect field of characteristic $p,$ $W(k)$ its ring of Witt vectors and $\frakS=\op{Spf}W(k).$ Following (\cite{DXU19} 2.5), we say that a $p$-adic formal $\frakS$-scheme $\frakX$ is \emph{flat over $\frakS$} or that $\frakX$ is a \emph{flat formal $\frakS$-scheme} if multiplication by $p$ on $\Ox_{\frakX}$ is injective. This is equivalent to the fact that, for every affine open formal subscheme $U$ of $\frakX,$ the algebra $\Gamma(U,\Ox_{\frakX})$ is flat over $W(k).$ In this sense, a $p$-adic formal $\frakS$-scheme $\frakX$ is flat over $\frakS$ if and only if $\frakX_n$ is flat over $\frakS_n$ for all positive integers $n,$ where $\frakX_n$ and $\frakS_n$ are obtained from $\frakX$ and $\frakS$ respectively by reduction modulo $p^n.$ Indeed, let $A$ be a $p$-adic $W(k)$-algebra and suppose that, for all integers $n\ge 1,$ $A_n:=A/p^nA$ is flat over $W_n(k):=W(k)/p^nW(k).$ Let $x\in A$ such that $px=0.$ By the flatness of $A_n$ over $W_n(k),$ the image of $x$ in $A_n$ belongs to $p^{n-1}A/p^nA.$ Since $A$ is separated, $x=0$ and so $A$ is flat over $W(k).$ 
\end{parag}

\begin{proposition}\label{era3proplogstr}
Let $f:X\ra Y$ and $g:Z\ra Y$ be morphisms of logarithmic schemes such that $g$ is strict. Suppose that there exists a morphism of schemes $h:X \ra Z$ such that $f=g\circ h$ in the category of schemes. Then there exists a unique morphism of sheaves of monoids $h^{\flat}:h^{-1}\calM_Z \ra \calM_X$ making $h$ a morphism of logarithmic schemes and such that $f=g\circ h$ in the category of logarithmic schemes. The result remains true if we replace schemes by formal schemes (see section 7 of the present document).
\end{proposition}

\begin{proof}
Since $g$ is strict, the canonical morphism $g^*\calM_Y \ra \calM_Z,$ where $g^*\calM_Y$ is the logarithmic structure pullback of $\calM_Y,$ is an isomorphism. Then
$$h^*\calM_Z=f^*\calM_Y.$$
The data of a morphism of prelogarithmic structures
$$h^{-1}\calM_Z \ra \calM_X$$
is equivalent to the data of a morphism of logarithmic structures
$$h^*\calM_Z \ra \calM_X.$$
Since $h^*\calM_Z=f^*\calM_Y,$ the only morphism satisfying the desired condition is
$$f^{\flat}:f^*\calM_Y \ra \calM_X.$$
\end{proof}

\begin{corollaire}\label{propstrcart}~ 
\begin{enumerate}
\item Let
$$
\begin{tikzcd}
Z\ar{r}{q} \ar[swap]{d}{p} & Y\ar{d}{g} \\
X\ar{r}{f} & S
\end{tikzcd}
$$
be a commutative diagram of logarithmic schemes. If it is cartesian in the category of schemes and $p$ and $g$ are strict, then it is cartesian in the category of logarithmic schemes.
\item Let $g:Y \ra S,$ $p:Z\ra X$ and $f:X\ra S$ be morphisms of logarithmic schemes  and $q:Z\ra Y$ a morphism of schemes such that $p$ and $g$ are strict and the diagram
\begin{equation}\label{era4cardiag1}
\begin{tikzcd}
Z\ar{r}{q} \ar[swap]{d}{p} & Y\ar{d}{g} \\
X\ar{r}{f} & S
\end{tikzcd}
\end{equation}
is commutative in the category of schemes. Then there exists a morphism of sheaves of monoids $q^{\flat}:q^{-1}\calM_Y \ra \calM_Z$ such that \eqref{era4cardiag1} is commutative in the category of logarithmic schemes. In addition, if \eqref{era4cardiag1} is cartesian in the category of schemes, then it is cartesian in the category of logarithmic schemes.
\end{enumerate}
These two assertions remain true if we consider formal schemes instead of schemes (see section 7 of the present document). 
\end{corollaire}

\begin{proof}
We only prove the result for schemes as the case of formal schemes is similar.
\begin{enumerate}
\item We start by proving the first assertion.
Let $h_1:T \ra X$ and $h_2:T\ra Y$ be morphisms of logarithmic schemes such that $f\circ h_1=g\circ h_2.$ There exists a unique morphism of schemes $h:T \ra Z$ such that $p\circ h=h_1$ and $q\circ h=h_2.$ By \ref{era3proplogstr}, there exists a morphism $h^{\flat}:h^{-1}\calM_Z \ra \calM_T$ making $h$ a morphism of logarithmic schemes such that $p\circ h=h_1.$ To prove that $q\circ h=h_2$ in the category of logarithmic schemes, it is sufficient to prove that the diagram
$$
\begin{tikzcd}
h_2^*\calM_Y \ar{r} \ar{dr} & h^*\calM_Z \ar{d} \\
& \calM_T
\end{tikzcd}
$$
is commutative.
This follows from the commutativity of the diagram
$$
\begin{tikzcd}
h_2^*g^*\calM_S \ar{r} \ar{dr} & h_2^*\calM_Y \ar{r}  & h^*\calM_Z \ar{d} \\
& h_2^*\calM_Y \ar{r} & \calM_T
\end{tikzcd}
$$
and the fact that $g^*\calM_S \ra \calM_Y$ is an isomorphism.
\item For the second assertion, $g$ is strict so $\calM_Y=g^*\calM_S$ and we define $q^{\flat}$ as the composition
$$q^{\flat}:q^{-1}\calM_Y \ra q^*\calM_Y = q^*g^*\calM_S=p^*f^*\calM_S\xrightarrow{p^*f^{\flat}}p^*\calM_X \xrightarrow{p^{\flat}} \calM_Z,$$
where the first arrow is the canonical one. The diagram (\ref{era4cardiag1}) is then commutative in the category of logarithmic schemes. We conclude by applying the first assertion.
\end{enumerate}
\end{proof}

\begin{lemma}\label{Mfs}
For a monoid $M,$ let $\sim$ be the equivalence relation defined, for any $x,y\in M,$ by $x\sim y$ if $px=py.$ We equip the quotient $\underline{M}:=M_{/\sim}$ with the natural monoid structure inherited from $M$ and call it the monoidal Frobenius image of $M.$
\begin{enumerate}
\item For any monoid $M,$ there exists a canonical isomorphism
$$\underline{M^{gp}} \xrightarrow{\sim} \underline{M}^{gp}.$$
\item If $M$ is an integral (resp. fine, resp. saturated, resp. fs) monoid then so is $\underline{M}.$
\item If $M$ is a saturated monoid then the canonical morphism $M \ra \underline{M}$ is exact and strict.
\end{enumerate}
\end{lemma}

\begin{proof}
The composition
$$M \ra \underline{M} \ra \underline{M}^{gp}$$
factors through $\underline{M^{gp}}$ yielding a morphism
$$u:\underline{M^{gp}} \ra \underline{M}^{gp}.$$
The composition
$$M \ra M^{gp} \ra \underline{M^{gp}}$$
factors through $\underline{M}^{gp}$ yielding a morphism
$$v:\underline{M}^{gp} \ra \underline{M^{gp}}.$$
The morphisms $u$ and $v$ are inverse to each other.

It is clear that, if $M$ is finitely generated, then so is $\underline{M}.$ Suppose that $M$ is integral and let $x,y,t\in M$ such that
$$\ov{x}+\ov{t}=\ov{y}+\ov{t}\in \underline{M}.$$
Then
$$px+pt=py+pt \in M$$
and by integrality of $M,$ $px=py$ so $\ov{x}=\ov{y} \in \underline{M}.$
Now suppose that $M$ is saturated and let $x\in M^{gp}$ and $\ov{x}$ its image in $\underline{M}^{gp}.$ If there exists a positive integer $n$ such that $n\ov{x}\in \underline{M}$ then there exists $y\in M$ such that $pnx=py\in M^{gp}.$ Since $M$ is saturated and $y\in M,$ $x\in M$ and so $\ov{x} \in \underline{M}.$ In addition, if the image of an element $x\in M^{gp}$ in $\underline{M}^{gp}$ belongs to $\underline{M},$ then $px\in M$ and so $x\in M.$ This proves that the square
$$
\begin{tikzcd}
M \ar{r} \ar{d} & \underline{M} \ar{d} \\
M^{gp} \ar{r} & \underline{M}^{gp}
\end{tikzcd}
$$ 
is cartesian, hence the exactness of $M \ra \underline{M}.$

Finally, suppose that $M$ is saturated and let us prove that the canonical morphism $M\ra \underline{M}$ is strict. It is surjective so it is sufficient to prove that, if $m,m'\in M$ are such that there exists a $u\in M$ satisfying $pm=pm'+pu$ and $\ov{u}$ is invertible in $\underline{M},$ then there exists $v\in M^{\times}$ such that $m=m'+v.$ Since $\ov{u}$ is invertible in $\underline{M},$ there exists $u'\in M$ such that $pu+pu'=0.$ So $u\in M^{\times}.$ In $M^{gp},$ $p(m-m')=pu\in M$ and $M$ is saturated so $m-m'\in M.$ Similarly $p(m'-m)=-pu\in M$ so $m'-m\in M,$ hence $m-m'\in M^{\times}.$ We conclude by taking $v=m-m'.$
\end{proof}

\begin{lemma}\label{Mfs2}
Let $u:M\ra N$ be a morphism of saturated monoids, $\underline{M}$ and $\underline{N}$ the monoidal Frobenius images of $M$ and $N$ respectively \eqref{Mfs} and $\underline{u}:\underline{M}\ra \underline{N}$ the morphism induced by $u.$ If $u$ is strict, then so is $\underline{u}.$
\end{lemma}

\begin{proof}
We have to prove that the morphism
$$\underline{M}/\underline{M}^{\times} \ra \underline{N}/\underline{N}^{\times},$$
induced by $\underline{u},$ is an isomorphism.
First, note that an element $x\in M$ is invertible if and only if its image $\ov{x}$ in $\underline{M}$ is invertible. This, along with the surjectivity of the projection $N\ra \underline{N},$ yields the surjectivity of $\underline{u}$ modulo the invertibles. Let $x,y\in M$ and $t\in N$ such that $\ov{t}$ is invertible and
$$\underline{u}(\ov{x})=\underline{u}(\ov{y})+\ov{t}.$$
Then
$$\ov{u(x)}=\ov{u(y)}+\ov{t}.$$
So
$$u(px)=u(py)+pt.$$
Since $u$ is strict, there exists $z\in M^{\times}$ such that $px=py+z.$ Set $w=x-y\in M^{gp}.$ Since $pw=z\in M^{\times}$ and $M$ is saturated, $w\in M^{\times}.$ We have $px=p(y+w)$ so $\ov{x}=\ov{y}+\ov{w}$ and $\underline{u}$ is injective modulo the invertibles.
\end{proof}

\begin{parag}\label{era3logimagedef}
Let $X$ be a logarithmic scheme of characteristic $p$ and $\underline{X}$ the scheme theoretic image of the absolute Frobenius $F_X$ i.e. the closed subscheme of $X$ defined by the ideal consisting of local sections of $\Ox_X$ whose $p^{th}$ power vanishes. The ideal of the closed immersion $\underline{X} \ra X$ is, by definition, a nilideal and $\underline{X} \ra X$ is hence a universal homeomorphism (\href{https://stacks.math.columbia.edu/tag/054M}{Lemma 054M}). We identify the small étale sites of $\underline{X}$ and $X$ via this immersion (\href{https://stacks.math.columbia.edu/tag/04DZ}{Theorem 04DZ}). For any étale $X$-scheme $U,$ let $\underline{\calM_X(U)}$ be the monoidal Frobenius image of the monoid $\calM_X(U)$ (\ref{Mfs}).
We denote by $\calM_{\underline{X}}$ the sheaf of monoids associated to the presheaf
$$U \mapsto \underline{\calM_X(U)}$$
and call it \emph{the monoidal Frobenius image of $\calM_X$}.
The canonical morphism $\calM_X \ra \calM_{\underline{X}}$ is an epimorphism of sheaves of sets.
For any local sections $m$ and $m'$ of $\calM_X,$ we have
$$\left (\alpha_X(m)-\alpha_X(m') \right )^p=\alpha_X(pm)-\alpha_X(pm').$$
This proves that the composition
$$\calM_X\xrightarrow{\alpha_X} \Ox_X \ra \Ox_{\underline{X}},$$
where the second arrow is the canonical projection, induces a morphism of monoids
$$\alpha_{\underline{X}} :\calM_{\underline{X}} \ra \Ox_{\underline{X}}.$$
We claim that the morphism of monoids
\begin{equation}\label{era4eq1}
\alpha_{\underline{X}}^{-1}\left (\Ox_{\underline{X}}^* \right ) \ra \Ox_{\underline{X}}^*,
\end{equation}
induced by $\alpha_{\underline{X}},$ is an isomorphism and hence that $(\calM_{\underline{X}},\alpha_{\underline{X}})$ is a logarithmic structure on $\underline{X}.$ Indeed, let $m$ and $m'$ be local sections of $\calM_X$ and $\ov{m}$ and $\ov{m'}$ their images in $\calM_{\underline{X}}.$ If $\alpha_{\underline{X}}(\ov{m})=\alpha_{\underline{X}}(\ov{m'})$ then $(\alpha_X(m)-\alpha_X(m'))^p=0$ and hence $\alpha_X(pm)=\alpha_X(pm').$ If, in addition, $\alpha_{\underline{X}}(\ov{m}),\alpha_{\underline{X}}(\ov{m'})\in \Ox_{\underline{X}}^*$ and since a local section of $\Ox_X$ is invertible if and only if its image in $\Ox_{\underline{X}}$ is invertible, then $\alpha_X(m),\alpha_X(m') \in \Ox_X^*$ and so $pm=pm'.$ This proves the injectivity of (\ref{era4eq1}). The surjectivity follows from the fact that a local section of $\Ox_X$ is invertible if and only if its image in $\Ox_{\underline{X}}$ is invertible.

We call the logarithmic scheme $(\underline{X},\calM_{\underline{X}}),$ the \emph{logarithmic scheme theoretic image of $F_X$}.
\end{parag}

\begin{proposition}\label{era4prop15}
Let $X$ be a logarithmic scheme of characteristic $p.$
\begin{enumerate}
\item The absolute Frobenius morphism $F_X:X\ra X$ factors uniquely through the logarithmic scheme $\underline{X}$ defined in \ref{era3logimagedef}. If we denote by $G_X:X \ra \underline{X}$ the resulting morphism, then the composition
$$
\underline{X} \ra X \xrightarrow{G_X} \underline{X},
$$
of the canonical immersion $\underline{X} \ra X$ with $G_X,$ is equal to the absolute Frobenius $F_{\underline{X}}$ of $\underline{X}.$
\item If $X$ is fs then the canonical immersion $\underline{X} \ra X$ is strict and inseparable (\cite{Ogus2018} III 2.4.2).
\end{enumerate}
\end{proposition}

\begin{proof}
The first assertion follows from the fact that $F_X^{\#}:\Ox_X \ra \Ox_X,\ x\mapsto x^p$ factors through $\Ox_{\underline{X}}$ and $F_X^{\flat}:\calM_X \ra \calM_X,\ m\mapsto pm$ factors through $\calM_{\underline{X}}$ and the commutativity of the diagram
$$
\begin{tikzcd}
\calM_X \ar{r}\ar[bend right=-30]{rrr}{m \mapsto pm} \ar[swap]{d}{\alpha_X} & \calM_{\underline{X}} \ar{rr}{\ov{m}\mapsto pm} \ar[swap]{d}{\alpha_{\underline{X}}} & & \calM_X \ar{d}{\alpha_X} \\
\Ox_X \ar{r} \ar[swap, bend right=30]{rrr}{x \mapsto x^p} & \Ox_{\underline{X}} \ar{rr}{\ov{x}\mapsto x^p} & & \Ox_X 
\end{tikzcd}
$$
The fact that the composition
$$
\underline{X} \ra X \xrightarrow{G_X} \underline{X}
$$
is equal to $F_{\underline{X}}$ is clear. 
Now suppose that $X$ is fs. The canonical morphism $\ov{\calM}_{X,\ov{x}} \ra \ov{\calM}_{\underline{X},\ov{x}}$ is an isomorphism for every geometric point $\ov{x}$ of $X,$ by \ref{Mfs}. It follows that the canonical immersion $\underline{X} \ra X$ is strict and hence exact. Since the canonical immersion $\underline{X} \ra X$ is clearly weakly inseparable (\cite{Ogus2018} III 2.4.2), this finishes the proof. 
\end{proof}

\begin{proposition}\label{era3functoriality}~ 
\begin{enumerate}
\item The correspondance $X\mapsto \underline{X}$ is functorial on logarithmic schemes of characteristic $p.$
\item Let $X\ra S$ and $Y\ra S$ be morphisms of logarithmic schemes of characteristic $p.$ There exists a canonical morphism
$$\underline{X\times_SY}\ra \underline{X}\times_{\underline{S}}\underline{Y}.$$
In addition, this morphism is a closed nilimmersion i.e. a closed immersion whose defining ideal is a nilideal.
\item Let $X\ra S$ and $Y\ra S$ be morphisms of fs logarithmic schemes of characteristic $p.$ There exists a canonical morphism
$$\underline{X\times_S^{\op{log}}Y}\ra \underline{X}\times_{\underline{S}}^{\op{log}}\underline{Y}.$$
In addition, this morphism is a closed nilimmersion.
\item Let $f:X\ra Y$ be a morphism of fs logarithmic schemes of characteristic $p$ and $\underline{f}:\underline{X}\ra \underline{Y}$ the morphism induced by $f.$ If $f$ is strict, then so is $\underline{f}.$
\end{enumerate}
\end{proposition}

\begin{proof}
Let $f:X\ra Y$ be a morphism of logarithmic schemes of characteristic $p.$ By \ref{era3logimagedef}, $\underline{X}$ (resp. $\underline{Y}$) is the closed subscheme of $X$ (resp. $Y$) defined by the ideal $\calI_X=\{x\in \Ox_X, x^p=0\}$ (resp. $\calI_Y=\{ y\in \Ox_Y, y^p=0\}$). Let $i_X:\underline{X}\ra X$ and $i_Y:\underline{Y}\ra Y$ be the canonical closed immersions. For any $y\in \calI_Y,$ $f^{\#}(y)^p=0,$ so there exists a morphism of schemes $\underline{f}:\underline{X}\ra \underline{Y}$ making the following diagram commutative in the category of schemes
$$
\begin{tikzcd}
\underline{X}\ar{r}{\underline{f}} \ar{d}{i_X} & \underline{Y}\ar{d}{i_Y} \\
X \ar{r}{f} & Y 
\end{tikzcd}
$$
Note that, as maps of topological spaces, $f=\underline{f}.$
By definition of $\calM_{\underline{X}}$ and $\calM_{\underline{Y}},$ the morphism $f^{\flat}:\calM_Y \ra f_*\calM_X$ induces a morphism
$$\underline{f}^{\flat}:\calM_{\underline{Y}} \ra \underline{f}_*\calM_{\underline{X}}.$$
This turns $\underline{f}$ into a morphism of logarithmic schemes. The correspondances $X\mapsto \underline{X}$ and $f\mapsto \underline{f}$ clearly define a functor. By functoriality, we get the canonical morphisms $\underline{X\times_SY}\ra \underline{X}\times_{\underline{S}}\underline{Y}$ and $\underline{X\times_S^{\op{log}}Y}\ra \underline{X}\times_{\underline{S}}^{\op{log}}\underline{Y}.$ The fact that they are closed nilimmersions follows from the commutativity of the diagrams
$$
\begin{tikzcd}
\underline{X\times_SY} \ar{r} \ar{d} & \underline{X}\times_{\underline{S}}\underline{Y} \ar{dl} & \underline{X\times_S^{\op{log}}Y} \ar{r} \ar{d} & \underline{X}\times_{\underline{S}}^{\op{log}}\underline{Y} \ar{dl} \\
X\times_SY & & X\times^{\op{log}}_SY, & 
\end{tikzcd}
$$
the fact that $\underline{X\times_SY} \ra X\times_SY$ and $\underline{X\times_S^{\op{log}}Y}\ra X\times_S^{\op{log}}Y$ are closed nilimmersions and $\underline{X}\times_{\underline{S}}\underline{Y} \ra X\times_SY,$ $\underline{X}\times_{\underline{S}}^{\op{log}}\underline{Y} \ra X\times^{\op{log}}_SY$ are affine morphisms, and (\cite{EGA1} 4.3.6). The assertion on strictness follows from \ref{Mfs2}.
\end{proof}

\begin{remark}\label{rem178}
Let $k$ be a perfect field of characteristic $p,$ $S=\Sp k$ equipped with the trivial logarithmic structure and $X$ a fine logarithmic scheme over $S.$ We have $\underline{S}=S.$ The canonical morphism
$$\underline{X}'=\underline{X}\times_{S,F_S}S \ra X\times_{S,F_S}S=X'$$
factors through $\underline{X'}.$ The resulting morphism $\underline{X}' \ra \underline{X'}$ is inverse to the canonical morphism $\underline{X'}=\underline{X\times_{S,F_S}S} \ra \underline{X}\times_{\underline{S},F_{\underline{S}}}\underline{S}= \underline{X}'$ (\ref{era3functoriality}).
From now on, we will identify $\underline{X}'$ and $\underline{X'}$ via these isomorphisms.
\end{remark}

\begin{proposition}\label{propfT/S}
Let $k$ be a perfect field of characteristic $p,$ $S=\Sp k$ equipped with the trivial logarithmic structure and $f:X\ra S$ a morphism of fs logarithmic schemes. Then the exact relative Frobenius $F:X\ra X'$ factors uniquely through the canonical closed immersion $i_X':\underline{X}' \ra X'$ \eqref{rem178}. We denote by $f_{X/S}:X\ra \underline{X}'$ the unique morphism satisfying $F=i_{X}'\circ f_{X/S}.$ This morphism $f_{X/S}$ is inseparable (\cite{Ogus2018} III 2.4.2) and an epimorphism in the category of logarithmic schemes.
\end{proposition}

\begin{proof}
The morphism $f$ is saturated so $F$ is equal to the relative Frobenius morphism $F_{X/S}$ (\cite{Ogus2018} III 2.5.4).
In addition, the logarithmic structure on $S$ is trivial so $F_S$ is strict and hence $X'=X\times_{S,F_S}S.$ Let $\underline{X}'=\underline{X}\times_{S,F_S}S.$
By \ref{era4prop15}, the absolute Frobenius $F_X$ factors uniquely through the canonical closed immersion $i_X:\underline{X} \ra X.$ Hence the existence of a morphism of logarithmic schemes $f_{X/S}:X\ra \underline{X}'$ fitting into the following commutative diagram
$$
\begin{tikzcd}
X \ar[swap, bend right=-30]{rrd} \ar{dr}{f_{X/S}} \ar[bend right =30]{ddr}{F} \ar[bend right=30]{dddr} &  & \\
 & \underline{X}' \ar{r} \ar{d}{i_X'} & \underline{X} \ar{d}{i_X} \\
 & X' \ar{r} \ar{d} & X\ar{d}{f} \\
 & S\ar{r}{F_S} & S
\end{tikzcd}
$$
The uniqueness follows from the fact that $i_X':\underline{X}' \ra X'$ is a closed immersion and closed immersions of logarithmic schemes are monomorphisms. The (exact) relative Frobenius morphism $F_{X/S}$ and $i_{X}'$ are inseparable (\ref{era4prop15}) so, by (\cite{Ogus2018} III 2.4.8 (3)), $f_{X/S}$ is inseparable.
We now prove that $f_{X/S}$ is an epimorphism. If $X$ is a affine, say $X=\Sp A$ for a $k$-algebra $A,$ then $\underline{X}=\Sp (A/I),$ where $I$ is the ideal of $A$ consisting of elements of $A$ whose $p^{th}$ power vanishes. Then $\underline{X}'=\Sp \left ( A/I)\otimes_{k,F_k}k\right ),$ where $F_k$ is the absolute Frobenius morphism of $k,$ and $f_{X/S}$ corresponds to
$$
g:\begin{array}[t]{clc}
(A/I)\otimes_{k,F_k}k & \ra & A \\
\ov{x} \otimes a & \mapsto & ax^p,
\end{array}
$$
for $x\in A$ and $a\in k.$ Since $k$ is perfect, $g$ is injective and hence the underlying morphism of schemes of $f_{X/S}$ is an epimorphism in the category of schemes. It remains to prove that
$$f_{X/S}^{\flat}:f_{X/S}^{-1}\calM_{\underline{X}'} \ra \calM_X$$
is injective. We identify the small étale sites of $X,$ $X'$ and $\underline{X}'$ via the universal homeomorphisms $X' \ra X$ and $\underline{X}'\ra X'.$ Let $\ov{x}$ be a geometric point of $X$ and $\ov{s}=f(\ov{x}).$
Then $\calM_{\underline{X}',\ov{x}}=\calM_{\underline{X},\ov{x}}\oplus_{\calM_{S,\ov{s}},F_{S,\ov{s}}^{\flat}}\calM_{S,\ov{s}}$ and $f_{X/S,\ov{x}}$ is given by
$$
f_{X/S,\ov{x}}:\begin{array}[t]{clc}
\calM_{\underline{X},\ov{x}}\oplus_{\calM_{S,\ov{s}},F_{S,\ov{s}}^{\flat}}\calM_{S,\ov{s}} & \ra & \calM_{X,\ov{x}} \\
(\ov{m},t) & \mapsto & pm+f_{\ov{x}}^{\flat}(t),
\end{array}
$$
for $m\in \calM_{X,\ov{x}}$ and $t\in\calM_{S,\ov{s}}.$ Let $m,m' \in \calM_{X,\ov{x}}$ and $t,t\in \calM_{S,\ov{s}}$ such that $f_{X/S,\ov{x}}^{\flat}(\ov{m},t)=f_{X/S,\ov{x}}^{\flat}(\ov{m'},t').$ Then
$$pm+f_{\ov{x}}^{\flat}(t)=pm'+f_{\ov{x}}^{\flat}(t').$$
By definition of $\calM_{\underline{X}}$ and amalgamated sums in the category of monoids, to prove that $(\ov{m},t)=(\ov{m'},t')$ in $\calM_{\underline{X},\ov{x}}\oplus_{\calM_{S,\ov{s}},F_{S,\ov{s}}^{\flat}}\calM_{S,\ov{s}},$ it is sufficient to prove that there exist $u,u'\in \calM_{S,\ov{s}}$ such that
\begin{equation}\label{era4eq4}
\begin{cases}
pm+f_{\ov{x}}^{\flat}(pu)=pm'+f_{\ov{x}}^{\flat}(pu') \\
t+pu'=t'+pu.
\end{cases}
\end{equation}
Since $S=\Sp k$ is equipped with the trivial logarithmic structure, $\calM_S=\Ox_S^*.$ Since $k$ is perfect, there exist $u,u'\in \calM_{S,\ov{s}}$ such that $t=pu$ and $t'=pu'.$ The condition (\ref{era4eq4}) is then automatically satisfied.
\end{proof}

\begin{proposition}\label{Xreduced}
Let $k$ be a perfect field of characteristic $p,$ $S=\Sp k$ equipped with the trivial logarithmic structure and $f:X\ra S$ a smooth morphism of fs logarithmic schemes. Then $X$ is reduced, hence $\underline{X}=X.$
\end{proposition}

\begin{proof}
Let $\ov{x}$ be a geometric point of $X$ and $\ov{s}=f(\ov{x}).$ Note that since $S$ is equipped with the trivial logarithmic structure, $\calM_{S,\ov{s}}=\Ox_{S,\ov{s}}^*.$ The morphism
$$
f_{\ov{x}}^{\flat}:\calM_{S,\ov{s}} \ra \calM_{X,\ov{x}}
$$
induces, modulo the invertibles, the morphism
$$
\ov{f}^{\flat}_{\ov{x}}: \{1\} \ra \ov{\calM}_{X,\ov{x}}.
$$
Since the direct sum of saturated morphisms is saturated, the morphism $\ov{f}^{\flat}_{\ov{x}}$ is saturated. It follows that $f$ is saturated (\cite{Ogus2018} III 2.5.1). The result then follows from (\cite{Tsuji} II.4.2).
\end{proof}

\begin{parag}\label{Shihoparag1}
Let $S$ be an fs logarithmic scheme of characteristic $p,$ $U$ and $T$ two fs logarithmic schemes over $S$ and $u:\underline{T} \ra U$ an $S$-morphism. Recall the diagram defining the exact relative Frobenius
\begin{equation}\label{diag51U}
\begin{tikzcd}
U\ar{r}\ar{dr}\ar[bend right=30]{rdd} & U'\ar{d}\ar{dr} & \\
 & U''\ar{d}\ar{r} & U\ar{d} \\
 & S \ar{r}{F_S} & S,
\end{tikzcd}
\end{equation}
where the square is cartesian in the category of fine logarithmic schemes and $U' \ra U''$ is log étale. By \ref{era4prop15}, the absolute Frobenius $F_T:T \ra T$ factors as
$$
\begin{tikzcd}
T \ar{r}{G_T} \ar[swap]{dr}{F_T} & \underline{T} \ar{d}{i} \\
 & T,
\end{tikzcd}
$$
where $i:\underline{T} \ra T$ is the canonical immersion, which is strict by \ref{era4prop15}. We then have a unique morphism $T \ra U''$ fitting into the commutative diagram
\begin{equation}\label{morphismShiho2}
\begin{tikzcd}
T\ar{r}{G_T} \ar{dr}\ar[bend right=30]{rdd} & \underline{T} \ar{dr}{u} \ar{r}{i} & T\ar[bend right=-30]{dd} \\
 & U''\ar{d}\ar{r} & U\ar{d} \\
 & S \ar{r}{F_S} & S.
\end{tikzcd}
\end{equation}
\end{parag}

\begin{proposition}\label{propShiho185}
Keep the notation of \ref{Shihoparag1}. The diagram
\begin{equation}\label{morphismShiho1}
\begin{tikzcd}
U \ar{rr} & & U' \ar{d} \\
\underline{T} \ar[swap]{u}{u} \ar{r}{i} & T \ar[dashed]{ur}{\Delta_{u}} \ar{r} & U'',
\end{tikzcd}
\end{equation}
where $T \ra U''$ is defined in \eqref{morphismShiho2}, is commutative and there exists a unique morphism $\Delta_u:T\ra U'$ fitting into it.
In addition, we have the following properties:
\begin{enumerate}
\item If $i_X:\underline{X} \ra X$ is the canonical immersion then $\Delta_{i_X}:X\ra X'$ is the exact relative Frobenius.
\item If $u$ is affine, then so is $\Delta_u.$
\item If $S=\op{Spec}k$ equipped with the trivial logarithmic structure and $f_{T/S}$ is the morphism defined in \ref{propfT/S}, then $U'=U''$ and $\Delta_u=u' \circ f_{T/S}.$
\end{enumerate}
\end{proposition}

\begin{proof}
The first point is clear.

The composition
$$
\underline{T} \xrightarrow{i} T \xrightarrow{G_T} \underline{T}
$$
is equal to $F_{\underline{T}}.$
On one hand, by \eqref{diag51U}, the composition
$$
\underline{T} \xrightarrow{u} U \ra U' \ra U'' \ra U,
$$
where $U'' \ra U$ is the canonical projection, is equal to
$$
\underline{T} \xrightarrow{u} U\xrightarrow{F_U} U.
$$
This is equal to
$$
\underline{T} \xrightarrow{F_{\underline{T}}} \underline{T} \xrightarrow{u} U.
$$
This is then equal to
$$
\underline{T} \xrightarrow{i} T \xrightarrow{G_T} \underline{T} \xrightarrow{u} U.
$$
By \eqref{morphismShiho2}, this is equal to
$$
\underline{T} \xrightarrow{i} T \ra U'' \ra U.
$$
On the other hand, by \eqref{diag51U}, the composition
$$
\underline{T} \xrightarrow{u} U \ra U' \ra U'' \ra S,
$$
where $U'' \ra S$ is the canonical projection, is equal to
$$
\underline{T} \xrightarrow{u} U\ra S.
$$
Since $u$ is an $S$-morphism, this is equal to
$$
\underline{T} \ra S.
$$
This is, by definition, equal to
$$
\underline{T} \xrightarrow{i} T \ra S.
$$
By definition of $T \ra U'',$ this is equal to
$$
\underline{T} \xrightarrow{i} T \ra U'' \ra S.
$$
Since $U''$ is, by definition, the fiber product of $U\ra S$ and $F_S$ in the category of fine logarithmic schemes, we get the commutativity of the outer diagram of
\begin{equation}
\begin{tikzcd}
U \ar{rr} & & U' \ar{d} \\
\underline{T} \ar[swap]{u}{u} \ar{r}{i} & T \ar[dashed]{ur}{\Delta_{u}} \ar{r} & U''.
\end{tikzcd}
\end{equation}
Then, the existence and uniqueness of $\Delta_{u}$ follows from the fact that $U' \ra U''$ is log étale and $\underline{T} \ra T$ is a strict nilpotent closed immersion. The fact that $\Delta_{u}$ is affine if so is $u,$ is clear.

Suppose now that $S=\op{Spec}k$ equipped with the trivial logarithmic structure. It follows, since $U$ is fs, that the morphism $U\ra S$ is saturated and so $U'=U''$ (\cite{Ogus2018} III 2.5.4). We have the diagram
$$
\begin{tikzcd}
T \ar{dr}{G_T} \ar[swap,bend right=70]{dd}{F_{T/S}} \ar[bend right=-70]{ddr}{F_T} \ar[swap]{d}{f_{T/S}} & \\
\underline{T}' \ar{d}{i'} \ar{r} & \underline{T} \ar{d}{i} \\
T' \ar{r} & T,
\end{tikzcd}
$$
where $T' \ra T$ and $\underline{T}' \ra \underline{T}$ are the canonical projections.
The square and the outer diagram are commutative. Since $i$ is a strict immersion hence a monomorphism, it follows that the diagram
$$
\begin{tikzcd}
T \ar{dr}{G_T} \ar[swap]{d}{f_{T/S}} & \\
\underline{T}' \ar{r} & \underline{T}
\end{tikzcd}
$$
is commutative. It follows that the diagram
$$
\begin{tikzcd}
T \ar{dr}{f_{T/S}} \ar[bend right=-30]{drr}{G_T} \ar[bend right=30]{dddrr} & & & \\
 & \underline{T}' \ar{r} \ar{dr}{u'} & \underline{T} \ar{dr}{u} & \\
 & & U' \ar{r} \ar{d} & U \ar{d} \\
 & & S \ar{r}{F_S} & S
\end{tikzcd}
$$
is commutative.
The morphism $\Delta_{u}$ is then equal to the composition
$$
T \xrightarrow{f_{T/S}} \underline{T}' \xrightarrow{u'} U'.
$$
\end{proof}

\begin{proposition}
Let $S$ be an fs logarithmic scheme of characteristic $p$ and $T_1\ra T_2$ a morphism of fs logarithmic schemes over $S$ and $u_1:\underline{T_1} \ra U_1$ and $u_2:\underline{T_2} \ra U_2$ two $S$-morphisms fitting into a commutative diagram
$$
\begin{tikzcd}
\underline{T_1} \ar{r} \ar[swap]{d}{u_1} & \underline{T_2} \ar{d}{u_2} \\
U_1 \ar{r} & U_2.
\end{tikzcd}
$$
The diagram
$$
\begin{tikzcd}
T_1 \ar{r} \ar[swap]{d}{\Delta_{u_1}} & T_2\ar{d}{\Delta_{u_2}} \\
U_1' \ar{r} & U_2'
\end{tikzcd}
$$
is then commutative.
\end{proposition}

\begin{proof}
Consider the notation of \ref{Shihoparag1}. We have the commutative diagram
\begin{equation}\label{diagShiho8}
\begin{tikzcd}
T_1 \ar{rr}{G_{T_1}} \ar{dd} \ar{dr} & & \underline{T_1} \ar[dashed]{dd} \ar{dr}{u_1} & \\
 & U_1'' \ar{rr} \ar{dd} & & U_1 \ar{dd} \\
 T_2 \ar[bend right=30]{ddr} \ar[dashed,pos=0.20]{rr}{G_{T_2}} \ar{dr} & & \underline{T_2} \ar[dashed]{dr}{u_2} & \\
& U_2'' \ar{rr} \ar{d} & & U_2 \ar{d} \\
 & S \ar{rr}{F_S} & & S.
\end{tikzcd}
\end{equation}
We also have the commutative diagrams
\begin{equation}\label{diagShiho6}
\begin{tikzcd}
 & U_2 \ar[bend right =-30]{rrd} &  & \\
U_1 \ar{rr} \ar{ur} & & U_1' \ar{d} \ar{r} & U_2' \ar{d} \\
\underline{T_1} \ar{r} \ar{u}{u_1} & T_1 \ar{ur}{\Delta_{u_1}} \ar{r} & U_1'' \ar{r} & U_2''
\end{tikzcd}
\end{equation}
and
\begin{equation}\label{diagShiho7}
\begin{tikzcd}
U_1 \ar{r}  & U_2 \ar{rr} &  & U_2' \ar{d} \\
\underline{T_1} \ar{r} \ar{u}{u_1} & \underline{T_2} \ar{u}{u_2} \ar{r} & T_2 \ar{ur}{\Delta_{u_2}} \ar{r}  & U_2''.
\end{tikzcd}
\end{equation}
The lower composition
$$
\underline{T_1} \ra \underline{T_2} \ra T_2 \ra U_2''
$$
of \eqref{diagShiho7}, is equal to
$$
\underline{T_1} \ra T_1 \ra T_2 \ra U_2''.
$$
This is equal, by \eqref{diagShiho8}, to
$$
\underline{T_1} \ra T_1 \ra U_1'' \ra U_2'',
$$
which is the lower composition of \eqref{diagShiho6}.
We conclude by the log étaleness of $U_2' \ra U_2''.$
\end{proof}

\newpage

\section{Logarithmic differential operators}

Let $S$ be a scheme of positive characteristic $p$ equipped with a fine logarithmic structure and a PD structure $\gamma$ and let $X\ra S$ be a smooth morphism of fine logarithmic schemes. We suppose that $\gamma$ extends to $X.$

\begin{parag}\label{P11}
We denote, for any integers $n,k \ge 0,$ by $P_{X/S}(n)$ the logarithmic PD-envelope of the diagonal immersion $X\ra X^{n+1},$ by $\ov{\calI}(n)$ the PD-ideal of $P_{X/S}(n)$ and by $P_{X/S}^k(n)$ the closed subscheme of $P_{X/S}(n)$ defined by the ideal $\ov{\calI}(n)^{[k+1]}.$ Since any element of a PD-ideal in characteristic $p$ is nilpotent, the canonical morphisms $X\ra P_{X/S}(n)$ and $P^k_{X/S}(n)\ra P_{X/S}(n)$ are thickenings (\cite{Ogus78} 3.31). So they induce universal homeomorphisms of the underlying topological spaces. We identify the underlying étale sites of the schemes $X,\ P_{X/S}(n)$ and $P^k_{X/S}(n)$ by these thickenings. We also denote by $p_1,p_2:P_{X/S}(1)\ra X$ and $p_{ij}:P_{X/S}(2)\ra X\times_SX,\ 1\le i<j\le 3$ (resp. $p_1^k,p_2^k:P_{X/S}^k(1)\ra X$ and $p_{ij}^k:P_{X/S}^k(2)\ra X\times_SX$) the canonical projections. When $n=1,$ we will drop the $(1)$ and write $P_{X/S}$ (resp. $P^k_{X/S},$ resp. $\ov{\calI}$) instead of $P_{X/S}(1)$ (resp. $P^k_{X/S}(1),$ resp. $\ov{\calI}(1)$). Recall that there exists a canonical isomorphism (\cite{Kat89} 5.8)
$$\ov{\calI}/\ov{\calI}^{[2]} \xrightarrow{\sim} \omega^1_{X/S}.$$
If $\calE$ is an $\Ox_X$-module, we denote by $\Ox_{P_{X/S}}\otimes_{\Ox_X}\calE$ (resp. $\Ox_{P_{X/S}^k}\otimes_{\Ox_X}\calE,$ resp. $\calE\otimes_{\Ox_X}\Ox_{P_{X/S}},$ resp. $\calE\otimes_{\Ox_X}\Ox_{P_{X/S}^k}$) the $\Ox_X$-module $p_{1*}p_2^*\calE$ (resp. $p^k_{1*}p_2^{k*}\calE,$ resp. $p_{2*}p_1^*\calE,$ resp. $p^k_{2*}p_1^{k*}\calE$).\\ 
Recall that there exists a canonical isomorphism $P_{X/S}\times_XP_{X/S}\overset{\sim}{\ra} P_{X/S}(2)$ where the left $P_{X/S}$ in the fiber product is considered as a scheme over $X$ via the second projection $p_2$ and the right $P_{X/S}$ is considered as a scheme over $X$ via the first projection $p_1.$ We denote by $\delta:\Ox_{P_{X/S}}\ra \Ox_{P_{X/S}}\otimes_{\Ox_X}\Ox_{P_{X/S}}$ the homomorphism of structural rings associated with the composition
$$P_{X/S}\times_XP_{X/S}\ra P_{X/S}(2)\xrightarrow{p_{13}}P_{X/S}.$$ 

Recall also that for any integers $n,m\ge 0,$ there exists a canonical morphism $P_{X/S}^n\times_XP_{X/S}^m\ra P_{X/S}^{n+m}(2)$ making the following diagram commutative
$$\begin{tikzcd}
X\ar[equal]{rrr}\ar{d} & & & X\ar{d}\ar{dr} & \\
P_{X/S}^n\times_XP_{X/S}^m\ar[bend right=20]{rrrr}\ar{r} & (X\times_SX)\times_X(X\times_SX)\ar{r}{\sim} & X^3\ar{r}{p_{13}}& X^2 & P_{X/S}^{n+m}\ar{l}
\end{tikzcd}$$
We denote by $\delta^{n,m}:\Ox_{P^{n+m}_{X/S}}\ra \Ox_{P^n_{X/S}}\otimes_{\Ox_X}\Ox_{P^m_{X/S}}$ the homomorphism of structural rings associated with the morphism
$$P_{X/S}^n\times_XP_{X/S}^m\ra P_{X/S}^{n+m}.$$  
\end{parag}

\begin{parag}\label{P1}
We denote by $\ov{\calI}$ the ideal of the exact closed immersion $i:X\ra P_{X/S}.$ We have a canonical exact sequence of abelian groups
\begin{equation}
1\ra i^{-1}(1+\ov{\calI})\xrightarrow{\lambda} i^{-1}\calM_{P_{X/S}}^{gp}\ra \calM_X^{gp}\ra 1.
\end{equation}
Similarly, if $\ov{\calI}_k$ is the ideal of the exact closed immersion $i_k:X\ra P^k_{X/S},$ then we have a canonical exact sequence of abelian groups
\begin{equation}1\ra i_k^{-1}(1+\ov{\calI}_k)\xrightarrow{\lambda_k}i_k^{-1}\calM^{gp}_{P^k_{X/S}}\ra \calM^{gp}_X\ra 1.\end{equation}
We define a morphism of sheaves $\mu:\calM^{gp}_X\ra i^{-1}\Ox_{P_{X/S}}$ (resp. $\mu_k:\calM^{gp}_X\ra i_k^{-1}\Ox_{P_{X/S}^k}$) as follows:
if $m$ is a local section of $\calM_X^{gp}$ then $p^{\flat}_2(m)$ et $p^{\flat}_1(m)$ (resp. $p_2^{k,\flat}(m)$ and $p_1^{k,\flat}(m)$) have the same image in $\calM^{gp}_X$ so there exists a unique local section $\mu(m)$ (resp. $\mu_k(m)$) of $i^{-1}(1+\ov{\calI})$ (resp. $i_k^{-1}(1+\ov{\calI}_k)$) such that $\lambda(\mu(m))=(i^{-1}p^{\flat}_2)(m)-(i^{-1}p^{\flat}_1)(m)$ (resp. $\lambda_k(\mu_k(m))=(i_k^{-1}p^{k,\flat}_2)(m)-(i_k^{-1}p^{k,\flat}_1)(m)$). It is clear that for any local sections $m$ and $m'$ of $\calM^{gp}_X,$ $\mu(m+m')=\mu(m)\mu(m')$ (resp. $\mu_k(m+m')=\mu_k(m)\mu_k(m')$).
We denote $\mu(m)-1$ by $\eta(m).$
\end{parag}

\begin{definition} \label{defmokrez1}
Let $\calE$ and $\F$ be two $\Ox_X$-modules and $k$ a nonnegative integer. A \emph{differential operator of order $\le k$ from $\calE$ to $\F$} is an $\Ox_X$-linear homomorphism 
\begin{equation}
\phi:\Ox_{P_{X/S}^k}\otimes_{\Ox_X}\calE\ra \F,
\end{equation} 
where $\Ox_{P_{X/S}^k}\otimes_{\Ox_X}\calE$ is considered as an $\Ox_X$-module via the first projection $p^k_1:P_{X/S}^k\ra X.$
For such an operator, we define 
\begin{equation}
\phi^{\flat}:\begin{array}[t]{clc}\calM_X\times \calE &\rightarrow &\F\\ (m,x)&\mapsto & \phi(\mu_k(m)\otimes x). \end{array}\end{equation} 

We denote by $\calD^k_{X/S}(\calE,\F)$ the sheaf of differential operators of order $\le k$ from $\calE$ to $\F$ and
$$\calD_{X/S}(\calE,\F)=\bigcup_{k\ge 0}\calD^k_{X/S}(\calE,\F).$$
For $\calE=\F=\Ox_X,$ we simply denote $\calD^k_{X/S}(\calE,\F)$ and $\calD_{X/S}(\calE,\F)$ by $\calD^k_{X/S}$ and $\calD_{X/S}$ respectively.

If $f:\Ox_{P_{X/S}^m}\otimes_{\Ox_X}\calE\ra \F$ and $g:\Ox_{P_{X/S}^n}\otimes_{\Ox_X}\F\ra \calG$ are two differential operators, then their composition $g\circ f$ if defined as follows
$$\Ox_{P_{X/S}^{n+m}}\otimes_{\Ox_X}\calE\xrightarrow{\delta^{n,m}\otimes \op{Id}}\Ox_{P_{X/S}^n}\otimes_{\Ox_X}\Ox_{P_{X/S}^m}\otimes_{\Ox_X}\calE\xrightarrow{\op{Id}\otimes f} \Ox_{P_{X/S}^n}\otimes_{\Ox_X}\F\xrightarrow{g} \calG.$$ 
The group $\calD_{X/S}^k(\calE,\F)$ is equipped with an $\Ox_{P_{X/S}}$-module structure by setting for any local sections $z$ and $f$ of $\Ox_{P_{X/S}}$ and $\calD_{X/S}^k$ respectively,
$$z\cdot f:x\mapsto f(zx).$$ 
\end{definition}

\begin{exmp}\label{exmp24}
Consider the split exact sequence of $\Ox_X$-modules
$$0 \ra \ov{\calI}/\ov{\calI}^{[2]} \ra \Ox_{P_{X/S}}/\ov{\calI}^{[2]} \ra \Ox_{P_{X/S}}/\ov{\calI} \ra 0.$$
We have canonical isomorphisms
$$\ov{\calI}/\ov{\calI}^{[2]} \xrightarrow{\sim} \omega^1_{X/S},\ \Ox_{P_{X/S}}/\ov{\calI}\xrightarrow{\sim} \Ox_X,$$
hence the split exact sequence
$$0 \ra \omega^1_{X/S} \ra \Ox_{P_{X/S}^1} \ra \Ox_X \ra 0.$$
Considering the splitting given by the first projection $p^1_1:P_{X/S}^1 \ra X,$ we get an isomorphism
$$\Ox_{P_{X/S}^1} \xrightarrow{\sim}\Ox_X\oplus \omega^1_{X/S}.$$
A logarithmic derivation $\partial:\Ox_X\times \calM_X^{gp} \ra \Ox_X$ corresponds to an $\Ox_X$-linear map
$$\partial:\omega^1_{X/S} \ra \Ox_X.$$
Composing with the projection
$$\Ox_{P_{X/S}^1}=\Ox_X\oplus \omega^1_{X/S}\ra \omega^1_{X/S},$$
we obtain a differential operator of degree $\le 1$
$$\partial: \Ox_{P_{X/S}^1} \ra \Ox_X.$$
Via this construction, we will often consider sections of the tangent sheaf $\calT_{X/S}=\mathscr{Hom}_{\Ox_X}(\omega^1_{X/S},\Ox_X)$ as sections of $\calD_{X/S}^1.$
\end{exmp}

\begin{definition} \label{def14logop}
Let $\calE$ and $\F$ be two $\Ox_X$-modules. A \emph{HPD-differential operator from $\calE$ to $\F$} is an $\Ox_X$-linear morphism 
$$\phi:\Ox_{P_{X/S}}\otimes_{\Ox_X}\calE\ra \F,$$
where $\Ox_{P_{X/S}}\otimes_{\Ox_X}\calE$ is considered as an $\Ox_X$-module via the first projection $p_1:P_{X/S}\ra X.$ 
For such an operator, we define 
$$\phi^{\flat}:\begin{array}[t]{clc}\calM_X\times \calE & \rightarrow & \F\\ (m,x) & \mapsto & \phi(\mu(m)\otimes x). \end{array}$$ 

We denote by $\widehat{\calD}_{X/S}(\calE,\F)$ the sheaf of HPD-differential operators from $\calE$ to $\F.$ If $\calE=\F=\Ox_X,$ we denote $\widehat{\calD}_{X/S}(\calE,\F)$ simply by $\widehat{\calD}_{X/S}.$

If $f:\Ox_{P_{X/S}}\otimes_{\Ox_X}\calE\ra \F$ and $g:\Ox_{P_{X/S}}\otimes_{\Ox_X}\F\ra \calG$ are two differential operators then their composition $g\circ f$ if defined as follows
$$\Ox_{P_{X/S}}\otimes_{\Ox_X}\calE\xrightarrow{\delta\otimes \op{Id}}\Ox_{P_{X/S}}\otimes_{\Ox_X}\Ox_{P_{X/S}}\otimes_{\Ox_X}\calE\xrightarrow{\op{Id}\otimes f} \Ox_{P_{X/S}}\otimes_{\Ox_X}\F\xrightarrow{g} \calG.$$ 
The ring $\widehat{\calD}_{X/S}(\calE,\F)$ is equipped with an $\Ox_{P_{X/S}}$-module structure by setting for any local sections $z$ and $f$ of $\Ox_{P_{X/S}}$ and $\widehat{\calD}_{X/S}$ respectively,
$$z\cdot f:x\mapsto f(zx).$$
\end{definition}

\begin{parag}\label{P2}
In what follows, we suppose that $X$ is smooth over $S.$ So any point of $X$ has an open neighborhood $U$ and local sections $m_1,\hdots,m_d\in \Gamma(U,\calM_X)$ such that $(\op{dlog}m_i)_{1\le i\le d}$ form a basis of $\omega^1_{U/S}.$ Such local sections $m_i$ will be refered to as \emph{local coordinates for $X$ over $S$}. For all $i=1,\hdots,d,$ let $\eta_i=\mu(m_i)-1\in \Gamma(U,\Ox_{P_{X/S}})$ (\ref{P1}). We have an isomorphism of $\Ox_U$-PD-algebras (\cite{Kat89} Proposition 6.5) given by
$$\begin{array}[t]{clc}\Ox_U\langle T_1,\hdots,T_d\rangle & \ra & \Ox_{P_{U/S}}\\ T_i& \mapsto &\eta_i\end{array}$$
where $\Ox_U\langle T_1,\hdots,T_d\rangle$ is the PD-polynomial algebra on the PD-ring $\Ox_U$ and $\Ox_{P_{U/S}}$ is an $\Ox_U$-algebra by the first projection $p_1:P_{U/S}\ra U$ (we still get an isomorphism if we use $p_2$ instead of $p_1$). We denote by $(\partial_I)_{I\in \mathbb{N}^d}$ the dual basis of $(T^{[I]})_{I\in\mathbb{N}^d},$ so that for any multi-index $I\in\mathbb{N}^d$
\begin{equation}\label{partialI}
\partial_I:\Ox_{P_{U/S}}\ra \Ox_U,\ T^{[J]}\mapsto \delta_{IJ}=\begin{cases}1\ \op{if}\ I=J\\ 0\ \op{else} \end{cases}\end{equation}
where $T^{[J]}=\prod_{i=1}^dT_i^{[J_i]}.$ The homomorphism $\partial_I$ is clearly an HPD-differential operator. The operator $\partial_I$ annihilates $\eta^{[J]}$ for any $J\in \N^d$ such that $|J|>n=|I|.$ Hence it induces a differential operator $\Ox_{P_{U/S}^n} \ra \Ox_U$ of order $\le n,$ that we will also denote by $\partial_I.$
Note that, if $\partial_i$ denotes the derivation
$$\partial_i:\begin{array}[t]{clc}\omega^1_{U/S} & \ra & \Ox_U\\
\ \omega=\sum_{j=1}^d\omega_j\op{dlog}m_j & \mapsto & \omega_i,
\end{array}$$
where $\omega_j\in \Ox_U$ for $1\le j\le d,$
then the differential operator defined by $\partial_i,$ as in \ref{exmp24}, is the same as $\partial_{\epsilon_i},$ where $\epsilon_i \in \N^d$ is the multi-index whose all coefficients are zero except for the $i^{th}$ which is equal to $1.$
\end{parag}

\begin{lemma}\label{lemma16}
Recall the notation introduced in \ref{P1}. For any local section $m$ of $\calM_X,$ we have
$$\delta(\mu(m))=\mu(m)\otimes \mu(m)\in \Ox_{P_{X/S}}\otimes_{\Ox_X} \Ox_{P_{X/S}},$$
$$\delta(\eta(m))=\eta(m)\otimes \eta(m)+\eta(m)\otimes 1+1\otimes \eta(m) \in \Ox_{P_{X/S}}\otimes_{\Ox_X} \Ox_{P_{X/S}}.$$
\end{lemma}

\begin{proof}
The exact closed immersion $j:X\ra P_{X/S}(2)$ of ideal $\calJ$ gives rise to the canonical exact sequence of abelian groups
$$1\ra j^{-1}(1+\calJ)\xrightarrow{\nu} j^{-1}\calM^{gp}_{P_{X/S}(2)}\ra \calM^{gp}_X\ra 1.$$
Then by the same argument as in \ref{P1}, for any local section $m\in \calM^{gp}_X$ and any $i,j\in \{1,2,3\},$ there exists a unique section $\mu_{ij}(m)\in j^{-1}(1+\calJ)$ such that
$$\nu(\mu_{ij}(m))=\pi_j^{\flat}(m)-\pi_i^{\flat}(m)$$
where $\pi_i:P_{X/S}(2)\ra X$ are the canonical projections. 
The sections $\mu_{ij}(m)$ clearly satisfy the cocycle condition
$$\mu_{ij}(m)\mu_{jk}(m)=\mu_{ik}(m).$$
Denote by $p_{ij}:P_{X/S}(2)\ra P_{X/S}$ the canonical projections and by $p_{ij}^{\#}$ the associated morphisms of structure sheaves. We have the following commutative diagram
$$\begin{tikzcd}[ampersand replacement=\&]
\Ox_{P_{X/S}}\ar{r}{p_{ij}^{\#}}\& \Ox_{P_{X/S}(2)}\\
\calM_{P_{X/S}}\ar{u}{\alpha_{P_{X/S}}}\ar{r}{p_{ij}^{\flat}}\& \calM_{P_{X/S}(2)}\ar[swap]{u}{\alpha_{P_{X/S}(2)}}
\end{tikzcd}.$$
It follows that
\begin{alignat*}
5
\nu(p_{ij}^{\#}(\mu(m)))&=p_{ij}^{\flat}(\lambda(\mu(m)))\\
&=p_{ij}^{\flat}(p_2^{\flat}m-p_1^{\flat}m)\\
&=\pi_j^{\flat}m-\pi_i^{\flat}m\\
&=\nu(\mu_{ij}(m)).
\end{alignat*}
So $p_{ij}^{\#}(\mu(m))=\mu_{ij}(m).$ In particular, $\delta(\mu(m))=\mu_{13}(m)=\mu_{12}(m)\mu_{23}(m)$ and the result follows.
\end{proof}

\begin{proposition}\label{prop17}
Suppose that $X\ra S$ has local coordinates $m_1,\hdots,m_d\in \Gamma(X,\calM_X).$ For any multi-index $I\in \mathbb{N}^d$ and any $1\le i\le d,$ denote by $\partial_I$ the differential operator defined in \ref{P2} and let $\epsilon_i=(0,\hdots,0,1,0\hdots,0)$ be the multi-index whose all coefficients are zero except for the $i$th which is equal to $1.$ Then
$$\partial_J\circ \partial_{\epsilon_i}=\partial_{J+\epsilon_i}+J_i\partial_J\ \forall J=(J_1,\hdots,J_d)\in\mathbb{N}^d\ \forall 1\le i\le d.$$
\end{proposition}

\begin{proof}
Let $\eta_i=\mu(m_i)-1.$
Recall that $\partial_J\circ \partial_{\epsilon_i}$ is the composition
$$\Ox_{P_{X/S}}\xrightarrow{\delta} \Ox_{P_{X/S}}\otimes_{\Ox_X}\Ox_{P_{X/S}}\xrightarrow{\op{Id}\otimes \partial_{\epsilon_i}}\Ox_{P_{X/S}}\xrightarrow{\partial_J} \Ox_X.$$
By lemma \ref{lemma16}, $\delta(\mu(m))=\mu(m)\otimes \mu(m)$ for any local section $m$ of $\calM_X.$ In particular, $$\delta(\eta_i)=\delta(\mu(m_i))-1=(\eta_i+1)\otimes(\eta_i+1)-1=\eta_i\otimes\eta_i+\eta_i\otimes 1+1\otimes\eta_i$$
\begin{alignat*}
5
\delta(\eta_i^{[n_i]})&=\delta(\eta_i)^{[n_i]}\\
&=\left (\eta_i\otimes\eta_i+\eta_i\otimes 1+1\otimes \eta_i \right )^{[n_i]}\\
&=\sum_{\alpha+\beta+\gamma=n_i}\alpha !\begin{pmatrix}\alpha+\beta \\ \alpha \end{pmatrix} \begin{pmatrix}\alpha+\gamma \\ \alpha\end{pmatrix}\eta_i^{[\alpha+\beta]}\otimes\eta_i^{[\alpha+\gamma]}
\end{alignat*}
We can check that this formula remains true if we use multi-indices instead. In other words, if $n=(n_1,\hdots,n_d)\in\mathbb{N}^d$ and $\eta^{[n]}=\prod_{i=1}^d\eta_i^{[n_i]}$ then
\begin{equation}\label{eq171}
\delta(\eta^{[n]})=\sum_{\alpha+\beta+\gamma=n}\alpha !\begin{pmatrix}\alpha+\beta \\ \alpha \end{pmatrix} \begin{pmatrix}\alpha+\gamma \\ \alpha\end{pmatrix}\eta^{[\alpha+\beta]}\otimes\eta^{[\alpha+\gamma]}.
\end{equation}
It follows that
\begin{alignat*}
5
(\partial_J\circ \partial_{\epsilon_i})(\eta^{[n]})&= \partial_J\left ((\op{Id}\otimes \partial_{\epsilon_i})( \delta(\eta^{[n]}))\right )\\
&=\partial_J(\eta^{[n-\epsilon_i]}+n_i\eta^{[n]})\\
&=(\partial_{J+\epsilon_i}+J_i\partial_J)(\eta^{[n]}).
\end{alignat*}
\end{proof}

\begin{corollaire}\label{cor18}
Suppose that $X\ra S$ has local coordinates $m_1,\hdots,m_d\in \Gamma(X,\calM_X)$ and consider the differential operators $\partial_I$ defined in \ref{P2}. For any multi-indices $I,J\in\mathbb{N}^d,$ any positive integer $k$ and $i,j\in \{1,\hdots,d\},$ we have the following formulae:
\begin{enumerate}
\item $$\partial_{k\epsilon_i}=\prod_{j=0}^{k-1}(\partial_{\epsilon_i}-j).$$
\item If $I_j=0$ then
$$\partial_I\circ \partial_{k\epsilon_j}=\partial_{I+k\epsilon_j}.$$
\item $$\partial_{\epsilon_i}\circ \partial_J=\partial_{J+\epsilon_i}+J_i\partial_J=\partial_J\circ \partial_{\epsilon_i}$$
\item $$\partial_I\circ \partial_J = \partial_J\circ \partial_I.$$
\end{enumerate}
\end{corollaire}

\begin{lemma}\label{lemiran1}
Suppose that $X\ra S$ has local coordinates $m_1,\hdots,m_d\in \Gamma(X,\calM_X)$ and consider the differential operators $\partial_I$ defined in \ref{P2}. For any integers $n,r\ge 0$ and $1\le i\le d,$
\begin{equation}\label{eqiran2}
\partial_{(pn+r)\epsilon_i}=(\partial_{p\epsilon_i})^n \circ \partial_{r\epsilon_i}.
\end{equation}
For any $I,J\in \N^d,$
\begin{equation}\label{eq3102}
\partial_I\circ \partial_{pJ}=\partial_{I+pJ}.
\end{equation}
\end{lemma}

\begin{proof}
By \ref{cor18}, since we work in characteristic $p,$
\begin{alignat*}{2}
\partial_{(pn+r)\epsilon_i} &= \prod_{k=0}^{pn+r-1}(\partial_{\epsilon_i}-k) \\
&= \prod_{k=0}^{pn-1}(\partial_{\epsilon_i}-k) \circ \prod_{k=pn}^{pn+r-1}(\partial_{\epsilon_i}-k) \\
&= \prod_{l=0}^{n-1} \prod_{k=pl}^{p(l+1)-1}(\partial_{\epsilon_i}-k) \circ \prod_{k=pn}^{pn+r-1}(\partial_{\epsilon_i}-k) \\
&= \prod_{l=0}^{n-1} \prod_{k=0}^{p-1}(\partial_{\epsilon_i}-k) \circ \prod_{k=0}^{r-1}(\partial_{\epsilon_i}-k) \\
&= (\partial_{p\epsilon_i})^n \circ \partial_{r\epsilon_i}.
\end{alignat*}
Now the second equality. We have
\begin{alignat*}{2}
\partial_{I+pJ} &= \prod_{k=1}^d\partial_{(I_k+pJ_k)\epsilon_k} \\
&= \prod_{k=1}^d\left (\partial_{p\epsilon_k}\right )^{J_k}\circ \partial_{I_k\epsilon_k} \\
&= \prod_{k=1}^d\left (\partial_{p\epsilon_k}\right )^{J_k} \circ \prod_{k=1}^d\partial_{I_k\epsilon_k} \\
&= \prod_{k=1}^d\partial_{pJ_k\epsilon_k} \circ \prod_{k=1}^d\partial_{I_k\epsilon_k} \\
&= \partial_{pJ}\circ \partial_I,
\end{alignat*}
where the second and forth lines result from \eqref{eqiran2} and the others from \ref{cor18}.
\end{proof}

\begin{lemma}\label{lem19}
For any local section $z$ of $\Ox_{P_{X/S}}$ and any multi-index $N\in\mathbb{N}^d$
$$z\cdot \partial_N=\sum_{I+J=N}\begin{pmatrix}N\\ I \end{pmatrix}\partial_I(z)\partial_J.$$ 
\end{lemma}

\begin{proof}
Since both expressions are $\Ox_X$-linear in $z,$ it is sufficient to prove the equality for $z=\eta^{[r]}$ where $r\in \mathbb{N}^d.$ In this case, for a multi-index $s\in \mathbb{N}^d$
\begin{alignat*}
5
(z\cdot \partial_N)(\eta^{[s]})&=\begin{pmatrix}r+s \\r \end{pmatrix}\partial_N (\eta^{[r+s]})\\
&=\begin{cases}\begin{pmatrix}N \\r \end{pmatrix}\ \op{if}\ r+s=N\\ 0\ \op{else} \end{cases}\\
&=\sum_{I+J=N}\begin{pmatrix}N \\I \end{pmatrix}\partial_I(z)\partial_J(\eta^{[s]}).
\end{alignat*} 
\end{proof}

\newpage

\section{Connections and stratifications}

We start by a basic review of $\lambda$-connections.
\begin{definition}
Let $f:X\ra S$ be a morphism of logarithmic schemes, $\calE$ an $\Ox_X$-module and $\lambda$ an integer. A \emph{logarithmic $\lambda$-connection on $\calE$} is an $f^{-1}\Ox_S$-linear morphism $\nabla:\calE\ra \calE\otimes_{\Ox_X}\omega^1_{X/S}$ satisfying the Leibniz rule
$$\nabla(ax)=a\nabla(x)+\lambda x\otimes da,$$
for all local sections $a$ and $x$ of $\Ox_X$ and $\calE$ respectively and where $d:\Ox_X\ra \omega^1_{X/S}$ denotes the universal derivation. 
Logarithmic $1$-connections are called \emph{logarithmic connections} and logarithmic $0$-connections are called \emph{logarithmic Higgs fields.} 
\end{definition} 

In the rest of this manuscript, we will always drop the term "logarithmic" and refer to logarithmic connections simply by connections.

\begin{proposition}\label{prop42}
Let $X\ra S$ be a morphism of logarithmic schemes, $\calE$ an $\Ox_X$-module and $\lambda$ an integer. A $\lambda$-connection $\nabla$ on $\calE$ induces, for any positive integer $i,$ a morphism
$$\nabla^i:\calE\otimes_{\Ox_X}\omega^i_{X/S}\ra \calE\otimes_{\Ox_X}\omega^{i+1}_{X/S}$$
defined for any local sections $x$ and $\omega$ of $\calE$ and $\omega^i_{X/S}$ respectively by
$$\nabla^i(x\otimes \omega)=\nabla(x)\wedge \omega+\lambda x\otimes d\omega$$
where $\nabla(x)\wedge \omega$ denotes the image of $\nabla(x)\otimes \omega$ by the canonical morphism $\calE\otimes_{\Ox_X} \omega^1_{X/S}\otimes_{\Ox_X}\omega^i_{X/S}\ra \calE\otimes_{\Ox_X}\omega^{i+1}_{X/S}.$
\end{proposition}

\begin{proof}
Let $i\ge 1.$ 
For any local sections $a,\ x,\ y,\ \omega$ and $\eta$ of $\Ox_X,\ \calE,\ \calE,\ \omega^i_{X/S}$ and $\omega^i_{X/S}$ respectively, we have the following equalities:
\begin{alignat*}
5
\nabla^i((x+y)\otimes \omega)&=\nabla^i(x\otimes \omega)+\nabla^i(y\otimes \omega).\\
\nabla^i(x\otimes (\omega+\eta))&=\nabla^i(x\otimes\omega)+\nabla^i(x\otimes\eta).\\
\nabla^i((ax)\otimes \omega)&=a\nabla(x)\wedge \omega+\lambda x\otimes(da\wedge \omega)+ax\otimes d\omega\\
&=\nabla^i(x\otimes (a\omega)).
\end{alignat*}
It follows that
$$\nabla^i:x\otimes \omega\mapsto \nabla(x)\wedge\omega+\lambda x\otimes d\omega$$ 
is well-defined and the proof is complete.
\end{proof}

\begin{definition}\label{definteglambdaconn}
Let $X\ra S$ be a morphism of logarithmic schemes, $\calE$ an $\Ox_X$-module, $\lambda$ an integer and $\nabla:\calE\ra \calE\otimes_{\Ox_X}\omega^1_{X/S}$ a $\lambda$-connection on $\calE.$ The \emph{curvature of $\nabla$,} denoted by $K(\nabla),$ is the composition $\nabla^1\circ \nabla$ (where $\nabla^1$ is defined in \ref{prop42}). We say that the $\lambda$-connection $\nabla$ is \emph{integrable} if $K(\nabla)=0.$ We denote by $\boldsymbol{\op{MIC}^{\lambda}(X/S)}$ the category of $\Ox_X$-modules equipped with an integrable $\lambda$-connection.
\end{definition}
If $\partial:\Ox_X\times\calM_X\ra \Ox_X$ is a logarithmic derivation, we denote by $\nabla(\partial)$ the composition
\begin{equation}\label{nablacomp}
\calE\xrightarrow{\nabla}\calE\otimes_{\Ox_X}\omega^1_{X/S}\xrightarrow{\op{Id}_{\calE}\otimes \partial}\calE\otimes_{\Ox_X}\Ox_X\xrightarrow{\sim}\calE
\end{equation}
where $\partial:\omega^1_{X/S}\ra \Ox_X$ is the $\Ox_X$-linear homomorphism induced by the logarithmic derivation $\partial$ and $\calE\otimes_{\Ox_X}\Ox_X\xrightarrow{\sim}\calE$ is the canonical isomorphism. It is clear that $\nabla(\partial)$ is $f^{-1}\Ox_S$-linear.

If $\partial_1,\partial_2:\Ox_X\times\calM_X\ra \Ox_X$ are two logarithmic derivations, we denote by $K(\nabla)(\langle \partial_1,\partial_2\rangle)$ the composition
$$\calE\xrightarrow{K(\nabla)}\calE\otimes_{\Ox_X}\omega^2_{X/S}\xrightarrow{\op{Id}_{\calE}\otimes (\langle\partial_1, \partial_2\rangle)} \calE\otimes_{\Ox_X}\Ox_X\xrightarrow{\sim}\calE$$
where $\langle\partial_1, \partial_2\rangle$ denotes the morphism
$$\langle\partial_1, \partial_2\rangle:\begin{array}[t]{clc}\omega^2_{X/S} & \ra & \Ox_X\\ \omega_1\wedge \omega_2 & \mapsto & \partial_1(\omega_1)\partial_2(\omega_2)-\partial_1(\omega_2)\partial_2(\omega_1). \end{array}$$

\begin{parag}\label{parag23}
The algebra $\op{Der}_{X/S}(\Ox_X)$ of logarithmic derivations $\Ox_X\times \calM_X\ra \Ox_X$ relative to $S$ has a Lie algebra structure given by the Lie bracket 
$$[\partial_1,\partial_2]=([D_1,D_2],D_1\delta_2-D_2\delta_1)$$ 
for all logarithmic derivations $\partial_i=(D_i,\delta_i):\Ox_X\times\calM_X\ra \Ox_X,\ i=1,2$ (\cite{Ogus2018} V 2.1.2). 
\end{parag}

\begin{proposition}
Let $X\ra S$ be a morphism of logarithmic schemes and $\nabla:\calE\ra \calE\otimes_{\Ox_X}\omega^1_{X/S}$ be a $\lambda$-connection on an $\Ox_X$-module $\calE.$ We denote by $K(\nabla)$ the curvature of $\nabla.$
For all logarithmic derivations $\partial_i=(D_i,\delta_i):\Ox_X\times\calM_X\ra \Ox_X,\ i=1,2$ 
$$
\lambda \nabla([\partial_1,\partial_2])=[\nabla(\partial_1),\nabla(\partial_2)]-K(\nabla)(\langle\partial_1, \partial_2\rangle).
$$
We say that $\nabla$ is integrable if $K(\nabla)=0.$
\end{proposition}

\begin{proof}
Let $x$ be a local section of $\calE.$ The sheaf $\omega^1_{X/S}$ is locally generated by sections of the form $\op{dlog}m,$ so there exist local sections $x_i$ and $x_{ij},$ $1\le i,j\le d,$ of $\calE$ and local sections $m_i,\ 1\le i\le d,$ of $\calM_X$ such that
$$\nabla(x)=\sum_{i=1}^dx_i\otimes \op{dlog}m_i$$
and
$$\nabla(x_i)=\sum_{j=1}^dx_{ij}\otimes \op{dlog}m_j\ \forall 1\le i\le d.$$
A straightforward computation then yields
\begin{alignat*}
5
\nabla([\partial_1,\partial_2])(x)&=\sum_{i=1}^d(D_1(\delta_2(m_i))-D_2(\delta_1(m_i)))x_i,\\
\nabla(\partial_1)\circ \nabla(\partial_2)(x)&=\sum_{1\le i,j\le d}\delta_2(m_i)\delta_1(m_j)x_{ij}+\lambda \sum_{i=1}^dD_1(\delta_2(m_i))x_i,\\
\nabla(\partial_2)\circ \nabla(\partial_1)(x)&=\sum_{1\le i,j\le d}\delta_1(m_i)\delta_2(m_j)x_{ij}+\lambda \sum_{i=1}^dD_2(\delta_1(m_i))x_i,\\
K(\nabla)(\langle\partial_1,\partial_2\rangle)&=\sum_{1\le i,j\le d}(\delta_1(m_j)\delta_2(m_i)-\delta_2(m_j)\delta_1(m_i))x_{ij},
\end{alignat*}
and the result follows.
\end{proof}

\begin{corollaire}\label{cor25}
Let $f:X\ra S$ be a morphism of logarithmic schemes and $\nabla:\calE\ra \calE\otimes_{\Ox_X}\omega^1_{X/S}$ be a $\lambda$-connection on an $\Ox_X$-module $\calE.$ If $\nabla$ is integrable then $\nabla:\op{Der}_{X/S}(\Ox_X)\ra \op{End}_{f^{-1}\Ox_S}(\calE)$ is a morphism of Lie algebras, i.e. for all logarithmic derivations $\partial_1,\partial_2:\Ox_X\times\calM_X\ra \Ox_X,$ 
$$\nabla([\partial_1,\partial_2])=[\nabla(\partial_1),\nabla(\partial_2)].$$
\end{corollaire}

\begin{parag}\label{P46}
In the remaining of this section, we fix a smooth morphism $f:X\ra S$ of logarithmic schemes of positive characteristic $p.$
For any logarithmic derivation $\partial=(D,\delta):\Ox_X\times\calM_X\ra \Ox_X,$ let
$$\partial^{(p)}=(D^p,F_X^{\#}\circ \delta+D^{p-1}\circ \delta)$$
where $F_X:X\ra X$ is the absolute Frobenius morphism of $X$ (\ref{Not3}) and $F_X^{\#}:\Ox_X\ra \Ox_X$ is the induced homomorphism of structural rings. Then, $\partial^{(p)}$ is also a logarithmic derivation and the $p$-operation $\partial\mapsto \partial^{(p)}$ defines a restricted Lie algebra structure on $\op{Der}_{X/S}(\Ox_X)$ (\cite{Ogus94} proposition 1.2.1). The differential operators $\partial_{\epsilon_i},$ defined in \ref{prop17}, satisfy $\partial_{\epsilon_i}^{(p)}=\partial_{\epsilon_i}$ (\cite{Ogus94} remark 1.2.2).
We view a logarithmic derivation $\partial:\Ox_X\times \calM_X \ra \Ox_X$ as a differential operator (\ref{exmp24}) and we denote by $\psi(\partial)$ the differential operator of degree $\le p$
$$\psi(\partial)=\partial^p-\partial^{(p)},$$
where $\partial^p$ denotes the composition of differential operators.
If $U$ is an étale $X$-scheme and $\partial:\omega^1_{U/S} \ra \Ox_{U}$ is a section of the tangent sheaf $\calT_{X/S}:=\mathscr{Hom}_{\Ox_{X}}(\omega^1_{X/S},\Ox_{X})$ over $U,$ then $\psi(\partial)$ is a section of $F_{X*}\calD_{X/S}$ over $U.$ So we have a morphism
$$
\psi:\begin{array}[t]{clc}
\calT_{X/S} & \ra & F_{X*}\calD_{X/S} \\
\partial & \mapsto & \psi(\partial).
\end{array}
$$
We define the $p$-curvature $\psi_{\nabla}$ of the connection $\nabla$ as in the usual case: 
\begin{equation}
\psi_{\nabla}:\begin{array}[t]{clc}\op{Der}_{X/S}(\Ox_X) & \ra & \op{End}_{f^{-1}\Ox_S}(\calE)\\ \partial & \mapsto & \nabla(\partial)^p-\nabla(\partial^{(p)}).\end{array}
\end{equation}
For any logarithmic derivation $\partial=(D,\delta):\Ox_X\times \calM_X\ra \Ox_X,$ the endomorphism $\psi_{\nabla}(\partial)$ is $\Ox_X$-linear. Indeed, for any local sections $a$ and $x$ of $\Ox_X$ and $\calE$ respectively,
\begin{alignat}{2}
\nabla(\partial)(ax) &= a\nabla(\partial)(x)+D(a)x \label{eq461}\\
\nabla(\partial^{(p)})(ax) &= a\nabla(\partial^{(p)})(x)+D^p(a)x.
\end{alignat} 
We can prove by induction that, for all positive integers $k,$
\begin{equation}\label{eq463}
\nabla(\partial)^k(ax)=\sum_{i=0}^k\begin{pmatrix}k\\ i\end{pmatrix}D^{k-i}(a)\nabla(\partial)^i(x).
\end{equation}
For $k=p,$ we get
$$\nabla(\partial)^p(ax)=a\nabla(\partial)^p(x)+D^p(a)x.$$
The $\Ox_X$-linearity of $\psi_{\nabla}(\partial)$ follows:
\begin{alignat*}{2}
\psi_{\nabla}(\partial)(ax) &= \nabla(\partial)^p(ax)-\nabla(\partial^{(p)})(ax) \\
&= a\nabla(\partial)^p(x)+D^p(a)x-a\nabla(\partial^{(p)})(x)-D^p(a)x\\
&= a\psi_{\nabla}(\partial)(x).
\end{alignat*}
\end{parag}

\begin{definition}
Let $\calE$ be an $\Ox_X$-module. A \emph{stratification on $\calE$} is a sequence of morphisms $(\varepsilon_n:\Ox_{P^n_{X/S}}\otimes_{\Ox_X}\calE\ra \calE\otimes_{\Ox_X}\Ox_{P^n_{X/S}})_{n\ge 0}$ satisfying the following conditions:
\begin{enumerate}
\item $\varepsilon_n$ is an $\Ox_{P^n_{X/S}}$-linear isomorphism for all $n.$
\item $\varepsilon_0=\op{Id}_{\calE}.$
\item The following diagram is commutative for all $n\ge m$
$$\begin{tikzcd}
\Ox_{P^n_{X/S}}\otimes_{\Ox_X}\calE\ar{r}{\varepsilon_n}\ar{d} & \calE\otimes_{\Ox_X}\Ox_{P^n_{X/S}}\ar{d} \\
\Ox_{P^m_{X/S}}\otimes_{\Ox_X}\calE\ar{r}{\varepsilon_m} & \calE\otimes_{\Ox_X}\Ox_{P^m_{X/S}}
\end{tikzcd}$$
where the vertical arrows are the canonical homomorphisms.
\item $p_{13}^{n*}\varepsilon_n=p_{23}^{n*}\varepsilon_n\circ p_{12}^{n*}\varepsilon_n$ for all $n$ (see \ref{P11} for the notation). 
\end{enumerate}
\end{definition}

\begin{definition}
Let $\calE$ be an $\Ox_X$-module. A \emph{hyperstratification on $\calE$} is a morphism 
$$\varepsilon:\Ox_{P_{X/S}}\otimes_{\Ox_X}\calE\ra \calE\otimes_{\Ox_X}\Ox_{P_{X/S}}$$ satisfying the following conditions:
\begin{enumerate}
\item $\varepsilon$ is an $\Ox_{P_{X/S}}$-linear isomorphism.
\item $\varepsilon$ is equal to the identity $\op{Id}_{\calE}$ modulo the ideal $\calI\subset \Ox_{P_{X/S}}$ of $X$ in $P_{X/S}.$
\item $p_{13}^*\varepsilon=p_{23}^*\varepsilon\circ p_{12}^*\varepsilon.$
\end{enumerate}
\end{definition}

\begin{proposition}\label{prop39}
Let $\calE$ be a quasi-coherent $\Ox_X$-module and suppose that $X\ra S$ is log smooth. The following are equivalent:
\begin{enumerate}
\item A stratification $(\varepsilon_n)_{n\ge 0}$ on $\calE.$
\item A sequence of $\Ox_X$-linear morphisms $\theta_n:\calE\ra \calE\otimes_{\Ox_X}\Ox_{P_{X/S}^n}$ satisfying the following conditions:
\begin{itemize}
\item $\theta_0=\op{Id}_{\calE}.$
\item For any integer $n\ge 0,$ the morphism $\theta_n$ is equal to the composition
$$\calE\xrightarrow{\theta_{n+1}}\calE\otimes_{\Ox_X}\Ox_{P^{n+1}_{X/S}}\ra \calE\otimes_{\Ox_X}\Ox_{P^n_{X/S}}$$
where the second arrow is the canonical morphism.
\item For all integers $n,m\ge 0$
$$\begin{tikzcd}
\calE\ar{r}{\theta_{n+m}}\ar{d}{\theta_m} & \calE\otimes_{\Ox_X}\Ox_{P^{n+m}_{X/S}}\ar{d}{\op{Id}\otimes \delta^{n,m}}\\
\calE\otimes_{\Ox_X}\Ox_{P^m_{X/S}}\ar{r}{\theta_n\otimes \op{Id}}& \calE\otimes_{\Ox_X}\Ox_{P^n_{X/S}}\otimes_{\Ox_X}\Ox_{P^m_{X/S}} 
\end{tikzcd}$$
\end{itemize}
\item A sequence of compatible $\Ox_{P_{X/S}}$-linear morphisms
$$\nabla_n:\calD_{X/S}^n(\Ox_X,\Ox_X)\ra \calD_{X/S}^n(\calE,\calE)$$
satisfying for all integers $n,m\ge 0$ and $(f,g)\in \calD_{X/S}^n(\Ox_X,\Ox_X)\times \calD_{X/S}^m(\Ox_X,\Ox_X)$ $$\nabla_{n+m}(f\circ g)=\nabla_n(f)\circ \nabla_m(g)$$ 
and $\nabla_0(\op{Id}_{\Ox_X})=\op{Id}_{\calE}.$
\item An integrable connection $\nabla:\calE\ra \calE\otimes_{\Ox_X}\omega^1_{X/S}.$
\end{enumerate}
\end{proposition}

\begin{proof}
The proof is analogous to the proof of theorem 4.8 in \cite{Ogus78} so we just outline it here.

Suppose we are given a stratification $(\varepsilon_n)$ on $\calE.$ We define for any integer $n\ge 0$ the morphism $\theta_n$ as the composition
$$\calE\ra \Ox_{P^n_{X/S}}\otimes_{\Ox_X}\calE \xrightarrow{\varepsilon_n}\calE\otimes_{\Ox_X}\Ox_{P^n_{X/S}}$$
where the first morphism is the canonical morphism $x\mapsto 1\otimes x.$ We thus get the data of (2) which clearly satisfies the desired conditions.

Suppose we are given the morphisms $(\theta_n).$ We define $\varepsilon_n$ as the linearization of $\theta_n.$ We thus get the equivalence between (1) and (2).  

Suppose we are given a stratification $(\varepsilon_n).$ For any differential operator $f:\Ox_{P^n_{X/S}}\ra \Ox_X$ of order $\le n$ we define $\nabla_n(f)$ as the composition
$$\Ox_{P^n_{X/S}}\otimes_{\Ox_X}\calE \xrightarrow{\varepsilon_n}\calE\otimes_{\Ox_X}\Ox_{P^n_{X/S}}\xrightarrow{\op{Id}\otimes f}\calE.$$
The compatibility of $(\varepsilon_n)$ implies that of $\nabla_n.$ The fact that $\epsilon_0=\op{Id}_{\calE}$ implies that $\nabla_0(\op{Id}_{\Ox_X})=\op{Id}_{\calE}.$ To prove that $\nabla_n$ are compatible with the composition of differential operators, we consider the following diagram
$$\begin{tikzcd}
\Ox_{P^{n+m}_{X/S}}\otimes \calE \ar{r}{\delta\otimes \op{Id}}\ar{ddd}{\varepsilon_{n+m}} & \Ox_{P^n_{X/S}}\otimes \Ox_{P^m_{X/S}}\otimes \calE\ar{rr}{\op{Id}\otimes \varepsilon_m} & & \Ox_{P^n_{X/S}}\otimes\calE\otimes\Ox_{P^m_{X/S}}\ar{d}{\op{Id}\otimes\op{Id}\otimes g}\ar{dddll}{\varepsilon_n\otimes \op{Id}} \\
 & & &\Ox_{P^n_{X/S}}\otimes\calE\ar{d}{\varepsilon_n}\\
 & & & \calE\otimes\Ox_{P^n_{X/S}}\ar{d}{\op{Id}\otimes f}\\
\calE\otimes \Ox_{P^{n+m}_{X/S}}\ar{r}{\op{Id}\otimes \delta^{n,m}} & \calE\otimes \Ox_{P^n_{X/S}}\otimes\Ox_{P^m_{X/S}}\ar{r}{\op{Id}\otimes\op{Id}\otimes g} & \calE\otimes\Ox_{P^n_{X/S}}\ar[equal]{ru}\ar{r}{\op{Id}\otimes f} & \calE
\end{tikzcd}$$
The right-down composition is $\nabla_n(f)\circ \nabla_m(g).$ The down-right composition is $\nabla_{n+m}(f\circ g).$ The commutativity of the left part of the diagram is a consequence of the cocycle condition of the stratification $(\varepsilon_n).$ The commutativity of the right part of the diagram is a result of the following straightforward computation: for all local sections $z, z'$ and $x$ of $\Ox_{P^n_{X/S}},\ \Ox_{P^m_{X/S}}$ and $\calE$ respectively, we have
\begin{alignat*}
5
\varepsilon_n\circ (\op{Id}\otimes\op{Id}\otimes g)(z\otimes x\otimes z') &= \varepsilon(g(z')z\otimes x)\\
&= g(z')z\varepsilon(1\otimes x)\\
(\op{Id}\otimes\op{Id}\otimes g)(\varepsilon_n\otimes \op{Id})(z\otimes x\otimes z')&=(\op{Id}\otimes\op{id}\otimes g)(\varepsilon(z\otimes x)\otimes z')\\
&=g(z')\varepsilon(z\otimes x)\\
&=g(z')z\varepsilon(1\otimes x)
\end{alignat*}
It remains to prove that $\nabla_n$ is $\Ox_{P_{X/S}}$-linear: let $f,z,z'$ et $x$ be local sections of $\calD_{X/S}^n(\Ox_X,\Ox_X),$ $\Ox_{P_{X/S}},$ $\Ox_{P_{X/S}}$ and $\calE$ respectively. There exist unique sections $x_I$ of $\calE$ for all multi-indices $I\in\mathbb{N}^d$ such that
$$\varepsilon_n(1\otimes x)=\sum_{|I|\le n}x_I\otimes \eta^{[I]}.$$
Then we have
\begin{alignat*}
5
\nabla_n(z\cdot f)(z'\otimes x)&= (\op{Id}\otimes (z\cdot f))(z'\otimes x)\\
&= \sum_{|I|\le n}f(zz'\eta^{[I]})x_I\\
&=\nabla_n(f)(zz'\otimes x).
\end{alignat*} 
So $\nabla_n(z\cdot f)=z\cdot \nabla_n(f).$

Suppose we are given the morphisms $\nabla_n$ of (3). We want to prove the existence of morphisms $\varepsilon_n$ making the following diagram commutative
$$\begin{tikzcd}
\Ox_{P^n_{X/S}}\otimes_{\Ox_X}\calE\ar{r}{ev}\ar{d}{\varepsilon_n} & \mathscr{Hom}_{\Ox_{X}}(\mathscr{Hom}_{\Ox_{X}}(\Ox_{P^n_{X/S}}\otimes_{\Ox_X}\calE,\calE),\calE)\ar{d}{\mathscr{Hom}_{\Ox_{X}}(\nabla_n,\calE)}\\
\calE\otimes_{\Ox_X}\Ox_{P^n_{X/S}}\ar{r} & \mathscr{Hom}_{\Ox_{X}}(\mathscr{Hom}_{\Ox_{X}}(\Ox_{P^n_{X/S}},\Ox_X),\calE)
\end{tikzcd}$$
where $ev$ is the evaluation morphism and the lower morphism is given by
$$x\otimes z\mapsto (f\mapsto f(z)x).$$
This lower morphism is an isomorphism since $\Ox_{P^n_{X/S}}$ is locally free, hence the existence of the morphism $\varepsilon_n.$ Thus we get the equivalence between (1) and (3).

Suppose we are given an integrable connection $\nabla$ on $\calE$ and let us construct the morphisms $\nabla_n$ of (3). We start by defining $\nabla_1.$ After eventual localization, we can suppose that there exist sections $m_1,\hdots,m_d \in \Gamma(X,\calM_X)$ such that $(\op{dlog}m_i)_{1\le i\le d}$ form a basis of $\omega^1_{X/S}.$ We know then that $\Ox_{P_{X/S}}\cong \Ox_X\langle \eta_1,\hdots,\eta_d\rangle,$ where $\eta_i=\mu(m_i)-1$ (see \ref{P2}). Denote by $\partial_I$ the dual basis of $\eta^{[I]}$ and let $\varepsilon_1:\Ox_{P^1_{X/S}}\otimes_{\Ox_X}\calE\ra \calE\otimes_{\Ox_X}\Ox_{P^1_{X/S}}$ be the $\Ox_{P_{X/S}}$-linear morphism defined by
$$\varepsilon_1(1\otimes x)=\nabla(x)+x\otimes 1.$$ 
Then for any differential operator $f:\Ox_X\ra \Ox_X$ of order $\le 1,$ let $\nabla_1(f)$ be the composition
$$\Ox_{P^1_{X/S}}\otimes_{\Ox_X}\calE\xrightarrow{\varepsilon_1} \calE\otimes_{\Ox_X}\Ox_{P^1_{X/S}}\xrightarrow{\op{Id}\otimes f}\calE.$$
The connection $\nabla$ is integrable so the morphisms $\nabla(\partial_{\epsilon_i})$ pairwise commute and we can define $\nabla_n(\partial_N)$ for any multi-index $N=(n_1,\hdots,n_d)\in\mathbb{N}^d$ of length $|N|\le n$ by
$$\nabla_n(\partial_N)=\prod_{i=1}^d\prod_{j=0}^{n_i-1}(\nabla_1(\partial_{i})-j),$$
see \ref{cor18} and \ref{cor25}.
The morphisms $\nabla_n$ are by definition compatible and satisfy $\nabla_{n+m}(f\circ g)=\nabla_n(f)\circ \nabla_m(g)$ and $\nabla_0(\op{Id}_{\Ox_X})=\op{Id}_{\calE}.$ All that remains is to prove that $\nabla_n$ is $\Ox_{P_{X/S}}$-linear i.e. that
$$\nabla_n(z\cdot f)=z\cdot \nabla_n(f)$$
for any differential operator $f\in \calD_{X/S}^n(\Ox_X,\Ox_X)$ and any local section $z$ of $\Ox_{P_{X/S}}.$
Since $\nabla_n$ is $\Ox_X$-linear, it is sufficient to prove the equality for $z=\eta^{[r]}$ and $f=\partial_N.$ We proceed by induction on $|N|:$

If $|N|=0$, the result is clear.

Suppose that the result is true for a multi-index $N\in \mathbb{N}^d$ and let us prove it is true for $N+\varepsilon_i,$ where $1\le i\le d.$ Let $z=\eta^{[r]}$ and $f=\partial_{N+\varepsilon_i}.$
We know that $\partial_{N+\varepsilon_i}=\partial_{\epsilon_i}\circ \partial_N-N_i\partial_N$ (\ref{prop17}) so for any local section $x$ of $\calE$
\begin{alignat*}
5
z\cdot \nabla(f)(x)&= \eta^{[r]}\cdot \nabla(\partial_{\epsilon_i})\circ \nabla(\partial_N)(x)-N_i\eta^{[r]}\cdot \nabla(\partial_N)(x).
\end{alignat*}
On the other hand
\begin{alignat*}
5
z\cdot \nabla(\partial_{\epsilon_i})\circ \nabla(\partial_N)(x)&= \nabla(\partial_{\epsilon_i})\circ (\op{Id}\otimes \nabla(\partial_N))\circ (P\otimes\op{Id})(\eta^{[r]}\otimes x)\\
&=\sum_{a+b+c=r}c!\begin{pmatrix}a+c\\c \end{pmatrix}\begin{pmatrix}b+c\\c \end{pmatrix}\nabla(\partial_{\epsilon_i})\circ (\op{Id}\otimes \nabla(\partial_N))\left (\eta^{[a+c]}\otimes \eta^{[b+c]}\otimes x\right )\\
&=\sum_{a+b+c=r}c!\begin{pmatrix}a+c\\c \end{pmatrix}\begin{pmatrix}b+c\\c \end{pmatrix}\nabla(\partial_{\epsilon_i})\left (\eta^{[a+c]}\otimes \nabla(\partial_N)(\eta^{[b+c]}\otimes x)\right )\\
&=\sum_{a+b+c=r}c!\begin{pmatrix}a+c\\c \end{pmatrix}\begin{pmatrix}b+c\\c \end{pmatrix}\nabla(\eta^{[a+c]}\cdot\partial_{\epsilon_i})\left (1\otimes \nabla(\eta^{[b+c]}\cdot \partial_N)(1\otimes x)\right )
\end{alignat*}
where the second line follows from \ref{eq171} and the last line follows from the induction hypothesis.
By lemma \ref{lem19}
$$\eta^{[a+c]}\cdot \partial_{\epsilon_i}=\sum_{I+J=\varepsilon_i}\begin{pmatrix} \varepsilon_i\\ I\end{pmatrix} \partial_I\left (\eta^{[a+c]}\right )\partial_J=\partial_0\left (\eta^{[a+c]}\right )\partial_{\epsilon_i}+\partial_{\epsilon_i}\left (\eta^{[a+c]}\right )\partial_0.$$
We deduce that
\begin{alignat*}
5
z\cdot \nabla(\partial_{\epsilon_i})\circ \nabla(\partial_N)(x) &= \nabla(\partial_{\epsilon_i})\left (1\otimes \nabla(\eta^{[r]}\cdot \partial_N)(1\otimes x)\right )\\
&\ \ +r_i\nabla(\partial_0)(1\otimes \nabla(\eta^{[r]}\cdot \partial_N)(1\otimes x))\\
&\ \ +\nabla(\partial_0) (1\otimes \nabla(\eta^{[r-\varepsilon_i]}\cdot \partial_N)(1\otimes x)).
\end{alignat*}
Since $\eta^{[r]}\cdot \partial_N=\begin{pmatrix}N\\r \end{pmatrix}\partial_{N-r}$ and
$\eta^{[r-\varepsilon_i]}\cdot \partial_N=\begin{pmatrix}N\\r-\varepsilon_i \end{pmatrix}\partial_{N-r+\varepsilon_i},$
we get
\begin{alignat*}{2}
z\cdot \nabla(\partial_{\epsilon_i})\circ \nabla(\partial_N)(x) &= \begin{pmatrix}N\\r \end{pmatrix}\nabla(\partial_{\epsilon_i})(1\otimes \nabla(\partial_{N-r}))\\
&+r_i\begin{pmatrix}N\\r \end{pmatrix}\nabla(\partial_{N-r})(x)+\begin{pmatrix}N\\r-\varepsilon_i \end{pmatrix}\nabla(\partial_{N-r+\varepsilon_i})(x).
\end{alignat*}
We also have
$$\nabla(z\cdot \partial_{N+\varepsilon_i})(x)=\begin{pmatrix}N+\varepsilon_i\\r \end{pmatrix}\nabla(\partial_{N-r+\varepsilon_i})(x).$$
We conclude that
\begin{alignat*}
5
&z\cdot \nabla(f)-\nabla(z\cdot f)\\
=&z\cdot\nabla(\partial_{\epsilon_i})\circ \nabla(\partial_N)-N_iz\cdot\nabla(\partial_N)-\nabla(z\cdot f)\\
=&\begin{pmatrix}N\\r \end{pmatrix}\nabla(\partial_{\epsilon_i}\circ \partial_{N-r})+(r_i-N_i)\begin{pmatrix}N\\r \end{pmatrix}\nabla(\partial_{N-r})\\
&+\begin{pmatrix}N\\r-\varepsilon_i \end{pmatrix}\nabla(\partial_{N-r+\varepsilon_i})-\begin{pmatrix}N+\varepsilon_i\\r \end{pmatrix}\nabla(\partial_{N-r+\varepsilon_i})\\
=&\begin{pmatrix}N\\r \end{pmatrix}\left (\nabla(\partial_{\epsilon_i}\circ \partial_{N-r}) +(r_i-N_i)\nabla(\partial_{N-r})-\nabla(\partial_{N-r+\varepsilon_i})\right )\\
=&\ 0
\end{alignat*}
where the last line results from the formula $\partial_{\epsilon_i}\circ \partial_I=\partial_{I+\varepsilon_i}+I_i\partial_I$ (\ref{prop17}).

\end{proof}

\begin{parag}\label{Shcomment}
Just like the classical case, we want a necessary and sufficient condition for a connection $\nabla_1:\calD_{X/S}^1(\Ox_X,\Ox_X)\ra \calD_{X/S}^1(\calE,\calE)$ to extend to a hyperstratification $\varepsilon:\Ox_{P_{X/S}}\otimes_{\Ox_X}\calE\ra \calE\otimes_{\Ox_X}\Ox_{P_{X/S}}.$ Suppose that such an extension $\varepsilon$ exists and denote $\nabla:\calD_{X/S}(\Ox_X,\Ox_X)\ra \calD_{X/S}(\calE,\calE)$ the corresponding morphism.
Let $x$ be a local section of $\calE.$ Then
$$\varepsilon(1\otimes x)=\sum_{N\in\mathbb{N}^d}\nabla(\partial_N)(1\otimes x)\otimes \eta^{[N]}.$$
The sum is finite so $\nabla(\partial_N)(1\otimes x)=0$ for $|N|$ sufficiently large, which is equivalent to 
$$\prod_{i=1}^d\prod_{j=0}^{n_i-1}\left (\nabla_1(\partial_{\epsilon_i})-j\right )(1\otimes x)=0$$ for $N=(n_1,\hdots,n_d)$ of sufficiently large length. This motivates the following definition:
\end{parag}

\begin{definition}\label{nil}
Let $X\ra S$ be a smooth morphism of logarithmic schemes and consider local coordinates $m_1,\hdots,m_d$ of $X$ with respect to $S.$ For all $1\le i\le d,$ denote by $\partial_i$ the logarithmic derivation associated to $m_i.$ Let $\nabla$ be an integrable connection on an $\Ox_X$-module $\calE.$ By \ref{cor25}, the morphisms $\nabla(\partial_i)$ and $\nabla(\partial_j)$ commute for all $1\le i,j\le d.$ Then we can define, for a multi-index $I=(I_1,\hdots,I_d)\in\mathbb{N}^d,$ a morphism $\nabla(\partial_I)$ as follows:
$$\nabla(\partial_I)=\prod_{k=1}^d\prod_{l=0}^{I_k-1}(\nabla(\partial_k)-l\op{Id}_{\calE}).$$
We say that $\nabla$ is \emph{quasi-nilpotent} if for any open subset $U\subset X$ and for any $x\in \Gamma(U,\calE),$ there exists, locally on $U,$ an integer $N$ such that $\nabla(\partial_I)(x)=0$ for all multi-index $I\in\mathbb{N}^d$ satisfying $|I|\ge N.$ We denote by $\boldsymbol{\op{MIC}^{qn}(X/S)}$ the full subcategory of $\boldsymbol{\op{MIC}(X/S)}$ consisting of $\Ox_X$-modules equipped with a quasi-nilpotent integrable connection. By \ref{Shcomment}, an integrable connection is quasi-nilpotent if and only if the corresponding stratification extends to hyperstratification. We see then that this definition is actually independant of the choice of the local system of coordinates, just as in the classical smooth case.
\end{definition}

We end this section with the following proposition that we will later need:
\begin{proposition}\label{prop412}
Let $\nabla:\calE\ra \calE\otimes_{\Ox_X}\omega^1_{X/S}$ be an integrable connection. Denote by $\psi_{\nabla}$ its $p$-curvature and consider the differential operators $\partial_I$ defined in (\ref{prop17}).
\begin{enumerate}
\item $\partial_{\epsilon_i}^p-\partial_{\epsilon_i}^{(p)}=\partial_{\epsilon_i}^p-\partial_{\epsilon_i}=\partial_{p\epsilon_i}.$
\item $\nabla(\partial_{p\epsilon_i})=\psi_{\nabla}(\partial_{\epsilon_i}).$
\item Let $I=(I_1,\hdots,I_d)\in \N^d$ be a multi-index. If $\psi_{\nabla}=0$ and there exists $1\le i\le d$ such that $I_i\ge p,$ then
$$\nabla(\partial_I)=0.$$
\end{enumerate}
\end{proposition}

\begin{proof}
We have the following equality of polynomials in $\mathbb{F}_p[x]:$
$$x(x-1)\hdots (x-p+1)=x^p-x.$$
By (\ref{cor18}), we get
$$\partial_{\epsilon_i}^p-\partial_{\epsilon_i}=\partial_{\epsilon_i}(\partial_{\epsilon_i}-1)\hdots (\partial_{\epsilon_i}-p+1)=\partial_{p\epsilon_i}.$$
Then, by \cite{Ogus94} remark 1.2.2, $\partial_{\epsilon_i}^{(p)}=\partial_{\epsilon_i}$ and so
$$\nabla(\partial_{p\epsilon_i})=\nabla(\partial_{\epsilon_i}^p)-\nabla(\partial_{\epsilon_i}^{(p)})=\psi_{\nabla}(\partial_{\epsilon_i}).$$
For the third point, by (\ref{cor18}), we have
\begin{alignat*}{2}
\nabla(\partial_I) &= \prod_{j\neq i}\nabla(\partial_{I_j\epsilon_j})\circ \nabla(\partial_{I_i\epsilon_i})\\
&= \prod_{j\neq i}\nabla(\partial_{I_j\epsilon_j})\circ \prod_{k=p}^{n-1}(\nabla(\partial_{\epsilon_i})-k)\circ \nabla(\partial_{p\epsilon_i}) \\
&= 0.
\end{alignat*}
\end{proof}

\newpage

\section{Frames on logarithmic schemes}

\begin{parag}\label{parag42}
We will review in this section the notion of a frame on a logarithmic scheme, that was introduced by Kato and Saito in (\cite{Saito04} \S 4). A monoid $P$ is said to be \emph{integral} if the canonical morphism $P\ra P^{gp},$ from $P$ to its associated group $P^{gp},$ is injective. A monoid $P$ is said to be \emph{saturated} if it is integral and for every $x\in P^{gp},$ if there exists a positive integer $n$ such that $nx\in P,$ then $x\in P.$
We will often consider finitely-generated and saturated monoids and we will refer to them as \emph{fs monoids}.
A logarithmic structure is finitely-generated (resp. integral, resp. saturated) if, étale locally, it has a chart given by a finitely-generated (resp. integral, resp. saturated) monoid. In this paper, we often consider finitely generated and saturated logarithmic schemes, we refer to them by fs logarithmic schemes and we denote the category of fs logarithmic schemes by $\boldsymbol{\op{L}}.$ We denote by $\boldsymbol{\widehat{\op{L}}}$ the category of presheaves of sets on $\boldsymbol{\op{L}}.$ For morphisms of fs logarithmic schemes $X\ra S$ and $Y\ra S,$ we denote by $X\times_S^{\op{log}}Y$ the fiber product of $X\ra S$ and $Y\ra S$ in the category $\boldsymbol{\op{L}}$ of fs logarithmic schemes. Recall that the Yoneda embedding $\boldsymbol{\op{L}}\ra \boldsymbol{\widehat{\op{L}}}$ commutes with representable projective limits and in particular with fiber products. We keep the notation $X\times_SY$ for the fiber product of $X\ra S$ and $Y\ra S$ in the category of logarithmic schemes.
\end{parag}

\begin{parag}\label{parag43}
Let $P$ be an fs monoid. We define a presheaf of sets $[P]$ on the category $\boldsymbol{\op{L}}$ of fs logarithmic schemes as follows :
\begin{equation}
[P]:\begin{array}[t]{clc} \boldsymbol{\op{L}} & \ra & \boldsymbol{\op{Sets}}\\ T & \mapsto & \op{Hom}\left (P,\Gamma(T,\ov{\calM}_T)\right )\end{array}
\end{equation}
where $\ov{\calM}_T$ is defined in \ref{Not5}. Consider the logarithmic scheme $A[P]$ defined in \ref{Not4}. We have a canonical morphism of presheaves $A[P]\ra [P]$ induced by the tautological map $P\ra \Gamma(A[P],\calM_{A[P]}).$
If $X$ is an fs logarithmic scheme, a morphism $X\ra [P]$ of presheaves of sets on $\boldsymbol{\op{L}},$ is said to be \emph{strict} if, for every geometric point $\ov{x}$ of $X,$ there exists an étale neighborhood $U$ of $\ov{x}$ such that the induced morphism $U\ra [P]$ factors
$$\begin{tikzcd}
U\ar{r}\ar{d} & \left [P\right ]\\
A\left [P\right ]\ar{ur} &
\end{tikzcd}$$
where $U\ra A\left [P\right ]$ is a strict morphism of logarithmic schemes and $A\left [P\right ]\ra [P]$ is the canonical morphism. A \emph{frame} on a logarithmic scheme $X$ is a strict morphism $X\ra [P].$ An fs logarithmic scheme $X$ equipped with a frame $X\ra [P]$ will be called a \emph{framed logarithmic scheme} and will be denoted by $(X,P).$ A morphism of framed logarithmic schemes $(X,Q)\ra (S,P)$ is a morphism of logarithmic schemes $X\ra Y$ and a morphism of monoids $P\ra Q,$ such that the diagram
$$\begin{tikzcd}
X\ar{r}\ar{d} & Y\ar{d} \\
\left [Q\right ]\ar{r} & \left [P\right ]
\end{tikzcd}$$
is commutative, where $[Q]\ra [P]$ is the morphism of presheaves induced by $P\ra Q.$

Finally, given an fs logarithmic scheme $X$ equipped with a map $X\ra [Q]$ and a morphism of monoids $Q\ra P,$ we denote by $X\times_{[Q]}^{\op{log}}[P]$ the presheaf of $\boldsymbol{\op{L}}$ defined for any fs logarithmic scheme $T$ by
$$(X\times_{[Q]}^{\op{log}}[P])(T)=X(T) \times _{[Q](T)} [P](T).$$
Also, given logarithmic schemes $X,\ Y$ and $S,$ a monoid $Q$ and morphisms $X\ra S\times [Q]$ and $Y\ra S\times [Q],$ we denote by $X\times_{S,[Q]}^{\op{log}}Y$ the presheaf of $\boldsymbol{\widehat{\op{L}}}$ defined for any fs logarithmic scheme $T$ by
$$(X\times_{S,[Q]}^{\op{log}}Y)(T)= X(T)\times_{S(T)\times [Q](T)}Y(T).$$ 
We now recall results about the representability of these presheaves. We start by recalling the following proposition, which is a slight improvement of (\cite{Saito04} 4.2.1). 
\end{parag}

\begin{proposition}[\cite{Saito04} 4.2.1]\label{prop49}
Let $P$ and $Q$ be fs monoids, $\theta:Q\ra P$ a morphism of monoids such that $\theta^{gp}$ is surjective, $X$ a logarithmic scheme and $X\ra [Q]$ a frame. Then $X\times_{[Q]}^{\op{log}}[P]$ is representable by a logarithmic scheme which affine and log étale over $X.$ In addition, if $X\ra A[Q]$ lifts $X\ra [Q],$ then $X\times_{[Q]}^{\op{log}}[P]$ is representable by $X\times^{\op{log}}_{A[Q]}A[\tilde{P}],$ where $\tilde{P}$ is the inverse image of $P$ by $\theta^{gp}.$
\end{proposition}

\begin{proposition}[\cite{Saito04} 4.2.3]\label{rep1}
Let $Q$ be an fs monoid, $X,$ $Y$ and $S$ fs logarithmic schemes and $X\ra S\times [Q],\ Y\ra S\times [Q]$ morphisms of presheaves. Then, the presheaf $X\times_{S,[Q]}^{\op{log}}Y$ is representable by the logarithmic scheme $\left (X\times_S^{\op{log}}Y\right )\times_{[Q \oplus Q]}^{\op{log}}[Q],$ which is log étale over $X\times_S^{\op{log}}Y.$
\end{proposition}

\begin{proposition}[\cite{Saito04} 4.1.6]\label{framelift}
Let $f:(X,Q) \ra (S,P)$ be a morphism of framed fs logarithmic schemes. Etale locally on $X$ and $S,$ there exists a chart $A[M] \ra A[P]$ of $f$ lifting the diagram
$$
\begin{tikzcd}
X \ar{r} \ar{d} & S \ar{d} \\
\left [Q \right ] \ar{r} & \left [P \right ].
\end{tikzcd}
$$
\end{proposition}

\begin{proposition}\label{etaleframelift}
Let $f:(X,Q) \ra (S,P)$ be a log smooth morphism of framed fs logarithmic schemes. Suppose that $S\ra [P]$ lifts to a chart $S\ra A[P].$ Etale locally on $X,$ there exists a chart $A[M] \ra A[P]$ of $f$ lifting the diagram
$$
\begin{tikzcd}
X \ar{r} \ar{d} & S \ar{d} \\
\left [Q \right ] \ar{r} & \left [P \right ],
\end{tikzcd}
$$
and satisfying the following conditions:
\begin{enumerate}
\item $P^{gp} \ra M^{gp}$ is injective and the torsion subgroup of its cokernel is of finite order invertible in $\Ox_X.$
\item The morphism $X \ra S\times_{A[P]}A[M],$ induced by $f$ and the chart $X\ra A[M],$ is étale and strict.
\end{enumerate}
\end{proposition}

\begin{proof}
Let $\ov{x} \ra X$ be a geometric point and $\ov{s}=f(\ov{x}).$ By (\cite{Ogus2018} IV 3.3.1), after restricting to an étale neighborhood of $\ov{x},$ there exists a chart $P\ra M$ of $f$ satisfying the conditions (1) and (2) and such that the chart $X \ra A[M]$ is exact at $\ov{x}.$ By (\cite{Ogus2018} II 2.3.1), the morphism $\gamma:\ov{M} \ra \ov{\calM}_{X,\ov{x}}$ is an isomorphism. We have the commutative diagram
\begin{equation}\label{diagtak1}
\begin{tikzcd}
X \ar{r} \ar{d} & S \ar{d} \\
A[M] \ar{r} & A[P].
\end{tikzcd}
\end{equation}
It is then sufficient to prove the existence of a morphism $A[M] \ra [Q]$ such that the diagram
\begin{equation}\label{diagtak2}
\begin{tikzcd}
A[M] \ar{r} \ar{d} & A[P] \ar{d} \\
\left [Q\right ] \ar{r} & \left [ P\right ]
\end{tikzcd}
\end{equation}
is commutative.
The chart $S\ra A[P]$ induces a morphism $\beta:\ov{P} \ra \ov{\calM}_{S,\ov{s}}.$ The frame $X\ra [Q]$ corresponds to a morphism $Q\ra \Gamma(X,\ov{\calM}_X).$ Consider the composition
$$
\delta:Q \ra \Gamma(X,\ov{\calM}_X) \ra \ov{\calM}_{X,\ov{x}} \xrightarrow{\gamma^{-1}} \ov{M}.
$$
The outer rectangle and the lower square of the following diagram
$$
\begin{tikzcd}
P \ar{r} \ar{d} & Q \ar{d}{\delta} \\
\ov{P} \ar{r} \ar{d} & \ov{M} \ar{d}{\gamma} \\
\ov{\calM}_{S,\ov{s}} \ar{r}{\ov{f^{\flat}_{\ov{x}}}} & \ov{\calM}_{X,\ov{x}}
\end{tikzcd}
$$
are commutative. Since $\gamma$ is an isomorphism,
the upper square is commutative. The result follows.
\end{proof}

\begin{proposition}[\cite{Saito04} 4.2.5.2]\label{rep3}
Let $X\ra S$ and $Y\ra S$ be morphisms of fs logarithmic schemes and $\theta:P\ra Q$ a morphism of fs monoids that fit into a commutative diagram of $\boldsymbol{\widehat{\op{L}}}:$
$$
\begin{tikzcd}
X \ar{r}\ar{d} & S \ar{d} & Y\ar{l}\ar{d} \\
\left [Q\right ] \ar{r} & \left [P\right ] & \left [Q\right ]\ar{l} 
\end{tikzcd}
$$
where both of the morphisms $[Q] \ra [P]$ are induced by $\theta.$
If $X\ra [Q]$ and $S\ra [P]$ are frames then the canonical projection $X\times_{S,[Q]}^{\op{log}}Y\ra Y$ is strict.
\end{proposition}

\begin{proposition}[\cite{Saito04} 4.2.6]\label{rep2}
Let $\theta : P\ra Q$ be a morphism of fs monoids and $X,$ $Y$ and $S$ be fs logarithmic schemes that fit into a commutative diagram 
$$\begin{tikzcd}
X \ar{r}\ar{d} & S\ar{d} & Y\ar{l}\ar{d} \\
A\left [Q\right ] \ar{r} & A\left [P\right ] & A\left [Q\right ]\ar{l}
\end{tikzcd}$$
where both of the morphisms $A[Q]\ra A[P]$ are induced by the morphism $\theta.$
If the vertical arrows are strict, then the presheaf $X\times_{S,[Q]}^{\op{log}}Y$ is representable by the scheme $X\times_SY\times_{A[Q\oplus Q]}A[M]$ equipped with the logarithmic structure pullback of that of $A[M],$ where $M$ is the inverse image of $Q$ by the canonical morphism $Q^{gp}\oplus_{p^{gp}}Q^{gp}\ra Q^{gp},\ (x,y)\mapsto x+y.$
\end{proposition}

\begin{proposition}[\cite{Saito04} 4.2.8]\label{prop45}
Let $(X,Q)\ra (S,P)$ be a morphism of framed fs logarithmic schemes such that $P$ and $Q$ are fs monoids and consider the factorization $$X\xrightarrow{\Delta}X\times_{S,[Q]}^{\op{log}}X\xrightarrow{g}X\times_S^{\op{log}}X$$
of the diagonal morphism $X\ra X\times_S^{\op{log}}X.$ Then $g$ is log étale and $\Delta$ is an exact immersion. Furthermore, if $X\ra U \ra X\times_{S,[Q]}^{\op{log}}X$ is a factorization of $\Delta$ into a closed immersion and an open one and if $\calI$ is the ideal of $X\ra U,$ then the canonical morphism of $\Ox_X$-modules
\begin{equation}\label{iso}
\Delta^{-1}\left (\calI/\calI^2\right )\ra \omega^1_{X/S}
\end{equation}
is an isomorphism.
\end{proposition}

\begin{proof}
Set $Y=X\times_{S,[Q]}^{\op{log}}X$ and denote by $p_1,p_2:X\times_S^{\op{log}}X\ra X$ and $q_1,q_2:Y\ra X$ the canonical projections. By \ref{rep1}, the morphism $g$ is log étale and by \ref{rep3} the immersion $\Delta$ is strict and hence exact. Consider $Y$ (resp. $X\times_S^{\op{log}}X$) as a scheme over $X$ by the first projection $q_1$ (resp. $p_1$). The conormal exact sequence associated to $X\xrightarrow{\Delta} Y\ra X$ yields a canonical isomorphism 
$$\Delta^{-1}\left (\calI/\calI^2\right )\xrightarrow{\sim}\Delta^*\omega^1_{Y/X}.$$
Since $g$ is log étale and $p_2$ induces an isomorphism $p_2^*\omega^1_{X/S}\xrightarrow{\sim}\omega^1_{X\times_S^{\op{log}}X/X},$ the morphism $q_2$ induces an isomorphism 
$$q_2^*\omega^1_{X/S}\xrightarrow{\sim}\omega^1_{Y/X}.$$
It follows that we have an isomorphism 
\begin{equation}
\Delta^{-1}(\calI/\calI^2)\xrightarrow{\sim}\omega^1_{X/S}.
\end{equation}
\end{proof}

\begin{parag}\label{parag46}
We keep the same assumptions of \ref{prop45} and we determine the image of $\op{dlog}m,$ for a local section $m$ of $\calM_X,$ by the inverse of (\ref{iso}).
Keep the same notations as in the proof of \ref{prop45}. Let $m$ be a local section of $\calM_X.$ The local section $\op{dlog}m$ of $\omega^1_{X/S}$ is sent by the isomorphism $\omega^1_{X/S}\xrightarrow{\sim}\Delta^*\omega^1_{Y/X}$ to $\op{dlog}q_2^{\flat}m.$
The exact immersion $\Delta$ gives rise to the following exact sequence of abelian groups
$$0\ra \Delta^{-1}(1+\calI)\ra \Delta^{-1}\calM_Y^{gp}\ra \calM_X^{gp}\ra 0.$$
Let $\ov{x}\ra X$ be a geometric point and $\ov{y}\ra Y$ its image. By the previous exact sequence, there exists a unique $\mu(m)\in 1+\calI_{\ov{y}}$ such that
$$\alpha_{Y,\ov{y}}^{-1}(\mu(m))=q^{\flat}_{2,\ov{y}}(m)-q_{1,\ov{y}}^{\flat}(m).$$
In $\omega^1_{Y/X,\ov{y}},$ we have
\begin{alignat*}{2}
\op{dlog}q_{2,\ov{y}}^{\flat}m &= \op{dlog}(q_{2,\ov{y}}^{\flat}m-q_{1,\ov{y}}^{\flat}m)\\
&= \op{dlog} \alpha_{Y,\ov{y}}^{-1}(\mu(m))\\
&= \mu(m)^{-1}d\mu(m).
\end{alignat*}
It follows that in $(\Delta^*\omega^1_{Y/X})_{\ov{x}}$ we have
$$1\otimes \op{dlog}q_{2,\ov{y}}^{\flat}m = d(\mu(m)-1),$$
so the section $\mu(m)-1$ is sent, by the isomorphism (\ref{iso}), to the section $\op{dlog}m.$
\end{parag}

\begin{remark}
Under the same assumptions of \ref{prop45}, we can prove without much difficulty that the same conclusions hold for the $r$-fold diagonal immersion $X\ra X_S^{r+1},$ where $X_S^{r+1}$ denotes the fiber product over $S,$ in $\boldsymbol{\op{L}},$ of $X$ by itself $r+1$ times. More precisely, denote by $X_{S,[Q]}^{r+1}$ the fiber product $X\times^{\op{log}}_{S,[Q]}X\times^{\op{log}}_{S,[Q]}\hdots \times^{\op{log}}_{S,[Q]}X$ ($r+1$ times). Then the canonical morphism $X^{r+1}_{S,[Q]}\ra X^{r+1}_S$ is étale and the canonical morphism $X\ra X^{r+1}_{S,[Q]}$ is a strict immersion.
\end{remark}

\begin{parag}\label{modQQ}
Let $f:X\ra S$ be a morphism of fs logarithmic schemes of characteristic $p.$ Recall the diagram of the exact relative Frobenius $F:X\ra X'$ \eqref{diag51}:
\begin{equation}
\begin{tikzcd}
X\ar{r}{F}\ar{dr}\ar[bend right=30]{rdd} & X'\ar{d}{G}\ar{dr}{\pi} & \\
 & X''\ar{d}\ar{r} & X\ar{d}{f} \\
 & S \ar{r}{F_S} & S
\end{tikzcd}
\end{equation}
Suppose that $f$ has a chart $\theta:P\ra Q$ and let
$$F_P:P \ra P,\ x\mapsto px$$
be the Frobenius morphism of $P.$ Consider the monoid $Q''=(Q\oplus_{P,F_P}P)^{int}$ and $Q'$ the inverse image of $Q$ by
$$v:Q''^{gp} \ra Q^{gp},\ (x,y) \mapsto px+\theta^{gp}(y).$$
Recall that the canonical morphism $X''\ra A_1[Q'']$ is a chart, $X'=X''\times_{A_1[Q'']}A_1[Q']$ and $v:Q' \ra Q$ is a chart of the exact relative Frobenius $F:X\ra X'.$

We now prove that if $X$ and $S$ are equipped with frames $X\ra [Q]$ and $S\ra [P]$ and $X\ra S$ is a morphism of framed logarithmic schemes which is of Cartier type (\cite{Kat89} (4.8)), then there exists a canonical frame $X'\ra [Q'].$
\end{parag}

\begin{theorem}\label{thm410}
Let $\theta:P\ra Q$ be a morphism of fs monoids and $(f,\theta):(X,Q)\ra (S,P)$ a morphism of framed logarithmic schemes of characteristic $p.$ Recall the exact relative Frobenius diagram \eqref{diag51} and the monoids $Q'$ and $Q''$ and the morphism $F_P:P\ra P$ defined in \ref{modQQ}.
Then there exists a canonical frame $X'\ra [Q'].$
\end{theorem}

\begin{proof}
Denote by $F:X\ra X'$ the exact relative Frobenius. 
The canonical projections $X''\ra X$ and the canonical morphism $X''\ra S$ induce morphisms
$$\Gamma \left (X,\ov{\calM}_X \right )\ra \Gamma \left (X'',\ov{\calM}_{X''} \right )\leftarrow \Gamma \left (S,\ov{\calM}_S \right ).$$
The frames $X\ra [Q]$ and $S\ra [P]$ correspond to morphisms
$$
Q\ra \Gamma(X,\ov{\calM}_X),\ P\ra \Gamma \left (S,\ov{\calM}_S \right ).
$$
We obtain a morphism
$$Q''=(Q\oplus_{P,F_P}P)^{int}\ra \Gamma \left (X'',\ov{\calM}_{X''} \right ),$$
and hence a morphism
$$X''\ra [Q''].$$
The exactness of $F^{\flat}:F^{-1}\calM_{X'} \ra \calM_X$ yields a cartesian square
$$
\begin{tikzcd}
\Gamma \left (X',\ov{\calM}_{X'} \right ) \ar{r} \ar[hook]{d} & \Gamma \left (X,\ov{\calM}_{X} \right ) \ar[hook]{d} \\
\Gamma \left (X',\ov{\calM}_{X'}^{gp} \right ) \ar{r} & \Gamma \left (X,\ov{\calM}_{X}^{gp} \right ).
\end{tikzcd}
$$
Similarly, the canonical morphism $Q' \ra Q$ is exact and hence the square
$$
\begin{tikzcd}
Q' \ar{r} \ar[hook]{d} & Q \ar[hook]{d} \\
Q'^{gp} \ar{r} & Q^{gp}
\end{tikzcd}
$$
is cartesian. Consider the composition
$$
Q'^{gp}=Q''^{gp} \ra \Gamma \left (X'',\ov{\calM}_{X''} \right ) \ra \Gamma \left (X',\ov{\calM}_{X'} \right ),
$$
where the second arrow is induced by the canonical morphism $X' \ra X''.$
There exists a unique morphism $Q' \ra \Gamma \left (X',\ov{\calM}_{X'} \right )$ fitting into the commutative diagram
$$
\begin{tikzcd}
Q' \ar{rr} \ar[hook]{dd} \ar{dr} & & Q \ar{dr} \ar[hook,dashed]{dd} & \\
 & \Gamma \left (X',\ov{\calM}_{X'} \right ) \ar{rr} \ar[hook]{dd} & & \Gamma \left (X,\ov{\calM}_{X} \right ) \ar[hook]{dd} \\
Q'^{gp} \ar[dashed]{rr} \ar{dr} & & Q^{gp} \ar[dashed]{dr} & \\
 & \Gamma \left (X',\ov{\calM}_{X'}^{gp} \right ) \ar{rr} & & \Gamma \left (X,\ov{\calM}_{X}^{gp} \right ).
\end{tikzcd}
$$
This morphism $Q' \ra \Gamma \left (X',\ov{\calM}_{X'} \right )$ yields the desired frame.
\end{proof}

\begin{corollaire}\label{cor411}
Let $\theta:P\ra Q$ be a morphism of fs monoids and $(f,\theta):(X,Q)\ra (S,P)$ a morphism of framed logarithmic schemes of characteristic $p.$ Recall the exact relative Frobenius diagram \eqref{diag51} and the monoid $Q'$ and the morphism $F_P:P\ra P$ defined in \ref{modQQ}. Let
\begin{equation}
\pi_{Q/P}:Q \ra Q',\ x\mapsto (x,0)
\end{equation}
and
\begin{equation}\label{FQP}
F_{Q/P}:Q' \ra Q,\ (x,y)\mapsto px+\theta^{gp}(y).
\end{equation}
Then $(F,F_{Q/P}):(X,Q)\ra (X',Q')$ and $(\pi,\pi_{Q/P}):(X',Q') \ra (X,Q)$ are morphisms of framed logarithmic schemes.
\end{corollaire}

\subsection*{Embedding $\calT_{X'/S}$ into $\calD_{X/S}$}

\begin{parag}
Let $f:(X,Q) \ra (S,P)$ be a log smooth morphism of framed fs logarithmic schemes of characteristic $p.$ Recall that the diagonal morphism $X \ra X\times^{\op{log}}_SX$ factors uniquely as
$$
\begin{tikzcd}
 & Y \ar{d}{g} \\
X \ar{ur}{\Delta} \ar{r} & X\times_S^{\op{log}}X
\end{tikzcd}
$$
where $\Delta:X \ra Y$ is an exact immersion of fs logarithmic schemes and $g$ is log étale. Furthermore, the normal sheaf of $\Delta$ is canonically isomorphic to $\omega^1_{X/S}.$ Let $P_{X/S}$ be the PD envelope of $\Delta$ and equip it with the logarithmic structure pullback of that of $Y.$ Note that this is the same as the logarithmic PD envelope of $X \ra X\times_S^{\op{log}}X$ as constructed by Kato in \cite{Kat89}. We identify the étale sites of $X$ and $P_{X/S}$ via the immersion $i:X \ra P_{X/S}.$
Let $\calI$ and $\ov{\calI}$ be the ideals of the immersions $X \ra Y$ and $X \ra P_{X/S}$ respectively. Let $\calJ=\ov{\calI}^{[p+1]}+\calI \calP_{X/S},$ where $\calP_{X/S}$ is the structural ring of $P_{X/S}.$ The $\Delta^{-1}\Ox_Y$-module $\ov{\calI}/\calJ$ is annihilated by $\calI$ and is thus an $\Ox_X$-module.
\end{parag}

\begin{proposition}\label{propDoma1}
The map
$$u:\ov{\calI} \ra \ov{\calI}/\calJ,\ x\mapsto x^{[p]}$$
induces an $\Ox_X$-linear isomorphism
$$F_X^*\omega^1_{X/S}=F_X^*\left (\ov{\calI}/\ov{\calI}^{[2]}\right ) \xrightarrow{\sim} \ov{\calI}/\calJ.$$
\end{proposition}

\begin{proof}
If $x,y\in \ov{\calI}$ then
$$u(xy)=(xy)^{[p]}=x^py^{[p]}=p!x^{[p]}y^{[p]}=0$$
and
$$u(x^{[2]})=\left (x^{[2]}\right )^{[p]}=\alpha x^{[2p]}=0,$$
where $\alpha=\frac{(2p)!}{p!2^p}.$ This proves that $u$ induces a map
$$u:\ov{\calI}/\ov{\calI}^{[2]}\ra \ov{\calI}/\calJ.$$
Let $x,y\in \ov{\calI}.$ By \cite{Ogus78} 3.20, there exists $a_i,b_i\in \calP_{X/S},$ local sections $x_i,y_i$ of the image of $\calI$ in $\calP_{X/S}$ and positive integers $n_i,$ $1\le i\le r$ such that
$$x=\sum_{i=1}^ra_ix_i,\ y=\sum_{i=1}^rb_iy_i.$$
We have
$$(x+y)^{[p]}=x^{[p]}+y^{[p]}+\sum_{k=1}^{p-1}x^{[k]}y^{[p-k]}.$$
If $n_i<p$ then $a_ix_i^{[n_i]}=n_i!a_ix_i^{n_i}\in \calI \calP_{X/S},$ so, modulo $\calJ,$ we can suppose that $n_i\ge p$ for any $1\le i\le r.$ For $1\le k,l\le p-1$ and $1\le i,j\le r,$ there exists an integer $c_{k,i,j}$ such that
$$\left (x_i^{[n_i]}\right )^k \left (y_j^{[n_j]}\right )^l=k!l! c_{k,i,j}x_i^{[kn_i]}y_j^{[ln_j]}\in \ov{\calI}^{[p+1]}.$$
It follows that, modulo $\calJ,$
$$(x+y)^{[p]}=x^{[p]}+y^{[p]}$$
and so $u$ is additive.
In addition, $u(ax)=(ax)^{[p]}=a^px^{[p]}=a^pu(x)$ for $a\in \Ox_X$ and $x\in \ov{\calI},$ so $u$ is $F_X$-linear so it induces an $\Ox_X$-linear morphism
$$u:F_X^*\left ( \ov{\calI}/\ov{\calI}^{[2]} \right ) \ra \ov{\calI}/\calJ.$$
To prove that $u$ is an isomorphism, we may work étale locally and thus suppose that we have local coordinates $m_1,\hdots,m_d \in \Gamma(X,\calM_X)$ for $X/S.$ Let $\eta_i \in \calI$ be a local lifting of $\op{dlog}m_i \in \omega^1_{X/S}=\Delta^{-1} \left (\calI/\calI^2\right ).$ Then $(\eta_i^{[p]})_{1\le i\le d}$ is a local basis for the $\Ox_X$-module $\ov{\calI}/\calJ.$ Indeed, let $I=(I_1,\hdots,I_d)\in \N^d$ and $\eta^{[I]}=\prod_{i=1}^d\eta_i^{[I_i]}\in \calP_{X/S}.$ If $|I|\ge p+1$ then $\eta^{[I]}\in \ov{\calI}^{[p+1]}\subset \calJ$ and if there exists $1\le i\le d$ such that $1\le I_i<p$ then $\eta_i^{[I_i]}=(I_i)!^{-1}\eta^{I_i}\in \calI \calP_{X/S}\subset \calJ.$ On the other hand, $(\eta_1,\hdots,\eta_d)$ is a local basis for the $\Ox_X$-module $F_X^*\left ( \ov{\calI}/\ov{\calI}^{[2]} \right ).$
\end{proof}

\begin{parag}\label{parDoma}
We are now ready to embed the tangent sheaf $\calT_{X'/S}$ into $F_{1*}\calD_{X/S},$ where $F_1:X \ra X'$ is the exact relative Frobenius. Let $U' \ra X'$ be an étale morphism of schemes and $U=U'\times_{X'}X.$ To a logarithmic derivation $\partial':\omega^1_{U'/S} \ra \Ox_{U'},$ we associate the differential operator of order $\le p,$
$$\calP_{U/S}^p=\calP_{U/S}/\ov{\calI}_{|U}^{[p+1]} \xrightarrow{\op{Id}-p_1^{\#}\circ i^{\#}} \ov{\calI}_{|U}/\ov{\calI}_{|U}^{[p+1]} \ra \ov{\calI}_{|U}/\calJ_{|U} \xrightarrow{\sim} F_X^*\omega^1_{U/S}=F_1^*\omega^1_{U'/S} \xrightarrow{F_1^*\partial'} \Ox_U,$$
where the second arrow is the canonical surjection, the third arrow is the isomorphism given in \ref{propDoma1}, $i:X \ra P_{X/S}$ is the canonical immersion and $p_1:P_{X/S}\ra X$ is the first projection.
This construction gives an $\Ox_{X'}$-linear morphism
\begin{equation}\label{Doma11}
\psi:\calT_{X'/S} \ra F_{1*}\calD_{X/S}.
\end{equation}
By (\cite{Ohkawa} 4.15), the image of $\psi$ is a submodule of the center of $\calD_{X/S}$ and so it induces a morphism of $\Ox_{X'}$-algebras
\begin{equation}\label{eqDoma1}
\psi:S^{\cdot}\calT_{X'/S} \ra F_{1*}\calD_{X/S}.
\end{equation}
We refer to $\psi$ as the \emph{$p$-curvature map}.
If $m_1,\hdots,m_d \in \calM_X$ are local coordinates for $X/S$ then $\pi^{\flat}m_1,\hdots,\pi^{\flat}m_d\in \calM_{X'}$ (where $\pi:X'\ra X$ is defined in \ref{diag51}) are local coordinates for $X'/S.$ For every multi-index $I\in\N^d,$ denote by $\partial_I$ the differential operator corresponding to $(m_i)_{1\le i\le d},$ as defined in \ref{P2}. Let $(\partial_i')_{1\le i\le d}$ be the dual basis of $(\op{dlog}\pi^{\flat}m_i)_{1\le i\le d}.$ Then, for all $1\le i\le d,$
\begin{equation}\label{era2Doma11}
\psi(\partial_i')=\partial_{p\epsilon_i},
\end{equation}
where $\epsilon_i$ is the multi-index whose all coefficients are zero except for the $i^{th}$ coefficient which is equal to $1.$ 
For all $J\in \N^d,$ set
$$
\partial'^J=\prod_{j=1}^d\partial_j'^{J_j} \in S^{\cdot}\calT_{X'/S}.
$$
Then,
\begin{alignat*}{2}
\psi\left (\partial'^J\right ) &= \prod_{j=1}^d\psi\left (\partial_j' \right )^{J_j} \\
&= \prod_{j=1}^d\partial_{p\epsilon_j}^{J_j} \\
&= \partial_{pJ_1\epsilon_1+\hdots+pJ_d\epsilon_d} \\
&= \partial_{pJ},
\end{alignat*}
where the third line is a consequence of \eqref{eq3102}.
In conclusion, for all $J\in \N^d,$ we have
\begin{equation}\label{era2Koko11}
\psi(\partial'^J)=\partial_{pJ}.
\end{equation}
\end{parag}

\begin{lemma}\label{era4tak2}
The composition
$$v:S^{\bullet}\calT_{X'/S} \xrightarrow{\psi} \calD_{X/S} \ra \widehat{\calD}_{X/S}=\mathscr{Hom}_{\Ox_X}(\calP_0,\Ox_X),$$
where the first arrow $\psi$ is the $p$-curvature map \eqref{eqDoma1} and the second arrow is induced by the canonical projections $\calP_0 \ra \calP_0/\ov{\calI}_0^{[n]}$ for $n\ge 1,$ induces
\begin{equation}\label{eq13191}
\widehat{\psi}:\widehat{S}^{\bullet}\calT_{X'/S} \ra \widehat{\calD}_{X/S}=\mathscr{Hom}_{\Ox_X}(\calP_0,\Ox_X).
\end{equation}
\end{lemma}

\begin{proof}
Suppose that we have local coordinates $m_1,\hdots,m_d \in \Gamma(X,\calM_X)$ and $\eta_1,\hdots,\eta_d \in \calP_0$ such that
$$\left (\eta^{[I]}:=\prod_{i=1}^d\eta_i^{[I_i]} \right )_{I\in \N^d}$$
is a basis for the $\Ox_X$-module $\calP_0,$ as in \ref{P2}.
First, note that
\begin{alignat*}{2}
\w{S}^{\bullet}\calT_{X'/S} &= \lim\limits_{\substack{\longleftarrow \\ n\ge 1}} S^{\bullet}\calT_{X'/S} / \bigoplus_{k\ge n}S^k\calT_{X'/S} \\
&= \lim\limits_{\substack{\longleftarrow \\ n\ge 0}} \bigoplus_{k\le n}S^k\calT_{X'/S} \\
&= \prod_{k\in \N}S^k\calT_{X'/S}.
\end{alignat*}
It is sufficient to define $\w{\psi}:\calP_0\ra \Ox_X$ by
$$\w{\psi}(t)=\sum_{k=0}^{\infty}v(t_k).$$
for every
$$t=(t_0,t_1,\hdots )\in \prod_{k\in \N}S^k\calT_{X'/S}.$$
Note that this infinite sum is well-defined since, for every $I\in \N^d,$ $v(t_k)\left (\eta^{[I]} \right )=0$ for sufficiently large $k$ by \eqref{era2Doma11}.
\end{proof}

\newpage

\section{Logarithmic formal schemes}

\begin{parag}\label{parag62}
The basic definitions and results for logarithmic schemes extend to formal logarithmic schemes. In particular, we have the notions of prelogarithmic and logarithmic structures on formal schemes, logarithmic structure associated to a prelogarithmic structure, as well as the notions of pullback structures, strict morphisms and charts (see \cite{Kat89}).
\end{parag}

\begin{definition}
An \emph{adic} (resp. \emph{$p$-adic}) \emph{logarithmic formal scheme} is a logarithmic formal scheme whose underlying formal scheme is adic (\cite{Ahmed2010} 2.1.24) (resp. $p$-adic i.e. an adic formal scheme such that the ideal $(p)$ generated by $p$ is an ideal of definition).
\end{definition}

\begin{parag}
Let $P$ be a monoid. We denote by $\Z_p\langle P\rangle$ the $p$-adic completion of the ring $\Z[P];$
\begin{equation}\label{BM}
\Z_p\langle P\rangle=\lim\limits_{\substack{\longleftarrow\\ n\ge 1}}\left (\Z/p^n\Z\right )[P].
\end{equation}
We denote by $B\langle P \rangle$ the $p$-adic formal scheme $\op{Spf}(\Z_p\langle P\rangle)$ equipped with the logarithmic structure induced by the canonical morphism $P\ra \Z_p\langle P\rangle.$ 
The following lemma is a variant of (\cite{Ogus2018} III 1.2.4) for $p$-adic logarithmic formal schemes:
\end{parag}

\begin{lemma}\label{Lemlog18}
Let $\frakX$ be a $p$-adic logarithmic formal scheme and $P$ a monoid. Then, we have a canonical bijection
$$\op{Hom}(\frakX,B\langle P \rangle)\stackrel{\sim}{\rightarrow} \op{Hom}_{\boldsymbol{\op{Mon}}}(P,\Gamma(\frakX,\calM_{\frakX})).$$
\end{lemma}

\begin{proof}
Consider the morphism $$\alpha:\op{Hom}(\frakX,B\langle P \rangle)\ra \op{Hom}_{\boldsymbol{\op{Mon}}}(P,\Gamma(\frakX,\calM_{\frakX}))$$ defined by composition with the canonical morphism
$$P\ra \Gamma(B\langle P \rangle,\calM_{B\langle P \rangle}).$$
Conversely, we define a morphism $$\beta:\op{Hom}_{\boldsymbol{\op{Mon}}}(P,\Gamma(\frakX,\calM_{\frakX}))\ra \op{Hom}(\frakX,B\langle P \rangle)$$ as follows: if $\theta:P\ra \Gamma(\frakX,\calM_{\frakX})$ is a morphism of monoids, we consider a covering $\frakX=\bigcup_{i\in I}\op{Spf}A_i$ of $\frakX$ by affine formal schemes, where $A_i$ is a complete $\Z_p$-algebra. For any $i\in I,$ let $\theta_i$ be the composition
$$P\ra \Gamma(\frakX,\calM_{\frakX})\ra \Gamma(\op{Spf}A_i,\Ox_{\frakX})=A_i.$$
Since $A_i$ is a complete $\Z_p$-algebra, the morphism $\theta_i$ induces a continuous morphism $\Z_p\langle P\rangle\ra A_i$ and so a morphism of formal schemes $f_i:\op{Spf}A_i\ra B\langle P \rangle.$ The morphism $P\xrightarrow{\theta} \Gamma(\frakX,\calM_{\frakX})\ra \Gamma(\op{Spf}A_i,\calM_{\frakX})$ induces a morphism $f_i^{\flat}:\calM_{B\langle P \rangle}\ra f_{i*}\calM_{\op{Spf}A_i}.$ The commutativity of the diagram
$$\begin{tikzcd}
P\ar{r}\ar{d} & \Gamma(\op{Spf}A_i,\calM_{\frakX})\ar{d} \\
\Z_p\langle P\rangle \ar{r} & A_i
\end{tikzcd}$$ 
implies that $(f_i,f_i^{\flat})$ is a morphism of logarithmic formal schemes. These morphisms glue into a morphism $\beta(\theta):\frakX\ra B\langle P \rangle$ and this morphism is independant of the choice of the covering. The maps $\alpha$ and $\beta$ are inverse to each other. 
\end{proof}

\begin{parag}\label{parag63}
We say that a $p$-adic logarithmic formal scheme $\frakX$ is \emph{quasi-coherent} (resp. \emph{coherent}, resp. \emph{fine}, resp. \emph{finitely-generated and saturated}, or \emph{fs} for short) if étale locally on $\frakX,$ there exists a strict morphism $\frakX\ra B\langle P \rangle$ for a monoid (resp. a finitely generated monoid, resp. a fine monoid, resp. an fs monoid) $P.$ We denote by $\boldsymbol{\op{LFS}}$ the category of fs $p$-adic logarithmic formal schemes.
\end{parag}

\begin{parag}\label{parag27}
Let $\frakX\ra \frakS$ and $\frakY\ra \frakS$ be morphisms in the category $\boldsymbol{\op{LFS}}.$ Following the conventions of \ref{parag42}, we denote by $\frakX\times_{\frakS}^{\op{log}}\frakY$ the fiber product of $\frakX \ra \frakS$ and $\frakY \ra \frakS$ in the category $\boldsymbol{\op{LFS}}$ and we keep the notation $\frakX \times_{\frakS} \frakY$ for the fiber product in the category of logarithmic formal schemes. The two are compared as follows: locally, there exist strict morphisms $\frakX\ra B\langle M \rangle,$ $\frakY\ra B\langle N \rangle$ and $\frakS\ra B\langle P \rangle$ for fs monoids $M,\ N$ and $P$ (\ref{parag63}), that fit into a commutative diagram (\cite{Ogus2018} III 1.2.7.3)
$$\begin{tikzcd}
\frakX\ar{r}\ar{d} & \frakS\ar{d} & \frakY\ar{l}\ar{d} \\
B\langle M \rangle\ar{r} & B\langle P \rangle & B\langle N \rangle\ar{l}
\end{tikzcd}$$ 
Let $Q=M\oplus_PN$ and $Q^{sat}$ its saturation. The commutative diagram above induces the following commutative diagram of monoids:
$$\begin{tikzcd}
P\ar{r}\ar{d} & M\ar{d} \\
N\ar{d} & \Gamma(\frakX,\calM_{\frakX})\ar{d} \\
\Gamma(\frakY,\calM_{\frakY})\ar{r} & \Gamma(\frakX\times_{\frakS}\frakY,\calM_{\frakX\times_{\frakS}\frakY})
\end{tikzcd}$$
and so we get a morphism $M\oplus_PN\ra \Gamma(\frakX\times_{\frakS}\frakY,\calM_{\frakX\times_{\frakS}\frakY}),$ which, by \ref{Lemlog18}, is equivalent to a morphism $B\langle M\oplus_PN\rangle \ra \frakX\times_{\frakS}\frakY.$ By definition of the logarithmic structure on a fiber product (\cite{Ogus2018} III 2.1), this morphism is a chart. 
Then the underlying formal scheme of $\frakX\times_{\frakS}^{\op{log}}\frakY$ is $\frakX\times_{\frakS}\frakY\times_{B\langle M\oplus_PN \rangle}B\langle Q^{sat} \rangle$ and its logarithmic structure is the pullback of that of $B\langle Q^{sat} \rangle.$
\end{parag}

\begin{parag}\label{parag17}
Let $(\frakX_n,\Ox_{\frakX_n},\calM_{\frakX_n})_{n\ge 1}$ be an inductive system of logarithmic schemes such that for every $n\ge 1,$ $\frakX_n$ is a $\Z/p^n\Z$-scheme, the diagram
$$\begin{tikzcd}
\frakX_n\ar{r}\ar{d}\ar{r} & \Sp \Z/p^n\Z \ar{d}\\
\frakX_{n+1}\ar{r} & \Sp \Z/p^{n+1}\Z
\end{tikzcd}$$
is cartesian in the category of schemes and $\frakX_n\ra\frakX_{n+1}$ is strict for all $n\ge 1.$  
Then the inductive limit $(\frakX,\Ox_{\frakX})$ of $((\frakX_n,\Ox_{\frakX_n}))_{n\ge 1}$ is a $p$-adic formal scheme. Let $f_n:\frakX_n\ra \frakX$ be the canonical morphisms. The morphisms $f_n$ are homeomorphisms of the underlying topological spaces and so we use them to identify the topological spaces $\frakX_n$ and $\frakX.$ Let $\calM_{\frakX}:=\lim\limits_{\longleftarrow}\calM_{\frakX_n}$ and $\alpha_{\frakX}:\calM_{\frakX}\ra \Ox_{\frakX}$ be the morphism induced by the morphisms $\alpha_{\frakX_n}.$ Since $\Ox_{\frakX}^*=\lim\limits_{\longleftarrow}\Ox_{\frakX_n}^*,$ the morphism $\alpha_{\frakX}:\calM_{\frakX}\ra \Ox_{\frakX}$ defines a logarithmic structure on $\frakX.$ Then $\frakX$ is a $p$-adic logarithmic formal scheme and $(\frakX,\calM_{\frakX})$ is the inductive limit of $(\frakX_n,\calM_{\frakX_n})$ in the category logarithmic formal schemes. 
Conversely, we have the following proposition:
\end{parag}

\begin{proposition}\label{propkey}
Let $\frakX$ be a $p$-adic integral logarithmic formal scheme and denote, for any integer $n\ge 1,$ by $\frakX_n$ the scheme $(\frakX,\Ox_{\frakX}/(p^n))$ and by $f_n:\frakX_n\ra \frakX$ the canonical morphism. Equip $\frakX_n$ with the logarithmic structure $f_n^*\calM_{\frakX}.$ Then the canonical morphism
\begin{equation}\label{canmo1}
\varphi:\calM_{\frakX}\ra \lim_{\substack{\longleftarrow\\ n\ge 1}}\calM_{\frakX_n},\ m\mapsto (f_n^{\flat}(m))_{n\ge 1}
\end{equation}
is an isomorphism.
\end{proposition}

\begin{proof}
Let $\alpha_n$ be the composition $\calM_{\frakX}\xrightarrow{\alpha_{\frakX}} \Ox_{\frakX}\ra \Ox_{\frakX_n},$ where $\Ox_{\frakX}\ra \Ox_{\frakX_n}$ is reduction modulo $p^n.$
By the definition of pullback of logarithmic structures (\cite{Ogus2018} III 1.1.5.2),
\begin{equation}
\calM_{\frakX_n}=f_n^*\calM_{\frakX}=\calM_{\frakX}\oplus_{\alpha_n^{-1}(\Ox_{\frakX_n}^*)}\Ox_{\frakX_n}^*.
\end{equation}
By this equality, the morphism (\ref{canmo1}) is given by
\begin{equation}
\varphi(m)=((m,1))_{n\ge 1}.
\end{equation}
Let $m$ and $m'$ be local sections of $\calM_{\frakX}$ such that $\varphi(m)=\varphi(m').$ For any positive integer $n,$ $(m,1)=(m',1)$ in $\calM_{\frakX_n}=\calM_{\frakX}\oplus_{\alpha_n^{-1}(\Ox_{\frakX}^*)}\Ox_{\frakX_n}^*.$ So, by (\cite{Kat89} (1.3)), there exists local sections $t_n,t_n'$ of $\alpha_n^{-1}(\Ox^*_{\frakX})$ such that
\begin{equation}\label{amsum}
\begin{cases}m+t_n=m'+t_n' \\ \alpha_n(t_n)=\alpha_n(t_n').\end{cases}
\end{equation} 
The local sections $\alpha_{\frakX}(t_n)$ and $\alpha_{\frakX}(t_n')$ are invertible modulo $p^n$ and so they are both invertible in $\Ox_{\frakX}.$ Since the morphism $\alpha_{\frakX}^{-1}(\Ox_{\frakX}^*)\rightarrow\Ox_{\frakX}^*$ induced by $\alpha_{\frakX}$ is an isomorphism, $t_n$ and $t_n'$ are invertible in $\calM_{\frakX}.$ Furthermore, $m+t_n-t_n'=m'$ and $\calM_{\frakX}$ is integral, so $t_n-t_n'=t_k-t'_k$ for all positive integers $n$ and $k.$ Set $d=t_1-t'_1.$ By (\ref{amsum}), $\alpha_{\frakX}(d)=1$ in $\Ox_{\frakX_n}$ for every positive integer $n$ and so $\alpha_{\frakX}(d)=1$ in $\Ox_{\frakX},$ $d=0$ and $t_1=t_1'.$ Again, since $\calM_{\frakX}$ is integral and by (\ref{amsum}), $m=m'.$ Hence (\ref{canmo1}) is injective.

For the surjectivity, let $\frakU=\op{Spf}A$ be an affine open formal subscheme of $\frakX$ and $(m_n)_{n\ge 1}$ a sequence of sections of $\Gamma \left ( \frakU, \calM_{\frakX} \right )$ and, for every $n\ge 1,$ $a_n\in A/(p^n)$ such that $((m_n,a_n))_{n\ge 1}$ yields a local section of $\lim\limits_{\substack{\longleftarrow\\ n\ge 1}}\calM_{\frakX_n}$ over $\frakU.$ For every positive integer $n\ge 1,$ let $b_n\in A$ such that, modulo $p^n,$ $b_n=a_n.$ First, note that $a_1$ is invertible in $A/(p)$ so $b_1$ is invertible in $A.$ Then there exists $t_1\in \Gamma \left (\frakU,\calM^*_{\frakX}\right )$ such that $\alpha_{\frakX}(t_1)=b_1.$ Setting $m=m_1+t_1,$ we get by (\cite{Kat89} (1.3))
\begin{equation}\label{eq683}
(m_1,a_1)=(m_1,\alpha_1(t_1))=(m,1).
\end{equation}
We now construct by induction on $n$ a sequence $(c_n)_{n\ge 1} \in \lim\limits_{\longleftarrow}\left (A/(p^n)\right )^*=A^*$ such that, in $\Gamma \left ( \frakU,\calM_{\frakX_n}\right ),$
$$(m_n,a_n)=(m,c_n)\ \forall n\ge 1.$$
By (\ref{eq683}), we set $c_1=1.$ Fix $n\ge 1$ and suppose $c_n$ is defined. Let $\ov{x} \ra \frakX_n$ be a geometric point. In $\calM_{\frakX_n,\ov{x}}=\calM_{\frakX,\ov{x}}\oplus_{\alpha_n^{-1}(\Ox_{\frakX_n}^*)_{\ov{x}}}\Ox_{\frakX_n,\ov{x}}^*,$ $$(m_{n+1},a_{n+1})=(m_n,a_n)=(m,c_n).$$
Then, again by (\cite{Kat89} (1.3)), there exist a neighborhood $V$ of $\ov{x}$ and local sections $t,t' \in \alpha_n^{-1}(\Ox_{\frakX_n}^*)$ over $V$ such that
$$
\begin{cases}
m_{n+1}+t=m+t' \\
a_{n+1}\alpha_n(t')=c_n\alpha_n(t)\ (\op{mod} p^n).
\end{cases}
$$
We set $c_{n+1}=a_{n+1}\alpha_{n+1}(t'-t).$ Note that $t'-t$ is independant of the choice of $t$ and $t'$ and we can hence suppose that $t'-t \in \Gamma\left (\frakU,\calM_{\frakX_n} \right ).$ Then
$$(m_{n+1},a_{n+1})=(m+t'-t,a_{n+1})=(m,a_{n+1}\alpha_{n+1}(t'-t))=(m,c_{n+1}).$$
The compatibility condition follows from the following equality modulo $p^n:$
$$c_{n+1}=a_{n+1}\alpha_n(t'-t)=c_n.$$
Let $c=(c_n)_{n\ge 1} \in \lim\limits_{\longleftarrow} \left (A/(p^n)\right )^*=A^*$ and $s\in \Gamma\left (\frakU,\calM_{\frakX} \right )$ such that $\alpha_{\frakX}(s)=c.$ Then $((m_n,a_n))_{n\ge 1}$ is the image of $m+s$ by $\varphi$ \eqref{canmo1}.
\end{proof}

\begin{corollaire}\label{corkey}
Under the same hypothesis of \ref{propkey}, the logarithmic formal scheme $\frakX$ is the inductive limit, in the category of logarithmic formal schemes, of $(\frakX_n)_{n\ge 1}.$
\end{corollaire}

\begin{parag}\label{parag68}
Let $f:\frakX\ra \frakY$ be a morphism of $p$-adic logarithmic formal schemes. Denote, for any integer $m\ge n\ge 1,$ by $\frakX_n$ (resp. $\frakY_n$) the scheme $(\frakX,\Ox_{\frakX}/p^n)$ (resp. $(\frakY,\Ox_{\frakY}/p^n)$) and by $i_n:\frakX_n\ra \frakX,$ $j_n:\frakY_n\ra \frakY,$ $i_{n,m}:\frakX_n \ra \frakX_m$ and $j_{n,m}:\frakY_n \ra \frakY_m$ the canonical morphisms. Equip $\frakX_n$ (resp. $\frakY_n$) with the logarithmic structure $i_n^*\calM_{\frakX}$ (resp. $j_n^*\calM_{\frakY}$). By (\cite{Ahmed2010}, 2.2.1), the morphism $f$ induces an inductive system of morphisms of schemes $(f_n:\frakX_n\ra \frakY_n).$ Then, by the commutativity of the diagram
$$\begin{tikzcd}
\frakX_n\ar{r}{f_n}\ar{d}{i_n} & \frakY_n\ar{d}{j_n}\\
\frakX\ar{r}{f} & \frakY
\end{tikzcd}$$ 
and since $i_n$ and $j_n$ are strict, the morphism $f^{\flat}:f^*\calM_{\frakY}\ra \calM_{\frakX}$ induces a morphism $f_n^*\calM_{\frakY_n}\ra \calM_{\frakX_n}$ for all $n\ge 1,$ making $f_n$ into a morphism of logarithmic schemes. Similarly, since $i_{n,m}^*\calM_{\frakX_m}=i_n^*\calM_{\frakX}$ (resp. $j_{n,m}^*\calM_{\frakY_m}=j_n^*\calM_{\frakY}$), the morphism $i_n^{\flat}$ (resp. $j_n^{\flat}$) turns $i_{n,m}$ (resp. $j_{n,m}$) into a morphism of logarithmic schemes. The map $f\mapsto (f_n)_{n\ge 1}$ is a bijection, as stated in the following proposition:
\end{parag}

\begin{proposition}\label{prop18}
Keep the same hypothesis of \ref{parag68} and assume $\calM_{\frakX}$ and $\calM_{\frakY}$ are integral. The map $f\mapsto (f_n)_{n\ge 1}$ defined in \ref{parag68} is a bijection between the set of morphisms of logarithmic formal schemes $f:\frakX\ra \frakY$ and the set of sequences of morphisms of logarithmic schemes $(f_n:\frakX_n\ra \frakY_n)$ making the following diagram commutative for all $n\ge m$
$$
\begin{tikzcd}
\frakX_m\ar{r}{f_m}\ar{d} & \frakY_m\ar{d}\\
\frakX_n\ar{r}{f_n} & \frakY_n
\end{tikzcd}$$
\end{proposition}

\begin{proof}
Let $(f_n:\frakX_n \ra \frakY_n)_{n\ge 1}$ a sequence of morphisms of logarithmic schemes such that
$$\begin{tikzcd}
\frakX_n\ar{r}{f_n}\ar{d}{i_{n,m}} & \frakY_n\ar{d}{j_{n,m}}\\
\frakX_m\ar{r}{f_m} & \frakY_m
\end{tikzcd}$$ 
is commutative for every integers $m\ge n\ge 1.$
By (\cite{Ahmed2010} 2.2.2), the morphisms $f_n$ induce a morphism of formal schemes. In addition, by \ref{propkey}, $\calM_{\frakY}=\lim\limits_{\longleftarrow} \calM_{\frakY_n}$ and $\calM_{\frakX}=\lim\limits_{\longleftarrow} \calM_{\frakX_n}.$ Then the compositions
$$f^{-1}\calM_{\frakY} \ra f^{-1}\calM_{\frakY_n} \xrightarrow {f_n^{\flat}} \calM_{\frakX_n}$$
induce a morphism $f^{-1}\calM_{\frakY} \ra \calM_{\frakX}$ which turns $f$ into a morphism of logarithmic formal schemes. This construction is inverse to that of \ref{parag68}.
\end{proof}

\begin{parag}\label{parag69}
Keep the same notations of \ref{parag68}.
For every integer $n\ge 1,$ denote by $i_{n,n+1}:\frakX_n\ra \frakX_{n+1}$ the canonical morphism and by $$(d_n,\delta_n):\Ox_{\frakX_n}\times \calM_{\frakX_n}\ra \omega^1_{\frakX_n/\frakY_n}$$
the universal derivation. Then the morphisms
$$\Ox_{\frakX_{n+1}}\xrightarrow{i_{n,n+1}^{\#}}\Ox_{\frakX_n}\xrightarrow{d_n}\omega^1_{\frakX_n / \frakY_n}$$
and
$$\calM_{\frakX_{n+1}}\xrightarrow{i_{n,n+1}^{\flat}}\calM_{\frakX_n}\xrightarrow{\delta_n}\omega^1_{\frakX_n / \frakY_n}$$
define a logarithmic derivation $\Ox_{\frakX_{n+1}}\times\calM_{\frakX_{n+1}}\ra \omega^1_{\frakX_n / \frakY_n}$ and thus a morphism
$$\omega^1_{\frakX_{n+1} / \frakY_{n+1}}\ra \omega^1_{\frakX_n / \frakY_n}.$$
These morphisms define a projective system $(\omega^1_{\frakX_n/\frakY_n})_{n\ge 1}.$ We set
\begin{equation}\label{eq681}
\omega^1_{\frakX/\frakY}=\lim_{\substack{\longleftarrow \\ n\ge 1}}\omega^1_{\frakX_n / \frakY_n},
\end{equation}
that we call \emph{the sheaf of logarithmic differentials of $\frakX$ over $\frakY$}.
\end{parag}

\begin{proposition}\label{prop69}
Let $\frakX$ be a $p$-adic logarithmic formal scheme and denote, for any integer $n\ge 1,$ by $\frakX_n$ the scheme $(\frakX,\Ox_{\frakX}/p^n)$ and, for any $1\le n\le k,$ denote by $f_n:\frakX_n\ra \frakX$ and $f_{n,k}:\frakX_{n}\ra \frakX_{k}$ the canonical morphisms. Equip $\frakX_n$ with the logarithmic structure $f_n^*\calM_{\frakX}.$ Then for any $1\le n\le k,$ the morphisms $f_n^{\flat}:\calM_{\frakX}\ra \calM_{\frakX_n}$ and $f_{n,k}^{\flat}:\calM_{\frakX_k}\ra \calM_{\frakX_n}$ are surjective.
\end{proposition}

\begin{proof}
Let $k\ge n\ge 1.$ Since $f_n^{\flat}=f_{n,k}^{\flat}\circ f_k^{\flat},$ it is sufficient to prove that $f_n^{\flat}$ is surjective. By definition of the logarithmic structure on $\frakX_n,$ the morphism $f_{n}$ is strict so $$\calM_{\frakX_n}=f_{n}^*\calM_{\frakX}=\calM_{\frakX}\oplus_{\alpha_n^{-1}(\Ox_{\frakX_n}^*)}\Ox_{\frakX_n}^*,$$
where $\alpha_n:\calM_{\frakX}\ra \Ox_{\frakX_n}$ is the composition
$\calM_{\frakX}\xrightarrow{\alpha_{\frakX}}\Ox_{\frakX}\xrightarrow{f_{n}^{\#}}\Ox_{\frakX_{n}}.$ Let $m'$ and $a_n$ be local sections of $\calM_{\frakX}$ and $\Ox_{\frakX_{n}}^*$ respectively. Set $m_n=(m',a_n),$ local section of $\calM_{\frakX_{n}}.$ Let $a$ be a local lifting of $a_n$ to $\Ox_{\frakX,\ov{x}}.$ Since $\Ox_{\frakX}=\lim\limits_{\substack{\longleftarrow\\ n\ge 1}}\Ox_{\frakX_{n}}$ and $a_n$ is invertible in $\Ox_{\frakX_n},$ it follows that $a$ is also invertible in $\Ox_{\frakX}$ and so $a=\alpha_{\frakX}(t)$ for a local section $t$ of $\calM_{\frakX}.$ It follows that
$$m_n=(m',a_n)=(m',\alpha_{n}(t))=(t+m',1)=f_{n}^{\flat}(t+m').$$
\end{proof}

\begin{parag}
Let $\Delta:\frakX \ra \frakY$ be an exact closed immersion of integral $p$-adic logarithmic locally Noetherian formal schemes with ideal $\calI.$ For any $x\in \frakX,$ the ideal $\calI_x$ is a proper ideal of $\Ox_{\frakY,x}$ so $1+\calI_x \subset \Ox_{\frakY,x}^*.$ It follows that $\Delta^{-1}(1+\calI)\subset \Delta^{-1}\Ox_{\frakY}^*.$ Denote by $\lambda:\Delta^{-1}(1+\calI) \ra \Delta^{-1}\calM_{\frakY}$ the composition of $\Delta^{-1}\alpha_{\frakY}^{-1}:\Delta^{-1}\Ox_{\frakY}^* \ra \Delta^{-1}\calM_{\frakY}$ with the canonical embedding $\Delta^{-1}(1+\calI) \hookrightarrow \Delta^{-1}\Ox_{\frakY}^*.$ Consider the sequence of monoids
\begin{equation}
0 \ra \Delta^{-1}(1+\calI) \xrightarrow{\lambda} \Delta^{-1}\calM_{\frakY} \xrightarrow{\Delta^{\flat}} \calM_{\frakX} \ra 0.
\end{equation}
This is exact by the following proposition:
\end{parag}

\begin{proposition}
Let $\Delta:\frakX \ra \frakY$ be an exact closed immersion of integral $p$-adic logarithmic locally Noetherian formal schemes with ideal $\calI.$ Then the sequence \eqref{era2exactseq}
\begin{equation}\label{era2exactseq}
0 \ra \Delta^{-1}(1+\calI) \xrightarrow{\lambda} \Delta^{-1}\calM_{\frakY} \xrightarrow{\Delta^{\flat}} \calM_{\frakX} \ra 0
\end{equation}
is an exact sequence of monoids i.e. $\lambda$ is injective, $\Delta^{\flat}$ is surjective and for any local sections $m$ and $m'$ of $\Delta^{-1}\calM_{\frakY}$ such that $\Delta^{\flat}(m)=\Delta^{\flat}(m'),$ there exists a local section $a$ of $\Delta^{-1}(1+\calI)$ such that
$$
m+\lambda(a)=m'.
$$
\end{proposition}

\begin{proof}
The morphism $\lambda$ (resp. $\Delta^{\flat}$) is injective (resp. surjective) by definition.
For a positive integer $n,$ let $\Delta_n:\frakX_n \ra \frakY_n$ be the immersion obtained from $\Delta$ by reduction modulo $p^n,$ as in \ref{prop69} and denote by $\calI_n$ its ideal. It is sufficient to prove that the sequence
\begin{equation}\label{era2exactseqprime}
0 \ra \Delta_n^{-1}(1+\calI_n) \xrightarrow{\lambda_n} \Delta_n^{-1}\calM_{\frakY_n} \xrightarrow{\Delta_n^{\flat}} \calM_{\frakX_n}
\end{equation}
is exact for all positive integers $n.$
Fix an integer $n\ge 1.$
The immersion $\Delta_n$ is exact so it is strict. It follows that $\calM_{\frakX_n}=\Delta_n^*\calM_{\frakY_n}$ is equal to the direct sum
$$
\begin{tikzcd}
\gamma^{-1}\Ox_{\frakX_n}^* \ar[hook]{r} \ar[swap]{d}{\gamma} & \Delta_n^{-1}\calM_{\frakY_n}\ar{d} \\
\Ox_{\frakX_n}^* \ar{r} & \Delta_n^{-1}\calM_{\frakY_n}\oplus_{\gamma^{-1}\Ox_{\frakX_n}^*}\Ox_{\frakX_n}^*,
\end{tikzcd}
$$
where $\gamma:\Delta_n^{-1}\calM_{\frakY_n} \xrightarrow{\Delta_n^{-1}\alpha_{\frakY_n}}\Delta_n^{-1}\Ox_{\frakY_n} \xrightarrow{\Delta_n^{\#}} \Ox_{\frakX_n}.$
The sequence \eqref{era2exactseqprime} becomes isomorphic to
$$
0 \ra \Delta_n^{-1}(1+\calI_n) \xrightarrow{\lambda_n} \Delta_n^{-1}\calM_{\frakY_n} \xrightarrow{\Delta_n^{\flat}} \Delta_n^{-1}\calM_{\frakY_n}\oplus_{\gamma^{-1}\Ox_{\frakX_n}^*}\Ox_{\frakX_n}^*,
$$
where $\Delta_n^{\flat}(m)=(m,1)$ for any local section $m$ of $\Delta_n^{-1}\calM_{\frakY_n}.$
Let $\ov{x} \ra \frakX_n$ be a geometric point and $m,m'\in \calM_{\frakY_n,\ov{x}}$ such that
$$
\Delta^{\flat}_{n,\ov{x}}(m)=\Delta^{\flat}_{n,\ov{x}}(m').
$$
By definition of amalgamated sums in the category of monoids, there exists $t,t'\in \left (\gamma^{-1}\Ox_{\frakX_n}^*\right )_{\ov{x}}$ such that
$$
\begin{cases}
m+t=m'+t' \\
\gamma(t)=\gamma(t').
\end{cases}
$$
Since $\gamma(t)\in \Ox_{\frakX_n}^*$ and $\Ox_{\frakX_n,\ov{x}}=\Ox_{\frakY_n,\ov{x}}/\calI_{\ov{x}},$ we get
$$
\alpha_{\frakY_n,\ov{x}}(t)\in \Ox_{\frakY_n,\ov{x}}^*.
$$
It follows that $t\in \calM_{\frakY_n,\ov{x}}^*$ and so
$$
m=m'+t'-t \in \calM_{\frakY_n,\ov{x}}.
$$
Since
$$
t-t'=\alpha_{\frakY_n,\ov{x}}^{-1} \left (\alpha_{\frakY_n,\ov{x}}(t)\alpha_{\frakY_n,\ov{x}}(t')^{-1} \right )
$$
and
$$
\Delta^{\#}_{\ov{x}} \left (\alpha_{\frakY_n,\ov{x}}(t)\alpha_{\frakY_n,\ov{x}}(t')^{-1} \right )=\gamma(t)\gamma(t')^{-1}=1,
$$
we get
$$
\alpha_{\frakY_n,\ov{x}}(t)\alpha_{\frakY_n,\ov{x}}(t')^{-1}\in 1+\calI_{\ov{x}}
$$
and
$$
t-t'=\lambda_{n,\ov{x}} \left (\alpha_{\frakY_n,\ov{x}}(t)\alpha_{\frakY_n,\ov{x}}(t')^{-1} \right ).
$$
\end{proof}

\begin{parag}
We end this section by defining log étaleness, log smoothness and log flatness in the context of $p$-adic logarithmic formal schemes and giving a useful property of infinitesimal liftings. 
\end{parag}

\begin{definition}\label{erafsmoothdef}
A morphism $f:\frakX\ra \frakY$ of $p$-adic logarithmic formal schemes is said to be \emph{log étale} (resp. \emph{log smooth}) if $f$ is locally of finite presentation (resp. locally of finite type) (\cite{Ahmed2010} 2.3.13 and 2.3.15), étale locally on $\frakX$ and $\frakY,$ there exists a chart $\theta:P\ra Q$ of $f$ such that
\begin{enumerate}
\item $\theta^{gp}$ is injective and $\op{coker}\theta^{gp}$ (resp. the torsion subgroup of $\op{coker}\theta^{gp}$) is finite and of order coprime with $p.$
\item The morphism $\frakX \ra \frakY\times_{B\langle P\rangle}B\langle Q\rangle,$ induced by $f$ and the chart $\frakX \ra B\langle Q\rangle$ is étale.
\end{enumerate}
It is clear that if $f$ is log étale (resp. log smooth), then, for every integer $n\ge 1,$ the induced morphism of logarithmic schemes $f_n:\frakX_n\ra \frakY_n$ (\ref{parag68}) is log étale (resp. log smooth).
\end{definition}

\begin{proposition}\label{proplift}
Let $f:\frakX \ra \frakY$ be a morphism of $p$-adic logarithmic formal schemes. Suppose we are given a commutative diagram
\begin{equation}\label{diag6161}
\begin{tikzcd}
 & & \frakX \ar{d}{f} \\
T \ar{r}\ar[bend right=-30]{urr} & \frakT \ar{r} \ar[dashed]{ur} & \frakY
\end{tikzcd}
\end{equation}
where $\frakT$ is a $p$-adic logarithmic formal scheme, $T$ is the special fiber of $\frakT$ equipped with the pullback structure of $\frakT$ and $T \ra \frakT$ is the canonical strict morphism. If $f$ is log étale (resp. log smooth) then the dashed arrow exists and is unique (resp. exists locally on $\frakT$).
\end{proposition}

\begin{proof}
For every integer $n\ge 1,$ we denote by $\frakX_n$ (resp. $\frakY_n,$ resp. $\frakT_n$) the logarithmic scheme obtained from $\frakX$ (resp. $\frakY,$ resp. $\frakT$) by reduction modulo $p^n,$ as in \ref{propkey}.
Suppose $f$ is log étale. We construct, by induction on $n,$ a sequence of compatible morphisms $(g_n:\frakT_n\ra \frakX_n).$ For $n=1,$ we take the morphism $g_1:T \ra \frakX_1$ given in (\ref{diag6161}). Let $n$ be a positive integer and suppose we are given a morphism $g_n:\frakT_n \ra \frakX_n$ such that the diagram
\begin{equation}
\begin{tikzcd}
 & & \frakX_n \ar{d}{f_n} \\
T \ar{r}\ar[bend right=-30]{urr} & \frakT_n \ar{r} \ar{ur}{g_n} & \frakY_n
\end{tikzcd}
\end{equation}
is commutative. We obtain the following solid commutative diagram
$$
\begin{tikzcd}
\frakX_n \ar{rr} & & \frakX_{n+1} \ar{d}{f_{n+1}} \\
\frakT_n \ar{u}{g_n} \ar{r} & \frakT_{n+1} \ar[dashed]{ur}{g_{n+1}} \ar{r} & \frakY_{n+1}
\end{tikzcd}
$$
The existence and uniqueness of the dashed arrow $g_{n+1}:\frakT_{n+1} \ra \frakX_{n+1}$ follows from the fact that $f_{n+1}$ is log étale and $\frakT_n \ra \frakT_{n+1}$ is a strict thickening. The morphisms $g_n$ are compatible by construction, so they induce a morphism
$$g:\lim_{\substack{\longrightarrow \\ n\ge 1}}\frakT_n \ra \lim_{\substack{\longrightarrow \\ n\ge 1}}\frakX_n.$$
And by \ref{propkey}, $$\lim_{\substack{\longrightarrow \\ n\ge 1}}\frakT_n=\frakT,\ \lim_{\substack{\longrightarrow \\ n\ge 1}}\frakX_n=\frakX.$$
The proof in the log étale case is complete.
The log smooth case can be proved in a similar way by restricting to affine formal subschemes of $\frakY,$ $\frakX$ and $\frakT.$ 
\end{proof}

\begin{definition}\label{logflatdef}
Let $f:\frakX \ra \frakY$ be a morphism of fine $p$-adic logarithmic formal schemes and, for all positive integers $n,$ $f_n:\frakX_n \ra \frakY_n$ the morphism obtained from $f$ by reduction modulo $p^n$ as in \ref{parag68}. The morphism $f$ is said to be \emph{log flat} if, for every positive integer $n,$ $f_n$ is log flat (\ref{loggfflat}).
\end{definition}

\begin{proposition}\label{Wflat}
Let $W(k)$ be the ring of Witt vectors of a perfect field $k$ of characteristic $p$ and $\frakS=\op{Spf}W$ equipped with the trivial logarithmic structure. If $\frakT$ is an fs logarithmic $p$-adic formal scheme log flat over $\frakS,$ then $\frakT$ is flat over $\frakS$ \eqref{dxuflat}. 
\end{proposition}

\begin{proof}
For any positive integer $n,$ $\frakT_n \ra \frakS_n$ is log flat. Let $n$ be a positive integer. Fppf locally on $\frakT_n,$ there exists a chart $\frakT_n \ra A_n[ M]$ such that $\frakT_n \ra \frakS_n\times_{\op{Spec}\Z/p^n\Z} B[ M]=\op{Spec}(W(k)/(p^n)[M])$ is flat. Since $W(k)/(p^n) \ra W(k)/(p^n)[ M]$ is flat, we get that $\frakT_n$ is flat over $\frakS_n.$ It follows that $\frakT$ is flat over $\frakS.$
\end{proof}

\begin{proposition}\label{logflatfiber}
Let $f:\frakX \ra \frakY$ and $g:\frakY \ra \frakZ$ be two morphisms of fine $p$-adic logarithmic formal schemes, $x\in \frakX,$ and $z=g\circ f(x)$ such that $f$ and $g$ are locally of finite presentation (\cite{Ahmed2010} 2.3.15) and $g\circ f$ and $f_z:\frakX_z\ra \frakY_z$ are log flat in a neighborhood of $x.$ Then $f$ is log flat in a neighborhood of $x.$
\end{proposition}

\begin{proof}
This is a consequence of (\cite{Ogus2018} IV 4.2.2).
\end{proof}

\newpage

\section{Frames on logarithmic formal schemes}
\begin{parag}\label{parag31}
In this section, we generalize the notion of frames, introduced by Kato and Saito in \cite{Saito04}, to logarithmic formal schemes. Let $P$ be an fs monoid. Recall (\ref{parag43}) that we have defined a presheaf of sets $[P]$ on the category $\boldsymbol{\op{L}}$ of fs logarithmic schemes by
$$[P](T)=\op{Hom}_{\boldsymbol{\op{Mon}}}(P,\Gamma(T,\ov{\calM}_{T})).$$
We consider $\boldsymbol{\op{L}}$ as a full subcategory of the category of fs logarithmic formal schemes in the canonical way (\cite{SP}\href{https://stacks.math.columbia.edu/tag/0AHY}{ 0AHY}).
By \ref{Lemlog18}, we have a canonical morphism of presheaves $B\langle P \rangle\ra [P].$ In this section, all logarithmic formal schemes are supposed to be fs. We start with the following proposition:
\end{parag}

\begin{proposition}
The functor
\begin{equation}
F:\boldsymbol{\op{LFS}} \ra \boldsymbol{\widehat{\op{L}}},\ \frakX \mapsto \op{Hom}(-,\frakX)
\end{equation}
is fully faithful.
\end{proposition}

\begin{proof}
Let $f,g:\frakX \ra \frakY$ be morphisms in $\boldsymbol{\op{LFS}}$ such that $F(f)=F(g)$ and consider the notations of \ref{parag68}. For any integer $n\ge 1,$ the equality $f\circ i_n=F(f)(i_n)=F(g)(i_n)=g\circ i_n$ is equivalent to $j_n \circ f_n=j_n \circ g_n.$ Since $i_n:\frakX_n \ra \frakX$ is equal to the identity on the underlying topological spaces, $f=g$ on the underlying topological spaces. On the level of structural sheaves, the compositions $$\Ox_{\frakY} \xrightarrow{j_n^{\#}} \Ox_{\frakY_n} \xrightarrow{f_n^{\#}} f_*\Ox_{\frakX_n},\ \Ox_{\frakY} \xrightarrow{j_n^{\#}} \Ox_{\frakY_n} \xrightarrow{g_n^{\#}} g_*\Ox_{\frakX_n}$$
are equal. And since $j_n^{\#}$ is surjective, $f_n^{\#}=g_n^{\#}.$
On the level of monoids, the compositions
$$\calM_{\frakY}\xrightarrow{j_n^{\flat}}\calM_{\frakY_n}\xrightarrow{f_n^{\flat}}f_*\calM_{\frakX_n},\ \calM_{\frakY}\xrightarrow{j_n^{\flat}}\calM_{\frakY_n}\xrightarrow{g_n^{\flat}}g_*\calM_{\frakX_n}$$
are equal. And since $j_n^{\flat}$ is surjective (\ref{prop69}), $f_n^{\flat}=g_n^{\flat}.$ We conclude that $f_n=g_n$ as morphisms of logarithmic schemes and it follows by \ref{prop18} that $f=g.$ The functor $F$ is thus faithful.

Now let $\varphi: \op{Hom}(-,\frakX) \ra \op{Hom}(-,\frakY)$ be a morphism of presheaves on $\boldsymbol{\op{L}}.$ For every integer $n\ge 1,$ let $h_n=\varphi(i_n):\frakX_n \ra \frakY.$ The morphisms $h_n$ are compatible so, by \ref{corkey}, they induce a morphism $h:\frakX \ra \frakY$ of logarithmic formal schemes. Let us prove that $\varphi=F(h).$ Let $T$ be an fs logarithmic scheme and $f:T\ra \frakX$ be a morphism. 
Let $\op{Spf}A$ be a formal affine subset of $\frakX$ and $\op{Spec}B$ an affine subset of $T$ such that $f(\op{Spec}B)\subset \op{Spf}A.$ Let $u:A\ra B$ be the corresponding continuous morphism of adic rings. Since $p^n\longrightarrow 0$ in $A,$ $u(p)^n\longrightarrow 0$ in $B.$ But since $0$ is an open ideal in $B,$ it follows that $u(p)^n=0$ for some positive integer $n.$ Hence the morphism of formal schemes $\op{Spec}B \xrightarrow{f} \frakX$ factors through $\frakX_n.$ And since $i_n:\frakX_n \ra \frakX$ is strict, this factorization holds for logarithmic formal schemes.
Suppose there exists a covering $(T_i\xrightarrow{f_i} T)_{i\in I}$ such that, for every $i\in I,$
\begin{equation}\label{eq721}
\varphi_{T_i}(f\circ f_i)=h \circ f\circ f_i
\end{equation}
Considering the commutative diagram
$$
\begin{tikzcd}
\op{Hom}(T,\frakX) \ar{r}{\varphi_T} \ar{d} & \op{Hom}(T,\frakY) \ar{d} \\
\op{Hom}(T_i,\frakX) \ar{r}{\varphi_{T_i}} & \op{Hom}(T_i,\frakY)
\end{tikzcd}
$$
where the vertical arrows are induced by $f_i.$
We get $\varphi_{T_i}(f\circ f_i)=\varphi_T(f) \circ f_i.$ By the hypothesis (\ref{eq721}), $h \circ f\circ f_i=\varphi_T(f)\circ f_i,$ and since $(f_i)_{i\in I}$ is a covering, we conclude that $h\circ f=\varphi_T(f).$ It is thus sufficient to work locally on $T$ and we may suppose that there exists $f_n:T \ra \frakX_n$ such that $f$ is equal to the composition
$$T\xrightarrow{f_n} \frakX_n \xrightarrow{i_n} \frakX.$$
Consider the commutative diagram
$$
\begin{tikzcd}
\op{Hom}(\frakX_n,\frakX) \ar{r}{\varphi_{\frakX_n}} \ar{d} & \op{Hom}(\frakX_n,\frakY) \ar{d} \\
\op{Hom}(T,\frakX) \ar{r}{\varphi_{T}} & \op{Hom}(T,\frakY)
\end{tikzcd}
$$
where the vertical arrows are induced by $f_n.$ We get $\varphi_T(f)=h_n\circ f_n=h\circ f.$ We conclude that $F$ is full.
\end{proof}

\begin{definition}\label{def72}
Let $\frakX$ be an fs logarithmic $p$-adic formal scheme and $\frakX\ra[P]$ a morphism of presheaves. We say that $\frakX\ra [P]$ is a \emph{frame} if, for any geometric point $\ov{x}\ra \frakX,$ there exists an étale neighborhood $\frakU$ of $\ov{x}$ such that $\frakU\ra [P]$ factors as a composition $\frakU\ra B\langle P\rangle \ra [P]$ with $\frakU\ra B\langle P\rangle$ strict and $B\langle P\rangle\ra [P]$ is the canonical morphism of presheaves (\ref{parag31}). A logarithmic formal scheme $\frakX$ equipped with a frame $\frakX\ra [P]$ will be denoted by $(\frakX,P)$ and called a \emph{framed logarithmic formal scheme}. A morphism $(\frakX,P)\ra (\frakS,Q)$ of framed logarithmic formal schemes is a pair of morphisms $(\frakX\ra \frakS,Q\ra P)$ such that the following diagram is commutative:
$$\begin{tikzcd}
\frakX\ar{r}\ar{d} & \frakS\ar{d}\\
\left [P\right ]\ar{r} & \left [Q\right ]
\end{tikzcd}$$
where $[P]\ra [Q]$ is the morphism induced by $Q\ra P.$ 
\end{definition}

\begin{parag}\label{parag73}
Let $(\frakX,Q)$ be a framed logarithmic formal scheme and $n\ge 1.$ Let $\frakX_n$ be the scheme $(\frakX,\Ox_{\frakX}/p^n)$ equipped with the logarithmic structure pullback of that of $\frakX.$ We get a frame on $\frakX_n$ by composing $\frakX\ra [Q]$ with the canonical morphism $\frakX_n\ra \frakX$ which is strict.
Let $\frakX\ra \frakS$ and $\frakY\ra \frakS$ be morphisms in $\boldsymbol{\op{LFS}},$ $Q\ra P$ a morphism of fs monoids and $\frakX\ra [Q]$ and $\frakY\ra [Q]$ frames. We denote by $\frakX\times_{[Q]}^{\op{log}}[P]$ and $\frakX\times^{\op{log}}_{\frakS,[Q]}\frakY$ the presheaves of sets
\begin{equation}\label{eq731}
\frakX\times_{[Q]}^{\op{log}}[P]:\begin{array}[t]{clc}\boldsymbol{\op{L}} & \ra & \boldsymbol{\op{Sets}}\\ T & \mapsto & \frakX(T)\times_{[Q](T)}[P](T),\end{array}
\end{equation}
\begin{equation}\label{eq732}
\frakX\times^{\op{log}}_{\frakS,[Q]}\frakY:\begin{array}[t]{clc}\boldsymbol{\op{L}} & \ra & \boldsymbol{\op{Sets}}\\ T & \mapsto & \frakX(T)\times_{\frakS(T)\times [P](T)}\frakY(T).\end{array}
\end{equation}
We have the following proposition, which is analogous to \ref{prop49}:
\end{parag}

\begin{proposition}\label{formalframelift}
Let $f:(\frakX,Q) \ra (\frakS,P)$ be a morphism of fs framed $p$-adic logarithmic formal schemes. Etale locally on $\frakX$ and $\frakS,$ there exists a chart $B\langle M\rangle \ra B\langle P \rangle$ of $f$ lifting the diagram
$$
\begin{tikzcd}
\frakX \ar{r} \ar{d} & \frakS \ar{d} \\
\left [Q \right ] \ar{r} & \left [P \right ].
\end{tikzcd}
$$
\end{proposition}

\begin{proposition}\label{formaletaleframelift}
Equip $\op{Spf}\Z_p$ with the trivial logarithmic structure. Let $f:(\frakX,Q) \ra (\frakS,P)$ be a log smooth morphism of fs framed $p$-adic logarithmic formal schemes such that $\frakS$ is log flat over $\op{Spf}\Z_p.$ Etale locally on $\frakX$ and $\frakS,$ there exists a chart $B\langle M \rangle \ra B\langle P \rangle$ of $f,$ fitting into a commutative diagram
$$
\begin{tikzcd}
\frakX \ar{r} \ar{d} & \frakS \ar{d} \\
B\langle M\rangle \ar{r} \ar{d} & B\langle P \rangle \ar{d} \\ 
\left [Q \right ] \ar{r} & \left [P\right ],
\end{tikzcd}
$$
and satisfying the following conditions:
\begin{enumerate}
\item $P^{gp} \ra M^{gp}$ is injective and the torsion subgroup of its cokernel is of finite order invertible in $\Ox_{\frakX}.$
\item The morphism $\frakX \ra \frakS\times_{B\langle P \rangle}B\langle M \rangle,$ induced by $f$ and the chart $\frakX \ra B\langle M\rangle,$ is étale and strict.
\end{enumerate}
\end{proposition}

\begin{proof}
Let $X$ and $S$ be the special fibers of $\frakX$ and $\frakS$ respectively and $f_1:X\ra S$ the morphism induced by $f.$ Let $\ov{x} \ra X$ be a geometric point. After shrinking $\frakS$ to an étale neighborhood of $f_1(\ov{x}),$ we can suppose that the frame $\frakS \ra [P]$ factors through a chart $\frakS \ra B\langle P \rangle.$ By \ref{etaleframelift}, there exists an étale neighborhood $U$ of $\ov{x}$ and a chart $A_1[M] \ra A_1[P]$ of the restriction $U\ra S$ of $f_1:X \ra S,$ fitting into a commutative diagram
$$
\begin{tikzcd}
U \ar{r}{f_1} \ar{d} & S \ar{d} \\
A_1[M] \ar{r} \ar{d} & A_1[P] \ar{d} \\
\left [Q \right ] \ar{r} & \left [ P \right ],
\end{tikzcd}
$$
such that the torsion subgroup of the cokernel of $P^{gp} \ra M^{gp}$ is of finite order invertible in $\Ox_{U}$ and the morphism $U \ra S\times_{A_1[P]}A_1[ M ],$ induced by $f_1$ and the chart $U\ra A_1[ M],$ is étale and strict. We can furthermore suppose that $U$ is affine.
Since $B\langle M\rangle \ra B\langle P \rangle$ is log smooth, the morphism $\frakS \times_{B\langle P \rangle}B\langle M\rangle \ra \frakS$ is log smooth modulo $p^n$ for every $n\ge 1.$ Let $\frakU \ra \frakX$ be an étale morphism such that $U=\frakU\times_{\frakX}X.$ There exists a morphism $\frakU\ra \frakS \times_{B\langle P \rangle} B\langle M\rangle$ fitting into the commutative diagram
$$
\begin{tikzcd}
S\times_{A_1[P]} A_1[M] \ar{rr} & & \frakS\times_{B\langle P\rangle} B\langle M\rangle \ar{d} \\
U\ar{u} \ar{r} & \frakU \ar{ur} \ar{r} & \frakS.
\end{tikzcd}
$$
Since the morphisms $U\ra \frakU,$ $U\ra S\times_{A_1[P]}A_1[M]$ and $S\times_{A_1[P]}A_1[M] \ra \frakS\times_{B\langle P\rangle}B\langle M\rangle$ are all strict, the morphism $\frakU \ra \frakS\times_{B\langle P \rangle}B\langle M\rangle$ is also strict. It remains to prove that the morphism $\frakU \ra \frakS\times_{B\langle P\rangle}B\langle M \rangle$ is étale. The morphism $B\langle M\rangle \ra B\langle P\rangle$ is log smooth hence log flat, thus so is $\frakS\times_{B\langle P\rangle} B\langle M\rangle \ra \frakS.$ Since $U \ra S\times_{A_1[P]}A_1[M]$ is étale, by (\cite{Ahmed2010} 2.4.11), there exists an étale strict morphism $\frakZ\ra \frakS\times_{B\langle P\rangle} B\langle M\rangle$ of $p$-adic formal schemes fitting into a cartesian diagram
\begin{equation}\label{diag771}
\begin{tikzcd}
U \ar{r} \ar{d} & \frakZ \ar{d} \\
S\times_{A_1[P]}A_1[M] \ar{r} & \frakS\times_{B\langle P\rangle} B\langle M \rangle.
\end{tikzcd}
\end{equation}
Then, by the étaleness of $\frakZ \ra \frakS\times_{B\langle P\rangle} B\langle M \rangle,$ there exists a unique morphism
$\frakU \ra \frakZ$ fitting into the commutative diagram
$$
\begin{tikzcd}
 & & \frakZ \ar{d} \\
U \ar{r} \ar[bend right=-20]{urr} & \frakU \ar{ur} \ar{r} & \frakS\times_{B\langle P\rangle} B\langle M \rangle.
\end{tikzcd}
$$
By \eqref{diag771}, $\frakU \ra \frakZ$ is an isomorphism modulo $p.$ By assumption, $\frakX$ and $\frakS$ are log flat over $\op{Spf}\Z_p.$ Since $\frakU$ and $\frakZ$ are log flat over $\frakX$ and $\frakS$ respectively, $\frakU$ and $\frakZ$ are log flat, hence flat \eqref{Wflat}, over $\op{Spf}\Z_p.$ We conclude by (\cite{DXU19} 7.2) that $\frakU \ra \frakZ$ is an isomorphism of formal schemes. Since it is clearly strict, it is an isomorphism of logarithmic formal schemes.
\end{proof}

\begin{proposition}\label{prop73}
Let $P$ and $Q$ be fs monoids, $\theta:Q\ra P$ a morphism of monoids such that $\theta^{gp}$ is surjective and $(\frakX,Q)$ a $p$-adic framed logarithmic formal scheme \eqref{parag73}. Then
\begin{enumerate}
\item The presheaf $\frakX\times_{[Q]}^{\op{log}}[P]$ \eqref{eq731} is representable by a logarithmic formal scheme which is affine and log étale over $\frakX.$
\item If $\frakX \ra [Q]$ lifts to a chart $\frakX\ra B\langle Q\rangle,$ then $\frakX\times_{[Q]}^{\op{log}}[P]$ is representable by $\frakX\times_{B\langle Q\rangle}^{\op{log}} B\langle \tilde{P}\rangle,$ where $\tilde{P}$ is the inverse image of $P$ by $Q^{gp}\ra P^{gp}.$
\end{enumerate}
\end{proposition}

\begin{proof}
The proof is similar to the case of schemes (\cite{Saito04} 4.2.1).
\end{proof}

\begin{corollaire}\label{prop74}
Let $(\frakX,Q)\ra (\frakS,P)\leftarrow (\frakY,Q)$ be morphisms of fs framed logarithmic $p$-adic formal schemes \eqref{def72}. The presheaf $\frakX\times^{\op{log}}_{\frakS,[Q]}\frakY$ \eqref{eq732} is representable by a logarithmic formal scheme which is log étale over $\frakX\times_{\frakS}^{\op{log}}\frakY.$
\end{corollaire}

\begin{proof}
We apply the previous proposition to the presheaf
$$\frakX\times_{\frakS,[Q]}^{\op{log}}\frakY=\frakX\times_{\frakS}^{\op{log}}\frakY\times_{[Q\oplus Q]}^{\op{log}}[Q].$$
\end{proof}

\begin{corollaire}\label{cor76}
Let $(\frakX,Q)\ra (\frakS,P)\leftarrow (\frakY,Q)$ be morphisms of fs framed logarithmic $p$-adic formal schemes \eqref{def72}. Then the projection $\frakX\times_{\frakS,[Q]}^{\op{log}}\frakY\ra \frakY$ is strict.
\end{corollaire}

\begin{proof}
By the proof of \ref{prop73}, the map $\frakX\times_{\frakS,[Q]}^{\op{log}}\frakY\ra [Q]$ is strict. Since $\frakY\ra [Q]$ is also strict, it follows that $\frakX\times_{\frakS,[Q]}^{\op{log}}\frakY\ra \frakY$ is strict. 
\end{proof}

\begin{proposition}\label{prop76}
Let $(\frakX,Q)\ra (\frakS,P)$ be a morphism of locally Noetherian fs framed logarithmic $p$-adic formal schemes \eqref{def72}. Then the diagonal immersion
$$\frakX \ra \frakX\times_{\frakS}^{\op{log}}\frakX$$
factors canonically into a strict immersion $\frakX\ra \frakX\times^{\op{log}}_{\frakS,[Q]}\frakX$ followed by a log étale morphism $\frakX\times^{\op{log}}_{\frakS,[Q]}\frakX\ra \frakX\times_{\frakS}^{\op{log}}\frakX.$
\end{proposition}

\begin{proof}
The fact that $\frakX\ra \frakX\times^{\op{log}}_{\frakS,[Q]}\frakX$ is strict and $\frakX\times^{\op{log}}_{\frakS,[Q]}\frakX\ra \frakX\times_{\frakS}^{\op{log}}\frakX$ is log étale is a consequence of \ref{cor76} and \ref{prop74}. 
\end{proof}

\begin{corollaire}\label{cor79}
Let $(\frakX,Q)\ra (\frakS,P)$ be a morphism of locally Noetherian fs framed logarithmic $p$-adic formal schemes. Denote by $\frakX_{\frakS,[Q]}^{r+1}$ the fiber product $$\frakX\times^{\op{log}}_{\frakS,[Q]}\frakX\times^{\op{log}}_{\frakS,[Q]}\hdots \times^{\op{log}}_{\frakS,[Q]}\frakX\ (r+1\ \op{times}).$$ Then the canonical morphism $\frakX^{r+1}_{\frakS,[Q]}\ra \frakX^{r+1}$ is étale and $\frakX\ra \frakX^{r+1}_{\frakS,[Q]}$ is a strict immersion, so that the diagonal immersion $\frakX\ra \frakX^{r+1}$ factors as the composition of an étale morphism and a strict immersion.
\end{corollaire}

\begin{lemma}\label{Glogflat}
Let $W$ be the ring of Witt vectors of a perfect field $k$ of characteristic $p$ and $(\frakX,Q)\ra (\frakS,P)$ a log smooth morphism of framed fs $p$-adic logarithmic formal schemes, log flat and locally of finite type over $\op{Spf}W.$ Note that this implies, by \ref{Wflat}, that $\frakX$ and $\frakS$ are flat over $\op{Spf}W$ \eqref{dxuflat}. Let $X$ and $S$ be the special fibers of $\frakX$ and $\frakS$ respectively and $F_1:X \ra X'$ the exact relative Frobenius morphism \eqref{PFrob}.
Recall the monoid $Q'$ and the morphism $F_{Q/P}:Q' \ra Q$ given in \ref{cor411}. Suppose that there exists a morphism $F:(\frakX,Q) \ra (\frakX',Q')$ of framed logarithmic formal schemes, over $(\frakS,P),$ lifting the exact relative Frobenius $F_1:X \ra X'$ and such that $\frakX'\ra \frakS$ is log smooth. Then $F$ and the morphism
$$G:\frakY=\frakX\times_{\frakS,[Q]}^{\op{log}}\frakX \ra \frakX'\times_{\frakS,[Q']}^{\op{log}}\frakX'=\frakY',$$
induced by $F,$ are log flat.
\end{lemma}

\begin{proof}
The lifting $F$ is log flat by \ref{thmlogflat} and \ref{logflatfiber}. Note that $\frakY$ and $\frakY'$ are log étale over $\frakX\times_{\frakS}^{\op{log}}\frakX$ and $\frakX'\times_{\frakS}^{\op{log}}\frakX'$ respectively \eqref{prop76}. The morphisms $\frakX \ra \frakS$ and $\frakX' \ra \frakS$ are log smooth hence so are the projections $\frakX\times_{\frakS}^{\op{log}}\frakX \ra \frakX$ and $\frakX'\times_{\frakS}^{\op{log}}\frakX' \ra \frakX'.$ Since the projections $\frakY \ra \frakX$ and $\frakY' \ra \frakX'$ are strict \eqref{cor76}, we get that the canonical morphisms $\frakY \ra \frakS$ and $\frakY' \ra \frakS$ are log smooth, hence $\frakY$ and $\frakY'$ are locally of finite presentation and log flat over $\op{Spf}W.$ Consider $X$ as a logarithmic scheme over $[Q']$ by the composition
$$X \ra [Q] \xrightarrow{[F_{Q/P}]} [Q'].$$
We have the commutative diagram
$$
\begin{tikzcd}
X \times_{S,[Q]}^{\op{log}}X \ar{r} \ar{dr} & X \times_{S,[Q']}^{\op{log}}X\ar{d} \\
 & X \times_{S}^{\op{log}}X.
\end{tikzcd}
$$
By \ref{rep1}, the morphisms $X \times_{S,[Q]}^{\op{log}}X \ra X \times_{S}^{\op{log}}X$ and $X \times_{S,[Q']}^{\op{log}}X \ra X \times_{S}^{\op{log}}X$ are log étale so $X \times_{S,[Q]}^{\op{log}}X \ra X \times_{S,[Q']}^{\op{log}}X$ is log étale. It remains to prove that
$$G_1=F_1\times_{S,[Q']} F_1:X \times_{S,[Q']}^{\op{log}}X \ra X' \times_{S,[Q']}^{\op{log}}X'$$
is log flat. This morphism $F_1\times_{S,[Q']} F_1$ decomposes as
$$X\times_{S,[Q']}^{\op{log}}X \xrightarrow{\op{Id} \times F_1} X\times_{S,[Q']}^{\op{log}}X' \xrightarrow{F_1\times \op{Id}} X'\times_{S,[Q']}^{\op{log}}X'.$$
By \ref{thmlogflat} and the cartesian squares
$$
\begin{tikzcd}
X\times_{S,[Q']}^{\op{log}}X \ar{r} \ar{d} & X\ar{d}{F_1} & X\times_{S,[Q']}^{\op{log}}X' \ar{r} \ar{d} & X\ar{d}{F_1} \\
X\times_{S,[Q']}^{\op{log}}X' \ar{r}{q_2} & X' & X'\times_{S,[Q']}^{\op{log}}X' \ar{r}{p_1} & X',
\end{tikzcd}
$$
where $p_1$ is the first projection and $q_2$ is the second projection, we deduce that $G_1:X \times_{S,[Q']}^{\op{log}}X \ra X' \times_{S,[Q']}^{\op{log}}X'$ is log flat. Since $\frakY$ is log flat over $\op{Spf}W,$ we conclude by \ref{logflatfiber}.
\end{proof}

\begin{parag}\label{parag77}
Let $(\frakX,Q)\ra (\frakS,P)$ be a morphism of locally Noetherian fs framed logarithmic $p$-adic formal schemes \eqref{def72}. Then, by \ref{prop76}, the diagonal immersion $\frakX \ra \frakX \times_{\frakS}^{\op{log}} \frakX$ factors into a strict closed immersion $\Delta:\frakX \ra \frakY$ followed by a log étale morphism $\frakY\ra \frakX \times_{\frakS}^{\op{log}} \frakX.$ For every integer $n\ge 1,$ denote by $\frakX_n$ (resp. $\frakY_n$) the logarithmic scheme obtained from $\frakX$ (resp. $\frakY$) by reduction modulo $p^n,$ as in \ref{prop69}, and denote by $\Delta_n:\frakX_n \ra \frakY_n$ the morphism induced by $\Delta$ and by $\calI$ (resp. $\calI_n$) the ideal of the immersion $\Delta$ (resp. $\Delta_n$). By \ref{prop45}, we have a canonical isomorphism $$\Delta_n^{-1} \left (\calI_n/\calI_n^2 \right ) \xrightarrow{\sim} \omega^1_{\frakX_n / \frakS_n}.$$
It follows, by (\ref{eq681}), that we have a canonical isomorphism
\begin{equation}\label{isoomega1}
\Delta^{-1} \left (\calI / \calI^2 \right ) \xrightarrow{\sim} \omega^1_{\frakX / \frakS}.
\end{equation}
Denote by $p_1,p_2:\frakY \ra \frakX$ and $p_{1,n},p_{2,n}:\frakY_n \ra \frakX_n$ the canonical projections. Note that $\Delta^{-1}(1+\calI) \subset \Delta^{-1}\Ox_{\frakY}^*.$
We have the following commutative diagram
$$
\begin{tikzcd}
0 \ar{r} & \Delta^{-1}(1+\calI) \ar{r}{\lambda} \ar{d} & \Delta^{-1} \calM_{\frakY} \ar{d}\ar{r}{\Delta^{\flat}} & \calM_{\frakX} \ar{d}\ar{r} & 0 \\
0 \ar{r} & \Delta_n^{-1}(1+\calI_n) \ar{r}{\lambda_n} & \Delta_n^{-1} \calM_{\frakY_n} \ar{r}{\Delta_n^{\flat}} & \calM_{\frakX_n}\ar{r} & 0
\end{tikzcd}
$$
with exact rows (\ref{era2exactseq}) and where $\lambda$ (resp. $\lambda_n$) is induced by $\alpha_{\frakY}^{-1}$ (resp. $\alpha_{\frakY_n}^{-1}$). Given a local section $m\in \Gamma(\frakU,\calM_{\frakX})$ (resp. $m\in \Gamma(U,\calM_{\frakX_n})$) over an étale $\frakX$-formal scheme $\frakU$ (resp. an étale $X$-scheme U), we denote by $\mu(m)$ (resp. $\mu_n(m)$) the unique section of $\Gamma(\frakU,\Delta^{-1}(1+\calI))$ (resp. $\Gamma(U,\Delta_n^{-1}(1+\calI_n))$) such that $(\Delta^{-1}p_1^{\flat})m+\lambda(\mu(m))=(\Delta^{-1}p_2^{\flat})m$ (resp. $(\Delta_n^{-1}p_{1,n}^{\flat})m+\lambda_n(\mu_n(m))=(\Delta_n^{-1}p_{2,n}^{\flat})m$) and $\eta(m)=\mu(m)-1\in \Delta^{-1}\calI$ (resp. $\eta_n(m)=\mu_n(m)-1\in \Delta_n^{-1}\calI_n$).
 
Suppose $\frakX\ra \frakS$ is smooth and suppose we are given local coordinates $m_{1,1},\hdots,m_{1,d}\in \Gamma(U,\calM_{\frakX_1})$ of $\frakX_1$ over $\frakS_1$ (\ref{P2}), where $U$ is an étale $\frakX_1$-scheme. By \ref{prop69}, after eventually shrinking $U,$ there exists an étale $\frakX$-scheme $\frakU$ and, for every $1\le i\le d,$ a local section $m_i\in \Gamma(\frakU,\calM_{\frakX})$ that lifts $m_{1,i}.$ Since $\eta_1(m_{1,1}),\hdots,\eta_1(m_{1,d}),$ modulo $\calI_1^2,$ form a local basis for $\omega^1_{\frakX_1 / \frakS_1},$ the sections $\eta(m_1),\hdots,\eta(m_d)$ modulo $\calI^2$ also form a local basis for $\omega^1_{\frakX / \frakS}.$ Such local sections $m_i$ will be refered to as \emph{local coordinates for $\frakX$ over $\frakS$}. Note that $\eta(m_1),\hdots,\eta(m_d)$ locally generate the ideal $\calI.$
\end{parag}

\begin{proposition}\label{era2prop611}
Keep the same hypothesis and notation of \ref{parag77}. We identify $\frakY \times_{\frakX}^{\op{log}}\frakY$ with $\frakY(2):=\frakX \times_{\frakS,[Q]}^{\op{log}} \frakX\times_{\frakS,[Q]}^{\op{log}} \frakX$ via the canonical isomorphism. Let $\w{\delta}$ be the morphism of structural rings associated with the $(1,3)$-projection $p_{13}:\frakY(2) \rightarrow \frakY.$ Let $\Delta(2): \frakX \ra  \frakY(2)$ be the diagonal embedding. Then, for any local section $m$ of $\calM_{\frakX},$
$$(\Delta(2)^{-1}\w{\delta})(\eta(m))=\eta(m)\otimes \eta(m)+\eta(m)\otimes 1+1\otimes \eta(m).$$
\end{proposition}

\begin{proof}
For $1\le i\le 3$ and $1\le j\le 2,$ denote by $\pi_i:\frakY(2) \ra \frakX$ and $p_j:\frakY \ra \frakX$ the projections on the $i$th and $j$th factor respectively and, for $1\le i<j\le 3,$ denote by $p_{ij}:\frakY(2) \ra \frakX$ the projection on the $(i,j)$th factor. 
Let $\calI$ and $\calI(2)$ be the ideals of the exact immersions $\Delta:\frakX \ra \frakY$ and $\Delta(2):\frakX \ra \frakY(2)$ respectively. By \eqref{era2exactseq}, we have exact sequences fitting into the commutative diagram
$$
\begin{tikzcd}
0 \ar{r} & \Delta^{-1}(1+\calI) \ar{r}{\lambda} \ar{d}{\Delta(2)^{-1}p_{ij}^{\#}} & \Delta^{-1}\calM_{\frakY} \ar{r}{\Delta^{\flat}} \ar{d}{\Delta(2)^{-1}p_{ij}^{\flat}} & \calM_{\frakX} \ar{r} \ar[equal]{d} & 0 \\
0 \ar{r} & \Delta(2)^{-1}(1+\calI(2)) \ar{r}{\lambda_2} & \Delta(2)^{-1}\calM_{\frakY(2)} \ar{r}{\Delta(2)^{\flat}} & \calM_{\frakX} \ar{r} & 0.
\end{tikzcd}
$$
For $1\le i<j\le 3,$ let $\eta_{ij}(m)$ be the unique local section of $\Delta(2)^{-1}(\calI(2))$ satisfying
$$
\lambda_2(1+\eta_{ij}(m))+\pi_i^{\flat}(m)=\pi_j^{\flat}(m).
$$
We drop the notation $\Delta(2)^{-1}$ and $\Delta^{-1}$ to lighten the notation.
We have
\begin{alignat*}{2}
\lambda_2(p_{13}^{\#}(1+\eta(m)))+\pi_1^{\flat}(m) &= p_{13}^{\flat}(\lambda(1+\eta(m)))+p_{13}^{\flat}p_1^{\flat}(m) \\
&= p_{13}^{\flat}(p_2^{\flat}(m)) \\
&= \pi_3^{\flat}(m) \\
&= \lambda_2(1+\eta_{23}(m))+\lambda_2(1+\eta_{12}(m))+\pi_1^{\flat}(m) \\
&= \lambda_2\left (1+\eta_{12}(m)+\eta_{23}(m)+\eta_{12}(m)\eta_{23}(m) \right )+\pi_1^{\flat}(m).
\end{alignat*}
Since $\calM_{\frakY}$ is integral and $\lambda_2$ is injective, we get
$$p_{13}^{\#}(\eta(m))=\eta_{12}(m)+\eta_{23}(m)+\eta_{12}(m)\eta_{23}(m).$$
We conclude by the following fact: let $q_1,q_2:\frakY\times_{\frakX}\frakY \ra \frakY$ be the canonical projections and $\Delta':\frakX \ra \frakY\times_{\frakX}\frakY$ the morphism induced by $\Delta:\frakX \ra \frakY.$ If we identify $\frakY\times_{\frakX}\frakY$ and $\frakY(2)$ via the canonical isomorphism then $\Delta'$ identifies with $\Delta(2),$ $\eta_{12}(m)$ identifies with $(\Delta'^{-1}q_1^{\#})(\eta(m))$ and $\eta_{23}(m)$ identifies with $(\Delta'^{-1}q_2^{\#})(\eta(m)).$
\end{proof}

\begin{proposition}\label{era2prop612}
Keep the same hypothesis and notation of \ref{parag77}. Let $\sigma:\frakY \ra \frakY$ be the morphism exchanging the $\frakX$ factors. For any local section $m$ of $\calM_{\frakX},$ we have
$$(\Delta^{-1}\sigma^{\#})(\eta(m))=(1+\eta(m))^{-1}-1.$$ 
\end{proposition}

\begin{proof}
We have the following commutative diagram
$$
\begin{tikzcd}
\Delta^{-1}(1+\calI) \ar{r}{\lambda} \ar[swap]{d}{\Delta^{-1}\sigma^{\#}} & \Delta^{-1}\calM_{\frakY} \ar{d}{\Delta^{-1}\sigma^{\flat}} \\
\Delta^{-1}(1+\calI) \ar{r}{\lambda} & \Delta^{-1}\calM_{\frakY}.
\end{tikzcd}
$$
We get
\begin{alignat*}{2}
\lambda\left (\Delta^{-1}\sigma^{\#}(\mu(m))\right ) + \left (\Delta^{-1}p_2^{\flat}\right )(m) &= \lambda\left (\Delta^{-1}\sigma^{\#}(\mu(m))\right ) + \left (\Delta^{-1} \sigma^{\flat} \right ) \left (\Delta^{-1}p_1^{\flat} \right )(m) \\
&= \left (\Delta^{-1}\sigma^{\flat} \right ) (\lambda(\mu(m))) + \left (\Delta^{-1} \sigma^{\flat} \right ) \left (\Delta^{-1}p_1^{\flat} \right )(m) \\
&= \left (\Delta^{-1}\sigma^{\flat} \right ) \left ( \Delta^{-1}p_2^{\flat} \right )(m) \\
&= \left (\Delta^{-1}p_1^{\flat} \right )(m)  \\
&= \lambda\left (\mu(m)^{-1} \right ) + \left (\Delta^{-1}p_2^{\flat}\right )(m).
\end{alignat*}
The morphism $\lambda$ is injective and $\Delta^{-1}\calM_{\frakY}$ is integral, it follows that
$$
\left (\Delta^{-1}\sigma^{\#}\right )(\mu(m))=\mu(m)^{-1}.
$$
\end{proof}

\newpage

\section{Indexed objects}

\subsection*{Indexed algebras and modules}

\begin{parag}\label{indpar1}
In this section, we review the theory of indexed objects introduced by Lorenzon in \cite{Lor2000}. We fix a site $C$ in which all fiber products are representable and consider the associated topos $X$ and a sheaf $\calI$ of $X.$ We denote by $C_{/\calI}$ the category of objects of $C$ over $\calI$ and consider the functor
\begin{equation}
j_{\calI}:\begin{array}[t]{clc} C_{/\calI} & \ra & C \\ (A\ra \calI) & \mapsto & A. \end{array}
\end{equation}
We equip $C_{/\calI}$ with the topology induced by $j_{\calI}$ (\cite{SGA43} exposé III 3.1). Then $j_{\calI}$ is continuous and cocontinuous (\cite{SGA43} exposé III 5.2) and so it induces three functors
\begin{equation}\label{indeq1}
j_{\calI !}:\widetilde{C_{/\calI}}\ra \widetilde{C},\ j_{\calI}^*:\widetilde{C}\ra \widetilde{C_{/\calI}},\ j_{\calI *}:\widetilde{C_{/\calI}}\ra \widetilde{C}
\end{equation}
such that $(j_{\calI *},j_{\calI}^*)$ is a morphism of topoi, and $j_{\calI !}$ is a left adjoint of $j_{\calI}^*$ (\cite{SGA43} exposé III 5.3). By (\cite{SGA43} exposé III 5.4), the functor $j_{\calI !}$ factors as
$$\widetilde{C_{/\calI}} \xrightarrow{e} \widetilde{C}_{/\calI} \ra \widetilde{C}$$
where $e$ is an equivalence of categories and $\widetilde{C}_{/\calI}\ra \widetilde{C}$ is the canonical functor. By identifying $\widetilde{C_{/\calI}}$ with $\widetilde{C}_{/\calI}$ via $e,$ the three functors given above (\ref{indeq1}) correspond to functors:
\begin{equation}\label{indeq2}
j_{\calI !}:\widetilde{C}_{/\calI}\ra \widetilde{C},\ j_{\calI}^*:\widetilde{C}\ra \widetilde{C}_{/\calI},\ j_{\calI *}:\widetilde{C}_{/\calI}\ra \widetilde{C}
\end{equation}
We denote simply by $j_{\calI}:\widetilde{C}_{/\calI} \ra \widetilde{C}$ the morphism $(j_{\calI*},j_{\calI}^*),$ called the localisation morphism of $\widetilde{C}$ at $\calI.$
By (\cite{SGA43} exposé III 5.4), the functor $j_{\calI}^*$ is given, for any sheaf $\F$ of $\widetilde{C},$ by 
\begin{equation}\label{indeq8}
j_{\calI}^*(\F)=(\F\times \calI \ra \calI),
\end{equation}
where $\F\times \calI \ra \calI$ is the canonical projection.

A morphism $f:\calJ\ra \calI$ in $\widetilde{C}$ defines an object of $\widetilde{C}_{/\calI}$ and so we have a functor
$$\begin{array}[t]{clc} (C_{/\calI})_{/f} & \ra & C_{/\calI}\\ (g\ra f) & \mapsto & g \end{array}$$
that we canonically identify with the functor abusively denoted by
$$
f:\begin{array}[t]{clc}
C_{/\calJ} & \ra & C_{/\calI} \\
g & \mapsto & f\circ g.
\end{array}
$$
Note that the topology on $C_{/\calJ}$ induced by $f$ is the same as the topology induced by the canonical functor $j_{\calJ}:C_{/\calJ} \ra C.$ Indeed, $j_{\calI}\circ f=j_{\calJ}.$ Since all fiber products in $C$ are representable, by (\cite{SGA43} Exposé II 1.3.1), any topology is determined by a pretopology. By (\cite{SGA43} Exposé III 3.3), for a family $F=(Y_i\ra Y)_{i\in I}$ of morphisms of $C_{/\calJ},$ $F$ is a covering family for the topology induced by $j_{\calJ}$ on $C_{/\calJ}$ if and only if $j_{\calJ}(F)=(j_{\calJ}(Y_i)\ra j_{\calJ}(Y))_{i\in I}$ is a covering if and only if $f(F)=(f(Y_i)\ra f(Y))_{i\in I}$ is a covering for the topology induced by $j_{\calJ}$ on $C_{/\calI}.$

We get three functors
\begin{equation}\label{indeq3}
f_!:\widetilde{C}_{/\calJ}\ra \widetilde{C}_{/\calI},\ f^*:\widetilde{C}_{/\calI} \ra \widetilde{C}_{/\calJ},\ f_*:\widetilde{C}_{/\calJ} \ra \widetilde{C}_{/\calI}
\end{equation}
such that $(f_*,f^*)$ is the morphism of localisation of $\widetilde{C}_{/\calI}$ at $f$ and $f_!$ is left adjoint to $f^*.$

If $U$ is an object of $C$ and $i\in\calI(U),$ then, by (\ref{indeq3}), we get three functors
\begin{equation}\label{indeq6}
i_!:\widetilde{C}_{/U}\ra \widetilde{C}_{/\calI},\ i^*:\widetilde{C}_{/\calI} \ra \widetilde{C}_{/U},\ i_*:\widetilde{C}_{/U} \ra \widetilde{C}_{/\calI}
\end{equation} 
From now on, for every object $\F$ of $\widetilde{C}_{/\calI},$ we denote $i^*\F$ by $\F_i.$
\end{parag}

\begin{parag}\label{indfiber73}
Suppose that $C$ is the small étale site on a scheme $S,$ $\calI$ an étale sheaf on $S$ and $i\in\calI(U)$ a section of $\calI$ over an étale $S$-scheme $U.$ Denote by $\calI_{\text{ét}}$ the topos $S_{\text{ét}/\calI.}$ We have a canonical equivalence of topoi $U_{\text{ét}}\xrightarrow{\sim} S_{\text{ét}/U}.$ Then, by (\ref{indeq6}), we get three morphisms
\begin{equation}\label{indeq4}
i_!:U_{\text{ét}}\ra \calI_{\text{ét}},\ i^*:\calI_{\text{ét}} \ra U_{\text{ét}},\ i_*:U_{\text{ét}} \ra \calI_{\text{ét}}
\end{equation}
If $\F\ra \calI$ is an object of $\calI_{\text{ét}},$ then $i^*(\F\ra \calI)=\F\times_{\calI}U.$ We will denote it by $\F_i$ or $\F_{|U}.$
\end{parag}

\begin{remark}\label{indremark73kraz} 
Let $X$ be a topos, $f:\calJ\ra \calI$ a morphism of sheaves of $X$ and equip $X$ with its canonical topology and $X_{/\calI}$ with the topology induced by the functor $j_{\calI}:X_{/\calI} \ra X.$ Then $\widetilde{X}$ is canonically equivalent to $X$ (\cite{SGA43} Exposé IV 1.2) and by applying the results of \ref{indpar1} to the site $C=X,$ we get functors
\begin{equation}
j_{\calI !}:X_{/\calI}\ra X,\ j_{\calI}^*:X\ra X_{/\calI},\ j_{\calI *}:X_{/\calI}\ra X
\end{equation}
and
\begin{equation}
f_!:X_{/\calJ}\ra X_{/\calI},\ f^*:X_{/\calI} \ra X_{/\calJ},\ f_*:X_{/\calJ} \ra X_{/\calI}.
\end{equation}
\end{remark}

\begin{parag}\label{indremark73}
Let $C$ be a $\mathbb{U}$-site, $X$ the associated topos and $\calI$ a sheaf of $X.$ Let $\calT$ be the split fibered $\left (X_{/\calI} \right )$-category defined, for every sheaf $\F$ of $X_{/\calI},$ by the localized category $\calT_{/\F}$ and, for every morphism $\varphi:\F \ra \calG,$ by the inverse image functor of the morphism of topoi
\begin{equation*}
\left (X_{/\calI} \right )_{/\F} \ra \left (X_{/\calI} \right )_{/\calG},
\end{equation*}
induced by $\varphi.$ By (\cite{Giraud} II 3.4.4), $\calT$ is a stack. If $i$ is a local section of $\calI$ over an object $U$ of $C,$ then we have a canonical isomorphism of topoi
$$
\left (X_{/\calI} \right )_{/i}\xrightarrow{\sim} X_{/U}.
$$
The data of a sheaf $\F$ of $X_{/\calI}$ is then, equivalent to the data, for every local section $i$ of $\calI,$ of a sheaf $\F_i$ of $X_{/U},$ and, for every morphism $g:V\ra U$ of $C$ and every local section $i$ of $\calI$ over $U,$ of an isomorphism
$$
c_{i,g}:g^{-1}\F_{i} \xrightarrow{\sim} \F_{i_{|V}},
$$
satisfying cocycle conditions and such that $c_{i,\op{Id}_U}=\op{Id}_{\F_i}.$ In particular, the data of a morphism $f:\F \ra \calG$ in $X_{/\calI}$ is equivalent to the data, for every section $i\in\calI(U)$ over an object $U$ of $C,$ of a morphism $f_i:\F_i \ra \calG_i$ satisfying the following compatibility condition: for all morphisms $U\ra V$ of $C$ and sections $i\in \calI(U)$ and $j\in \calI(V)$ such that $j_{|U}=i,$ the following diagram is commutative:
$$
\begin{tikzcd}
\F_{j|U} \ar{r}{f_{j|U}} \ar[swap,sloped]{d}{\sim} & \calG_{j|U}\ar[sloped]{d}{\sim} \\
\F_i\ar{r}{f_i} & \calG_i
\end{tikzcd}
$$
where the vertical arrows are the restriction morphisms, which are isomorphisms.
\end{parag}

\begin{parag}\label{indpar2}
Let $(X,\Ox)$ be a ringed topos and $f:\calJ\ra \calI$ a morphism of $X$ and consider the functors
\begin{equation}\label{indeq5}
j_{\calI}^*:X\ra X_{/\calI},\ j_{\calI *}:X_{/\calI}\ra X,
\end{equation}
\begin{equation}\label{indeq7}
f^*:X_{/\calI} \ra X_{/\calJ},\ f_*:X_{/\calJ} \ra X_{/\calI},
\end{equation}
defined in \ref{indremark73kraz}. Then $(X,j_{\calI}^*\Ox)$ and $(X,j_{\calJ}^*\Ox)$ are ringed topoi and $(j_{\calI *},j_{\calI}^*)$ and $(f_*,f^*)$ are morphisms of ringed topoi. Note that there is no difference between $j_{\calI}^*$ and $j_{\calI}^{-1}$ (resp. $f^*$ and $f^{-1}$). In the rest of this section, we denote $j_{\calI}^*\Ox$ by $\Ox_{\calI}.$ Denote by $\boldsymbol{\op{Mod}}_{\Ox}$ (resp. $\boldsymbol{\op{Mod}}_{\Ox_{\calI}}$) the category of $\Ox$-modules (resp. $\Ox_{\calI}$-modules). The functor $j_{\calI}^*:\boldsymbol{\op{Mod}}_{\Ox}\ra \boldsymbol{\op{Mod}}_{\Ox_{\calI}}$ (resp. $f^*:\boldsymbol{\op{Mod}}_{\Ox_{\calJ}}\ra \boldsymbol{\op{Mod}}_{\Ox_{\calI}}$) has a left adjoint, that we denote by $j_{\calI !}$ (resp. $f_!$). Note that $j_{\calI !}$ and $f_!$ are different from the functors $j_{\calI !}$ and $f_!$ defined in \ref{indpar1}. Otherwise stated, we will only use the notions for modules. 

Consider the canonical projections $p_1:\calI \times \calJ \ra \calI$ and $p_2:\calI\times \calJ \ra \calJ.$ Let $\calE$ (resp. $\F$) be an $\Ox_{\calI}$-module (resp. $\Ox_{\calJ}$-module). We denote by $\calE \boxtimes \F$ \emph{the exterior tensor product of $\calE$ and $\F$}, i.e. the $\Ox_{\calI\times \calJ}$-module defined as follows:
\begin{equation}
\calE \boxtimes \F := p_1^*\calE \otimes_{\Ox_{\calI\times \calJ}}p_2^*\F.
\end{equation}
\end{parag}

\begin{proposition}\label{indpropshriek}
Let $C$ be a $\mathbb{U}$-site, $X$ the associated topos, $\Ox$ a sheaf of rings on a $C$ and $f:\calJ\ra \calI$ a morphism of sheaves of $X.$ Then:
\begin{enumerate}
\item For any $\Ox_{\calI}$-module $\F,$ $j_{\calI !}\F$ is the sheaf associated to the presheaf
$$U \mapsto \bigoplus_{u\in \op{Hom}(U,\calI)}\F(u),$$
where we view $\F$ as a sheaf on the localised site $C_{/\calI}.$
\item For any $\Ox_{\calJ}$-module $\F,$ $f_!\F$ is the sheaf associated to the presheaf
$$(i:U\ra \calI) \mapsto \bigoplus_{\substack{j\in \op{Hom}(U,\calJ)\\ f\circ j=i}}\F(j)$$
where $U$ is an object of $C$ and $\F$ is viewed as a sheaf on the localised site $C_{/\calJ}.$
\item Suppose that $\calI$ is a sheaf of monoids and let $\sigma:\calI^2 \ra \calI$ be the addition map and $\F$ an $\Ox_{\calI^2}$-module. Then $\sigma_!\F$ is the sheaf associated to the presheaf
$$(s:U \ra \calI) \mapsto \bigoplus_{\substack{u,v\in \Gamma(U,\calI)\\ u+v=s}}\F((u,v)).$$
\end{enumerate}
\end{proposition}

\begin{proof}
We just prove the first assertion as the second is a special case of the first and the third is an application of the second.
For any $\Ox_{\calI}$-module $\F$ of $X_{/\calI},$ denote by $v(\F)$ the sheaf associated to the presheaf
$$U \mapsto \bigoplus_{u\in \op{Hom}_C(U,\calI)}\F(u).$$
We want to prove that $v$ is left adjoint to $j_{\calI}^*.$ Let $\F$ (resp. $\calG$) be an $\Ox_{\calI}$-module (resp. $\Ox$-module). Then the data of a morphism $g:v(\F)\ra \calG$ is equivalent to the data of compatible morphisms $$g_U:v(\F)(U) \ra \calG(U)$$
for every object $U$ of $C.$ This is equivalent to the data of compatible morphisms
$$\F(u)\ra \calG(U)$$
for every object $U$ of $C$ and every morphism $u:U \ra \calI.$ Since $(j_{\calI}^*\calG)(u)=\calG(U)$ (\ref{indeq8}), we conclude that this is equivalent to the data of compatible morphisms
$$\F(u) \ra (j_{\calI}^*\calG)(u)$$
for every object $u$ of $C_{/\calI}.$
\end{proof}

\begin{proposition}\label{indpropshriekoh}
Let $C$ be a $\mathbb{U}$-site, $X$ the associated topos and $\Ox$ a sheaf of rings on a $C.$
\begin{enumerate} 
\item Let $\calI$ and $\calJ$ be sheaves of $X.$ Denote by $p:\calI\times \calJ \ra \calI$ and $q:\calI\times \calJ \ra \calJ$ the canonical projections. Then, for any $\Ox_{\calI}$-module $\F$ of $X_{/\calI},$ there exists a canonical isomorphism
$$q_!p^*\F\xrightarrow{\sim}j_{\calJ}^*j_{\calI !}\F.$$
\item Let $f:\calJ \ra \calI$ and $g:\calK \ra \calI$ be morphisms of sheaves of $X.$ Denote by $p:\calJ\times_{\calI} \calK \ra \calJ$ and $q:\calJ\times_{\calI} \calK \ra \calK$ the canonical projections. Then, for any $\Ox_{\calJ}$-module $\F$ of $X_{/\calJ},$ there exists a canonical isomorphism
$$q_!p^*\F\xrightarrow{\sim}g^*f_!\F.$$
\end{enumerate}
\end{proposition}

\begin{proof}
We have the following commutative diagram of ringed topoi
$$
\begin{tikzcd}
X_{/\calI\times\calJ}\ar{r}{p}\ar{d}{q}\ar{dr}{j_{\calI\times\calJ}} & X_{/\calI}\ar{d}{j_{\calI}} \\
X_{/\calJ}\ar{r}{j_{\calJ}} & X
\end{tikzcd}
$$
The adjunction morphism
$$p_!p^*\F \ra \F$$
induces a morphism
$$j_{\calJ !}q_! p^*\F=j_{\calI !}p_!p^*\F\ra j_{\calI !}\F.$$
By adjunction, it yields a morphism
\begin{equation}
q_! p^*\F=j_{\calI !}p_!p^*\F\ra j_{\calJ}^*j_{\calI !}\F.
\end{equation}
We check that it is an isomorphism using \ref{indpropshriek}. On one hand, $q_!p^*\F$ is the sheaf on $C_{/\calJ}$ associated to the presheaf
$$(j:U\ra \calJ) \mapsto \bigoplus\limits_{\substack{s:U\ra \calI\times \calJ \\ q\circ s=j}}(p^*\F)(s)=\bigoplus\limits_{i:U\ra \calI}\F(i).$$
On the other hand, $j_{\calJ}^*j_{\calI !}\F$ is the sheaf
$$(j:U\ra \calJ) \mapsto (j_{\calI !}\F)(U)$$ and thus it is associated to the presheaf
$$(j:U\ra \calJ) \mapsto \bigoplus\limits_{i:U\ra \calI}\F(i).$$
The second assertion is proved similary by considering the commutative diagram of ringed topoi
$$
\begin{tikzcd}
X_{/\calJ\times_{\calI}\calK}\ar{r}{p}\ar{d}{q}\ar{dr}{j_{\calJ\times_{\calI}\calK}} & X_{/\calJ}\ar{d}{f} \\
X_{/\calK}\ar{r}{g} & X_{/\calI}
\end{tikzcd}
$$
\end{proof}

\begin{corollaire}\label{indpropshriek2}
Let $C$ be a $\mathbb{U}$-site, $X$ the associated topos, $\Ox$ a sheaf of rings on a $C$ and $f:\calJ \ra \calI$ a morphism of sheaves of $X.$ For any module $\F$ of $X_{/\calJ}$ and any section $k\in \Gamma(U,\calI)$ over an object $U$ of $C,$ the sheaf $(j_{\calI !}\F)_k$ of $X_{/U}$ is associated to the presheaf
$$(i:V\ra U) \mapsto \bigoplus\limits_{\substack{j\in \Gamma(V,\calJ)\\ f\circ j=k\circ i}}\F(j),$$
where $V$ is an object of $C.$ 
\end{corollaire}

\begin{proof}
By the second assertion of \ref{indpropshriekoh} applied to $g=k,$ we have a canonical isomorphism
$$q_!p^*\F\xrightarrow{\sim}(f_!\F)_k,$$
where $p:U\times_{\calI}\calJ \ra \calJ$ and $q:U\times_{\calI}\calJ\ra U$ are the canonical projections. By \ref{indpropshriek}, the sheaf $q_!p^*\F$ is associated to the presheaf
$$(i:V\ra U) \mapsto \bigoplus\limits_{\substack{s:V\ra U\times_{\calI}\calJ\\ p\circ s=i}}(p^*\F)(s).$$
\end{proof}

\begin{parag}
Keep the same notations of \ref{indpar2} and suppose $\calI$ is a sheaf of commutative monoids. Consider the canonical projections 
$$p_{ij}:\calI^3\ra \calI^2,\ 1\le i<j\le 3,\ q_i:\calI^3\ra \calI,\ 1\le i\le 3.$$
For any $\Ox_{\calI}$-modules $\calE,$ $\F$ and $\calG,$ we have canonical isomorphisms of $\Ox_{\calI^3}$-modules
\begin{equation}
(\calE \boxtimes \F )\boxtimes \calG \xrightarrow{\sim} q_1^*\calE \otimes_{\Ox_{\calI^3}}q_2^*\F \otimes_{\Ox_{\calI^3}}q_3^*\calG \xleftarrow{\sim} \calE \boxtimes (\F \boxtimes \calG).
\end{equation}
We identify the triple exterior tensor products $(\calE \boxtimes \F) \boxtimes \calG$ and $\calE \boxtimes (\F \boxtimes \calG)$ via this canonical isomorphism and we denote them simply by $\calE\boxtimes \F \boxtimes \calG.$
We can now define $\calI$-indexed algebras and $\calI$-indexed modules:
\end{parag}

\begin{definition}\label{inddef}
Let $(X,\Ox)$ be a ringed topos, $\calI$ a sheaf of commutative monoids of $X,$ $p_1,p_2:\calI^2\ra \calI,$ $q_1,q_2,q_3:\calI^3\ra \calI$ and $p_{ij}:\calI^3\ra \calI^2,\ 1\le i<j\le 3$ the canonical projections, $\sigma:\calI^2\ra \calI$ and $\tau:\calI^3\ra \calI$ the addition maps, $e$ the final object of $X$ and $0:e\ra \calI$ the zero section and $\nu:\calI \ra e$ the canonical morphism.
\begin{enumerate}
\item An \emph{$\calI$-indexed $\Ox$-algebra} is an $\Ox_{\calI}$-module $\calA$ equipped with a $\Ox_{\calI^2}$-linear map
$$\pi:\calA \boxtimes \calA \ra \sigma^*\calA$$
and a global section $1_{\calA}:e\ra \calA_0,$ where $\calA_0$ is the $\Ox$-module $0^*\calA,$ such that the following diagrams are commutative:
\begin{enumerate}
\item Associativity:
$$\begin{tikzcd}
\calA \boxtimes \calA \boxtimes \calA \ar{r}{\pi \boxtimes \op{Id}_{\calA}} \ar[swap]{d}{\op{Id}_{\calA} \boxtimes \pi} & (\sigma \times\op{Id}_{\calI})^*(\calA \boxtimes \calA) \ar{d}{(\sigma \times\op{Id}_{\calI})^*\pi} \\
(\op{Id}_{\calI}\times \sigma)^*(\calA \boxtimes \calA) \ar{r}{(\op{Id}_{\calI}\times \sigma)^*\pi} & \tau^*\calA
\end{tikzcd}
$$
\item Unity:
$$
\begin{tikzcd}
\calA \ar{r}{\op{Id}_{\calA}\otimes \nu^*1_{\calA}} \ar{d}[swap]{\nu^*1_{\calA}\otimes \op{Id}_{\calA}} \ar{dr}{\op{Id}_{\calA}} & \calA \otimes_{\Ox_{\calI}}\nu^*\calA_0=(\op{Id}_{\calI}\times 0)^*(\calA \boxtimes \calA) \ar{d}{(\op{Id}_{\calI}\times 0)^*\pi} \\  \nu^*\calA_0\otimes_{\Ox_{\calI}}\calA=(0\times \op{Id}_{\calI})^*(\calA \boxtimes \calA) \ar[swap]{r}{(0\times \op{Id}_{\calI})^*\pi} & \calA
\end{tikzcd}
$$
where the identification $\calA \otimes_{\Ox_{\calI}}\nu^*\calA_0=(\op{Id}_{\calI}\times 0)^*(\calA \boxtimes \calA)$ follows from the commutative diagram
$$
\begin{tikzcd}
\calI=\calI\times e \ar{r}{\nu} \ar{d}{\op{Id}_{\calI}\times 0} \ar[swap,bend right=50]{dd}{\op{Id}_{\calI}} & e\ar{d}{0} \\
 \calI^2\ar{r}{p_2} \ar{d}{p_1} & \calI \\ \calI &  
\end{tikzcd}
$$
and $\nu^*\calA_0\otimes_{\Ox_{\calI}}\calA=(0\times \op{Id}_{\calI})^*(\calA \boxtimes \calA)$ follows from a similar diagram.
\end{enumerate}
We say that $\calA$ is \emph{commutative} if the diagram
$$
\begin{tikzcd}
\calA \boxtimes \calA \ar{r}{\pi} \ar{d}{s} & \sigma^*\calA \\
\calA \boxtimes \calA \ar{ur}{\pi} &
\end{tikzcd}
$$
where $s$ is the isomorphism $x\boxtimes y\mapsto y\boxtimes x,$ is commutative.
\item Let $\calA$ be an $\calI$-indexed $\Ox$-algebra. The center of $\calA$ is the subsheaf $\calC$ of $\calA$ consisting of local sections $a$ of $\calA$ such that $\pi(p_1^*a\otimes p_2^*b)=\pi(p_1^*b\otimes p_2^*a)$ for every local section $b$ of $\calA.$ It is immediate that $\calC$ is an $\calI$-indexed $\Ox$-subalgebra of $\calA.$
\item Let $\calA$ and $\calB$ be two $\calI$-indexed $\Ox$-algebras. A morphism of $\calI$-indexed $\Ox$-algebras $f:\calA \ra \calB$ is a $\Ox_{\calI}$-linear morphism $f:\calA \ra \calB$ such that the following diagram is commutative:
$$
\begin{tikzcd}
\calA \boxtimes \calA \ar{r} \ar{d}{f\boxtimes f} & \sigma^*\calA \ar{d}{\sigma^*f} \\
\calB\boxtimes \calB \ar{r} & \sigma^*\calB
\end{tikzcd}
$$
where the horizontal arrows are the multiplication morphisms.
\item If $f:\calA \ra \calB$ is a morphism of $\calI$-indexed $\Ox$-algebras, we say that $\calB$ is an \emph{$\calI$-indexed $\calA$-algebra.}
\item Let $\calA$ be an $\calI$-indexed $\Ox$-algebra. An \emph{$\calI$-indexed left $\calA$-module} is an $\Ox_{\calI}$-module $\calE$ equipped with a $\Ox_{\calI^2}$-linear map
$$\rho:\calA \boxtimes \calE \ra \sigma^*\calE,$$
such that the following diagrams are commutative
\begin{enumerate}
\item Associativity:
$$
\begin{tikzcd}
\calA\boxtimes \calA \boxtimes \calE \ar{r}{\op{Id}_{\calA}\boxtimes \rho} \ar{d}{\pi \boxtimes \op{Id}_{\calE}} & \calA \boxtimes \sigma^*\calE = (\op{Id}_{\calI}\times \sigma)^*(\calA \boxtimes \calE) \ar{d}{(\op{Id}_{\calI}\times \sigma)^*\rho} \\ 
\sigma^*\calA \boxtimes \calE = (\sigma \times \op{Id}_{\calI})^*(\calA \boxtimes \calE) \ar{r}{(\sigma \times \op{Id}_{\calI})^*\rho} & \tau^*\calE
\end{tikzcd}
$$
\item Unity:
$$
\begin{tikzcd}
\calE \ar{dr}{\op{Id}_{\calE}} \ar[swap]{d}{\nu^*1_{\calA}\otimes \op{Id}_{\calE}} & \\
\nu^*\calA_0\otimes_{\Ox_{\calI}}\calE = (0\times \op{Id}_{\calI})^* (\calA \boxtimes \calE) \ar[swap]{r}{(0\times \op{Id}_{\calI})^*\rho} & \calE
\end{tikzcd}
$$
\end{enumerate}
We also define a \emph{$\calI$-indexed right $\calA$-module} as an $\Ox_{\calI}$-module $\calE$ equipped with a $\Ox_{\calI^2}$-linear map
$$\calE \boxtimes \calA \ra \sigma^*\calE$$
satisfying similar conditions. Unless explicitly stated, we only consider indexed left module structures.

\item Let $\calA$ be an $\calI$-indexed $\Ox$-algebra and $\calE$ and $\F$ two $\calI$-indexed $\calA$-modules. A \emph{morphism of $\calI$-indexed $\calA$-modules $f:\calE \ra \F$} is a $\Ox_{\calI}$-linear morphism
$$f:\calE \ra \F$$
such that the following diagram is commutative
$$\begin{tikzcd}
\calA \boxtimes \calE \ar{r}{\op{Id}_{\calA} \boxtimes f} \ar{d} & \calA \boxtimes \F \ar{d} \\
\sigma^*\calE \ar{r}{\sigma^*f} & \sigma^* \F
\end{tikzcd}$$ 
where the vertical arrows are the maps defining the $\calA$-module structures on $\calE$ and $\F.$

For an $\calI$-indexed $\Ox$-algebra $\calA,$ we denote by $\boldsymbol{\op{Mod}}_{\calA}$ the category of $\calI$-indexed $\calA$-modules. The category $\boldsymbol{\op{Mod}}_{\calA}$ is an abelian subcategory of $\boldsymbol{\op{Mod}}_{\Ox_{\calI}}.$
\end{enumerate}
\end{definition}

\begin{remark}\label{era2indrem}
Let $(X,\Ox)$ be a ringed topos on a $\mathbb{U}$-site $C,$ $e$ the final object of $X$ and $\calI$ a sheaf of monoids of $X.$ Recall the notation \eqref{indfiber73}. For an object $U$ of $C,$ we denote by $\Ox_U$ the restriction of $\Ox$ to $U.$ By \ref{indremark73}, an $\calI$-indexed algebra structure on an $\Ox_{\calI}$-module $\calA$ is equivalent to the data of a global section $1:e\ra \calA_0$ and, for sections $s,t\in \Gamma(U,\calI)$ of $\calI$ over an object $U$ of $C,$ of an $\Ox_U$-linear morphism
$$\calA_s\otimes_{\Ox_U}\calA_t\ra \calA_{s+t},$$
such that, for any sections $r,s,t\in \Gamma(U,\calI)$ over an object $U$ of $C,$ the following properties are satisfied:
\begin{enumerate}
\item The diagram
$$
\begin{tikzcd}
\calA_r\otimes_{\Ox_U}\calA_s\otimes_{\Ox_U}\calA_t \ar{r}\ar{d} & \calA_{r+s}\otimes_{\Ox_U}\calA_t \ar{d} \\
\calA_r\otimes_{\Ox_U}\calA_{s+t} \ar{r} & \calA_{r+s+t}
\end{tikzcd}
$$
is commutative.
\item The composed morphisms
$$\calA_s \xrightarrow{x\mapsto x\otimes 1} \calA_s\otimes_{\Ox_U}\calA_0 \ra \calA_s$$
$$\calA_s \xrightarrow{x\mapsto 1\otimes x} \calA_0\otimes_{\Ox_U}\calA_s \ra \calA_s$$
are both equal to the identity.
\end{enumerate}
\end{remark}

\begin{proposition}\label{krazzzz}
Let $(X,\Ox)$ be a ringed topos on a site $C,$ $\calI$ a sheaf of monoids of $X,$ $\sigma:\calI^2 \ra \calI$ the addition morphism, $\calA$ an $\calI$-indexed algebra, $\calE$ an $\calI$-indexed $\calA$-module and $\F$ an $\Ox_{\calI}$-module. Then, there is a canonical $\calI$-indexed $\calA$-algebra structure on $\sigma_!\left (\calE \boxtimes \F\right ).$
\end{proposition}

\begin{proof}
For a morphism $V\ra U$ of $C$ and local sections $s,t\in \Gamma(U,\calI),$ we have a canonical morphism
$$
\calA_s(V) \otimes_{\Ox(V)} \bigoplus_{\substack{u,v\in \calI(V) \\ u+v=t}}\calE_u(V) \otimes_{\Ox(V)}\F_v(V) \ra \bigoplus_{\substack{\alpha,\beta \in \calI(V) \\ \alpha+\beta=s+t}} \calE_{\alpha}(V)\otimes_{\Ox(V)}\F_{\beta}(V),
$$
induced by the multiplications
$$
\calA_s\otimes_{\Ox_U}\calE_u \ra \calE_{s+u}.
$$
By \ref{indremark73} and \ref{indpropshriek} (3), these morphisms define the desired structure.
\end{proof}

\subsection*{Localization of indexed objects}

\begin{parag}\label{defresind}
Let $X$ be a topos, $\calI$ and $S$ sheaves of $X$ and $\F$ be a sheaf of $X_{/\calI}.$ Consider the commutative diagram of topoi
$$
\begin{tikzcd}
(X_{/S})_{/\calI_{|S}}\ar{r}{\sim}\ar[swap]{dr}{j_{\calI_{|S}}} & X_{/S\times\calI} \ar{r}{j_{S\times \calI/\calI}} \ar{d}{j_{S\times\calI/S}} & X_{/\calI} \ar{d}{j_{\calI}} \\
 & X_{/S} \ar{r}{j_S} & X
\end{tikzcd}
$$
We call \emph{localization of $\F$ at $S$} the sheaf $j_{S\times \calI/\calI}^*\F$ of $(X_{/S})_{\calI_{|S}}$ and we denote it by $\F_{|S}.$
We will also denote the composition $(X_{/S})_{\calI_{|S}} \xrightarrow{\sim} X_{/S\times\calI} \ra X_{/\calI}$ by $j_{S\times \calI/\calI}.$
\end{parag}

\begin{proposition}\label{Locprop1}
Let $(X,\Ox)$ be a ringed topos, $\calI$ a sheaf of monoids of $X,$ $S$ a sheaf of sets of $X,$ $\calA$ an $\calI$-indexed $\Ox$-algebra and $\calE$ an $\calI$-indexed $\calA$-module. Then $\calA_{|S}$ (resp. $\calE_{|S}$) is canonically equipped with an $\calI_{|S}$-indexed $\Ox_S$-algebra structure (resp. an $\calI_{|S}$-indexed $\calA_{|S}$-module structure).
\end{proposition}

\begin{proof}
Let $p_1,p_2:\calI^2\ra \calI$ and $p_{1,S},p_{2,S}:\calI_{|S}^2\ra \calI_{|S}$ be the canonical projections and $\sigma:\calI^2\ra \calI$ and $\sigma_S:\calI_{|S}^2\ra \calI_{|S}$ the addition maps. We have the following commutative diagrams of topoi
$$
\begin{tikzcd}
(X_{/S})_{/\calI^2_{|S}} \ar{r}{\sim} \ar[bend right=-30]{rr}{j_{S\times \calI^2/\calI^2}} \ar{d}{p_{i,S}} & X_{/S\times \calI^2} \ar{r} \ar{d} & X_{/\calI^2} \ar{d}{p_i} \\
(X_{/S})_{/\calI_{|S}} \ar{r}{\sim} \ar[bend right=30]{rr}{j_{S\times\calI/\calI}} & X_{S\times \calI} \ar{r} & X_{/\calI}
\end{tikzcd}
$$
$$
\begin{tikzcd}
(X_{/S})_{/\calI^2_{|S}} \ar{r}{\sim} \ar[bend right=-30]{rr}{j_{S\times \calI^2/\calI^2}} \ar{d}{\sigma_{S}} & X_{/S\times \calI^2} \ar{r} \ar{d} & X_{/\calI^2} \ar{d}{\sigma} \\
(X_{/S})_{/\calI_{|S}} \ar{r}{\sim} \ar[bend right=30]{rr}{j_{S\times\calI/\calI}} & X_{S\times \calI} \ar{r} & X_{/\calI}
\end{tikzcd}
$$
Applying $j_{S\times \calI^2/\calI^2}^*$ to the morphism
$$\calA\boxtimes \calA \ra \sigma^*\calA$$
defining the $\calI$-indexed algebra structure on $\calA,$ we get
$$\calA_{|S}\boxtimes \calA_{|S}\ra \sigma_{S}^*\calA_{|S}.$$
Similarly, applying $j_{S\times \calI^2/\calI^2}^*$ to the morphism
$$\calA\boxtimes \calE \ra \sigma^*\calE$$
defining the $\calI$-indexed $\calA$-module structure on $\calE,$ we get
$$\calA_{|S}\boxtimes \calE_{|S} \ra \sigma_S^*\calE_{|S}.$$
These two morphisms define the desired structures.
\end{proof}

\begin{parag}\label{inddef2}
Let $(X,\Ox)$ be a ringed topos on a $\mathbb{U}$-site $C,$ $\calI$ a sheaf of commutative monoids of $X,$ $\calA$ an $\calI$-indexed $\Ox$-algebra and $\calE$ and $\F$ two $\calI$-indexed $\calA$-modules. Denote by $p_1,p_2:\calI^2 \ra \calI$ the canonical projections and by $\sigma:\calI^2 \ra \calI$ the addition map. Equip $p_{1*}\mathscr{Hom}_{\Ox_{\calI^2}}(p_2^*\calE,\sigma^*\F)$ with the $\Ox_{\calI}$-module structure induced by the canonical morphism of rings $\Ox_{\calI} \ra p_{1*}\Ox_{\calI^2}.$
We denote by $\mathscr{Hom}_{\calA}(\calE,\F)$ the sub-$\Ox_{\calI}$-module of $p_{1*}\mathscr{Hom}_{\Ox_{\calI^2}}(p_2^*\calE,\sigma^*\F)$ defined as follows: let $s\in \Gamma(U,\calI)$ be a section of $\calI$ over an object $U$ of $C$ and the morphism $\sigma_s:\calI_{|U} \ra \calI_{|U},\ x\mapsto x+s.$ The $\Ox(U)$-module $p_{1*}\mathscr{Hom}_{\Ox_{\calI^2}}(p_2^*\calE,\sigma^*\F)(s)=\mathscr{Hom}_{\Ox_{\calI^2}}(p_2^*\calE,\sigma^*\F)(s\times \op{Id}_{\calI})$ is the set of $\Ox_U$-linear morphisms $f:(s\times \op{Id}_{\calI})^*p_2^*\calE \ra (s\times \op{Id}_{\calI})^*\sigma^*\F.$
Consider the commutative diagrams of topoi:
$$
\begin{tikzcd}
X_{/U\times \calI} \ar{r}{s\times \op{Id}_{\calI}}  \ar[swap]{dr}{j_{U\times \calI/\calI}} & X_{/\calI^2} \ar{d}{p_2} & & X_{U\times\calI} \ar{r}{s\times \op{Id}_{\calI}} \ar{d}{\sigma_s} & X_{/\calI^2}\ar{d}{\sigma} \\
 & X_{/\calI} & & X_{U\times \calI} \ar{r}{j_{U\times \calI/\calI}} & X_{/\calI}
\end{tikzcd}
$$ 
Define $\mathscr{Hom}_{\calA}(\calE,\F)(s)$ to be the subset of $\left (p_{1*}\mathscr{Hom}_{\Ox_{\calI^2}}(p_2^*\calE,\sigma^*\F)\right )(s)$ consisting of the morphisms of $\calI_{|U}$-indexed $\calA_{|U}$-modules $\calE_{|U} \ra \sigma_s^*\F_{|U}.$ This is clearly a subsheaf of $\Ox_{\calI}$-modules of $p_{1*}\mathscr{Hom}_{\Ox_{\calI^2}}(p_2^*\calE,\sigma^*\F).$

If $\calE=\F,$ we denote $\mathscr{Hom}_{\calA}(\calE,\F)$ simply by $\mathscr{End}_{\calA}(\calE).$
\end{parag}

\begin{parag}\label{indremtakriz1}
Keep the same hypothesis of \ref{inddef2}.
\begin{enumerate}
\item Let $U$ be an object of $C$ and $s\in \Gamma(U,\calI).$ Then $\mathscr{Hom}_{\calA}(\calE,\F)_s$ is the $\Ox_U$-module of $X_{/U}$ defined as follows: it sends an object $v:V\ra U$ of $C_{/U}$ to the set of morphisms of $\calI_{|V}$-indexed $\calA_{|V}$-modules $\calE_{|V} \ra \sigma_{s_{|V}}^*\F_{|V}.$

\item Let $U$ an object of $C,$ $s\in \Gamma(U,\calI)$ and $a\in \Gamma(U,\calA_s).$ The section $a$ can be considered as an element of $\mathscr{Hom}_{\calA}(\F,\F)(s)$ as follows: for any object $V\ra U$ of $C_{/U}$ and any section $t\in \Gamma(V,\calI_{|U})=\Gamma(V,\calI),$ consider the composition
$$
\begin{array}[t]{clclc}
\F_t & \ra & \calA_s\otimes_{\Ox_U}\F_t & \ra & \F_{t+s}\\
x& \mapsto & a\otimes x& \mapsto & ax.
\end{array}
$$
By \ref{era2indrem}, these compositions define a morphism
\begin{equation}\label{era2a}
a_{\F}:\F_{|U}\ra \sigma_s^*\F_{|U}.
\end{equation}

\item $\mathscr{Hom}_{\calA}(\calE,\F)$ has an $\calI$-indexed $\calA$-module structure defined as follows: Let $s,t\in \Gamma(U,\calI)$ be sections of $\calI$ over an object $U$ of $C.$ We will define a morphism
$$\calA_s\otimes_{\Ox_U}\mathscr{Hom}_{\calA}(\calE,\F)_t \ra \mathscr{Hom}_{\calA}(\calE,\F)_{s+t}.$$
Let $V \ra U$ be an object of $C_{/U},$ $a\in \Gamma(V,\calA_s)$ and $f\in \Gamma(V,\mathscr{Hom}_{\calA}(\calE,\F)_t).$ We set
$$af:\calE_{|V} \xrightarrow{f}\sigma_{t_{|V}}^*\F_{|V} \xrightarrow{\sigma_{t_{|V}}^*a_{\F}}\sigma_{(s+t)_{|V}}^*\F_{|V},$$
where $a$ is defined in (\ref{era2a}). Note that, by the $\calA_{|V}$-linearity of $f,$ the morphism $af$ is equal to the composition
$$\calE_{|V} \xrightarrow{a_{\calE}} \sigma_{s_{|V}}^* \calE_{|V} \xrightarrow{\sigma_{s_{|V}}^*f}\sigma^*_{(s+t)_{|V}}\F_{|V}.$$ 

\item The operation (\ref{era2compshift}) defines an $\calI$-indexed algebra structure on $\mathscr{End}_{\calA}(\calE).$
\end{enumerate}
\end{parag}

\begin{definition}\label{inddef33}
Let $(X,\Ox)$ be a ringed topos, $\calI$ a sheaf of commutative monoids of $X,$ $\calA$ an $\calI$-indexed algebra and $\calE$ and $\F$ two $\calI$-indexed $\calA$-modules. We call \emph{internal Hom functor of $\calI$-indexed $\calA$-modules} the sheaf $\mathscr{Hom}_{\calA}(\calE,\F)$ defined in \ref{inddef2}, which is naturally equipped with an $\calI$-indexed $\calA$-module structure in \ref{indremtakriz1} (3).
\end{definition}

\begin{parag}\label{algstrendleft}
The following is a generalization of \ref{indremtakriz1} (3). Let $(X,\Ox)$ be a ringed topos on a $\mathbb{U}$-site $C,$ $\calI$ a sheaf of monoids of $X,$ $p_1,p_2:\calI^2 \ra \calI$ the canonical projections, $\sigma:\calI^2\ra \calI$ the addition map, $\calA \ra \calB$ a morphism of $\calI$-indexed algebras and $\calE$ and $\F$ two $\calI$-indexed $\calB$-modules. Let $s,t\in \Gamma(U,\calI)$ be two sections of $\calI$ over an object $U$ of $C$ and $\sigma_s:\calI_{|U} \ra \calI_{|U},\ x\mapsto s+x.$ As in \ref{indremtakriz1}, a section $b\in \calB(s)$ defines morphisms
$$b_{\calE}:\calE_{|U} \ra \sigma_s^*\calE_{|U},$$
$$b_{\F}:\F_{|U} \ra \sigma_s^*\F_{|U}.$$
We have two morphisms
$$
\varphi^g(s,t):\begin{array}[t]{clc}
\calB(s) \otimes_{\Ox(U)}\mathscr{Hom}_{\calA}(\calE,\F)(t) & \ra & \mathscr{Hom}_{\calA}(\calE,\F)(s+t) \\
b \otimes f & \mapsto & f \circ b_{\calE}
\end{array}
$$
and
$$
\varphi^d(s,t):\begin{array}[t]{clc}
\calB(s) \otimes_{\Ox(U)}\mathscr{Hom}_{\calA}(\calE,\F)(t) & \ra & \mathscr{Hom}_{\calA}(\calE,\F)(s+t) \\
b \otimes f & \mapsto & b_{\F} \circ f,
\end{array}
$$
where the compositions are defined as follows:
\begin{alignat*}{2}
f\circ b_{\calE} &: \calE_{|U} \xrightarrow{b_{\calE}} \sigma_s^*\calE_{|U} \xrightarrow{\sigma_s^*f} \sigma_{s+t}^*\F_{|U}, \\
b_{\F}\circ f &: \calE_{|U} \xrightarrow{f} \sigma_t^*\F_{|U} \xrightarrow{\sigma_t^*b_{\F}} \sigma_{s+t}^*\F_{|U}.
\end{alignat*}
These morphisms, for every $(s,t),$ define two morphisms
$$
\varphi^g,\varphi^d:\calB \boxtimes \mathscr{Hom}_{\calA}(\calE,\F) \ra \sigma^*\mathscr{Hom}_{\calA}(\calE,\F).
$$
These define two structures of $\calI$-indexed $\calB$-modules on $\mathscr{Hom}_{\calA}(\calE,\F).$ We call them \emph{the left} (resp. \emph{right}) \emph{$\calI$-indexed $\calB$-module structure on $\mathscr{Hom}_{\calA}(\calE,\F).$}
\end{parag}

\begin{parag}\label{indtensprod}
Let $(X,\Ox)$ be a ringed topos over a $\mathbb{U}$-site $C,$ $\calI$ a sheaf of commutative monoids of $X,$ $\calA$ a commutative $\calI$-indexed algebra, $\calE$ and $\F$ two $\calI$-indexed $\calA$-modules (\ref{inddef}), $p_1,p_2:\calI^2\ra \calI,$ $q_1,q_2,q_3:\calI^3\ra \calI$ and $p_{ij}:\calI^3\ra \calI^2,\ 1\le i<j\le 3$ the canonical projections and $\sigma:\calI^2\ra \calI$ and $\tau:\calI^3\ra \calI$ the addition maps. The $\calA$-module structures on $\calE$ and $\F$ induce morphisms
$$\calA\boxtimes \calE \boxtimes \F \ra (\sigma^*\calE) \boxtimes \F$$
and
$$\calA\boxtimes \calE \boxtimes \F \ra \calE \boxtimes (\sigma^*\F).$$
The commutativity of the diagrams
$$
\begin{tikzcd}
\calI^3 \ar{r}{p_{12}} \ar[swap]{d}{\sigma \times \op{Id}_{\calI}} & \calI^2\ar{d}{\sigma} & & \calI^3\ar{dr}{q_3}\ar[swap]{d}{\sigma\times \op{Id}_{\calI}} & \\
\calI^2\ar{r}{p_1} & \calI & & \calI^2\ar{r}{p_2} & \calI
\end{tikzcd}
$$
implies that
$$(\sigma^*\calE) \boxtimes \F=(\sigma \times \op{Id}_{\calI})^*(\calE \boxtimes \F).$$
Similarly, we have
$$\calE \boxtimes (\sigma^*\F)=(\op{Id}_{\calI}\times\sigma)^*(\calE \boxtimes \F).$$
Thus we have morphisms
\begin{alignat}{2}
\calA\boxtimes \calE \boxtimes \F \ra (\sigma \times \op{Id}_{\calI})^*(\calE \boxtimes \F), \label{indeqqq1} \\
\calA\boxtimes \calE \boxtimes \F \ra (\op{Id}_{\calI}\times\sigma)^*(\calE \boxtimes \F). \label{indeqqq2}
\end{alignat}
Since $\tau=\sigma \circ (\sigma \times \op{Id}_{\calI})=\sigma \circ (\op{Id}_{\calI}\times \sigma),$ we get
$$\tau_!=\sigma_! \circ (\sigma \times \op{Id}_{\calI})_!=\sigma_! \circ (\op{Id}_{\calI}\times \sigma)_!.$$
By adjunction and applying $\sigma_!$ on both sides of (\ref{indeqqq1}) and (\ref{indeqqq2}), we get two morphisms
\begin{equation}\label{indeq10}
\tau_!(\calA \boxtimes \calE \boxtimes \F) \rightrightarrows \sigma_!(\calE \boxtimes \F).
\end{equation}
We denote by $\calE \circledast_{\calA} \F$ the coequalizer of these two morphisms. We then have an exact sequence of $\Ox_{\calI}$-modules
\begin{equation}\label{exseqindtens}
\tau_!(\calA \boxtimes \calE \boxtimes \F) \ra \sigma_!(\calE \boxtimes \F) \ra \calE \circledast_{\calA}\F \ra 0.
\end{equation}
In the next proposition, we prove that the $\calI$-indexed $\calA$-module structures on $\calE$ and $\F$ induce an $\calI$-indexed $\calA$-module structure on $\calE\circledast_{\calA}\F.$
\end{parag}

\begin{proposition}
Keep the same hypothesis of \ref{indtensprod}. Then the $\calI$-indexed $\calA$-module structures on $\calE$ and $\F$ induce an $\calI$-indexed $\calA$-module structure on $\calE\circledast_{\calA}\F.$
\end{proposition}

\begin{proof}
By \ref{indpropshriek} and \ref{indpropshriek2}, the two morphisms (\ref{indeq10}) are given explicitly as follows: for a section $l\in \calI(U)$ over an object $U$ of $C,$ the morphisms
$$\tau_!(\calA \boxtimes \calE \boxtimes \F)_l \rightrightarrows \sigma_!(\calE \boxtimes \F)_l$$
are induced by
\begin{equation}\label{indeq14}
\begin{array}[t]{clc}
\bigoplus\limits_{\substack{i,j,k\in \calI(V) \\ i+j+k=l}}\calA_i(V)\otimes_{\Ox(V)}\calE_i(V)\otimes_{\Ox(V)}\F_k(V) & \ra & \bigoplus\limits_{\substack{r,s\in \calI(V) \\ r+s=l}}\calE_r(V)\otimes_{\Ox(V)} \F_s(V) \\
a \otimes x \otimes y & \mapsto & (ax) \otimes y \\
a \otimes x \otimes y & \mapsto & x \otimes (ay)
\end{array}
\end{equation}
where $V\ra U$ is an object of $C_{/U},$ $a \in \calA_i(V),$ $x\in \calE_j(V)$ and $y\in \F_k(V).$

The adjunction morphism
$$\calE \boxtimes \F \ra \sigma^*\sigma_!(\calE\boxtimes \F)$$
induces a morphism
$$p_{23}^*(\calE\boxtimes \F) \ra  p_{23}^*\sigma ^* \sigma_!(\calE\boxtimes \F)=(\op{Id}_{\calI}\times \sigma)^*p_2^*\sigma_!(\calE\boxtimes \F).$$
Again, by adjunction we obtain a morphism
\begin{equation}\label{indeq9}
(\op{Id}_{\calI}\times \sigma)_!p_{23}^*(\calE\boxtimes \F) \ra p_2^*\sigma_!(\calE\boxtimes \F).
\end{equation}
By \ref{indpropshriekoh}, it is an isomorphism.
Its inverse
$$p_2^*\sigma_!(\calE\boxtimes \F) \ra (\op{Id}_{\calI}\times \sigma)_!p_{23}^*(\calE\boxtimes \F)$$
induces a morphism
$$\calA \boxtimes \sigma_!(\calE \boxtimes \F) \ra p_1^*\calA \otimes_{\Ox_{\calI^2}} (\op{Id}_{\calI}\times \sigma)_!p_{23}^*(\calE \boxtimes \F).$$
By (\cite{SGA43} Exposé IV 12.11.b), there exists a canonical isomorphism
\begin{alignat*}{2}
p_1^*\calA \otimes_{\Ox_{\calI^2}} (\op{Id}_{\calI}\times \sigma)_!p_{23}^*(\calE \boxtimes \F) &\xrightarrow{\sim}  (\op{Id}_{\calI}\times \sigma)_!\left ( (\op{Id}_{\calI}\times \sigma)^*p_1^*\calA \otimes_{\Ox_{\calI^3}}p_{23}^*(\calE \boxtimes \F)\right ) \\
&\xrightarrow{\sim} (\op{Id}_{\calI}\times \sigma)_!(\calA \boxtimes \calE \boxtimes \F).
\end{alignat*}
It follows that we have a morphism
$$\calA \boxtimes \sigma_!(\calE\boxtimes \F) \ra (\op{Id}_{\calI}\times \sigma)_!(\calA\boxtimes \calE\boxtimes \F)$$
and then a morphism
$$\sigma_!(\calA \boxtimes \sigma_!(\calE\boxtimes \F)) \ra \tau_!(\calA \boxtimes \calE\boxtimes \F).$$
Composing with one of the morphisms (\ref{indeq10}) and then the canonical projection, we get a morphism
$$\sigma_!(\calA \boxtimes \sigma_!(\calE\boxtimes \F)) \ra \calE \circledast_{\calA}\F$$
which induces by adjunction
\begin{equation}\label{indeq13}
\calA \boxtimes \sigma_!(\calE \boxtimes \F)\ra \sigma^*(\calE \circledast_{\calA}\F).
\end{equation}
For sections $k,l\in \calI(U)$ over an object $U$ of $C,$ the fiber of (\ref{indeq13}) at $(k,l)$ is induced by
$$
\begin{array}[t]{clc}
\calA_k(U) \otimes_{\Ox(U)} \bigoplus\limits_{i+j=l}\calE_i(U) \otimes_{\Ox(U)}\F_j(U) & \ra & \bigoplus\limits_{r+s=k+l}\calE_r(U)\otimes_{\Ox(U)}\F_s(U) / \sim \\
a \otimes x\otimes y & \mapsto & (ax)\otimes y=x\otimes (ay)
\end{array}
$$
where $a\in \calA_k(U),$ $x\in \calE_i(U)$ and $y\in \F_j(U)$ and $\sim$ is the equalizer of the morphisms (\ref{indeq14}).
This morphism factors through the quotient of
$$\calA_k(U) \otimes_{\Ox(U)} \bigoplus\limits_{i+j=l}\calE_i(U) \otimes_{\Ox(U)}\F_j(U)$$
by the submodule generated by sections
$$a \otimes ((bx)\otimes y-x\otimes(by)),$$
where $a\in \calA_k(U),$ $b\in \calA_{\alpha}(U),$ $x\in \calE_{\beta}(U),$ $y\in \F_{\gamma}(U)$ and $\alpha,\beta,\gamma \in \calI(U)$ such that $\alpha+\beta+\gamma=i+j.$
It follows that the morphism (\ref{indeq13}) factors through
$$\calA \boxtimes (\calE\circledast_{\calA}\F).$$
The resulting morphism satisfies the conditions of (\ref{inddef} 3) and thus gives the desired structure.
The $\calI$-indexed $\calA$-module $\calE\circledast_{\calA}\F$ is called \emph{the tensor product of $\calE$ and $\F$ over $\calA.$} 
\end{proof}

\begin{proposition}\label{indproptensprod}
Keep the same hypothesis of \ref{indtensprod} and consider a section $k\in \Gamma(U,\calI)$ of $\calI$ over an object $U$ of $C.$ By construction, $(\calE \circledast_{\calA}\F)_k$ is the sheaf associated to the the presheaf
$$
\begin{array}[t]{clc}
C_{/U} & \mapsto & \op{Sets} \\
V & \mapsto & \bigoplus\limits_{\substack{i,j\in \calI(V) \\ i+j=k}}\calE_i(V)\otimes_{\Ox(V)}\F_j(V)/\sim,
\end{array}
$$
where $(ax)\otimes y\sim x\otimes (ay)$ for every section $a\in\calA_r(V),$ $x\in \calE_s(V),$ $y\in \F_t(V)$ and $r,s,t\in \calI(V)$ such that $r+s+t=k.$
\end{proposition}

\begin{remark}\label{indrem1}
Keep the same hypothesis of \ref{indtensprod}. 
\begin{enumerate}
\item Let $f:\calA \ra \calA'$ be a morphism of $\calI$-indexed algebras. The composition
$$\calA\boxtimes \calA' \xrightarrow{f\boxtimes \op{Id}_{\calA'}} \calA' \boxtimes \calA' \xrightarrow{\pi_{\calA'}} \sigma^*\calA'$$
defines an $\calI$-indexed $\calA$-module structure on $\calA'.$
The morphism 
$$\pi_{\calA'}:\calA' \boxtimes \calA' \ra \sigma^*\calA'$$
defining the indexed algebra structure on $\calA,$ induces a morphism
$$\sigma_!(\calA' \boxtimes \calA') \ra \calA'$$
and then an $\calA$-linear morphism
\begin{equation}\label{indrem11}
v:\calA' \circledast_{\calA}\calA' \ra \calA'.
\end{equation}
We call $v$ the \emph{$\calI$-indexed codiagonal map.}
\item Let $\calA \ra \calA'$ be a morphism of $\calI$-indexed algebras. Denote by $e$ the final object of $X,$ $0:e\ra \calI$ the zero section of $\calI$ and $\nu:\calI \ra e$ the canonical morphism. The morphisms
$$\calA' \xrightarrow{\op{Id}_{\calA'}\times 1_{\calA'}} \calA' \otimes_{\Ox_{\calI}} \nu^*\calA'_0=(\op{Id}_{\calI}\times 0)^*(\calA' \boxtimes \calA')$$
$$\calA' \xrightarrow{1_{\calA'}\times \op{Id}_{\calA'}} \nu^*\calA'_0\otimes_{\Ox_{\calI}} \calA'=(0\times \op{Id}_{\calI})^*(\calA' \boxtimes \calA')$$
induce two morphisms
$$(\op{Id}_{\calI}\times 0)_! \calA' \ra \calA' \boxtimes \calA'$$
$$(0\times \op{Id}_{\calI})_! \calA' \ra \calA' \boxtimes \calA'$$
and then two morphisms
\begin{equation}\label{indeq17}
\calA'=\sigma_!(\op{Id}_{\calI}\times 0)_!\calA' \rightrightarrows \sigma_!(\calA' \boxtimes \calA').
\end{equation}
Composing with the canonical projection
$$\sigma_!(\calA' \boxtimes \calA') \ra \calA' \circledast_{\calA}\calA',$$
we get two $\calA$-linear morphisms
\begin{equation}\label{indrem12}
\pi_1,\pi_2:\calA' \ra \calA' \circledast_{\calA} \calA'.
\end{equation}
These morphisms are induced by the morphisms of presheaves
$$
\begin{array}[t]{clclc}
\calA'_k(U) & \ra & \calA'_k(U)\otimes_{\Ox(U)}\calA'_0(U) & \subset & \bigoplus\limits_{\substack{i,j\in \calI(U) \\ i+j=k}}\calA'_i(U) \otimes_{\Ox(U)}\calA'_j(U) \\
x & \mapsto & x \otimes 1 & &
\end{array}
$$
and
$$
\begin{array}[t]{clclc}
\calA'_k(U) & \ra & \calA'_0(U)\otimes_{\Ox(U)}\calA'_k(U) & \subset & \bigoplus\limits_{\substack{i,j\in \calI(U) \\ i+j=k}}\calA'_i(U) \otimes_{\Ox(U)}\calA'_j(U) \\
x & \mapsto & 1 \otimes x & &
\end{array}
$$
\end{enumerate}
\end{remark}

\begin{lemma}\label{lemtop15}
Let $u:(T',\Ox_{T'})\ra (T,\Ox_T)$ be a morphism of ringed topoi, $I$ a sheaf of $T,$ $I'=u^{-1}I$ and $j_I$ and $j_{I'}$ the localization morphisms. Equip $T_{/I}$ and $T'_{/I'}$ with the rings $j_I^{-1}\Ox_T$ and $j_{I'}^{-1}\Ox_{T'}$ respectively. Consider the commutative diagram of ringed topoi
$$
\begin{tikzcd}
T'_{/I'} \ar{r}{u'} \ar[swap]{d}{j_{I'}} & T_{/I} \ar{d}{j_I} \\
T' \ar{r}{u} & T,
\end{tikzcd}
$$
where $u'$ is induced by $u.$
\begin{enumerate}
\item For every sheaf $\F'$ of $T',$ there exists a canonical isomorphism
$$j_I^*u_*\F' \xrightarrow{\sim} u'_*j_{I'}^*\F'.$$
\item For every sheaf $\F$ of $T_{/I},$ there exists a canonical isomorphism
$$j_{I'!}u'^*\F \xrightarrow{\sim} u^*j_{I!}\F.$$
\end{enumerate}
The proposition is true if we consider topoi instead of ringed topoi and we replace $u^*,$ $u'^*,$ $j_{\calI}^*$ and $j_{\calI'}^*$ by $u^{-1},$ $u'^{-1},$ $j_{\calI}^{-1}$ and $j_{\calI'}^{-1}$ respectively.
\end{lemma}

\begin{proof}
The adjunction morphism
$$\F' \ra j_{I'*}j_{I'}^*\F'$$
induces
$$u_*\F' \ra u_*j_{I'*}j_{I'}^*\F'=j_{I*}u'_*j_{I'}^*\F'.$$
We then get a canonical morphism
\begin{equation}\label{eqtop11}
j_I^*u_*\F' \ra u'_*j_{I'}^*\F'.
\end{equation}
This morphism is an isomorphism. Indeed, if $U$ is an object of $T$ then
$$
\left (j_I^*u_*\F' \right )\left (U_{/I} \right )=\left (u_*\F'\right )(U)=\F'(u^{-1}U),
$$
$$
\left (u'_*j_{I'}^*\F' \right )\left (U_{/I} \right )=\left (j_{I'}^*\F' \right )\left ( (u^{-1}U)_{/I'}\right )=\F'(u^{-1}U).
$$
The second assertion follows from the first as follows: for any object $\calG$ of $T',$ we have
\begin{alignat*}{2}
\op{Hom}_{T'}(j_{I'!}u'^*\F,\calG) &= \op{Hom}_{T_{/I}}(\F,u'_*j_{I'}^*\calG )\\
&= \op{Hom}_{T_{/I}}(\F,j_I^*u_*\calG) \\
&= \op{Hom}_{T'}(u^*j_{I!}\F,\calG).
\end{alignat*}
The proof for topoi instead of ringed topoi is similar.
\end{proof}

\begin{proposition}\label{ren3}
Let $(X,\Ox)$ be a ringed topos, $\calI$ a sheaf of monoids of $X,$ $S$ a sheaf of $X$, $\calA$ an $\calI$-indexed algebra and $\calE$ and $\F$ two $\calI$-indexed $\calA$-modules. There exists a canonical isomorphism
$$
\left (\calE_{|S} \right ) \circledast_{\calA_{|S}} \left (\F_{|S} \right ) \xrightarrow{\sim} \left (\calE \circledast_{\calA}\F \right )_{|S}  .
$$
\end{proposition}

\begin{proof}
Denote by $\sigma:\calI^2 \ra \calI$ and $\tau:\calI^3 \ra \calI$ the addition maps and $p_1,p_2:\calI^2 \ra \calI$ the canonical projections. We have the commutative diagram of topoi
$$
\begin{tikzcd}
\left (X_{/S} \right )_{/\calI_{|S}^2} \ar{r}{j_{\calI^2\times S/\calI^2}} \ar[swap]{d}{\sigma_{|S}} & X_{/\calI^2} \ar{d}{\sigma} \\
\left (X_{/S} \right )_{/\calI_{|S}} \ar{r}{j_{\calI\times S/\calI}} & X_{/\calI}.
\end{tikzcd}
$$
By \ref{lemtop15}, there exists a canonical isomorphism
$$
\left (\sigma_{|S} \right )_! \left ( \left ( \calE \boxtimes \F \right )_{|S} \right ) \xrightarrow{\sim} \left ( \sigma_!(\calE \boxtimes \F) \right )_{|S}.
$$
By the commutative diagram
$$
\begin{tikzcd}
\left (X_{/S} \right )_{/\calI_{|S}^2} \ar{r}{j_{\calI^2\times S/\calI^2}} \ar[swap]{d}{p_{i|S}} & X_{/\calI^2} \ar{d}{p_i} \\
\left (X_{/S} \right )_{/\calI_{|S}} \ar{r}{j_{\calI\times S/\calI}} & X_{/\calI},
\end{tikzcd}
$$
we get a canonical isomorphism
$$
\left (\calE_{|S} \right ) \boxtimes \left (\F_{|S} \right ) \xrightarrow{\sim} \left ( \calE \boxtimes \F \right )_{|S}.
$$
It follows that we have a canonical isomorphism
$$
\left (\sigma_{|S} \right )_! \left (\calE_{|S}  \boxtimes \F_{|S}  \right ) \xrightarrow{\sim} \left ( \sigma_!(\calE \boxtimes \F) \right )_{|S}.
$$
Similarly, we have a canonical isomorphism
$$
\left (\tau_{|S} \right )_! \left ( \calA_{|S} \boxtimes \calE_{|S} \boxtimes \F_{|S} \right ) \xrightarrow{\sim} \left ( \tau_!(\calA \boxtimes \calE \boxtimes \F) \right )_{|S}.
$$
By definition of $\calE \circledast_{\calA} \F$ and the right exactness of $\left (\tau_{|S} \right )_!$ and the restriction functor $\calM \mapsto \calM_{|S,}$ these isomorphisms fit into a commutative diagram
$$
\begin{tikzcd}
\left (\tau_{|S} \right )_! \left ( \calA_{|S} \boxtimes \calE_{|S} \boxtimes \F_{|S} \right ) \ar[sloped]{d}{\sim} \ar{r} & \left (\sigma_{|S} \right )_! \left ( \calE_{|S} \boxtimes \F_{|S} \right ) \ar{r} \ar[sloped]{d}{\sim} & \left (\calE_{|S} \right ) \circledast_{\calA_{|S}} \left (\F_{|S} \right ) \ar{r}  \ar[sloped,dashed]{d}{\sim} & 0 \\
\left ( \tau_!(\calA \boxtimes \calE \boxtimes \F) \right )_{|S} \ar{r} & \left ( \sigma_!(\calE \boxtimes \F) \right )_{|S} \ar{r} & \left (\calE \circledast_{\calA} \F\right )_{|S} \ar{r} & 0,
\end{tikzcd}
$$
where the rows are induced by \eqref{exseqindtens} and hence exact.
Hence the existence of the third vertical isomorphism.
\end{proof}

\begin{proposition}
Let $C$ be a $\mathbb{U}$-site, $(X,\Ox)$ a ringed topos on $C,$ $\calI$ a sheaf of monoids of $X,$ $\calB \ra \calA$ a morphism of $\calI$-indexed algebras and $\calE,$ $\F$ and $\calG$ three $\calI$-indexed $\calA$-modules. Denote by $p_1,p_2:\calI^2 \ra \calI$ the canonical projections and by $\sigma:\calI^2 \ra \calI$ the addition map. Consider the morphism
\begin{equation}\label{takcompadj}
\mathscr{Hom}_{\calA}(\calE,\F) \boxtimes \mathscr{Hom}_{\calA}(\F,\calG) \ra \sigma^*\mathscr{Hom}_{\calA}(\calE,\calG)\end{equation}
defined, for every sections $s,t\in \Gamma(U,\calI)$ of $\calI$ over an object $U$ of $C,$ by
\begin{equation}
\begin{array}[t]{clc}
\mathscr{Hom}_{\calA}(\calE,\F)(s) \otimes_{\Ox(U)}\mathscr{Hom}_{\calA}(\F,\calG)(t) & \ra & \mathscr{Hom}_{\calA}(\calE,\calG)(s+t) \\
f\otimes g & \mapsto & g\circ f,
\end{array}
\end{equation}
where $g\circ f$ is the following composition
\begin{equation}\label{era2compshift}
g\circ f:\calE_{|U} \xrightarrow{f} \sigma_s^*\F_{|U} \xrightarrow{\sigma_s^*g} \sigma_{s+t}^*\calG_{|U}.
\end{equation}
The morphism \eqref{takcompadj} induces a morphism
\begin{equation}\label{era3comptakriz}
\mathscr{Hom}_{\calA}(\calE,\F) \circledast_{\calB} \mathscr{Hom}_{\calA}(\F,\calG) \ra \mathscr{Hom}_{\calA}(\calE,\calG).\end{equation}
\end{proposition}

\begin{proof}
The morphism \eqref{takcompadj} induces, by adjunction,
\begin{equation}\label{eqa}
\sigma_!\left (\mathscr{Hom}_{\calA}(\calE,\F) \boxtimes \mathscr{Hom}_{\calA}(\F,\calG) \right ) \ra \mathscr{Hom}_{\calA}(\calE,\calG).
\end{equation}
By \ref{indpropshriek} (3), $\sigma_!\left (\mathscr{Hom}_{\calA}(\calE,\F) \boxtimes \mathscr{Hom}_{\calA}(\F,\calG) \right )$ is the sheaf associated to the presheaf
$$(s:U\ra \calI) \mapsto \bigoplus_{\substack{u,v\in \Gamma(U,\calI) \\ u+v=s}}\mathscr{Hom}_{\calA}(\calE,\F)(u) \otimes_{\Ox(U)} \mathscr{Hom}_{\calA}(\F,\calG)(v).$$
The morphism \eqref{eqa} is then induced, for every sections $u,v,s\in \Gamma(U,\calI)$ such that $u+v=s,$ by
$$
\varphi_{u,v}:\begin{array}[t]{clc}
\mathscr{Hom}_{\calA}(\calE,\F)(u) \otimes_{\Ox(U)} \mathscr{Hom}_{\calA}(\F,\calG)(v) & \ra & \mathscr{Hom}_{\calA}(\calE,\calG)(s) \\
f \otimes g & \mapsto & g\circ f.
\end{array}
$$
By \ref{indremtakriz1} (3), for any sections $u,v,w\in \Gamma(U,\calI),$ $f\in \mathscr{Hom}_{\calA}(\calE,\F)(u),$ $g\in\mathscr{Hom}_{\calA}(\F,\calG)(v)$ and $a\in \calB(w),$ we have
$$\varphi_{u+w,v}((af)\otimes g)=\varphi_{u,v+w}(f\otimes (ag)).$$
The result then follows by \ref{indproptensprod}.
\end{proof}

\begin{proposition}\label{indprop710}
Let $(X,\Ox)$ be a ringed topos and $\calI$ a sheaf of monoids of $X.$
\begin{enumerate}
\item For every $\calI$-indexed algebra $\calA$ and $\calI$-indexed $\calA$-module $\calE,$ there exists a canonical isomorphism
$$\calA\circledast_{\calA}\calE \xrightarrow{\sim} \calE.$$
\item For every $\Ox$-algebra $\calA$ over $X$ and $\calI$-indexed algebra $\calB,$ the $\Ox_{\calI}$-module $j_{\calI}^*\calA\otimes_{\Ox_{\calI}}\calB$ is naturally an $\calI$-indexed algebra and the canonical morphism
$$\calB \ra j_{\calI}^*\calA \otimes_{\Ox_{\calI}}\calB$$
is a morphism of $\calI$-indexed algebras.
\item For every $\calI$-indexed algebra $\calA,$ $\Ox$-module $\calE$ and $\calI$-indexed $\calA$-modules $\F$ and $\calG,$ the $\Ox_{\calI}$-module $(j_{\calI}^*\calE)\otimes_{\Ox_{\calI}}\F$ is naturally an $\calI$-indexed $\calA$-module and there exists a canonical isomorphism
\begin{equation}\label{indiso1}
\left ( j_{\calI}^*\calE \otimes _{\Ox_{\calI}}\F \right ) \circledast_{\calA} \calG \xrightarrow{\sim} j_{\calI}^*\calE \otimes_{\Ox_{\calI}} \left ( \F \circledast_{\calA}\calG \right ).
\end{equation}
\end{enumerate}
\end{proposition}

\begin{proof}
For the first assertion, the morphism
$$\calA \boxtimes \calE \ra \sigma^*\calE$$
defining the indexed $\calA$-module structure on $\calE,$ induces by adjunction
\begin{equation}\label{indeq16}
\sigma_!(\calA \boxtimes \calE) \ra \calE.
\end{equation}
By the description given in \ref{indproptensprod}, (\ref{indeq16}) induces a morphism
\begin{equation}\label{indeq20}
\calA \circledast_{\calA}\calE \rightarrow \calE.
\end{equation}
Let $k\in \calI(U)$ be a section over an object $U$ of $C.$ The morphism
$$\left (\calA \circledast_{\calA}\calE \right )_k\ra \calE_k,$$
induced by (\ref{indeq20}), is defined by the morphism of presheaves $\varphi$ given, for every object $V$ of $C_{/U},$ by
$$
\varphi(V):\begin{array}[t]{clc}
\bigoplus\limits_{\substack{i,j\in \calI(V)\\ i+j=k}}\calA_i(V)\otimes \calE_j(V) /\sim& \ra & \calE_k(V) \\
a\otimes x & \mapsto & ax
\end{array}$$
for $a\in \calA_i(V),\ x\in \calE_j(V)$ and where $\sim$ is the relation defined in \ref{indproptensprod}.
The morphism $\varphi$ is clearly surjective. Now let $n$ be a positive integer, $i_1,\hdots,i_n,j_1,\hdots,j_n \in \calI(V)$ such that $i_1+j_1=i_2+j_2=\hdots=i_n+j_n=k,$ $a_{l}\in \calA_{i_l}(V),$ $x_{l}\in \calE_{j_l}(V)$ for all $1\le l\le n,$ such that 
$$\varphi(a_1\otimes x_1+\hdots+a_n\otimes x_n)=0.$$
Then
$$\sum_{l=1}^na_{l}x_{l}=0.$$
But
$$\sum_{l=1}^na_{l}\otimes x_{l} \sim 1\otimes \left (\sum_{l=1}^n a_{l}x_{l} \right ).$$
This proves that $\varphi$ is injective.

For the second assertion, the indexed algebra structure on $j_{\calI}^*\calA\otimes_{\Ox_{\calI}}\calB$ is given by the composition
\begin{alignat*}{2}
(j_{\calI}^*\calA\otimes_{\Ox_{\calI}}\calB)\boxtimes (j_{\calI}^*\calA\otimes_{\Ox_{\calI}}\calB) & \xrightarrow{\sim} j_{\calI^2}^*(\calA\otimes_{\Ox}\calA)\otimes_{\Ox_{\calI^2}}(\calB\boxtimes \calB) \\
& \xrightarrow{j_{\calI^2}^*\pi_{\calA}\otimes \pi_{\calB}} (j_{\calI^2}^*\calA)\otimes_{\Ox_{\calI^2}}\sigma^*\calB\\
& \xrightarrow{\sim} \sigma^*(j_{\calI}^*\calA\otimes_{\Ox_{\calI}}\calB),
\end{alignat*}
where the first and third arrows are the canonical isomorphisms and $\pi_{\calA}$ and $\pi_{\calB}$ are the multiplication maps defining the algebra structures on $\calA$ and $\calB$ respectively. The fact that the canonical morphism
$$\calB \ra j_{\calI}^*\calA \otimes_{\Ox_{\calI}}\calB$$
is a morphism of $\calI$-indexed algebras follows from the obvious commutativity of the following diagram, for any local sections $i$ and $j$ of $\calI$ over an object $U,$
$$
\begin{tikzcd}
\calB_i \otimes_{\Ox_U}\calB_j \ar{r}\ar{d} & (\calA_{|U}\otimes_{\Ox_U}\calB_i)\otimes_{\Ox_U}(\calA_{|U}\otimes_{\Ox_U}\calB_j)\ar{d} \\
\calB_{i+j}\ar{r} & \calA_{|U}\otimes_{\Ox_U}\calB_{i+j}
\end{tikzcd}
$$
For the third assertion, the $\calI$-indexed $\calA$-module structure on $j_{\calI}^*\calE \otimes_{\Ox_{\calI}}\F$ is given by the composition
\begin{alignat*}{2}
\calA \boxtimes (j^*_{\calI}\calE \otimes_{\Ox_{\calI}}\F) &\xrightarrow{\sim} (\calA\boxtimes \F)\otimes_{\Ox_{\calI^2}}j_{\calI^2}^*\calE \\
& \xrightarrow{\pi \otimes \op{Id}_{j_{\calI^2}^*\calE}} (\sigma^*\F) \otimes_{\Ox_{\calI^2}}j_{\calI^2}^*\calE \\
& \xrightarrow{\sim} \sigma^*\left (j_{\calI}^*\calE \otimes_{\Ox_{\calI}}\F \right ),
\end{alignat*}
where the first and third arrows are the canonical isomorphisms and $\pi$ is the multiplication map defining the indexed module structure on $\F.$
The canonical isomorphism
$$(j_{\calI}^*\calE \otimes_{\Ox_{\calI}}\F)\boxtimes \calG \xrightarrow{\sim} j_{\calI^2}^*\calE \otimes_{\Ox_{\calI^2}}(\F \boxtimes \calG)$$
induces an isomorphism
$$\sigma_!\left ((j_{\calI}^*\calE \otimes_{\Ox_{\calI}}\F)\boxtimes \calG \right )\xrightarrow{\sim} \sigma_!\left (j_{\calI^2}^*\calE \otimes_{\Ox_{\calI^2}}(\F \boxtimes \calG)\right ).$$
Composing it with the canonical isomorphism (\cite{SGA43} Exposé IV 12.11.b)
$$\sigma_!\left (j_{\calI^2}^*\calE \otimes_{\Ox_{\calI^2}}(\F \boxtimes \calG)\right ) \xrightarrow{\sim} j_{\calI}^*\calE \otimes_{\Ox_{\calI}} \sigma_!(\F \boxtimes \calG),$$
we get a canonical isomorphism
$$\sigma_!\left ((j_{\calI}^*\calE \otimes_{\Ox_{\calI}}\F)\boxtimes \calG \right )\xrightarrow{\sim}j_{\calI}^*\calE \otimes_{\Ox_{\calI}} \sigma_!(\F \boxtimes \calG)$$
which induces the desired isomorphism.
\end{proof}

\begin{proposition}\label{propisoAkaza}
Let $(X,\Ox)$ be a ringed topos, $\calI$ a sheaf of monoids of $X,$ $\calA$ an $\Ox$-algebra, $\calA'$ an $\calA$-algebra, $\calB$ an $\calI$-indexed $\Ox$-algebra and $\calE$ an $\calI$-indexed $(\calA_{\calI}\otimes_{\Ox_{\calI}}\calB)$-module. Then $\calE$ is canonically equipped with an $\calA_{\calI}$-module structure and there exists a canonical isomorphism
\begin{equation}\label{propisoAkazahh}
\calA'_{\calI}\otimes_{\calA_{\calI}}\calE \xrightarrow{\sim} (\calA'_{\calI}\otimes_{\Ox_{\calI}}\calB)\circledast_{\calA_{\calI}\otimes_{\Ox_{\calI}}\calB}\calE.
\end{equation}
\end{proposition}

\begin{proof}
Let $p_1,p_2:\calI^2\ra \calI$ be the canonical projections, $\sigma:\calI^2\ra \calI$ the addition map and $i_2:\calI \ra \calI^2,\ s\mapsto (0,s).$ The canonical map
$$(\calA_{\calI}\otimes_{\Ox_{\calI}}\calB)\boxtimes \calE \ra \sigma^*\calE$$
induces
\begin{equation}\label{eq7241}
i_2^*p_1^*(\calA_{\calI}\otimes_{\Ox_{\calI}}\calB)\otimes_{\Ox_{\calI}}i_2^*p_2^*\calE \ra i_2^*\sigma^*\calE.
\end{equation}
Let $e$ be the final object of $X.$ By the commutativity of the diagrams
$$
\begin{tikzcd}
\calI \ar{r}{i_2} \ar{d} & \calI^2\ar{d}{p_1} & & \calI \ar{r}{i_2} \ar[swap]{dr}{\op{Id}_{\calI}} & \calI^2 \ar{d}{p_2} & & \calI \ar{r}{i_2} \ar[swap]{dr}{\op{Id}_{\calI}} & \calI^2 \ar{d}{\sigma} \\
e \ar{r}{0} & \calI & & & \calI & & & \calI
\end{tikzcd}
$$
the equality \eqref{eq7241} is equivalent to
$$(\calA_{\calI} \otimes_{\Ox_{\calI}}\calB_{0,\calI})\otimes_{\Ox_{\calI}}\calE \ra \calE.$$
Composing this morphism with the morphism induced by
$$\calA_{\calI} \ra \calA_{\calI}\otimes_{\Ox_{\calI}} \calB_{0,\calI},\ x\mapsto x\otimes 1,$$
we get the $\calA_{\calI}$-module structure on $\calE.$

For any section $s\in  \Gamma(U,\calI)$ over an object $U,$ consider the composition of
$$\begin{array}[t]{clc}
\calA'_U\otimes_{\Ox_U}\calE_s & \ra & \calA'_U\otimes_{\Ox_U}\calB_0\otimes_{\Ox_U}\calE_s \\
a\otimes x & \mapsto & a\otimes 1 \otimes x
\end{array}$$
with the canonical morphism
$$\calA'_U\otimes_{\Ox_U}\calB_0\otimes_{\Ox_U}\calE_s \ra \left ((\calA'_{\calI}\otimes_{\Ox_{\calI}}\calB)\circledast_{\calA_{\calI}\otimes_{\Ox_{\calI}}\calB}\calE \right )_s.$$
These compositions induce a morphism
\begin{equation}\label{eq1Akaza}
\calA'_{\calI}\otimes_{\calA_{\calI}}\calE \ra (\calA'_{\calI}\otimes_{\Ox_{\calI}}\calB)\circledast_{\calA_{\calI}\otimes_{\Ox_{\calI}}\calB}\calE.
\end{equation}
Conversely, the canonical morphism
$$(\calA_{\calI}\otimes_{\Ox_{\calI}}\calB)\boxtimes \calE \ra \sigma^*\calE$$
induces a morphism
$$\calA'_{\calI^2}\otimes_{\calA_{\calI^2}}((\calA_{\calI}\otimes_{\Ox_{\calI}}\calB)\boxtimes \calE)= (\calA'_{\calI}\otimes_{\Ox_{\calI}}\calB)\boxtimes \calE \ra \calA'_{\calI^2}\otimes_{\calA_{\calI^2}}\sigma^*\calE=\sigma^*(\calA'_{\calI}\otimes_{\calA_{\calI}}\calE).$$
By adjunction, we get
$$\sigma_!\left ( (\calA'_{\calI}\otimes_{\Ox_{\calI}}\calB)\boxtimes \calE\right ) \ra \calA'_{\calI}\otimes_{\calA_{\calI}}\calE.$$
By checking on fibers, this morphism yields
\begin{equation}\label{eq2Akaza}
(\calA'_{\calI}\otimes_{\Ox_{\calI}}\calB)\circledast_{\calA_{\calI}\otimes_{\Ox_{\calI}}\calB}\calE \ra \calA'_{\calI}\otimes_{\calA_{\calI}}\calE.
\end{equation}
Again, by checking on fibers, we prove that (\ref{eq1Akaza}) and (\ref{eq2Akaza}) are inverse to each other.
\end{proof}

\begin{proposition}\label{kraz827}
Let $(X,\Ox)$ be a ringed topos, $\calI$ a sheaf of monoids of $X,$ $\calA$ an $\Ox$-algebra, $\calB$ an $\calI$-indexed $\Ox$-algebra, $\calB'$ an $\calI$-indexed $\calB$-algebra and $\calE$ an $\calI$-indexed $\calA_{\calI}\otimes_{\Ox_{\calI}}\calB$-module. Then $\calE$ is canonically equipped with an $\calI$-indexed $\calB$-module structure and with an $\calA_{\calI}$-module structure and there exists a canonical isomorphism
$$\calB'\circledast_{\calB}\calE \xrightarrow{\sim} (\calB'\otimes_{\Ox_{\calI}}\calA_{\calI})\circledast_{\calB\otimes_{\Ox_{\calI}}\calA_{\calI}}\calE.$$
\end{proposition}

\begin{proof}
Let $p_1,p_2:\calI^2\ra \calI$ be the canonical projections and $\sigma:\calI^2\ra \calI$ the addition map.
The fact that $\calE$ is naturally equipped with an $\calA_{\calI}$-module structure follows from proposition \ref{propisoAkaza}. The $\calI$-indexed $\calB$-module structure on $\calE$ comes from the composition
$$\calB \boxtimes \calE \ra (\calA_{\calI}\otimes_{\Ox_{\calI}}\calB)\boxtimes \calE \ra \sigma^*\calE,$$
where the first arrow is induced by 
$$\calB \ra \calA_{\calI}\otimes_{\Ox_{\calI}} \calB,\ b\mapsto 1\otimes b$$
and the second one comes from the $\calI$-indexed $(\calA_{\calI}\otimes_{\Ox_{\calI}}\calB)$-module structure on $\calE.$

The canonical morphism
$$\calA_{\calI}\otimes_{\Ox_{\calI}}\calE \ra \calE$$
induces a morphism
$$\calB' \boxtimes (\calA_{\calI} \otimes_{\Ox_{\calI}}\calE)=(\calB'\otimes_{\Ox_{\calI}}\calA_{\calI})\boxtimes \calE \ra \calB' \boxtimes \calE$$
and then
$$\sigma_!((\calB'\otimes_{\Ox_{\calI}}\calA_{\calI})\boxtimes \calE) \ra \sigma_!(\calB' \boxtimes \calE).$$
By checking on fibers, we deduce a morphism
\begin{equation}\label{eq3Akaza}
(\calB'\otimes_{\Ox_{\calI}}\calA_{\calI})\circledast_{\calB\otimes_{\Ox_{\calI}}\calA_{\calI}} \calE \ra \calB' \circledast_{\calB} \calE.
\end{equation}
Conversely, the morphism
$$\calB' \ra \calB'\otimes_{\Ox_{\calI}}\calA_{\calI},\ b\mapsto b\otimes 1$$
induces
$$\sigma_! (\calB' \boxtimes \calE ) \ra \sigma_!\left ( (\calB' \otimes_{\Ox_{\calI}}\calA_{\calI})\boxtimes \calE \right )$$
and then
\begin{equation}\label{eq4Akaza}
\calB' \circledast_{\calB} \calE \ra (\calB'\otimes_{\Ox_{\calI}}\calA_{\calI})\circledast_{\calB\otimes_{\Ox_{\calI}}\calA_{\calI}} \calE.
\end{equation}
By checking on fibers, we prove that \eqref{eq3Akaza} and \eqref{eq4Akaza} are inverse to each other.
\end{proof}

\begin{proposition}\label{BCmatrix}
Let $(X,\Ox)$ be a ringed topos, $\calI$ a sheaf of monoids of $X,$ $\calA \ra \calB$ a morphism of $\calI$-indexed $\Ox$-algebras and $\calE$ an $\calI$-indexed $\calA$-module. There exists a canonical morphism of $\calI$-indexed $\calB$-algebras
\begin{equation}\label{erafisoprime}
\calB\circledast_{\calA}\mathscr{End}_{\calA}(\calE) \ra \mathscr{End}_{\calB}(\calB \circledast_{\calA} \calE).
\end{equation}
\end{proposition}

\begin{proof}
Let $\sigma:\calI^2 \ra \calI$ be the addition morphism and $p_1,p_2:\calI^2\ra \calI$ the projection morphisms. Let $s,t\in \Gamma(U,\calI)$ be sections of $\calI$ over an object $U.$ We define a morphism
$$\varphi_{(s,t)}:\calB_s\otimes_{\Ox_U}\mathscr{End}_{\calA}(\calE)_t \ra \mathscr{End}_{\calB}\left (\calB\circledast_{\calA}\calE\right )_{s+t}$$ as follows:
Let $V \ra U$ be a morphism in the underlying site of $X,$ $b\in \Gamma(V,\calB_s)$ and $f\in \Gamma \left (V,\mathscr{End}_{\calA}(\calE)_t\right ).$ By \ref{indremtakriz1}, $f$ is a morphism of $\calI_{|V}$-indexed $\calA_{|V}$-modules
$$f:\calE_{|V} \ra \sigma_{t_{|V}}^*\calE_{|V},$$
where $\sigma_{t_{|V}}:\calI_{|V} \ra \calI_{|V}$ is addition by $t_{|V}.$
Also, by \ref{indremtakriz1}, the morphism $b$ induces a morphism of $\calI_{|V}$-indexed $\calB_{|V}$-modules
$$b:\calB_{|V} \ra \sigma_{s_{|V}}^*\calB_{|V}.$$
We then get an $\Ox_{\calI_{|V}^2}$-linear morphism
$$
b\boxtimes f:\calB_{|V} \boxtimes \calE_{|V} \ra \left (\sigma_{s_{|V}}^*\calB_{|V} \right )\boxtimes \left (\sigma_{t_{|V}}^*\calE_{|V} \right )=\sigma_{\left (s_{|V},t_{|V} \right )}^*\left ( \calB_{|V}\boxtimes \calE_{|V}\right ).
$$
We then get an $\Ox_{\calI_{|V}}$-linear morphism
\begin{equation}\label{ren1}
\sigma_{|V!} \left (b\boxtimes f\right ): \sigma_{|V!} \left (\calB_{|V} \boxtimes \calE_{|V}\right ) \ra \sigma_{|V!} \left (\sigma_{\left (s_{|V},t_{|V} \right )}^*\left ( \calB_{|V}\boxtimes \calE_{|V}\right ) \right ).
\end{equation}
By the commutativity of diagram
$$
\begin{tikzcd}
\left (X_{/U} \right )_{/\calI_{|U}^2} \ar{r}{\sigma_{(s,t)}} \ar[swap]{d}{\sigma_{|U}} & \left (X_{/U} \right )_{/\calI_{|U}^2} \ar{d}{\sigma_{|U}} \\
\left (X_{/U} \right )_{/\calI_{|U}} \ar{r}{\sigma_{s+t}} & \left (X_{/U} \right )_{/\calI_{|U}},
\end{tikzcd}
$$
we get a canonical morphism
$$
\sigma_{|V!} \left (\sigma_{\left (s_{|V},t_{|V} \right )}^*\left ( \calB_{|V}\boxtimes \calE_{|V}\right ) \right ) \ra \sigma_{s_{|V}+t_{|V}}^*\sigma_{|V!}\left ( \calB_{|V}\boxtimes \calE_{|V}\right ).
$$
Composing this with \eqref{ren1}, we get an $\Ox_{\calI_{|V}}$-linear morphism
\begin{equation}\label{ren2}
\sigma_{|V!} \left (b\boxtimes f\right ): \sigma_{|V!} \left (\calB_{|V} \boxtimes \calE_{|V}\right ) \ra \sigma_{s_{|V}+t_{|V}}^*\left (\sigma_{|V!}\left ( \calB_{|V}\boxtimes \calE_{|V}\right ) \right ).
\end{equation}
By \ref{indproptensprod}, \eqref{ren2} induces a $\calB_{|V}$-linear morphism
$$
b\circledast_{\calA_{|V}}f: \calB_{|V}\circledast_{\calA_{|V}}\calE_{|V} \ra \sigma_{s_{|V}+t_{|V}}^*\left (\calB_{|V}\circledast_{\calA_{|V}}\calE_{|V}\right ).
$$
By \ref{ren3}, we get
$$
\left (b\circledast_{\calA_{|V}}f:\left (\calB\circledast_{\calA}\calE\right )_{|V} \ra \sigma_{s_{|V}+t_{|V}}^*\left (\calB\circledast_{\calA}\calE\right )_{|V} \right ) \in \mathscr{End}_{\calB}\left (\calB\circledast_{\calA}\calE\right )_{s+t}(V).
$$
We set $\varphi_{(s,t)}(b\otimes f)=b\circledast_{\calA_{|V}}f.$ By \ref{indremark73}, the morphisms $\varphi_{(s,t)}$ define a morphism
$$
\calB\boxtimes \mathscr{End}_{\calA}(\calE) \ra \sigma^*\mathscr{End}_{\calB}(\calB\circledast_{\calA}\calE),
$$
and then, by adjunction,
\begin{equation}\label{ren4}
\sigma_!\left ( \calB\boxtimes \mathscr{End}_{\calA}(\calE) \right ) \ra \mathscr{End}_{\calB}(\calB\circledast_{\calA}\calE).
\end{equation}
Again, by \ref{indremark73},  \eqref{ren4} induces
$$
\varphi:\calB\circledast_{\calA}\mathscr{End}_{\calA}(\calE) \ra \mathscr{End}_{\calB}(\calB \circledast_{\calA} \calE).
$$
\end{proof}

\begin{proposition}\label{Koko829}
Let $f:Y\ra X$ be a morphism of ringed topoi, $\calI$ a sheaf of monoids of $X,$ $\calA$ an $\calI$-indexed algebra of $X$ and $\calE$ and $\F$ two $\calI$-indexed $\calA$-modules. There exist a canonical isomorphism of $f^{-1}\calI$-indexed $f_{/\calI}^{-1}\calA$-modules
$$
f_{/\calI}^{-1}\calE \circledast_{f_{/\calI}^{-1}\calA}f_{/\calI}^{-1}\F \xrightarrow{\sim} f_{/\calI}^{-1}\left ( \calE \circledast_{\calA}\F \right )
$$
and a canonical isomorphism of $f^{-1}\calI$-indexed $f^*\calA$-modules
$$
f_{/\calI}^*\calE \circledast_{f_{/\calI}^*\calA}f_{/\calI}^*\F \xrightarrow{\sim} f_{/\calI}^*\left ( \calE \circledast_{\calA}\F \right ).
$$
\end{proposition}

\begin{proof}
Denote by $\sigma:\calI^2 \ra \calI$ and $\tau:\calI^3\ra \calI$ the addition morphisms. By \ref{lemtop15} and the commutative diagram
$$
\begin{tikzcd}
Y_{/f^{-1}\calI^2} \ar{r}{f_{/\calI^2}} \ar[swap]{d}{f^{-1}\sigma} & X_{/\calI^2} \ar{d}{\sigma} \\
Y_{/f^{-1}\calI} \ar[swap]{r}{f_{/\calI}} & X_{/\calI},
\end{tikzcd}
$$
we get a canonical isomorphism
$$
\left (f^{-1}\sigma\right )_!f_{/\calI^2}^{-1} \left (\calE \boxtimes \F \right ) \xrightarrow{\sim} f_{/\calI}^{-1}\sigma_!\left ( \calE\boxtimes \F\right ).
$$
We similarly have a canonical isomorphism
$$
\left (f^{-1}\tau\right )_!f_{/\calI^3}^{-1} \left (\calA \boxtimes \calE \boxtimes \F \right ) \xrightarrow{\sim} f_{/\calI}^{-1}\tau_!\left ( \calA \boxtimes \calE\boxtimes \F\right ).
$$
By definition of $\calE \circledast_{\calA} \F$ and the exactness of $f_{/\calI}^{-1},$ these isomorphisms fit into a commutative diagram
$$
\begin{tikzcd}
\left (f^{-1}\tau\right )_!f_{/\calI^3}^{-1} \left (\calA \boxtimes \calE \boxtimes \F \right ) \ar[sloped]{d}{\sim} \ar{r} & \left (f^{-1}\sigma\right )_!f_{/\calI^2}^{-1} \left (\calE \boxtimes \F \right ) \ar{r} \ar[sloped]{d}{\sim} & \left (f_{/\calI}^{-1}\calE \right ) \circledast_{f_{/\calI}^{-1}\calA} \left (f_{/\calI}^{-1}\F \right )   \ar[sloped,dashed]{d}{\sim} \ra 0 \\
f^{-1}_{/\calI}\tau_!(\calA \boxtimes \calE \boxtimes \F)  \ar{r} & f^{-1}_{/\calI}\sigma_!(\calE \boxtimes \F)  \ar{r} & f^{-1}_{/\calI}\left (\calE \circledast_{\calA} \F\right ) \ra 0,
\end{tikzcd}
$$
where the rows are exact and induced by \eqref{exseqindtens}.
Hence the existence of the third vertical isomorphism. The proof of the second isomorphism is similar.
\end{proof}

\begin{parag}
We end this subsection by defining local freeness for indexed modules. We introduce the following notation: let $(X,\Ox)$ be a ringed topos on a $\mathbb{U}$-site $C,$ $\calI$ a sheaf of monoids of $X,$ $\calA$ an $\calI$-indexed algebra and $\calE$ an $\calI$-indexed $\calA$-module. Let $m\in \Gamma(U,\calI)$ be a section of $\calI$ over an object $U$ of $C$ and consider the commutative diagram of topoi
\begin{equation}\label{diagtakriz22}
\begin{tikzcd}
\left (X_{/U} \right )_{/\calI_{|U}} \ar{dr}{j_{\calI_{|U}}} \ar{r}{\sim} & X_{/U\times \calI} \ar{r}{j_{U\times\calI/\calI}} \ar{d}{j_{\calI\times U/U}} & X_{/\calI} \ar{d}{j_{\calI}} \\
 & X_{/U} \ar{r}{j_U} & X
\end{tikzcd}
\end{equation}
We denote by $\sigma_m:\calI_{|U} \ra \calI_{|U},\ x\mapsto x+m$ and $\calE(m)=\sigma_m^*\calE_{|U},$ where $\calE_{|U}$ is the $\calI_{|U}$-indexed $\calA_{|U}$-module defined in \ref{Locprop1}. The sheaf $\calE(m)$ is clearly an $\calI_{|U}$-indexed $\calA_{|U}$-module.
\end{parag}

\begin{definition}\label{deflocfreeind}
Let $X$ be a scheme, $\calI$ a sheaf of abelian groups on $X,$ $\sigma:\calI^2 \ra \calI$ the addition map, $\calA$ an $\calI$-indexed algebra and $\calE$ an $\calI$-indexed $\calA$-module.
\begin{enumerate}
\item We say that $\calE$ is a \emph{free $\calI$-indexed $\calA$-module} (resp. \emph{free $\calI$-indexed $\calA$-module of finite rank}) if there exists sections (resp. finitely many sections) $m_{i}\in \Gamma(X,\calI), i\in I$ satisfying
$$\calE \cong \bigoplus_{i\in I} \calA(m_{i}).$$
In this case, a \emph{local basis} of the $\calI$-indexed $\calA$-module $\calE$ is a family of local sections $(e_i)_{i\in I}$ of $\prod_{i\in I}\calE_{s_i},$ where the $s_i$'s are local sections of $\calI$ over an étale $X$-scheme $U,$ satisfying: for all $V\ra U$ étale and all local sections $s\in \Gamma (U,\calI)$ and $x\in \Gamma(V,\calE_s),$ there exists unique local sections $(a_i)_{i\in I}$ of $\prod_{i\in I}\Gamma(V,\calA_{s-s_i})$ (all equal to zero but finitely many) such that $x=\sum_{i\in I}a_ie_i,$ where $a_ie_i$ denotes the image of $a_i\otimes e_i$ by $\calA_{s-s_i} \otimes_{\Ox_U} \calE_{s_i} \ra \calE_s.$
\item We say that $\calE$ is a \emph{locally free $\calI$-indexed $\calA$-module} (resp. \emph{locally free $\calI$-indexed $\calA$-module of finite rank}) if there exists an étale covering $(U_i)_{i\in I}$ of $X$ such that for every $i\in I,$ the $\calI_{|U_i}$-indexed $\calA_{|U_i}$-module $\calE_{|U_i}$ is free (resp. free of finite rank).
\end{enumerate}
\end{definition}

\begin{proposition}\label{propfreetakriz}
Let $X$ be a scheme, $\calI$ a sheaf of monoids on $X,$ $\calA \ra \calA'$ a morphism of $\calI$-indexed algebras and $\calE$ an $\calI$-indexed $\calA$-module. If $\calE$ is a locally free $\calI$-indexed $\calA$-module then $\calA' \circledast_{\calA}\calE$ is a locally free $\calI$-indexed $\calA'$-module. In addition, if $(e_i)_{i\in I}$ is a basis for the $\calI$-indexed $\calA$-module $\calE,$ where $e_i$ is over a section $s_i$ of $\calI,$ then $(1\otimes e_i)_{i\in I}$ is a basis for the $\calI$-indexed $\calA'$-module $\calA'\circledast_{\calA}\calE.$ Here, $1\otimes e_i$ is considered as a section of $\calA'\circledast_{\calA}\calE$ via the morphism $\calA'_0\otimes_{\Ox_X}\calE_{s_i} \ra (\calA' \circledast_{\calA}\calE)_{s_i}$ (\ref{indproptensprod}).
\end{proposition}

\begin{proof}
It is sufficient to prove that, for any section $m\in \Gamma(X,\calI),$
$$\calA(m)\circledast_{\calA}\calA' \cong \calA'(m).$$
This is a special case of the lemma \ref{ren5} below.
\end{proof}

\begin{lemma}\label{ren5}
Let $X$ be a scheme, $\calI$ a sheaf of monoids on $X,$ $\calA$ an $\calI$-indexed algebra, $\calE$ an $\calI$-indexed $\calA$-module and $m\in \Gamma(X,\calI).$ Then there exists a canonical isomorphism
$$\calA(m)\circledast_{\calA}\calE \xrightarrow{\sim} \calE(m).$$
\end{lemma}

\begin{proof}
Let $p_1,p_2:\calI^2 \ra \calI$ be the canonical projections, $\sigma:\calI^2 \ra \calI$ the addition map, $\sigma_m:\calI \ra \calI,\ x\mapsto x+m$ the addition by $m,$ $\sigma_{(m,0)}:\calI^2 \ra \calI^2,\ x\mapsto x+(m,0)$ the addition by $(m,0)$ and $\pi:\calA \boxtimes \calE \ra \sigma^*\calE$ the canonical map. We have the following commutative diagrams
$$
\begin{tikzcd}
\calI^2 \ar{r}{\sigma_{(m,0)}} \ar{d}{p_1} & \calI^2 \ar{d}{p_1} & \calI^2 \ar{r}{\sigma_{(m,0)}} \ar{dr}{p_2} & \calI^2 \ar{d}{p_2} & \calI^2 \ar{r}{\sigma_{(m,0)}} \ar{d}{\sigma} & \calI^2 \ar{d}{\sigma} \\
\calI \ar{r}{\sigma_m} & \calI & & \calI & \calI \ar{r}{\sigma_m} & \calI
\end{tikzcd}
$$
We then obtain a morphism
$$\sigma_{(m,0)}^*\pi: \calA(m)\boxtimes \calE \ra \sigma^*\calE(m).$$
By adjunction, we get a morphism
$$\sigma_!(\calA(m) \boxtimes \calE) \ra \calE(m).$$
By checking on fibers, we prove it induces an isomorphism
$$\calA(m)\circledast_{\calA}\calE \xrightarrow{\sim} \calE(m).$$
\end{proof}

\begin{proposition}
Let $X$ be a scheme, $\calI$ a sheaf of abelian groups of $X_{\text{ét}},$ $\calA \ra \calB$ a morphism of $\calI$-indexed algebras and $\calE$ a locally free $\calI$-indexed $\calA$-module of finite rank \eqref{deflocfreeind}. The canonical morphism \eqref{erafisoprime} is then an isomorphism:
\begin{equation}\label{erafiso}
\calB\circledast_{\calA}\mathscr{End}_{\calA}(\calE) \xrightarrow{\sim} \mathscr{End}_{\calB}(\calB \circledast_{\calA} \calE).
\end{equation}
\end{proposition}

\begin{proof}
If $U\ra X$ is an étale morphism, then $\left (\calB \circledast_{\calA}\mathscr{End}_{\calA}(\calE) \right )_{|U}=\calB_{|U}\circledast_{\calA_{|U}}\mathscr{End}_{\calA}(\calE)_{|U}$ by \ref{ren3} and $\mathscr{End}_{\calA}(\calE)_{|U}=\mathscr{End}_{\calA_{|U}}\left (\calE_{|U} \right )$ by definition of $\mathscr{End}_{\calA}$ (\ref{inddef2}). We can then suppose that $\calE$ is a free $\calI$-indexed $\calA$-module. It follows that there exists $s_1,\hdots,s_d\in \Gamma(X,\calI)$ and an isomorphism
$$
\calE\xrightarrow{\sim} \bigoplus_{i=1}^d\calA(s_i).
$$
The result then follows from \ref{ren5} and the canonical isomorphisms
$$
\bigoplus_{i,j=1}^d\calA(s_j-s_i) \xrightarrow{\sim} \mathscr{End}_{\calA}\left (\bigoplus_{i=1}^d\calA(s_i) \right ),\ \bigoplus_{i,j=1}^d\calB(s_j-s_i) \xrightarrow{\sim} \mathscr{End}_{\calB}\left (\bigoplus_{i=1}^d\calB(s_i) \right ).
$$
\end{proof}

\subsection*{Indexed connections}

\begin{parag}\label{indparag77}
In this subsection, we fix a morphism $f:X\ra S$ of logarithmic schemes or logarithmic formal schemes and an étale sheaf $\calI$ on $X.$ Applying the results of the previous section to the small étale topos $X_{\text{ét}},$ we get three functors
\begin{equation}\label{indeq71}
j_{\calI !}:\calI_{\text{ét}}\ra X_{\text{ét}},\ j_{\calI}^*:X_{\text{ét}}\ra \calI_{\text{ét}},\ j_{\calI *}:\calI_{\text{ét}}\ra X_{\text{ét}}
\end{equation}
In the rest of this manuscript and for any étale sheaf $\F$ on $X$ (resp. any morphism $f:\F \ra \calG$ of étale sheaves on $X$), we denote $j_{\calI}^*\F$ (resp. $j_{\calI}^*f$) by $\F_{\calI}$ (resp. $f_{\calI}$). In particular, we denote $j_{\calI}^*\Ox_X$ by $\Ox_{X,\calI}$ and, for any positive integer $i,$ we denote the sheaf $j_{\calI}^*\left (\omega^i_{X/S}\right )$ by $\omega^i_{X/S,\calI}.$
\end{parag}

\begin{definition}\label{indprop42}
Let $\calE$ be an $\Ox_{X,\calI}$-module and $\lambda$ a global section of $\Ox_{X,\calI}.$ A \emph{logarithmic $\lambda$-connection on $\calE$} is an $(f^{-1}\Ox_S)_{\calI}$-linear morphism $\nabla:\calE\ra \calE\otimes_{\Ox_{X,\calI}}\omega^1_{X/S,\calI}$ satisfying the Leibniz rule
$$\nabla(ax)=a\nabla(x)+\lambda x\otimes d_{\calI}a,$$
for all local sections $a$ and $x$ of $\Ox_{X,\calI}$ and $\calE$ respectively and where $d:\Ox_X\ra \omega^1_{X/S}$ denotes the universal derivation. 
Logarithmic $1$-connections are called \emph{logarithmic connections.}
\end{definition} 

In the rest of this manuscript, we will always drop the term "logarithmic" and refer to logarithmic connections and logarithmic Higgs fields simply by connections and Higgs fields.

\begin{proposition}
Let $\calE$ be an $\Ox_{X,\calI}$-module and $\lambda$ a global section of $\Ox_{X,\calI}.$ A $\lambda$-connection $\nabla$ on $\calE$ induces, for any positive integer $i,$ a morphism
$$\nabla^i:\calE\otimes_{\Ox_{X,\calI}}\omega^i_{X/S,\calI}\ra \calE\otimes_{\Ox_{X,\calI}}\omega^{i+1}_{X/S,\calI}$$
defined for any local sections $x$ and $\omega$ of $\calE$ and $\omega^i_{X/S,\calI}$ respectively by
$$\nabla^i(x\otimes \omega)=\nabla(x)\wedge \omega+\lambda x\otimes d_{\calI}\omega$$
where $\nabla(x)\wedge \omega$ denotes the image of $\nabla(x)\otimes \omega$ by the canonical morphism $$\calE\otimes_{\Ox_{X,\calI}} \omega^1_{X/S,\calI}\otimes_{\Ox_{X,\calI}}\omega^i_{X/S,\calI}\ra \calE\otimes_{\Ox_{X,\calI}}\omega^{i+1}_{X/S,\calI}.$$
\end{proposition}

\begin{proof}
Same as \ref{prop42}
\end{proof}

\begin{definition}\label{indhiggs}
Let $\calE$ be an $\calI$-indexed $\Ox_{X}$-module, $\lambda$ a global section of $\Ox_{X,\calI}$ and $\nabla:\calE\ra \calE\otimes_{\Ox_{X,\calI}}\omega^1_{X/S,\calI}$ a $\lambda$-connection on $\calE.$ The \emph{curvature of $\nabla$,} denoted by $K(\nabla),$ is the composition $\nabla^1\circ \nabla$ (where $\nabla^1$ is defined in \ref{indprop42}). The curvature $K(\nabla)$ is $\Ox_{X,\calI}$-linear and we say that the $\lambda$-connection $\nabla$ is \emph{integrable} if $K(\nabla)=0.$ Integrable $0$-connections are called \emph{Higgs fields.} We denote by $\boldsymbol{\op{MIC^{ind}}(X/S)}$ (resp. $\boldsymbol{\op{HIG^{ind}}(X/S)}$) the category of $\Ox_{X,\calI}$-modules equipped with an integrable connection (resp. Higgs field).
\end{definition}

\begin{remark}
\begin{enumerate}
\item Let $\calE$ be an $\Ox_{X,\calI}$-module and $\lambda$ a global section of $\Ox_{X,\calI}.$ By \ref{indremark73}, the data of a $\lambda$-connection $$\nabla : \calE \ra \calE \otimes_{\Ox_{X,\calI}} \omega^1_{X/S,\calI}$$
is equivalent to the data of a compatible family of $\lambda$-connections
$$\nabla_i:\calE_i \ra \calE_i \otimes_{\Ox_U}\omega^1_{U/S}$$
for every section $i\in \calI(U).$ Furthermore, $K(\nabla)_i=K(\nabla_i)$ and $\nabla$ is integrable if and only if $\nabla_i$ is integrable for every section $i$ of $\calI.$

\item Let $(\calE_1,\nabla_1)$ and $(\calE_2,\nabla_2)$ be two $\Ox_{X,\calI}$-modules equipped with $\lambda$-connections, where $\lambda$ is a global section of $\Ox_{X,\calI}.$ Then $\calE=\calE_1\otimes_{\Ox_{X,\calI}}\calE_2$ is canonically equipped with a $\lambda$-connection $\nabla$ defined by
$$
\nabla(x_1\otimes x_2)=\nabla_1(x_1)\otimes x_2+x_1\otimes \nabla_2(x_2),
$$
for local sections $x_1$ and $x_2$ of $\calE_1$ and $\calE_2$ respectively.
\end{enumerate}
\end{remark}

\subsection*{Indexed algebra associated to a logarithmic structure}

\begin{parag}\label{paragAX}
Let $X$ be a logarithmic scheme or logarithmic formal scheme.
We denote by $\calI$ the sheaf $\ov{\calM}^{gp}_X$ and by $\sigma:\calI^2\ra \calI$ the addition map and consider the functors introduced in (\ref{indeq71}):
\begin{equation}
j_{\calI !}:\calI_{\text{ét}}\ra X_{\text{ét}},\ j_{\calI}^*:X_{\text{ét}}\ra \calI_{\text{ét}},\ j_{\calI *}:\calI_{\text{ét}}\ra X_{\text{ét}}.
\end{equation}
The exact sequence of abelian groups
$$0 \ra \Ox_X^* \xrightarrow{\alpha_X^{-1}} \calM_X^{gp} \ra \calI \ra 0$$
turns $\calM_X^{gp}$ into an $\Ox_{X,\calI}^*$-torsor of $\calI_{\text{ét}}.$ To this torsor corresponds the invertible $\Ox_{X,\calI}$-module
\begin{equation}\label{Agroupe}
\calA_X=\Ox_{X,\calI}\times_{\calI}\calM_X^{gp}/\sim,
\end{equation}
where $\sim$ is the relation
\begin{equation}\label{eq8422c}
(a,m)\sim (a',m') \iff \exists m''\in \calM_X^{\times},\begin{cases} m+m''=m' \\ a'\alpha_X(m'')=a\end{cases}
\end{equation}
for any local sections $a,a'$ of $\Ox_{\calI,\ov{m}}$ and $m,m'$ of $\calM_X^{gp}$ such that $\ov{m}=\ov{m'}.$

For a local section $m\in\calM_X^{gp}(U)$ over an étale $X$-scheme $U,$ the $\Ox_U$-module
\begin{equation}\label{eq8422b}
\calA_{X,\ov{m}}=\Ox_U\times \calM_{X,\ov{m}}^{gp}/\sim
\end{equation}
admits a basis
\begin{equation}\label{eqAbase}
e_m=(1,m).
\end{equation}
For any local sections $m,m'\in\calM_X^{gp}(U),$ we consider the morphism
\begin{equation}
\begin{array}[t]{clc}
\calA_{X,\ov{m}}\otimes_{\Ox_U}\calA_{X,\ov{m'}} & \ra & \calA_{X,\ov{m}+\ov{m'}}\\
a e_m\otimes a'e_{m'} & \mapsto & aa'e_{m+m'}
\end{array}
\end{equation}
These morphisms are well-defined and compatible and so, by \ref{era2indrem}, they define a morphism
\begin{equation}\label{Agroupe2}
\pi:\calA_X\boxtimes \calA_X \ra \sigma^*\calA_X
\end{equation}
which defines an $\calI$-indexed $\Ox_X$-algebra structure on $\calA_X$ (\ref{inddef} (2)).
\end{parag}

\begin{proposition}[\cite{Lor2000} 1.6 and 1.7]
Let $X\ra S$ be a morphism of logarithmic schemes or logarithmic formal schemes. We denote by $\calI$ the sheaf $\ov{\calM}^{gp}_X.$ Then, the $\calI$-indexed $\Ox_X$-algebra $\calA_X$ \eqref{Agroupe} admits a unique connection \eqref{indprop42}
$$d_{\calA_X}:\calA_X \ra \calA_X \otimes_{\Ox_{X,\calI}}\omega^1_{X/S,\calI}$$
given by 
\begin{equation}\label{eq8431Koko}
d_{\calA_X}(e_m)=e_m\otimes \op{dlog}m,
\end{equation}
for any local section $m$ of $\calM_X^{gp}.$ In addition, $d_{\calA_X}$ is integrable.
\end{proposition}

\begin{definition}\label{Bgroupe}
Let $X\ra S$ be a morphism of logarithmic schemes or logarithmic formal schemes, $\calI$ the sheaf $\ov{\calM}^{gp}_X,$ $p_1,p_2:\calI^2\ra \calI$ the canonical projections and $\sigma:\calI^2\ra \calI$ the addition map. 
A connection $\nabla$ on an $\calI$-indexed $\calA_X$-module $\calE$ (\ref{indprop42}) is said to be \emph{admissible} if the map $\calA_X\boxtimes \calE \ra \sigma^*\calE$ defining the $\calA_X$-module structure on $\calE$ is horizontal with respect to the connections
$$
\begin{array}[t]{clc}
\calA_X \boxtimes \calE & \ra & (\calA_X\boxtimes \calE) \otimes_{\Ox_{X,\calI^2}} \omega^1_{X/S,\calI^2}\\
a \otimes x & \mapsto & \left ( p_1^*d_{\calA_X} \right )(a)\otimes x+a\otimes \left (p_2^*\nabla \right )(x)
\end{array}
$$
and $\sigma^*\nabla.$
\end{definition}

\begin{parag}\label{Bgroupe2}
Let $X\ra S$ be a morphism of logarithmic schemes or logarithmic formal schemes and $\calI$ the sheaf $\ov{\calM}^{gp}_X.$ Then the canonical connection $d_{\calA_X}$ is admissible. Indeed, for any local sections $s,t:U\ra \calM_X^{gp}$ over an étale $X$-scheme $U,$
$$
\left (\pi \otimes \op{Id}_{\omega^1_{U/S}} \right ) \left (d_{\calA_X}(e_s)\otimes e_t+e_s\otimes d_{\calA_X}(e_t) \right )=d_{\calA_X}(e_{s+t}),$$
where $\pi$ is the multiplication morphism \eqref{Agroupe2}. So, for any local section $a,b\in \Gamma(U,\Ox_X),$
\begin{alignat*}{2}
\left (\pi \otimes \op{Id}_{\omega^1_{U/S}} \right ) \left (d_{\calA_X}(ae_s)\otimes (be_t)+(ae_s)\otimes d_{\calA_X}(be_t) \right ) &= e_{s+t}\otimes \left ( d(ab)+ab \op{dlog}(s+t) \right ) \\
&= d_{\calA_X} \left (abe_{s+t} \right ).
\end{alignat*}
and so $\pi:\calA_X \boxtimes \calA_X \ra \sigma^*\calA_X$ is horizontal with respect to the connections $\sigma^*d_{\calA_X}$ and $\op{Id}_{\calA_X}\boxtimes d_{\calA_X}+d_{\calA_X}\boxtimes \op{Id}_{\calA_X}.$
It follows that the kernel of the connection $d_{\calA_X}$ is an $\calI$-indexed $\Ox_X$-subalgebra of $\calA_X.$ We denote it by $\calB_{X/S}.$
\end{parag}

\begin{proposition}\label{lemBres}
Let $X\ra S$ be a morphism of logarithmic schemes, $U\ra X$ an étale and strict morphism, $\calI=\ov{\calM}_X^{gp}$ and $\calI_U=\calI_{|U}.$
Let $u:U_{\text{ét}/\cal_U} \ra X_{\text{ét}/\calI}$ the morphism of ringed topoi induced by $U\ra X.$ Then, there exists a canonica isomorphism of $\calI$-indexed algebras
$$
u^*\calB_{X/S} \xrightarrow{\sim} \calB_{U/S}.
$$
\end{proposition}

\begin{proof}
The morphism $U\ra X$ is étale and strict so we have a canonical isomorphism
$$
u^*\omega^1_{X/S,\calI} \xrightarrow{\sim} \omega^1_{U/S,\calI_U}.
$$
We have the following commutative diagram of $\Ox_{U',\calI_U}$-modules:
$$
\begin{tikzcd}
0 \ar{r} & u^*\calB_{X/S} \ar{r} \ar[dashed]{d} & u^*\calA_X \ar[sloped]{d}{\sim} \ar{r}{u^*d_{\calA_X}} & u^*\calA_{X} \otimes_{\Ox_{U,\calI_U}} \omega^1_{U/S,\calI} \ar[sloped]{d}{\sim} \\
0 \ar{r} & \calB_{U/S} \ar{r} & \calA_U \ar{r}{d_{\calA_U}} & \calA_U \otimes_{\Ox_{U,\calI_U}} \omega^1_{U/S,\calI}.
\end{tikzcd}
$$
The lines are exact by definition of $\calB_{X/S}$ and $\calB_{U/S}$ and the left exactness of $u$ (which comes from the flatness of $U\ra X$). The vertical isomorphisms are given in (\cite{Lor2000} I 3.1 page 273). Hence the existence of the dashed arrow and the fact that it is an isomorphism.
\end{proof}

We end this section by recalling the following result of Lorenzon: 
\begin{theorem}[\cite{Lor2000} 2.6]\label{thmlor}
Let $f:X\ra S$ be a log smooth morphism of fine logarithmic schemes of characteristic $p.$ Denote by $\calI$ the sheaf $\ov{\calM}^{gp}_X.$ Then $\calA_X$ is an $\calI$-indexed $\calB_{X/S}$-module locally free of finite rank. More precisely, if $m_1,\hdots,m_d \in \Gamma(X,\calM_X)$ are local coordinates for $X/S$ and, for every $I=(I_1,\hdots,I_d)\in \Z^d,$ $m_I=\sum_{i=1}^dI_im_i\in \Gamma(X,\calM_X^{gp})$ then $(e_{m_I})_{I\in \{0,\hdots p-1\}^d}$ is a basis for the $\calI$-indexed $\calB_{X/S}$-module $\calA_X.$
\end{theorem}

\newpage

\section{A logarithmic Shiho functor}

\begin{parag}\label{parag48}
In this section, we fix a log smooth morphism $f:\frakX\ra \frakS$ of fs $p$-adic logarithmic locally Noetherian formal schemes flat over $\Z_p$ (\ref{parag63}). For any integer $n\ge 1,$ we denote by $\frakX_n$ and $\frakS_n$ the logarithmic schemes obtained from $\frakX$ and $\frakS$ respectively by reduction modulo $p^n,$ as in \ref{prop69}. Set $X=\frakX_1$ and $S=\frakS_1.$ We use the same notations as in \ref{PFrob}. Namely, we denote by $F_S:S\ra S$ the absolute Frobenius morphism of $S$ \eqref{Not3}, by $F_{X/S}:X\ra X''$ the relative Frobenius morphism of $X$ with respect to $S,$ by $F_1:X\ra X'$ the exact relative Frobenius of $X$ with respect to $S$ and by $F_1^{\#}:\Ox_{X'}\ra F_{1*}\Ox_X$ the associated morphism of structural rings. Denote also by $G:X'\ra X''$ and $\pi:X'\ra X$ the morphisms defined in \eqref{diag51}.
\end{parag}

\begin{parag}
In the rest of this section, we suppose that we have an $\frakS$-lifting $F:\frakX\ra \frakX'$ of the exact relative Frobenius $F_1:X\ra X'.$ In other words, we suppose that there exists an fs $p$-adic logarithmic formal scheme $\frakX'$ over $\frakS$ and an $\frakS$-morphism $F:\frakX\ra \frakX'$ that fit into cartesian squares
$$\begin{tikzcd}
X'\ar{r}\ar{d} & S\ar{d} & & X\ar{r}{F_1}\ar{d} & X'\ar{d}\\
\frakX'\ar{r} & \frakS & & \frakX\ar{r}{F} & \frakX'
\end{tikzcd}$$
Note that, by (\cite{Ahmed2010} 2.1.36), $\frakX'$ is locally Noetherian.
We also suppose that $F$ is log flat.
For any integer $n\ge 1,$ we denote by $F_n:\frakX_n\ra \frakX_n'$ the base change of $F:\frakX\ra \frakX'$ by $\frakS_n\ra \frakS.$
\end{parag}

\begin{lemma}\label{lem12}
Consider the composition $X' \xrightarrow{\pi} X \ra \frakX$ where $\pi$ is given in \eqref{diag51}. Let $\frakU \ra \frakX$ be an étale morphism of $p$-adic formal schemes and $\frakU' \ra \frakX'$ the unique étale morphism such that $\frakU \times_{\frakX} X'=\frakU' \times_{\frakX'}X'.$ Let $m$ be a local section of $\calM_{\frakX}$ (resp. $\calM_{\frakX_n}$) over $\frakU$ and $m_1$ its image in $\Gamma \left (\frakU\times_{\frakX}X,\calM_X \right ).$ By \ref{prop69}, after eventually reducing $\frakU,$ there exists a lifting $m'$ of $\pi^{\flat}m_1\in \Gamma \left (\frakU\times_{\frakX}X',\calM_{X'} \right )$ to $\Gamma \left (\frakU\times_{\frakX}X',\calM_{\frakX'} \right )$ (resp. $\Gamma \left (\frakU\times_{\frakX}X',\calM_{\frakX'_n} \right )$). Then there exists an invertible local section $u$ of $\calM_{\frakX}$ (resp. $\calM_{\frakX_n}$) such that $u+pm=F^{\flat}(m')$ (resp. $u+pm=F_n^{\flat}(m')$) and $\alpha_{\frakX}(u)=1+pb$ (resp. $\alpha_{\frakX_n}(u)=1+pb$) for a local section $b$ of $\Ox_{\frakX}$ (resp. $\Ox_{\frakX_n}$).
\end{lemma}

\begin{proof}
We prove the lemma for $\frakX.$ The proof for $\frakX_n$ follows.
Let $\gamma:\calM_{\frakX}\ra \Ox_X$ be the composition $\calM_{\frakX}\xrightarrow{\alpha_{\frakX}} \Ox_{\frakX}\ra \Ox_X,$ where $\Ox_{\frakX}\ra \Ox_X$ is reduction modulo $p.$ The morphism $X\ra \frakX$ is strict (\ref{prop69}) so $\calM_X=\calM_{\frakX}\oplus_{\gamma^{-1}(\Ox_X^*)}\Ox_X^*.$ The commutative diagram
$$\begin{tikzcd}
X\ar{r}\ar{d}{F_1} & \frakX\ar{d}{F}\\
X'\ar{r} & \frakX'
\end{tikzcd}$$
gives rise to the following commutative diagram of sheaves of monoids
$$\begin{tikzcd}
\calM_X & \calM_{\frakX}\ar{l}\\
\calM_{X'}\ar{u}{F_1^{\flat}} & \calM_{\frakX'}\ar{u}{F^{\flat}}\ar{l}.
\end{tikzcd}$$
It follows that, in $\calM_X,$ we have the equality
$$(pm,1)=(F^{\flat}(m'),1).$$
By the definition of amalgamated sums in the category of monoids (\cite{Kat89}  (1.3)), there exist local sections $v$ and $v'$ of $\gamma^{-1}(\Ox_X^*)$ such that $F^{\flat}(m')+v'=pm+v$ and $\gamma(v')=\gamma(v).$ Since $\gamma(v')=\gamma(v),$ there exists a local section $c$ of $\Ox_{\frakX}$ such that $\alpha_{\frakX}(v)=\alpha_{\frakX}(v')+pc.$ The local sections $\gamma(v)$ and $\gamma(v')$ are both invertible so $\alpha_{\frakX}(v)$ and $\alpha_{\frakX}(v')$ are also invertible and so $v$ and $v'$ are invertible in $\calM_{\frakX}.$ It follows that we can take $u=v-v'$ and $\alpha_{\frakX}(u)=\alpha_{\frakX}(v)\alpha_{\frakX}(v')^{-1}=1+pb$ where $b=\alpha_{\frakX}(v')^{-1}c.$
\end{proof}

\begin{lemma}\label{lemcalc}
Suppose that $\frakX$ and $\frakS$ are equipped with frames $\frakX \ra [Q]$ and $\frakS\ra [P]$ and that $f$ underlies a morphism of framed logarithmic formal schemes. Let $\Delta:\frakX \ra \frakY:=\frakX\times_{\frakS,[Q]}^{\op{log}}\frakX$ and $\Delta':\frakX' \ra \frakY':=\frakX' \times_{\frakS,[Q']}^{\op{log}}\frakX'$ be the strict diagonal immersions \eqref{prop76} and $\frakX \ra \frakU \ra \frakY$ and $\frakX' \ra \frakU' \ra \frakY'$ factorizations of $\Delta$ and $\Delta'$ respectively into a closed immersion followed by an open one. Denote by $\calI$ and $\calI'$ the respective ideals of $\frakX \ra \frakU$ and $\frakX' \ra \frakU'$ and by $p_1,p_2:\frakY\ra \frakX$ and $p_1',p_2':\frakY'\ra \frakX'$ the canonical projections. Suppose that there exists a morphism $F:(\frakX,Q) \ra (\frakX',Q')$ of framed logarithmic formal schemes, over $(\frakS,P),$ lifting the exact relative Frobenius $F_1:X \ra X'.$ Denote by $G:\frakY\ra \frakY'$ the morphism induced by $F:(\frakX,Q)\ra (\frakX',Q').$ Let $m_1$ be a local section of $\calM_X,$ $m'_1$ its image in $\calM_{X'}$ by $\pi^{\flat}:\calM_X\ra \pi_*\calM_{X'}$ \eqref{diag51}, $m$ (resp. $m'$) a local lifting of $m_1$ (resp. $m'_1$) to $\calM_{\frakX}$ (resp. $\calM_{\frakX'}$).
Let $\eta(m)\in \Delta^{-1}\calI,\ \eta(m')\in \Delta'^{-1}\calI',$ $\mu(m)=\eta(m)+1$ and $\mu(m')=\eta(m')+1$ as defined in \ref{parag77}. By \ref{lem12}, there exists a local invertible section $u$ of $\calM_{\frakX}$ and a local section $b$ of $\Ox_{\frakX}$ such that $F^{\flat}(m')=pm+u$ and $\alpha_{\frakX}(u)=1+pb.$
Then
\begin{equation}\label{eq1321}
\Delta^{-1}G^{\#}\left ( \eta(m')\right ) = \left (\eta(m)^p+\sum_{k=1}^{p-1}\begin{pmatrix}p\\k \end{pmatrix}\eta(m)^k+1 \right )\alpha_{\frakY}(p_2^{\flat}u-p_1^{\flat}u)-1
\end{equation}
and
\begin{equation}\label{eq1322}
\alpha_{\frakY}(p_2^{\flat}u-p_1^{\flat}u)=\frac{1+pp_2^{\#}(b)}{1+pp_1^{\#}(b)}.
\end{equation}
\end{lemma}

\begin{proof}
By \ref{parag77}, we have the following commutative diagram with exact rows
$$
\begin{tikzcd}
0 \ar{r} & \Delta'^{-1}(1+\calI') \ar{r}{\lambda'} \ar{d}{\Delta^{-1}G^{\#}} & \Delta'^{-1}\calM_{\frakY'} \ar{d}{\Delta^{-1}G^{\flat}} \ar{r} & \calM_{\frakX'} \ar{d}{F^{\flat}} \ar{r} & 0 \\
0 \ar{r} & \Delta^{-1}(1+\calI) \ar{r}{\lambda}  & \Delta^{-1}\calM_{\frakY} \ar{r} & \calM_{\frakX} \ar{r} & 0.
\end{tikzcd}
$$
In the rest of this proof, we drop $\Delta^{-1}$ and $\Delta'^{-1}$ to lighten the notation. By definition, we have
$$p_1'^{\flat}(m')+\lambda'(\mu(m'))=p_2'^{\flat}(m').$$
Applying $G^{\flat},$ we get
$$p_1^{\flat}F^{\flat}(m')+\lambda(G^{\#}(\mu(m'))=p_2^{\flat}F^{\flat}(m').$$
So we have
$$pp_1^{\flat}(m)+p_1^{\flat}(u)+\lambda(G^{\#}(\mu(m')))=pp_2^{\flat}(m)+p_2^{\flat}(u).$$
By definition,
$$p_1^{\flat}(m)+\lambda(\mu(m))=p_2^{\flat}(m).$$
It follows that
$$\lambda(G^{\#}(\mu(m')))=\lambda(\mu(m)^p)+p_2^{\flat}u-p_1^{\flat}u.$$
But $\lambda(\alpha_{\frakY}(p_2^{\flat}u-p_1^{\flat}u))=p_2^{\flat}u-p_1^{\flat}u$ by definition of $\lambda.$ The equality \eqref{eq1321} then follows from the injectivity of $\lambda.$ For \eqref{eq1322}, we have
\begin{equation*}
\alpha_{\frakY}(p_2^{\flat}u-p_1^{\flat}u) = \frac{p_2^{\#}\alpha_{\frakX}(u)}{p_1^{\#}\alpha_{\frakX}(u)}=\frac{1+pp_2^{\#}(b)}{1+pp_1^{\#}(b)}.
\end{equation*}
\end{proof}

\begin{parag}
We have the following commutative diagram
$$\begin{tikzcd}
F_1^*G^*\omega^1_{X''/S}\ar{rd}{dF_{X/S}}\ar{d}{F_1^*dG} & \\
F_1^*\omega^1_{X'/S}\ar{r}{dF_1} & \omega^1_{X/S}
\end{tikzcd}.$$
Since $dF_{X/S}=0$ and $dG$ is an isomorphism (because $G$ is étale), we deduce that $dF_1=0.$
Let $n\ge 1$ and denote by $F_n:\frakX_n\ra \frakX'_n$ the reduction of $F$ modulo $p^n.$ By the following commutative diagram
$$\begin{tikzcd}
F_n^*\omega^1_{\frakX_n'/\frakS_n}\ar{d}\ar{r}{dF_n} & \omega^1_{\frakX_n/\frakS_n} \ar{d}\\
F_1^*\omega^1_{X'/S}\ar{r}{dF_1} & \omega^1_{X/S}
\end{tikzcd}$$
(where the vertical arrows are reduction modulo $p$) and the fact that $dF_1=0,$ we get that $\op{Im}(dF_n)\subset p\omega^1_{\frakX_n/\frakS_n}.$
Now by the flatness of $\frakX$ over $\op{Spf} \Z_p$ and the exact sequence
$$0\ra \Z/p^n\Z\xrightarrow{\times p} \Z/p^{n+1}\Z\rightarrow \Z/p\Z\ra 0,$$
we get the exact sequence
$$0\ra \Ox_{\frakX_n}\xrightarrow{\times p} \Ox_{\frakX_{n+1}}\rightarrow \Ox_{\frakX_1} \ra 0.$$
Since $f:\frakX\ra \frakS$ is smooth, the $\Ox_{\frakX_{n+1}}$-module $\omega^1_{\frakX_{n+1}/\frakS_{n+1}}$ is locally free and so we get the exact sequence
$$0\ra \omega^1_{\frakX_n/\frakS_n}\xrightarrow{\times p}\omega^1_{\frakX_{n+1}/\frakS_{n+1}}\ra \omega^1_{\frakX_1/\frakS_1}\ra 0.$$
Hence the morphism of multiplication by $p;$ $\omega^1_{\frakX_n/\frakS_n}\xrightarrow{\times p}\omega^1_{\frakX_{n+1}/\frakS_{n+1}}$ is injective. We deduce the existence of a unique $\Ox_{\frakX_n}$-linear morphism $p^{-1}dF_{n+1}:F_n^*\omega^1_{\frakX_n'/\frakS_n}\ra \omega^1_{\frakX_n/\frakS_n}$ making the following diagram commutative
\begin{equation}\label{surp}
\begin{tikzcd}
F_{n+1}^*\omega^1_{\frakX_{n+1}'/\frakS_{n+1}}\ar{d}\ar{r}{dF_{n+1}} & p\omega^1_{\frakX_{n+1}/\frakS_{n+1}} \\
F_n^*\omega^1_{\frakX'_n/\frakS_n}\ar{r}{p^{-1}dF_{n+1}} & \omega^1_{\frakX_n/\frakS_n}\ar{u}{\times p}
\end{tikzcd}
\end{equation}
Let us give an explicit expression of $p^{-1}dF_{n+1}$ in light of \ref{lem12}.
Let $U\ra \frakX_{n+1}$ be an étale morphism, $\tilde{m}\in \Gamma(U,\calM_{\frakX_{n+1}}),$ $m$ and $\ov{m}$ the images of $\tilde{m}$ in $\Gamma(U\times_{\frakX_{n+1}}\frakX_n,\calM_{\frakX_n})$ and $\Gamma(U\times_{\frakX_{n+1}}X,\calM_X)$ respectively. Let $m'$ and $\tilde{m}'$ be compatible local liftings of $\pi^{\flat}\ov{m}$ to $\calM_{\frakX_n'}$ and $\calM_{\frakX_{n+1}'}$ respectively (\ref{prop69}). Set $\tilde{u}=F_{n+1}^{\flat}(\tilde{m}')-p\tilde{m}.$ By \ref{lem12}, $\tilde{u}$ is invertible in $\calM_{\frakX_{n+1}}$ and there exists a local section $\tilde{b}$ of $\Ox_{\frakX_{n+1}}$ such that $\alpha_{\frakX_{n+1}}(\tilde{u})=1+p\tilde{b}.$ Let $b$ (resp. $u$) be the image of $\tilde{b}$ (resp. $\tilde{u}$) in $\Ox_{\frakX_n}$ (resp. $\calM_{\frakX_n}$). 
By the commutativity of (\ref{surp}), we deduce that in $\omega^1_{\frakX_{n+1}/\frakS_{n+1}}:$
\begin{alignat*}{2}
p\times (p^{-1}dF_{n+1})(F_n^*\op{dlog}m') &= dF_{n+1}(F_{n+1}^*\op{dlog}\tilde{m}') \\
&=  \op{dlog}(F_{n+1}^{\flat}\tilde{m}') \\
&= \op{dlog}(p\tilde{m}+\tilde{u}) \\
&= p\op{dlog}\tilde{m}+\op{dlog}\tilde{u}\\
&= p\op{dlog}\tilde{m}+\frac{d\alpha_{\frakX_{n+1}}(\tilde{u})}{\alpha_{\frakX_{n+1}}(\tilde{u})}\\
&= p\op{dlog}\tilde{m}+p\frac{d\tilde{b}}{1+p\tilde{b}}. 
\end{alignat*}
By the injectivity of $\omega^1_{\frakX_n/\frakS_n}\xrightarrow{\times p}\omega^1_{\frakX_{n+1}/\frakS_{n+1}},$ we conclude that, in $\omega^1_{\frakX_n/\frakS_n},$ we have:
\begin{equation}\label{surp2}
(p^{-1}dF_{n+1})(F_n^*\op{dlog}m') = \op{dlog}m+\frac{db}{1+pb}= \op{dlog}m+(1-pb+p^2b^2-\hdots )db.
\end{equation}
\end{parag}

\begin{lemma}\label{lem85}
Let $n\ge 1$ be an integer and suppose that $\frakX_n \ra \frakS_n$ is log smooth. Let $(\calE',\nabla')$ be an $\Ox_{\frakX_n'}$-module equipped with a $p$-connection and $\zeta$ the composition
$$\zeta:F_n^*(\calE'\otimes_{\Ox_{\frakX_n'}}\omega^1_{\frakX_n'/\frakS_n})\xrightarrow{\sim} F_n^*(\calE')\otimes_{\Ox_{\frakX_n}}F_n^*(\omega^1_{\frakX_n'/\frakS_n})\xrightarrow{\op{Id}_{F_n^*(\calE')}\otimes p^{-1}dF_{n+1}}F_n^*(\calE')\otimes_{\Ox_{\frakX_n}}\omega^1_{\frakX_n/\frakS_n},$$
where the first arrow is the canonical isomorphism.
Consider the morphism 
\begin{equation}\label{conn}
\nabla:\begin{array}[t]{lc}F_n^*\calE'\ra F_n^*\calE'\otimes_{\Ox_{\frakX_n}}\omega^1_{\frakX_n/\frakS_n}\\ x\otimes a\mapsto a\zeta( F_n^*\nabla'(x))+F_n^*x\otimes da,\end{array}
\end{equation}
where $x$ and $a$ are local sections $\calE'$ and $\Ox_{\frakX_n}$ respectively.
Then $\nabla$ a well defined connection on $F_n^*\calE'.$ In addition, if $\nabla'$ is integrable then $\nabla$ is also integrable.
\end{lemma}

\begin{proof}
It is clear that for any local sections $x$ and $y$ of $\calE',$ $a$ and $b$ of $\Ox_{\frakX_n}$ and $\alpha$ of $\Ox_{\frakX_n'},$
$\nabla((x+y)\otimes a)=\nabla(x\otimes a)+\nabla(y\otimes a)$ and 
$\nabla(x\otimes (a+b))=\nabla(x\otimes a)+\nabla(x\otimes b).$
In addition
\begin{equation}
\nabla((\alpha x)\otimes b)=F_n^{\#}(\alpha) b\zeta(F_n^*\nabla'(x))+F_n^*x\otimes bdF_n^{\#}(\alpha)+F_n^*(\alpha x)\otimes db
\end{equation}
and
\begin{alignat}{2}
\nabla(x\otimes(F_n^{\#}(\alpha) b)) &= F_n^{\#}(\alpha) b\zeta(F_n^*\nabla'(x))+F_n^*x\otimes d(F_n^{\#}(\alpha) b)\\
&=\nabla((\alpha x)\otimes b).
\end{alignat}
It follows that $\nabla$ is well defined and by (\ref{conn}), it is a connection on the $\Ox_{\frakX_n}$-module $F_n^*\calE'.$ 

Now suppose that $\nabla'$ is integrable and let us check that $\nabla$ is also integrable. Let $x$ be a local section of $\calE'.$ Let $(\ov{m}_i)_{1\le i\le d}$ be local coordinates for $X\ra S$ (\ref{P2}). For $1\le i\le d,$ let $m_i$ (resp. $m_i'$) be a local lifting of $\ov{m}_i$ (resp. $\pi^{\flat}\ov{m}_i$) to $\frakX_n$ (resp. $\frakX_n'$). Then $(\op{dlog}m_i')_{1\le i\le d}$ generates the $\Ox_{\frakX_n'}$-module $\omega^1_{\frakX_n'/\frakS_n}.$ Let $x_i,\ x_{ij},\ i,j=1,\hdots,d$ be local sections of $\calE'$ such that $$\nabla'(x)=\sum_{i=1}^dx_i\otimes \op{dlog}m_i'$$
and 
$$\nabla'(x_i)=\sum_{j=1}^dx_{ij}\otimes \op{dlog}m_j'.$$ It follows that the curvature $K(\nabla')$ of $\nabla'$ is given by $$K(\nabla')(x)=\sum_{i,j}x_{ij}\otimes \op{dlog}m_i'\wedge \op{dlog}m_j'.$$ On the other hand, by the definition \ref{conn}, we have
$$\nabla(aF_n^*x)=\sum_{i=1}^daF_n^*x_i\otimes p^{-1}F_{n+1}(\op{dlog}m_i')+F_n^*x\otimes da,$$
for any local sections $a$ and $x$ respectively of $\Ox_{\frakX_n}$ and $\calE'.$ It follows that the curvature $K(\nabla)$ of $\nabla$ is given, for any local sections $a$ and $x$ respectively of $\Ox_{\frakX_n}$ and $\calE',$ by
\begin{alignat*}{2}
K(\nabla)(aF_n^*x)&=\nabla^1\circ \nabla(aF_n^*x)\\
&=\sum_{i=1}^d\nabla^1(aF_n^*x_i\otimes (p^{-1}F_{n+1})(\op{dlog}m_i'))+\nabla^1(F_n^*x\otimes da)\\
&=\sum_{i=1}^d\nabla(aF_n^*x_i)\wedge (p^{-1}F_{n+1})(\op{dlog}m_i')+aF_n^*x_i\otimes d((p^{-1}F_{n+1})(\op{dlog}m_i'))\\ 
&\ +\nabla(F_n^*x)\wedge da.
\end{alignat*}
Finally, we get:
\begin{alignat}{2}
K(\nabla)(aF_n^*x)&=\sum_{i,j}aF_n^*x_{ij}\otimes (p^{-1}dF_{n+1})(\op{dlog}m_j')\wedge (p^{-1}dF_{n+1})(\op{dlog}m_i') \label{856} \\
&\ +\sum_iF_n^*x_i\otimes da\wedge (p^{-1}dF_{n+1})(\op{dlog}m_i') \label{857}\\
&\ +\sum_iaF_n^*x_i\otimes d((p^{-1}dF_{n+1})(\op{dlog}m_i')) \label{858}\\
&\ +\nabla(F_n^*x)\wedge da. \label{859}
\end{alignat} 
The first sum is simply the image of $K(\nabla')(x)$ by the morphism $F_n^*\otimes p^{-1}dF_{n+1}\wedge p^{-1}dF_{n+1}$ and is thus equal to zero since $\nabla'$ is integrable. We also have
\begin{equation*}
\nabla(F_n^*x)\wedge da = \sum_iF_n^*x_i\otimes (p^{-1}dF_{n+1})(\op{dlog}m'_i)\wedge da
\end{equation*}
so (\ref{857}) and (\ref{859}) cancel out. It remains to prove that (\ref{858}) vanishes. By (\ref{surp2}), there exists a local section $b_i$ of $\Ox_{\frakX_n}$ such that 
$$(p^{-1}dF_{n+1})(F_{n+1}^*\op{dlog}m_i')=\op{dlog}m_i+(1-pb_i+p^2b_i^2-\hdots )db_i.$$
It follows that $(p^{-1}dF_{n+1})(\op{dlog}m_i')$ is a closed form and so (\ref{858}) vanishes.
We conclude that $K(\nabla)(x)=0$ and so $\nabla$ is integrable. 
\end{proof}

\begin{parag}
Keep the same hypothesis of \ref{lem85}. By \ref{lem85}, we have a functor
\begin{equation}\label{krazphin}
\Phi_n:\begin{array}[t]{clc}
p\op{-MIC}(\frakX_n'/\frakS_n) & \ra & \op{MIC}(\frakX_n/\frakS_n)\\ 
(\calE',\nabla') & \mapsto & (F_n^*\calE',\nabla)
\end{array}
\end{equation}
from the category $p\op{-MIC}(\frakX_n'/\frakS_n)$ of $\Ox_{\frakX_n'}$-modules with an integrable $p$-connection to the category $\op{MIC}(\frakX_n/\frakS_n)$ of $\Ox_{\frakX_n}$-modules with an integrable connection.
\end{parag}

\newpage

\section{Dilatations}

In this section, we fix a perfect field $k$ of characteristic $p$ and denote by $W(k)$ its ring of Witt vectors. Let $\frakS=\op{Spf}W(k)$ equipped with the trivial logarithmic structure. For any logarithmic $p$-adic formal scheme $\frakY$ and any positive integer $n,$ we denote by $\frakY_n$ the logarithmic scheme obtained from $\frakY$ by reduction modulo $p^n$ (as in \ref{prop69}).

\begin{parag}
Let $n$ be a positive integer, $\frakX$ a logarithmic $p$-adic formal scheme flat and locally of finite type over $\frakS$ (\ref{dxuflat}) and $\calI$ an open ideal of finite type of $\Ox_{\frakX}$ containing $p^n$ (\cite{Ahmed2010} 2.1.19). Let $\frakY \ra \frakX$ be the admissible blow-up of $\calI$ in $\frakX$ (\cite{Ahmed2010} 3.1.2). We equip it with the logarithmic structure pullback of that of $\frakX.$ By (\cite{Ahmed2010} 3.1.4), the ideal $\calI\Ox_{\frakY}$ is invertible. We denote by $\frakX_{(\calI/p^n)}$ the dilatation of $\calI$ with respect to $p^n$ i.e. the largest open formal subscheme of $\frakY$ such that the restriction of $\calI \Ox_{\frakY}$ on $\frakX_{(\calI/p^n)}$ is generated by $p^n$ (\cite{Ahmed2010} 3.2.3.4 and 3.2.7). The formal scheme $\frakX_{(\calI/p^n)}$ is flat over $\frakS$ (\ref{dxuflat}).
We denote by $\calI^{\#}$ the ideal generated by $p$ and the $p$th powers of the sections of $\calI$ and by $\frakX^{\#}_{(\calI/p)}$ the dilatation of $\calI^{\#}$ with respect to $p.$ 
If $\frakX=\op{Spf}A$ and $I$ is the open ideal of $A$ corresponding to $\calI,$ and $(a_0,\hdots,a_r)$ is a set of generators of $I$ such that $a_0=p^n,$ then $\frakX_{(\calI/p^n)}=\op{Spf}B$ where $B$ is the $p$-adic completion of the ring $A_0/(p^n\text{-tor}),$ where
\begin{equation}\label{dilaff}
A_0=A\left [x_1,\hdots,x_r\right ]/(p^nx_1-a_1,\hdots,p^nx_r-a_r)
\end{equation}
and $(p^n\text{-tor})$ is the $p^n$-torsion ideal of $A_0.$ In particular, we see that the canonical morphism $\frakX_{(\calI/p^n)} \ra \frakX$ is affine.
\end{parag}

\begin{proposition}\label{erafdil1}
Let $n$ be a positive integer, $\frakY$ a logarithmic $p$-adic formal scheme flat and locally of finite type over $\frakS$ and $i:T \ra \frakY_n$ a strict immersion. There exists a strict morphism of logarithmic $p$-adic formal schemes flat over $\frakS,$ $g:\frakY_{(T/p^n)} \ra \frakY,$ unique up to a canonical isomorphism, satisfying the following conditions:
\begin{enumerate}
\item The morphism $\left (\frakY_{(T/p^n)}\right )_n \ra \frakY_n$ factors uniquely through $T \ra \frakY_n$ and the resulting morphism $\left (\frakY_{(T/p^n)}\right )_n \ra T$ is affine.
\item If $\frakZ$ is a logarithmic $p$-adic formal scheme flat over $\frakS$ and $f:\frakZ \ra \frakY$ is an $\frakS$-morphism such that $\frakZ_n \ra \frakY_n$ factors through $T\ra \frakY_n$ then there exists a unique $\frakS$-morphism $f':\frakZ \ra \frakY_{(T/p^n)}$ such that $f=g\circ f'.$ In addition, if $T \ra \frakY$ and $f$ are closed immersions, then so is $f'.$
\end{enumerate}
\end{proposition}

\begin{proof}
If the immersion $i:T \ra \frakY_n$ is closed, we denote by $\calI$ the ideal of $T \ra \frakY.$ We take $\frakY_{(T/p^n)}$ to be the dilatation of $\calI$ with respect to $p^n,$ we equip it with the logarithmic structure pull back of that of $\frakY$ and take $g:\frakY_{(T/p^n)} \ra \frakY$ to be the canonical morphism. By \ref{era3proplogstr}, to prove that the morphism of logarithmic formal schemes $\left (\frakY_{(T/p^n)} \right )_n \ra \frakY_n$ factors through $i,$ it is sufficient to prove the factorization for the underlying formal schemes, which is clear. For the second assertion, again by \ref{era3proplogstr}, it is sufficient to prove the existence of $f'$ on the underlying formal schemes. The rest of the proof is the same as in (\cite{DXU19} 3.5).
\end{proof}

\begin{proposition}\label{erafdil2}
Let $\frakY$ be a logarithmic $p$-adic formal scheme flat and locally of finite type over $\frakS$ and $i:T \ra \frakY_1$ a strict immersion. There exists a strict morphism of logarithmic $p$-adic formal schemes $g:\frakY_{(T/p)}^{\#} \ra \frakY,$ unique up to a canonical isomorphism, satisfying the following conditions:
\begin{enumerate}
\item The morphism $\underline{\left (\frakY_{(T/p)}^{\#}\right )_1}\ra \left (\frakY_{(T/p)}^{\#}\right )_1 \ra \frakY_1$ factors uniquely through $T \ra \frakY_1$ and the resulting morphism $\underline{\left (\frakY_{(T/p)}^{\#}\right )_1} \ra T$ is affine.
\item If $\frakZ$ is a logarithmic $p$-adic formal scheme flat over $\frakS$ and $f:\frakZ \ra \frakY$ is an $\frakS$-morphism such that $\underline{\frakZ_1} \ra \frakY_1$ factors through $T\ra \frakY_1$ then there exists a unique $\frakS$-morphism $f':\frakZ \ra \frakY_{(T/p)}^{\#}$ such that $f=g\circ f'.$
\end{enumerate}
\end{proposition}

\begin{proof}
The proof is similar to \ref{erafdil1}.
\end{proof}

\begin{proposition}\label{erafdil3}
Let $n$ be a positive integer, $\frakY$ a logarithmic $p$-adic formal scheme flat and locally of finite type over $\frakS$ and $i:T \ra \frakY_n$ a strict immersion. For any integer $k\ge n,$ there exists a canonical morphism
$$\frakY_{(T/p^{k+1})} \ra \frakY_{(T/p^{k})}.$$
\end{proposition}

\begin{proof}
We have the commutative diagram on the left:
$$
\begin{tikzcd}
\frakY_{(T/p^{k+1}),k} \ar{r} \ar{d} & \frakY_k \ar{d} \\
\frakY_{(T/p^{k+1}),k+1} \ar{r} \ar{d} & \frakY_{k+1} \\
T \ar{ur} & 
\end{tikzcd}
\begin{tikzcd}
\frakY_{(T/p^{k+1}),k} \ar{r} \ar{d} & \frakY_k \\
\frakY_{(T/p^{k+1}),k+1} \ar{d} &  \\
T \ar{uur} & 
\end{tikzcd}
$$
Since $T \ra \frakY_{k+1}$ decomposes as $T \ra \frakY_k \ra \frakY_{k+1}$ and $\frakY_k \ra \frakY_{k+1}$ is an immersion, hence a monomorphism, we obtain the commutative diagram on the right.
The desired morphism is then obtained by the flatness of $\frakY_{(T/p^{k+1})}$ over $\frakS$ and by applying the universal property of $\frakY_{(T/p^k)}$ to the canonical morphism
$$\frakY_{(T/p^{k+1})} \ra \frakY.$$
\end{proof}

\begin{proposition}[\cite{Oyama} 1.1.4 and \cite{DXU19} 3.7]\label{propdiletale}
Let $n$ be a positive integer, $\frakX$ and $\frakY$ two logarithmic $p$-adic formal schemes flat and locally of finite type over $\frakS$ and $f:\frakX \ra \frakY$ a log étale $\frakS$-morphism. Suppose that there exists two strict $\frakS_n$-immersions $i:T \ra \frakX_n$ and $j:T \ra \frakY_n$ such that $f_n\circ i=j.$ Then $f$ induces a canonical isomorphism of logarithmic formal schemes
$$\frakX_{(T/p^n)} \xrightarrow{\sim} \frakY_{(T/p^n)}.$$
If $n=1,$ then $f$ also induces a canonical isomorphism of logarithmic formal schemes
$$\frakX_{(T/p)}^{\#} \xrightarrow{\sim} \frakY_{(T/p)}^{\#}.$$
\end{proposition}

\begin{proof}
Denote by $g_{\frakX}:\frakX_{(T/p^n)} \ra \frakX$ the canonical morphism. Since the formal scheme $\frakX_{(T/p^n)}$ is flat over $\frakS,$
the universal property of $\frakY_{(T/p^n)}$ applied to the morphism
$$\frakX_{(T/p^n)} \xrightarrow{g_{\frakX}} \frakX \xrightarrow{f} \frakY$$
implies the existence of a unique morphism $f':\frakX_{(T/p^n)} \ra \frakY_{(T/p^n)}$ such that the diagram
$$
\begin{tikzcd}
\frakX_{(T/p^n)} \ar{r}{g_{\frakX}} \ar{d}{f'} & \frakX \ar{r}{f} & \frakY \\
\frakY_{(T/p^n)} \ar[swap]{urr}{g_{\frakY}} & & 
\end{tikzcd}
$$
is commutative.
The log étaleness of $f:\frakX \ra \frakY$ implies the existence of a unique morphism $h:\frakY_{(T/p^n)} \ra \frakX$ fitting into the commutative diagram
$$
\begin{tikzcd}
T \ar{r}{i} & \frakX_n \ar{r} & \frakX \ar{d}{f} \\
\left (\frakY_{(T/p^n)} \right )_n \ar{r} \ar{u} & \frakY_{(T/p^n)} \ar{r} \ar{ur}{h} & \frakY
\end{tikzcd}
$$
The universal property of $\frakX_{(T/p^n)}$ applied to $h:\frakY_{(T/p^n)} \ra \frakX$ implies the existence of a unique morphism $f'':\frakY_{(T/p^n)} \ra \frakX_{(T/p^n)}$ such that the diagram
$$
\begin{tikzcd}
\frakY_{(T/p^n)} \ar{r}{h} \ar{d}{f''} & \frakX \\
\frakX_{(T/p^n)} \ar[swap]{ur}{g_{\frakX}} & 
\end{tikzcd}
$$
is commutative. Consider the commutative diagram
$$
\begin{tikzcd}
\frakY_{(T/p^n)} \ar{r} \ar{d}{f''} & \frakX \ar{r}{f} & \frakY \\
\frakX_{(T/p^n)} \ar{d}{f'} \ar{ur} & & \\
\frakY_{(T/p^n)} \ar{uurr} & & 
\end{tikzcd}
$$
The universal property of $\frakY_{(T/p^n)}$ implies that
$$f' \circ f'' = \op{Id}_{\frakY_{(T/p^n)}}.$$
Now let $u=f''\circ f'.$ The commutative diagram
$$
\begin{tikzcd}
\frakX_{(T/p^n),n} \ar{dr} \ar{r}{f'_n} & \frakY_{(T/p^n),n} \ar{d} \ar{r}{h_n} & \frakX_n \\
 & T\ar[hook]{ur}{i} & 
\end{tikzcd}
$$
implies that $h_n \circ f'_n=g_{\frakX,n}.$
We then have the following commutative diagram:
$$
\begin{tikzcd}
 & \frakX_n \ar[hook]{rr} & & \frakX \ar{d}{f} \\
\left ( \frakX_{(T/p^n)} \right )_n \ar{ur}{g_{\frakX,n}} \ar{r} & \frakX_{(T/p^n)} \ar{urr}{h\circ f'} \ar{r} & \frakX \ar{r}{f} & \frakY
\end{tikzcd}
$$
By the log étaleness of $f,$ we deduce that $h \circ f'=g_{\frakX}.$ Since $h=g_{\frakX}\circ f'',$ we get
$$g_{\frakX}\circ  f''\circ f'=g_{\frakX}.$$
We deduce by the universal property of $\frakX_{(T/p^n)}$ that
$$f'' \circ f'=\op{Id}_{\frakX_{(T/p^n)}}.$$
If $n=1,$ the isomorphism $\frakX_{(T/p)}^{\#} \xrightarrow{\sim} \frakY_{(T/p)}^{\#}$ is established by a similar argument.
\end{proof}

\begin{proposition}[\cite{Oyama} 1.1.5 and \cite{DXU19} 3.8]\label{propdilflat}
Let $n$ be a positive integer, $\frakX$ and $\frakY$ two logarithmic $p$-adic formal schemes flat and locally of finite type over $\frakS,$ $f:\frakX \ra \frakY$ an $\frakS$-morphism which is flat on the underlying formal schemes, $T \ra \frakY_n$ a strict immersion and $S=\frakX\times_{\frakY}T.$ Then $f$ induces a canonical isomorphism
$$\frakY_{(T/p^n)}\times_{\frakY}\frakX \xrightarrow{\sim} \frakX_{(S/p^n)}.$$
If $n=1,$ then $f$ also induces a canonical isomorphism
$$\frakY_{(T/p)}^{\#}\times_{\frakY}\frakX \xrightarrow{\sim} \frakX_{(S/p)}^{\#}.$$
\end{proposition}

\begin{proof}
Denote by $g_{\frakX}:\frakX_{(S/p^n)} \ra \frakX$ and $g_{\frakY}:\frakY_{(T/p^n)} \ra \frakY$ the canonical morphisms.
The morphism $f:\frakX \ra \frakY$ is flat on the underlying formal schemes, hence so is the projection $\frakY_{(T/p^n)}\times_{\frakY}\frakX \ra \frakY_{(T/p^n)}.$ We deduce that $\frakY_{(T/p^n)}\times_{\frakY}\frakX$ is flat over $\frakS.$
By the commutativity of the diagram
$$
\begin{tikzcd}
\left (\frakX_{(S/p^n)} \right )_n \ar{r} \ar{d} & \frakX_n \ar{r} & \frakY_n \\
S \ar{ur} \ar{r} & T\ar{ur} &
\end{tikzcd}
$$
and the universal property of $\frakY_{(T/p^n)},$ there exists a unique morphism $\frakX_{(S/p^n)} \ra \frakY_{(T/p^n)}$ fitting into the commutative diagram
$$
\begin{tikzcd}
\frakX_{(S/p^n)} \ar{r} \ar[swap]{d}{g_{\frakX}} & \frakY_{(T/p^n)} \ar{d}{g_{\frakY}} \\
\frakX \ar{r}{f} & \frakY
\end{tikzcd}
$$
It follows that there exists a unique morphism $f':\frakX_{(S/p^n)} \ra \frakY_{(T/p^n)}\times_{\frakY}\frakX$ making the following diagram commutative:
$$
\begin{tikzcd}
\frakX_{(S/p^n)} \ar{dr}{f'} \ar[bend right=-20]{drr} \ar[swap, bend right=30]{ddr}{g_{\frakX}}  & & \\
 & \frakY_{(T/p^n)}\times_{\frakY}\frakX \ar{r}{\alpha} \ar{d}{\beta} & \frakY_{(T/p^n)} \ar{d}{g_{\frakY}} \\
 & \frakX \ar{r}{f} & \frakY,
\end{tikzcd}
$$
where $\alpha$ and $\beta$ are the canonical projection.
Applying the universal property of $\frakX_{(S/p^n)}$ to the projection
$$\beta:\frakY_{(T/p^n)}\times_{\frakY}\frakX \ra \frakX$$
proves the existence of a unique morphism
$$f'':\frakY_{(T/p^n)}\times_{\frakY}\frakX \ra \frakX_{(S/p^n)}$$
such that the following diagram is commutative
$$
\begin{tikzcd}
\frakY_{(T/p^n)}\times_{\frakY}\frakX \ar{r}{\beta} \ar{d}{f''} & \frakX \\
\frakX_{(S/p^n)} \ar[swap]{ur}{g_{\frakX}} & 
\end{tikzcd}
$$
The commutativity of the diagram
$$
\begin{tikzcd}
\frakX_{(S/p^n)} \ar{r}{f'} \ar{d}{g_{\frakX}} & \frakY_{(T/p^n)} \times_{\frakY} \frakX \ar{r}{f''} & \frakX_{(S/p^n)} \ar{dll}{g_{\frakX}} \\
\frakX & & 
\end{tikzcd}
$$
and the universal property of $\frakX_{(S/p^n)}$ proves that
$$f'' \circ f' =\op{Id}_{\frakX_{(S/p^n)}}.$$
The commutativity of the diagram
$$
\begin{tikzcd}
\frakY_{(T/p^n)}\ar{dr} \times_{\frakY} \frakX \ar{r}{f''} \ar[swap]{dd}{\alpha} & \frakX_{(S/p^n)} \ar{d} \ar[bend right=-30]{rr} \ar{r}{f'} & \frakY_{(T/p^n)} \times_{\frakY} \frakX \ar{r}{\alpha} & \frakY_{(T/p^n)} \ar{dd}{g_{\frakY}} \\
& \frakX \ar{drr} & & \\
\frakY_{(T/p^n)} \ar{rrr}{g_{\frakY}} & & & \frakY
\end{tikzcd}
$$
and the universal property of $\frakY_{(T/p^n)}$ imply that
$$\alpha \circ f' \circ f''=\alpha.$$
We also have
$$\beta \circ f' \circ f''=g_{\frakX} \circ f''=\beta.$$
This proves that
$$f' \circ f'' =\op{Id}_{\frakY_{(T/p^n)} \times_{\frakY} \frakX}.$$
If $n=1,$ the isomorphism $\frakY_{(T/p)}^{\#}\times_{\frakY}\frakX \xrightarrow{\sim} \frakX_{(S/p)}^{\#}$ is established by a similar argument.
\end{proof}

\newpage

\section{Groupoids}

\begin{parag}\label{paragJapon1}
In this section, we fix a perfect field $\kappa$ of characteristic $p$ and denote by $W$ its ring of Witt vectors. We equip $\op{Spf}W$ with the trivial logarithmic structure. For a logarithmic $p$-adic formal scheme $\frakX$ over $\op{Spf}W$ and a positive integer $n,$ we denote by $\frakX_n$ the logarithmic scheme obtained from $\frakX$ by reduction modulo $p^n.$ We use the corresponding roman letter $X$ for $\frakX_1.$ Unless otherwise stated, all formal $\op{Spf}W$-schemes considered are $p$-adic and all logarithmic structures are assumed to be fs.
\end{parag}

\begin{parag}\label{parag86}
In this section, let $P \ra Q$ be a morphism of fs monoids (\cite{Ogus2018} I 4.6.2), $f:(\frakX,Q)\ra (\frakS,P)$ a log smooth morphism of framed logarithmic formal schemes \eqref{def72}, log flat and locally of finite type over $\op{Spf}W.$ Note that this implies, by \ref{Wflat}, that $\frakX$ and $\frakS$ are flat over $\op{Spf}W$ \eqref{dxuflat}. Let $r$ be a nonnegative integer. By \ref{cor79}, we have a factorization
$$\begin{tikzcd}
 & \frakX_{\frakS,\left [Q\right ]}^{r+1}\ar{d} \\
\frakX \ar{r} \ar{ur}{\Delta(r)}&  \frakX^{r+1}_{\frakS}
\end{tikzcd}$$
where $\frakX^{r+1}_{\frakS}$ denotes the fiber product, in the category fs logarithmic $p$-adic formal schemes, of $\frakX$ over $\frakS$ with itself $r+1$ times, $\Delta(r):\frakX\ra \frakX_{\frakS,[Q]}^{r+1}$ is a strict immersion and $\frakX_{\frakS,[Q]}^{r+1}\ra \frakX_{\frakS}^{r+1}$ is log étale. Set $\frakY(r)=\frakX_{\frakS,[Q]}^{r+1}.$ The canonical morphism $\frakY(r) \ra \frakS$ is log smooth hence locally of finite type, so $\frakY(r)$ is locally of finite type over $\op{Spf}W.$ The projections $\frakY(r) \ra \frakX$ are strict (\ref{cor76}) and log smooth hence smooth. It follows that $\frakY(r)$ is flat over $\op{Spf}W.$ 

Let $\calI_{\frakX/\frakS}(r)$ be the ideal of the immersion $\Delta(r):\frakX \ra \frakY(r).$ For any integer $n\ge 1,$ denote by $R_{\frakX,n}(r)$ (resp. $Q_{\frakX}(r)$) the dilatation $\frakY(r)_{(\frakX_n/p^n)}$ (resp. $\frakY(r)_{(X/p)}^{\#}$) defined in \ref{erafdil1} (resp. \ref{erafdil2}). By the universal property, there exists a unique strict immersion of logarithmic formal schemes $\frakX\ra R_{\frakX,n}(r)$ (resp. $\frakX\ra Q_{\frakX}(r)$) making the following diagrams commutative
$$\begin{tikzcd}
 & R_{\frakX,n}(r)\ar{d} \\
 & \frakY(r)\ar{d} \\
\frakX\ar{uur}\ar{ur}\ar{r} & \frakX^{r+1}_{\frakS}
\end{tikzcd} 
\begin{tikzcd}
 & Q_{\frakX}(r)\ar{d} \\
 & \frakY(r)\ar{d} \\
\frakX\ar{uur}\ar{ur}\ar{r} & \frakX^{r+1}_{\frakS}
\end{tikzcd}
$$
We extend the notation to $n=0$ by setting $R_{\frakX,0}(r)=\frakY(r).$ For every positive integer $k,$ we denote by $\frakX_k$ and $R_{\frakX,n}(r)_k$ the logarithmic schemes obtained from $\frakX$ and $R_{\frakX,n}(r)$ respectively, by reduction modulo $p^k$ (as in \ref{prop69}), by $P_{\frakX/\frakS,n}(r)_k$ the PD-envelope of the strict immersion $\frakX_k\ra R_{\frakX,n}(r)_k$ and by $\ov{\calI}_{\frakX/\frakS,n,k}(r)$ the PD-ideal of $P_{\frakX/\frakS,n}(r)_k.$ We equip $P_{\frakX/\frakS,n}(r)_k$ with the logarithmic structure pullback of that of $R_{\frakX,n}(r)_k.$ By extension of scalars (\cite{Ber74} I 2.8.2 and \cite{Ogus78} 3.20.8), the canonical morphism
$$
P_{\frakX/\frakS,n}(r)_k \ra P_{\frakX/\frakS,n}(r)_{k+1}\times_{\frakS_{k+1}}\frakS_k
$$
is an isomorphism for any positive integer $k.$
The schemes $(P_{\frakX/\frakS,n}(r)_k)_{k\ge 1}$ form then a $p$-adic system of schemes. We denote by $P_{\frakX/\frakS,n}(r)$ its limit, equip it with the logarithmic structure pullback of that of $R_{\frakX,n}(r).$ We obtain the commutative diagram:
$$
\begin{tikzcd}
 & P_{\frakX/\frakS,n}(r)\ar{d} \\
 & R_{\frakX,n}(r)\ar{d} \\
 & \frakY(r)\ar{d} \\
\frakX \ar[bend right =-20]{uuur} \ar{uur}\ar{ur}\ar{r} & \frakX^{r+1}_{\frakS}
\end{tikzcd}
$$
The canonical morphism $X\ra P_{\frakX/\frakS,n}(r)_1,$ is a universal homeomorphism on the underlying topological spaces and hence so is $\frakX\ra P_{\frakX/\frakS,n}(r).$ We identify the underlying small étale sites via this homeomorphism. We denote by $\calP_{\frakX/\frakS,n}(r)$ the structural ring of $P_{\frakX/\frakS,n}(r)$ and, for any integers $k\ge 1$ and $l\ge 0,$ $\calP_{\frakX/\frakS,n}^{\{ l \}}(r)_k=\calP_{\frakX/\frakS,n}(r)_k/\ov{\calI}_{\frakX/\frakS,n,k}(r)^{[l+1]}$ and $\ov{\calI}^{\{ l \}}_{\frakX/\frakS,n,k}(r)=\ov{\calI}_{\frakX/\frakS,n,k}(r)/\ov{\calI}_{\frakX/\frakS,n,k}(r)^{[l+1]}.$

From now on, for $r=1,$ we will drop $(r)$ from the notation we just introduced. Note that $\calP_{\frakX/\frakS,n,k}^{\{ 0 \}}=\Ox_{\frakX_k}$ for every positive integer $k.$

By Kato's construction of the PD envelope, given in (\cite{Kat89} 5.6), and in view of \ref{prop45}, for all integers $k,r\ge 1,$ Kato's logarithmic PD-envelope of $\frakX_k \ra (\frakX_k)^{r+1}$ and the logarithmic scheme $P_{\frakX/\frakS,0}(r)_k$ are the same.
It follows, by (\cite{Kat89} 5.8.1), that, for any integer $k\ge 1,$ there exists a canonical isomorphism
\begin{equation}\label{eqtakrizIIsquare}
\ov{\calI}^{\{1\}}_{\frakX/\frakS,0,k} \xrightarrow{\sim} \omega^1_{\frakX_k/\frakS_k}.
\end{equation}
\end{parag}

\begin{parag}
Let $k$ be a positive integer. Since $f:\frakX \ra \frakS$ is log smooth and $\frakY=\frakX\times_{\frakS,[Q]}^{\op{log}}\frakX \ra \frakX\times_{\frakS}^{\op{log}}\frakX$ is log étale (\ref{prop74}), the projections $\frakY \ra \frakX$ are log smooth. Since they are also strict (\ref{cor76}), they are smooth. It follows that the exact diagonal immersion $\frakX_k \ra \frakY_k$ is regular. Let $\calI_k$ be its ideal.
The canonical morphism of $\Ox_{\frakX_k}$-modules
\begin{equation}\label{grsym}
S^n\left (\calI_k/\calI_k^2\right ) \xrightarrow{\sim} \calI_k^n/\calI_k^{n+1}
\end{equation}
is then an isomorphism for all $n\ge 1.$
\end{parag}

\begin{proposition}\label{QU}
Let $\frakU \ra \frakX$ be a strict étale morphism of formal logarithmic $\frakS$-schemes. Equip $\frakU$ with the frame $\frakU \ra \frakX \ra [Q]$ and consider the dilatation $Q_{\frakU}=\left (\frakU\times_{\frakS,[Q]}^{\op{log}}\frakU \right )_{(U/p)}^{\#}$ \eqref{erafdil2}. Then, we have a canonical isomorphism
$$
Q_{\frakU} \xrightarrow{\sim} Q_{\frakX} \times_{\frakX}\frakU,
$$
where $Q_{\frakX}$ is considered over $\frakX$ by the second projection.
\end{proposition}

\begin{proof}
We have the commutative diagram
$$
\begin{tikzcd}
\frakU \ar{d} \ar[equal]{r} & \frakU \ar{r} \ar{d} &  \frakX \ar{d}  \\
\frakU\times_{\frakS,[Q]}^{\op{log}}\frakU  \ar{r} & \frakX\times_{\frakS,[Q]}^{\op{log}}\frakU \ar{r} & \frakX\times_{\frakS,[Q]}^{\op{log}}\frakX,
\end{tikzcd}
$$
where the lower morphisms and the middle vertical arrow are induced by the strict étale morphism $\frakU \ra \frakX$ and the first and third vertical arrows are the strict diagonal immersions. By \ref{propdiletale} and \ref{propdilflat}, we get canonical isomorphisms
$$
Q_{\frakU} \xrightarrow{\sim} \left (\frakX\times_{\frakS,[Q]}^{\op{log}}\frakU \right )_{(U/p)}^{\#} \xrightarrow{\sim} \left (\frakX\times_{\frakS,[Q]}^{\op{log}}\frakX \right )_{(X/p)}^{\#} \times_{\frakX\times_{\frakS,[Q]}^{\op{log}} \frakX}\left (\frakX\times_{\frakS,[Q]}^{\op{log}}\frakU \right ) \xrightarrow{\sim} Q_{\frakX} \times_{\frakX}\frakU.
$$
\end{proof}

\begin{proposition}\label{erafprop13}
For any positive integer $n,$ let $R_n=\left (R_{\frakX,n}\right )_1$ and $Q_{1}=\left (Q_{\frakX}\right )_1.$ Denote by $\underline{Q_1}$ the logarithmic scheme-theoretic image of the Frobenius of $Q_1$ \eqref{era3logimagedef}. The canonical morphism $R_n \ra \frakY_n$ (resp. the composition
$$\underline{Q_1} \ra Q_1 \ra Y,$$
of the canonical morphism $Q_1 \ra Y$ and the canonical closed immersion $\underline{Q_1} \ra Q_1$) factors uniquely through the strict diagonal immersion $\frakX_n \ra \frakY_n$ (resp. $X\ra Y$). In addition, the morphisms $R_n \ra \frakX_n$ and $\underline{Q_1} \ra X$ are affine.
\end{proposition}

\begin{proof}
This follows from the definitions of $R_{\frakX,n}$ and $Q_{\frakX},$ \ref{erafdil1} and \ref{erafdil2}.
\end{proof}

\begin{proposition}\label{qF}
Suppose that $X\ra S$ is of Cartier type. Let $Q_1$ be the special fiber of $Q_{\frakX}$ and $\underline{Q_1}$ the logarithmic scheme-theoretic image of the Frobenius of $Q_1$ \eqref{era3logimagedef}. Denote by $q_1,q_2:Q_1\ra X$ the canonical projections and $F_1:X\ra X'$ the exact relative Frobenius (which is equal to the relative Frobenius since $X\ra S$ is of Cartier type). Then
$$
F_1\circ q_1=F_1\circ q_2.
$$
\end{proposition}

\begin{proof}
Since $X\ra S$ is of Cartier type, $X'$ is the base change of $X$ by the Frobenius $F_S$ of $S.$ Since $Q_1 \xrightarrow{q_i} X\ra X' \ra S,$ where $X'\ra S$ is the canonical projection, is equal to the structuralmorphism $Q_1 \ra S$ for $i=1,2,$ it is sufficient to prove that $Q_1 \xrightarrow{q_i} X\ra X' \ra X$ is the same for $i=1,2,$ where $X'\ra X$ is the canonical projection. This is equal to $Q_1 \xrightarrow{q_i} X\xrightarrow{F_X} X.$ We have the following commutative diagram
$$
\begin{tikzcd}
Q_1 \ar{r}{q_i} \ar{dr}{F_{Q_1}} \ar[swap]{d}& X \ar{r}{F_X} & X \\
\underline{Q_1} \ar[hook]{r} & Q_1 \ar{ur}{q_i} & 
\end{tikzcd}
$$
By \ref{erafprop13}, the compositions $\underline{Q_1} \hookrightarrow Q_1 \xrightarrow{q_1} X$ and $\underline{Q_1} \hookrightarrow Q_1 \xrightarrow{q_2} X$ are equal. The result follows.
\end{proof}

\begin{remark}\label{rem87}
Let $r$ and $s$ be nonnegative integers. We have a canonical isomorphism
\begin{equation}\label{eq871}
\frakY(r)\times_{\frakX}\frakY(s)\xrightarrow{\sim} \frakY(r+s),
\end{equation}
where $\frakY(r)$ (resp. $\frakY(s)$) is considered as a logarithmic formal scheme over $\frakX$ via the $(r+1)^{th}$ projection (resp. first projection). For every integer $n\ge 1,$ the isomorphism (\ref{eq871}) induces isomorphisms $$R_{\frakX,n}(r)\times_{\frakX}R_{\frakX,n}(s)\xrightarrow{\sim} R_{\frakX,n}(r+s),\ P_{\frakX/\frakS,n}(r)\times_{\frakX}P_{\frakX/\frakS,n}(s)\xrightarrow{\sim} P_{\frakX/\frakS,n}(r+s),$$
$$Q_{\frakX}(r) \times_{\frakX} Q_{\frakX}(s) \xrightarrow{\sim} Q_{\frakX}(r+s),$$
where $R_{\frakX,n}(r),$ $Q_{\frakX}(r)$ and $P_{\frakX/\frakS,n}(r)$ (resp. $R_{\frakX,n}(s),$ $Q_{\frakX}(s)$ and $P_{\frakX/\frakS,n}(s)$) are considered as a logarithmic formal scheme over $\frakX$ via the $(r+1)^{th}$ projection (resp. first projection). These isomorphisms fit into the following commutative diagrams
$$
\begin{tikzcd}
P_{\frakX/\frakS,n}(r)\times_{\frakX}P_{\frakX/\frakS,n}(s) \ar{r}{\sim} \ar{d} & P_{\frakX/\frakS,n}(r+s) \ar{d} \\
R_{\frakX,n}(r)\times_{\frakX}R_{\frakX,n}(s) \ar{r}{\sim} \ar{d} & R_{\frakX,n}(r+s) \ar{d} \\
\frakY(r)\times_{\frakX}\frakY(s) \ar{r}{\sim} & \frakY(r+s)
\end{tikzcd}
$$
$$
\begin{tikzcd}
Q_{\frakX}(r)\times_{\frakX}Q_{\frakX}(s) \ar{r}{\sim} \ar{d} & Q_{\frakX}(r+s) \ar{d} \\
\frakY(r)\times_{\frakX}\frakY(s) \ar{r}{\sim} & \frakY(r+s)
\end{tikzcd}
$$
We now introduce the groupoid structures on $R_{\frakX,n},$ $Q_{\frakX}$ and $P_{\frakX/\frakS,n}.$
\end{remark}

\begin{proposition}\label{Kokoprop126}
Let $n$ be a positive integer. Consider the following morphisms:
\begin{enumerate}
\item The morphism
\begin{alignat*}{2}
\alpha & :P_{\frakX/\frakS,n}\times_{\frakX}P_{\frakX/\frakS,n} \xrightarrow{\sim} P_{\frakX/\frakS,n}(2) \ra P_{\frakX/\frakS,n} \\
(\text{resp.}\ \alpha & :R_{\frakX,n}\times_{\frakX}R_{\frakX,n} \xrightarrow{\sim} R_{\frakX,n}(2) \ra R_{\frakX,n}, \\
\text{resp.}\ \alpha & :Q_{\frakX}\times_{\frakX}Q_{\frakX} \xrightarrow{\sim} Q_{\frakX}(2) \ra Q_{\frakX}),
\end{alignat*}
induced by the $(1,3)$-projection $\frakY(2)\ra \frakY.$
\item The morphism $\iota:\frakX\ra P_{\frakX/\frakS,n}$ (resp. $\iota:\frakX\ra R_{\frakX,n},$ resp. $\iota:\frakX\ra Q_{\frakX}$).
\item The morphism
\begin{alignat*}{2}
\eta & :P_{\frakX/\frakS,n} \ra P_{\frakX/\frakS,n} \\
(\text{resp.}\ \eta & :R_{\frakX,n} \ra R_{\frakX,n}, \\
\text{resp.}\ \eta & :Q_{\frakX} \ra Q_{\frakX}),
\end{alignat*}
induced by the morphism $\frakY\ra \frakY$ that exchanges factors.
\end{enumerate} 
Then the morphisms $\alpha,\ \iota$ and $\eta$ define a $p$-adic $\frakX$-groupoid structure on $P_{\frakX/\frakS,n}$ (resp. $R_{\frakX,n},$ resp. $Q_{\frakX}$) (\cite{DXU19} 4.7).
\end{proposition}

\begin{proof}
We just check that the canonical morphisms $P_{\frakX/\frakS,n} \ra \frakY,$ $R_{\frakX,n} \ra \frakY$ and $Q_{\frakX} \ra \frakY$ factor, as maps of topological spaces, through the exact diagonal immersion $\frakX \ra \frakY.$ If $n \ge 1,$ the assertion for $P_{\frakX/\frakS,n}$ and $R_{\frakX,n}$ follows from \ref{erafprop13}. So does the assertion for $Q_{\frakX}.$ It remains to check it for $P_{\frakX/\frakS,0}.$ Let $P_0$ be the logarithmic scheme $\left (P_{\frakX/\frakS,0} \right )_1$ and $\underline{P_0}$ the logarithmic scheme-theoretic image of the Frobenius of $P_0$ \eqref{era3logimagedef}. Since $P_0,$ $\underline{P_0}$ and $P_{\frakX/\frakS,0}$ have the same underlying topological spaces, it is sufficient to prove that $\underline{P_0} \ra Y$ factors, as a map of topological spaces, through the exact diagonal immersion $X\ra Y.$ If $x$ is a local section of the ideal of the canonical immersion $X \ra P_0,$ then
$$x^p=p!x^{[p]}=0.$$
It follows that the canonical immersion $\underline{P_0} \ra P_0$ factors through $X \ra P_0.$
\end{proof}

\begin{parag}\label{paragomega}
Let $\frakG$ be a $p$-adic $\frakX$-groupoid (\cite{DXU19} 4.7) and $\omega:\frakG_{\text{zar}} \ra \frakX_{\text{zar}}$ the morphism of topoi induced by the factorization of the morphism of underlying topological spaces of $\frakG \ra \frakX^2_{\frakS}$ through the diagonal. Let $q_1,q_2:\frakG \ra \frakX$ be the projections. Then
$$\omega_*\Ox_{\frakG}=q_{1*}\Ox_{\frakG}=q_{2*}\Ox_{\frakG}.$$
In this way, we see $\Ox_{\frakG}$ as a bialgebra of $\frakX_{\text{zar}}.$ The groupoid structure on $\frakG$ induces a formal Hopf algebra structure on $\Ox_{\frakG}.$ For any positive integers $r$ and $n,$ we denote by $\calR_{\frakX,n}(r)$ and $\calQ_{\frakX}(r)$ the formal Hopf algebras $\Ox_{R_{\frakX,n}(r)}$ and $\Ox_{Q_{\frakX}(r)}$ respectively. 
\end{parag}

\begin{proposition}\label{locdescRQ}
Let $n$ be a positive integer and consider $R_{\frakX,n}$ as a logarithmic formal scheme over $\frakX$ via the first projection $q_1: R_{\frakX,n} \ra \frakX.$ Suppose that there exists a chart $\alpha:P \ra M$ of $f:\frakX \ra \frakS,$ lifting the frame $[Q] \ra [P],$ satisfying the conditions of \ref{erafsmoothdef} i.e.
\begin{enumerate}
\item $\alpha^{gp}$ is injective and the torsion subgroup of $\op{coker}\alpha^{gp}$ is finite and of order coprime with $p.$
\item The morphism $g:\frakX \ra \frakT=\frakS \times_{B\langle P\rangle}B\langle M\rangle,$ induced by $f$ and the chart $\frakX \ra B\langle M\rangle,$ is étale and strict.
\end{enumerate}
Suppose that there exist $m_1,\hdots,m_d \in \Gamma\left (\frakX,\calM_{\frakX}\right )$ lifting local coordinates of $X\ra S.$ Let $\eta_1,\hdots,\eta_d \in \Gamma \left (\frakX,\Delta^{-1}\calI_{\frakX/\frakS} \right )$ be as defined in \ref{parag77} and consider them as sections of $\Ox_{R_{\frakX,n}}$ (resp. $\Ox_{Q_{\frakX,n}}$). Then, there exists an isomorphism of $\Ox_{\frakX}$-algebras:
\begin{equation}\label{isoR}
\Ox_{\frakX}\{x_1,\hdots,x_d\} \xrightarrow{\sim} q_{1*}\Ox_{R_{\frakX,n}}
\end{equation}
\begin{equation}\label{isoQ}
\left (\text{resp. } \Ox_{\frakX}\{x_1,\hdots,x_d,y_1,\hdots,y_d\}/(y_1^p-px_1,\hdots,y_d^p-px_d) \xrightarrow{\sim} q_{1*}\Ox_{Q_{\frakX}}\right ),
\end{equation}
sending $x_i$ to $\frac{\eta_i}{p^n}$ (resp. $x_i$ to $\frac{\eta_i^p}{p}$ and $y_i$ to $\eta_i$).
\end{proposition}

\begin{proof}
By \ref{prop73} and \ref{prop74}, we have
$$\frakY = \left (\frakX \times_{\frakS}^{\op{log}}\frakX \right )\times _{B\langle \left (M\oplus_PM\right )^{sat}\rangle }B\langle N \rangle,$$
$$\frakT \times_{\frakS,[M]}^{\op{log}}\frakT = \left (\frakT \times_{\frakS}^{\op{log}}\frakT \right )\times _{B\langle \left (M\oplus_PM\right )^{sat}\rangle }B\langle N \rangle=\frakS \times_{B\langle P\rangle} B\langle N\rangle,$$
where $N$ is the inverse image of $M$ by
$$\left (M\oplus_PM\right )^{gp} \ra M^{gp},\ (x,y)\mapsto x+y.$$
The morphism $g:\frakX \ra \frakT$ induces an étale and strict morphism
$$\frakY \ra \left (\frakX \times_{\frakS}^{\op{log}}\frakT \right )\times _{B\langle \left (M\oplus_PM\right )^{sat}\rangle }B\langle N \rangle= \frakX \times_{\frakS,[M]}^{\op{log}}\frakT,$$
fitting into a commutative diagram
$$
\begin{tikzcd}
\frakX \ar[swap]{d}{\Delta} \ar{dr} & \\
\frakY \ar{r} & \frakX \times_{\frakS,[M]}^{\op{log}}\frakT.
\end{tikzcd}
$$
By \ref{propdiletale}, we get an isomorphism
$$R_{\frakX,n}=\frakY_{(\frakX_n/p^n)} \xrightarrow{\sim} \left (\frakX \times_{\frakS,[M]}^{\op{log}}\frakT\right )_{(\frakX_n/p^n)}.$$
The strict and étale morphism $g$ induces a strict and étale, hence flat, morphism 
$$h:\frakX \times_{\frakS,[M]}^{\op{log}}\frakT \ra \frakT \times_{\frakS,[M]}^{\op{log}}\frakT$$
fitting into the commutative diagram
\begin{equation}\label{diagest1}
\begin{tikzcd}
\frakX \ar{r}{g} \ar{d} & \frakT \ar{d} \\
\frakX \times_{\frakS,[M]}^{\op{log}}\frakT \ar{r}{h} \ar{d} & \frakT \times_{\frakS,[M]}^{\op{log}}\frakT \ar{d} \\
\frakX \ar{r}{g} & \frakT,
\end{tikzcd}
\end{equation}
where the upper vertical arrows are the exact diagonal morphisms and the lower vertical arrows are the projections on the first factor. Since the lower square and the outer rectangle of \eqref{diagest1} are cartesian, so is the upper square.
Then, by \ref{propdilflat}, $h$ induces an isomorphism
\begin{alignat*}{2}
\left (\frakX \times_{\frakS,[M]}^{\op{log}}\frakT\right )_{(\frakX_n/p^n)}  & \xrightarrow{\sim} \left (\frakT \times_{\frakS,[M]}^{\op{log}}\frakT\right )_{(\frakT_n/p^n)} \times_{\frakT \times_{\frakS,[M]}^{\op{log}}\frakT}\left (\frakX \times_{\frakS,[M]}^{\op{log}}\frakT\right ) \\
& \xrightarrow{\sim} \frakX \times_{\frakT} \left (\frakT \times_{\frakS,[M]}^{\op{log}}\frakT\right )_{(\frakT_n/p^n)},
\end{alignat*}
where $\left (\frakT \times_{\frakS,[M]}^{\op{log}}\frakT\right )_{(\frakT_n/p^n)}$ is considered as a logarithmic scheme over $\frakT$ via the first projection.
By this and similar isomorphisms for $Q_{\frakX},$ we can reduce to the case $\frakX=\frakT=\frakS \times_{B\langle P\rangle}B\langle M\rangle.$ We can also suppose that $\frakS=\op{Spf}A$ for a $p$-adic algebra $A.$ In this case,
\begin{alignat*}{2}
\frakX &= \frakS \times_{B\langle P\rangle}B\langle M\rangle = \op{Spf} \left ( A \widehat{\otimes}_{\Z_p\langle P \rangle} \Z_p \langle M\rangle \right ) = \op{Spf} C, \\
\frakY &= \frakS \times_{B\langle P \rangle} B\langle N \rangle = \op{Spf} \left ( A\widehat{\otimes}_{\Z_p\langle P\rangle } \Z_p\langle N\rangle \right ) = \op{Spf} D,
\end{alignat*}
where $\widehat{\otimes}$ is the $p$-adically completed tensor product. Let $I$ be the kernel of the exact diagonal $D\ra C$ and $a_1,\hdots,a_d\in I$ correspond to $\eta_1,\hdots,\eta_d.$ Then, by definition of $R_{\frakX,n}$ and $Q_{\frakX}$ and \eqref{dilaff}, we have
\begin{alignat*}{2}
R_{\frakX,n} &= \op{Spf}\left ( \left ( \frac{D\{ x_1,\hdots,x_d\}}{\left (p^nx_1-a_1,\hdots,p^nx_d-a_d \right )} \right )_{/p^n\text{-tor}} \right ), \\
Q_{\frakX} &= \op{Spf}\left ( \left ( \frac{D\{ x_1,\hdots,x_d\}}{\left (px_1-a_1^p,\hdots,px_d-a_d^p \right )} \right )_{/p\text{-tor}} \right ),
\end{alignat*}
where $p^n$-tor denotes the $p^n$-torsion ideal.
The result then follows from the lemmas \ref{lemR} and \ref{lemQ} below.
\end{proof}

\begin{lemma}\label{liftbase}
Keep the assumptions of \ref{parag86} and suppose that $\frakS=\op{Spf}A$ for a $p$-adic algebra $A$ and $\frakX=\frakS \times_{B\langle P\rangle}B\langle M\rangle=\op{Spf} C.$ The morphism $\frakX \ra \frakS$ corresponds then to a continuous morphism $\xi:A \ra C.$ In this case, $\frakY=\frakS \times_{B\langle P \rangle} B\langle N \rangle=\op{Spf} D,$ where $N$ is the inverse image of $M$ by
$$\left (M\oplus_PM\right )^{gp} \ra M^{gp},\ (x,y)\mapsto x+y,$$
and
$$
D=A\widehat{\otimes}_{\Z_p\langle P\rangle } \Z_p\langle N\rangle.
$$
Let $I$ be the ideal of the exact diagonal $D \ra C$ \eqref{BM}, $I_k$ its reduction modulo $p^k$ for $k\ge 1$ and $a_1,\hdots,a_d\in I$ whose images in $I_1/I_1^2$ are a basis of the $\left ( C/pC \right )$-module $I_1/I_1^2.$ Then, for every positive integer $k,$ the images of $a_1,\hdots,a_d$ in $I_k/I_k^2$ form a basis of the $\left ( C/p^kC \right )$-module $I_k/I_k^2.$
\end{lemma}

\begin{proof}
Let $\Omega=I/I^2$ and consider the morphism
\begin{equation}\label{eq1295}
C^d \ra \Omega,\ (\beta_1,\hdots,\beta_d) \mapsto \sum_{i=1}^d\beta_ia_i.
\end{equation}
Let $x\in I$ such that $p\ov{x}=0$ in $\Omega$ and $k$ a positive integer. Since $\frakX_k \ra \frakS_k$ is log smooth and by \eqref{iso}, $I_k/I_k^2$ is a free $C/p^kC$-module. It follows that $x\in p^{k-1}I/I^2.$ This being true for every $k$ and $I/I^2$ being separated, we deduce that $x=0.$ It follows that $I/I^2$ is flat over $\Z_p$ and so, by (\cite{Bourbaki} III §5 theorem 1), the canonical morphism
\begin{equation}\label{eq1296}
\left ((p^k)/(p^{k+1})\right ) \otimes_{\mathbb{F}_p} \left (\Omega/p\Omega\right ) \xrightarrow{\sim} p^k\Omega/p^{k+1}\Omega
\end{equation}
is an isomorphism.
Since \eqref{eq1295} is an isomorphism modulo $p,$ we conclude by \eqref{eq1296} and (\cite{Bourbaki} III §2, 8, corollary 3). 
\end{proof}

\begin{lemma}\label{lemR}
Keep the assumptions of \ref{liftbase} and recall the notation \ref{Not2}. Let $r$ be a positive integer. The morphism of topological $C$-algebras
$$\varphi:C\{x_1,\hdots,x_d\} \ra \left ( \frac{D\{ x_1,\hdots,x_d\}}{\left (p^rx_1-a_1,\hdots,p^rx_d-a_d \right )} \right )_{/p^r\text{-tor}},$$
sending $x_i$ to $x_i$ and $C$ to $D$ via the first projection, is an isomorphism.
\end{lemma}

\begin{proof}
Denote by $\pi:C \ra D$ the morphism corresponding to the first projection $\frakY \ra \frakX$ and $\Delta:D \ra C$ the morphism corresponding to the exact diagonal $\frakX \ra \frakY.$ Note that $\Delta \circ \pi=\op{Id}_{C}.$
For $g\in D \{x_1,\hdots,x_d\}$ and $J=(J_1,\hdots,J_d)\in \N^d,$ set
$$x^J=\prod_{j=1}^dx_j^{J_j},\ a^J=\prod_{j=1}^da_j^{J_j}$$
and denote by $g_J$ the coefficient of $x^J$ in $g.$
Let $f=\sum_{J\in \N^d}\alpha_Jx^J\in \op{Ker}\varphi$ with $\alpha_J\in C$ for all $J\in \N^d.$ Then there exists a positive integer $l$ and $g_1,\hdots,g_d \in D \{x_1,\hdots,x_d\}$ such that
\begin{equation}\label{eq1281}
\sum_{J\in \N^d}p^l\pi(\alpha_J)x^J=\sum_{i=1}^dg_i(p^rx_i-a_i).
\end{equation}
We will prove by induction on $n$ that $p^l\alpha_{J}=0$ for all $J\in \N^d$ such that $|J|=n$ and that
$$
\sum_{\substack{1\le i\le d \\ |J|=n}}g_{i,J}a^{J+\epsilon_i}=0.
$$
Taking the constant coefficients of both sides in \eqref{eq1281}, we get
$$
p^l\pi (\alpha_0)=-\sum_{i=1}^dg_{i,0}a_i.
$$
Applying $\Delta$ to both sides, we get
$$p^l\Delta\circ \pi (\alpha_0)=p^l\alpha_0=0.$$
It follows that
$$\sum_{i=1}^dg_{i,0}a_i=0.$$
The first step of the induction is thus proven.
Let $n\ge 0$ be an integer and suppose that $p^l\alpha_{J}=0$ for all $J\in \N^d$ such that $|J|=n$ and that
$$
\sum_{\substack{1\le i\le d \\ |J|=n}}g_{i,J}a^{J+\epsilon_i}=0.
$$
Taking the monomials of degree $n+1$ in \eqref{eq1281}, we get
$$
\sum_{|J|=n+1}p^l\pi(\alpha_J)x^J=\sum_{\substack{1\le i\le d\\ |J|=n}}p^rg_{i,J}x^{J+\epsilon_i}-\sum_{\substack{1\le i\le d\\ |J|=n+1}}g_{i,J}a_ix^J.
$$
Taking $x_i=a_i$ for all $1\le i\le d$ and using the induction hypothesis, we get
\begin{equation}\label{eq1284}
\sum_{|J|=n+1}p^l\pi(\alpha_J)a^J=-\sum_{\substack{1\le i\le d\\ |J|=n+1}}g_{i,J}a^{J+\epsilon_i}.
\end{equation}
It follows that, in $I^{n+1}/I^{n+2},$
$$\sum_{|J|=n+1}p^l\pi(\alpha_J)a^J=0.$$
Let $k$ be a positive integer. In the $\left (C/p^kC\right )$-module $I_k^{n+1}/I_k^{n+2},$ we have
$$\sum_{|J|=n+1}p^l\alpha_Ja^J=0.$$
By \ref{liftbase}, the family $(a_1,\hdots,a_d)$ is a basis of the $\left (C/p^kC\right )$-module $I_k/I_k^2.$ Then, by \eqref{grsym}, $(a^J)_{|J|=n+1}$ is free in the $\left (C/p^kC\right )$-module $I_k^{n+1}/I_k^{n+2}.$ It follows that
$p^l\alpha_J=0$ for $|J|=n+1.$
Then, by \eqref{eq1284}, we get
$$\sum_{\substack{1\le i\le d\\ |J|=n+1}}g_{i,J}a^{J+\epsilon_i}=0.$$
This concludes the induction. Since $C$ is flat over $\Z_p,$ we get $\alpha_J=0$ for all $J\in \N^d.$ This concludes the injectivity of $\varphi.$ The surjectivity follows from the fact that, for $y\in D,$
$$y=\pi \circ \Delta(y)+y-\pi \circ \Delta(y)$$
and $y-\pi \circ \Delta(y)\in I.$
\end{proof}

\begin{lemma}\label{lemQ}
Keep the assumptions of \ref{liftbase} and recall the notation \ref{Not6}. The morphism of $C$-modules
$$\varphi: \frac{C\{x_1,\hdots,x_d,y_1,\hdots,y_d\}}{(px_1-y_1^p,\hdots,px_r-y_r^p)} \ra \left ( \frac{D\{ x_1,\hdots,x_d\}}{\left (px_1-a_1^p,\hdots,px_d-a_d^p \right )} \right )_{/p\text{-tor}},$$
sending $x_i$ to $x_i,$ $y_i$ to $a_i$ and $C$ to $D$ via the first projection, is an isomorphism.
\end{lemma}

\begin{proof}
We keep the notations introduced in the proof of \ref{lemR}.
Let us prove that $\varphi$ is injective. Let
$$f=\sum_{\substack{K\in \N^d \\ J\in \llbracket 0,p-1\rrbracket^d}}\alpha_{K,J}x^Ky^J\in C\{x_1,\hdots,x_d,y_1,\hdots,y_d\}$$
such that $\varphi(f)=0.$
There exists a positive integer $l$ and $g_1,\hdots,g_d\in D \{x_1,\hdots,x_d\}$ such that
\begin{equation}\label{eq1291}
p^l\sum_{\substack{K\in \N^d \\ J\in \llbracket 0,p-1\rrbracket^d}}\pi(\alpha_{K,J})a^Jx^K=\sum_{i=1}^dg_i(px_i-a_i^p).
\end{equation}
From here on, we work modulo $p^k$ for a fixed positive integer $k.$
We prove by induction on $n$ that $p^l\alpha_{K,J}=0$ for all $(K,J)\in \N^d \times \llbracket 0,p-1\rrbracket^d$ such that $|K|=n$ and that 
$$
\sum_{\substack{1\le i\le d \\ |K|=n}}g_{i,K}a^{pK+p\epsilon_i}=0.
$$
Taking the constant coefficients of both sides in \eqref{eq1291}, we get
\begin{equation}\label{eq1294}
p^l\sum_{J\in \llbracket 0,p-1\rrbracket^d}\pi(\alpha_{0,J})a^J=-\sum_{i=1}^dg_{i,0}a_i^p.
\end{equation}
We rewrite this as
\begin{equation}\label{eq1297}
\sum_{J\in \N^d}\beta_Ja^J=0,
\end{equation}
where
$$\beta_J=\begin{cases} p^l\pi(\alpha_{0,J})\ \text{if}\ J\in \llbracket 0,p-1\rrbracket^d,\\ g_{i,0}\ \text{if}\ J=p\epsilon_i,\\  0\ \text{otherwise.}\end{cases}$$ 
Proceeding by induction on $|J|,$ reducing modulo $I_k,$ $I_k^2,\hdots$ and using \eqref{grsym}, we prove that $\Delta(\beta_J)=0$ for all $J.$ Indeed, applying $\Delta$ to \eqref{eq1297}, we get
$$\Delta(\beta_0)=0.$$
Let $m\ge 0$ be an integer and suppose that $\Delta(\beta_J)=0$ for $|J|\le m.$ In the $\left (C/p^kC\right )$-module $I_k^{m+1}/I_k^{m+2},$ we have
$$\sum_{|J|=m+1}\Delta (\beta_J)a^J=0.$$
By \ref{liftbase}, the family $(a_1,\hdots,a_d)$ is a basis of the $\left (C/p^kC\right )$-module $I_k/I_k^2.$ Then, by \eqref{grsym}, $(a^J)_{|J|=m+1}$ is free in the $\left (C/p^kC\right )$-module $I_k^{m+1}/I_k^{m+2}$ so $\Delta(\beta_J)=0$ for $|J|=m+1.$
In particular, if $J\in \llbracket 0,p-1\rrbracket^d,$ then
$$\Delta(\beta_J)=p^l\Delta\circ \pi(\alpha_{0,J})=p^l\alpha_{0,J}=0.$$
By \eqref{eq1294}, we deduce that
$$\sum_{i=1}^dg_{i,0}a_i^p=0.$$
This concludes the first step of the induction.

Now taking the monomials of degree $1$ in \eqref{eq1291}, we get
$$p^l\sum_{\substack{1\le i\le d \\ J\in \llbracket 0,p-1 \rrbracket^d}} \pi(\alpha_{\epsilon_i,J})a^Jx_i=\sum_{i=1}^dpg_{i,0}x_i-\sum_{1\le i,j\le d}g_{i,\epsilon_j}a_i^px_j.$$
Taking $x_i=a_i^p,$ we get
$$p^l\sum_{\substack{1\le i\le d \\ J\in \llbracket 0,p-1 \rrbracket^d}} \pi(\alpha_{\epsilon_i,J})a^{J+p\epsilon_i}=-\sum_{1\le i,j\le d}g_{i,\epsilon_j}a^{p\epsilon_i+p\epsilon_j}.$$
Just as in the previous case, we prove that $p^l\alpha_{\epsilon_i,J}=0$ for all $1\le i\le d$ and $J.$
Let $n$ be a positive integer and suppose that $p^l\alpha_{K,J}=0$ for all $K\in \N^d$ and $J\in \llbracket 1,p-1\rrbracket^d$ such that $|K|\le n$ and that
\begin{equation}\label{eq1293}
\sum_{\substack{|K|= n \\ 1\le i\le d}}g_{i,I}a^{pK+p\epsilon_i}=0.
\end{equation}
Taking the monomials of degree $n+1$ in \eqref{eq1291}, we get
$$p^l\sum_{\substack{|K|=n+1 \\ J\in \llbracket 0,p-1 \rrbracket^d}} \pi(\alpha_{K,J})a^Jx^K=\sum_{\substack{|K|=n \\ 1\le i\le d}}pg_{i,K}x^{K+\epsilon_i}-\sum_{\substack{|K|=n+1 \\ 1\le i\le d}}g_{i,K}a_i^px^K.$$
Taking $x_i=a_i^p$ for all $1\le i\le d,$ we get
$$p^l\sum_{\substack{|K|=n+1 \\ J\in \llbracket 0,p-1 \rrbracket^d}} \pi(\alpha_{K,J})a^{J+pK}=\sum_{\substack{|K|=n \\ 1\le i\le d}}pg_{i,K}a^{pI+p\epsilon_i}-\sum_{\substack{|K|=n+1 \\ 1\le i\le d}}g_{i,K}a^{pK+p\epsilon_i}.$$
By \eqref{eq1293}, we get
$$ p^l\sum_{\substack{|K|=n+1 \\ J\in \llbracket 0,p-1 \rrbracket^d}} \pi(\alpha_{K,J})a^{J+pK}=-\sum_{\substack{|K|=n+1 \\ 1\le i\le d}}g_{i,K}a^{pK+p\epsilon_i}.$$
Again, using the same argument as in the case $n=0,$ we prove that $p^l\alpha_{K,J}=0$ for $|K|=n+1$ and $J\in \llbracket 0,p-1 \rrbracket^d$ and then
$$\sum_{\substack{|K|= n+1 \\ 1\le i\le d }}g_{i,K}a^{pK+p\epsilon_i}=0.
$$
This concludes the induction, and since this is true modulo $p^k$ for all $k,$ this concludes the injectivity of $\varphi.$ The surjectivity follows from the fact that, for $y\in D,$
$$y=\pi \circ \Delta(y)+y-\pi \circ \Delta(y)$$
and $y-\pi \circ \Delta(y)\in I.$
\end{proof}

\begin{parag}\label{par1212}
Keep the assumptions and notation of \ref{parag86}. By \ref{formaletaleframelift}, there exists an étale covering $\left ( \frakS_i \ra \frakS \right )_{i\in I},$ and, for every $i\in I,$ an étale covering $\left (\frakX_{ij} \ra \frakX_i:=\frakX\times_{\frakS}\frakS_i \right )_{j\in I_i}$ and charts $P_i \ra M_{ij}$ of $\frakX_{ij} \ra \frakS_i,$ fitting into a commutative diagram
$$
\begin{tikzcd}
\frakX_{ij} \ar{r} \ar{d} & \frakS_i \ar{d} \\
B\langle M_{ij} \rangle \ar{r} \ar{d} & B\langle P_i \rangle \ar{d} \\
\left [Q\right ] \ar{r} & \left [ P\right ], 
\end{tikzcd}
$$
and satisfying:
\begin{enumerate}
\item $P_i^{gp} \ra M_{ij}^{gp}$ is injective and the torsion subgroup of its cokernel has a finite order coprime with $p.$
\item The morphism $\frakX_{ij} \ra \frakS_i\times_{B\langle P_i \rangle} B\langle M_{ij} \rangle,$ induced by $\frakX_{ij}\ra \frakS_i$ and the chart $\frakX_{ij} \ra B\langle M_{ij}\rangle,$ is étale and strict.
\end{enumerate}
For any $i\in I$ and $j\in I_i,$ let $X_{ij}$ and $S_i$ be the special fibers of $\frakX_{ij}$ and $\frakS_i$ respectively. After eventually shrinking $X_{ij},$ we can suppose that $X_{ij} \ra S_i$ has local coordinates $m_1,\hdots,m_d\in \Gamma \left (X_{ij},\calM_{X_{ij}} \right ).$ By \ref{prop69}, the morphism $\calM_{\frakX_{ij}}\ra \calM_{X_{ij}}$ is surjective. Hence, after eventually a furthermore shrinking of $X_{ij},$ there exists $\widetilde{m}_1,\hdots,\widetilde{m}_d\in \Gamma\left (\frakX_{ij},\calM_{\frakX_{ij}}\right )$ that lift $m_1,\hdots,m_d.$ We can then define $\widetilde{\eta}_1,\hdots,\widetilde{\eta}_d\in \Gamma\left (\frakU_i,\Delta^{-1}\calI_{\frakX/\frakS} \right )$ from $\widetilde{m}_d,\hdots,\widetilde{m}_d$ as in \ref{parag77}. This proves that the isomorphisms \eqref{isoR} and \eqref{isoQ} exist étale locally on $\frakX$ and $\frakS.$
\end{parag}

\begin{proposition}\label{Qlogflat}
For any positive integer $n,$ the logarithmic formal schemes $Q_{\frakX}$ and $R_{\frakX,n}$ are log flat over $\op{Spf}W$ \eqref{logflatdef} and hence flat over $\op{Spf}W$ (\ref{Wflat} and \ref{dxuflat}).
\end{proposition}

\begin{proof}
By \ref{locdescRQ} and \ref{par1212}, for every positive integer $k,$ $R_{\frakX,n,k} \ra \frakX_k$ and $Q_{\frakX,k}\ra \frakX_k$ are flat. Since they are also strict, they are log flat. We conclude by the log flatness of $\frakX \ra \op{Spf}W.$
\end{proof}

\begin{proposition}\label{propfree}
Let $n$ be a positive integer and consider $P_{\frakX/\frakS,n}$ as a logarithmic formal scheme over $\frakX$ by the first projection $q_1:P_{\frakX/\frakS,n} \ra \frakX.$
Under the assumptions of \ref{locdescRQ}, we have an isomorphism of $\Ox_{\frakX}$-algebras
$$\Ox_{\frakX}\langle \langle x_1,\hdots,x_d \rangle \rangle \xrightarrow{\sim} \calP_{\frakX/\frakS,n}$$
sending $x_i$ to $\frac{\eta_i}{p^n},$ where $\Ox_{\frakX}\langle \langle x_1,\hdots,x_d \rangle \rangle$ is the $p$-adic completion of the PD algebra $\Ox_{\frakX}\langle x_1,\hdots,x_d \rangle.$
\end{proposition}

\begin{proof}
Let $k$ be a positive integer and $\eta_{i,k}$ the image of $\frac{\eta_i}{p^n}$ in $\calP_{\frakX/\frakS,n,k}.$ By the local description of $R_{\frakX,n}$ given in \ref{locdescRQ}, there exists an isomorphism of $\Ox_{\frakX_k}$-algebras
$$
\Ox_{\frakX_k}\langle x_1,\hdots,x_d \rangle \xrightarrow{\sim} \calP_{\frakX/\frakS,n,k}
$$
sending $x_i$ to $\eta_{i,k}.$
The result then follows by taking the projective limit for $k\ge 1.$
\end{proof}

\begin{proposition}\label{propflat}
For any nonnegative integers $n$ and $r,$ the formal scheme $P_{\frakX / \frakS,n}(r)$ is log flat over $\op{Spf}W$ \eqref{logflatdef} and hence flat over $\op{Spf}W$ (\ref{Wflat} and \ref{dxuflat}).
\end{proposition}

\begin{proof}
By hypothesis, $\frakX$ is log flat over $\op{Spf}W.$
By \ref{propfree}, the projections $P_{\frakX/\frakS,n,k}(r) \ra \frakX_k$ are flat for every positive integer $k.$ Since they are also strict, they are log flat. It follows that $P_{\frakX / \frakS,n,k}(r)$ is log flat over $\op{Spec}\left (W/(p^k)\right )$ for every positive integer $k.$ The result follows \eqref{logflatdef}.
\end{proof}

\begin{parag}\label{loccoord}
We will often make the following hypothesis:
suppose that we have local coordinates $m_1,\hdots,m_d\in \Gamma(X,\calM_X)$ for $X\ra S$ (\ref{P2}) that lift to local sections $\widetilde{m}_1,\hdots,\widetilde{m}_d\in \Gamma(\frakX,\calM_{\frakX}).$ For every $1\le i\le d,$ let $m_i'=\pi^{\flat}m_i$ \eqref{diag51} and $\widetilde{\eta}_i=\eta(\widetilde{m}_i) \in \Gamma(\frakX,\Delta^{-1}\calI_{\frakX/\frakS})$ be as defined in \ref{parag77} and $\eta_i$ the reduction of $\widetilde{\eta}_i$ modulo $p.$ Let $n$ be a positive integer. By \ref{locdescRQ} and  \ref{propfree}, we have isomorphisms:
\begin{alignat}{2}
\calR_{\frakX,n} & \xrightarrow{\sim} \Ox_{\frakX}\left \{ \frac{\widetilde{\eta}_1}{p^n},\hdots,\frac{\widetilde{\eta}_d}{p^n} \right \}, \\
\calQ_{\frakX} & \xrightarrow{\sim} \Ox_{\frakX}\{ \widetilde{\eta}_1,\hdots,\widetilde{\eta}_d,y_1,\hdots,y_d\}/(py_1-\widetilde{\eta}_1^p,\hdots,py_d-\widetilde{\eta}_d^p), \\
\calP_{\frakX/\frakS,n} & \xrightarrow{\sim} \Ox_{\frakX}\left \langle \left \langle \frac{\widetilde{\eta}_1}{p^n},\hdots,\frac{\widetilde{\eta}_d}{p^n} \right \rangle \right \rangle. \label{eqPkraz6}
\end{alignat}
Let $\eta_{i(r),n}$ be the image of $\frac{\widetilde{\eta}_i}{p^n}$ in $\calR_n:=\calR_{\frakX,n}/p\calR_{\frakX,n}$ and $\eta_{i(q)}$ the image of $\frac{\widetilde{\eta}_i^p}{p}$ in $\calQ_1:=\calQ_{\frakX}/p\calQ_{\frakX}.$ For $n=1,$ we denote $ \eta_{i(r),1}$ simply by $\eta_{i(r)}.$ For any $I\in \N^d,$ we set
$$
\eta^I=\prod_{i=1}^d\eta_i^{I_i},\ \eta_{(q)}^I=\prod_{i=1}^d\eta_{i(q)}^{I_i},\ \eta_{(r)}^I=\prod_{i=1}^d\eta_{i(r)}^{I_i}.
$$
Then $\left (\eta^I\eta_{(q)}^J \right )_{\substack{I\in \llbracket 0,p-1 \rrbracket^d \\ J\in \N^d}}$ is a basis of the $\Ox_X$-module $\calQ_1.$ We denote by $\left ( \varphi_{I,J} \right )$ its dual basis.
If $X'$ lifts to a log smooth logarithmic $\frakS$-scheme $\frakX',$ we define $\eta_i',$ $\widetilde{\eta}_i',$ $\eta_{i(r)}'$ and $\eta_{i(q)}'$ in a similar way from $m_i'.$ We also denote by $(\partial_1',\hdots,\partial_d')$ the dual basis of $(\op{dlog}m_1',\hdots,\op{dlog}m_d').$
\end{parag}

\begin{proposition}\label{HopffrakR}
Let
\begin{alignat*}{2}
\delta:&\calR_{\frakX,n} \ra \calR_{\frakX,n}\otimes_{\Ox_{\frakX}}\calR_{\frakX,n}\\
\pi:&\calR_{\frakX,n} \ra \Ox_{\frakX}\\
\sigma:&\calR_{\frakX,n} \ra \calR_{\frakX,n}
\end{alignat*}
the morphisms defining the Hopf algebra structure on $\calR_{\frakX,n}.$ Under the hypothesis \ref{loccoord},
we have
\begin{alignat*}{2}
\delta(\widetilde{\eta}_i) &= 1\otimes \widetilde{\eta}_i+\widetilde{\eta}_i\otimes 1+\widetilde{\eta}_i\otimes \widetilde{\eta}_i, \\
\pi(\widetilde{\eta}_i) &= 0,\\
\sigma(\widetilde{\eta}_i) &= (1+\widetilde{\eta}_i)^{-1}-1=-\widetilde{\eta}_i+\widetilde{\eta}_i^2-\hdots
\end{alignat*}
\end{proposition}

\begin{proof}
This is a result of \ref{era2prop611} and \ref{era2prop612}. Note that the series $-\widetilde{\eta}_i+\widetilde{\eta}_i^2-\hdots$ is convergent since $\widetilde{\eta}_i=p^n\frac{\widetilde{\eta}_i}{p^n}\in p\calR_{\frakX,n}.$
\end{proof}

\begin{proposition}\label{HopfR}
Let $n$ be a positive integer, $\calR_n=\calR_{\frakX,n}/p\calR_{\frakX,n}$ and
\begin{alignat*}{2}
\delta:&\calR_n \ra \calR_n\otimes_{\Ox_X}\calR_n\\
\pi:&\calR_n \ra \Ox_X\\
\sigma:&\calR_n \ra \calR_n
\end{alignat*}
the morphisms defining the Hopf algebra structure on $\calR_n.$ Under the hypothesis \ref{loccoord},
we have
\begin{alignat*}{2}
\delta\left (\eta_{i(r),n} \right ) &= 1\otimes \eta_{i(r),n}+\eta_{i(r),n}\otimes 1, \\
\pi\left (\eta_{i(r),n} \right ) &= 0,\\
\sigma\left (\eta_{i(r),n} \right ) &= -\eta_{i(r),n}.
\end{alignat*}
\end{proposition}

\begin{proof}
This results from \ref{HopffrakR} and the fact that $\widetilde{\eta}_i=p^n\frac{\widetilde{\eta}_i}{p^n}\in p\calR_{\frakX,n}.$
\end{proof}

\begin{proposition}\label{HopffrakQ}
Let
\begin{alignat*}{2}
\widetilde{\delta}:&\calQ_{\frakX} \ra \calQ_{\frakX}\otimes_{\Ox_{\frakX}}\calQ_{\frakX}\\
\widetilde{\pi}:&\calQ_{\frakX} \ra \Ox_{\frakX}\\
\widetilde{\sigma}:&\calQ_{\frakX} \ra \calQ_{\frakX}
\end{alignat*}
the morphisms defining the Hopf algebra structure on $\calQ_{\frakX}.$ Under the hypothesis \ref{loccoord},
we have
\begin{alignat*}{2}
\widetilde{\delta}(\widetilde{\eta}_i) &= 1\otimes \widetilde{\eta}_i+\widetilde{\eta}_i\otimes 1+\widetilde{\eta}_i\otimes \widetilde{\eta}_i, \\
\widetilde{\pi}(\widetilde{\eta}_i) &= 0,\\
\widetilde{\sigma}(\widetilde{\eta}_i) &= (1+\widetilde{\eta}_i)^{-1}-1.
\end{alignat*}
\end{proposition}

\begin{proof}
This is a result of \ref{era2prop611} and \ref{era2prop612}.
\end{proof}

\begin{proposition}\label{HopfQ}
Let $\calQ_1=\calQ_{\frakX}/p\calQ_{\frakX}$ and
\begin{alignat*}{2}
\delta:&\calQ_1 \ra \calQ_1\otimes_{\Ox_X}\calQ_1\\
\pi:&\calQ_1 \ra \Ox_X\\
\sigma:&\calQ_1 \ra \calQ_1
\end{alignat*}
the morphisms defining the Hopf algebra structure on $\calQ_1.$ Under the hypothesis \ref{loccoord},
we have
\begin{alignat*}{2}
\delta\left (\eta_{i(q)} \right ) &= 1\otimes \eta_{i(q)}+\sum_{0<b+c<p} \frac{(-1)^{b+c}}{b+c}\begin{pmatrix}b+c \\ b \end{pmatrix} \eta_i^{b+c}\otimes \eta_i^{p-b} +\eta_{i(q)}\otimes 1, \\
\delta(\eta_i) &= 1\otimes \eta_i+\eta_i\otimes 1+\eta_i\otimes \eta_i, \\
\pi\left (\eta_{i(q)} \right ) &= \pi(\eta_i)= 0,\\
\sigma\left (\eta_{i(q)} \right ) &= \frac{-\eta_{i(q)}}{(1+\eta_i)^p}.
\end{alignat*}
\end{proposition}

\begin{proof}
Let
\begin{alignat*}{2}
\widetilde{\delta}:&\calQ_1 \ra \calQ_1\otimes_{\Ox_X}\calQ_1\\
\widetilde{\pi}:&\calQ_1 \ra \Ox_X\\
\widetilde{\sigma}:&\calQ_1 \ra \calQ_1
\end{alignat*}
be the morphisms defining the Hopf algebra structure on $\calQ_{\frakX}.$ By \ref{HopffrakQ}, we have
\begin{alignat*}{2}
p\widetilde{\delta}\left (\frac{\widetilde{\eta}_i^p}{p}\right ) =& \widetilde{\delta}(\widetilde{\eta}_i)^p \\
=& \left (1\otimes \widetilde{\eta}_i + \widetilde{\eta}_i\otimes 1 + \widetilde{\eta}_i\otimes \widetilde{\eta}_i\right )^p\\
=& \sum_{a+b+c=p} \begin{pmatrix}p \\ b+c \end{pmatrix} \begin{pmatrix} b+c \\ b\end{pmatrix} \widetilde{\eta}_i^{b+c} \otimes \widetilde{\eta}_i^{a+c} \\
=& 1\otimes \widetilde{\eta}_i^p + \sum_{0<b+c<p} \begin{pmatrix}p \\ b+c \end{pmatrix} \begin{pmatrix} b+c \\ b\end{pmatrix} \widetilde{\eta}_i^{b+c} \otimes \widetilde{\eta}_i^{p-b} + \sum_{b+c=p}\begin{pmatrix}p \\ b \end{pmatrix}\widetilde{\eta}_i^p\otimes \widetilde{\eta}_i^c \\
=& 1\otimes \widetilde{\eta}_i^p + \sum_{0<b+c<p} \begin{pmatrix}p \\ b+c \end{pmatrix} \begin{pmatrix} b+c \\ b\end{pmatrix} \widetilde{\eta}_i^{b+c} \otimes \widetilde{\eta}_i^{p-b} + \sum_{b=1}^{p-1}\begin{pmatrix}p \\ b \end{pmatrix}\widetilde{\eta}_i^p\otimes \widetilde{\eta}_i^{p-b} + \widetilde{\eta}_i^p\otimes 1 + \widetilde{\eta}_i^p\otimes \widetilde{\eta}_i^p \\
=& p \left (1\otimes \frac{\widetilde{\eta}_i^p}{p}\right )  + \sum_{0<b+c<p} \begin{pmatrix}p \\ b+c \end{pmatrix} \begin{pmatrix} b+c \\ b\end{pmatrix} \widetilde{\eta}_i^{b+c} \otimes \widetilde{\eta}_i^{p-b} + \sum_{b=1}^{p-1}\begin{pmatrix}p \\ b \end{pmatrix}\widetilde{\eta}_i^p\otimes \widetilde{\eta}_i^{p-b} + p\left (\frac{\widetilde{\eta}_i^p}{p}\otimes 1 \right )\\
& + p^2\left (\frac{\widetilde{\eta}_i^p}{p}\otimes \frac{\widetilde{\eta}_i^p}{p} \right ).
\end{alignat*}
By the flatness of $Q_{\frakX}$ over $\op{Spf}W$ (\ref{Qlogflat} and \ref{dxuflat}), we get
\begin{alignat*}{2}
\widetilde{\delta}\left (\frac{\widetilde{\eta}_i^p}{p} \right )=&  1\otimes \frac{\widetilde{\eta}_i^p}{p}  + \sum_{0<b+c<p} \frac{(p-1)!}{(b+c)!(p-b-c)!} \begin{pmatrix} b+c \\ b\end{pmatrix} \widetilde{\eta}_i^{b+c} \otimes \widetilde{\eta}_i^{p-b} \\
& + \sum_{b=1}^{p-1}\frac{(p-1)!}{b!(p-b)!}\widetilde{\eta}_i^p\otimes \widetilde{\eta}_i^{p-b} + \frac{\widetilde{\eta}_i^p}{p}\otimes 1 \\
& + p\left (\frac{\widetilde{\eta}_i^p}{p}\otimes \frac{\widetilde{\eta}_i^p}{p} \right ).
\end{alignat*}
Since
$$
\sum_{b=1}^{p-1}\frac{(p-1)!}{b!(p-b)!}\widetilde{\eta}_i^p\otimes \widetilde{\eta}_i^{p-b}=\sum_{b=1}^{p-1}\begin{pmatrix}p\\ b \end{pmatrix} \frac{\widetilde{\eta}_i^p}{p}\otimes \widetilde{\eta}_i^{p-b},
$$
this sum vanishes modulo $p.$
By \ref{HopffrakQ},
$$
(1+\widetilde{\eta}_i)\widetilde{\sigma}(\widetilde{\eta}_i)=-\widetilde{\eta}_i
$$
so
$$
p(1+\widetilde{\eta}_i)^p\widetilde{\sigma}\left (\frac{\widetilde{\eta}_i^p}{p}\right )=(1+\widetilde{\eta}_i)^p\widetilde{\sigma}(\widetilde{\eta}_i^p)=-\widetilde{\eta}_i^p=-p\frac{\widetilde{\eta}_i^p}{p}.
$$
It follows, by the flatness of $Q_{\frakX}$ over $\op{Spf}W$ (\ref{Qlogflat} and \ref{dxuflat}), that
$$
(1+\widetilde{\eta}_i)^p\widetilde{\sigma} \left (\frac{\widetilde{\eta}_i^p}{p} \right )=-\frac{\widetilde{\eta}_i^p}{p}.
$$
The result then follows from \ref{HopffrakQ} and the fact that, for any $1\le k\le p-1,$ we have the equality $\frac{(p-1)!}{k!(p-k)!}=\frac{(-1)^k}{k}$ in $\mathbb{F}_p.$
\end{proof}

\begin{proposition}\label{HopffrakP}
Let
\begin{alignat*}{2}
\delta:&\calP_{\frakX / \frakS,n} \ra \calP_{\frakX / \frakS,n}\otimes_{\Ox_{\frakX}}\calP_{\frakX / \frakS,n} \\
\pi:&\calP_{\frakX / \frakS,n} \ra \Ox_{\frakX}\\
\sigma:&\calP_{\frakX / \frakS,n} \ra \calP_{\frakX / \frakS,n}
\end{alignat*}
the morphisms defining the Hopf algebra structure on $\calP_{\frakX / \frakS,n}.$ Under the hypothesis \ref{loccoord},
we have
\begin{alignat*}{2}
\delta(\widetilde{\eta}_i) &= 1\otimes \widetilde{\eta}_i+\widetilde{\eta}_i\otimes 1+\widetilde{\eta}_i\otimes \widetilde{\eta}_i, \\
\pi(\widetilde{\eta}_i) &= 0,\\
\sigma(\widetilde{\eta}_i) &= (1+\widetilde{\eta}_i)^{-1}-1=-\widetilde{\eta}_i+\widetilde{\eta}_i^2-\hdots
\end{alignat*}
\end{proposition}

\begin{proof}
This is a result of \ref{era2prop611} and \ref{era2prop612}. Note that the series $-\widetilde{\eta}_i+\widetilde{\eta}_i^2-\hdots$ is convergent since $\widetilde{\eta}_i=p^n\frac{\widetilde{\eta}_i}{p^n}\in p\calP_{\frakX / \frakS,n}.$
\end{proof}

\begin{definition}\label{defhpdstrat}
Let $k\ge 1$ and $n\ge 0$ be integers and $\calE$ an $\Ox_{\frakX_k}$-module. Let $P_{n,k}$ be the logarithmic scheme obtained from $P_{\frakX/\frakS,n}$ by reduction modulo $p^k.$ 
An \emph{$n$-stratification on $\calE$} (resp. \emph{$n$-HPD-stratification on $\calE$}) is a stratification (resp. HPD stratification) (\cite{DXU19} 5.4) on $\calE$ with respect to the Hopf algebra corresponding to the groupoid $P_{n,k}$ (\ref{Kokoprop126} and \ref{paragomega}).
We denote by $n$-$\op{MHS}(\frakX_k/\frakS_k)$ the category of $\Ox_{\frakX_k}$-modules equipped with $n$-HPD-stratifications.
\end{definition}

\newpage

\section{\texorpdfstring{$Q_{\frakX}$}{Q} %
     and \texorpdfstring{$R_{\frakX,n}$}%
     {R}-stratifications}

\begin{parag}
In this section, we keep the set-up of the previous one i.e. we keep the notation and assumption of \ref{paragJapon1} and \ref{parag86}. In other words, we consider a perfect field of positive characteristic $p,$ denote by $W$ its ring of Witt vectors and equip $\op{Spf}W$ woth the trivial log structure. We consider a log smooth morphism of framed logarithmic $p$-adic formal schemes $f:(\frakX,Q)\ra (\frakS,0),$ such that $\frakS$ is log flat and locally of finite type over $\op{Spf}W.$ Note that this implies, by \ref{Wflat}, that the formal schemes $\frakX$ and $\frakS$ are flat over $\op{Spf}W$ \eqref{dxuflat}. We also consider the formal groupoids $Q_{\frakX},$ $R_{\frakX,n}$ and $P_{\frakX/\frakS,n}$ defined in \ref{parag86}.
\end{parag}

\begin{proposition}\label{wiw1}
There exists a unique morphism of logarithmic formal schemes
\begin{equation}\label{wiwkraz}
P_{\frakX/\frakS,0} \ra Q_{\frakX},
\end{equation}
such that the diagram
$$
\begin{tikzcd}
 & Q_{\frakX} \ar{d} \\
P_{\frakX/\frakS,0} \ar{r} \ar{ur} & \frakY
\end{tikzcd}
$$
is commutative. In addition, this morphism is compatible with the groupoid structures of $P_{\frakX/\frakS,0}$ and $Q_{\frakX}.$
\end{proposition}

\begin{proof}
If $x$ is a local section of the ideal of the immersion $X=\frakX_1 \ra \frakY_1=Y$ then its image in $P_0= \left (P_{\frakX/\frakS,0}\right )_1$ satisfies
$$x^p=p! x^{[p]}=0.$$
It follows that the composition
$$\underline{P_0} \ra P_0 \ra Y$$
factors through $X.$
Since $P_{\frakX/\frakS,0}$ is flat over $\op{Spf}W$ \eqref{propflat}, we conclude by the universal property of $Q_{\frakX}$ (\ref{erafdil2}).
\end{proof} 

\begin{parag}
Let $P_0$ and $Q_1$ be the special fibers of $P_{\frakX/\frakS,0}$ and $Q_{\frakX},$ and let $\calP_0$ and $\calQ_1$ be the Hopf algebras defined by $P_0$ and $Q_1$ respectively.
The morphism $P_{\frakX/\frakS,0} \ra Q_{\frakX}$ \eqref{wiw1} induces a morphism of groupoids $P_0 \ra Q_1$ and then a morphism of Hopf algebras
\begin{equation}\label{Muzanu1}
u:\calQ_1 \ra \calP_0.
\end{equation}
By taking the dual, we obtain a morphism of $\Ox_X$-algebras
\begin{equation}\label{checku}
\check{u}:\widehat{\calD}_{X/S} \ra \calQ_1^{\vee},
\end{equation}
where $\widehat{\calD}_{X/S}:=\mathscr{Hom}_{\Ox_X}(\calP_0,\Ox_X)$ is the sheaf of hyper PD-differential operators.
Under the hypothesis \ref{loccoord}, since $(p-1)!\equiv -1\ (\text{mod}\ p),$ we clearly have, for all $1\le i\le d,$
\begin{equation}\label{ueta}
u(\eta_i)=\eta_i,\ u(\eta_{i(q)})=-\eta_i^{[p]}.
\end{equation}
For an upcoming proposition, we recall the canonical Hopf algebra structure on $S^{\bullet}\omega^1_{X/S}:$
\begin{equation}\label{HopfalgstrS}
\begin{alignedat}{2}
\delta : &\begin{array}[t]{clc}
S^{\bullet}\omega^1_{X/S} & \ra & S^{\bullet}\omega^1_{X/S} \otimes_{\Ox_X} S^{\bullet}\omega^1_{X/S} \\
x & \mapsto & 1\otimes x+x\otimes 1,
\end{array} \\
\sigma : &\begin{array}[t]{clc}
S^{\bullet}\omega^1_{X/S} & \ra & S^{\bullet}\omega^1_{X/S} \\
x & \mapsto & -x,
\end{array} \\
\pi : &\begin{array}[t]{clc}
S^{\bullet}\omega^1_{X/S} & \ra & \Ox_X \\
x=(x_0,x_1,\hdots ) & \mapsto & x_0,
\end{array}
\end{alignedat}
\end{equation}
where $x_i \in S^{i}\omega^1_{X/S}.$
\end{parag}

\begin{proposition}\label{dualAkaza}
For a Hopf $\Ox_X$-algebra $\calA,$ let $\mathscr{Hom}_{\Ox_X}(\calA,\Ox_X)$ be the sheaf of $\Ox_X$-linear morphisms of $\Ox_X$-modules, equipped with the ring structure induced by the Hopf algebra structure on $\calA.$
Denote by $\widehat{S}^{\bullet}\calT_{X/S}$ and $\widehat{\Gamma}^{\bullet}\calT_{X/S}$ respectively the completions of $S^{\bullet}\calT_{X/S}$ and $\Gamma^{\bullet}\calT_{X/S}$ with respect to the ideals
$$S^{\ge 1}\calT_{X/S}:=\bigoplus_{n\ge 1}S^n\calT_{X/S}$$
and
$$\Gamma^{\ge 1}\calT_{X/S}:=\bigoplus_{n\ge 1}\Gamma^n\calT_{X/S}.$$
There exist canonical isomorphisms of $\Ox_X$-algebras
\begin{equation}\label{isodualAkaza}
s:\widehat{S}^{\bullet}\calT_{X/S}\xrightarrow{\sim} \mathscr{Hom}_{\Ox_X}(\Gamma^{\bullet}\omega^1_{X/S},\Ox_X),
\end{equation}
\begin{equation}\label{isodualAkaza2}
s':\widehat{\Gamma}^{\bullet}\calT_{X/S}\xrightarrow{\sim} \mathscr{Hom}_{\Ox_X}(S^{\bullet}\omega^1_{X/S},\Ox_X),
\end{equation}
such that the following diagram is commutative
$$
\begin{tikzcd}
\calT_{X/S} \ar[equal]{r} \ar{d} & \mathscr{Hom}_{\Ox_X}(\omega^1_{X/S},\Ox_X) \ar{d} \\
\widehat{S}^{\bullet}\calT_{X/S} \ar{r}{s} & \mathscr{Hom}_{\Ox_X}(\Gamma^{\bullet}\omega^1_{X/S},\Ox_X)
\end{tikzcd}
$$
where the left vertical arrow is the canonical one and the right vertical arrow sends $f:\omega^1_{X/S} \ra \Ox_X$ to the morphism induced by $f$ and the zero maps $0:\Gamma^{\ge 2}\omega^1_{X/S} \ra \Ox_X$ and $0:\Ox_X \ra \Ox_X.$
\end{proposition}

\begin{proof}
The $\Ox_X$-module $\omega^1_{X/S}$ is locally free of finite rank so, by (\cite{Ogus78} A10), there exists, for any integer $n\ge 0,$ canonical isomorphisms
$$s_n:S^{n}\calT_{X/S} \xrightarrow{\sim}\mathscr{Hom}_{\Ox_X}(\Gamma^n\omega^1_{X/S},\Ox_X),$$
$$S^{n}\omega^1_{X/S} \xrightarrow{\sim}\mathscr{Hom}_{\Ox_X}(\Gamma^{n}\calT_{X/S},\Ox_X).$$
By taking the dual of the second arrow, we obtain an isomorphism
$$s_n':\Gamma^{n}\calT_{X/S} \xrightarrow{\sim}\mathscr{Hom}_{\Ox_X}(S^{n}\omega^1_{X/S},\Ox_X).$$
Then
\begin{alignat*}{2}
\mathscr{Hom}_{\Ox_X}(\Gamma^{\bullet}\omega^1_{X/S},\Ox_X) & \xrightarrow{\sim} \prod_{n\ge 0}\mathscr{Hom}_{\Ox_X}(\Gamma^n\omega^1_{X/S},\Ox_X) \\
&\xrightarrow{\sim} \prod_{n\ge 0} S^n\calT_{X/S} \\
&\xrightarrow{\sim} \lim_{\substack{\longleftarrow \\ n\ge 1}} \bigoplus_{k<n}S^k\calT_{X/S} \\
&\xrightarrow{\sim} \lim_{\substack{\longleftarrow \\ n\ge 0}} \left (S^{\bullet}\calT_{X/S} / \bigoplus_{k\ge n}S^k\calT_{X/S} \right ) \\
&\xrightarrow{\sim} \w{S}^{\bullet}\calT_{X/S},
\end{alignat*}
where the first arrow is induced by the canonical embeddings 
$$\Gamma^n\omega^1_{X/S} \hookrightarrow \Gamma^{\bullet}\omega^1_{X/S},$$
the second arrow is induced by $(s_n^{-1})_{n\ge 0},$ the third arrow sends $(t_n)\in \prod_{n\ge 0}S^n\calT_{X/S}$ to the section $(y_n)$ defined, for any $n\ge 1,$ by
$$y_n=(t_0,t_1,\hdots,t_{n-1}),$$
and the two last arrows are the canonical ones. The fact that this isomorphism is an isomorphism of algebras follows from (\cite{OgusVol} 5.19). The isomorphism $s'$ is constructed in a similar way and the remaining assertion follows from the construction of $s.$
\end{proof}

\begin{proposition}\label{erafprop115}
Let $\calR=\calR_{\frakX,1}/p\calR_{\frakX,1}.$ There exists a canonical isomorphism of $\Ox_X$-algebras
\begin{equation}\label{eq1181}
S^{\bullet}\omega^1_{X/S} \xrightarrow{\sim} \calR,
\end{equation}
compatible with the Hopf algebra structures. In particular, there exists a canonical isomorphism of $\Ox_X$-algebras
\begin{equation}\label{eq1182}
\calR^{\vee} \xrightarrow{\sim} \widehat{\Gamma}^{\bullet}\calT_{X/S}.
\end{equation}
In addition, under the hypothesis \ref{loccoord}, \eqref{eq1181} sends $\op{dlog}m_i$ to $\eta_{i(r)}.$
\end{proposition}

\begin{proof}
Let $g:R_{\frakX,1} \ra \frakY$ be the canonical morphism, $\Delta:\frakX \ra \frakY$ the exact diagonal immersion and $\calI$ the ideal of $\Delta.$ The morphism of topoi $g:\left (R_{\frakX,1} \right )_{\text{zar}} \ra \frakY_{\text{zar}}$ factors as
$$
\begin{tikzcd}
\left (R_{\frakX,1} \right )_{\text{zar}} \ar{r}{\omega} \ar{d}{g} & \frakX_{\text{zar}} \ar{dl}{\Delta} \\
\frakY_{\text{zar}} &
\end{tikzcd}
$$
The morphism
$$
g^{-1}\calI \ra \Ox_{R_{\frakX,1}},\ x\mapsto \frac{x}{p}
$$
induces
$$\Delta^{-1}\calI \ra \omega_*\Ox_{R_{\frakX,1}} = \calR_{\frakX,1}.$$
Since $\calR$ is of characteristic $p,$ the composition with the canonical projection $\calR_{\frakX,1}\ra \calR,$
$$\Delta^{-1}\calI \ra \calR_{\frakX,1} \ra \calR,$$
vanishes on $\Delta^{-1}\calI^2.$ By (\ref{isoomega1}), we get an $\Ox_X$-linear morphism
$$\omega^1_{X/S} \ra \calR$$
and then a morphism of $\Ox_X$-algebras
$$S^{\bullet}\omega^1_{X/S} \ra \calR.$$
It is an isomorphism by \eqref{isoR}. The compatibility with the Hopf algebra structures comes from \ref{HopfR}. By taking the dual and \ref{dualAkaza}, we get an isomorphism
$$\calR^{\vee} \xrightarrow{\sim} \widehat{\Gamma}^{\bullet}\calT_{X/S}.$$
\end{proof}

\begin{parag}\label{Q'recall}
Recall the monoid $Q'$ and the morphisms $F_{Q/P}:Q'\ra Q,\ (x,t)\mapsto px+\theta^{gp}(t)$ and $\pi_{Q/P}:Q\ra Q',\ x\mapsto (x,0)$ defined in \ref{cor411}.
\end{parag}

\begin{proposition}\label{groupprop17}
Keep the notation \ref{Q'recall}. Suppose that there exists a morphism $F:(\frakX,Q) \ra (\frakX',Q')$ of framed logarithmic formal schemes, over $(\frakS,P),$ lifting the exact relative Frobenius $F_1:X \ra X'$ and such that $\frakX'\ra \frakS$ is log smooth. Let $G:\frakY \ra \frakY':=\frakX'\times_{\frakS,[Q']}^{\op{log}}\frakX'$ be the morphism induced by $F.$ Then the projection
$$\frakY\times_{\frakY'}R_{\frakX',1} \ra \frakY$$
factors through $Q_{\frakX}.$
\end{proposition}

\begin{proof}
Let $R'_1$ and $Y'$ be the special fibers of $R_{\frakX',1}$ and $\frakY'$ respectively\footnote{Not to confuse $Y'$ with the logarithmic scheme appearing in the exact relative Frobenius of $Y/S.$}.
By \ref{Glogflat}, $\frakY\times_{\frakY'}R_{\frakX',1} \ra R_{\frakX',1}$ is log flat. By \ref{locdescRQ}, $R_{\frakX,1} \ra \frakX$ is flat and strict hence log flat. We deduce that $\frakY\times_{\frakY'}R_{\frakX',1}$ is log flat over $\frakX$ and then over $\op{Spf}W.$ Then, by \ref{Wflat}, $\frakY\times_{\frakY'}R_{\frakX',1}$ is flat over $\op{Spf}W.$ To apply the universal property of $Q_{\frakX}$ (\ref{erafdil2}), it remains to prove that
$$\underline{Y\times_{Y'}R'_1} \ra Y\times_{Y'}R'_1 \ra Y$$
factors through the exact diagonal immersion $X \ra Y.$ For that it is sufficient to prove that
$$Y\times_{Y'}R'_1 \ra Y \xrightarrow{F_Y} Y$$
factors through $X \ra Y.$ We may work étale locally on $X$ and suppose that the frames $X\ra [Q]$ and $X' \ra [Q']$ lift to charts $X\ra A_1[Q]$ and $X' \ra A_1[Q'].$ Let $\widetilde{Q}$ and $\widetilde{Q}'$ be the inverse images of $Q$ and $Q'$ by the addition maps
$$Q^{gp} \oplus Q^{gp} \ra Q^{gp},\ Q'^{gp} \oplus Q'^{gp} \ra Q'^{gp}$$
respectively. By \ref{prop49} and \ref{rep1}, we have
$$
Y=X\times_{S,[Q]}^{\op{log}}X=\left (X\times_S^{\op{log}}X\right )\times_{A_1[Q\oplus Q]}^{\op{log}}A_1[\widetilde{Q}],
$$
$$
Y'=X'\times_{S,[Q']}^{\op{log}}X'=\left (X'\times_S^{\op{log}}X'\right )\times_{A_1[Q'\oplus Q']}^{\op{log}}A_1[\widetilde{Q}'].
$$
The absolute Frobenius $F_Y$ is then induced by the absolute Frobenius morphisms $F_X$ and $F_{A_1[\widetilde{Q}]}.$
Recall that the morphism $\pi:X' \ra X$ \eqref{diag51} is underlying a morphism of framed logarithmic schemes \eqref{cor411} and satisfies $\pi \circ F_1=F_X.$ The morphism $\pi$ then induces a morphism $\rho:Y' \ra Y.$ The commutative diagram
$$
\begin{tikzcd}
X \ar{r}{F_1} \ar[swap]{dr}{F_X} & X' \ar{d}{\pi}  \\
 & X 
\end{tikzcd}
$$
imply that $F_Y=\rho \circ G.$ Let $\Delta':X' \ra Y'$ and $\Delta:X\ra Y$ be the exact diagonal immersions and $p_1,p_2:Y \ra X$ the canonical projections. The compositions
\begin{alignat*}{1}
X' \xrightarrow{\Delta'} Y' \xrightarrow{\rho} Y \xrightarrow{p_i} X \\
X' \xrightarrow{\pi} X \xrightarrow{\Delta} Y \xrightarrow{p_i} X
\end{alignat*}
are both equal to $\pi$ for $i\in \{1,2\}.$ This implies that $\rho \circ \Delta'=\Delta \circ \pi.$
The result then follows from the commutative diagram
$$
\begin{tikzcd}
Y\times_{Y'}R'_1 \ar{r} \ar{d} & Y \ar{r}{F_Y} \ar{d}{G} & Y \\
R'_1 \ar{r} \ar{dr} & Y' \ar[swap]{ru}{\rho} & X \ar[swap]{u}{\Delta} \\
 & X' \ar{u}{\Delta'} \ar{ur}{\pi}, &
\end{tikzcd}
$$
where the upper left square is cartesian.
\end{proof}

\begin{parag}\label{erafparag121}
Keep the assumptions of \ref{groupprop17}. Let $Q_1$ be the special fiber of $Q_{\frakX}.$ By \ref{erafprop13}, we have a canonical morphism of logarithmic schemes $\omega':R'_1 \ra X'.$ Let $\omega:Q_1 \ra X$ be the map of topological spaces coming from the groupoid structure on $Q_1.$ We have a commutative diagram
$$
\begin{tikzcd}
Q_1\ar{rr} & & Q_{\frakX}\ar{d}{g} \\
X\times_{X'}R'_1 \ar{u}{v} \ar{r}{\delta} \ar{urr}{\widetilde{v}} & \frakY \times_{\frakY'}R_{\frakX',1} \ar{ur}{\chi} \ar{r}{\alpha} \ar{d}{\beta} & \frakY \ar{d}{G} \\
 & R_{\frakX',1} \ar{r}{h} & \frakY',
\end{tikzcd}
$$
where the lower square is cartesian, $g$ and $h$ are the canonical morphisms, $\chi:\frakY \times_{\frakY'}R_{\frakX',1} \ra Q_{\frakX}$ is given in \ref{groupprop17} and $\delta$ is induced by $X \ra \frakX \xrightarrow{\Delta}\frakY,$ $X' \ra \frakX' \xrightarrow{\Delta'} \frakY'$ and $R_1' \ra R_{\frakX',1}.$
Let us provide local descriptions for $\widetilde{v}$ and $v.$ Let $p_1,p_2:\frakY \ra \frakX$ be the canonical projections.
Consider the hypothesis \ref{loccoord}. Then, by \ref{lemcalc}, there exists $b\in \Ox_{\frakX}$ such that
$$G^{\#}\left ( \widetilde{\eta}_i'\right ) = \left (\widetilde{\eta}_i^p+\sum_{k=1}^{p-1}\begin{pmatrix}p\\k \end{pmatrix}\widetilde{\eta}_i^k+1 \right )\frac{1+pp_2^{\#}(b)}{1+pp_1^{\#}(b)}-1.$$
Set
$$
c=\frac{p_2^{\#}(b)-p_1^{\#}(b)}{1+pp_1^{\#}(b)},
$$
so that
$$
1+pc=\frac{1+pp_2^{\#}(b)}{1+pp_1^{\#}(b)}
$$
and
$$G^{\#}\left ( \widetilde{\eta}_i'\right ) = \left (\widetilde{\eta}_i^p+\sum_{k=1}^{p-1}\begin{pmatrix}p\\k \end{pmatrix}\widetilde{\eta}_i^k+1 \right )(1+pc)-1.$$
Then
\begin{alignat}{2}\label{eq12301}
\widetilde{\eta}_i^p &= G^{\#}(\widetilde{\eta}_i')-\sum_{k=1}^{p-1}\begin{pmatrix}p\\k \end{pmatrix}\widetilde{\eta}_i^k -pc\widetilde{\eta}_i^p-pc \left (\sum_{k=1}^{p-1}\begin{pmatrix}p\\k \end{pmatrix}\widetilde{\eta}_i^k\right )-pc.
\end{alignat}
Applying $g^{\#}$ to \eqref{eq12301}, since $g^{\#}(\widetilde{\eta}_i)^p\in p\Ox_{Q_{\frakX}},$ we get $g^{\#}G^{\#}(\widetilde{\eta}_i') \in p\Ox_{Q_{\frakX}}.$ Using the flatness of $Q_{\frakX}$ over $\Z_p,$ we get
\begin{equation}\label{eq12302}
\frac{g^{\#}(\widetilde{\eta}_i)^p}{p} = \frac{g^{\#}G^{\#}(\widetilde{\eta}_i')}{p}-\sum_{k=1}^{p-1}\frac{(p-1)!}{k!(p-k)!}g^{\#}(\widetilde{\eta}_i)^k -g^{\#}(c\widetilde{\eta}_i^p)- \left (\sum_{k=1}^{p-1}\begin{pmatrix}p\\k \end{pmatrix}g^{\#}(c\widetilde{\eta}_i^k)\right )-g^{\#}(c).
\end{equation}
Since $\widetilde{\eta}_i$ and $c=\frac{p_2^{\#}(b)-p_1^{\#}(b)}{1+pp_1^{\#}(b)}$ belong to the ideal of $\Delta:\frakX \ra \frakY,$ we get
\begin{alignat}{2}
\widetilde{v}^{\#}g^{\#}(\widetilde{\eta}_i) & =\delta ^{\#}\alpha ^{\#}(\widetilde{\eta}_i)=0, \label{eq12303}\\
\widetilde{v}^{\#}g^{\#}(c) & =\delta ^{\#}\alpha ^{\#}(c)=0. \label{eq12304}
\end{alignat}
We also have
$$
\chi^{\#}g^{\#}G^{\#}(\widetilde{\eta}_i') = \beta^{\#}h^{\#}(\widetilde{\eta}_i')=p\beta^{\#}\left (\frac{\widetilde{\eta}_i'}{p} \right )=p\left (1\otimes \frac{\widetilde{\eta}_i'}{p} \right ).
$$
The morphism $G:\frakY \ra \frakY'$ is log flat \eqref{Glogflat} and hence so is $\beta.$ Since $R_{\frakX',1}$ is log flat over $\op{Spf}W$ \eqref{Qlogflat}, we deduce that $\frakY \times_{\frakY'}R_{\frakX',1}$ is log flat hence flat over $\op{Spf}W.$ We then get
$$
\chi^{\#}\left ( \frac{g^{\#}G^{\#}(\widetilde{\eta}_i')}{p} \right ) = \beta^{\#}\left (\frac{h^{ \#}(\widetilde{\eta}_i')}{p} \right )=1\otimes \frac{\widetilde{\eta}_i'}{p}.
$$
It follows that
\begin{equation}\label{eq12305}
\widetilde{v}^{\#}\left (\frac{g^{\#}G^{\#}(\widetilde{\eta}_i')}{p}\right )=1\otimes \eta_{i(r)}'.
\end{equation}
By \eqref{eq12302}, \eqref{eq12303}, \eqref{eq12304} and \eqref{eq12305}, dropping the notation $g^{\#},$ we get
$$\widetilde{v}^{\#}\left (\frac{\widetilde{\eta}_i^p}{p} \right )=1\otimes \eta_{i(r)}',\ \widetilde{v}^{\#}(\widetilde{\eta}_i)=0,$$
and so
\begin{equation}\label{eq1211v}
v^{\#}\left (\eta_{i(q)} \right ) = 1\otimes \eta_{i(r)}',\ v^{\#}(\eta_i)=0.
\end{equation}
This proves that $v$ is independant of the choice of the lifting $F,$ so we have a global morphism
$$v:X\times_{X'}R'_1 \ra Q_1.$$
Let $\calR'_1=\omega'_*\Ox_{R'_1}$ and $\calQ_1=\omega_*\Ox_{Q_1}$ be the respective Hopf algebras defined by $R'_1$ and $Q_1.$ (\ref{paragomega}). 
Consider the cartesian square
$$
\begin{tikzcd}
X\times_{X'}R'_1 \ar{r}{\alpha'} \ar{d}{\beta'} & X \ar{d}{F_1} \\
R'_1 \ar{r}{\omega'} & X'.
\end{tikzcd}
$$
Since $F_1$ and $\omega'$ are affine, there exists a canonical isomorphism
\begin{equation}\label{eq12307}
F_1^*\calR'_1=F_1^*\omega'_*\Ox_{R'_1} \xrightarrow{\sim} \alpha'_*\Ox_{X\times_{X'}R'_1}.
\end{equation}
We have a diagram of topological spaces (recall the notation \ref{Not8}):
$$
\begin{tikzcd}
 & & \left|Q_1\right| \ar[swap]{dl}{\omega} \ar{d} \\
\left|X\times_{X'}R'_1\right| \ar{r}{\alpha'} \ar[bend right=-30]{urr}{v} & \left|X\right| \ar{r}{\Delta} & \left|Y\right|.
\end{tikzcd}
$$
Since the outer diagram is commutative and $\Delta$ is injective as a map of sets, we get $\omega \circ v=\alpha'$ as maps of topological spaces.
Applying $\omega_*$ to $v^{\#}:\Ox_{Q_1} \ra v_*\Ox_{X\times_{X'}R'_1},$ by \eqref{eq12307}, we get a morphism that we abusively denote by $v:$
\begin{equation}
v:\calQ_1 \ra F_1^*\calR'_1.
\end{equation}
We prove that $v$ is a morphism of Hopf algebras. Let $(\delta_q,\pi_q,\sigma_q)$ and $(\delta_r,\pi_r,\sigma)$ be the Hopf algebra structures of $\calQ_1$ and $\calR_1'$ respectively. We check that the diagram
\begin{equation}\label{diag12309}
\begin{tikzcd}
\calQ_1 \ar{r}{\delta_q} \ar[swap]{d}{v} & \calQ_1 \otimes_{\Ox_X} \calQ_1 \ar{d}{v\otimes v} \\
F_1^*\calR_1' \ar{r}{F_1^*\delta_r} & F_1^*\calR_1' \otimes_{\Ox_X} F_1^*\calR_1'.
\end{tikzcd}
\end{equation}
By \ref{HopfQ}, we have
\begin{alignat*}{2}
\delta_q\left (\eta_{i(q)} \right ) &= 1\otimes \eta_{i(q)}+\sum_{0<b+c<p} \frac{(-1)^{b+c}}{b+c}\begin{pmatrix}b+c \\ b \end{pmatrix} \eta_i^{b+c}\otimes \eta_i^{p-b} + \eta_{i(q)}\otimes 1, \\
\delta_q(\eta_i) &= 1\otimes \eta_i+\eta_i\otimes 1+\eta_i\otimes \eta_i.
\end{alignat*}
By \eqref{eq1211v}, we get
\begin{alignat*}{2}
(v\otimes v)\circ \delta_q(\eta_{i(q)}) &= 1\otimes F_1^*\eta_{i(r)}' + F_1^*\eta_{i(r)}' \otimes 1, \\
(v\otimes v) \circ \delta_q(\eta_i) &= 0.
\end{alignat*}
By \eqref{eq1211v} and \ref{HopfR}, we have
\begin{alignat*}{2}
(F_1^*\delta_r)\circ v(\eta_{i(q)}) &= 1\otimes F_1^*\eta_{i(r)}' + F_1^*\eta_{i(r)}'\otimes 1, \\
(F_1^*\delta_r)\circ v(\eta_{i}) &= 0.
\end{alignat*}
This proves the commutativity of \eqref{diag12309}. Similarly, we prove the commutativity of the diagrams
$$
\begin{tikzcd}
\calQ_1 \ar{r}{\sigma_q} \ar[swap]{d}{v} & \calQ_1 \ar{d}{v} & \calQ_1 \ar[swap]{d}{v} \ar{r}{\pi_q} & \Ox_X \\
F_1^*\calR_1' \ar{r}{F_1^*\sigma} & F_1^*\calR_1' & F_1^*\calR_1' \ar[swap]{ur}{F_1^*\pi_r} & 
\end{tikzcd}
$$
This proves that $v$ is a morphism of Hopf algebras.
By taking the dual of $v,$ we obtain a morphism of $\Ox_X$-algebras
\begin{equation}\label{checkv}
\check{v}:F_1^*\calR_1'^{\vee} \ra \calQ_1^{\vee}.
\end{equation}
\end{parag}

\begin{lemma}\label{Kokoprop1234}
Embed $S^{\bullet}\calT_{X'/S}$ into $\calD_{X/S}$ via the $p$-curvature map $\psi$ \eqref{eqDoma1}.
There exists a canonical isomorphism of $\Ox_X$-algebras:
\begin{equation}
\calD_{X/S} \otimes_{S^{\bullet}\calT_{X'/S}} \widehat{S}^{\bullet}\calT_{X'/S} \xrightarrow{\sim} \widehat{\calD}_{X/S},
\end{equation}
where $\widehat{\calD}_{X/S}$ is defined in \ref{def14logop}.
\end{lemma}

\begin{proof}
By construction of the morphism $\widehat{\psi}:\widehat{S}^{\bullet}\calT_{X'/S} \ra \widehat{\calD}_{X/S}$ \eqref{eq13191}, the diagram
$$
\begin{tikzcd}
S^{\bullet}\calT_{X'/S} \ar{r}{\psi} \ar{d} & \calD_{X/S} \ar{d} \\ 
\widehat{S}^{\bullet}\calT_{X'/S} \ar{r}{\widehat{\psi}} & \widehat{\calD}_{X/S} 
\end{tikzcd}
$$
is commutative.
We deduce a morphism
\begin{equation}\label{ren8}
\calD_{X/S} \otimes_{S^{\bullet}\calT_{X'/S}} \widehat{S}^{\bullet}\calT_{X'/S} \ra \widehat{\calD}_{X/S}.
\end{equation}
To prove it is an isomorphism, we may work étale locally on $X$ and suppose that the hypothesis \ref{loccoord} is satisfied. We keep the same notations of \ref{erafparag121}. Let $(\partial_1',\hdots,\partial_d')\in \left (\calT_{X'/S}\right )^d$ be the dual basis of $(\op{dlog}m_1',\hdots,\op{dlog}m_1')$ and consider the differential operators $\partial_I$ defined in \ref{P2}. For $J\in \N^d,$ set
$$
\partial'^J=\prod_{j=1}^d \partial_j'^{J_j}\in S^{\bullet}\calT_{X'/S}.
$$
A local section of $\calD_{X/S} \otimes_{S^{\bullet}\calT_{X'/S}} \widehat{S}^{\bullet}\calT_{X'/S}$ is uniquely written as an infinite sum
\begin{equation}\label{eq12311}
\sum_{\substack{I\in \llbracket 0,p-1\rrbracket^d \\ J\in \N^d}} a_{I,J}\partial_I\otimes \partial'^J,
\end{equation}
with $a_{I,J}\in \Ox_X.$
By definition, the morphism $\widehat{S}^{\bullet}\calT_{X'/S} \ra \widehat{\calD}_{X/S}$ \eqref{eq13191} is induced by the $p$-curvature map \eqref{eqDoma1}. By \eqref{era2Koko11}, for all $J\in \N^d,$
\begin{alignat*}{2}
\widehat{\psi}\left (\partial'^J\right ) &= \partial_{pJ}.
\end{alignat*}
It follows that \eqref{eq12311} is sent by \eqref{ren8} to
$$
\sum_{\substack{I\in \llbracket 0,p-1\rrbracket^d \\ J\in \N^d}}a_{I,J}\partial_{I}\circ \partial_ {pJ},
$$
which is equal, by \ref{lemiran1}, to
$$
\sum_{\substack{I\in \llbracket 0,p-1\rrbracket^d \\ J\in \N^d}}a_{I,J}\partial_{I+pJ}.
$$
The result follows.
\end{proof}

\begin{proposition}\label{prop1232}
Consider the morphisms $\check{u}:\widehat{\calD}_{X/S} \ra \calQ_1^{\vee}$ \eqref{checku}, $\check{v}:F_1^*\calR_1'^{\vee} \ra \calQ_1^{\vee}$ \eqref{checkv}, the antipode morphism $\sigma:\calR'_1 \ra \calR'_1$ \eqref{HopfR} and its dual $\check{\sigma}:\calR_1'^{\vee} \ra \calR_1'^{\vee}.$ Recall that we identify the small étale sites of $X$ and $X'$ via $F_1:X\ra X'.$
Consider the canonical morphisms $\calR_1'^{\vee}\ra F_1^*\calR_1'^{\vee},\ x\mapsto 1\otimes x$ and $\widehat{S}^{\bullet}\calT_{X'/S} \ra \widehat{\Gamma}^{\bullet}\calT_{X'/S}$ and the morphism $\widehat{\psi}:\widehat{S}^{\bullet}\calT_{X'/S} \ra \widehat{\calD}_{X/S}$ \eqref{eq13191}.
We identify $\widehat{\Gamma}^{\bullet}\calT_{X'/S}$ with $\calR_1'^{\vee}$ via the isomorphism $\xi:\widehat{\Gamma}^{\bullet}\calT_{X'/S} \xrightarrow{\sim} \calR_1'^{\vee}$ \eqref{eq1182}.
Then, the diagram
$$
\begin{tikzcd}
\widehat{S}^{\bullet}\calT_{X'/S} \ar{d} \ar{rrrr}{\widehat{\psi}} & & & & \widehat{\calD}_{X/S} \ar{d}{\check{u}} \\
\widehat{\Gamma}^{\bullet}\calT_{X'/S} \ar{r}{\xi} & \calR_1'^{\vee} \ar{r}{\check{\sigma}} & \calR_1'^{\vee} \ar{r} & F_1^*\calR_1'^{\vee} \ar{r}{\check{v}} & \calQ_1^{\vee}
\end{tikzcd}
$$
is commutative.
\end{proposition}

\begin{proof}
We may suppose that the hypothesis \ref{loccoord} is satisfied. We keep the same notations of \ref{erafparag121}. Let $(\partial_1',\hdots,\partial_d')\in \left (\calT_{X'/S} \right )^d$ be the dual basis of $(\op{dlog}m_1',\hdots,\op{dlog}m_d').$ We consider the differential operators $\partial_I,$ local sections of $\w{\calD}_{X/S},$ defined from the local coordinates $m_1,\hdots,m_d$ as in \ref{P2}. For $1\le i\le d,$ denote by $\epsilon_i$ the multi-index of $\N^d$ whose all coefficients are zero except for the $i$th which is equal to $1.$ For any $I=(I_1,\hdots,I_d)\in \N^d,$ set
$$\eta^{[I]}=\prod_{i=1}^d\eta_i^{[I_i]} \in \calP_0,\ \eta_{(r)}'^I=\prod_{i=1}^d\eta_{i(r)}'^{I_i} \in \calR'_1.$$
By \ref{propfree}, $\left (\eta^{[I]}\right )_{I\in \N^d}$ is a basis for the $\Ox_X$-module $\calP_0,$ where $\calP_0$ is considered as an $\Ox_X$-module via the first projection $P_0 \ra X.$ By definition and \eqref{era2Doma11}, $\check{u}(\partial_i')=\partial_{p\epsilon_i}\circ u.$ Then, by \eqref{ueta}, we have
\begin{alignat*}{2}
\left (\check{u}(\partial_i')\right )(\eta_j) &= \partial_{p\epsilon_i}(\eta_j)=0.\\
\left (\check{u}(\partial_i')\right )(\eta_{j(q)}) &= \partial_{p\epsilon_i}\left (-\eta_j^{[p]}\right )=-\delta_{ij}.
\end{alignat*}
By \ref{locdescRQ}, $F_1^*\calR_1' =\Ox_X\left [\eta_{1(r)}',\hdots,\eta_{d(r)}'\right ].$ Let $(\varphi_I)_{I\in \N^d}$ be the dual basis of $(\eta_{(r)}'^I)_{I\in \N^d}.$
By \ref{erafprop115}, the local section $\partial_i'\in F_1^*\widehat{\Gamma}^{\bullet}\calT_{X'/S}$ corresponds, by $F_1^*\xi,$ to $\varphi_{\epsilon_i} \in F_1^*\calR_1'^{\vee}.$ It follows that $(\check{v}\circ\check{\sigma})(\partial_i')=\varphi_{\epsilon_i}\circ (F_1^*\sigma) \circ v.$ Then, by \eqref{eq1211v},
\begin{alignat*}{2}
\left (\check{v}\circ \check{\sigma}(\partial_i')\right )(\eta_j) &= \varphi_{\epsilon_i}\circ (F_1^*\sigma)\circ v(\eta_j)=0.\\
\left (\check{v}\circ \check{\sigma}(\partial_i')\right )(\eta_{j(q)}) &= \varphi_{\epsilon_i}\circ (F_1^*\sigma)\circ v\left (\eta_{j(q)}\right )=\varphi_{\epsilon_i}(-\eta_{j(r)}')=-\delta_{ij}.
\end{alignat*}
This proves that
$$\check{u}(\partial_i')=\check{v}\circ \check{\sigma}(\partial_i')$$
hence the proposition.
\end{proof}

\begin{proposition}
Keep the notation of \ref{prop1232}.
\begin{enumerate}
\item For any $\partial \in \widehat{\calD}_{X/S}$ and $D'\in \widehat{\Gamma}^{\bullet}\calT_{X'/S},$ we have the following equality in $\calQ_1^{\vee}:$
$$\check{u}(\partial)\cdot (\check{v}\circ \check{\sigma})(D')=(\check{v}\circ \check{\sigma})(D')\cdot \check{u}(\partial),$$
where the algebra structure on $\calQ_1^{\vee}$ comes from the Hopf algebra structure on $\calQ_1.$
\item The map
\begin{equation}\label{isoDQ}
\xi:\begin{array}[t]{clc}
\widehat{\calD}_{X/S} \otimes_{\widehat{S}^{\bullet}\calT_{X'/S}} \widehat{\Gamma}^{\bullet}\calT_{X'/S} & \ra & \calQ_1^{\vee} \\
x\otimes y & \mapsto & \check{u}(x) \cdot \check{v}(y),
\end{array}
\end{equation}
is a well-defined isomorphism of $\Ox_X$-algebras.
\end{enumerate}
\end{proposition}

\begin{proof}
For the first assertion, it is enough to prove the equality for $\partial=\partial_{\alpha}$ and $D'=\partial'^{[\beta]}:=\prod_{i=1}^d\partial_i'^{[\beta_i]},$ for $\alpha,\beta\in \N^d.$
Just as in the proof of \ref{prop1232}, we have
$$
(\check{v} \circ \check{\sigma})(\partial'^{[\beta]})=\varphi_{\beta} \circ (F_1^*\sigma) \circ v
$$
and
$$
\check{u}(\partial_{\alpha})=\partial_{\alpha}\circ u.
$$
By definition of the ring structure on $\calQ_1^{\vee},$ the product $\check{u}(\partial_{\alpha}) \cdot (\check{v} \circ \check{\sigma})(\partial'^{[\beta]})$ is then equal to the composition
\begin{equation*}
\calQ_1 \xrightarrow{\delta} \calQ_1\otimes_{\Ox_X}\calQ_1 \xrightarrow{\op{Id}\otimes v} \calQ_1\otimes_{\Ox_X}F_1^*\calR_1' \xrightarrow{\op{Id}\otimes F_1^*\sigma} \calQ_1\otimes_{\Ox_X}F_1^*\calR_1' \xrightarrow{\op{Id}\otimes \varphi_{\beta}} \calQ_1 \xrightarrow{u} \calP_0 \xrightarrow{\partial_{\alpha}} \Ox_X,
\end{equation*}
where $\delta$ is given in \ref{HopfQ}.
This is equal to the composition
\begin{equation*}
\calQ_1 \xrightarrow{\delta} \calQ_1\otimes_{\Ox_X}\calQ_1 \xrightarrow{u\otimes v} \calP_0\otimes_{\Ox_X} F_1^*\calR_1' \xrightarrow{\partial_{\alpha}\otimes (\varphi_{\beta}\circ F_1^*\sigma)}  \Ox_X.
\end{equation*}
For $J,K\in \N^d,$ set
$$
\eta_{(q)}^J=\prod_{j=1}^d\eta_{j(q)}^{J_j} \in \calQ_1,\ 
\eta^K=\prod_{k=1}^d\eta_k^{K_k} \in \calQ_1.
$$
By \eqref{isoQ}, the $\Ox_X$-module $\calQ_1$ has a basis
$$\left (\eta^K\eta_{(q)}^J\right )_{K\in \llbracket 0,p-1 \rrbracket^d,\ J\in \N^d}.$$
Let $K\in \llbracket 0,p-1 \rrbracket^d$ and $J\in \N^d.$ By \ref{HopfQ},
\begin{equation}\label{12332}
\begin{alignedat}{2}
\delta(\eta^K\eta_{(q)}^J) =& \prod_{1\le k\le d} \delta(\eta_k)^{K_k} \prod_{1\le j\le d}\delta(\eta_{j(q)})^{J_j} \\
=& \prod_{1\le k\le d} (1\otimes\eta_k+\eta_k\otimes 1+\eta_k\otimes\eta_k)^{K_k} \\
&\times \prod_{1\le j\le d} \left (1\otimes \eta_{j(q)}+\sum_{0<b+c<p} \frac{(-1)^{b+c}}{b+c}\begin{pmatrix}b+c \\ b \end{pmatrix} \eta_j^{b+c}\otimes \eta_j^{p-b} + \eta_{j(q)}\otimes 1 \right )^{J_j}.
\end{alignedat}
\end{equation}
Recall the notation \ref{Not9}. Developing the right-hand side of \eqref{12332} with the binomial formula and by \eqref{eq1211v} and \eqref{ueta}, we get
\begin{alignat*}{2}
(u\otimes v)\circ \delta(\eta^K\eta_{(q)}^J) =& \left (u \left (\eta^K \right )\otimes 1\right ) \left (\sum_{L\le J} \begin{pmatrix}J\\L \end{pmatrix} u \left (\eta_{(q)}^L\right )\otimes v\left (\eta_{(q)}^{J-L} \right ) \right ) \\
=& \left (\eta^K\otimes 1\right )\left (\sum_{L\le J} \begin{pmatrix} J\\ L \end{pmatrix} \left (-\eta^{[p]} \right )^L\otimes \eta_{(r)}'^{J-L}\right ) \\
=& \sum_{L\le J} \begin{pmatrix} J\\ L \end{pmatrix} \left (\eta^K\left (-\eta^{[p]} \right )^L\right )\otimes \eta_{(r)}'^{J-L}.
\end{alignat*}
We then get
$$
(\op{Id}_{\calP_0}\otimes F_1^*\sigma) \circ (u\otimes v)\circ \delta(\eta^K\eta_{(q)}^J) = (-1)^{|J|}\sum_{L\le J} \begin{pmatrix} J\\ L \end{pmatrix} \left (\eta^K\left (\eta^{[p]} \right )^L\right )\otimes \eta_{(r)}'^{J-L},
$$
and so
\begin{equation}
\left ((\check{v} \circ \check{\sigma})(\partial'^{[\beta]}) \cdot \check{u}(\partial_{\alpha}) \right ) (\eta^K\eta_{(q)}^J) = (-1)^{|J|} \begin{pmatrix} J \\ \beta\end{pmatrix} \partial_{\alpha}\left (\eta^K\left (\eta^{[p]} \right )^{J-\beta} \right ).
\end{equation}
By basic PD structure properties,
$$\left (\eta^{[p]} \right )^{J-\beta}=\frac{(p(J-\beta))!}{p!^{J-\beta}}\eta^{[p(J-\beta)]}.$$
Note that, for any integer $n\ge 0,$
$$v_p\left ( (pn)! \right )=\sum_{k\ge 1} \left \lfloor \frac{pn}{p^k} \right \rfloor \ge n,$$
so $\frac{(p(J-\beta))!}{p!^{J-\beta}}\in \N.$
It follows that
\begin{equation}\label{1222v}
 (\check{v} \circ \check{\sigma})(\partial'^{[\beta]}) \cdot \left ( \check{u}(\partial_{\alpha}) \right ) (\eta^K\eta_{(q)}^J) = (-1)^{|J|} \begin{pmatrix} J \\ \beta\end{pmatrix} \frac{(K+p(J-\beta))!}{p!^{J-\beta}} \partial_{\alpha}\left (\eta^{[K+p(J-\beta)]} \right ).
\end{equation}
The product $\left ((\check{v} \circ \check{\sigma})(\partial'^{[\beta]})\right ) \cdot  \check{u}(\partial_{\alpha})$ is equal to the composition
$$
\calQ_1 \xrightarrow{\delta} \calQ_1\otimes_{\Ox_X}\calQ_1 \xrightarrow{\op{Id} \otimes u} \calQ_1\otimes_{\Ox_X}\calP_0 \xrightarrow{\op{Id} \otimes \partial_{\alpha}} \calQ_1 \xrightarrow{v} F_1^*\calR_1' \xrightarrow{F_1^*\sigma} F_1^*\calR_1' \xrightarrow{\varphi_{\beta}} \Ox_X.
$$
This is equal to the composition
\begin{equation*}
\calQ_1 \xrightarrow{\delta} \calQ_1\otimes_{\Ox_X}\calQ_1 \xrightarrow{v\otimes u} \calP_0\otimes_{\Ox_X} F_1^*\calR_1' \xrightarrow{(\varphi_{\beta}\circ F_1^*\sigma)\otimes \partial_{\alpha}}  \Ox_X.
\end{equation*}
A similar computation for $\left ((\check{v} \circ \check{\sigma})(\partial'^{[\beta]}) \right ) \circ \check{u}(\partial_{\alpha})(\eta^K\eta_{(q)}^J)$ proves that
$$(\check{v} \circ \check{\sigma})(\partial'^{[\beta]}) \cdot \check{u}(\partial_{\alpha})=\check{u}(\partial_{\alpha}) \cdot (\check{v} \circ \check{\sigma})(\partial'^{[\beta]}).$$
This concludes the proof of the first assertion, which in turn proves that $\xi$ \eqref{isoDQ} is well-defined. It remains to prove that it is an isomorphism. A local section of $\widehat{\calD}_{X/S} \otimes_{\widehat{S}^{\bullet}\calT_{X'/S}} \widehat{\Gamma}^{\bullet}\calT_{X'/S}$ can be uniquely written as an infinite sum
$$
\sum_{\substack{\alpha \in \llbracket0,p-1\rrbracket^d \\ \beta\in \N^d}} c_{\alpha\beta} \partial_{\alpha}\otimes \partial'^{[\beta]},\ c_{\alpha \beta}\in \Ox_X.
$$
For $\alpha\in \llbracket 0,p-1 \rrbracket^d$ and $\beta,J,K\in \N^d,$ we have, by \eqref{1222v},
\begin{equation}\label{12315}
\xi(\partial_{\alpha}\otimes \partial'^{[\beta]})(\eta^K\eta_{(q)}^J)=(-1)^{|\beta|}\alpha! \delta_{\alpha K}\delta_{\beta J}.
\end{equation}
This finishes the proof.
\end{proof}

\begin{corollaire}
Let
\begin{equation}\label{eq12331bis}
\calD^{\gamma}_{X/S}=\calD_{X/S}\otimes_{S^{\bullet} \calT_{X'/S}} \widehat{\Gamma}^{\bullet}\calT_{X'/S}.
\end{equation}
There exists an isomorphism of $\Ox_X$-algebras
\begin{equation}\label{eq12331}
\calD^{\gamma}_{X/S} \xrightarrow{\sim} \calQ_1^{\vee}.
\end{equation}
In addition, under the assumption \ref{loccoord}, the local section $(-1)^{|J|} I! \partial_I \otimes \partial'^{[J]}$ of $\calD_{X/S}^{\gamma}$ corresponds to $\varphi_{I,J},$ for $I\in \llbracket 0,p-1 \rrbracket^d$ and $J\in \N^d,$ where $\left (\varphi_{I,J}\right )_{\substack{I\in \llbracket 0,p-1 \rrbracket^d \\ J\in \N^d}}$ is the dual basis of $\left (\eta^I\eta_{(q)}^J \right )_{\substack{I\in \llbracket 0,p-1 \rrbracket^d \\ J\in \N^d}}.$
\end{corollaire}

\begin{proof}
This results from \eqref{isoDQ}, \ref{Kokoprop1234} and \eqref{12315}.
\end{proof}

\begin{definition}\label{locPDnil}
Let $\mathfrak{I}$ be the PD-ideal of $\Gamma^{\bullet}\calT_{X/S}.$
\begin{enumerate}
\item A $\widehat{\Gamma}^{\bullet}\calT_{X/S}$-module $\calE$ is said to be \emph{locally PD-nilpotent} if any local section $x$ of $\calE$ is locally annihilated by $\mathfrak{I}^{[m]}$ for some positive integer $m.$
\item A $\calD_{X/S}^{\gamma}$-module $\calE$ is said to be \emph{locally PD-nilpotent} if it is so when considered as a $\w{\Gamma}^{\bullet}\calT_{X/S}$-module via the canonical morphism
$$\w{\Gamma}^{\bullet}\calT_{X/S} \ra \calD_{X/S}^{\gamma}.$$
\end{enumerate}
\end{definition}

\begin{theorem}\label{thm1237}
Let $\calR_1=\calR_{\frakX,1}/p\calR_{\frakX,1},$ $\calQ_1=\calQ_{\frakX}/p\calQ_{\frakX}$ and
$$\calD^{\gamma}_{X/S}=\calD_{X/S}\otimes_{S^{\bullet} \calT_{X'/S}} \widehat{\Gamma}^{\bullet}\calT_{X'/S}.$$
The following tensor categories are canonically equivalent:
\begin{enumerate}
\item The category of $\Ox_X$-modules equipped with a stratification relative to $\calR_1$ (resp. $\calQ_1$).
\item The category of locally PD-nilpotent $\widehat{\Gamma}^{\bullet}\calT_{X/S}$-modules (resp. locally PD-nilpotent $\calD^{\gamma}_{X/S}$-modules).
\end{enumerate}
\end{theorem}

\begin{proof}
The proof being similar to (\cite{Oyama} 1.2.10), we just outline the idea of the equivalences. 

Let $\calE$ be an $\Ox_X$-module equipped with an $\calR_1$-stratification
$$\epsilon:\calR_1\otimes_{\Ox_X}\calE \ra \calE \otimes_{\Ox_X}\calR_1.$$
Since we have an isomorphism of $\Ox_X$-algebras \eqref{eq1182}
$$\calR^{\vee}_1 \xrightarrow{\sim} \widehat{\Gamma}^{\bullet}\calT_{X/S},$$
it is sufficient to define a $\calR^{\vee}_1$-module structure on $\calE.$
This structure is obtained as follows: for local sections $x$ and $\varphi$ of $\calE$ and $\calR^{\vee}_1$ respectively, the action $\varphi \cdot x$ of $\varphi$ on $x$ is the image of $x$ by the composition
$$\begin{array}[t]{clclc}
\calE & \ra & \calE\otimes_{\Ox_X}\calR_1 & \xrightarrow{\op{Id}_{\calE}\otimes \varphi }  & \calE \\
x & \mapsto & \epsilon(1\otimes x) & &
\end{array}
$$
We check that the $\widehat{\Gamma}^{\bullet}\calT_{X/S}$-module structure obtained is PD-nilpotent by reducing to local coordinates: suppoe that the hypothesis \ref{loccoord} is satisfied. We have a basis $(\op{dlog}m_i)_{1\le i\le d}$ of $\omega^1_{X/S}.$ Let $(D_i)_{1\le i\le d}$ be its dual basis and consider $D_i$ as a section of $\widehat{\Gamma}^{\bullet}\calT_{X/S}$ via the canonical morphism $\calT_{X/S} \ra \widehat{\Gamma}^{\bullet}\calT_{X/S}.$ For $I=(I_1,\hdots,I_d)\in \N^d,$ set
$$\eta_{(r)}^I=\prod_{i=1}^d\eta_{i(r)}^{I_i} \in \calR_1.$$
By \ref{HopfR}, the family $\left (\eta_{(r)}^I\right )_{I\in \N^d}$ is a basis for the $\Ox_X$-module $\calR_1.$ Denote by $\left (\varphi_I \right )_{I\in \N^d}$ its dual base. By \ref{erafprop115}, the section $D_i\in \widehat{\Gamma}^{\bullet}\calT_{X/S}$ corresponds to $\varphi_{\epsilon_i} \in \calR_1^{\vee},$ where $\epsilon_i\in \N^d$ is such that all its coefficients are zero except for the $i$th which is equal to $1.$ Let $x$ be a local section of $\calE.$ There exists a positive integer $n$ and local sections $x_I\in \calE$ such that
$$
\epsilon(1\otimes x)=\sum_{\substack{I\in \N^d \\ |I|\le n}} x_I\otimes \eta_{(r)}^I \in \calE \otimes_{\Ox_X}\calR_1.
$$
It follows that, for $I\in \N^d$ such that $|I|>n,$
$$
\varphi_I\cdot x=0.
$$
Conversly, let $\calE$ be a locally PD-nilpotent $\widehat{\Gamma}^{\bullet}\calT_{X/S}$-module and
$$\alpha:\calR^{\vee}_1\otimes_{\Ox_X}\calE \xrightarrow{\sim} \widehat{\Gamma}^{\bullet}\calT_{X/S}\otimes_{\Ox_X}\calE \ra \calE$$
the structure morphism.
Consider the morphism
$$
\Lambda:\mathscr{Hom}_{\Ox_X}\left (\calE , \calE\otimes_{\Ox_X}\calR_1 \right ) \ra \mathscr{Hom}_{\Ox_X}\left (\calR^{\vee}_1\otimes_{\Ox_X}\calE , \calE \right )
$$
sending a morphism $f:\calE \ra \calE \otimes_{\Ox_X}\calR_1$ to
$$
\begin{array}[t]{clc}
\calR^{\vee}_1 \otimes_{\Ox_X} \calE & \ra & \calE \\
\varphi \otimes x & \mapsto & (\op{Id}_{\calE} \otimes \varphi)\circ f(x).
\end{array}
$$
The $\Ox_X$-module $\calR_1$ is locally free so $\Lambda$ is injective. Suppose that there exists an étale covering $(U_i \ra X)_{i\in I}$ and morphisms $\theta_i:\calE_{|U_i} \ra \calE_{|U_i}\otimes_{\Ox_{U_i}}\calR_{1|U_i}$ such that $\Lambda(\theta_i)=\alpha_{|U_i}$ for all $i\in I.$ The injectivity of $\Lambda$ implies that the morphisms $\theta_i$ glue together into a morphism $\theta:\calE \ra \calE\otimes_{\Ox_X}\calR_1$ such that $\Lambda(\theta)=\alpha.$ It is hence sufficient to construct, étale locally on $X,$ a morphism
$$\theta:\calE \ra \calE \otimes_{\Ox_X}\calR_1$$
sent, by $\Lambda,$ to $\alpha.$ We can thus suppose that the hypothesis \ref{loccoord} is satisfied. We have a basis $\left (\eta_{(r)}^I\right )_{I\in \N^d}$ of the $\Ox_X$-module $\calR_1.$ Let $(\alpha_I)_{i\in \N^d}$ be the dual basis. We set
$$
\theta:\begin{array}[t]{clc}
\calE & \ra & \calE\otimes_{\Ox_X}\calR_1 \\
x & \mapsto & \sum_{I\in \N^d}(\alpha_I\cdot x)\otimes \eta_{(r)}^I.
\end{array}
$$
Note that the sum is finite since $\calE$ is locally PD-nilpotent.

Let us now construct the second equivalence. Let $\calE$ be an $\Ox_X$-module equipped with a $\calQ_1$-stratification $\epsilon.$ By \eqref{eq12331}, it is sufficient to define a $\calQ_1^{\vee}$-module structure on $\calE.$ We define the action $\varphi\cdot x$ of $\varphi \in \calQ_1^{\vee}$ on $x\in \calE$ as the image of $x$ by the composition
$$
\begin{array}[t]{clclc}
\calE & \ra & \calE\otimes_{\Ox_X}\calQ_1 & \xrightarrow{\op{Id}_{\calE}\otimes \varphi }  & \calE \\
x & \mapsto & \epsilon(1\otimes x) & &
\end{array}
$$
To check that this $\calD_{X/S}^{\gamma}$-module structure is PD-nilpotent, we suppose that the hypothesis \ref{loccoord} is satisfied. For $I=(I_1,\hdots,I_d)\in \N^d,$ set
$$\eta_{(q)}^I=\prod_{i=1}^d\eta_{i(q)}^{I_i} \in \calQ_1,\ \eta^I=\prod_{i=1}^d\eta_{i}^{I_i} \in \calQ_1.$$
Then $\left (\eta^I\eta_{(q)}^J\right )_{\substack{I\in \llbracket 0,p-1 \rrbracket^d \\ J\in \N^d}}$ is a basis for the $\Ox_X$-module $\calQ_1.$ Denote by $\left (\varphi_{I,J}\right )$ its dual basis. Let $(\partial'_1,\hdots,\partial'_d) \in \left (\calT_{X'/S}\right )^d$ be the dual basis of $(\op{dlog}m_1',\hdots,\op{dlog}m_d').$ For $\beta\in  \N^d,$ set
$$
\partial'^{[\beta]}:=\prod_{i=1}^d\partial_i'^{[\beta_i]} \in \Gamma^{\bullet}\calT_{X'/S}.
$$
By \eqref{12315}, the section $\partial_I\otimes \partial'^{[\beta]}\in \calD_{X/S}^{\gamma}$ corresponds, by \eqref{eq12331}, to $(-1)^{|\beta|}I!\varphi_{I,\beta}\in \calQ_1^{\vee}.$
For $x\in \calE,$ there exist $x_{I,J}\in \calE$ such that
$$
\epsilon(1 \otimes x)=\sum_{\substack{I\in \llbracket 0,p-1 \rrbracket^d \\ J\in \N^d}} x_{I,J}\otimes \eta^I\eta_{(q)}^J
$$
and $x_{I,J}=0$ except for finitely many multi-indices $I$ and $J.$
It follows that $\varphi_{I,\beta}\cdot x=0$ for sufficiently large $|\beta|.$ The rest of the equivalence is similar to the previous case.
\end{proof}

\begin{parag}\label{paragfinalMuzan1}
We end this section by proving that stratifications with respect to the groupoid $P_{\frakX/\frakS,n}$ are equivalent to integrable quasi-nilpotent $p^n$-connections.
For integers $n\ge 0$ and $k\ge 1,$ let $\ov{\calI}_{\frakX/\frakS,n}$ be the PD-ideal of $P_{\frakX/\frakS,n}$ and $\ov{\calI}_{n,k}$ its reduction modulo $p^k.$ For all positive integers $k,$ the universal property of PD-envelopes implies the existence of a morphism
$$
\left ( P_{\frakX/\frakS,n} \right )_k \ra \left ( P_{\frakX/\frakS,0} \right )_k
$$
fitting into the commutative diagram
$$
\begin{tikzcd}
\frakX_k \ar[equal]{r} \ar{d} & \frakX_k \ar{d} \ar{dr} & \\
\left ( P_{\frakX/\frakS,n} \right )_k \ar[bend right=30]{rr} \ar{r} & \frakY_k & \left ( P_{\frakX/\frakS,0} \right )_k, \ar{l}
\end{tikzcd}
$$
where $\frakX_k \ra \frakY_k$ is the exact diagonal immersion and the other morphisms are the canonical ones.
By taking the inductive limit, these morphisms induce a morphism of formal groupoids
\begin{equation}\label{finalMuzan1}
P_{\frakX/\frakS,n} \ra P_{\frakX/\frakS,0}.
\end{equation}
Under the assumption \ref{loccoord}, the corresponding morphism of rings is
$$
\calP_{\frakX/\frakS,0} \ra \calP_{\frakX/\frakS,n},\ \eta_i \mapsto \eta_i=p^n\frac{\eta_i}{p^n}.
$$
The image of $\ov{\calI}_{\frakX / \frakS,0}$ by this morphism is $p^n\ov{\calI}_{\frakX / \frakS , n}.$ It induces, by flatness of $P_{\frakX / \frakS,n}$ over $\op{Spf}\Z_p,$ an isomorphism
$
\ov{\calI}_{\frakX / \frakS,0} \xrightarrow{\sim} \ov{\calI}_{\frakX / \frakS , n}.
$
Reducing modulo $p^k$ and using \eqref{eqtakrizIIsquare}
\begin{equation}\label{finalMuzan3}
\ov{\calI}_{0,k}^{ \{ 1 \} } := \ov{\calI}_{0,k}/\ov{\calI}_{0,k}^{[2]} \xrightarrow{\sim} \omega^1_{\frakX_k / \frakS_k },
\end{equation}
we get an isomorphism of $\Ox_{\frakX_k}$-modules
$$\beta: \omega^1_{\frakX_k/\frakS_k} \xrightarrow{\sim} \ov{\calI}_{n,k}^{\{1\}} := \ov{\calI}_{n,k}/\ov{\calI}_{n,k}^{[2]}.$$
Under the assumption \ref{loccoord}, if $\widehat{m}_i$ and $\widehat{\eta}_i$ are the reductions of $\widetilde{m}_i$ and $\frac{\widetilde{\eta}_i}{p^n}$ modulo $p^k$ respectively, then
\begin{equation}\label{KhamineiBeta}
\beta(\op{dlog}\widehat{m}_i)=\widehat{\eta}_i.
\end{equation}
Let $p_{1,n},p_{2,n}:P_{\frakX / \frakS , n} \ra \frakX$ be the canonical projections. We abusively denote by $p_{1,n}$ and $p_{2,n}$ their reductions modulo $p^k.$ For a local section $a$ of $\Ox_{\frakX_k},$ the differential $da$ corresponds, by \eqref{finalMuzan3}, to $p_{2,0}^{\#}a-p_{1,0}^{\#}a.$ Consider the derivation
$$
d':\begin{array}[t]{clc}
\Ox_{\frakX_k} & \ra & \ov{\calI}_{n,k}^{\{1\}} \\
a & \mapsto & p_{2,n}^{\#}(a)-p_{1,n}^{\#}(a).
\end{array}
$$
By construction of \eqref{finalMuzan1}, the diagram
$$
\begin{tikzcd}
P_{\frakX / \frakS , n} \ar{dr}{p_{i,n}} \ar{d} & \\
P_{\frakX / \frakS , 0} \ar[swap]{r}{p_{i,0}} & \frakX
\end{tikzcd}
$$
is commutative. It follows that, for any local section $a$ of $\Ox_{\frakX_k},$ we have
\begin{equation}\label{finalMuzan4}
p^n\beta(da) = p^n\beta \left ( p_{2,0}^{\#} a - p_{1,0}^{\#}a \right ) = p_{2,n}^{\#} a - p_{1,n}^{\#} a =d'(a).
\end{equation}
\end{parag}

\begin{proposition}\label{propfinalMuzan1}
Keep the notation of \ref{paragfinalMuzan1}. The data of a $p^n$-connection on an $\Ox_{\frakX_k}$-module $\calE$
$$
\nabla:\calE \ra \calE \otimes_{\Ox_{\frakX_k}}\omega^1_{\frakX_k/\frakS_k}
$$
is equivalent to the data of an additive morphism
$$\nabla':\calE \ra \calE \otimes_{\Ox_{\frakX_k}}\ov{\calI}_{n,k}^{\{1\}}$$
satisfying the Leibniz rule
$$\nabla'(ax)=a\nabla'(x)+x\otimes d'(a),$$
for all local sections $a$ and $x$ of $\Ox_{\frakX_k}$ and $\calE$ respectively. This equivalence is given by
$$
\nabla \mapsto \nabla'=\left ( \op{Id}_{\calE} \otimes \beta \right ) \circ \nabla.
$$
\end{proposition}

\begin{proof}
The Leibniz rule satisfied by $\nabla'$ follows from \eqref{finalMuzan4}.
\end{proof}

\begin{definition}
Let $\nabla$ be a $p^n$-connection on an $\Ox_{\frakX_k}$-module $\calE$ and $\nabla':\calE\ra \calE\otimes_{\Ox_{\frakX_k}}\ov{\calI}_{n,k}^{\{1\}}$ the corresponding additive morphism by \ref{propfinalMuzan1}. We say that $\nabla'$ is \emph{integrable} if $\nabla$ is integrable.
\end{definition}

\begin{parag}\label{compdiffopn}
Let $k\ge 1$ and $n\ge 0$ be integers and $\calE$ an $\Ox_{\frakX_k}$-module. Let $P_{n,k}$ be the logarithmic scheme obtained from $P_{\frakX/\frakS,n}$ by reduction modulo $p^k,$ $\ov{\calI}_{n,k}$ the PD-ideal of $P_{n,k},$ $\calP_{n,k}$ the structural ring of $P_{n,k}$ and $\calP_{n,k}^{\{l\}}=\calP_{n,k}/\ov{\calI}_{n,k}^{[l+1]}$ for any integer $l\ge 0.$ We recall how the composition of differential operators
$$
f:\calP_{n,k}^{\{l'\}} \otimes_{\Ox_{\frakX_k}} \calE \ra \calE,\ g:\calP_{n,k}^{\{l\}} \otimes_{\Ox_{\frakX_k}} \calE \ra \calE
$$
is defined. Let 
$$
\delta^{l,l'}:\calP^{\{l+l'\}}_{n,k} \ra \calP^{\{l\}}_{n,k} \otimes_{\Ox_{\frakX_k}} \calP^{\{l'\}}_{n,k}
$$
be the morphism induced by the comultiplication map of the Hopf algebra $\calP_{n,k}.$
The composition $g \circ f$ is defined as the composition
\begin{equation}\label{compjapon}
g\circ f:\calP^{\{l+l'\}}_{n,k}\otimes \calE \xrightarrow{\delta^{l,l'}} \calP^{\{l\}}_{n,k} \otimes \calP^{\{l'\}}_{n,k} \otimes \calE \xrightarrow{\op{Id} \otimes f} \calP^{\{l\}}_{n,k} \otimes \calE \xrightarrow{g} \calE.
\end{equation}
\end{parag}

\begin{lemma}\label{lemjapon}
Let $k\ge 1$ and $n\ge 0$ be integers and suppose \ref{loccoord} is satisfied. For any $1\le i\le d,$ let $\xi_i$ be the image of $\frac{\widetilde{\eta}_i}{p^n} \in \calP_{\frakX / \frakS , n}$ in $\calP_{n,k}:=\calP_{\frakX/\frakS,n}/(p^k)$ and set, for any $I=(I_1,\hdots,I_d)\in \N^d,$
$$\xi^{[I]}=\prod_{i=1}^d\xi_i^{[I_i]}.$$
Let $\delta:\calP_{n,k} \ra \calP_{n,k}\otimes_{\Ox_{\frakX_k}}\calP_{n,k}$ be the comultiplication map of the Hopf algebra $\calP_{n,k}.$ . Then $\left (\xi^{[I]}\right )_{I\in \N^d}$ is a basis for the $\Ox_{\frakX_k}$-module $\calP_{n,k}.$ Let $(\partial_I)_{I\in \N^d}$ be its dual basis. Then, for any $k\in \N,$ $1\le i\le d$ and $I=(I_1,\hdots,I_d)\in \N^d,$ we have
\begin{enumerate}
\item $\partial_I\circ \partial_{\epsilon_i}=\partial_{I+\epsilon_i}+p^nI_i\partial_I,$
\item $\partial_{k\epsilon_i}=\prod_{j=0}^{k-1}(\partial_{\epsilon_i}-p^nj),$
\item $\partial_I=\prod_{i=1}^d\prod_{j=0}^{I_i-1}(\partial_{\epsilon_i}-p^nj),$
\end{enumerate}
where the composition is the composition of differential operators \eqref{compjapon}.
\end{lemma}

\begin{proof}
Using \ref{lem116f}, the proof of the first formula is similar to \ref{prop17}. The other two follow from the first.
\end{proof}

\begin{parag}\label{paragJapon2}
Let $n\ge 0$ and $k\ge 1$ be integers and consider an integrable $p^n$-connection $\calE \ra \calE \otimes_{\Ox_{\frakX_k}}\omega^1_{\frakX_k/\frakS_k}$ on an $\Ox_{\frakX_k}$-module $\calE.$ Let
$$\nabla:\calE \ra \calE \otimes_{\Ox_{\frakX_k}}\ov{\calI}_{n,k}^{\{1\}}$$
be the corresponding connection given by \ref{propfinalMuzan1}. Suppose that the hypothesis \ref{loccoord} is satisfied and consider the notation introduced in \ref{lem116f}. We have a canonical isomorphism
$$\calP_{n,k}\xrightarrow{\sim} \Ox_{\frakX_k}\langle \xi_1,\hdots,\xi_d \rangle.$$
Denote by $\left ( \partial_I \right )_{I\in \N^d}$ the dual basis of $\left ( \xi^I \right )_{I\in \N^d}.$ Let $\epsilon_1$ be the $\calP_{n,k}$-linear morphism defined by
$$
\epsilon_1:\begin{array}[t]{clc}
\calP_{n,k}^{\{1\}} \otimes_{\Ox_{\frakX_k}}\calE & \ra & \calE\otimes_{\Ox_{\frakX_k}}\calP_{n,k}^{\{1\}} \\
1\otimes x & \mapsto & \nabla(x)+x\otimes 1,
\end{array}
$$
and, for $i\le 1\le d,$ the differential operator
$$
\nabla \left ( \partial_{\epsilon_i} \right ): \calP_{n,k}^{ \{ 1 \} } \otimes_{\Ox_{\frakX_k}} \calE \xrightarrow{\epsilon_1} \calE \otimes_{\Ox_{\frakX_k}} \calP_{n,k}^{ \{ 1 \} } \xrightarrow{\op{Id}_{\calE} \otimes \partial_{\epsilon_i}} \calE.
$$
The integrability of $\nabla$ implies that the differential operators $\nabla (\partial_{\epsilon_i})$ pairwise commute (here, we consider composition as differential operators \eqref{compdiffopn}) and we can define $\nabla (\partial_N),$ for any multi-index $N=(n_1,\hdots,n_d)\in\mathbb{N}^d,$ by
$$
\nabla (\partial_N)=\prod_{i=1}^d\prod_{j=0}^{n_i-1} \left ( \nabla (\partial_{\epsilon_i})-p^nj \right ),
$$
where $p^nj$ denotes $p^nj \op{Id}_{\calE},$ considered as a differential operator.
\end{parag}

\begin{definition}\label{pnqn}
Keep the notation of the proof of \ref{paragJapon2}. We say that an integrable $p^n$-connection $\nabla$ on a $\Ox_{\frakX_k}$-module $\calE$ is \emph{quasi-nilpotent} if, for every local section $x$ of $\calE,$ there exists, locally on $\frakX_k,$ an integer $m$ such that
$$
\nabla\left (\partial_I \right )(1 \otimes x)=0
$$
for all $I\in \N^d$ such that $|I|\ge m.$
\end{definition}

\begin{proposition}\label{equivstrat}
Let $k\ge 1$ and $n\ge 0$ be integers and $\calE$ an $\Ox_{\frakX_k}$-module. Let $P_{n,k}$ be the logarithmic scheme obtained from $P_{\frakX/\frakS,n}$ by reduction modulo $p^k,$ $\ov{\calI}_{n,k}$ the PD-ideal of $P_{n,k},$ $\calP_{n,k}$ the structural ring of $P_{n,k}$ and $\calP_{n,k}^{\{l\}}=\calP_{n,k}/\ov{\calI}_{n,k}^{[l+1]}$ for any integer $l\ge 0.$ The following data are equivalent:
\begin{enumerate}
\item A $\calP_{n,k}$-stratification $(\varepsilon_l)_{l\ge 0}$ on $\calE.$
\item A sequence of $\Ox_{\frakX_k}$-linear morphisms $\theta_l:\calE\ra \calE\otimes_{\Ox_{\frakX_k}}\calP_{n,k}^{\{l\}}$ satisfying the following conditions:
\begin{itemize}
\item $\theta_0=\op{Id}_{\calE}.$
\item For any integer $l\ge 0,$ the morphism $\theta_l$ is equal to the composition
$$\calE\xrightarrow{\theta_{l+1}}\calE\otimes_{\Ox_{\frakX_k}}\calP^{\{l+1\}}_{n,k}\ra \calE\otimes_{\Ox_{\frakX_k}}\calP^{\{l\}}_{n,k}$$
where the second arrow is the canonical projection.
\item For all integers $l,l'\ge 0$
\begin{equation}\label{diagstrconn1}
\begin{tikzcd}
\calE\ar{r}{\theta_{l+l'}}\ar{d}{\theta_{l'}} & \calE\otimes_{\Ox_{\frakX_k}}\calP^{\{l+l'\}}_{n,k}\ar{d}{\op{Id}\otimes \delta^{l,l'}}\\
\calE\otimes_{\Ox_{\frakX_k}}\calP^{\{l'\}}_{n,k}\ar{r}{\theta_l\otimes \op{Id}}& \calE\otimes_{\Ox_{\frakX_k}}\calP^{\{l\}}_{n,k}\otimes_{\Ox_{\frakX_k}}\calP^{\{l'\}}_{n,k},
\end{tikzcd}
\end{equation}
where $\delta^{l,l'}:\calP^{\{l+l'\}}_{n,k} \ra \calP^{\{l\}}_{n,k} \otimes_{\Ox_{\frakX_k}} \calP^{\{l'\}}_{n,k}$ is induced by the comultiplication map of the Hopf algebra $\calP_{n,k}.$
\end{itemize}
\item An integrable $p^n$-connection $\nabla:\calE\ra \calE\otimes_{\Ox_{\frakX_k}}\omega^1_{\frakX_k/\frakS_k}.$
\end{enumerate}
There exists also an equivalence between the data of a stratification
$$
\calP_{n,k} \otimes_{\Ox_{\frakX_k}} \calE \ra \calE \otimes_{\Ox_{\frakX_k}} \calP_{n,k}
$$
and an integrable quasi-nilpotent connection
$$
\nabla:\calE \ra \calE\otimes_{\Ox_{\frakX_k}} \omega^1_{\frakX_k/\frakS_k}.
$$
The connection is obtained from the stratification as follows:
$$
\nabla : \begin{array}[t]{clclclc}
\calE & \ra & \calE \otimes \ov{\calI}_{n,k} & \ra & \calE\otimes \ov{\calI}_{n,k}^{\{1\}} & \xrightarrow{\sim} & \omega^1_{\frakX_k/\frakS_k}\\
x & \mapsto & \epsilon(1\otimes x)-x\otimes 1 & & & &
\end{array},
$$
where the second arrow is the canonical projection.
\end{proposition}

\begin{proof}
See appendix.
\end{proof}

\newpage

\section{A stratified interpretation of the logarithmic Shiho functor}

\begin{parag}\label{parag91}

In this section, we consider a perfect field of characteristic $p>0.$ We denote by $W$ its ring of Witt vectors and we equip $\op{Spf}W$ with the trivial logarithmic structure. We consider a morphism $\theta:P\ra Q$ of fs monoids and a log smooth morphism of framed logarithmic $p$-adic formal schemes $(f,\theta):(\frakX,Q)\ra (\frakS,P),$ such that $\frakS$ is log flat and locally of finite type over $\op{Spf}W$ (\cite{Ahmed2010} 2.3.13). Note that this implies, by \ref{Wflat}, that the formal schemes $\frakX$ and $\frakS$ are flat over $\op{Spf}W$ \eqref{dxuflat}. We also consider the formal groupoids $Q_{\frakX},$ $R_{\frakX,n}$ and $P_{\frakX/\frakS,n}$ defined in \ref{parag86}.
Denote, for every integers $n\ge 0$ and $k\ge 1,$ by $\frakX_k,$ $\frakS_k$ and $P_{\frakX/\frakS,n,k}$ the logarithmic schemes obtained from $\frakX,$ $\frakS$ and $P_{\frakX/\frakS,n}$  by reduction modulo $p^k,$ as in \ref{prop69}. Set $X=\frakX_1,$ $S=\frakS_1,$ $Y=\frakY_1,$ $P_n=P_{\frakX/\frakS,n,1}$ and $P'_n=P_{\frakX'/\frakS,n,1}.$ Denote by $\calP_n$ and $\calP'_n$ the structural rings of $P_n$ and $P'_n$ respectively.

Denote by $F_1:X\ra X'$ the exact relative Frobenius of $X$ over $S.$ Let $Q',$ $F_P:P\ra P$ and $F_{Q/P}:Q' \ra Q$ be as defined in \ref{cor411}.
By \ref{thm410}, the logarithmic scheme $X'$ is naturally equipped with a frame $X'\ra [Q']$ such that $(F_1,F_{Q/P})$ is a morphism of framed logarithmic schemes. We suppose that $(X',Q')$ (resp. $F_1:X\ra X'$) lifts to a $p$-adic framed logarithmic formal scheme (resp. to an $\frakS$-morphism $F:(\frakX,Q)\ra (\frakX',Q')$) such that $\frakX'$ is log smooth over $\frakS.$
Set (see \ref{prop74})
$$\frakY=\frakX\times_{\frakS,[Q]}^{\op{log}}\frakX,\ \frakY'=\frakX'\times_{\frakS,[Q']}^{\op{log}}\frakX'.$$ 
\end{parag}

\begin{lemma}\label{lem92}
Denote by $p_1,p_2:\frakY\ra \frakX,$ $p_1',p_2':\frakY'\ra \frakX'$ and $q_1,q_2:Q_{\frakX} \ra \frakX$ the canonical projections and by $\Delta:\frakX \ra \frakY$ and $\Delta':\frakX' \ra \frakY'$ the strict diagonal immersions. Let $G:\frakY\ra \frakY'$ be the morphism induced by $F:(\frakX,Q)\ra (\frakX',Q').$ Then $G$ induces morphisms of formal groupoids:
\begin{equation}\label{phi}
\varphi:P_{\frakX/\frakS,0}\ra P_{\frakX'/\frakS,1},
\end{equation}
\begin{equation}\label{nu}
\nu:Q_{\frakX}\ra R_{\frakX',1}.
\end{equation}
fitting into the commutative diagram
$$
\begin{tikzcd}
P_{\frakX/\frakS,0} \ar{r}{\varphi} \ar{d} & P_{\frakX'/\frakS,1} \ar{d} \\
Q_{\frakX} \ar{r}{\nu} & R_{\frakX',1}
\end{tikzcd}
$$
where the left vertical arrow is the morphism defined in \ref{wiw1}.
In addition, if the hypothesis \ref{loccoord} is satisfied and $\nu_1:\calR_1' \ra \calQ_1$ is the reduction modulo $p$ of the morphism of Hopf algebras induced by $\nu,$ then, for $1\le i\le d,$ there exist local sections $c_1,\hdots,c_d$ of $\Ox_X$ such that
\begin{equation}\label{eqKoko1323}
\nu_1 \left (\eta_{i(r)} \right ) = \eta_{i(q)} + \sum_{k=1}^{p-1}\frac{(-1)^{k+1}}{k}\eta_i^k + q_2^{\#}c_i-q_1^{\#}c_i.
\end{equation}
\end{lemma}

\begin{proof}
We have the following commutative diagram
$$\begin{tikzcd}
\frakY\ar{r}{G}\ar{d}\ar[swap, bend right=30]{dd}{p_i} & \frakY'\ar{d}\ar[bend right=-30]{dd}{p'_i} \\
\frakX^2_{\frakS}\ar{r}{F^2}\ar{d} & \frakX'^2_{\frakS}\ar{d} \\
\frakX\ar{r}{F} & \frakX'
\end{tikzcd}$$
Let $m$ be a local section of $\calM_X$ and $m'$ its image in $\calM_{X'}$ by $\pi^{\flat}:\calM_X\ra \pi_*\calM_{X'}$ (\ref{diag51}). Let $\widetilde{m}$ (resp. $\widetilde{m}'$) be a local lifting of $m$ (resp. $m'$) to $\calM_{\frakX}$ (resp. $\calM_{\frakX'}$). Let $\eta=\eta(m)\in \Delta^{-1}\calI,\ \eta'=\eta(m')\in \Delta'^{-1}\calI'$ be as defined in \ref{parag77}.
To prove that $P_{\frakX/\frakS,0}\ra \frakY\xrightarrow{G} \frakY'$ and $Q_{\frakX} \ra \frakY \xrightarrow{G} \frakY'$ factor through $R_{\frakX',1},$ and since the local sections $\eta(m')$ locally generate the ideal of $\frakX' \ra \frakY'$ (\ref{parag77}), it is sufficient to prove that the image of $\eta(m')$ in $P_{\frakX/\frakS,0}$ (resp. $Q_{\frakX}$) belongs to $p\Ox_{P_{\frakX/\frakS,0}}$ (resp. $p\Ox_{Q_{\frakX}}$) and that $P_{\frakX/\frakS,0}$ (resp. $Q_{\frakX}$) is flat over $\op{Spf}\Z_p.$ Denote by $G^{\#}$ the morphism of structure sheaves associated to $G.$ By \ref{lem12}, there exists a local invertible section $u$ of $\calM_{\frakX}$ and a local section $b$ of $\Ox_{\frakX}$ such that $F^{\flat}(m')=pm+u$ and $\alpha_{\frakX}(u)=1+pb.$ By \ref{lemcalc}, we have
\begin{alignat}{2}
 & G^{\#}\left (\eta(m')\right ) \nonumber\\
=& \alpha_{\frakY}(p_2^{\flat}F^{\flat}m'-p_1^{\flat}F^{\flat}m')-1 \label{eq921}\\
=& \left ((\eta(m))^p+\sum_{k=1}^{p-1}\begin{pmatrix}p\\k \end{pmatrix}(\eta(m))^k+1 \right )\alpha_{\frakY}(p_2^{\flat}u-p_1^{\flat}u)-1. \nonumber
\end{alignat}
In $P_{\frakX/\frakS,0},$ this is equal to
$$\left (p!(\eta(m))^{[p]}+\sum_{k=1}^{p-1}\begin{pmatrix}p\\k \end{pmatrix}(\eta(m))^k\right )\alpha_{\frakY}(p_2^{\flat}u-p_1^{\flat}u)+\alpha_{\frakY}(p_2^{\flat}u-p_1^{\flat}u)-1.$$
In $Q_{\frakX},$ it is equal to
$$\left (p\frac{(\eta(m))^p}{p}+\sum_{k=1}^{p-1}\begin{pmatrix}p\\k \end{pmatrix}(\eta(m))^k\right )\alpha_{\frakY}(p_2^{\flat}u-p_1^{\flat}u)+\alpha_{\frakY}(p_2^{\flat}u-p_1^{\flat}u)-1.$$
All we need is to prove that
$$\alpha_{\frakY}(p_2^{\flat}u-p_1^{\flat}u)-1$$
is a section of $p\Ox_{\frakY},$ which is indeed the case:
\begin{alignat*}{2}
\alpha_{\frakY}(p_2^{\flat}u-p_1^{\flat}u)-1 &= \frac{\alpha_{\frakY}(p_2^{\flat}u)-\alpha_{\frakY}(p_1^{\flat}u)}{\alpha_{\frakY}(p_1^{\flat}u)} \\
&= \frac{p_2^{\#}(\alpha_{\frakX}(u))-p_1^{\#}(\alpha_{\frakX}(u))}{\alpha_{\frakY}(p_1^{\flat}u)} \\
&= p \frac{p_2^{\#}b-p_1^{\#}b}{\alpha_{\frakY}(p_1^{\flat}u)}.
\end{alignat*}
The flatness of $P_{\frakX/\frakS,0}$ and $Q_{\frakX}$ over $\op{Spf}\Z_p$ follows from \ref{propflat} and \ref{Qlogflat}, hence the formula \eqref{eqKoko1323}.
We also deduce that there exist unique morphisms $P_{\frakX/\frakS,0}\ra R_{\frakX',1}$ and $Q_{\frakX} \ra R_{\frakX',1}$ making the following diagrams commutative
$$\begin{tikzcd}
P_{\frakX/\frakS,0}\ar{r} \ar{d} & R_{\frakX',1}\ar{d} & & Q_{\frakX} \ar{r} \ar{d} & R_{\frakX',1} \ar{d} \\
\frakY\ar{r} & \frakY' & & \frakY\ar{r} & \frakY'
\end{tikzcd}$$
Then by the universal property of the PD-envelope, there exists a unique morphism $\varphi:P_{\frakX / \frakS,0}\ra P_{\frakX' / \frakS,1}$ making the following diagram commutative
$$\begin{tikzcd}
\frakX\ar{r}\ar{d} & \frakX'\ar{d}\ar{dr} & \\
P_{\frakX / \frakS,0}\ar{r}\ar[bend right=30,swap]{rr}{\varphi} & R_{\frakX',1}& P_{\frakX'/\frakS,1}\ar{l}.
\end{tikzcd}$$
\end{proof}

\begin{parag}\label{remarktakriz}
Denote by $p_1,p_2:\frakY\ra \frakX$ and $p_1',p_2':\frakY'\ra \frakX'$ the canonical projections. Let $G:\frakY\ra \frakY'$ be the morphism induced by $F:(\frakX,Q)\ra (\frakX',Q'),$ where $\frakX^2_{\frakS}=\frakX\times_{\frakS}^{\op{log}}\frakX$ and $\frakX'^2_{\frakS}=\frakX'\times_{\frakS}^{\op{log}}\frakX'.$
By a similar proof to that of \ref{lem92}, there exists a unique morphism of PD-logarithmic formal schemes $\varphi':P_{\frakX/\frakS,0} \ra P_{\frakX'/\frakS,0}$ making the following diagram commutative
$$
\begin{tikzcd}
\frakX \ar{d}{\iota} \ar{r}{F} & \frakX' \ar{d}{\iota'}\\
P_{\frakX/\frakS,0} \ar{r}{\varphi'} \ar{d} & P_{\frakX'/\frakS,0}\ar{d} \\
\frakY \ar{d} \ar{r}{G} & \frakY' \ar{d} \\
\frakX\times^{\op{log}}_{\frakS }\frakX \ar{r}{F^2} & \frakX'\times_{\frakS}^{\op{log}}\frakX'
\end{tikzcd}
$$
In addition, $\varphi'^{\#}(\ov{\calI}_{\frakX'/\frakS,0})\subset p\calP_{\frakX/\frakS,0},$ where we denote by $\varphi'^{\#}:\calP_{\frakX'/\frakS,0} \ra \calP_{\frakX/\frakS,0}$ the morphism of structural rings associated with $\varphi'$ and by $\ov{\calI}_{\frakX'/\frakS,0}$ the PD ideal of $\calP_{\frakX'/\frakS,0}.$

Moreover, suppose that we have local coordinates $m_{1,1},\hdots, m_{1,d}\in \Gamma(X,\calM_X)$ of $X$ with respect to $S$ and let $m_{1,1}',\hdots,m_{1,d}'$ be their images in $\calM_{X'}$ by $\pi^{\flat}:\calM_X\ra \pi_*\calM_{X'}$ \eqref{diag51}. Let $m_1,\hdots,m_d$ (resp. $m_1',\hdots,m_d'$) be local liftings of $m_{1,1},\hdots,m_{1,d}$ (resp. $m'_{1,1},\hdots,m'_{1,d}$) to $\calM_{\frakX}$ (resp. $\calM_{\frakX'}$). By \ref{lem12}, for every $1\le i\le d,$ there exists a local invertible section $u_i$ of $\calM_{\frakX}$ and a local section $b_i$ of $\Ox_{\frakX}$ such that $F^{\flat}(m_i')=pm_i+u_i$ and $\alpha_{\frakX}(u_i)=1+pb_i.$ Let $\widetilde{\eta}'_i$ and $\widetilde{\eta}_i$ be as defined in \ref{parag77} and $\eta'_i$ and $\eta_i$ their reduction modulo $p$ respectively. For any $1\le i\le d,$ a proof similar to that of \ref{lem92} shows that
\begin{alignat*}{1}
\varphi'^{\#}\left ( \widetilde{\eta}'_i \right ) = \left (p!\widetilde{\eta}_i^{[p]}+\sum_{k=1}^{p-1}\begin{pmatrix}p\\k \end{pmatrix}\widetilde{\eta}_i^k\right )\alpha_{\frakY}(p_2^{\flat}u_i-p_1^{\flat}u_i) +p \frac{p_2^{\#}b_i-p_1^{\#}b_i}{\alpha_{\frakY}(p_1^{\flat}u_i)}.
\end{alignat*}
For any $1\le i\le d$ and $1\le j\le 2,$
$$\alpha_{\frakY}(p_j^{\flat}u_i)=p_j^{\#}(\alpha_{\frakX}(u_i))=1+p p_j^{\#}(b_i).$$
Denote by $d:p\calP_{\frakX / \frakS,0} \ra \calP_{\frakX / \frakS,0}$ the inverse map of the isomorphism
$\calP_{\frakX / \frakS,0} \ra p \calP_{\frakX / \frakS,0},\ x\mapsto px$ and by $\pi:\calP_{\frakX / \frakS,0} \ra \calP_{\frakX / \frakS,0}/p\calP_{\frakX / \frakS,0}$ the canonical projection.
Modulo $p,$ we have, for any $1\le k\le p-1,$
$$d\left (\begin{pmatrix} p \\ k \end{pmatrix} \right )=\frac{(-1)^{k-1}}{k}.$$
We get
\begin{equation}
\pi \circ d\circ \varphi'^{\#}(\widetilde{\eta}'_i)=-\eta_i^{[p]}+\sum_{k=1}^{p-1}\frac{(-1)^{k-1}}{k}\eta_i^k+\pi(p_2^{\#}b_i-p_1^{\#}b_i).
\end{equation}
Note that, for every $1\le i\le d,$ the image of $b_i$ in $\Ox_X$ is equal to the local section $c_i$ given in \eqref{eqKoko1323}. If we abusively denote by $p_1,p_2:Y \ra X$ the canonical projections, then
\begin{equation}\label{era2eq10111}
\pi \circ d\circ \varphi'^{\#}(\widetilde{\eta}'_i)=-\eta_i^{[p]}+\sum_{k=1}^{p-1}\frac{(-1)^{k-1}}{k}\eta_i^k+p_2^{\#}c_i-p_1^{\#}c_i,
\end{equation}
where the sections $c_i$ are given in \eqref{eqKoko1323}.
Also, we have a relation between $\varphi:P_{\frakX/\frakS,0} \ra P_{\frakX'/\frakS,1}$ and $\varphi':P_{\frakX/\frakS,0} \ra P_{\frakX'/\frakS,0}$ given as follows: the universal property of PD envelopes yields a unique morphism of PD logarithmic formal schemes
$$P_{\frakX'/\frakS,1} \ra P_{\frakX'/\frakS,0}$$
that fits into a commutative diagram
$$
\begin{tikzcd}
\frakX' \ar[equal]{r} \ar{d} & \frakX' \ar{d} \ar{dr} & \\
P_{\frakX'/\frakS,1} \ar[bend right=30]{rr} \ar{r} & \frakY' & P_{\frakX'/\frakS,0} \ar{l}
\end{tikzcd}
$$
The uniqueness of $\varphi'$ proves that the following diagram is commutative
$$
\begin{tikzcd}
P_{\frakX/\frakS,0} \ar{r}{\varphi} \ar{dr}{\varphi'} & P_{\frakX'/\frakS,1} \ar{d} \\
& P_{\frakX'/\frakS,0}
\end{tikzcd}
$$ 
\end{parag}

\begin{remark}\label{rem93}
Using the canonical isomorphisms (\ref{rem87})
$$P_{\frakX / \frakS,0}\times_{\frakX}P_{\frakX / \frakS,0}\xrightarrow{\sim}P_{\frakX / \frakS,0}(2),\ P_{\frakX' / \frakS,1}\times_{\frakX'}P_{\frakX' / \frakS,1}\xrightarrow{\sim}P_{\frakX' / \frakS,1}(2)$$
and the morphism $\varphi:P_{\frakX / \frakS,0} \ra P_{\frakX' / \frakS,1}$ defined in \ref{lem92}, we deduce the existence of a morphism
\begin{equation}\label{eq931}
\varphi(2):P_{\frakX/\frakS,0}(2)\ra P_{\frakX'/\frakS,1}(2).
\end{equation}
\end{remark}

\begin{proposition}\label{Khaminei145}
The morphism $\varphi : P_{\frakX/\frakS,0} \ra P_{\frakX'/\frakS,1}$ \eqref{phi} is log flat.
\end{proposition}

\begin{proof}
Consider the canonical projections $p_1:P_{\frakX/\frakS,0} \ra \frakX$ and $p_1':P_{\frakX'/\frakX,1} \ra \frakX'.$ The morphism $\varphi$ factors as
$$
\varphi:P_{\frakX/\frakS,0} \xrightarrow{\widetilde{\varphi}} P_{\frakX'/\frakS,1} \times_{\frakX'}\frakX \ra P_{\frakX'/\frakS,1}.
$$
The second morphism is log flat since so is $F:\frakX\ra \frakX',$ it is hence sufficient to prove that the first morphism is log flat. Since it is strict, it is sufficient to prove that it is flat. We may suppose that the hypothesis \ref{loccoord} is satisfied. The morphism of structural rings of $\widetilde{\varphi}$ identifies then with a morphism
$$
\Ox_{\frakX} \left \langle \left \langle \frac{\widetilde{\eta}_1'}{p},\hdots ,\frac{\widetilde{\eta}_d'}{p} \right \rangle \right \rangle \ra \Ox_{\frakX} \left \langle \left \langle \widetilde{\eta}_1,\hdots,\widetilde{\eta}_d\right \rangle \right \rangle.
$$
By \ref{lemcalc}, keeping the same notation of \ref{lemcalc}, this morphism sends $\frac{\widetilde{\eta}_i'}{p}$ to a section of the form
\begin{equation}\label{Khaminei1451}
\left ( (p-1)!\widetilde{\eta}_i^{[p]}+\sum_{k=1}^{p-1} \frac{(p-1)!}{k!(p-k)!}\widetilde{\eta}_i^k \right ) a+xy,
\end{equation}
where $y=p_2^{\#}b-p_1^{\#}b$ is in the ideal $I$ generated by $\widetilde{\eta}_1,\hdots,\widetilde{\eta}_d,$ $x=(1+pp_1^{\#}b)^{-1}$ and $a-1=pxy.$
Since $(\widetilde{\eta}_1,\hdots,\widetilde{\eta}_d)$ is a basis of the $\Ox_{\frakX}$-module $I/I^2$ and by \eqref{grsym}, the section \eqref{Khaminei1451} is equal to
$$
-\widetilde{\eta}_i^{[p]} + \sum_{\substack{k\in \llbracket 0,p-1\rrbracket ^d\\ k\neq 0}} u_{i,k}\widetilde{\eta}^k + pv,
$$
for local sections $u_{i,k}$ of $\Ox_{\frakX}$ and $v$ of $\Ox_{\frakX} \left \langle \left \langle \widetilde{\eta}_1,\hdots,\widetilde{\eta}_d\right \rangle \right \rangle.$
The rest of the proof is the same as in (\cite{Shiho} 3.3).
\end{proof}

\begin{lemma}\label{lem96}
Let $\frakX\ra \frakU \ra \frakX\times^{\op{log}}_{\frakX',[Q]}\frakX$ a factorization of the diagonal immersion $\Delta_{\frakX/\frakX'}:\frakX\ra \frakX\times^{\op{log}}_{\frakX',[Q]}\frakX$ into a closed immersion and an open one. The ideal $\calI_{\frakX/\frakX'}$ of $\frakX \ra \frakU$ has a unique PD structure and the canonical morphism $\frakX\times_{\frakX',[Q]}^{\op{log}}\frakX \ra \frakX\times_{\frakS,[Q]}^{\op{log}}\frakX=\frakY$ induces a morphism $\psi:\frakX\times_{\frakX',[Q]}^{\op{log}}\frakX\ra P_{\frakX / \frakS,0}.$
\end{lemma}

\begin{proof}
Denote by $q_1,q_2:\frakX\times_{\frakX',[Q]}^{\op{log}}\frakX \ra \frakX$ the canonical projections.
For any local section $m$ of $\calM_{\frakX},$ consider the local section $\mu(m)$ of $\Delta_{\frakX/\frakX'}^{-1} \left (1+\calI_{\frakX/\frakX'} \right )$ defined as in \ref{parag77} and $\eta(m)=\mu(m)-1.$ Recall that we have an exact sequence \eqref{era2exactseq}
$$
0 \ra \Delta_{\frakX/\frakX'}^{-1} \left (1+\calI_{\frakX/\frakX'} \right ) \xrightarrow{\lambda}  \Delta_{\frakX/\frakX'}^{-1}\calM_{\frakX \times_{\frakX'}^{\op{log}}\frakX} \ra \calM_{\frakX}  \ra 0,
$$
and that
$
\Delta_{\frakX/\frakX'}^{-1}q_1^{\flat}m +\lambda \left ( \mu(m) \right ) = \Delta_{\frakX/\frakX'}^{-1}q_2^{\flat}m.
$
It follows that
$$
\left (\Delta_{\frakX/\frakX'}^{-1}\alpha_{\frakX\times_{\frakX',[Q]}^{\op{log}}\frakX} \right ) \left (\Delta_{\frakX/\frakX'}^{-1}q_1^{\flat}m \right )+ \mu(m) = \left (\Delta_{\frakX/\frakX'}^{-1}\alpha_{\frakX\times_{\frakX',[Q]}^{\op{log}}\frakX} \right ) \left ( \Delta_{\frakX/\frakX'}^{-1}q_2^{\flat}m \right ).
$$
From here on out, we drop $\Delta_{\frakX/\frakX'}^{-1}$ to lighten the notation. The ideal $\calI_{\frakX/\frakX'}$ is locally generated by sections of the form $\eta(m).$ 
By \ref{thmlogflat} and \ref{logflatfiber}, the lifting $F:\frakX\ra \frakX'$ is log flat. Then the base change $\frakX \times_{\frakX'}^{\op{log}}\frakX \ra \frakX$ is also log flat. Composing with the log étale morphism $\frakX\times_{\frakX',[Q]}^{\op{log}}\frakX \ra \frakX\times_{\frakX'}^{\op{log}}\frakX,$ we obtain the log flat projection $\frakX\times_{\frakX',[Q]}^{\op{log}}\frakX \ra \frakX.$ This projection is also strict and thus flat on the underlying formal schemes. Since $\frakX$ is flat over $\op{Spf}\Z_p,$ the formal scheme $\frakX \times_{\frakX',[Q]}^{\op{log}}\frakX$ is also flat over $\op{Spf}\Z_p.$ It is thus sufficient to prove that for any
local section $m$ of $\calM_{\frakX}$ and any positive integer $k,$
$\eta(m)^k\in k!\Ox_{\frakX\times_{\frakX',[Q]}^{\op{log}}\frakX}.$
We start by proving that for any local section $x$ of $\Ox_{\frakX}$ and any positive integer $k,$
$$(q_2^{\#}x-q_1^{\#}x)^k \in k!\Ox_{\frakX\times_{\frakX',[Q]}^{\op{log}}\frakX}.$$
We proceed by induction on $k:$
for $1\le k\le p-1,$ the result is obvious since $k$ is invertible in $\Ox_{\frakX\times_{\frakX'}\frakX}.$ We now check the case $k=p:$
let $x_1$ be the image of $x$ in $\Ox_X.$ Consider the morphism
$$\pi:X'\ra X$$
defined in \eqref{diag51}. Locally, $\pi^{\#}x_1$ admits a lifting $x'$ to $\Ox_{\frakX'}.$ There exists a local section $a$ of $\Ox_{\frakX}$ such that $F^{\#}(x')=x^p-pa.$
\begin{alignat*}{2}
0 &= q_2^{\#} (F^{\#}(x'))-q_1^{\#}(F^{\#}(x')) \\
&= q_2^{\#}(x)^p-q_1^{\#}(x)^p-p(q_2^{\#}(a)-q_1^{\#}(a)) \\
&= (q_2^{\#}(x)-q_1^{\#}(x))^p+\sum_{k=1}^{p-1}\begin{pmatrix}p\\k \end{pmatrix}q_2^{\#}(x)^k(q_2^{\#}(x)-q_1^{\#}(x))^k-p(q_2^{\#}(a)-q_1^{\#}(a)) \\
&= (q_2^{\#}(x)-q_1^{\#}(x))^p-pb(q_2^{\#}(x)-q_1^{\#}(x))-p(q_2^{\#}(a)-q_1^{\#}(a)).
\end{alignat*} 
So
$$(q_2^{\#}(x)-q_1^{\#}(x))^p=pb(q_2^{\#}(x)-q_1^{\#}(x))+p(q_2^{\#}(a)-q_1^{\#}(a)).$$
Now let $k>p$ and suppose the result is true for all integers $i\le k.$ Consider the largest integer $r$ such that $p^r \le k.$
\begin{alignat*}{2}
(q_2^{\#}(x)-q_1^{\#}(x))^{p^r} &= (pb(q_2^{\#}(x)-q_1^{\#}(x))+p(q_2^{\#}(a)-q_1^{\#}(a)))^{p^{r-1}} \\
&= \sum_{l=0}^{p^{r-1}}\begin{pmatrix}p^{r-1}\\ l \end{pmatrix}p^{p^{r-1}}b^l(q_2^{\#}(x)-q_1^{\#}(x))^l(q_2^{\#}(a)-q_1^{\#}(a))^{p^{r-1}-l}.
\end{alignat*}
By the induction hypothesis, the $l^{\text{th}}$ term of this sum is a local section of
$$
p^{p^{r-1}}\begin{pmatrix}p^{r-1}\\ l \end{pmatrix}l!(p^{r-1}-l)!\Ox_{\frakX\times_{\frakX'}^{\op{log}}\frakX}=p^{p^{r-1}}(p^{r-1})!\Ox_{\frakX\times_{\frakX'}^{\op{log}}\frakX}=\left (p^r\right )!\Ox_{\frakX\times_{\frakX'}^{\op{log}}\frakX}.
$$
If $k=p^r,$ we get that $\left (q_2^{\#}(x)-q_1^{\#}(x)\right )^k$ is a local section of $k!\Ox_{\frakX\times_{\frakX'}^{\op{log}}\frakX}.$ If $k>p^r,$ then
$$
\left (q_2^{\#}(x)-q_1^{\#}(x)\right )^k=\left (q_2^{\#}(x)-q_1^{\#}(x)\right )^{p^r}\left (q_2^{\#}(x)-q_1^{\#}(x)\right )^{k-p^r}
$$
is a local section of $\left (p^r\right )! \left (k-p^r \right )!\Ox_{\frakX\times_{\frakX'}^{\op{log}}\frakX}=k!\Ox_{\frakX\times_{\frakX'}^{\op{log}}\frakX}.$

Now let $m,\ m',\ u$ and $b$ as in \ref{lem12} and $\mu(F^{\flat}m'),$ $\mu(m)$ and $\eta(m)$ be as in \ref{parag77}. Since $F \circ q_1 =F \circ q_2$, we have $\mu(F^{\flat}m')=1.$
Then
$$
1=\mu(F^{\flat}m')=\mu(pm+u)=\mu(m)^p\mu(u).
$$
Since $u$ is invertible in $\calM_{\frakX},$ we also have
\begin{alignat*}{2}
\mu(u) &= \alpha_{\frakX \times_{\frakX'}^{\op{log}} \frakX} \left (q_2^{\flat}u-q_1^{\flat}u \right ) \\
&= \frac{\alpha_{\frakX \times_{\frakX'}^{\op{log}} \frakX} \left (q_2^{\flat}u \right )}{\alpha_{\frakX \times_{\frakX'}^{\op{log}} \frakX} \left (q_1^{\flat}u \right )} \\
&= \frac{q_2^{\#} \alpha_{\frakX} \left (u \right )}{q_1^{\#}\alpha_{\frakX} \left (u \right )} \\
&= \frac{1+pq_2^{\#}b}{1+pq_1^{\#}b}.
\end{alignat*}
We deduce that
\begin{alignat*}{2}
1 &= (\eta(m)+1)^p \frac{1+pq_2^{\#}b}{1+pq_1^{\#}b}.
\end{alignat*}
Then
\begin{alignat*}{2}
(\eta(m)+1)^p&=(1+pq_1^{\#}b)(1-pq_2^{\#}b+p^2(q_2^{\#}b)^2-\hdots )\\
&=1-p(q_2^{\#}b-q_1^{\#}b)\sum_{n=1}^{\infty}(-pq_2^{\#}b)^{n-1}\\
&=1+pc(q_2^{\#}b-q_1^{\#}b).
\end{alignat*}
It follows that
\begin{alignat*}{2}
\eta(m)^p &= -\sum_{k=1}^{p-1}\begin{pmatrix}p\\k \end{pmatrix}\eta(m)^k+pc(q_2^{\#}b-q_1^{\#}b) \\
&= pd\eta(m)+pc(q_2^{\#}b-q_1^{\#}b).
\end{alignat*}
We then prove by induction on $k,$ in the same way we did above, that
$$\eta(m)^k \in k!\Ox_{\frakX\times_{\frakX',[Q]}^{\op{log}}\frakX}.$$
The canonical morphism $$\frakX\times_{\frakX',[Q]}^{\op{log}}\frakX\ra \frakY=\frakX\times_{\frakS,[Q]}^{\op{log}}\frakX$$
and the fact that the ideal of the diagonal immersion $\frakX\ra \frakX\times_{\frakX'}^{\op{log}}\frakX$ has a unique PD-structure, imply the existence of a PD-morphism $\psi:\frakX\times_{\frakX',[Q]}^{\op{log}}\frakX\ra P_{\frakX / \frakS,0}$ fitting into the following commutative diagram
$$
\begin{tikzcd}
 & & P_{\frakX / \frakS,0} \ar{d} \\
 & \frakX\times^{\op{log}}_{\frakX',[Q]}\frakX \ar{ur}{\psi} \ar{r} \ar{d} & \frakY\ar{d} \\
\frakX \ar{ur} \ar{r} & \frakX\times_{\frakX'}^{\op{log}}\frakX \ar{r} & \frakX \times_{\frakS}^{\op{log}}\frakX
\end{tikzcd}
$$
\end{proof}

\begin{lemma}\label{lem95}
The diagonal immersion $X\ra X \times_{X'}^{\op{log}} X$ is exact and so $X\times_{X',[Q]}^{\op{log}}X=X\times^{\op{log}}_{X'}X.$
\end{lemma}

\begin{proof}
let $Z$ be the fiber product of $F_1:X \ra X'$ by itself in the category of fine logarithmic schemes.
Since the class of exact morphisms is stable by base change in the category of fine logarithmic schemes (\cite{Ogus2018} III 2.2.1.3), the second projection $p:Z\ra X$ is exact. Etale locally on $X$ and $X',$ the exact relative Frobenius $F_1:X\ra X'$ admits a chart
$$\varphi:Q'\ra Q,\ (x,y)\mapsto px+\theta^{gp} (y).$$
So, étale locally, the second projection $p_2:Z \ra X$ (resp. the canonical morphism $X\times_{X'}^{\op{log}}X\ra Z$) has a chart given by
\begin{alignat*}{2}
q_1 &: Q\ra (Q\oplus_{Q'}Q)^{int},\ x\mapsto (x,0) \\
( \text{resp.}\ i &: (Q\oplus_{Q'}Q)^{int}\ra (Q\oplus_{Q'}Q)^{sat} ).
\end{alignat*}
Let
$$\Delta:(Q\oplus_{Q'}Q)^{sat}\ra Q,\ (x,y)\mapsto x+y$$
be a local chart of the diagonal immersion $X\ra X\times_{X'}^{\op{log}}X.$
We thus have the following commutative diagram of monoids:
$$\begin{tikzcd}
(Q\oplus_{Q'}Q)^{sat}\ar{r}{\Delta} & Q\\
(Q\oplus_{Q'}Q)^{int}\ar{u}{i} & \\
Q\ar{u}{p_2}\ar[swap]{uur}{\op{Id}_Q} &
\end{tikzcd}$$
corresponding to
$$
\begin{tikzcd}
X\ar{r}\ar{ddr} & X\times_{X'}^{\op{log}}X\ar{d} & \\
 & Z\ar{d}{p_2} \\
 & X
\end{tikzcd}
$$
The composition $\Delta\circ i\circ p_2$ is exact and $(Q\oplus_{Q'}Q)^{sat}$ is saturated, so, by (\cite{Ogus2018} I 4.2.1.2), it is sufficient to prove that $\op{coker}(i\circ p_2)^{gp}$ is a torsion group. The morphism $(i\circ p_2)^{gp}$ is
$$Q^{gp}\ra Q^{gp}\oplus_{Q'^{gp}}Q^{gp},\ x\mapsto (0,x).$$ 
Let $(x,y)\in Q^{gp}\oplus_{Q'^{gp}}Q^{gp}.$ Then
$$p(x,y)=(\varphi(x,0),py)=(0,py-\varphi(x,0))$$
and that concludes the proof.
\end{proof}

\begin{corollaire}\label{cor1014}
The diagonal immersion $\frakX \ra \frakX \times_{\frakX'}^{\op{log}} \frakX$ is exact and so $\frakX \times_{\frakX',[Q]}^{\op{log}}\frakX = \frakX \times_{\frakX'}^{\op{log}} \frakX.$
\end{corollaire}

\begin{proof}
The exactness of the diagonal immersion $\frakX \ra \frakX \times_{\frakX'}^{\op{log}} \frakX$ follows from \ref{lem95} and the fact that the canonical morphisms $X \ra \frakX$ and $X\times_{X'}^{\op{log}}X \ra \frakX \times_{\frakX'}^{\op{log}}\frakX$ are strict.
\end{proof}

\begin{lemma}\label{lem98}
The composition $\varphi\circ\psi:\frakX\times_{\frakX',[Q]}^{\op{log}}\frakX\ra P_{\frakX' / \frakS,1},$ of the morphisms $\varphi$ and $\psi$ defined in \ref{lem92} and \ref{lem96}, factors through the canonical immersion $\frakX'\ra P_{\frakX' / \frakS,1}.$
\end{lemma}

\begin{proof}
Let $r_1,r_2:P_{\frakX'/\frakS,1}\ra \frakX'$ and $q_1,q_2:\frakX\times_{\frakX',[Q]}^{\op{log}}\frakX\ra \frakX$ be the canonical projections and $\Delta_{\frakX/\frakX'}:\frakX \ra \frakX \times_{\frakX',[Q]}^{\op{log}} \frakX$ the exact diagonal immersion. We denote by $\calI_{\frakX/\frakX'}$ its ideal. By \ref{lem96}, $\Delta_{\frakX/\frakX'}$ is a universal homeomorphism.
By \ref{parag86}, the formal scheme $P_{\frakX'/\frakS,1}$ is constructed as the inductive limit of $(P_{\frakX'/\frakS,1,k})_{k\ge 1},$ where $P_{\frakX'/\frakS,1,k}$ is the PD-envelope of $\frakX'_k\ra R_{\frakX',1,k}.$ Since $\Ox_{P_{\frakX'/\frakS,1,k}}$ is killed by $p^k,$ the ideal $\ov{\calI}_{\frakX'/\frakS,1,k}$ of the canonical immersion $\frakX'\ra P_{\frakX'/\frakS,1,k}$ is a nilideal and so $\frakX'_k\ra P_{\frakX'/\frakS,1,k}$ is a universal homeomorphism, hence so is $\iota':\frakX'\ra P_{\frakX'/\frakS,1}.$ To prove the lemma, it is sufficient to prove that the image of the ideal $\ov{\calI}_{\frakX'/\frakS,1}$ of $\iota',$ by $(\varphi\circ \psi)^{\#},$ vanishes. 
We have a commutative diagram with exact rows \eqref{era2exactseq}
$$
\begin{tikzcd}
0 \ar{r} & F^{-1}\iota'^{-1}\left (1+\ov{\calI}_{\frakX'/\frakS,1} \right ) \ar{r}{F^{-1}\lambda'} \ar{d}{\Delta_{\frakX/\frakX'}^{-1}\left (\varphi\circ \psi \right )^{\#}} & F^{-1}\iota'^{-1} \calM_{P_{\frakX'/\frakS,1}} \ar{d}{\Delta_{\frakX/\frakX'}^{-1}\left (\varphi \circ \psi \right )^{\flat}} \ar{r} & F^{-1}\calM_{\frakX'} \ar{d}{F^{\flat}} \ar{r} & 0 \\
0 \ar{r} & \Delta_{\frakX/\frakX'}^{-1} \left (1+\calI_{\frakX/\frakX'} \right ) \ar{r}{\lambda} & \Delta_{\frakX/\frakX'}^{-1}\calM_{\frakX \times_{\frakX'}^{\op{log}}\frakX} \ar{r} & \calM_{\frakX}  \ar{r} & 0
\end{tikzcd}
$$
We identify the étale sites of $\frakX$ and $\frakX'$ via $F,$ those of $\frakX$ and $\frakX\times_{\frakX',[Q]}^{\op{log}}\frakX$ via $\Delta_{\frakX/\frakX'}$ and those of $\frakX$ and $P_{\frakX'/\frakS,1}$ via $\iota'.$ We hence drop $F^{-1},$ $\iota'^{-1}$ and $\Delta_{\frakX/\frakX'}^{-1}$ from our future notations. Let $m'$ be a local section of $\calM_{\frakX'}$ and $\mu(m')$ be as defined in \ref{parag77}. We just have to prove that
$$
\left ( \varphi \circ \psi \right )^{\#} (\mu(m')-1)=0.
$$
For that, we just have to prove that
$$
\lambda \left (\left ( \varphi \circ \psi \right )^{\#} (\mu(m')) \right )=0.
$$
We have
$$
r_1^{\flat}m'+\lambda'(\mu(m'))=r_2^{\flat}m'.
$$
Applying $\left ( \varphi \circ \psi \right )^{\flat},$ by the commutativity of
$$
\begin{tikzcd}
\frakX \times_{\frakX'}^{\op{log}}\frakX \ar{r}{\varphi \circ \psi} \ar{d}{q_i} & P_{\frakX'/\frakS,1} \ar{d}{r_i} \\
\frakX \ar{r}{F} & \frakX',
\end{tikzcd}
$$
we get
$$
q_1^{\flat}F^{\flat}m'+\left ( \varphi \circ \psi \right )^{\flat}\lambda'(\mu(m'))=q_2^{\flat}F^{\flat}m'.
$$
Since $q_1^{\flat}F^{\flat}=q_2^{\flat}F^{\flat},$ we deduce that
$$
\left ( \varphi \circ \psi \right )^{\flat}\lambda'(\mu(m'))=0.
$$
It follows that
$$
\lambda \left ( (\varphi \circ \psi )^{\#}(\mu(m')) \right )=\left ( \varphi \circ \psi \right )^{\flat}\lambda'(\mu(m'))=0.
$$
We conclude by the injectivity of $\lambda$ that
$$
\left ( \varphi \circ \psi \right )^{\#}(\mu(m'))=1.
$$
\end{proof}

\begin{parag}
Lemmas \ref{lem92}, \ref{lem96} and \ref{lem98} and corollary \ref{cor1014}  imply the existence of the following commutative diagram
\begin{equation}\label{totdiag1}
\begin{tikzcd}
 & & \frakX'\ar{d}{\iota'} \\
 & P_{\frakX/\frakS,0} \ar{r}{\varphi} \ar{d} & P_{\frakX'/\frakS,1}\ar{d} \\
\frakX\times_{\frakX'}^{\op{log}}\frakX \ar{ur}{\psi} \ar[bend right=-30]{uurr} \ar{r} & \frakX\times^{\op{log}}_{\frakS,[Q]}\frakX \ar{r}{F^2} & \frakX'\times_{\frakS,[Q']}^{\op{log}}\frakX'
\end{tikzcd}
\end{equation}
The reduction of the previous diagram modulo $p$ yields the following commutative diagram
\begin{equation}\label{totdiag2}
\begin{tikzcd}
 & & X'\ar{d}{\iota_1'} \\
 & P_0 \ar{r}{\varphi_1} \ar{d} & P'_1\ar{d} \\
X\times_{X'}^{\op{log}}X \ar{ur}{\psi_1} \ar[bend right=-30]{uurr} \ar{r} & X\times^{\op{log}}_{S,[Q]}X \ar{r}{F_1^2} & X'\times_{S,[Q']}^{\op{log}}X'
\end{tikzcd}
\end{equation}

\end{parag}

\begin{parag}
Let $p_1,p_2:P_{\frakX / \frakS,0}\ra \frakX$ and $p_1',p_2':P_{\frakX'/\frakS,1}\ra \frakX'$ be the canonical projections. For $1\le i<j\le 3,$ let $p_{ij}:P_{\frakX/\frakS,0}(2)\ra P_{\frakX/\frakS,0}$ (resp. $p_{ij}':P_{\frakX'/\frakS,1}(2)\ra P_{\frakX'/\frakS,1}$) be the morphism induced by the $(i,j)$-projection $\frakY(2)\ra \frakY$ (resp. $\frakY'(2)\ra \frakY'$). Also let $\iota:\frakX\ra P_{\frakX / \frakS,0}$ and $\iota':\frakX'\ra P_{\frakX' / \frakS,1}$ be the canonical immersions.
Consider the morphism $\varphi:P_{\frakX / \frakS,0}\ra P_{\frakX' / \frakS,1}$ defined in \ref{lem92}.
Let $(\calE',\epsilon')$ be an object of $1\text{-}\op{MHS}(\frakX'/\frakS).$ The $1$-HPD-stratification
$$\epsilon':p_2'^*\calE'\ra p_1'^*\calE'$$
induces an isomorphism
\begin{equation}\label{eq981}
\varphi^*\epsilon':\varphi^*p_2'^*\calE' \xrightarrow{\sim} \varphi^*p_1'^*\calE'.
\end{equation}
By the commutative diagram
$$
\begin{tikzcd}
P_{\frakX / \frakS,0}\ar{r}{\varphi} \ar{d}{p_i} & P_{\frakX'/ \frakS,1} \ar{d}{p_i'} \\
\frakX\ar{r}{F} & \frakX'
\end{tikzcd}
$$
the isomorphism (\ref{eq981}) is equal to
\begin{equation}\label{eq982}
\varphi^*\epsilon':p_2^*F^*\calE' \xrightarrow{\sim} p_1^*F^*\calE'.
\end{equation}
By the commutative diagram
$$
\begin{tikzcd}
 & \frakX\ar{r}{F} \ar{d}{\iota} \ar[swap]{dl}{\op{Id}_{\frakX}} & \frakX'\ar{d}{\iota'}\ar{dr}{\op{Id}_{\frakX'}} & \\
\frakX & P_{\frakX / \frakS,0}\ar{l}{p_i} \ar{r}{\varphi} & P_{\frakX' / \frakS,1} \ar[swap]{r}{p_i'} & \frakX'
\end{tikzcd}
$$
and the fact that $\iota'^*\epsilon'=\op{Id}_{\calE'},$ we have
$$\iota^*\varphi^*\epsilon'=\op{Id}_{F^*\calE}.$$
Consider the morphism $\varphi(2):P_{\frakX / \frakS,0}(2) \ra P_{\frakX' / \frakS,1}(2)$ defined in \ref{rem93}.
The $1$-HPD-stratification $\epsilon'$ satisfies the cocycle condition
$$p_{13}'^*\epsilon'=p_{23}'^*\epsilon'\circ p_{12}'^*\epsilon',$$
which implies
\begin{equation}\label{eq983}
\varphi(2)^*p_{13}'^*\epsilon'=\varphi(2)^*p_{23}'^*\epsilon'\circ \varphi(2)^*p_{12}'^*\epsilon'.
\end{equation}
By the commutativity of the diagram
$$
\begin{tikzcd}
P_{\frakX / \frakS,0}(2) \ar{r}{p_{ij}} \ar{d}{\varphi(2)} & P_{\frakX / \frakS,0}\ar{d}{\varphi} \\
P_{\frakX' / \frakS,1}(2) \ar{r}{p'_{ij}} & P_{\frakX' / \frakS,1}
\end{tikzcd}
$$
the equality (\ref{eq983}) becomes
$$p_{13}^*\varphi^*\epsilon'=p_{23}^*\varphi^*\epsilon'\circ p_{12}^*\varphi^*\epsilon'.$$
We conclude that $\varphi^*\epsilon'$ is an HPD-stratification on $F^*\calE'$ and so we have a functor
\begin{equation}\label{eq994}
\Psi:\begin{array}[t]{clc}
1\text{-}\op{MHS}(\frakX'/\frakS) & \ra & \op{MHS}(\frakX/\frakS) \\
(\calE',\epsilon') & \mapsto & (F^*\calE',\varphi^*\epsilon').
\end{array}
\end{equation}
\end{parag}

\begin{parag}\label{parwiw1}
Consider the morphism
$$\varphi':P_{\frakX/\frakS,0} \ra P_{\frakX'/\frakS,0}$$
defined in \ref{remarktakriz}.
The composition
$$\ov{\calI}_{\frakX'/\frakS,0} \xrightarrow{\varphi'^{\#}} p\calP_{\frakX/\frakS,0} \xrightarrow{d} \calP_{\frakX/\frakS,0} \ra \calP_0,$$
where $d$ is the division by $p$ morphism (recall that, by \ref{propflat}, $P_{\frakX/\frakS,0}$ is flat over $\op{Spf}\Z_p)$ and the last arrow is reduction modulo $p,$ factors through $\ov{\calI}_{\frakX'/\frakS,0}/\ov{\calI}_{\frakX'/\frakS,0}^{[2]}.$ Indeed, for any local sections $x,y$ of $\ov{\calI}_{\frakX'/\frakS,0},$ there exists local sections $x',y'\in \calP_{\frakX/\frakS,0}$ such that $\varphi'^{\#}(x)=px'$ and $\varphi'^{\#}(y)=py'.$ Then
$$\varphi'^{\#}(xy)=p^2x'y'$$
$$\varphi'(x^{[2]})=(px')^{[2]}=p^2x'^{[2]}.$$
Since $\calP_0$ is of characteristic $p,$ and by \ref{eqtakrizIIsquare}, it then induces an $\Ox_{X'}$-linear morphism
$$\omega^1_{X'/S} \ra \calP_0.$$
Since $\varphi'$ is a morphism of PD formal schemes, the image of this morphism is included in the PD ideal of $\calP_0.$ It then induces a morphism of PD-$\Ox_{X'}$-algebras
$$\Gamma^{\bullet}\omega^1_{X'/S} \ra \calP_0.$$
By $\Ox_X$-linearization, we get an $\Ox_X$-linear morphism
\begin{equation}\label{uAkaza}
u:F_1^*\Gamma^{\bullet}\omega^1_{X'/S} \ra \calP_0.
\end{equation}
Suppose that we have local coordinates $m_1,\hdots,m_d \in \Gamma(X,\calM_X)$ and $m_i'=\pi^{\flat}m_i$ (where $\pi$ is defined in \ref{diag51}). Just as in \ref{remarktakriz}, let $\widetilde{m}_1,\hdots,\widetilde{m}_d$ (resp. $\widetilde{m}'_1,\hdots,\widetilde{m}'_d$) be local liftings of $m_1,\hdots,m_d$ (resp. $m'_1,\hdots,m'_d$) to $\calM_{\frakX}$ (resp. $\calM_{\frakX'}$) and set $\widetilde{\eta}'_i=\alpha_{\frakY'}(p_2'^{\flat}\widetilde{m}_i'-p_1'^{\flat}\widetilde{m}_i')-1 \in \ov{\calI}_{\frakX'/\frakS,0},$ $\widetilde{\eta}_i=\alpha_{\frakY}(p_2^{\flat}\widetilde{m}_i-p_1^{\flat}\widetilde{m}_i)-1 \in \ov{\calI}_{\frakX/\frakS,0}$  and $\eta'_i$ and $\eta_i$ respectively their reductions modulo $p.$ By the isomorphism (\ref{eqtakrizIIsquare})
$$\ov{\calI}^{\{1\}}_{\frakX'/\frakS,0}/p\ov{\calI}^{\{1\}}_{\frakX'/\frakS,0} \xrightarrow{\sim}\omega^1_{X'/S},$$
the image of $\eta_i'$ in $\ov{\calI}^{\{1\}}_{\frakX'/\frakS,0}/p\ov{\calI}^{\{1\}}_{\frakX'/\frakS,0}$ corresponds to $\op{dlog}m_i'.$
Then, by (\ref{era2eq10111}), there exists a section $c_i$ of $\Ox_X$ such that
\begin{equation}\label{era2equ}
u(F_1^*\op{dlog}m'_i)=-\eta_i^{[p]}+\sum_{k=1}^{p-1}\frac{(-1)^{k-1}}{k}\eta_i^k+p_2^{0\#}c_i-p_1^{0\#}c_i,
\end{equation}
where $p^0_1,p^0_2:P_0 \ra X$ are the canonical projections. Note that $c_i$ is the same local section $c_i$ given in \eqref{eqKoko1323}.
\end{parag}

\begin{proposition}\label{propAkaza}
Let $\theta$ be the composition
$$F_1^*\Gamma^{\bullet}\omega^1_{X'/S}\xrightarrow{\Delta} F_1^*\Gamma^{\bullet}\omega^1_{X'/S} \otimes_{\Ox_X}F_1^*\Gamma^{\bullet}\omega^1_{X'/S} \xrightarrow{\op{Id}\otimes u}F_1^*\Gamma^{\bullet}\omega^1_{X'/S} \otimes_{\Ox_X}\calP_0,$$
where $\Delta$ is the comultiplication of $F_1^*\Gamma^{\bullet}\omega^1_{X'/S}.$ Then the $\calP_0$-linearization of $\theta$ is a HPD stratification on $F_1^*\Gamma^{\bullet}\omega^1_{X'/S}.$ 
\end{proposition}

\begin{proof}
The composition of $\theta$ with $F_1^*\Gamma^{\bullet}\omega^1_{X'/S}\otimes_{\Ox_X}\calP_0 \ra F_1^*\Gamma^{\bullet}\omega^1_{X'/S}\otimes_{\Ox_X}\left (\calP_0/\ov{\calI}_0\right )=F_1^*\Gamma^{\bullet}\omega^1_{X'/S},$ where $\ov{\calI}_0$ is the PD ideal of $\calP_0$, is equal to the identity. Indeed, if $x$ is a local section of $\omega^1_{X'/S}$ then
$$(\op{Id}\otimes u)\circ\Delta(x)=1\otimes u(x)+x\otimes 1,$$
by definition of $u,$ $u(x)\in \ov{\calI}_0.$
Now we just have to check the commutativity of the diagram
\begin{equation}\label{diagAkaza2}
\begin{tikzcd}
F_1^*\Gamma^{\bullet}\omega^1_{X'/S} \ar{r}{\theta} \ar{d}{\theta} & F_1^*\Gamma^{\bullet}\omega^1_{X'/S}\otimes_{\Ox_X}\calP_0 \ar{d}{\op{Id}\otimes \delta} \\
F_1^*\Gamma^{\bullet}\omega^1_{X'/S}\otimes_{\Ox_X}\calP_0 \ar{r}{\theta\otimes \op{Id}} & F_1^*\Gamma^{\bullet}\omega^1_{X'/S}\otimes_{\Ox_X} \calP_0 \otimes_{\Ox_X}\calP_0
\end{tikzcd}
\end{equation}
The composition
$$F_1^*\Gamma^{\bullet}\omega^1_{X'/S} \xrightarrow{\theta} F_1^*\Gamma^{\bullet}\omega^1_{X'/S}\otimes_{\Ox_X}\calP_0 \xrightarrow{\theta \otimes \op{Id}} F_1^*\Gamma^{\bullet}\omega^1_{X'/S}\otimes_{\Ox_X} \calP_0 \otimes_{\Ox_X} \calP_0$$
is equal to
\begin{alignat*}{2}
F_1^*\Gamma^{\bullet}\omega^1_{X'/S} \xrightarrow{\Delta}  F_1^*\Gamma^{\bullet}\omega^1_{X'/S}\otimes_{\Ox_X} F_1^*\Gamma^{\bullet}\omega^1_{X'/S} \xrightarrow{\Delta\otimes \op{Id}}& F_1^*\Gamma^{\bullet}\omega^1_{X'/S}\otimes_{\Ox_X}F_1^*\Gamma^{\bullet}\omega^1_{X'/S}\otimes_{\Ox_X}F_1^*\Gamma^{\bullet}\omega^1_{X'/S} \\
\xrightarrow{\op{Id}\otimes u\otimes u} & F_1^*\Gamma^{\bullet}\omega^1_{X'/S}\otimes_{\Ox_X} \calP_0 \otimes_{\Ox_X}\calP_0.
\end{alignat*}
It follows that (\ref{diagAkaza2}) is equivalent to the commutativity of the diagram
$$
\begin{tikzcd}
F_1^*\Gamma^{\bullet}\omega^1_{X'/S} \ar{r}{\Delta} \ar{d}{\Delta} & F_1^*\Gamma^{\bullet}\omega^1_{X'/S}\otimes_{\Ox_X}F_1^*\Gamma^{\bullet}\omega^1_{X'/S} \ar{d}{\op{Id}\otimes u} \\
F_1^*\Gamma^{\bullet}\omega^1_{X'/S}\otimes_{\Ox_X}F_1^*\Gamma^{\bullet}\omega^1_{X'/S} \ar{d}{\Delta\otimes \op{Id}} & F_1^*\Gamma^{\bullet}\omega^1_{X'/S}\otimes_{\Ox_X}\calP_0 \ar{d}{\op{Id}\otimes \delta} \\
F_1^*\Gamma^{\bullet}\omega^1_{X'/S}\otimes_{\Ox_X}F_1^*\Gamma^{\bullet}\omega^1_{X'/S}\otimes_{\Ox_X}F_1^*\Gamma^{\bullet}\omega^1_{X'/S} \ar{r}{\op{Id}\otimes u \otimes u} & F_1^*\Gamma^{\bullet}\omega^1_{X'/S}\otimes_{\Ox_X}\calP_0\otimes_{\Ox_X}\calP_0
\end{tikzcd}
$$
Since
$$
\begin{tikzcd}
F_1^*\Gamma^{\bullet}\omega^1_{X'/S} \ar{r}{\Delta} \ar{d}{\Delta} & F_1^*\Gamma^{\bullet}\omega^1_{X'/S}\otimes_{\Ox_X}F_1^*\Gamma^{\bullet}\omega^1_{X'/S} \ar{d}{\op{Id}\otimes\Delta} \\
F_1^*\Gamma^{\bullet}\omega^1_{X'/S}\otimes_{\Ox_X}F_1^*\Gamma^{\bullet}\omega^1_{X'/S} \ar{r}{\Delta \otimes \op{Id}} & F_1^*\Gamma^{\bullet}\omega^1_{X'/S}\otimes_{\Ox_X}F_1^*\Gamma^{\bullet}\omega^1_{X'/S}\otimes_{\Ox_X}F_1^*\Gamma^{\bullet}\omega^1_{X'/S}
\end{tikzcd}
$$
is commutative, it is sufficient to prove the commutativity of
\begin{equation}\label{diagAkaza}
\begin{tikzcd}
F_1^*\Gamma^{\bullet}\omega^1_{X'/S} \ar{r}{u} \ar{d}{\Delta} & \calP_0\ar{d}{\delta} \\
F_1^*\Gamma^{\bullet}\omega^1_{X'/S}\otimes_{\Ox_X}F_1^*\Gamma^{\bullet}\omega^1_{X'/S}\ar{r}{u\otimes u} & \calP_0\otimes_{\Ox_X}\calP_0
\end{tikzcd}
\end{equation}
To prove the commutativity of (\ref{diagAkaza}), and given that
$$F_1^*\Gamma^{\bullet}\omega^1_{X'/S} \otimes_{\Ox_X} F_1^*\Gamma^{\bullet}\omega^1_{X'/S}=F_1^*\left (\Gamma^{\bullet}\omega^1_{X'/S} \otimes_{\Ox_{X'}} \Gamma^{\bullet}\omega^1_{X'/S} \right )=F_1^*\Gamma^{\bullet}(\omega^1_{X'/S}\oplus \omega^1_{X'/S}),$$
it is sufficient to prove the commutativity of the diagram
$$
\begin{tikzcd}
\omega^1_{X'/S} \ar{rrrr}{u} \ar{d}{\Delta'} & & & & \calP_0 \ar{d}{\delta} \\
\omega^1_{X'/S} \oplus \omega^1_{X'/S} \ar{rrrr}{(x,y) \mapsto 1\otimes u(y)+ u(x)\otimes 1} & & & & \calP_0\otimes_{\Ox_X}\calP_0
\end{tikzcd}
$$
where $\Delta'$ is the diagonal embedding.
We may work étale locally and so we can suppose that we have local coordinates $m_1,\hdots,m_d \in \Gamma(X,\calM_X).$ For any $1\le i\le d,$ let $m'_i,$ $\widetilde{\eta}_i$ and $\widetilde{\eta}'_i$ as defined in \ref{parwiw1}.
Consider the morphism $\widetilde{\delta}:\calP_{\frakX/\frakS,0} \ra \calP_{\frakX/\frakS,0}\otimes_{\Ox_{\frakX}}\calP_{\frakX/\frakS,0}$ (resp. $\widetilde{\delta'}:\calP_{\frakX'/\frakS,0} \ra \calP_{\frakX'/\frakS,0}\otimes_{\Ox_{\frakX'}}\calP_{\frakX'/\frakS,0}$) induced by the $(1,3)$-projection $\frakX \times_{\frakS}^{\op{log}} \frakX \times_{\frakS}^{\op{log}} \frakX \ra \frakX\times_{\frakS}^{\op{log}}\frakX$ (resp. $\frakX' \times_{\frakS}^{\op{log}} \frakX' \times_{\frakS}^{\op{log}} \frakX' \ra \frakX' \times_{\frakS}^{\op{log}}\frakX'$). To lighten notation, denote $\calP_{\frakX / \frakS,0}\otimes_{\Ox_{\frakX}}\calP_{\frakX / \frakS,0}$ by $\calP(2).$ By \ref{propflat} and \ref{rem87}, $\calP(2)$ is flat over $\Z_p.$ By construction of $\varphi',$ the following diagram is commutative
$$
\begin{tikzcd}
P_{\frakX/\frakS,0}\times_{\frakX} P_{\frakX/\frakS,0} \ar{r}{\varphi'\times \varphi'} \ar{d} & P_{\frakX'/\frakS,0}\times_{\frakX'}P_{\frakX'/\frakS,0} \ar{d} \\
P_{\frakX/\frakS,0} \ar{r}{\varphi'} & P_{\frakX'/\frakS,0}
\end{tikzcd}
$$
It follows that the left square of the following diagram is commutative:
\begin{equation}\label{era2diaghh}
\begin{tikzcd}
\ov{\calI}_{\frakX'/\frakS,0} \ar{r}{\varphi'^{\#}} \ar{d}{\widetilde{\delta'}} & p\calP_{\frakX / \frakS,0} \ar{d}{\widetilde{\delta}} \ar{r}{d} & \calP_{\frakX / \frakS,0} \ar{d}{\widetilde{\delta}} \ar{r} & \calP_0 \ar{d}{\delta} \\
\widetilde{\delta'}(\ov{\calI}_{\frakX'/\frakS,0}) \ar{r}{\varphi'^{\#2}} & p \calP(2) \ar{r}{d'} &  \calP(2) \ar{r} & \calP_0 \otimes_{\Ox_X} \calP_0
\end{tikzcd}
\end{equation}
where upper and lower right arrows are the canonical projections, $d$ and $d'$ are the division by $p$ isomorphisms and $\varphi'^{\#2}$ is induced by
$$\begin{array}[t]{clc}
\calP_{\frakX'/\frakS}\otimes_{\Ox_{\frakX'}}\calP_{\frakX'/\frakS} & \ra & p\calP_{\frakX/\frakS}\otimes_{\Ox_{\frakX}}\calP_{\frakX/\frakS}\\
x\otimes y & \mapsto & \varphi'(x)\otimes \varphi'(y).
\end{array}
$$
The morphism $u$ is induced by the upper line and the two right squares are clearly commutative. 
Since $(\varphi'^{\#}\otimes\varphi'^{\#})(\widetilde{\eta}_i'\otimes \widetilde{\eta}_i')\in p^2\calP(2)$ and $\widetilde{\delta'}(\widetilde{\eta}'_i)=1\otimes \widetilde{\eta}_i'+\widetilde{\eta}_i'\otimes 1+\widetilde{\eta}_i'\otimes \widetilde{\eta}_i'$ (\ref{era2prop611}), we get
$$\delta(u(\op{dlog}m_i'))=1\otimes u(\op{dlog}m_i')+u(\op{dlog}m'_i)\otimes 1.$$
\end{proof}

\begin{parag}
The HPD stratification on $F_1^*\Gamma^{\bullet}\omega^1_{X'/S},$ defined in \ref{propAkaza}, is equivalent to a $\calD_{X/S}$-module structure on $F_1^*\Gamma^{\bullet}\omega^1_{X'/S}.$ By duality and \ref{dualAkaza}, this yields a $\calD_{X/S}$-module structure on $F_1^*\widehat{S}^{\bullet}\calT_{X'/S},$ where $\widehat{S}^{\bullet}\calT_{X'/S}$ is the completion of $S^{\bullet}\calT_{X'/S}$ with respect to the ideal
$$S^{\ge 1}\calT_{X'/S}:=\bigoplus_{n\ge 1}S^n\calT_{X'/S}.$$
In the following proposition and corollary, we give explicit formulae for these actions.
\end{parag}

\begin{proposition}\label{propwiw1025}
Suppose that we have local coordinates $m_1,\hdots,m_d \in \Gamma(X,\calM_X)$ (\ref{P2}) and set $m_i'=\pi^{\flat}m_i,$ where $\pi$ is defined in \ref{diag51}. Let $n\ge 0,$ $u:F_1^*\Gamma^{\bullet}\omega^1_{X'/S} \ra \calP_0$ the morphism \eqref{uAkaza} and $p_1,p_2:P_0 \ra X$ and $\pi_n:\calP_0 \ra \calP_0^n:=\calP_0/\ov{\calI}_0^{[n+1]}$ the canonical projections. Consider the action of $\calD_{X/S}$ on $F_1^*\Gamma^{\bullet}\omega^1_{X'/S}$ given by the stratification defined in \ref{propAkaza}.
Then, the local sections $c_i \in \Ox_X,$ given in \eqref{eqKoko1323}, satisfy
\begin{alignat}{2} \label{era2action}
\partial_{\epsilon_j}\cdot F_1^*\op{dlog}m'_i &= \begin{cases} \partial_{\epsilon_j}(p_2^{\#}c_i-p_1^{\#}c_i)\ \op{if}\ i\neq j\\
1+\partial_{\epsilon_i}(p_2^{\#}c_i-p_1^{\#}c_i)\ \op{if}\ i=j.\end{cases}\\
\partial_{p\epsilon_i}\cdot F_1^*\op{dlog}m'_i &= -1+\partial_{p\epsilon_i}(p_2^{\#}c_i-p_1^{\#}c_i).
\end{alignat}
where $\partial_{\epsilon_j}$ and $\partial_{p\epsilon_i}$ are the differential operators defined in \ref{P2}.
\end{proposition}

\begin{proof}
By definition, the action $\partial \cdot F_1^*\op{dlog}m'_i$ of $\partial$ on $F_1^*\op{dlog}m'_i$ is the image of $F_1^*\op{dlog}m'_i$ by the composition
$$F_1^*\Gamma^{\bullet} \omega^1_{X'/S} \xrightarrow{(\op{Id}\otimes u)\circ \Delta}F_1^*\Gamma^{\bullet}\omega^1_{X'/S} \otimes_{\Ox_X}\calP_0 \xrightarrow{\op{Id}\otimes \pi_n} F_1^*\Gamma^{\bullet}\omega^1_{X'/S} \otimes_{\Ox_X}\calP_0^n \xrightarrow{\op{Id}\otimes \partial} F_1^*\Gamma^{\bullet}\omega^1_{X'/S},$$
where $\Delta:F_1^*\Gamma^{\bullet} \omega^1_{X'/S} \ra  (F_1^*\Gamma^{\bullet} \omega^1_{X'/S}) \otimes_{\Ox_X} (F_1^*\Gamma^{\bullet} \omega^1_{X'/S})$ is the comultiplication map. For a local section $x$ of $F_1^*\Gamma^{\bullet} \omega^1_{X'/S}$ of degree $1,$ we have
$$(\op{Id}\otimes u)\circ \Delta(x)=1\otimes u(x)+x\otimes 1.$$
The formulae \eqref{era2action} follow from \eqref{era2equ}.
\end{proof}

\begin{corollaire}\label{era2parag113}
Let $m_1,\hdots,m_d \in \Gamma(X,\calM_X),$ $m_1',\hdots,m_d'\in \Gamma(X',\calM_{X'}),$ $(\op{dlog}m')^{[I]}$ and $p_1,p_2:P_0 \ra X$ as in \ref{propwiw1025}. Let $(\partial'_i)_{1\le i\le d}$ be the dual basis of $(\op{dlog}m'_i)_{1\le i\le d}.$ We also consider the differential operators $\partial_{\epsilon_i}\in \calD_{X/S}$ defined in \ref{P2}. For every $1\le i\le d,$ the action of $\partial_{\epsilon_i}$ on $1\in F_1^*\w{S}^{\bullet}\calT_{X'/S}$ is given by
\begin{equation}\label{era2action3}
\partial_{\epsilon_i}\cdot 1=F_1^*\partial'_{i}+\sum_{j=1}^d\partial_{\epsilon_i}(p_2^{\#}c_j-p_1^{\#}c_j)F_1^*\partial'_j \in F_1^*\w{S}^{\bullet}\calT_{X'/S},
\end{equation} 
where the local sections $c_i$ are given in \eqref{eqKoko1323}.
\end{corollaire}

\begin{proof}
The action of $\calD_{X/S}$ on $\mathscr{Hom}_{\Ox_X}(F_1^*\Gamma^{\bullet}\omega^1_{X'/S},\Ox_X),$ defined by duality and the action of $\calD_{X/S}$ on $F_1^*\Gamma^{\bullet}\omega^1_{X'/S}$ is given as follows: if $f\in \mathscr{Hom}_{\Ox_X}(F_1^*\Gamma^{\bullet}\omega^1_{X'/S},\Ox_X)$ then the action of $\partial_{\epsilon_i}$ on $f$ is
$$F_1^*\Gamma^{\bullet}\omega^1_{X'/S} \ra \Ox_X,\ x\mapsto f(\partial_{\epsilon_i}\cdot x).$$
By the isomorphism (\ref{isodualAkaza}), the section $1\in F_1^*\w{S}^{\bullet}\calT_{X'/S}$ corresponds to the $\Ox_X$-linear morphism $\alpha:F_1^*\Gamma^{\bullet}\omega^1_{X'/S} \ra \Ox_X$ given by $\alpha(a)=a$ for any $a\in \Ox_X$ and $\alpha(a)=0$ for any $a\in F_1^*\Gamma^{\ge 1}\omega^1_{X'/S}.$
The action of $\partial_{\epsilon_i}$ on $\alpha$ is the composition
$$
\partial_{\epsilon_i}\cdot \alpha : F_1^*\Gamma^{\bullet}\omega^1_{X'/S} \xrightarrow{\Delta} F_1^*\Gamma^{\bullet}\omega^1_{X'/S}\otimes_{\Ox_X} F_1^*\Gamma^{\bullet}\omega^1_{X'/S} \xrightarrow{\op{Id} \otimes (\partial_{\epsilon_i}\circ u)} F_1^*\Gamma^{\bullet}\omega^1_{X'/S} \xrightarrow{\alpha} \Ox_X,
$$
where $u$ is given in \eqref{uAkaza} and $\Delta$ is the comultiplication map.
For simplicity, we set $x_i=F_1^*\op{dlog}m_i'$ and $x^I=F_1^* \left (\op{dlog}m'\right )^{[I]}$ for $1\le i\le d$ and $I\in \N^d.$ The image of $x^{[I]}$ by $\Delta$ is
$$
\sum_{0\le J\le I}x^{[I-J]}\otimes x^{[J]}.
$$
This is sent by $\op{Id} \otimes (\partial_{\epsilon_i}\circ u)$ to
\begin{equation}\label{eqKhaminei14151}
\sum_{0\le J\le I}\partial_{\epsilon_i}\left (u(x)^{[J]} \right )x^{[I-J]}.
\end{equation}
By \eqref{era2equ}, $\partial_{\epsilon_i}\left (u(x)^{[J]} \right )=0$ for $|J|\ge 2$ and $J=0,$ so \eqref{eqKhaminei14151} is equal to
\begin{equation}\label{eqKhaminei14152}
\sum_{j=1}^d \partial_{\epsilon_i}\left (u(x_j) \right ) x^{[I-\epsilon_j]},
\end{equation}
with the convention $x^{[I-\epsilon_j]}=0$ if $I-\epsilon_j$ has a negative coefficient. If $|I|\ge 2,$ the image of \eqref{eqKhaminei14152} by $\alpha$ vanishes. It follows that
$$
\left ( \partial_{\epsilon_i} \cdot \alpha \right )(x^{[I]})=0\quad \forall |I|\neq 1.
$$
In addition, by \eqref{era2equ}, the local sections $c_j\in \Ox_X,$ appearing in \eqref{eqKoko1323}, satisfy
$$
(\partial_{\epsilon_i}\cdot \alpha)(x_j)=\delta_{ij}+\partial_{\epsilon_i}(p_2^{\#}c_j-p_1^{\#}c_j).
$$
The family $(\op{dlog}m'^{[I]})_{I\in \N^d}$ is a basis for the $\Ox_{X'}$-module $\Gamma^{\bullet}\omega^1_{X'/S}.$ Denote by $(\alpha_I)$ the dual basis.
Then
$$\partial_{\epsilon_i}\cdot \alpha=F_1^*\alpha_{\epsilon_i}+\sum_{j=1}^d\partial_{\epsilon_i}(p_2^{\#}c_j-p_1^{\#}c_j)F_1^*\alpha_{\epsilon_j} \in F_1^*\mathscr{Hom}_{\Ox_{X'}}(\Gamma^{\bullet}\omega^1_{X'/S},\Ox_{X'}).$$
By the isomorphism (\ref{isodualAkaza}), the desired action on $1\in F_1^*\w{S}^{\bullet}\calT_{X'/S}$ is
\begin{equation*}
\partial_{\epsilon_i}\cdot 1=F_1^*\partial'_{i}+\sum_{j=1}^d\partial_{\epsilon_i}(p_2^{\#}c_j-p_1^{\#}c_j)F_1^*\partial'_j \in F_1^*\w{S}^{\bullet}\calT_{X'/S}.
\end{equation*}
\end{proof}

\begin{parag}
The action of $\calD_{X/S}$ on $F_1^*\widehat{S}^{\bullet}\calT_{X'/S}$ and the $p$-curvature map $\psi:S^{\bullet}\calT_{X'/S}\ra \calD_{X/S}$ (\ref{eqDoma1}) induce an action of $S^{\bullet}\calT_{X'/S}$ on $F_1^*\widehat{S}^{\bullet}\calT_{X'/S}.$ In the next proposition, we check that this action extends to an action of $\widehat{S}^{\bullet}\calT_{X'/S}$ on $F_1^*\widehat{S}^{\bullet}\calT_{X'/S}.$
\end{parag}

\begin{proposition}\label{era4tak1}
The action of $S^{\bullet}\calT_{X'/S}$ on $F_1^*\widehat{S}^{\bullet}\calT_{X'/S},$ induced by the action of $\calD_{X/S}$ on $F_1^*\widehat{S}^{\bullet}\calT_{X'/S}$ and the $p$-curvature map $\psi:S^{\bullet}\calT_{X'/S}\ra \calD_{X/S}$ (\ref{eqDoma1}), extends to an action of $\widehat{S}^{\bullet}\calT_{X'/S}$ on $F_1^*\widehat{S}^{\bullet}\calT_{X'/S}.$
\end{proposition}

\begin{proof}
Denote by $\ov{\calI}_0$ the PD-ideal of $\calP_0.$
Consider the composition
\begin{equation}\label{era4tak5}
v:S^{\bullet}\calT_{X'/S} \xrightarrow{\psi} \calD_{X/S} \ra \mathscr{Hom}_{\Ox_X}(\calP_0,\Ox_X),
\end{equation}
where the second arrow is induced by the canonical projections $\calP_0 \ra \calP_0/\ov{\calI}_0^{[n]}$ for $n\ge 1.$

Given a local section $x$ of $F_1^*\Gamma^{\bullet}\omega^1_{X'/S},$ an integer $n\ge 0$ and a differential operator $\partial:\calP_0^n=\calP_0/\ov{\calI}_0^{[n+1]}\ra \Ox_X$ of $\calD^n_{X/S},$ the action $\partial\cdot x$ of $\partial$ on $x$ by the stratification defined in \ref{propAkaza} is the image of $x$ by the composition
$$
F_1^*\Gamma^{\bullet}\omega^1_{X'/S} \xrightarrow{\Delta} F_1^*\Gamma^{\bullet}\omega^1_{X'/S}\otimes_{\Ox_X} F_1^*\Gamma^{\bullet}\omega^1_{X'/S} \xrightarrow{\op{Id} \otimes (\partial_{\epsilon_i}\circ u)} F_1^*\Gamma^{\bullet}\omega^1_{X'/S},
$$
where $u$ is given in \eqref{uAkaza} and $\Delta$ is the comultiplication map.
The action of $S^{\bullet}\calT_{X'/S}$ on $F_1^*\Gamma^{\bullet}\omega^1_{X'/S},$ induced by $\psi:S^{\bullet}\calT_{X'/S} \ra \calD_{X/S}$ and the action of $\calD_{X/S}$ on $F_1^*\Gamma^{\bullet}\omega^1_{X'/S},$ is then given by
\begin{equation}\label{era4takaction}
\begin{array}[t]{clc}
S^{\bullet}\calT_{X'/S} \times F_1^*\Gamma^{\bullet}\omega^1_{X'/S} & \ra & F_1^*\Gamma^{\bullet}\omega^1_{X'/S} \\
(\partial',x) & \mapsto & \psi(\partial')\cdot x=(\op{Id} \otimes v(\partial')\circ u)\circ \Delta(x).
\end{array}
\end{equation}
Since $v$ extends to $\widehat{S}^{\bullet}\calT_{X'/S}$ (\ref{era4tak2}), the action \eqref{era4takaction} extends to an action of $\widehat{S}^{\bullet}\calT_{X'/S}$ on $F_1^*\Gamma^{\bullet}\omega^1_{X'/S}.$ The result then follows by duality.
\end{proof}

\newpage

\section[Full faithfulness]{Full faithfulness of $\Psi$}

\begin{lemma}\label{dirsum}
Let $u:M\ra N$ be a Kummer morphism of fs monoids \eqref{defkummer}. Then $\Z[N \backslash u(M)]$ is a $\Z[M]$-submodule of $\Z[N].$
\end{lemma}

\begin{proof}
We have to prove that for any $a\in M$ and $b\in N,$ if $u(a)+b\in u(M)$ then $b\in u(M).$ So let $a,t\in M$ and $b\in N$ such that $u(a)+b=u(t).$ Since $u$ is Kummer, there exists a positive integer $n$ and $x\in M$ such that $nb=u(x).$ Then $u(x)=u(n(t-a))$ and so $x=n(t-a)\in M.$ It follows that $t-a\in M^{sat}=M$ and so $b=u(t-a)\in u(M).$
\end{proof}

\begin{lemma}\label{dirsum2}
Let $u:M\ra N$ be a Kummer morphism of fs monoids. Then the morphism
$$v:\begin{array}[t]{clc}
N^{gp}\oplus N^{gp} & \ra & N^{gp} \oplus (N^{gp}/u^{gp}(M^{gp})) \\
(x,y) & \mapsto & (x+y,\ov{y})
\end{array}$$
induces an isomorphism
$$(N\oplus_MN)^{sat}\xrightarrow{\sim}N\oplus (N^{gp}/u^{gp}(M^{gp})),$$
whose inverse is induced by
$$w:\begin{array}[t]{clc}
N^{gp}\oplus N^{gp} & \ra & N^{gp}\oplus_{M^{gp}} N^{gp} \\
(x,y) & \mapsto & (x-y,y).
\end{array}$$
\end{lemma}

\begin{proof}
For any $x\in N^{gp}$ and $t\in M^{gp},$ $v(u^{gp}(t),-u^{gp}(t))=0$ and $w(x,u^{gp}(t))=(x,0)$ so $v$ and $w$ induce morphisms
$$v:\begin{array}[t]{clc}
N^{gp}\oplus_{M^{gp}} N^{gp} & \ra & N^{gp} \oplus (N^{gp}/u^{gp}(M^{gp})) \\
(x,y) & \mapsto & (x+y,\ov{y}).
\end{array}$$
$$w:\begin{array}[t]{clc}
N^{gp}\oplus (N^{gp}/u^{gp}(M^{gp})) & \ra & N^{gp}\oplus_{M^{gp}} N^{gp} \\
(x,y) & \mapsto & (x-y,y).
\end{array}$$
Now let $(x,y)\in N^{gp}\oplus N^{gp}.$ If $\ov{(x,y)}\in (N\oplus_MN)^{sat}\subset N^{gp}\oplus_{M^{gp}}N^{gp}$ then there exists a positive integer $n,$ $t\in M^{gp}$ and $z,s \in N$ such that
$$n(x,y)=(z,s)+(u^{gp}(t),-u^{gp}(t)) \in N^{gp}\oplus N^{gp}.$$
Then
$$n(x+y)=z+s\in N.$$
Since $N$ is saturated, $x+y\in N$ and so $v$ induces a morphism
$$v:\begin{array}[t]{clc}
(N\oplus_{M} N)^{sat} & \ra & N \oplus (N^{gp}/u^{gp}(M^{gp})).\end{array}$$
Let $x\in N$ and $y\in N^{gp}.$ There exists a positive integer $n$ such that $ny\in u^{gp}(M^{gp}).$ Then $n(x-y,y)=(nx,0)$ in $N^{gp}\oplus_{M^{gp}}N^{gp}$ and so $w$ induces a morphism
$$w:\begin{array}[t]{clc}
N\oplus (N^{gp}/u^{gp}(M^{gp})) & \ra & (N\oplus_{M} N)^{sat}.
\end{array}$$
The morphisms $v$ and $w$ are clearly inverse to each other.
\end{proof}

\begin{lemma}\label{exactseq}
Let $u:M\ra N$ be a Kummer morphism of fs monoids. Consider the morphism
$$d:\begin{array}[t]{clc}\Z [N] & \ra & \Z[N]\otimes_{\Z} \Z\left [ N^{gp}/u^{gp}(M^{gp})\right ] \\
e^n & \mapsto & e^n\otimes e^{\ov{n}}-e^n\otimes 1.\end{array}$$
Then the sequence of $\Z[M]$-modules
\begin{equation}
0 \ra \Z[M] \xrightarrow{\Z [u]} \Z[N] \xrightarrow{d} \Z[N]\otimes_{\Z} \Z\left [ N^{gp}/u^{gp}(M^{gp})\right ]
\end{equation}
is homotopic to zero i.e. there exists morphisms
$$s:\Z[N] \ra \Z[M]$$
and
$$s':\Z[N]\otimes_{\Z} \Z\left [ N^{gp}/u^{gp}(M^{gp})\right ] \ra \Z[N]$$
such that $s \circ \Z[u]=\op{Id}_{\Z[M]}$ and $s' \circ d+\Z[u] \circ s=\op{Id}_{\Z[N]}.$
\end{lemma}

\begin{proof}
By \ref{dirsum}, the $\Z[M]$-module $\Z[N]$ decomposes as
$$\Z[N]=\Z[M]\oplus \Z[N\backslash M].$$
Let
$$s:\Z[N] \ra \Z[M]$$
be the projection on $\Z[M]$ and let
$$s':\Z[N]\otimes_{\Z} \Z\left [ N^{gp}/u^{gp}(M^{gp})\right ] \ra \Z[N]$$
be the base change of
$$\Z\left [N^{gp}/u^{gp}(M^{gp}) \right ] \ra \Z,\ e^{\ov{n}} \mapsto \begin{cases}1\ \op{if}\ \ov{n}\neq 0\\ 0\ \op{else} \end{cases}$$
by $\Z \hookrightarrow \Z[N].$ 
Then
$$s\circ \Z[u]= \op{Id}_{\Z[M]}$$
and
$$s' \circ d+d\circ s= \op{Id}_{\Z[N]}.$$
\end{proof}

\begin{proposition}[\cite{Kat19} 3.4.1]\label{Katolftop}
Let $u:M\ra N$ be a Kummer morphism of fs monoids and $T$ an fs logarithmic scheme equipped with a chart $T\ra A[M]$ and $T'=T\times_{A[M]}A[N].$
Let $\calE$ be a quasi-coherent $\Ox_T$-module, $T''=T'\times_{T}^{\op{log}}T'$ and $p:T' \ra T$ and $q_1,q_2:T'' \ra T'$ the canonical projections. Then the sequence of $\Ox_{T}$-modules
$$0 \ra \calE \ra p_*p^*\calE \xrightarrow{p_*(q_2^{\#}-q_1^{\#})} (p\circ q_1)_*(p\circ q_1)^*\calE$$
is homotopic to zero.
\end{proposition}

\begin{proof}
We can suppose that $T$ is affine. Since $p$ is affine, $T'$ and $T''$ are both affine. Let $E$ be the module corresponding to $\calE.$ It is sufficient to prove that the sequence
$$0 \ra E \ra \Gamma(U,\calE)\otimes_{\Z[M]} \Z[N] \ra E\otimes_{\Z[M]}\Z[(N\oplus_MN)^{sat}]$$
is homotopic to $0.$ It is thus sufficient to prove that the sequence
$$0 \ra \Z[M] \xrightarrow{\Z[u]} \Z[N] \xrightarrow{w} \Z[(N\oplus_MN)^{sat}]$$
is homotopic to zero, where
$$
w:\begin{array}[t]{clc}
\Z[N] & \ra & \Z[(N\oplus_MN)^{sat}] \\
e^n & \mapsto & e^{(0,n)}-e^{(n,0)}.
\end{array}$$
By \ref{dirsum2}, there exists an isomorphism
$$v:\Z[(N\oplus_MN)^{sat}]  \xrightarrow{\sim} \Z [N\oplus N^{gp}/u^{gp}(M^{gp})] \xrightarrow{\sim} \Z[N] \otimes_{\Z} \Z[N^{gp}/u^{gp}(M^{gp})]$$
such that
$$v\circ w(e^n)=e^n\otimes e^{\ov{n}}-e^n\otimes 1 \ \forall n\in N.$$
The result then follows from \ref{exactseq}.
\end{proof}

\begin{proposition}[\cite{INT} 1.3]\label{INTlog}
Let $f:X\ra S$ be a morphism of fs logarithmic schemes which is log flat, of Kummer type and locally of finite presentation, $x\in X$ and $y=f(x).$ Then, fppf locally around $x$ and $y,$ there exists a chart $\theta:P\ra Q$ of $f$ such that
\begin{enumerate}
\item $\theta$ is a Kummer morphism.
\item The morphism $f_1:X \ra T=S\times_{A[P]} A[Q],$ induced by $f$ and the chart $X \ra A[Q],$ is flat, surjective and locally of finite presentation.
\end{enumerate}
\end{proposition}

\begin{theorem}\label{logflatdescent}
Let $f:X \ra S$ be a quasi-compact surjective morphism of fs logarithmic schemes satisfying the following condition: 
For any $x\in X$ and $y=f(x),$ there exists, fppf locally around $x$ and $y,$ a chart $\theta:P\ra Q$ of $f$ satisfying:
\begin{enumerate}
\item $\theta$ is a Kummer morphism.
\item The morphism $f_1:X \ra T=S\times_{A[P]} A[Q],$ induced by $f$ and the chart $X \ra A[Q],$ is faithfully flat.
\end{enumerate}
Let $q_1,q_2:X\times_S^{\op{log}}X \ra X$ be the canonical projections and $g=f\circ q_1=f\circ q_2$ the structural morphism. For any quasi-coherent $\Ox_{S}$-module $\calE,$ the sequence
$$0 \ra \calE \ra f_*f^*\calE \rightarrow g_*g^*\calE,$$
where the second arrow is the difference of the morphisms induced by $q_1$ and $q_2,$ is exact.
\end{theorem}

\begin{proof}
We prove that for any $y\in S,$ the sequence
$$0 \ra \calE_y \ra (f_*f^*\calE)_y \ra (g_*g^*\calE)_y$$
is exact.
Let $y\in S.$ Since $f$ is surjective, there exists $x\in X$ such that $f(x)=y.$ By \ref{lemwiwfppf} we can suppose that $f$ has a chart $\theta:P\ra Q$ such that $\theta:P \ra Q$ is a Kummer morphism of fs monoids and the morphism $f_1:X \ra T=S\times_{A[P]} A[Q],$ induced by $f$ and the chart $X \ra A[Q],$ is faithfully flat and locally of finite presentation.
Let $f_2:T \ra S,$ $p_1,p_2:X\times_{T}^{\op{log}}X \ra X$ and $r_1,r_2:T \times_{S}^{\op{log}}T \ra T$ be the canonical projections and $g_1:X\times_T^{\op{log}}X\ra T$ and $g_2:T\times_S^{\op{log}}T \ra S$ the structural morphisms. Note that $f_1:X\ra T$ is strict so $X\times_T^{\op{log}}X=X\times_TX.$
By \ref{Katolftop}, the sequence
\begin{equation}\label{exseq1}
0 \ra \calE \ra f_{2*}f_2^*\calE \ra g_{2*}g_2^*\calE,
\end{equation}
where the second arrow is the difference between the morphisms induced by $r_1$ and $r_2,$ is exact. The morphism $f_1$ is quasi-compact and faithfully flat, so, by (\cite{Raynaud71} Exposé VIII 1.7), the sequence
\begin{equation}\label{exseq2}
0 \ra f_2^*\calE \ra f_{1*}f^*\calE \ra g_{1*}g_1^*f_2^*\calE
\end{equation}
is exact.
The morphism $f_1\times_{S}f_1:X\times_{S}X \ra T\times_{S}T$ is strict and faithfully flat so $f_1\times_{S}^{\op{log}}f_1:X\times_{S}^{\op{log}}X \ra T\times_{S}^{\op{log}}T$ is also faithfully flat (see \ref{parag42} for the notations $\times_S^{\op{log}}$ and $\times_S$). It follows that the canonical morphism
\begin{equation}\label{exmor1}
g_2^*\calE \ra (f_1\times_S^{\op{log}}f_1)_*(f_1\times_S^{\op{log}}f_1)^*g_2^*\calE=(f_1\times_S^{\op{log}}f_1)_*g^*\calE
\end{equation}
is injective. Now considering the commutative diagram
$$
\begin{tikzcd}
X\times_T^{\op{log}}X \ar{r}{h} \ar{d}{g_1} & X\times_S^{\op{log}}X \ar{d}{g} \\
T\ar{r}{f_2} & S,
\end{tikzcd}
$$
we have a canonical morphism
\begin{equation}\label{exmor2}
g^*\calE \ra h_*h^*g^*\calE=h_*g_1^*f_2^*\calE.
\end{equation}
We have the following commutative diagram:
$$
\begin{tikzcd}
 & X\times_T^{\op{log}}X \ar{dl}{h} \ar{d} \ar[bend right =-60]{ddr}{g_1}  & \\
X\times_S^{\op{log}}X \ar{r} \ar{d} \ar{dr}{f_1\times_S^{\op{log}}f_1} & T\times_S^{\op{log}}X \ar{r} \ar{d} & X \ar[swap]{d}{f_1} \ar[bend right=-30]{dd}{f} \\
X\times_S^{\op{log}}T \ar{r} \ar{d} & T\times_S^{\op{log}}T \ar{r} \ar{d} \ar[swap]{dr}{g_2} & T\ar[swap]{d}{f_2} \\
X \ar[swap]{r}{f_1} & T \ar[swap]{r}{f_2} & S
\end{tikzcd}
$$
It proves that the exact sequences \eqref{exseq1} and \eqref{exseq2} and the morphisms \eqref{exmor1} and \eqref{exmor2} fit into the commutative diagram
$$
\begin{tikzcd}
 & & 0\ar{d} & 0\ar{d} \\
0 \ar{r} & \calE \ar{r} \ar{dr} & f_{2*}f_2^*\calE \ar{r} \ar{d} & g_{2*}g_2^*\calE \ar{d} \\
 & & f_*f^*\calE \ar{r} \ar{d} & g_*g^*\calE \ar{dl} \\
& & f_{2*}g_{1*}g_1^*f_2^*\calE &
\end{tikzcd}
$$
The assertion follows then from a simple diagram chase.
\end{proof}

\begin{lemma}\label{lemwiwfppf}
Keep the hypothesis of \ref{logflatdescent}. Let $(S_i \xrightarrow{u_i}S)_{i\in I}$ be an fppf covering and consider, for any $i\in I,$ the cartesian square
$$
\begin{tikzcd}
X_i \ar{r}{v_i} \ar{d}{f_i} & X \ar{d}{f} \\
S_i \ar{r}{u_i} & S.
\end{tikzcd}
$$
Let $g_i:X_i\times_{S_i}^{\op{log}}X_i \ra S_i$ be the structural morphism for any $i\in I.$
If, for any $i\in I,$ the sequence of $\Ox_{S_i}$-modules
$$0 \ra u_i^*\calE \ra f_{i*}f_i^*u_i^*\calE \ra g_{i*}g_i^*u_i^*\calE$$
is exact, then the sequence of $\Ox_S$-modules
$$0 \ra \calE \ra f_*f^*\calE \ra g_*g^*\calE$$
is exact.
\end{lemma}

\begin{proof}
First, note that the following squares are cartesian in the category of logarithmic schemes:
$$
\begin{tikzcd}
X_i \ar{r}{v_i} \ar{d}{f_i} & X \ar{d}{f} & X_i\times_{S_i}^{\op{log}}X_i \ar{r}{w_i} \ar{d}{g_i} & X\times_S^{\op{log}}X \ar{d}{g} \\
S_i \ar{r}{u_i} & S & S_i\ar{r}{u_i} & S
\end{tikzcd}
$$
The exact sequence
$$0 \ra u_i^*\calE \ra f_{i*}f_i^*u_i^*\calE \ra g_{i*}g_i^*u_i^*\calE$$
becomes equal to
$$0 \ra u_i^*\calE \ra f_{i*}v_i^*f^*\calE \ra g_{i*}w_i^*g^*\calE.$$
By base change, we obtain the exact sequence
$$0 \ra u_i^*\calE \ra u_i^*f_*f^*\calE \ra u_i^*g_*g^*\calE.$$
The result then follows by faithfully flat descent.
\end{proof}

\begin{definition}[\cite{Kat19} 2.3]\label{logflattop}
Let $T$ be an fs logarithmic scheme. A family of morphisms $(T_i \xrightarrow{f_i}T)_{i\in I}$ of fs logarithmic schemes is said to be \emph{a covering for the log flat topology}, if the following conditions are satisfied:
\begin{enumerate}
\item $f_i$ is log flat, locally of finite presentation and of Kummer type for any $i\in I.$
\item Set theoretically,
$$T=\bigcup_{i\in I}f_i(T_i).$$
\end{enumerate}
This defines a pretopology on the category of fs logarithmic schemes.
\end{definition}

\begin{exmp}\label{Flogflatcover}
Let $X \ra S$ be a log smooth morphism of fine logarithmic schemes of characteristic $p,$ $F:X\ra X'$ the exact relative Frobenius and $g:X' \ra S$ the canonical morphism (\ref{PFrob}). By \ref{FKummer}, $F$ is of Kummer type. By \ref{thmlogflat}, $F$ is log flat. Since $f=g\circ F$ and $f$ and $g$ are locally of finite presentation, $F$ is also locally of finite presentation. Finally, $F$ is a homeomorphism on the underlying topological spaces. We conclude that $(F:X\ra X')$ is a log flat covering.
\end{exmp}

\begin{theorem}\label{eraflogflatdescent}
Let $T$ be an fs logarithmic scheme, $(T_i \xrightarrow{f_i} T)_{i\in I}$ a log flat covering \eqref{logflattop} and $\calE$ a quasi-coherent $\Ox_T$-module. For any $i,j\in I,$ let $f_{ij}:T_i\times_T^{\op{log}}T_j \ra T$ be the canonical morphism. The sequence of $\Ox_T$-modules
$$0 \ra \calE \ra \prod_{i\in I}f_{i*}f_i^*\calE \ra \prod_{i,j\in I}f_{ij*}f_{ij}^*\calE,$$
where the last arrow is the difference between the morphisms induced by the projections $T_i\times_T^{\op{log}}T_j \ra T_i$ and $T_i\times_T^{\op{log}}T_j \ra T_j,$ is exact.
\end{theorem}

\begin{proof}
Let $X=\coprod_{i\in I}T_i$ and $f:X\ra T$ the morphism induced by $(f_i)_{i\in I}.$ Since the assertion is local on $T,$ and $f$ is open by (\cite{Kat19} 2.5), we can suppose that $I$ is finite. We then apply \ref{logflatdescent} to $f.$
\end{proof}

\begin{parag}
For the remaining of this section, we keep the set-up of section 12 i.e. we keep the notation and assumption of \ref{paragJapon1} and \ref{parag86}. We consider a perfect field of positive characteristic $p,$ denote its ring of Witt vectors by $W,$ and consider a log smooth morphism of framed logarithmic $p$-adic formal schemes $f:(\frakX,Q)\ra (\frakS,P),$ such that $\frakS$ is log flat and locally of finite type over $\op{Spf}W.$ Note that this implies, by \ref{Wflat}, that the formal schemes $\frakX$ and $\frakS$ are flat over $\op{Spf}W$ \eqref{dxuflat}.
We denote by $F_1:X \ra X'$ the exact relative Frobenius and suppose that it lifts to a morphism of framed fs $p$-adic logarithmic formal schemes $F:(\frakX,Q) \ra (\frakX',Q')$ over $(\frakS,P),$ such that $\frakX'$ is log smooth over $\frakS.$
\end{parag}

\begin{proposition}\label{propFnlogtop}
For any positive integer $n,$ the morphism $F_n:\frakX_n \ra \frakX'_n$ forms a log flat covering \eqref{logflattop}.
\end{proposition}

\begin{proof}
The lifting $F_n$ is clearly a homeomorphism on the underlying topological spaces.
We have the commutative diagram
$$
\begin{tikzcd}
F_n^{-1}\ov{\calM}_{\frakX'_n} \ar{r}{\sim} \ar[swap]{d}{F_n^{\flat}} & F_1^{-1}\ov{\calM}_{X'} \ar{d}{F_1^{\flat}} \\
\ov{\calM}_{\frakX_n} \ar{r}{\sim} & \ov{\calM}_X.
\end{tikzcd}
$$
Since $F_1$ is of Kummer type \eqref{FKummer}, so is $F_n.$ The fact that $F_n$ is locally of finite presentation follows from the fact that $F_n$ is an $\frakS_n$-morphism and $\frakX_n$ and $\frakX'_n$ are locally of finite presentation over $\frakS_n.$ Finally, the log flatness of $F_n$ follows from \ref{logflatfiber}.
\end{proof}

\begin{theorem}\label{THMX1}
Let $n$ be a positive integer. Denote by $1\text{-}\op{MHS}^{\text{qcoh}}(\frakX_n'/\frakS_n)$ and $\op{MHS}^{\text{qcoh}}(\frakX_n/\frakS_n)$ the full subcategories of $1\text{-}\op{MHS}(\frakX_n'/\frakS_n)$ and $\op{MHS}(\frakX_n/\frakS_n)$ consisting of quasi-coherent modules. Then, the functor
$$
\Psi_n:\begin{array}[t]{clc}
1\text{-}\op{MHS}^{\text{qcoh}}(\frakX_n'/\frakS_n) & \ra & \op{MHS}^{\text{qcoh}}(\frakX_n/\frakS_n)
\end{array}
$$
induced by \eqref{eq994}, is fully faithful.
\end{theorem}

\begin{proof}
Recall the diagram (\ref{totdiag1})
\begin{equation}
\begin{tikzcd}
 & & \frakX'\ar{d}{\iota'} \\
 & P_{\frakX/\frakS,0} \ar{r}{\varphi} \ar{d} & P_{\frakX'/\frakS,1}\ar{d} \\
\frakX\times_{\frakX'}^{\op{log}}\frakX \ar{ur}{\psi} \ar[bend right=-30]{uurr} \ar{r} & \frakX\times^{\op{log}}_{\frakS,[Q]}\frakX \ar{r}{F^2} & \frakX'\times_{\frakS,[Q']}^{\op{log}}\frakX'
\end{tikzcd}
\end{equation}
The functor $\Psi_n$ is then defined by
$$\Psi_n(\calE',\epsilon')=(F_n^*\calE',\varphi_n^*\epsilon')$$
for any object $(\calE',\epsilon')$ of $1\text{-}\op{MHS}^{\text{qcoh}}(\frakX_n'/\frakS_n).$
Let $(\calE'_1,\epsilon'_1)$ and $(\calE'_2,\epsilon'_2)$ be two objects of $1\text{-}\op{MHS}^{\text{qcoh}}(\frakX_n'/\frakS_n),$ $(\calE_1,\epsilon_1)$ and $(\calE_2,\epsilon_2)$ their images by $\Psi_n$ and
$$u:(\calE_1,\epsilon_1) \ra (\calE_2,\epsilon_2)$$
a morphism of $\op{MHS}^{\text{qcoh}}(\frakX_n/\frakS_n).$
Let $p_1,p_2:\left (P_{\frakX/\frakS,0}\right )_n \ra \frakX_n$ and $q_1,q_2:\frakX_n \times_{\frakX_n'}^{\op{log}} \frakX_n \ra \frakX_n$ be the canonical projections.
Since $u$ is a morphism of stratified modules, the diagram to the left is commutative and that implies the commutativity of the diagram to the right:
$$
\begin{tikzcd}
p_2^*\calE_1 \ar{r}{\epsilon_1} \ar[swap]{d}{p_2^*u} & p_1^*\calE_1 \ar{d}{p_1^*u} & & q_2^*\calE_1 \ar[swap]{d}{q_2^*u} \ar{r}{\psi_n^*\epsilon_1} & q_1^*\calE_1 \ar{d}{q_1^*u} \\
p_2^*\calE_2 \ar{r}{\epsilon_2} & p_1^*\calE_2 & & q_2^*\calE_2 \ar{r}{\psi_n^*\epsilon_2} & q_1^*\calE_2
\end{tikzcd}
$$
We deduce the commutativity of
$$
\begin{tikzcd}
0 \ar{r} & \calE_1' \ar{r} \ar[dashed]{d}{u'} & \calE_1  \ar{r} \ar{d}{u} & \calE'_1\otimes_{\Ox_{\frakX_n'}}\Ox_{\frakX_n \times_{\frakX_n'}^{\op{log}}\frakX_n} \ar{d}{q_2^*u-q_1^*u} \\
0 \ar{r} & \calE_2' \ar{r} & \calE_2 \ar{r} & \calE'_2\otimes_{\Ox_{\frakX_n'}}\Ox_{\frakX_n \times_{\frakX_n'}^{\op{log}}\frakX_n}
\end{tikzcd}
$$
The rows are exact by \ref{eraflogflatdescent} and \ref{propFnlogtop}. We deduce the existence of the morphism $u':\calE'_1 \ra \calE'_2.$ It remains to prove that $u'$ is compatible with the stratifications $\epsilon_1'$ and $\epsilon_2'.$ Since $\varphi$ is log flat \eqref{Khaminei145}, it is sufficient, by \ref{eraflogflatdescent}, to prove that $F^*u'=u$ is compatible with $\varphi^*\epsilon'_i=\epsilon_i.$ This is true by definition of $u$ and the full faithfulness of $\Psi_n$ follows.
\end{proof}

\begin{theorem}\label{THMX2}
Let $n$ be a positive integer. The functor $\Phi_n$ \eqref{krazphin} preserves quasi-nilpotent connections, the diagram
$$
\begin{tikzcd}
1\text{-}\op{MHS}^{qcoh}(\frakX_n'/\frakS_n) \ar{rr}{\Psi_n} \ar[swap,sloped]{d}{\sim} & & \op{MHS}^{qcoh}(\frakX_n/\frakS_n) \ar[sloped]{d}{\sim} \\
p\text{-}\op{MIC}^{qcoh,qn}(\frakX_n'/\frakS_n) \ar{rr}{\Phi_n} & & \op{MIC}^{qcoh,qn}(\frakX_n/\frakS_n),
\end{tikzcd}
$$
where qn denotes quasi-nilpotent objects \eqref{krazphin}, is commutative and $\Phi_n$ is fully faithful.
\end{theorem}

\begin{proof}
It is sufficient to prove that the diagram
$$
\begin{tikzcd}
1\text{-}\op{MHS}^{qcoh}(\frakX_n'/\frakS_n) \ar{rr}{\Psi_n} \ar[hook]{d} & & \op{MHS}^{qcoh}(\frakX_n/\frakS_n) \ar[hook]{d} \\
p\text{-}\op{MIC}^{qcoh}(\frakX_n'/\frakS_n) \ar{rr}{\Phi_n} & & \op{MIC}^{qcoh}(\frakX_n/\frakS_n)
\end{tikzcd}
$$
is commutative.
Let $(\calE',\epsilon')$ be an object of the upper left category. Its image by $\Psi_n$ is
$$
(\calE,\epsilon)=\Psi_n(\calE',\epsilon')=(F_n^*\calE',\varphi_n^*\epsilon').
$$
Suppose the hypothesis \ref{loccoord} is satisfied. By \ref{lem12}, for every $1\le i\le d,$ there exists a local invertible section $u_i$ of $\calM_{\frakX}$ and a local section $b_i$ of $\Ox_{\frakX}$ such that $F^{\flat}(\widetilde{m}_i')=p\widetilde{m}_i+u_i$ and $\alpha_{\frakX}(u_i)=1+pb_i.$
For $1\le i\le d,$ by \eqref{eqKoko1323}, we have
$$
\varphi^{\#}\left (\frac{\widetilde{\eta}_i'}{p} \right )=\left ((p-1)!\widetilde{\eta}_i^{[p]}+\sum_{k=1}^{p-1}\frac{(p-1)!}{k!(p-k)!}\widetilde{\eta}_i^k \right )\alpha_{\frakY}(p_2^{\flat}u_i-p_1^{\flat}u_i)+\frac{p_2^{\#}b_i-p_1^{\#}b_i}{1+pp_1^{\#}b_i},
$$
where $p_1,p_2:\frakY \ra \frakX$ are the canonical projections.
Denote by $\widehat{\eta}_i$ (resp. $\widehat{\eta}_i'$) the class of $\widetilde{\eta}_i$ (resp. $\frac{\widetilde{\eta}_i'}{p}$) modulo $p^n,$ by $\calP'$ the structural ring of $P_{\frakX'/\frakS,1},$ by $\calP'_n$ its reduction modulo $p^n$ and by $\calI'_n$ its PD ideal. Similarly, let $\calI$ be the PD ideal of $P_{\frakX/\frakS,0}$ and $\calI_n$ its reduction modulo $p^n.$ Denote by $\calP$ the structure sheaf of $P_{\frakX/\frakS,0}$ and by $\calP_n$ its reduction modulo $p^n.$ Then $\left (\widehat{\eta}'^{[I]} \right )_{I\in \N^d}$ is a basis of the $\Ox_{\frakX'_n}$-module $\calP'_n.$ Denote by $(\partial'_I)$ its dual basis. It follows that, for a local section $x'$ of $\calE',$
$$
\epsilon':\calP'_n \otimes_{\Ox_{\frakX'_n}}\calE' \ra \calE' \otimes_{\Ox_{\frakX'_n}}\calP'_n,\ 1\otimes x'\mapsto \sum_{I\in\N^d}(\partial_I'\cdot x')\otimes \widehat{\eta}'^{[I]},
$$
$$
\varphi_n^*\epsilon': \calP_n \otimes_{\Ox_{\frakX_n}}(F_n^*\calE') \ra (F_n^*\calE')\otimes_{\Ox_{\frakX_n}} \calP_n,\ 1\otimes 1\otimes x' \mapsto \sum_{I \in \N^d}1\otimes \left (\partial'_I\cdot x' \right )\otimes \varphi_n^{\#}(\widehat{\eta}'^{[I]}).
$$
Composing $\varphi_n^*\epsilon'$ with the canonical projection
$$
\pi:\left (F_n^*\calE'\right ) \otimes _{\Ox_{\frakX_n}} \calP_n \ra \left (F_n^*\calE'\right ) \otimes_{\Ox_{\frakX_n}} \calP_n/\calI_n^{[2]},
$$
and since the image of $\alpha_{\frakY}(p_2^{\flat}u_i-p_1^{\flat}u_i)-1$ in $\calP$ belongs to $\calI,$ we get
$$
\left (\pi\circ \varphi_n^*\epsilon'\right ) \left (1\otimes 1\otimes x' \right ) =  1\otimes x'\otimes 1+\sum_{i=1}^d1\otimes \left (\partial_{\epsilon_i}'\cdot x' \right )\otimes \left (\widehat{\eta}_i+\frac{p_{2}^{\#}b_i-p_1^{\#}b_i}{1+pp_1^{\#}b_i} \right ).
$$
Identifying $\calI_n/\calI_n^{[2]}$ with $\omega^1_{\frakX_n/\frakS_n}$ via the canonical isomorphism \eqref{eqtakrizIIsquare}, we get 
$$
\left (\pi\circ \varphi_n^*\epsilon'\right ) \left (1\otimes 1\otimes x' \right ) -  1\otimes x'\otimes 1 = \sum_{i=1}^d1\otimes \left (\partial_{\epsilon_i}'\cdot x' \right )\otimes \left (\op{dlog}m_i+\frac{db_i}{1+pb_i} \right ).
$$
Let $\nabla'$ be the $p$-connection corresponding to $\epsilon'.$ Denote by $\calI'$ be the PD-ideal of $\calP_{\frakX'/\frakS,1}$ and $\calI'_n$ its reduction modulo $p^n.$ Then, by \eqref{KhamineiBeta}, we have
$$
\nabla':\begin{array}[t]{clclc}
\calE' & \ra & \calE'\otimes_{\Ox_{\frakX_n}} \left (\calI'_n/\calI_n'^{[2]} \right ) & \xrightarrow{\sim} & \calE' \otimes_{\Ox_{\frakX_n}} \omega^1_{\frakX_n'/\frakS_n} \\
x' & \mapsto & \sum_{i=1}^d(\partial_{\epsilon_i}'\cdot x') \otimes \ov{\widehat{\eta}_i'} & \mapsto & \sum_{i=1}^d(\partial_{\epsilon_i}'\cdot x') \otimes \op{dlog}m_i'.
\end{array}
$$
Set
$$
(\calE,\nabla)=\Phi_n(\calE',\nabla').
$$
If $x'$ is a local section of $\calE',$ then, by definition of $\Phi_n$ and \eqref{surp2},
$$
\nabla(1\otimes x')=\sum_{i=1}^d1\otimes (\partial_{\epsilon_i}'\cdot x') \otimes \left (\op{dlog}m_i +\frac{db_i}{1+pb_i} \right ).
$$
\end{proof}

\newpage

\section{An indexed Shiho functor}

\begin{parag}\label{indparag91}
In this section, we will provide an indexed version of the Shiho functor. However, we will only do it in characteristic $p$ for the time being. We fix a perfect field $k$ of positive characteristic $p$ and a log smooth morphism $f:\frakX\ra \frakS$ of fs $p$-adic logarithmic formal schemes flat over the ring of Witt vectors $W$ of $k.$ We keep the notation introduced in section 9. We suppose that we have a lifting $F:\frakX\ra \frakX'$ of the exact relative Frobenius $F_1:X\ra X'.$ We also denote by $\calI$ the sheaf $\ov{\calM}^{gp}_X.$ Consider the functors introduced in \eqref{indeq71}:
\begin{alignat}{2}
& j_{\calI !}:\calI_{\text{ét}}\ra X_{\text{ét}},\ j_{\calI}^*:X_{\text{ét}}\ra \calI_{\text{ét}},\ j_{\calI *}:\calI_{\text{ét}}\ra X_{\text{ét}}
\end{alignat}
We identify the étale sites $X_{\text{ét}}$ and $X'_{\text{ét}}$ via the universal homeomorphism $F_1:X\ra X'.$
Following the conventions of \ref{indpar2}, we denote $j_{\calI}^*\Ox_X$ (resp. $j_{\calI}^*\omega^1_{X/S},$ $j_{\calI}^*\Ox_{X'},$ $j_{\calI}^*\omega^1_{X'/S}$) by $\Ox_{X,\calI}$ (resp. $\omega^1_{X/S,\calI},$ $\Ox_{X',\calI},$ $\omega^1_{X'/S,\calI}$).

Recall the definition of the $\calI$-indexed algebra $\calA_X$ (\ref{Agroupe}) and the $\calI$-indexed subalgebra $\calB_{X/S}$ (\ref{Bgroupe}). The $\Ox_{X,\calI}$-module structure on $\calA_X$ and the morphism $\Ox_{X',\calI} \ra \Ox_{X,\calI},$ induced by $F_1:X\ra X',$ equip $\calB_{X/S}$ with an $\Ox_{X',\calI}$-module structure.

We denote by $\op{HIG}(\calB_{X/S})$ the category of $\calI$-indexed $\calB_{X/S}$-modules $\calE'$ (\ref{inddef} 3 and \ref{Bgroupe2}) equipped with a $\calB_{X/S}$-linear Higgs field (\ref{indhiggs})$$\theta:\calE' \ra \calE'\otimes_{\Ox_{X',\calI}}\omega^1_{X'/S,\calI}.$$ We also denote by $\op{MIC}(\calA_X)$ the category of $\calI$-indexed $\calA_X$-modules $\calE$ (\ref{inddef} 3 and \ref{Agroupe}) equipped with an admissible integrable connection (\ref{Bgroupe}) $$\nabla:\calE \ra \calE \otimes_{\Ox_{X,\calI}}\omega^1_{X/S,\calI}.$$
Let $(\calE',\theta)$ be an object of $\op{HIG}(\calB_{X/S}).$ We denote by $p_1,p_2:\calI^2\ra \calI$ the canonical projections and by 
$$p_i^*:\boldsymbol{\op{Mod}_{\Ox_{X,\calI}}} \ra \boldsymbol{\op{Mod}_{\Ox_{X,\calI^2}}}$$
$$p_i'^*:\boldsymbol{\op{Mod}_{\Ox_{X',\calI}}} \ra \boldsymbol{\op{Mod}_{\Ox_{X',\calI^2}}}$$
the functors they induce.
Let $\zeta$ be the composition
\begin{alignat}{2}
& p_1'^*(\calE' \otimes_{\Ox_{X',\calI}}\omega^1_{X'/S,\calI}) \\
\xrightarrow{\sim} & \ (p_1'^*\calE') \otimes_{\Ox_{X',\calI^2}}\omega^1_{X'/S,\calI^2}  \label{eq914} \\
\ra & \ (p_1'^*\calE') \otimes_{\Ox_{X',\calI^2}}\Ox_{X,\calI^2} \otimes_{\Ox_{X',\calI^2}}\omega^1_{X'/S,\calI^2}  \label{eq915} \\
\xrightarrow{\sim} & \ (p_1'^*\calE') \otimes_{\Ox_{X',\calI^2}}(F_1^*\omega^1_{X'/S})_{|\calI^2}
 \label{eq916} \\
\rightarrow &\ (p_1'^*\calE') \otimes_{\Ox_{X',\calI^2}}\omega^1_{X/S,\calI^2} \label{eq917}
\end{alignat}
where \eqref{eq914} and \eqref{eq916} are the canonical isomorphisms, \eqref{eq915} is $x\otimes \omega \mapsto x\otimes 1\otimes \omega,$ \eqref{eq917} is $\op{Id}_{p_1'^*\calE'}\otimes (p^{-1}dF_2)_{|\calI^2}$ and $p^{-1}dF_2$ is the morphism defined in \eqref{surp}.
For sections $\ov{s}$ and $\ov{t}$ of $\ov{\calM}^{gp}_X$ over an étale $X$-scheme $U,$
\begin{equation}
\zeta_{(\ov{s},\ov{t})}:\begin{array}[t]{clc}
\calE'_{\ov{s}}\otimes_{\Ox_{U'}}\omega^1_{U'/S} & \ra & \calE'_{\ov{s}}\otimes_{\Ox_{U'}}\omega^1_{U/S} \\
x\otimes \omega & \mapsto & x\otimes (p^{-1}dF_2)(\omega).
\end{array}
\end{equation}
\end{parag}

\begin{lemma}\label{indlem92}
Keep the same hypothesis and notations of the previous paragraph and assume that $\frakX_2\ra \frakS_2$ is log smooth. Then the morphism
\begin{equation}
\nabla_1:\begin{array}[t]{clc}\calE' \boxtimes \calA_X & \ra & (\calE' \boxtimes \calA_X) \otimes_{\Ox_{X,\calI^2}} (\omega^1_{X/S,\calI^2}) \\ 
x\otimes a & \mapsto & \zeta \left ( (p_1'^*\theta)(x)\right )\otimes a+x\otimes d_{\calA_X}(a)
\end{array}
\end{equation}
is well-defined and is an integrable connection on the $\Ox_{X,\calI^2}$-module $\calE' \boxtimes \calA_X.$ In addition, it induces an integrable admissible connection on the $\calI$-indexed $\calA_X$-module $\calE' \circledast_{\calB_{X/S}} \calA_X$
\begin{equation}\label{eq922}
\nabla: \calE' \circledast_{\calB_{X/S}} \calA_X \ra  (\calE' \circledast_{\calB_{X/S}} \calA_X) \otimes_{\Ox_{X,\calI}} \omega^1_{X/S,\calI}
\end{equation}
\end{lemma}

\begin{proof}
To prove the first part, it is sufficient to prove that for any sections $s,t$ of $\calM_X^{gp}$ over an étale $X$-scheme $U,$ with images $\ov{s}$ and $\ov{t}$ in $\calI=\ov{\calM}_X^{gp},$ the morphism
$$
\nabla_{1,(\ov{s},\ov{t})}:\begin{array}[t]{clc}
\calE'_{\ov{s}}\otimes_{\Ox_{U'}}\calA_{X,\ov{t}} & \ra & \calE'_{\ov{s}} \otimes_{\Ox_{U'}} \calA_{X,\ov{t}}\otimes_{\Ox_U}\omega^1_{U/S} \\
x\otimes a & \mapsto & \zeta_{(\ov{s},\ov{t})}(\theta_{\ov{s}}(x))\otimes a+x\otimes d_{\calA_X,\ov{t}}(a)
\end{array}
$$
is well-defined and is integrable, which is similar to that of \ref{lem85}. For the second part, denote by $\sigma:\calI^2\ra \calI$ the addition map. The adjunction morphism
$$\calE' \boxtimes \calA_X \ra \sigma^* \sigma_!(\calE' \boxtimes \calA_X)$$
induces a mophism
$$(\calE' \boxtimes \calA_X) \otimes_{\Ox_{X,\calI^2}} \omega^1_{X/S,\calI^2} \ra \sigma^* \sigma_!(\calE' \boxtimes \calA_X)\otimes_{\Ox_{X,\calI^2}} \omega^1_{X/S,\calI^2}.$$
Composing it with the connection $\nabla_1,$ we get a connection
$$\calE' \boxtimes \calA_X \ra \sigma^* \sigma_!(\calE' \boxtimes \calA_X)\otimes_{\Ox_{X,\calI^2}} \omega^1_{X/S,\calI^2}.$$
We have a canonical isomorphism
$$\sigma^* \sigma_!(\calE' \boxtimes \calA_X)\otimes_{\Ox_{X,\calI^2}} \omega^1_{X/S,\calI^2} \xrightarrow{\sim} \sigma^* \left ( \sigma_!(\calE' \boxtimes \calA_X)\otimes_{\Ox_{X,\calI}} \omega^1_{X/S,\calI} \right )$$
and so we get a connection
$$\nabla_2:\sigma_!(\calE' \boxtimes \calA_X) \ra \sigma_!(\calE' \boxtimes \calA_X)\otimes_{\Ox_{X,\calI}} \omega^1_{X/S,\calI}.$$
Since the Higgs field $\theta$ is $\calB_{X/S}$-linear and $d_{\calA_X}(\calB_{X/S})=0,$ the composition of $\nabla_2$ with the projection
$$\sigma_!(\calE' \boxtimes \calA_X)\otimes_{\Ox_{X,\calI}} \omega^1_{X/S,\calI} \ra (\calE' \circledast_{\calB_{X/S}} \calA_X)\otimes_{\Ox_{X,\calI}} \omega^1_{X/S,\calI}$$
factors through $\calE' \circledast_{\calB_{X/S}} \calA_X$ and we get a connection
$$\nabla:\calE' \circledast_{\calB_{X/S}} \calA_X \ra (\calE' \circledast_{\calB_{X/S}} \calA_X) \otimes_{\Ox_{X,\calI}} \omega^1_{X/S,\calI}.$$
More precisely, for a section $\ov{r}$ of $\calI$ over an étale $X$-scheme $U,$ the morphism $\nabla_{\ov{r}}$ is induced by
\begin{equation}\label{explicitconnimage}
\begin{array}[t]{clc}
\bigoplus\limits_{\ov{s}+\ov{t}=\ov{r}}\calE'_{\ov{s}}\otimes_{\Ox_{U'}}\calA_{X,\ov{t}} & \ra & \bigoplus\limits_{\ov{s}+\ov{t}=\ov{r}}\calE'_{\ov{s}}\otimes_{\Ox_{U'}}\calA_{X,\ov{t}}\otimes_{\Ox_{U}}\omega^1_{U/S} \\
x\otimes a & \mapsto & \nabla_{1,(\ov{s},\ov{t})}(x\otimes a).
\end{array}
\end{equation}
The properties of integrability and admissibility of $\nabla$ are inherited from $\nabla_1.$
\end{proof}

\begin{parag}
Keep the assumptions and notations of \ref{indparag91}. By \ref{indlem92}, we have a functor
\begin{equation}\label{Phiind}
\Phi^{ind}_1:\begin{array}[t]{clc} \op{HIG}(\calB_{X/S}) & \ra & \op{MIC}(\calA_X) \\
(\calE',\theta') & \mapsto & (\calE' \circledast_{\calB_{X/S}}\calA_X, \nabla) \end{array}
\end{equation}
where $\nabla$ is the connection defined in (\ref{eq922}).
\end{parag}

\newpage

\section{The ring of logarithmic differential operators as an Azumaya algebra}

\begin{parag}
In this section, we consider a log smooth morphism $f_1:(X,Q)\ra (S,P)$ of framed fs logarithmic schemes of characteristic $p$ and equip $S$ with the trivial PD structure on the zero ideal. Set
\begin{equation}
\upmu = \ov{\calM}_X^{gp}.
\end{equation}
Denote by $F_1:X\ra X'$ the exact relative Frobenius and by $P_0$ (resp. $P_0(2)$) the PD envelope of $X\ra X\times_S^{\op{log}}X$ (resp. $X \ra X\times_S^{\op{log}}X\times_S^{\op{log}}X$) (\cite{Kat89} 5.3). We identify the étale sites of $X,$ $X',$ $P_0$ and $P_0(2)$ via the universal homeomorphism $F_1$ and the nilimmersions $X\ra P_0$ and $X\ra P_0(2).$ Denote by $\calP_0$ the structural ring of $P_0$ and by $\ov{\calI}_{0}$ the PD ideal of $\calP_0.$ Following \ref{defmokrez1}, set $\calP^n_{0}=\calP_0/\ov{\calI}_{0}^{[n+1]},$ $\calD^n_{X/S}=\mathscr{Hom}_{\Ox_X}(\calP_0^n,\Ox_X)$ for any integer $n\ge 0$ and
$$\calD_{X/S}=\lim\limits_{\substack{\longrightarrow \\ n\ge 0}} \calD^n_{X/S},$$
where $\calP^n_{0}$ is considered as an $\Ox_X$-module via the first projection $P_0 \ra X.$
Recall that there exists a canonical isomorphism \eqref{P11}
$$\ov{\calI}_{0}/\ov{\calI}_{0}^{[2]} \xrightarrow{\sim} \omega^1_{X/S}.$$
Consider the $p$-curvature morphism (\ref{eqDoma1})
\begin{equation}\label{Takriz2}
\psi:S^{\bullet}\calT_{X'/S} \ra \calD_{X/S},
\end{equation}
where $S^{\bullet} \calT_{X'/S}$ is the symmetric algebra of $\calT_{X'/S}=\mathscr{Hom}_{\Ox_{X'}}(\omega^1_{X'/S},\Ox_{X'}).$
We denote the image of this morphism by $\calC_{X/S},$ which is a sub $\Ox_{X'}$-algebra of the center of $\calD_{X/S}.$
Finally, for any $\Ox_X$-module $\F,$ denote by $\F_{\upmu}$ the image of $\F$ by the functor $j^*_{\upmu}:X_{\text{ét}} \ra X_{\text{ét}/\upmu}$ defined in \ref{indpar2} and set
$$\widetilde{\calD}^n_{X/S}=\calA_X \otimes_{\Ox_{X,\upmu}}\calD^n_{X/S,\upmu}$$
\begin{equation}\label{Dtilda1713}
\widetilde{\calD}_{X/S}=\calA_X \otimes_{\Ox_{X,\upmu}}\calD_{X/S,\upmu},
\end{equation}
where $\calA_X$ is defined in \ref{paragAX}. In the next paragraph, we will equip $\widetilde{\calD}_{X/S}$ with a $\upmu$-indexed algebra structure.
Since the functors $j_{\upmu}^*$ and $(\calA_X\otimes_{\Ox_{X,\upmu}}-)$ commute with inductive limits, we have
$$\widetilde{\calD}_{X/S}=\lim\limits_{\substack{\longrightarrow \\ n\ge 0}}\widetilde{\calD}^n_{X/S}.$$
\end{parag}

\begin{parag}\label{era3parag142}
Since $f_1:X\ra S$ is log smooth, the PD envelope $\calP_0$ is a locally free $\Ox_X$-module. So, for every integer $n\ge 0,$ we have a canonical isomorphism
$$\widetilde{\calD}^n_{X/S}\xrightarrow{\sim} \mathscr{Hom}_{\Ox_{X,\upmu}}(\calP^n_{0,\upmu},\calA_X).$$
Let $p_1,p_2:\upmu^2 \ra \upmu$ be the canonical projections and $\sigma:\upmu^2 \ra \upmu$ the addition map.
Let $g_1:\calP_{0,\upmu^2}^n\ra p_1^*\calA_X$ and $g_2:\calP_{0,\upmu^2}^m\ra p_2^*\calA_X$ be sections of $p_1^*\mathscr{Hom}_{\Ox_{X,\upmu}}(\calP^n_{0,\upmu},\calA_X)=\mathscr{Hom}_{\Ox_{X,\upmu^2}}(\calP^n_{0,\upmu^2},p_1^*\calA_X)$ and $p_2^*\mathscr{Hom}_{\Ox_{X,\upmu}}(\calP^m_{0,\upmu},\calA_X)=\mathscr{Hom}_{\Ox_{X,\upmu^2}}(\calP_{0,\upmu^2}^m,p_2^*\calA_X)$ respectively, over an object $U \ra \upmu^2$ of $X_{\text{ét}/\upmu^2},$ where $U$ is an étale $X$-scheme. Let
$$
\delta^{m,n}:\calP_0^{m+n} \ra \calP_0^m\otimes_{\Ox_X}\calP_0^n
$$
be the morphism induced by the composition
$$
P_0\times_XP_0 \xrightarrow{\sim} P_0(2) \ra P_0,
$$
where the second morphism is the $(1,3)$-projection and the first isomorphism is induced by the canonical isomorphism
$$
(X\times^{\op{log}}_SX)\times^{\op{log}}_X(X\times^{\op{log}}_SX)\xrightarrow{\sim}X\times^{\op{log}}_SX\times^{\op{log}}_SX.
$$
Recall that $\calA_X$ is equipped with an indexed stratification (\cite{Mon} 4.1.1 page 46) which is a morphism of indexed algebras
$$
\epsilon_{\calA_X}:\calP_{0,\upmu}\otimes_{\Ox_{X,\upmu}} \calA_X \ra \calA_X \otimes_{\Ox_{X,\upmu}}\calP_{0,\upmu}.
$$
We define $g_2\circ g_1$ to be the following composition
\begin{alignat*}{2}
\calP_{0,\upmu^2}^{n+m} \xrightarrow{\delta_{\upmu^2}^{m,n}}\calP_{0,\upmu^2}^m\otimes_{\Ox_{U,\upmu^2}}\calP_{0,\upmu^2}^n & \xrightarrow{\op{Id}\otimes g_1} \calP_{0,\upmu^2}^m\otimes_{\Ox_{U,\upmu^2}}p_1^*\calA_U \\
& \xrightarrow{p_1^*\epsilon_{\calA_U}} p_1^*\calA_U\otimes_{\Ox_{U,\upmu^2}}\calP_{0,\upmu^2}^m \\
&\xrightarrow{\op{Id}\otimes g_2} \calA_U\boxtimes \calA_U \\
& \ra \sigma^*\calA_U.
\end{alignat*}
Then
\begin{equation}\label{lci}
g_2\circ g_1:\calP_{0,\upmu^2}^{n+m} \ra \sigma^*\calA_U
\end{equation}
is a local section of
$$\mathscr{Hom}_{\Ox_{X,\upmu^2}}(\calP_{0,\upmu^2}^{n+m},\sigma^*\calA_X)=\sigma^*\mathscr{Hom}_{\Ox_{X,\upmu}}(\calP_{0,\upmu}^{n+m},\calA_X)\xrightarrow{\sim} \sigma^*\widetilde{\calD}^{n+m}_{X/S}.$$
\end{parag}

\begin{proposition}\label{era2stralgD}
The operation \eqref{lci} defines a $\upmu$-indexed algebra structure on the $\Ox_{X,\upmu}$-module $\widetilde{\calD}_{X/S}.$
\end{proposition}

\begin{proof}
The hard part is to check the associativity of this operation. Let $q_1,q_2,q_3:\upmu^3\ra \upmu$ be the canonical projections, $n_1,n_2,n_3\ge 0$ integers and
$$g_i\in q_i^*\mathscr{Hom}_{\Ox_X}(\calP_{0,\upmu}^{n_i},\calA_X),\ 1\le i\le 3$$
be local sections. We want to check that
$$g_3\circ (g_2 \circ g_1)=(g_3\circ g_2)\circ g_1.$$
To lighten notations in the following diagrams, we consider $g_i$ as a section of
$$q_i^*\mathscr{Hom}_{\Ox_X}(\calP_{0,\upmu},\calA_X)$$
by composing with the canonical projection $\calP_{0,\upmu^2} \ra \calP_{0,\upmu^2}^{n_i}.$
Let $r,s,t\in \Gamma(U,\upmu)$ be sections of $\upmu$ over an étale $X$-scheme $U.$ To lighten notation, we use $\calP_0$ instead of $\calP_{0|U}$ and $\epsilon_r$ instead of $\epsilon_{\calA_X,r}.$ Consider the following diagrams, where all tensor products are over $\Ox_U:$
$$
\begin{tikzcd}
\calP_0\otimes \calP_0 \ar{d}{\op{Id}\otimes \delta} \ar{rr}{\op{Id}\otimes (g_2\circ g_1)_{r,s}} & & \calP_0\otimes \calA_{X,r+s} \ar{r}{\epsilon_{r+s}} & \calA_{X,r+s}\otimes \calP_0 \\
\calP_0\otimes \calP_0\otimes \calP_0 \ar{d}{\op{Id}\otimes \op{Id}\otimes g_{1,r}} & & & \calA_{X,r}\otimes \calA_{X,s} \otimes \calP_0 \ar{u} \\
\calP_0\otimes \calP_0\otimes \calA_{X,r} \ar{r}{\op{Id}\otimes \epsilon_r} & \calP_0\otimes \calA_{X,r} \otimes \calP_0 \ar{r}{\op{Id}\otimes \op{Id}\otimes g_{2,s}} & \calP_0\otimes \calA_{X,r}\otimes \calA_{X,s}\ar{r}{\epsilon_r\otimes \op{Id}} \ar{uu} & \calA_{X,r}\otimes \calP_0\otimes \calA_{X,s} \ar{u}{\op{Id}\otimes \epsilon_s}
\end{tikzcd}
$$
The left square of the above diagram is commutative by the definition of $g_2\circ g_1.$ The right square is commutative by the fact that $\epsilon_{\calA_X}$ is a morphism of $\upmu$-indexed algebras. Note that $(g_{3,t}\circ (g_2\circ g_1)_{s,t})$ is equal to the composition of the upper line in the above diagram with $\delta$ to the left and
$$\calA_{X,r+s}\otimes \calP_0 \xrightarrow{\op{Id}\otimes g_{3,t}} \calA_{X,r+s}\otimes \calA_{X,t} \ra \calA_{X,r+s+t}$$
to the right.
It follows that the equality $(g_3\circ (g_2\circ g_1))_{r,s,t}=((g_3\circ g_2)\circ g_1)_{r,s,t}$ is equivalent to the commutativity of the following diagram:
$$
\begin{tikzcd}
\calP_0 \ar{d}{\delta} \ar{r}{\delta} & \calP_0\otimes \calP_0 \ar{r}{\op{Id}\otimes \delta} & \calP_0\otimes \calP_0\otimes \calP_0 \ar{d}{\op{Id}\otimes \op{Id}\otimes g_{1,r}} \\
\calP_0\otimes \calP_0  \ar{dd}{\op{Id}\otimes g_{1,r}} & & \calP_0\otimes \calP_0\otimes \calA_{X,r} \ar{d}{\op{Id} \otimes \epsilon_{r}} \\
 & & \calP_0\otimes \calA_{X,r} \otimes \calP_0 \ar{d}{\op{Id}\otimes \op{Id} \otimes g_{2,s}} \ar[swap]{dl}{\epsilon_{r}\otimes\op{Id}} \\
\calP_0\otimes \calA_{X,r} \ar{d}{\epsilon_{r}} & \calA_{X,r} \otimes \calP_0\otimes \calP_0 \ar[swap]{dr}{\op{Id}\otimes \op{Id} \otimes g_{2,s}} & \calP_0\otimes \calA_{X,r} \otimes \calA_{X,s} \ar{d}{\epsilon_{r}\otimes \op{Id}} \\
\calA_{X,r}\otimes \calP_0 \ar{d}{\op{Id}\otimes \delta} & & \calA_{X,r} \otimes \calP_0 \otimes \calA_{X,s} \ar{d}{\op{Id}\otimes \epsilon_{s}} \\
\calA_{X,r}\otimes \calP_0 \otimes \calP_0 \ar{d}{\op{Id}\otimes \op{Id} \otimes g_{2,s}} \ar[equal]{uur} & & \calA_{X,r} \otimes \calA_{X,s} \otimes \calP_0 \ar{d}{\op{Id} \otimes \op{Id} \otimes g_{3,t}} \\
\calA_{X,r} \otimes \calP_0 \otimes \calA_{X,s} \ar{r}{\op{Id}\otimes \epsilon_{s}} \ar[equal]{uurr} & \calA_{X,r} \otimes \calA_{X,s} \otimes \calP_0\ar{r}{\op{Id} \otimes \op{Id} \otimes g_{3,t}} &  \calA_{X,r} \otimes \calA_{X,s} \otimes \calA_{X,t}
\end{tikzcd}
$$
It is sufficient to prove that the upper rectangle in the following diagram is commutative:
$$
\begin{tikzcd}
\calP_0 \ar{d}{\delta} \ar{r}{\delta} & \calP_0\otimes \calP_0\ar{r}{\op{Id}\otimes \delta} & \calP_0\otimes \calP_0\otimes \calP_0 \ar{d}{\op{Id}\otimes \op{Id} \otimes g_{1,r}} \\
\calP_0\otimes \calP_0 \ar{dd}{\op{Id}\otimes g_{1,r}} & & \calP_0\otimes \calP_0\otimes \calA_{X,r} \ar{d}{\op{Id} \otimes \epsilon_{r}} \\
& & \calP_0\otimes \calA_{X,r} \otimes \calP_0 \ar{d}{\epsilon_{r}\otimes\op{Id}} \\
\calP_0\otimes \calA_{X,r} \ar{r}{\epsilon_{r}} \ar{d}{\delta\otimes \op{Id}}& \calA_{X,r}\otimes \calP_0\ar{r}{\op{Id}\otimes \delta} & \calA_{X,r} \otimes \calP_0\otimes \calP_0 \\
\calP_0 \otimes \calP_0\otimes \calA_{X,r} \ar{rr}{\op{Id}\otimes \epsilon_r} & & \calP_0\otimes \calA_{X,r} \otimes \calP_0 \ar{u}{\epsilon_r\otimes \op{Id}}
\end{tikzcd}
$$
The lower rectangle is commutative by the cocycle condition satisfied by the HPD stratification $\epsilon_r.$ It is thus sufficient to prove the commutativity of the following diagram:
\begin{equation}\label{diagAkaza1}
\begin{tikzcd}
\calP_0 \ar{d}{\delta} \ar{r}{\delta} & \calP_0\otimes \calP_0\ar{d}{\op{Id}\otimes \delta} \\
\calP_0 \otimes \calP_0 \ar{d}{\op{Id}\otimes g_{1,r}} & \calP_0\otimes \calP_0 \otimes \calP_0 \ar{d}{\op{Id}\otimes \op{Id} \otimes g_{1,r}} \\
\calP_0 \otimes \calA_{X,r} \ar{r}{\delta \otimes \op{Id}} & \calP_0\otimes \calP_0\otimes \calA_{X,r}
\end{tikzcd}
\end{equation}
Since
$$
\begin{tikzcd}
\calP_0\otimes \calP_0 \ar{r}{\delta \otimes \op{Id}} \ar{d}{\op{Id}\otimes g_{1,r}} & \calP_0\otimes \calP_0\otimes \calP_0\ar{d}{\op{Id}\otimes \op{Id} \otimes g_{1,r}} \\
\calP_0\otimes \calA_{X,r} \ar{r}{\delta \otimes \op{Id}} & \calP_0\otimes \calP_0\otimes \calA_{X,r}
\end{tikzcd}
$$
is commutative, the diagram \eqref{diagAkaza1} is equal to
$$
\begin{tikzcd}
\calP_0 \ar{rr}{\delta} \ar{d}{\delta} & & \calP_0\otimes \calP_0\ar{d}{\op{Id}\otimes \delta} \\
\calP_0 \otimes \calP_0 \ar{d}{\delta\otimes \op{Id}} & & \calP_0 \otimes \calP_0 \otimes \calP_0\ar{d}{\op{Id}\otimes \op{Id} \otimes g_{1,r}} \\
\calP_0 \otimes \calP_0 \otimes \calP_0 \ar{rr}{\op{Id}\otimes \op{Id} \otimes g_{1,r}} & & \calP_0 \otimes \calP_0 \otimes \calA_{X,r}
\end{tikzcd}
$$
It is thus sufficient to prove the commutativity of
$$
\begin{tikzcd}
\calP_0 \ar{r}{\delta}\ar{d}{\delta} & \calP_0\otimes \calP_0 \ar{d}{\op{Id} \otimes \delta} \\
\calP_0\otimes \calP_0 \ar{r}{\delta \otimes \op{Id}} & \calP_0\otimes \calP_0\otimes \calP_0
\end{tikzcd}
$$
This follows from the commutativity of the diagram
$$
\begin{tikzcd}
X\times_S X\times_S X\times_S X \ar{rr}{p_{13}\times_S \op{Id}_X} \ar{d}{\op{Id}_X \times_S p_{13}} & & X\times_SX\times_SX \ar{d}{p_{13}} \\
X\times_SX\times_SX \ar{rr}{p_{13}} & & X\times_SX
\end{tikzcd}
$$
\end{proof}

\begin{proposition}\label{algstrD}
The indexed algebra structure on $\widetilde{\calD}_{X/S}$ satisfies the following properties:
\begin{enumerate}
\item $(a\otimes 1)(b\otimes 1)=(ab)\otimes 1$ for any local sections $a$ and $b$ of $\calA_X.$
\item $(1\otimes \partial_1)(1\otimes \partial_2)=1\otimes (\partial_1\circ \partial_2)$ in $\calA_{X,0}\otimes \calD_{X/S}$ for any PD differential operators $\partial_1$ and $\partial_2.$
\item $(a\otimes 1)(1\otimes \partial)=a\otimes \partial$ for any PD differential operator $\partial$ and any local section $a$ of $\calA_X.$
\item If we have $d$ local coordinates of $X$ over $S,$ then, with the notation of \ref{P2}, for any multi-index $I\in \mathbb{N}^d$ and any local section $a$ of $\calA_X,$
$$(1\otimes \partial_I)(a\otimes 1)=\sum_{J\le I}\begin{pmatrix} I \\ J\end{pmatrix} (\partial_J\cdot a)\otimes \partial_{I-J},$$
where on the right hand side, we consider the action of $\calD_{X/S,\upmu}$ on $\calA_X$ via the connection $d_{\calA_X}$ (\ref{prop39}) 
\end{enumerate}
\end{proposition}

\begin{proof}
See (\cite{Mon} 4.1.2).
\end{proof}

\begin{proposition}\label{Dcenter}~ 
\begin{enumerate}
\item The center of $\widetilde{\calD}_{X/S}$ is $\calB_{X/S}\otimes_{\Ox_{X',\upmu}}\calC_{X/S,\upmu},$ where $\calC_{X/S}$ is the image of \eqref{Takriz2}.
\item The centralizer of $\calA_X$ in $\widetilde{\calD}_{X/S}$ i.e. the subsheaf of $\widetilde{\calD}_{X/S}$ consisting of local sections $x\in \widetilde{\calD}_{X/S}$ such that $x(a\otimes 1)=(a\otimes 1)x$ for every local section $a$ of $\calA_X,$ is equal to $\calA_X\otimes_{\Ox_{X',\upmu}}\calC_{X/S,\upmu}.$
\item The indexed subalgebra $\calA_X\otimes_{\Ox_{X',\upmu}}\calC_{X/S,\upmu}$ is commutative and the canonical morphism
$$\calA_X\otimes_{\Ox_{X',\upmu}}\calC_{X/S,\upmu} \ra \widetilde{\calD}_{X/S}$$
is a morphism of $\upmu$-indexed algebras.
\end{enumerate}
\end{proposition}

\begin{proof}
We can work étale locally and so we can suppose that we have local coordinates $m_1,\hdots,m_d\in \Gamma(X,\calM_X)$ of $X$ over $S$ (\ref{P2}). Consider the corresponding differential operators $\partial_I$ for multi-indices $I\in \mathbb{N}^d$ (\ref{partialI}). The sections $(\op{dlog}\pi^{\flat}m_i)_{1\le i\le d}$ form a basis for the $\Ox_{X'}$-module $\omega^1_{X'/S}.$ It follows from (\ref{era2Doma11}) that the $\Ox_X$-algebra $\calC_{X/S}$ is generated by $\partial_{\epsilon_i}^p-\partial_{\epsilon_i}=\partial_{p\epsilon_i},\ 1\le i\le d.$ Let $s\in \Gamma(U,\upmu)$ be a section of $\upmu$ over an étale $X$-scheme $U$ and $w$ a section of $\widetilde{\calD}_{X/S}$ over $s.$ First, let us suppose that $w$ is of the form
$$w=b\otimes \partial_{p\epsilon_i},$$
where $b$ is a section of $\calB_{X/S}$ over $s.$
Then, by \ref{algstrD}, for any local section $a$ of $\calA_X$ and any multi-index $I\in \N^d,$
\begin{alignat*}{2}
(b\otimes \partial_{p\epsilon_i})(a\otimes \partial_I) &= \sum_{k=0}^p \begin{pmatrix}p\\ k \end{pmatrix}(b \partial_{k\epsilon_i}\cdot a)\otimes \partial_{(p-k)\epsilon_i}\circ \partial_I \\
&= (ab)\otimes (\partial_{p\epsilon_i}\circ \partial_I) + (b \partial_{p\epsilon_i}\cdot a)\otimes \partial_I\\
&= (ab)\otimes (\partial_{p\epsilon_i}\circ \partial_I),
\end{alignat*}
where the last equality follows from the fact that the connection $d_{\calA_X}$ has vanishing $p$-curvature. Indeed, if we denote by $(d_{A_X,n})_{n\ge 0}$ the morphisms defined by the connection $d_{A_X},$ as in \ref{prop39} (3), then $\partial_{p\epsilon_i}\cdot a=d_{\calA_X,p}(\partial_{p\epsilon_i})(a)=d_{\calA_X,p}(\partial_{\epsilon_i}^p-\partial_{\epsilon_i}^{(p)})(a)=0.$ On the other hand, since $b$ is a section of $\calB_{X/S},$ $\partial_J\cdot b=0$ for any multi-index $J\in \N^d$ such that $|J|\ge 1.$ Then, by \ref{algstrD} (4),
\begin{alignat*}{2}
(a\otimes \partial_I)(b\otimes \partial_{p\epsilon_i}) &= \sum_{J\le I} \begin{pmatrix} I \\ J \end{pmatrix} (a \partial_J\cdot b)\otimes (\partial_{I-J}\circ \partial_{p\epsilon_i}) \\
&= (ab) \otimes (\partial_I\circ \partial_{p\epsilon_i}).
\end{alignat*}
By \ref{cor18},
$$\partial_I\circ \partial_{p\epsilon_i}=\partial_{p\epsilon_i}\circ \partial_I.$$
So $w$ is in the center of $\widetilde{\calD}_{X/S}.$

Conversely, suppose that $w=\sum_{I\in \mathbb{N}^d}a_I\otimes \partial_I$ is in the center of $\widetilde{\calD}_{X/S}.$ Denote by $e_{m_i}$ the $\Ox_X$-basis of $\calA_{X,m_i}$ defined in \ref{paragAX} and let $j\in \{1,\hdots, d\}.$
\begin{alignat*}{2}
(1\otimes \partial_{\epsilon_j})w &= \sum_{I \in \mathbb{N}^d}a_I\otimes \partial_{\epsilon_j}\circ \partial_I+(\partial_{\epsilon_j}\cdot a_I)\otimes \partial_I \\
w (1\otimes \partial_{\epsilon_j}) &= \sum_{I \in\mathbb{N}^d} a_I \otimes \partial_I \circ \partial_{\epsilon_j}.
\end{alignat*}
By \ref{cor18}, $\partial_I\circ \partial_{\epsilon_j}=\partial_{\epsilon_j}\circ \partial_I,$ so
$$\sum_{I\in \mathbb{N}^d}\partial_{\epsilon_j}\cdot a_I\otimes \partial_I=0.$$
It follows that $\partial_{\epsilon_j}\cdot a_I=0$ for all $I\in\mathbb{N}^d$ and $1\le j\le d$ so $a_I$ is a section of $\calB_{X/S}$ for any multi-index $I\in \mathbb{N}^d.$ On the other hand, $\partial_{\epsilon_i}\cdot e_{m_i}=e_{m_i}$ and $\partial_{2\epsilon_i}\cdot e_{m_i}=((\partial_{\epsilon_i}-1)\circ \partial_{\epsilon_i}) \cdot e_{m_i}=0.$ So
\begin{alignat*}{2}
(e_{m_i}\otimes 1)w &= \sum_{I\in\mathbb{N}^d}(a_Ie_{m_i})\otimes \partial_I.\\
w(e_{m_i}\otimes 1) &= \sum_{I\in \mathbb{N}^d}\sum_{J\le I} \begin{pmatrix} I\\J \end{pmatrix} (a_I \partial_J\cdot e_{m_i})\otimes \partial_{I-J} \\
&= \sum_{I\in\mathbb{N}^d} a_Ie_{m_i}\otimes \partial_I + \sum_{\substack{I\in\mathbb{N}^d \\ I_i\ge 1}} I_ia_Ie_{m_i}\otimes \partial_{I-\epsilon_i}.
\end{alignat*}
It follows that
$a_I=0$ if $I \notin p\mathbb{N}^d.$ By \ref{cor18},
$$\partial_{pI}=\prod_{i=1}^d\partial_{pI_i\epsilon_i}$$
and
\begin{alignat*}{2}
\partial_{pI_i\epsilon_i} &= \prod_{k=0}^{pI_i-1}(\partial_{\epsilon_i}-k) \\
&= \prod_{l=0}^{I_i-1} \prod_{k=pl}^{p(l+1)-1}(\partial_{\epsilon_i}-k) \\
&= \prod_{l=0}^{I_i-1} \prod_{k=0}^{p-1}(\partial_{\epsilon_i}-k) \\
&= \left (\partial_{p\epsilon_i} \right )^{I_i}  \in \calC_{X/S}.
\end{alignat*}
This concludes the proof of the first assertion. The second assertion is proved by similar computations and the third is an immediate consequence of the first and the second.
\end{proof}

\begin{parag}
We consider $\widetilde{\calD}_{X/S}$ as a $\upmu$-indexed $\widetilde{\calD}_{X/S}$-module by multiplication to the right. Let $p_1,p_2:\upmu^2 \ra \upmu$ be the canonical projections and $\sigma:\upmu^2\ra \upmu$ the addition map.
Note that
$$\calA_X\otimes_{\Ox_{X,\upmu}}(F_1^*\calC_{X/S})_{\upmu}=\calA_X\otimes_{\Ox_{X',\upmu}}\calC_{X/S,\upmu}.$$
Let $U$ be an étale $X$-scheme, $s,t\in \Gamma(U,\upmu)$ and
$$\sigma:\upmu_{|U} \ra \upmu_{|U},\ x\mapsto x+s.$$
By \ref{indremtakriz1} (2), a local section $\partial \in \widetilde{\calD}_{X/S}(s)$ defines a local section $\partial_g$ of $\mathscr{End}_{\calA_X\otimes_{\Ox_{X',\upmu}}\calC_{X/S,\upmu}}(\widetilde{\calD}_{X/S})(s)$ (\ref{inddef2}) given by
$$
\partial_g:\begin{array}[t]{clc}
\left ( \widetilde{\calD}_{X/S}\right )_{|U} & \ra & \sigma_s^*\left (\widetilde{\calD}_{X/S} \right )_{|U} \\
x & \mapsto & \partial x,
\end{array}
$$
where $\partial x$ is the product of $\partial$ and $x$ in the $\upmu$-indexed algebra $\widetilde{\calD}_{X/S}$ (\ref{era2stralgD}).
Similarly, a local section $w\in \left (\calA_X\otimes_{\Ox_{X',\upmu}}\calC_{X/S,\upmu}\right )(t)$ defines a local section of $w_d\in \mathscr{End}_{\calA_X\otimes_{\Ox_{X',\upmu}}\calC_{X/S,\upmu}}(\widetilde{\calD}_{X/S})(t)$ given by
$$
w_d:\begin{array}[t]{clc}
\left (\widetilde{\calD}_{X/S}\right )_{|U} & \ra & \sigma_t^*\left (\widetilde{\calD}_{X/S} \right )_{|U} \\
x & \mapsto & xw.
\end{array}
$$
The $\left (\calA_X\otimes_{\Ox_{X',\upmu}}\calC_{X/S,\upmu} \right )$-linearity of $w_d$ follows from \ref{Dcenter} (3) and the fact that the $\upmu$-indexed $\left (\calA_X\otimes_{\Ox_{X',\upmu}}\calC_{X/S,\upmu} \right )$-module structure on $\widetilde{\calD}_{X/S}$ is given by multiplication to the right.

Hence we obtain two morphisms
$$
\varphi_1:\begin{array}[t]{clc}
\widetilde{\calD}_{X/S} & \rightarrow & \mathscr{End}_{\calA_X\otimes_{\Ox_{X',\upmu}}\calC_{X/S,\upmu}}(\widetilde{\calD}_{X/S}) \\
\partial & \mapsto & \partial_g
\end{array}
$$
$$
\varphi_2:\begin{array}[t]{clc}
\calA_X\otimes_{\Ox_{X',\upmu}}\calC_{X/S,\upmu} & \rightarrow & \mathscr{End}_{\calA_X\otimes_{\Ox_{X',\upmu}}\calC_{X/S,\upmu}}(\widetilde{\calD}_{X/S}) \\
w & \mapsto & w_d
\end{array}
$$
If we equip $\mathscr{End}_{\calA_X\otimes_{\Ox_{X',\upmu}}\calC_{X/S,\upmu}}(\widetilde{\calD}_{X/S})$ with the $\upmu$-indexed $\widetilde{\calD}_{X/S}$-module structure given by
$$
\begin{array}[t]{clc}
\widetilde{\calD}_{X/S} \boxtimes \mathscr{End}_{\calA_X\otimes_{\Ox_{X',\upmu}}\calC_{X/S,\upmu}}(\widetilde{\calD}_{X/S}) & \ra & \sigma^*\mathscr{End}_{\calA_X\otimes_{\Ox_{X',\upmu}}\calC_{X/S,\upmu}}(\widetilde{\calD}_{X/S}) \\
\partial \otimes f & \mapsto & (\partial\cdot f:x\mapsto f(\partial x)),
\end{array}
$$
then
$$\varphi_1(\partial_1\partial_2)=\partial_2\cdot \varphi_1(\partial_1),$$
for any local sections $\partial_1$ and $\partial_2$ of $\widetilde{\calD}_{X/S}.$
So the morphism $\varphi_1$ is $\widetilde{\calD}_{X/S}$-linear.

Similarly, if we equip $\mathscr{End}_{\calA_X\otimes_{\Ox_{X',\upmu}}\calC_{X/S,\upmu}}(\widetilde{\calD}_{X/S})$ with the $\upmu$-indexed $\left (\calA_X\otimes_{\Ox_{X',\upmu}}\calC_{X/S,\upmu} \right )$-module structure given by
$$
\begin{array}[t]{clc}
\left (\calA_X\otimes_{\Ox_{X',\upmu}}\calC_{X/S,\upmu} \right ) \boxtimes \mathscr{End}_{\calA_X\otimes_{\Ox_{X',\upmu}}\calC_{X/S,\upmu}}(\widetilde{\calD}_{X/S}) & \ra & \sigma^*\mathscr{End}_{\calA_X\otimes_{\Ox_{X',\upmu}}\calC_{X/S,\upmu}}(\widetilde{\calD}_{X/S}) \\
w \otimes f & \mapsto & (w\cdot f:x\mapsto f(xw)=f(x)w),
\end{array}
$$
then,
$$
\varphi_2(w_1w_2)=w_2\cdot \varphi_2(w_1),
$$
for any local sections $w_1$ and $w_2$ of $\calA_X\otimes_{\Ox_{X',\upmu}}\calC_{X/S,\upmu}$ and so $\varphi_2$ is $\left (\calA_X\otimes_{\Ox_{X',\upmu}}\calC_{X/S,\upmu}\right )$-linear.

We define a morphism a $\left (\calB_{X/S}\otimes_{\Ox_{X',\upmu}}\calC_{X/S,\upmu}\right )$-linear morphism
\begin{equation}\label{phitakriz}
\varphi:\widetilde{\calD}_{X/S} \circledast_{\calB_{X/S}\otimes_{\Ox_{X',\upmu}}\calC_{X/S,\upmu}}(\calA_X\otimes_{\Ox_{X',\upmu}}\calC_{X/S,\upmu}) \rightarrow \mathscr{End}_{\calA_X\otimes_{\Ox_{X',\upmu}}\calC_{X/S,\upmu}}(\widetilde{\calD}_{X/S})
\end{equation}
as the composition of
$$\varphi_1\circledast_{\calB_{X/S}\otimes_{\Ox_{X',\upmu}}\calC_{X/S,\upmu}} \varphi_2$$
and the morphism
\begin{alignat*}{2}
\mathscr{End}_{\calA_X\otimes_{\Ox_{X',\upmu}}\calC_{X/S,\upmu}}(\widetilde{\calD}_{X/S})\circledast_{\calB_{X/S}\otimes_{\Ox_{X',\upmu}}\calC_{X/S,\upmu}} \mathscr{End}_{\calA_X\otimes_{\Ox_{X',\upmu}}\calC_{X/S,\upmu}}(\widetilde{\calD}_{X/S}) \\
\ra \mathscr{End}_{\calA_X\otimes_{\Ox_{X',\upmu}}\calC_{X/S,\upmu}}(\widetilde{\calD}_{X/S})
\end{alignat*}
defined in \eqref{era3comptakriz}.
When considering the appropriate module structures given above, $\varphi$ is both $\widetilde{\calD}_{X/S}$-linear and $\left (\calA_X\otimes_{\Ox_{X',\upmu}}\calC_{X/S,\upmu}\right )$-linear.
\end{parag}

\begin{proposition}\label{AzumayaAkaza}
The $\left (\calB_{X/S}\otimes_{\Ox_{X',\upmu}}\calC_{X/S,\upmu}\right )$-linear morphism \eqref{phitakriz}
$$\varphi:\widetilde{\calD}_{X/S} \circledast_{\calB_{X/S}\otimes_{\Ox_{X',\upmu}}\calC_{X/S,\upmu}}(\calA_X\otimes_{\Ox_{X',\upmu}}\calC_{X/S,\upmu}) \rightarrow \mathscr{End}_{\calA_X\otimes_{\Ox_{X',\upmu}}\calC_{X/S,\upmu}}(\widetilde{\calD}_{X/S})$$
is an isomorphism.
\end{proposition}

\begin{proof}
We can work étale locally on $X,$ so we may suppose that we have local coordinates $m_1,\hdots,m_d \in \Gamma(X,\calM_X).$
For any $I\in \Z^d,$ we denote by $m_I=\sum_{i=1}^dI_im_i \in \Gamma(X,\calM_X^{gp})$ and by $e_{m_I}$ the corresponding basis for $\calA_{X,\ov{m}_I}$ (\ref{paragAX}). Let $M=\{0,\hdots,p-1\}^d.$ By \ref{thmlor}, $(e_{m_I})_{I\in M}$ is a basis for the $\upmu$-indexed $\calB_{X/S}$-module $\calA_X.$ Then, by \ref{propfreetakriz}, $(\theta_I)_{I\in M}=(e_{m_I}\otimes 1)_{I\in M}$ is a basis for the $\upmu$-indexed $(\calB_{X/S}\otimes_{\Ox_{X',\upmu}}\calC_{X/S,\upmu})$-module $\calA_X\otimes_{\Ox_{X',\upmu}}\calC_{X/S,\upmu}=\calA_X\otimes_{\Ox_{X',\upmu}}\calC_{X/S,\upmu}.$ It follows, again by \ref{propfreetakriz}, that
$$(1\otimes \theta_I)_{I\in M}$$
is a basis for the left $\widetilde{\calD}_{X/S}$-module $\widetilde{\calD}_{X/S} \circledast_{\calB_{X/S}\otimes_{\Ox_{X',\upmu}}\calC_{X/S,\upmu}}(\calA_X\otimes_{\Ox_{X',\upmu}}\calC_{X/S,\upmu}).$ Note that $\theta_I$ is invertible in the $\upmu$-indexed algebra $\calA_X\otimes_{\Ox_{X',\upmu}}\calC_{X/S,\upmu}$ with inverse
$$\theta_I^{-1}=e_{-m_I}\otimes 1.$$
On the other hand, $(\partial_I)_{I\in M}$ is a basis for the $F_{1}^*\calC_{X/S}$-module $\calD_{X/S}$ (as shown in \ref{era33lem}) and so $(1\otimes \partial_I)_{I\in M}$ is a basis for the $\upmu$-indexed $(\calA_X\otimes_{\Ox_{X',\upmu}}\calC_{X/S,\upmu})$-module $\widetilde{\calD}_{X/S}.$ For any $I\in M,$ denote by $\alpha_I$ the global section of
$$\mathscr{End}_{\calA_X\otimes_{\Ox_{X',\upmu}}\calC_{X/S,\upmu}}(\widetilde{\calD}_{X/S})$$
defined by
$$\alpha_I(1\otimes \partial_J)=\delta_{IJ}\ \forall I,J\in M.$$
The family $(\alpha_I)_{I\in M}$ is then a basis for the $\upmu$-indexed right $\widetilde{\calD}_{X/S}$-module $\mathscr{End}_{\calA_X\otimes_{\Ox_{X',\upmu}}\calC_{X/S,\upmu}}(\widetilde{\calD}_{X/S})$ (\ref{algstrendleft}). By definition of $\varphi$ and by \ref{algstrD} (4), we have the following equalities in $\widetilde{\calD}_{X/S}:$
\begin{alignat*}{2}
\varphi(1\otimes \theta_I)(1\otimes \partial_J) &= (1\otimes \partial_J)\theta_I \\
&= (1\otimes \partial_J)(e_{m_I}\otimes 1) \\
&= \sum_{K\le J} \begin{pmatrix}J\\K \end{pmatrix}(\partial_K\cdot e_{m_I})\otimes \partial_{J-K}.
\end{alignat*}
For $1\le i\le d$ and $I=\epsilon_i,$ we have $\partial_{\epsilon_i}\cdot e_{m_i}=e_{m_i}$ and $\partial_{2\epsilon_i}\cdot e_{m_i}=((\partial_{\epsilon_i}-1)\circ \partial_{\epsilon_i}) \cdot e_{m_i}=0$ (\ref{cor18}). So $\partial_K\cdot e_{m_i}=0$ if $K\neq 0$ or $K\neq \epsilon_i.$ We get
\begin{alignat*}{2}
\varphi(1\otimes \theta_{\epsilon_i})(1\otimes \partial_J) &= (1\otimes \partial_J)\theta_{\epsilon_i} \\
&= (1\otimes \partial_J)(e_{m_i}\otimes 1) \\
&= \sum_{K\le J} \begin{pmatrix}J\\K \end{pmatrix}(\partial_K\cdot e_{m_i})\otimes \partial_{J-K}\\
&= e_{m_i}\otimes (\partial_J+J_i\partial_{J-\epsilon_i}).
\end{alignat*}
If $J_i\ge 1,$
$$\varphi(\theta_{\epsilon_i}^{-1}\otimes \theta_{\epsilon_i}-1\otimes 1)(1\otimes \partial_J)=J_i (1\otimes \partial_{J-\epsilon_i})$$
and, if $J_i=0,$
$$\varphi(\theta_{\epsilon_i}^{-1}\otimes \theta_{\epsilon_i}-1\otimes 1)(1\otimes \partial_J)=0.$$
For any $I=(I_1,\hdots,I_d)\in \N^d,$ consider the local section of $\widetilde{\calD}_{X/S} \circledast_{\calB_{X/S}\otimes_{\Ox_{X',\upmu}}\calC_{X/S,\upmu}}(\calA_X\otimes_{\Ox_{X',\upmu}}\calC_{X/S,\upmu})$
$$\xi_I=\prod_{i=1}^d (\theta_{\epsilon_i}^{-1}\otimes \theta_{\epsilon_i}-1\otimes 1)^{I_i}.$$
A simple induction yields $\varphi(\xi_I)(1\otimes \partial_J)= \frac{J!}{(J-I)!}(1\otimes \partial_{J-I})$ if $J\ge I$ and $\varphi(\xi_I)(1\otimes \partial_J)=0$ otherwise. This, along with the definition of $\xi_I,$ imply the two following formulas:
$$\xi_I=\sum_{J\le I}(-1)^{|I-J|}\begin{pmatrix}I\\ J \end{pmatrix}\theta_J^{-1}\otimes \theta_J$$
$$\varphi(\xi_I)=\sum_{\substack{J\ge I \\ J\in M}}\frac{J!}{(J-I)!}(1\otimes \partial_{J-I})\alpha_J.$$
With a suitable ordering of $M\subset \N^d,$ we see that the matrix of $(\xi_I)_{I\in M}$ in the basis $(1\otimes \theta_I)_{I\in M}$ and the matrix of $(\varphi(\xi_I))_{I\in M}$ in the basis $(\alpha_I)_{I\in M}$ are both triangular with diagonal coefficients equal to $1.$
We deduce that $(\xi_I)_{I\in M}$ is a basis for the left $\widetilde{\calD}_{X/S}$-module $\widetilde{\calD}_{X/S} \circledast_{\calB_{X/S}\otimes_{\Ox_{X',\upmu}}\calC_{X/S,\upmu}}(\calA_X\otimes_{\Ox_{X',\upmu}}\calC_{X/S,\upmu})$ and that $\varphi$ sends it to a basis, it is thus an isomorphism.
\end{proof}

\begin{corollaire}\label{krazcor158}
$\widetilde{\calD}_{X/S}$ is a $\upmu$-indexed Azumaya algebra over its center $\calB_{X/S}\otimes_{\Ox_{X',\upmu}}\calC_{X/S,\upmu}.$
\end{corollaire}

\begin{lemma}\label{era33lem}
Suppose that $X\ra S$ has local coordinates $m_1,\hdots,m_d \in \Gamma(X,\calM_X)$ \eqref{P2}. The $F_1^*\calC_{X/S}$-module $\calD_{X/S}$ has a basis $(\partial_I)_{I\in \{0,\hdots,p-1\}^d}.$
\end{lemma}

\begin{proof}
First note that, by \ref{lemiran1}, for any integers $n,r\ge 0,$
\begin{alignat*}{2}
\partial_{(pn+r)\epsilon_i} &= \prod_{k=0}^{pn+r-1}(\partial_{\epsilon_i}-k) &= (\partial_{p\epsilon_i})^n \circ \partial_{r\epsilon_i}.
\end{alignat*}
We have a basis $(\partial_I)_{I\in \N^d}$ for the $\Ox_X$-module $\calD_{X/S}.$ Let $I=(I_1,\hdots,I_d)\in \N^d$ and $I_i=pq_i+r_i,\ 0\le r_i<p.$ Then
$$\partial_I=\prod_{i=1}^d\partial_{I_i\epsilon_i}=\prod_{i=1}^d (\partial_{p\epsilon_i})^{q_i}\circ \partial_{r_i\epsilon_i}.$$
We deduce that the family $(\partial_I)_{I\in \{0,\hdots,p-1\}^d}$ is generating.

Now let $(b_I)_{I\in \{0,\hdots,p-1 \}^d}$ be local sections of $F_1^*\calC_{X/S}$ such that
$$\sum_{i\in \{0,\hdots,p-1 \}^d }b_I\partial_I=0.$$
By (\ref{era2Doma11}), for any $I\in \{0,\hdots,p-1\}^d,$ there exists a finite subset $\Gamma_I\subset \N^d$ and local sections $a_{I,J}$ of $\Ox_X$ for $J\in \Gamma_I$ such that $b_I=\sum_{J\in \Gamma_I}a_{I,J}\partial_{pJ}.$ Hence
$$
\sum_{I\in \{0,\hdots,p-1 \}^d,J\in \Gamma_I}a_{I,J}\partial_{pJ}\circ \partial_I=0.
$$
By \ref{lemiran1},
$$\partial_{pJ}\circ \partial_I=\partial_{pJ+I}.$$
We deduce that
$$\sum_{I\in \{0,\hdots,p-1 \}^d,J\in \Gamma_I}a_{I,J}\partial_{pJ+I}=0.$$
Since $pK+I \neq pL+J$ for $K,L\in \N^d$ and $I,J\in \{0,\hdots,p-1 \}^d$ such that $I\neq J,$ it follows that for any $I\in \{0,\hdots,p-1 \}^d$ and $J\in \Gamma_I,$ $a_{I,J}=0.$
\end{proof}

\subsection*{A splitting of the Azumaya algebra $\widetilde{\calD}_{X/S}$}

\begin{parag}
In the remaining of this section, we consider a perfect field $k$ of positive characteristic $p$ and denote its ring of Witt vectors by $W.$ We equip $\op{Spf}W$ with the trivial logarithmic structure. We suppose that $f_1:X\ra S$ lifts to a log smooth morphism $f:(\frakX,Q) \ra (\frakS,P)$ of framed fs $p$-adic logarithmic formal schemes such that $\frakS$ is log flat and locally of finite type over $\op{Spf}W.$ Recall that $X'$ is naturally equipped with a frame $X' \ra [Q']$ and that the exact relative Frobenius $F_1:X\ra X'$ is a morphism of framed logarithmic schemes \eqref{cor411}. We suppose that $F_1:X \ra X'$ lifts to an $(\frakS,P)$-morphism of framed fs $p$-adic logarithmic formal schemes $F:(\frakX,Q) \ra (\frakX',Q')$ such that $\frakX'$ is log smooth over $\frakS.$ Let $\frakY=\frakX\times_{\frakS,[Q]}^{\op{log}}\frakX,$ $\frakY'=\frakX'\times_{\frakS,[Q']}^{\op{log}}\frakX',$ $P_{\frakX/\frakS,0}$ (resp. $P_{\frakX'/\frakS,0}$) the PD envelope of $\frakX \ra \frakY$ (resp. $\frakX' \ra \frakY'$) as defined in \ref{parag86}. Denote by $\iota:\frakX \ra P_{\frakX/\frakS,0}$ and $\iota':\frakX' \ra P_{\frakX'/\frakS,0}$ the canonical immersions and by $\calP_{\frakX/\frakS,0}$ and $\calP_{\frakX'/\frakS,0}$ the structural rings of $P_{\frakX/\frakS,0}$ and $P_{\frakX'/\frakS,0}$ respectively. We identify the étale sites on $\frakX,$ $\frakX',$ $P_{\frakX/\frakS,0}$ and $P_{\frakX'/\frakS,0}$ via the universal homeomorphisms $F,\ \iota$ and $\iota'.$ We finally denote by $G:\frakY \ra \frakY'$ the morphism induced by $F$ and by $p_1,p_2:\frakY \ra \frakX,$ $p_1',p_2':\frakY' \ra \frakX'$ the canonical projections.
Recall that we have a morphism (\ref{uAkaza}) of PD-$\Ox_{X'}$-algebras:
\begin{equation}\label{era2uAkaza}
u:F_1^*\Gamma^{\bullet}\omega^1_{X'/S} \ra \calP_0.
\end{equation}
We have also constructed an HPD stratification on $F_1^*\Gamma^{\bullet}\omega^1_{X'/S}$ (\ref{propAkaza}) given by the $\calP_0$-linearization of
$$F_1^*\Gamma^{\bullet}\omega^1_{X'/S}\xrightarrow{\Delta} F_1^*\Gamma^{\bullet}\omega^1_{X'/S} \otimes_{\Ox_X}F_1^*\Gamma^{\bullet}\omega^1_{X'/S} \xrightarrow{\op{Id}\otimes u}F_1^*\Gamma^{\bullet}\omega^1_{X'/S} \otimes_{\Ox_X}\calP_0,$$
where $\Delta$ is the comultiplication of $F_1^*\Gamma^{\bullet}\omega^1_{X'/S}.$ This HPD stratification on $F_1^*\Gamma^{\bullet}\omega^1_{X'/S}$ and \ref{dualAkaza} induce an HPD stratification on $F_1^*\w{S}^{\bullet}\calT_{X'/S}.$ We thus get an action of $\calD_{X/S,\upmu}$ on $(F_1^*\widehat{S}^{\bullet}\calT_{X'/S})_{\upmu}.$

Let $S=(S^{\bullet}\calT_{X'/S})_{\upmu},$ $\widehat{S}=(\widehat{S}^{\bullet}\calT_{X'/S})_{\upmu}$ and $F_1^*\widehat{S}=(F_1^*\widehat{S}^{\bullet}\calT_{X'/S})_{\upmu}.$
The actions of $\calD_{X/S,\upmu}$ on $\calA_X,$ given by the connection $d_{\calA_X},$ and on $F_1^*\widehat{S}$ define a canonical action of $\calD_{X/S,\upmu}$ on $\calA_X\otimes_{\Ox_{X,\upmu}}F_1^*\widehat{S}.$
By \ref{era4tak1}, this action extends to an action of $\widehat{S}\otimes_{S}\calD_{X/S,\upmu}$ on $\calA_X\otimes_{\Ox_{X,\upmu}}F_1^*\widehat{S}$ and then to an action of $\widehat{S}\otimes_{S}\widetilde{\calD}_{X/S}$ on $\calA_X\otimes_{\Ox_{X,\upmu}}F_1^*\widehat{S}.$
By \ref{Dcenter} (1), $\calB_{X/S}\otimes_{\Ox_{X',\upmu}}\widehat{S}$ is a sub-indexed-algebra of the center of $\widehat{S}\otimes_{S}\widetilde{\calD}_{X/S}.$ We thus obtain a morphism of $\upmu$-indexed algebras
\begin{equation}\label{era2chi}
\chi: \widehat{S}\otimes_{S}\widetilde{\calD}_{X/S}  \ra \mathscr{End}_{\calB_{X/S}\otimes_{\Ox_{X',\upmu}}\widehat{S}}\left ( \calA_X\otimes_{\Ox_{X,\upmu}}F_1^*\widehat{S}\right ).
\end{equation}
\end{parag}

\begin{proposition}\label{erafPsi}
The composition
\begin{alignat*}{2}
\psi:(\calB_{X/S}\otimes_{\Ox_{X',\upmu}}\widehat{S})\circledast_{\calB_{X/S}\otimes_{\Ox_{X',\upmu}}S}\widetilde{\calD}_{X/S} \xrightarrow{\sim}&
\widehat{S}\otimes_{S}\widetilde{\calD}_{X/S} \\
 \xrightarrow{\chi} & \mathscr{End}_{\calB_{X/S}\otimes_{\Ox_{X',\upmu}}\widehat{S}}\left ( \calA_X\otimes_{\Ox_{X,\upmu}}F_1^*\widehat{S}\right ),
\end{alignat*}
where the first arrow is the canonical isomorphism \eqref{propisoAkaza} and $\chi$ is given in \eqref{era2chi}, is an isomorphism of $\upmu$-indexed algebras.
\end{proposition}

\begin{proof}
Let $d$ be the rank of $\omega^1_{X/S}.$
By \ref{AzumayaAkaza}, $\widetilde{\calD}_{X/S}$ is an Azumaya algebra of rank $p^{2d}$ over $\calB_{X/S}\otimes_{\Ox_{X',\upmu}}S$ so $(\calB_{X/S}\otimes_{\Ox_{X',\upmu}}\widehat{S})\circledast_{\calB_{X/S}\otimes_{\Ox_{X',\upmu}}S}\widetilde{\calD}_{X/S}$ is an Azumaya algebra of rank $p^{2d}$ over $\calB_{X/S}\otimes_{\Ox_{X',\upmu}}\widehat{S}.$ On the other hand, by \ref{thmlor}, $\calA_X\otimes_{\Ox_{X,\upmu}}F_1^*\widehat{S}= \calA_X\otimes_{\Ox_{X',\upmu}}\widehat{S}$ is a locally free $(\calB_{X/S}\otimes_{\Ox_{X',\upmu}}\widehat{S})$-module of rank $p^d.$ So, by (\cite{Schepler} 2.5), $\psi$ is an isomorphism of $\upmu$-indexed algebras.
\end{proof}

\begin{corollaire}
The functor $\Psi^{az}$ from the category of $\upmu$-indexed $(\calB_{X/S}\otimes_{\Ox_{X',\upmu}}\widehat{S})$-modules to the category of $\upmu$-indexed left $(\widehat{S}\otimes_{S}\widetilde{\calD}_{X/S})$-modules, given by
$$\Psi^{az}(\calE')=(\calA_X\otimes_{\Ox_{X,\upmu}}F_1^*\widehat{S}) \circledast_{\calB_{X/S}\otimes_{\Ox_{X',\upmu}}\widehat{S}}  \calE'=\calA_X\circledast_{\calB_{X/S}}\calE'$$
is an equivalence of categories.
\end{corollaire}

\begin{corollaire}\label{cor1713}
Let $\widehat{\Gamma}^{\bullet}\calT_{X'/S}$ be the completed PD-algebra of $\calT_{X'/S}$ and $\widehat{\Gamma}=\left (\widehat{\Gamma}^{\bullet}\calT_{X'/S} \right )_{\upmu}.$ There exists a canonical isomorphism of $\upmu$-indexed algebras
\begin{equation}\label{iso17131}
(\calB_{X/S}\otimes_{\Ox_{X',\upmu}}\widehat{\Gamma})\circledast_{\calB_{X/S}\otimes_{\Ox_{X',\upmu}}S}\widetilde{\calD}_{X/S} \xrightarrow{\sim} \mathscr{End}_{\calB_{X/S}\otimes_{\Ox_{X',\upmu}}\widehat{\Gamma}}\left ( \calA_X\otimes_{\Ox_{X,\upmu}}F_1^*\widehat{\Gamma}\right ).
\end{equation}
\end{corollaire}

\begin{proof}
It is sufficient to base change the isomorphism $\psi$ of \ref{erafPsi} by the canonical morphism
$$\calB_{X/S}\otimes_{\Ox_{X',\upmu}}\widehat{S} \ra \calB_{X/S}\otimes_{\Ox_{X',\upmu}}\widehat{\Gamma}$$
and then apply \eqref{erafiso} and \ref{propisoAkaza}.
\end{proof}

\begin{definition}\label{deflocnilind}
Let $\widehat{\Gamma}$ be as defined in \ref{cor1713}. We say that a $\widehat{\Gamma}$-module $\calE$ is \emph{locally PD-nilpotent} if, for every local section $s\in \Gamma(U,\upmu)$ over an étale $X$-scheme $U,$ the $\widehat{\Gamma}^{\bullet}\calT_{U'/S}$-module $\calE_s$ is locally PD-nilpotent (\ref{locPDnil}). If $\calA$ is a $\widehat{\Gamma}$-algebra then a $\calA$-module is said to be \emph{locally PD-nilpotent} if it is so as a $\widehat{\Gamma}$-module.
\end{definition}

\begin{proposition}\label{locPDnilprop}
Let $\widehat{\Gamma}$ be as defined in \ref{cor1713}. Let $\calE$ be a $\widehat{\Gamma}$-module. Suppose that there exists an étale covering $(U_i \ra X)_{i\in I}$ such that the restriction $\calE_{|U_i}$ \eqref{defresind} is locally PD-nilpotent for all $i\in I.$ Then $\calE$ is locally PD-nilpotent.
\end{proposition}

\begin{proof}
Let $s\in \Gamma(U,\upmu)$ be a section over an étale $X$-scheme $U$ and $s_i=s_{|U\times_XU_i}.$ The diagram
$$
\begin{tikzcd}
\left (U_i\times_XU \right )_{\text{ét}} \ar{rr} \ar{d}{s_i} \ar{drr}{s_i} & & U_{\text{ét}} \ar{d}{s} \\
\left (\left (U_i\times_XU \right )_{\text{ét}} \right )_{/\upmu_{|U_i\times_XU}} \ar[swap]{rr}{j_{\upmu\times (U_i\times_XU)/\upmu}} & & X_{/\upmu},
\end{tikzcd}
$$
where the upper arrow is the canonical morphism and the lower morphism is given in \ref{defresind}, is commutative so
$$
\left (\calE_s\right )_{|U_i\times_XU}=\left (\calE_{|U_i\times_XU} \right )_{s_i}=\calE_{s_i}.
$$
The result follows.
\end{proof}

\begin{parag}
Let $\widehat{\Gamma}$ be as defined in \ref{cor1713}. By (\cite{Schepler} 2.2), the isomorphism \eqref{iso17131} induces an equivalence of categories
\begin{equation}\label{ren6}
\begin{array}[t]{clc}
\op{Mod}\left (\calB_{X/S}\otimes_{\Ox_{X',\upmu}}\widehat{\Gamma} \right ) & \xrightarrow{\sim} & \op{Mod} \left ( (\calB_{X/S}\otimes_{\Ox_{X',\upmu}}\widehat{\Gamma})\circledast_{\calB_{X/S}\otimes_{\Ox_{X',\upmu}}S}\widetilde{\calD}_{X/S} \right ) \\
\calE & \mapsto & \left ( \calA_X\otimes_{\Ox_{X,\upmu}}F_1^*\widehat{\Gamma}\right ) \circledast_{\calB_{X/S}\otimes_{\Ox_{X',\upmu}}\widehat{\Gamma}} \calE.
\end{array} 
\end{equation}
The inverse functor is given by
\begin{equation}\label{ren7}
\F \mapsto \mathscr{Hom}_{(\calB_{X/S}\otimes_{\Ox_{X',\upmu}}\widehat{\Gamma})\circledast_{\calB_{X/S}\otimes_{\Ox_{X',\upmu}}S}\widetilde{\calD}_{X/S}} \left (\calA_X\otimes_{\Ox_{X,\upmu}}F_1^*\widehat{\Gamma},\F \right ).
\end{equation}
\end{parag}

\begin{proposition}\label{propKoko1717}
The equivalence \eqref{ren6} induces and equivalence between the full subcategories of locally PD-nilpotent modules:
$$
\Phi_{\upmu}:\begin{Bmatrix}\mathrm{Locally\ PD}\text{-}\mathrm{nilpotent}\\
\upmu\text{-}\mathrm{indexed\ }\left (\calB_{X/S}\otimes_{\Ox_{X',\upmu}}\widehat{\Gamma}\right )\text{-}\\ \mathrm{modules} \end{Bmatrix} \xrightarrow{\sim} \begin{Bmatrix} \mathrm{Locally\ PD}\text{-}\mathrm{nilpotent\ }\upmu\text{-}\mathrm{indexed}\\ \left ((\calB_{X/S}\otimes_{\Ox_{X',\upmu}}\widehat{\Gamma})\circledast_{\calB_{X/S}\otimes_{\Ox_{X',\upmu}}S}\widetilde{\calD}_{X/S} \right )\text{-} \\ \mathrm{modules}\end{Bmatrix}.
$$
\end{proposition}

\begin{proof}
It is sufficient to prove that the functors \eqref{ren6} and \eqref{ren7} preserve locally PD-nilpotent modules.
To simplify notations, we set
$$
\calB_{\Gamma}=\calB_{X/S}\otimes_{\Ox_{X',\upmu}}\widehat{\Gamma},\ \calB_S=\calB_{X/S}\otimes_{\Ox_{X',\upmu}}S,\ \calA_{\Gamma}=\calA_X\otimes_{\Ox_{X,\upmu}}F_1^*\widehat{\Gamma}.
$$
Let $\calE$ be a locally PD-nilpotent object of $\op{Mod}\left (\calB_{\Gamma}\right ).$ Since the property of being locally PD-nilpotent is local \eqref{locPDnilprop} and $\calA_{\Gamma}$ is a locally free $\upmu$-indexed $\calB_{\Gamma}$-module of finite rank \eqref{thmlor}, we can suppose that there exists sections $m_1,\hdots,m_d\in \Gamma(X,\upmu)$ such that
$$
\calA_{\Gamma}=\bigoplus_{i=1}^d\calB_{\Gamma}(m_i).
$$
Then, by \ref{ren5}, we have
$$
\calA_{\Gamma}\circledast_{\calB_{\Gamma}}\calE=\bigoplus_{i=1}^d\calE(m_i).
$$
It is then clear that \eqref{ren6} preserves locally PD-nilpotent modules.

Suppose now that $\F$ is a locally PD-nilpotent module of $\op{Mod}\left (\calB_{\Gamma} \circledast_{\calB_S}\widetilde{\calD}_{X/S} \right ).$ Let $s\in \Gamma(U,\upmu)$ be a section over an étale $X$-scheme, $V$ an étale $U$-scheme and
$$
f\in \mathscr{Hom}_{\calB_{\Gamma}\circledast_{\calB_{S}}\widetilde{\calD}_{X/S}} \left (\calA_{\Gamma},\F \right )_s(V).$$
By replacing $s$ with $s_{|V},$ we can suppose that $V=U.$ By \ref{indremtakriz1} (1), $f$ is a morphism of $\upmu_{|U}$-indexed $\left (\calB_{\Gamma}\circledast_{\calB_{S}}\widetilde{\calD}_{X/S}\right )_{|U}$-modules
$$
f:\calA_{\Gamma|U} \ra \sigma_s^*\F_{|U},
$$
where $\sigma_s:\upmu_{|U} \ra \upmu_{|U}$ is addition by $s.$ After shrinking $U,$ we can suppose that $\calA_{\Gamma|U}$ is free of finite rank over $\calB_{\Gamma|U}.$ Take a basis $(e_i)_{1\le i\le d}$ of the $\upmu_{|U}$-indexed $\calB_{\Gamma|U}$-module $\calA_{\Gamma|U}$ \eqref{deflocfreeind}. There exists local sections $s_i$ of $\upmu$ such that $e_i\in \calA_{\Gamma,s_i}$ for all $i.$ Since $f(e_i)$ is locally annahilated by some power of the augmentation ideal of $\widehat{\Gamma}^{\bullet}\calT_{X'/S},$ we deduce that $\mathscr{Hom}_{\calB_{\Gamma}\circledast_{\calB_{S}}\widetilde{\calD}_{X/S}} \left (\calA_{\Gamma},\F \right )_s$ is locally PD-nilpotent.
\end{proof}

\begin{remark}
Let $(\calE',\theta')$ be an $\upmu$-indexed $\calB_{X/S}$-module equipped with a $\calB_{X/S}$-linear Higgs field $\theta':\calE' \ra \calE'\otimes_{\Ox_{X',\upmu}}\omega^1_{X'/S,\upmu}$ such that the morphism of $\Ox_{X'}$-algebras
$$S=(S^{\bullet}\calT_{X'/S})_{\upmu} \ra \mathscr{End}_{\Ox_{X',\upmu}}(\calE'),$$
induced by $\theta',$ extends to $\w{\Gamma}=(\w{\Gamma}^{\bullet}\calT_{X'/S})_{\upmu}.$ We thus have a $\w{\Gamma}$-module structure on $\calE'.$ Let $\calE=(\calA_X\otimes_{\Ox_{X,\upmu}}F_1^*\widehat{\Gamma}) \circledast_{\calB_{X/S}\otimes_{\Ox_{X',\upmu}}\widehat{\Gamma}}  \calE'=\calA_X\circledast_{\calB_{X/S}}\calE'.$ Suppose that we have local coordinates $m_1,\hdots,m_d \in \Gamma(X,\calM_X)$ and $m_i'=\pi^{\flat}m_i$ (where $\pi$ is defined in \eqref{diag51}). Let $(\partial_i')_{1\le i\le d}$ be the dual basis of $(\op{dlog}m_i')_{1\le i\le d}.$ Let $s$ and $t$ be local sections of $\calI$ over an étale $X$-scheme $U$ and $a\in \calA_{X,s},$ $x'\in \calE'_t.$ By \ref{indproptensprod}, $a\otimes 1\otimes x'\in (\calA_{X,s}\otimes_{\Ox_U} F_1^*\widehat{S}^{\bullet}\calT_{U'/S})\otimes_{\Ox_{U'}} \calE'_t$ can be seen as a local section of
$$\left ((\calA_X\otimes_{\Ox_{X,\upmu}}F_1^*\widehat{\Gamma}) \circledast_{\calB_{X/S}\otimes_{\Ox_{X',\upmu}}\widehat{\Gamma}}  \calE'\right )_{s+t}.$$
Let us explicit the actions of $\partial_{\epsilon_i} \in \calD_{X/S}$ and $\partial_i'^{[k]} \in \Gamma^{\bullet} \calT_{X'/S}$ on $a\otimes 1\otimes x'.$ By \eqref{era2action3}, the local sections $c_1,\hdots,c_d$ of $\Ox_X,$ given in \eqref{eqKoko1323}, satisfy
\begin{alignat*}{2}
\partial_{\epsilon_i} \cdot (a\otimes 1\otimes x') &= (\partial_{\epsilon_i} \cdot a) \otimes 1 \otimes x'+ a\otimes (\partial_{\epsilon_i}\cdot 1)\otimes x' \\
&= (\partial_{\epsilon_i}\cdot a) \otimes 1\otimes x' +a\otimes \left (F_1^*\partial'_{i}+\sum_{j=1}^d\partial_{\epsilon_i}(p_2^{\#}c_j-p_1^{\#}c_j)F_1^*\partial'_j\right )\otimes x',
\end{alignat*}
where $p_1,p_2:P_0 \ra X$ are the canonical projections.
Identifying
$$(\calA_X\otimes_{\Ox_{X,\upmu}}F_1^*\widehat{\Gamma}) \circledast_{\calB_{X/S}\otimes_{\Ox_{X',\upmu}}\widehat{\Gamma}}  \calE'=(\calA_X\otimes_{\Ox_{X',\upmu}}\widehat{\Gamma}) \circledast_{\calB_{X/S}\otimes_{\Ox_{X',\upmu}}\widehat{\Gamma}}  \calE'$$
with $\calA_X\circledast_{\calB_{X/S}}  \calE'$ via the isomorphism \eqref{propisoAkazahh}, we have
\begin{equation*}
\partial_{\epsilon_i} \cdot (a\otimes x') = (\partial_{\epsilon_i}\cdot a) \otimes x' +a\otimes \theta(\partial'_{i})(x')+\sum_{j=1}^d\partial_{\epsilon_i}(p_2^{\#}c_j-p_1^{\#}c_j) a \otimes \theta(\partial'_j)(x').
\end{equation*}
For $1\le i\le d,$ the local section $\theta(\partial_i')(x')$ is equal to the action $\partial_i'\cdot x'$ of $\partial_i'$ on $x'.$ We thus have
\begin{equation}\label{actionpartialax}
\begin{alignedat}{2}
\partial_{\epsilon_i} \cdot (a\otimes x') &= (\partial_{\epsilon_i}\cdot a) \otimes x' +a\otimes (\partial'_{i}\cdot x')+\sum_{j=1}^d\partial_{\epsilon_i}(p_2^{\#}c_j-p_1^{\#}c_j) a \otimes (\partial'_j\cdot x'), \\
\partial_i'^{[k]} \cdot (a \otimes x') &= a \otimes \left ( \partial_i'^{[k]} \cdot x' \right ).
\end{alignedat}
\end{equation}
 
\end{remark}

\newpage

\section{Logarithmic Oyama topoi and the Cartier transform}

\begin{parag}\label{intro18}
In this section, we fix a perfect field $k$ of positive characteristic $p$ and denote its ring of Witt vectors by $W(k).$ We equip $\op{Spf}W(k)$ with the trivial logarithmic structure. Unless explicitly stated, all formal schemes considered are $p$-adic and all logarithmic structures are fs. If a gothic letter $\frakT$ denotes a $p$-adic logarithmic formal scheme, the corresponding roman letter $T$ will denote the logarithmic scheme obtained from $\frakT$ by reduction modulo $p$ \eqref{propkey}. We also denote by $\frakT_n,$ for any positive integer $n,$ the logarithmic scheme obtained from $\frakT$ by reduction modulo $p^n$ \eqref{propkey}.
Let $\frakS$ be a logarithmic formal scheme, log flat and locally of finite type over $\op{Spf}W(k),$ and $X\ra S$ a log smooth morphism of logarithmic schemes.
\end{parag}

\begin{definition}\label{defforintro}
We define the categories $\calE(X/\frakS)$ and $\underline{\calE}(X/\frakS)$ as follows:
\begin{enumerate}
\item An object of $\calE(X/\frakS)$ (resp. $\underline{\calE}(X/\frakS)$) is a triple $(U,\frakT,u)$ consisting of an étale strict morphism of logarithmic schemes $U \ra X,$ a log flat fs $p$-adic formal logarithmic $\frakS$-scheme $\frakT$ (i.e. log flat modulo $p^n$ for all $n\ge 1$ (\ref{logflatdef})) and an $S$-morphism of logarithmic schemes $u:T \ra U$ (resp. $u:\underline{T}\ra U$) which is affine as a morphism of schemes, where $\underline{T}$ is defined in \ref{era3logimagedef}.
\item A morphism $(U_1,\frakT_1,u_1) \ra (U_2,\frakT_2,u_2)$ in $\calE(X/\frakS)$ (resp. $\underline{\calE}(X/\frakS)$) is a pair $(f,g)$ consisting of an $\frakS$-morphism $f:\frakT_1 \ra \frakT_2$ and an $X$-morphism $g:U_1 \ra U_2$ such that $u_2\circ f_1=g\circ u_1$ (resp. $u_2 \circ \underline{f_1}=g\circ u_1$), where $f_1:T_1 \ra T_2$ is the morphism induced by $f$ by reduction modulo $p$ and $\underline{f_1}:\underline{T_1} \ra \underline{T_2}$ is obtained from $f_1$ by functoriality (\ref{era3functoriality}).
\end{enumerate}
A morphism $(U_1,\frakT_1,u_1) \ra (U_2,\frakT_2,u_2)$ in $\calE(X/\frakS)$ or $\underline{\calE}(X/\frakS)$ is said to be \emph{log flat} if $\frakT_1 \ra \frakT_2$ is log flat \eqref{logflatdef}.
\end{definition}

\begin{definition}
We denote by $\calE_{\mathrm{strict}}(X/\frakS)$ (resp. $\underline{\calE}_{\mathrm{strict}}(X/\frakS)$) the subcategory of $\calE(X/\frakS)$ (resp. $\underline{\calE}(X/\frakS)$) consisting of objects $(U,\frakT,u)$ such that $u:T\ra U$ (resp. $u:\underline{T} \ra U$) is strict.
\end{definition}

\begin{parag}\label{era3fiberparag}
Fiber products by log flat morphisms in the categories $\calE(X/\frakS),$ $\underline{\calE}(X/\frakS),$ $\calE_{\mathrm{strict}}(X/\frakS)$ and $\underline{\calE}_{\mathrm{strict}}(X/\frakS)$ are representable. More precisely, given morphisms $(U_1,\frakT_1,u_1) \ra (U,\frakT,u)$ and $(U_2,\frakT_2,u_2) \ra (U,\frakT,u)$ in $\calE(X/\frakS)$ (resp. $\underline{\calE}(X/\frakS)$) with one of them being log flat, the fiber product $(U_1,\frakT_1,u_1)\times_{(U,\frakT,u)}(U_2,\frakT_2,u_2)$ is representable by $(U_1\times_UU_2,\frakT_1\times_{\frakT}^{\op{log}}\frakT_2,v)$ where $v$ is the morphism $T_1\times_T^{\op{log}}T_2 \ra U_1\times_UU_2$ induced by $u_1$ and $u_2$ (resp. is the composition of the mophism $\underline{T_1}\times_{\underline{T}}^{\op{log}}\underline{T_2} \ra U_1\times_UU_2$ induced by $u_1$ and $u_2$ with the canonical morphism $\underline{T_1\times_T^{\op{log}}T_2} \ra \underline{T_1}\times_{\underline{T}}^{\op{log}}\underline{T_2}$ (\ref{era3functoriality})), which is affine. Indeed, since one of the morphisms $\frakT_1 \ra \frakT$ and $\frakT_2 \ra \frakT$ is log flat, $\frakT_1\times_{\frakT}^{\op{log}}\frakT_2 \ra \frakS$ is also log flat (\cite{Ogus2018} IV 4.1.2 (3)). Note that $U_1\times_UU_2=U_1\times_U^{\op{log}}U_2$ since $U\ra X,$ $U_1\ra X$ and $U_2 \ra X$ are all strict. In addition, if $u_1$ and $u_2$ are strict, then so is $v.$
\end{parag}

\begin{definition}\label{era3defcart}~ 
\begin{enumerate}
\item A morphism $(U_1,\frakT_1,u_1) \ra (U_2,\frakT_2,u_2)$ in $\calE(X/\frakS)$ is said to be \emph{cartesian} if the canonical morphism
$$T_1 \ra T_2\times_{U_2}U_1$$
is an isomorphism.
\item A morphism $(U_1,\frakT_1,u_1) \ra (U_2,\frakT_2,u_2)$ in $\underline{\calE}(X/\frakS)$ is said to be \emph{cartesian} if the morphism $\frakT_1 \ra \frakT_2$ is strict and étale and the square
$$
\begin{tikzcd}
T_1 \ar{r} \ar[swap]{d}{\Delta_{u_1}} & T_2 \ar{d}{\Delta_{u_2}} \\
U_1' \ar{r} & U_2',
\end{tikzcd}
$$
where the horizontal arrows are induced by $\frakT_1 \ra \frakT_2$ and $U_1 \ra U_2$ and the vertical arrows are given in \eqref{morphismShiho1}, is cartesian.
\end{enumerate}
\end{definition}

\begin{proposition}\label{era3stablebc}
The class of cartesian morphisms in $\calE(X/\frakS)$ (resp. $\underline{\calE}(X/\frakS),$ resp. $\calE_{\mathrm{strict}}(X/\frakS),$ resp. $\underline{\calE}_{\mathrm{strict}}(X/\frakS)$) is stable under composition and base change by any morphism.
\end{proposition}

\begin{proof}
Stability under composition in all categories is clear and so is stability by base change for $\calE(X/\frakS)$ and $\calE_{\mathrm{strict}}(X/\frakS).$
Let $(U_1,\frakT_1,u_1) \ra (U,\frakT,u)$ and $(U_2,\frakT_2,u_2) \ra (U,\frakT,u)$ be morphisms in $\underline{\calE}(X/\frakS)$ such that $(U_1,\frakT_1,u_1) \ra (U,\frakT,u)$ is cartesian. By \ref{era3fiberparag}, the fiber product is given by
$$(U_1,\frakT_1,u_1)\times_{(U,\frakT,u)}(U_2,\frakT_2,u_2)=(V,\frakZ,v),$$
where $V=U_1\times_UU_2,$ $\frakZ=\frakT_1\times_{\frakT}^{\op{log}}\frakT_2$ and $v:\underline{T_1\times_T^{\op{log}}T_2}\ra U_1\times_UU_2$ is the composition of the canonical morphism $\underline{T_1\times^{\op{log}}_TT_2} \ra \underline{T_1}\times_{\underline{T}}^{\op{log}}\underline{T_2}$ with the morphism $\underline{T_1}\times_{\underline{T}}^{\op{log}}\underline{T_2} \ra U_1\times_UU_2$ induced by $u_1$ and $u_2.$ Considering the following commutative diagram
$$
\begin{tikzcd}
 & T_1\times_T^{\op{log}}T_2 \ar{rr} \ar[dashed]{ddd}{\Delta_{v}} \ar{dl} & & T_2 \ar{dl} \ar{ddd}{\Delta_{u_2}} \\
T_1 \ar{rr} \ar[swap]{ddd}{\Delta_{u_1}} & & T \ar[swap]{ddd}{\Delta_{u}} & \\
 &  & &  \\
 & U_1'\times_{U'} U_2' \ar[dashed]{rr} \ar[dashed]{dl} & & U_2' \ar{dl} \\
 U_1' \ar{rr} & & U', &
\end{tikzcd}
$$
we get that the diagram
$$
\begin{tikzcd}
T_1\times_T^{\op{log}}T_2 \ar{r} \ar[swap]{d}{\Delta_{v}} & T_2 \ar{d}{\Delta_{u_2}} \\
U_1'\times_{U'}U_2' \ar{r} & U_2'
\end{tikzcd}
$$
is cartesian. The stability for $\underline{\calE}_{\mathrm{strict}}(X/\frakS)$ follows immediately.
\end{proof}

\begin{parag}\label{era3paralpha}
Let $(U,\frakT,u)$ be an object of $\calE(X/\frakS)$ (resp. $\underline{\calE}(X/\frakS)$). We define a functor
\begin{equation}\label{eqKoko1861}
\alpha_{(U,\frakT,u)}:\begin{array}[t]{cll}
\text{ét}_{/U} & \ra & \calE(X/\frakS)\ (\op{resp.} \underline{\calE}(X/\frakS))\\ 
V & \mapsto & (V,\frakT_V,u_V)
\end{array}
\end{equation}
as follows: let $V$ be an étale $U$-scheme. We equip $V$ with the logarithmic structure pullback of that of $U,$ so that $V\ra U$ becomes a strict morphism. First, we construct $\alpha_{(U,\frakT,u)}$ for $\calE(X/\frakS).$ Let $T_V=T\times_UV$ and $u_v:T_V \ra V$ the canonical projection, which is affine as a base change of the affine morphism $u:T \ra U.$ Since $V\ra U$ is an étale strict morphism, so is $T_V \ra T.$ Then there exists a unique $p$-adic formal scheme $\frakT_V$ that fits into a cartesian square of formal schemes
$$
\begin{tikzcd}
T_V \ar{r} \ar{d} & \frakT_V \ar{d} \\
T \ar{r} & \frakT
\end{tikzcd}
$$
We equip $\frakT_V$ with the logarithmic structure pullback of that of $\frakT.$ This makes the previous square cartesian in the category of logarithmic formal schemes by \ref{propstrcart}. In addition, $\frakT_V \ra \frakT$ is étale and strict so $\frakT_V$ is log flat over $\frakS.$ We get an object $(V,\frakT_V,u_V)$ of $\calE(X/\frakS)$ and we set
\begin{equation}\label{eqKoko1862}
\alpha_{(U,\frakT,u)}(V)=(V,\frakT_V,u_V).
\end{equation}
We now prove the result for $\underline{\calE}(X/\frakS).$ Let $T_V=T\times_{U'}V',$ where $T \ra U'$ is the morphism $\Delta_{u}$ \eqref{morphismShiho1}. Since $V' \ra U'$ is strict, $T_V\ra T$ is also strict and so $T_V$ is fs. The canonical projection $T_V\ra T$ being étale, we construct $\frakT_V$ the same way we did in the previous case. 
Let $T \ra \underline{T} \ra T$ and $T_V \ra \underline{T_V} \ra T_V$ be the factorizations of $F_T$ and $F_{T_V}$ through the canonical immersions $\underline{T} \ra T$ and $\underline{T_V} \ra T_V.$ The morphism $T_V \ra \underline{T_V}$ is inseparable. By \ref{propShiho185}, the diagram
$$
\begin{tikzcd}
T \ar{r}{\Delta_{u}} \ar{d} & U' \ar{d} \\
\underline{T} \ar{r}{u} & U
\end{tikzcd}
$$ 
is commutative. In addition, so are the diagrams
$$
\begin{tikzcd}
V' \ar{r} \ar{d} & V \ar{d} & T_V \ar{r} \ar{d} & \underline{T_V} \ar{d} \\
U' \ar{r} & U & T \ar{r} & \underline{T}.
\end{tikzcd}
$$
It follows that the diagram
$$
\begin{tikzcd}
V' \ar{rrr} & & & V \ar{d} \\
T_V \ar{u} \ar{r} \ar{drr} & \underline{T_V} \ar[dashed]{urr}{u_V} \ar{r} & \underline{T} \ar{r}{u} & U \\
 & & T \ar{u} \ar{r}{\Delta_{u}} & U' \ar{u}
\end{tikzcd}
$$
is commutative. The existence of the dashed arrow $u_V$ follows then from the fact that $V \ra U$ is étale and strict, the fact that $T_V \ra \underline{T_V}$ is inseparable and (\cite{Ogus2018} IV 3.3.7).
If $U$ and $V$ are affine then so is $u_V.$ Note that, if $u$ is strict, then so is $u_V.$
We set
\begin{equation}\label{eqKoko1863}
\alpha_{(U,\frakT,u)}(V)=(V,\frakT_V,u_V).
\end{equation}
If $V_1 \ra V_2$ is a morphism of étale and strict $U$-schemes then we have cartesian squares
$$
\begin{tikzcd}
T_{V_1} \ar{r} \ar{d} & T_{V_2} \ar{r} \ar{d} & T \ar{d}{\Delta_u} \\
V_1' \ar{r} & V_2' \ar{r} & U'.
\end{tikzcd}
$$
By the étaleness and strictness of $\frakT_{V_2} \ra \frakT,$ we have a unique morphism $\frakT_{V_1} \ra \frakT_{V_2}$ fitting into the commutative diagram
$$
\begin{tikzcd}
T_{V_2} \ar{rr} & & \frakT_{V_2} \ar{d} \\
T_{V_1} \ar{u} \ar{r} & \frakT_{V_1} \ar{r} \ar{ur} & \frakT.
\end{tikzcd}
$$
By construction of $u_{V_1},$ we have a commutative diagram
$$
\begin{tikzcd}
V_1' \ar{rrr} & &  & V_1 \ar{r} \ar{d} & V_2 \ar{dl} \\
T_{V_1} \ar{u} \ar{r}{G_{T_{V_1}}} & \underline{T_{V_1}} \ar{r} \ar{urr}{u_{V_1}} & \underline{T} \ar{r}{u} & U. &
\end{tikzcd}
$$
We also have the commutative diagram
$$
\begin{tikzcd}
 & V_2' \ar{rrr} & & & V_2 \ar{dd} \\
V_1' \ar{ur} & T_{V_2} \ar{u} \ar{r} & \underline{T_{V_2}} \ar{urr}{u_{V_2}} \ar{dr} & &  \\
 & T_{V_1} \ar{ul} \ar{u} \ar{r} & \underline{T_{V_1}} \ar{u} \ar{r} & \underline{T} \ar{r}{u} & U.
\end{tikzcd}
$$
The morphism $V_2 \ra U$ is étale and strict so the diagram
$$
\begin{tikzcd}
\underline{T_{V_1}} \ar[swap]{d}{u_{V_1}} \ar{r} & \underline{T_{V_2}} \ar{d}{u_{V_2}} \\
V_1 \ar{r} & V_2
\end{tikzcd}
$$
is commutative.
We hence get a morphism
$$
\alpha_{(U,\frakT,u)}(V_1)=\left ( V_1,\frakT_{V_1},u_{V_1} \right ) \ra \left ( V_2,\frakT_{V_2},u_{V_2} \right )=\alpha_{(U,\frakT,u)}(V_2).
$$
This defines the functor $\alpha_{(U,\frakT,u)}.$
\end{parag}

\begin{remark}\label{remShiho1}~ 
\begin{enumerate}
\item If $(U,\frakT,u)$ is an object of $\calE(X/\frakS)$ (resp. $\underline{\calE}(X/\frakS)$) and if $V$ is an étale an strict $U$-scheme and $(V,\frakT_V,u_V)=\alpha_{(U,\frakT,u)}(V),$ then the construction of $u_V$ shows that $u_V$ is strict if $u$ is strict. It follows that, if $(U,\frakT,u)$ is an object of $\calE_{\mathrm{strict}}(X/\frakS)$ (resp. $\underline{\calE}_{\mathrm{strict}}(X/\frakS)$), then $\alpha_{(U,\frakT,u)}$ factors through $\calE_{\mathrm{strict}}(X/\frakS)$ (resp. $\underline{\calE}_{\mathrm{strict}}(X/\frakS)$).
\item Suppose that $\frakS=\op{Spf}W(k)$ equipped with the trivial log structure. Let $(U,\frakT,u)$ be an object of $\underline{\calE}(X/\frakS).$ Here, we provide a simpler definition of the functor
$$
\alpha_{(U,\frakT,u)}:\text{ét}_{/U} \ra \underline{\calE}(X/\frakS).
$$
Recall that, in this case, $\Delta_{u}=u' \circ f_{T/S}.$ Let $V$ be an étale $U$-scheme. Set $T_V=T\times_{U'}V',$ where $T \ra U'$ is the composition $u' \circ f_{T/S}.$ Since $V' \ra U'$ is strict, $T_V$ is fs. The canonical projection $T_V\ra T$ being étale, we construct $\frakT_V$ the same way we did in \ref{era3paralpha}. The morphism $f_{T_V/S}$ is inseparable \eqref{propfT/S} and $V' \ra U'$ is étale and strict. It follows, by (\cite{Ogus2018} IV 3.3.7), that there exists a unique morphism $v:\underline{T_V}' \ra V'$ fitting into the following commutative diagram
$$
\begin{tikzcd}
T_V \ar[swap, bend right=70]{dd} \ar{r} \ar[swap]{d}{f_{T_V/S}} & T  \ar{d}{f_{T/S}} \\
 \underline{T_V}' \ar{r} \ar[dashed,swap]{d}{v} & \underline{T}' \ar{d}{u'} \\
V' \ar{r} & U'.
\end{tikzcd}
$$
If $V$ is affine then $T_V=V'\times_{U',u'\circ f_{T/S}}T$ is affine and so so is $\underline{T_V}'.$ It follows that $v$ is affine.
There exists a unique morphism of logarithmic schemes $u_V:\underline{T_V} \ra V$ such that $v=u_V',$ which is affine since $v$ is affine.
We set
\begin{equation}
\alpha_{(U,\frakT,u)}(V)=(V,\frakT_V,u_V).
\end{equation}
\end{enumerate}
\end{remark}

\begin{proposition}\label{era3paralphaKoko}
Let $(U,\frakT,u)$ be an object of $\calE(X/\frakS)$ (resp. $\underline{\calE}(X/\frakS)$). The functor $\alpha_{(U,\frakT,u)}$ \eqref{eqKoko1861} satisfies the following properties
\begin{enumerate}
\item For any étale $U$-scheme $V,$ the canonical morphism $(V,\frakT_V,u_V)=\alpha_{(U,\frakT,u)}(V) \ra (U,\frakT,u)$ is cartesian and $\frakT_V \ra \frakT$ is étale and strict.
\item For any étale $U$-scheme $V$ and any morphism $a:(W,\frakZ,w) \ra \alpha_{(U,\frakT,u)}(V)$ in $\calE(X/\frakS)$ (resp. $\underline{\calE}(X/\frakS)$), there exists a unique morphism $b:(W,\frakZ,w) \ra \alpha_{(U,\frakT,u)}(W)$ such that the diagram
$$
\begin{tikzcd}
(W,\frakZ,w) \ar{r}{a} \ar{d}{b} & \alpha_{(U,\frakT,u)}(V) \\
\alpha_{(U,\frakT,u)}(W) \ar{ur} & 
\end{tikzcd}
$$
is commutative. In addition, $a$ is cartesian if and only if $b$ is cartesian.
\item For any $U$-morphism of étale $U$-schemes $V_1 \ra V_2,$ the morphism $\alpha_{(U,\frakT,u)}(V_1) \ra \alpha_{(U,\frakT,u)}(V_2)$ is cartesian.
\end{enumerate}
\end{proposition}

\begin{proof}
The first assertion is clear from the construction of $\alpha_{(U,\frakT,u)}.$ The existence of $b$ in the second assertion follows, in the case $\underline{\calE}(X/\frakS),$ from the following fact: let $(V,\frakT_V,u_V)=\alpha_{(U,\frakT,u)}(V)$ and $(W,\frakT_W,u_W)=\alpha_{(U,\frakT,u)}(W).$ By definition, we have cartesian squares
\begin{equation}\label{diagKoko1871}
\begin{tikzcd}
Z \ar[swap]{dr}{\Delta_{(W,\frakZ,w)}} \ar{r} \ar[bend right=-30]{rrr} & T_W \ar{r} \ar{d} & T_V \ar{d} \ar{r} & T \ar{d}{\Delta_{u}} \\
 & W' \ar{r} & V' \ar{r} & U',
\end{tikzcd}
\end{equation}
$$
\begin{tikzcd}
Z\ar{d} \ar[bend right=30]{ddd} \ar{r} & \frakZ \ar{d} \\
T_W \ar{r} \ar{d} & \frakT_W \ar{d} \\
T_V \ar{r} \ar{d} & \frakT_V \ar{d} \\
T \ar{r} & \frakT.
\end{tikzcd}
$$
The case $\calE(X/\frakS)$ is simpler.
The commutativity of \eqref{diagKoko1871} proves that $a$ is cartesian if and only if $b$ is cartesian. The third assertion follows from the first.
\end{proof}

\begin{parag}\label{era3parbeta}
Let $(f,g):(U_1,\frakT_1,u_1) \ra (U_2,\frakT_2,u_2)$ be a morphism in $\calE(X/\frakS)$ (resp. $\underline{\calE}(X/\frakS)$) and
\begin{equation}\label{jg}
j_g:\begin{array}[t]{clc}
\text{ét}_{/U_1} & \ra & \text{ét}_{/U_2} \\
(V\ra U_1) & \mapsto & (V \ra U_1 \xrightarrow[g]{} U_2).
\end{array}
\end{equation}
The morphism $(f,g)$ induces a morphism
\begin{equation}\label{defbeta115}
\beta_{(f,g)}:\alpha_{(U_1,\frakT_1,u_1)} \ra \alpha_{(U_2,\frakT_2,u_2)}\circ j_g
\end{equation}
as follows (we will only give the construction for $\underline{\calE}(X/\frakS)$): let $V$ be an étale $U_1$-scheme and $(V,\frakT_{i_V},v_i)=\alpha_{(U_i,\frakT_i,u_i)}(V).$
The morphism $T_1 \ra T_2$ induces
$$T_{1_V}=T_1\times_{U_1',\Delta_{u_1}}V' \ra T_2 \times_{U_2',\Delta_{u_2}}V'=T_{2_V}$$
fitting into the commutative diagram
$$
\begin{tikzcd}
T_{1_V} \ar{rr} \ar{dd} \ar{dr} & & T_{2_V} \ar{dr} \ar[dashed]{dd} & \\
 & T_1 \ar{rr}\ar{dd} & & T_2\ar{dd} \\
V' \ar{dr} \ar[dashed, equal]{rr} & & V' \ar{dr} & \\
 & U_1' \ar{rr} & & U_2'.
\end{tikzcd}
$$
This morphism extends, by étaleness and strictness of $\frakT_{2_V} \ra \frakT_2,$ to a morphism of $p$-adic logarithmic formal schemes $h:\frakT_{1_V} \ra \frakT_{2_V}$ over $\frakT_2.$
We hence obtain a morphism $(h,\op{Id}_V):(V,\frakT_{1_V},v_1) \ra (V,\frakT_{2_V},v_2)$ over $(f,g).$
We easily prove that the morphisms $\beta_{(f,g)}$ satisfy the following properties:
\begin{enumerate}
\item For every object $(U,\frakT,u)$ of $\calE(X/\frakS)$ (resp. $\underline{\calE}(X/\frakS)$),
$$\beta_{\op{Id}_{(U,\frakT,u)}}=\op{Id}_{\alpha_{(U,\frakT,u)}}.$$
\item For all composable morphisms $(f_1,g_1):(U_1,\frakT_1,u_1) \ra (U_2,\frakT_2,u_2)$ and $(f_2,g_2):(U_2,\frakT_2,u_2) \ra (U_3,\frakT_3,u_3)$ in $\calE(X/\frakS)$ (resp. $\underline{\calE}(X/\frakS)$), the following diagram
$$
\begin{tikzcd}
\alpha_{(U_1,\frakT_1,u_1)} \ar{r}{\beta_{(f_1,g_1)}} \ar[swap]{dr}{\beta_{(f_2,g_2) \circ (f_1,g_1)}} & \alpha_{(U_2,\frakT_2,u_2)} \circ j_{g_1} \ar{d}{j_{g_1}^*\beta_{(f_2,g_2)}} \\
& \alpha_{(U_3,\frakT_3,u_3)} \circ j_{g_2\circ g_1},
\end{tikzcd}
$$
where $j_{g_1}^*\beta_{(f_2,g_2)}$ is the morphism induced by $\beta_{(f_2,g_2)},$ is commutative up to a canonical isomorphism
\item If $(f,g)$ is a cartesian morphism, then $\beta_{(f,g)}$ is an isomorphism.
\end{enumerate}
We just check the third point: if $(f,g)$ is a cartesian morphism in $\underline{\calE}(X/\frakS)$ then
$$T_{1_V}=T_1\times_{U_1'}V'=T_2\times_{U_2'}V'=T_{2_V}.$$
It follows that $\beta_{(f,g)}$ is an isomorphism. This is also true in the category $\calE(X/\frakS).$
\end{parag}

\begin{proposition}\label{propalphacomfiber}
Let $(U,\frakT,u)$ be an object of $\calE(X/\frakS)$ (resp. $\underline{\calE}(X/\frakS)$). The functor $\alpha_{(U,\frakT,u)}$ \eqref{eqKoko1861} commutes with fiber products.
\end{proposition}

\begin{proof}
For any étale $U$-scheme $V,$ we set $(V,\frakT_V,u_V)=\alpha_{(U,\frakT,u)}(V).$ 
Let $V_1 \ra V$ and $V_2\ra V$ be étale morphisms of $U$-schemes. We want to prove that
$$\alpha_{(U,\frakT,u)}(V_1\times_VV_2)=\alpha_{(U,\frakT,u)}(V_1)\times_{\alpha_{(U,\frakT,u)}(V)}\alpha_{(U,\frakT,u)}(V_2).$$
Since
$$(V,\frakT_V,u_V)=\alpha_{(V,\frakT_V,u_V)}(V),$$
we can suppose that $U=V$ and hence $(V,\frakT_V,u_V)=(U,\frakT,u).$
We prove the result for $\underline{\calE}(X/\frakS)$ since the other case is simpler.
Consider the following commutative diagram with cartesian squares
$$
\begin{tikzcd}
T_{V_1\times_UV_2} \ar{rr} \ar{dd} \ar{dr} & & T_{V_2} \ar{dr} \ar[dashed]{dd}& \\
 & T_{V_1} \ar{rr} \ar{dd} & & T\ar{dd}{\Delta_u} \\
V_1'\times_{U'}V_2' \ar[dashed]{rr} \ar{dr} & & V_2' \ar[dashed]{dr} & \\
 & V_1' \ar{rr} & & U'.
\end{tikzcd}
$$
It follows that
\begin{equation}\label{17181}
T_{V_1\times_UV_2}=T_{V_1}\times_T^{\op{log}}T_{V_2}
\end{equation}
and so, by definition of $\alpha_{(U,\frakT,u)},$
$$\frakT_{V_1\times_UV_2}=\frakT_{V_1}\times_{\frakT}^{\op{log}}\frakT_{V_2}.$$
By definition, $u_{V_1\times_UV_2}$ is the unique morphism fitting in the following commutative diagram
\begin{equation}\label{diagS1}
\begin{tikzcd}
\left ( V_1\times_UV_2 \right )' \ar{rrr} & & & V_1\times_{U}V_2\ar{d} \\
T_{V_1\times_UV_2} \ar{u}{\Delta_{\left ( V_1\times_UV_2,T_{V_1\times_UV_2},u_{V_1\times_UV_2}\right )}}  \ar{r} & \underline{T_{V_1\times_UV_2}} \ar{urr}{u_{V_1\times_UV_2}} \ar{r} & \underline{T} \ar{r}{u} & U,
\end{tikzcd}
\end{equation}
where $T_{V_1\times_UV_2} \ra \underline{T_{V_1\times_UV_2}}$ is obtained from the factorization of $F_{T_{V_1\times_UV_2}}$ through $\underline{T_{V_1\times_UV_2}}.$ Similarly, $u_{V_i}$ fits into the commutative diagram
$$
\begin{tikzcd}
V'_i \ar{rrr} & & & V_i \ar{d} \\
T_{V_i} \ar{u} \ar{r}  & \underline{T_{V_i}} \ar[dashed]{urr}{u_{V_i}} \ar{r} & \underline{T} \ar{r}{u} & U.
\end{tikzcd}
$$
By \eqref{17181},
$$\underline{T_{V_1\times_UV_2}}=\underline{T_{V_1}\times_T^{\op{log}}T_{V_2}}.$$
We deduce that $u_{V_1} \times u_{V_2}$ fits also into the commutative diagram \eqref{diagS1}. Hence
$$
u_{V_1\times_UV_2}=u_{V_1} \times u_{V_2}.
$$
This finishes the proof.
\end{proof}

\begin{proposition}\label{era3propdescentdata}
A presheaf $\F$ on $\calE(X/\frakS)$ (resp. $\underline{\calE}(X/\frakS)$) is equivalent to the following data:
\begin{enumerate}
\item For every object $(U,\frakT,u)$ of $\calE(X/\frakS)$ (resp. $\underline{\calE}(X/\frakS)$), a presheaf $\F_{(U,\frakT,u)}$ on $\text{ét}_{/U}.$
\item For every morphism $(f,g):(U_1,\frakT_1,u_1) \ra (U_2,\frakT_2,u_2)$ in $\calE(X/\frakS)$ (resp. $\underline{\calE}(X/\frakS)$), a morphism of presheaves on $\text{ét}_{/U_1},$
$$\gamma_{\F,(f,g)}:g^{-1}\F_{(U_2,\frakT_2,u_2)} \ra \F_{(U_1,\frakT_1,u_1)},$$
\end{enumerate}
such that
\begin{enumerate}[(i)]
\item For every object $(U,\frakT,u)$ of $\calE(X/\frakS)$ (resp. $\underline{\calE}(X/\frakS)$),
$$\gamma_{\F,\op{Id}_{(U,\frakT,u)}}=\op{Id}_{\F_{(U,\frakT,u)}}.$$
\item For all composable morphisms $(f_1,g_1):(U_1,\frakT_1,u_1) \ra (U_2,\frakT_2,u_2)$ and $(f_2,g_2):(U_2,\frakT_2,u_2) \ra (U_3,\frakT_3,u_3)$ in $\calE(X/\frakS)$ (resp. $\underline{\calE}(X/\frakS)$),
$$\gamma_{\F,(f_1,g_1)} \circ g_1^{-1}\gamma_{\F,(f_2,g_2)}=\gamma_{\F,(f_2,g_2)\circ (f_1,g_1)}.$$
\item If $(f,g)$ is a cartesian morphism, then $\gamma_{\F,(f,g)}$ is an isomorphism.
\end{enumerate}
This equivalence is given as follows: for a presheaf $\F$ on $\calE(X/\frakS)$ (resp. $\underline{\calE}(X/\frakS)$),
$$\F_{(U,\frakT,u)}=\F \circ \alpha_{(U,\frakT,u)}$$
and $\gamma_{\F,(f,g)}$ is induced by $\beta_{(f,g)}$ \eqref{defbeta115}. Conversely, given a data $(\F_{(U,\frakT,u)},\gamma_{\F,(f,g)})$ as above, the presheaf is defined by
$$\F(U,\frakT,u)=\F_{(U,\frakT,u)}(U)$$
and, for a morphism $(f,g):(U_1,\frakT_1,u_1) \ra (U_2,\frakT_2,u_2),$ the morphism
$$\F(U_2,\frakT_2,u_2) \ra \F(U_1,\frakT_1,u_1)$$
is equal to $\gamma_{\F,(f,g)}(U_1).$
\end{proposition}

\begin{proof}
Suppose we are given the data $(\F_{(U,\frakT,u)},\gamma_{\F,(f,g)}).$ Conditions $(i)$ and $(ii)$ imply that the correspondance
$$(U,\frakT,u) \mapsto \F_{(U,\frakT,u)}(U)$$
is a presheaf.
Conversely, if $\F$ is a presheaf on $\calE(X/\frakS)$ (resp. $\underline{\calE}(X/\frakS)$) then, for any object $(U,\frakT,u)$ of $\calE(X/\frakS)$ (resp. $\underline{\calE}(X/\frakS)$), $\F_{(U,\frakT,u)}=\F \circ \alpha_{(U,\frakT,u)}$ is clearly a presheaf on $\text{ét}_{/U}.$ Conditons $(i),$ $(ii)$ and $(iii)$ are satisfied by the fact that $\F$ is a presheaf and \ref{era3parbeta}.

We check that the constructions are inverse to each other: let $\F$ be a presheaf on $\calE(X/\frakS)$ and $\calG$ the presheaf corresponding to $\left (\F_{(U,\frakT,u)}\right ).$ For any object $(U,\frakT,u)$ of $\calE(X/\frakS),$ we have
$$\calG(U,\frakT,u)=\F_{(U,\frakT,u)}(U)=\F(U,\frakT,u)$$
so $\calG=\F.$

Conversely, let $(\calG_{(U,\frakT,u)})$ be data satisfying conditions $(i),$ $(ii)$ and $(iii)$ and let $\F$ be the corresponding presheaf. Let $(U,\frakT,u)$ be an object of $\calE(X/\frakS),$ $V$ an étale $U$-scheme and $(V,\frakT_V,u_V)=\alpha_{(U,\frakT,u)}(V).$ Then
$$\F_{(U,\frakT,u)}(V)=\F(V,\frakT_V,u_V)=\calG_{(V,\frakT_V,u_V)}(V).$$
By condition $(iii),$ we have
$$\calG_{(V,\frakT_V,u_V)}(V)=\calG_{(U,\frakT,u)}(V)$$
so $\F_{(U,\frakT,u)}=\calG_{(U,\frakT,u)}.$
\end{proof}

\begin{remark}\label{rem1813}
Proposition \ref{era3propdescentdata} remains true if we replace $\calE(X/\frakS)$ and $\underline{\calE}(X/\frakS)$ with $\calE_{\mathrm{strict}}(X/\frakS)$ and $\underline{\calE}_{\mathrm{strict}}(X/\frakS)$ respectively.
\end{remark}

\begin{definition}\label{era3defdescentdata}
For a presheaf of sets $\F$ on $\calE(X/\frakS)$ (resp. $\underline{\calE}(X/\frakS),$ resp. $\calE_{\mathrm{strict}}(X/\frakS),$ resp. $\underline{\calE}_{\mathrm{strict}}(X/\frakS)$), we call \emph{descent data associated with $\F$} the data $(\F_{(U,\frakT,u)},\gamma_{\F,(f,g)})$ given in \ref{era3propdescentdata} and \ref{rem1813}.
\end{definition}

\begin{parag}\label{parettop}
For any object $(U,\frakT,u)$ of $\calE(X/\frakS)$ (resp. $\underline{\calE}(X/\frakS),$ resp. $\calE_{\mathrm{strict}}(X/\frakS),$ resp. $\underline{\calE}_{\mathrm{strict}}(X/\frakS)$), we denote by $\op{Cov}((U,\frakT,u))$ the collection of families of cartesian morphisms $((U_i\frakT_i,u_i) \ra (U,\frakT,u))$ such that $(U_i\ra U)$ is an étale covering of $U.$ By \ref{era3stablebc}, these collections define a pretopology on $\calE(X/\frakS)$ (resp. $\underline{\calE}(X/\frakS),$ $\calE_{\mathrm{strict}}(X/\frakS)$ and $\underline{\calE}_{\mathrm{strict}}(X/\frakS)$). We call the associated topology the \emph{étale topology of $\calE(X/\frakS)$} (resp. $\underline{\calE}(X/\frakS),$ resp. $\calE_{\mathrm{strict}}(X/\frakS),$ resp. $\underline{\calE}_{\mathrm{strict}}(X/\frakS)$) and we denote by $\widetilde{\calE}(X/\frakS)$ (resp. $\widetilde{\underline{\calE}}(X/\frakS),$ resp. $\widetilde{\calE}_{\mathrm{strict}}(X/\frakS),$ resp. $\widetilde{\underline{\calE}}_{\mathrm{strict}}(X/\frakS)$) the category of sheaves of sets on $\calE(X/\frakS)$ (resp. $\underline{\calE}(X/\frakS),$ resp. $\calE_{\mathrm{strict}}(X/\frakS),$ resp. $\underline{\calE}_{\mathrm{strict}}(X/\frakS)$).
\end{parag}

\begin{proposition}\label{era4prop1723}
Let $\F$ be a presheaf of sets on $\calE(X/\frakS)$ (resp. $\underline{\calE}(X/\frakS)$) and
$$(\F_{(U,\frakT,u)},\gamma_{\F,(f,g)})$$
the associated descent data \eqref{era3defdescentdata}. Then $\F$ is a sheaf for the étale topology \eqref{parettop} if and only if $\F_{(U,\frakT,u)}$ is a sheaf of $U_{\text{ét}}$ for every object $(U,\frakT,u)$ of $\calE(X/\frakS)$ (resp. $\underline{\calE}(X/\frakS)$).
\end{proposition}

\begin{proof}
Suppose that $\F$ is a sheaf for the étale topology and let $(U,\frakT,u)$ be an object of $\calE(X/\frakS)$ (resp. $\underline{\calE}(X/\frakS)$) and $(V_i\ra V)_{i\in I}$ an étale covering of $U.$ We have seen that the functor $\alpha_{(U,\frakT,u)}$ \eqref{eqKoko1861} sends morphisms to cartesian morphisms. It follows that
$$\alpha_{(U,\frakT,u)}(V_i) \ra \alpha_{(U,\frakT,u)}(V)$$
is a covering for the étale topology and so the sequence
$$0 \ra \F \circ \alpha_{(U,\frakT,u)}(V) \ra \prod_{i\in I}\F \circ \alpha_{(U,\frakT,u)}(V_i) \rightrightarrows \prod_{i,j\in I} \F \circ \alpha_{(U,\frakT,u)}(V_i\times_VV_j)$$
is exact. By \ref{propalphacomfiber}, we get the exactness of the sequence
$$0 \ra \F_{(U,\frakT,u)}(V) \ra \prod_{i\in I}\F_{(U,\frakT,u)}(V_i) \rightrightarrows \prod_{i,j\in I} \F_{(U,\frakT,u)}(V_i\times_VV_j),$$
and so $F_{(U,\frakT,u)}$ is a sheaf.

Conversely, suppose every $\F_{(U,\frakT,u)}$ is a sheaf and let $((U_i,\frakT_i,u_i) \xrightarrow{(f_i,g_i)} (U,\frakT,u))_{i\in I}$ be a covering for the étale topology and $(U_{ij},\frakT_{ij},u_{ij})=(U_{i},\frakT_{i},u_{i})\times_{(U,\frakT,u)}(U_{j},\frakT_{j},u_{j}).$ Since $\F_{(U,\frakT,u)}$ is a sheaf, we get the exact sequence
$$0 \ra \F_{(U,\frakT,u)}(U) \ra \prod_{i\in I} \F_{(U,\frakT,u)}(U_i) \rightrightarrows \prod_{i,j\in I} \F_{(U,\frakT,u)}(U_i\times_UU_j).$$
By \ref{era3propdescentdata} (iii), we get the exact sequence
$$0 \ra \F_{(U,\frakT,u)}(U) \ra \prod_{i\in I} \F_{(U_i,\frakT_i,u_i)}(U_i) \rightrightarrows \prod_{i,j\in I} \F_{(U_{ij},\frakT_{ij},u_{ij})}(U_i\times_UU_j),$$
which is equal to
$$0 \ra \F(U,\frakT,u) \ra \prod_{i\in I} \F(U_i,\frakT_i,u_i) \rightrightarrows \prod_{i,j\in I} \F(U_{ij},\frakT_{ij},u_{ij}).$$
\end{proof}

\begin{remark}
Proposition \ref{era4prop1723} remains true if we replace $\calE(X/\frakS)$ and $\underline{\calE}(X/\frakS)$ with $\calE_{\mathrm{strict}}(X/\frakS)$ and $\underline{\calE}_{\mathrm{strict}}(X/\frakS)$ respectively.
\end{remark}

\begin{proposition}\label{propalphacocont}
Let $(U,\frakT,u)$ be an object of $\calE(X/\frakS)$ (resp. $\underline{\calE}(X/\frakS)$) and $V$ an étale $U$-scheme. Consider the functor $\alpha_{(U,\frakT,u)}$ \eqref{eqKoko1861} and set $\left (V,\frakT_V,u_V\right )=\alpha_{(U,\frakT,u)}(V).$ For any covering family $\left ((U_i,\frakT_i,u_i) \xrightarrow{(f_i,g_i)} \left (V,\frakT_V,u_V \right ) \right )_{i\in I}$ for the étale topology \eqref{parettop}, there exists an étale covering $\left (U_i \ra V\right )_{i\in I},$ sent by $\alpha_{(U,\frakT,u)}$ to $(f_i,g_i).$
\end{proposition}

\begin{proof}
We first prove the proposition for $\calE(X/\frakS).$ For any $i\in I,$ the morphism $(U_i,\frakT_i,u_i) \ra \left (V,\frakT_V,u_V\right )$ is cartesian \eqref{era3defcart}. It follows that the left square in the diagram
$$
\begin{tikzcd}
T_i \ar{r} \ar[swap]{d}{u_i} & T_V \ar{r} \ar{d}{u_V} & T\ar{d}{u} \\
U_i \ar{r} & V \ar{r} & U
\end{tikzcd}
$$
is cartesian. The right square is cartesian by definition.
Set
$$
\left (U_i,\frakT_{U_i},u_{U_i}\right )=\alpha_{(U,\frakT,u)}(U_i).
$$
By definition, the square
$$
\begin{tikzcd}
T_{U_i} \ar{r} \ar[swap]{d}{u_{U_i}} & T \ar{d}{u} \\
U_i \ar{r} & U
\end{tikzcd}
$$
is cartesian. We deduce a canonical isomorphism
$$
\varphi_i:T_i \xrightarrow{\sim} T_{U_i}
$$
such that
$$
u_{U_i}\circ \varphi_i=u_i.
$$
By definition, we have cartesian squares
$$
\begin{tikzcd}
T_{U_i} \ar{r} \ar{d} & \frakT_{U_i} \ar{d} & T_i \ar{r} \ar{d} & \frakT_i \ar{d} \\
T \ar{r} & \frakT & T \ar{r} & \frakT
\end{tikzcd}
$$
where the left vertical arrow in each square is étale and strict. By \ref{era3paralphaKoko}, $\frakT_{U_i} \ra \frakT$ is étale and strict, hence the existence of a unique morphism $\widetilde{\varphi}_i:\frakT_i \ra \frakT_{U_i}$ fitting into the commutative diagram
$$
\begin{tikzcd}
T_{U_i} \ar{rr} & & \frakT_{U_i} \ar{d} \\
T_i \ar{u}{\varphi_i} \ar{r} & \frakT_i \ar{r} \ar{ur}{\widetilde{\varphi}_i} & \frakT.
\end{tikzcd}
$$
The morphism $\varphi$ is an isomorphism so, by (\cite{DXU19} 7.2), $
\widetilde{\varphi}_i:\frakT_i \ra \frakT_{U_i}
$
is also an isomorphism.
We have the commutative diagram
$$
\begin{tikzcd}
T_i \ar{rr} \ar{dr}{\varphi_i} \ar{dd} & & \frakT_i \ar[dashed]{dd} \ar{dr}{\widetilde{\varphi}_i} & \\
 & T_{U_i} \ar{rr}\ar{dd} & & \frakT_{U_i} \ar{dd} \\
T \ar[dashed]{rr} \ar[equal]{dr} & & \frakT \ar[dashed]{dr} & \\
 & T \ar{rr} &  & \frakT. 
\end{tikzcd}
$$
We deduce an isomorphism
$$
(U_i,\frakT_i,u_i) \xrightarrow{\sim} \left (U_i,\frakT_{U_i},u_{U_i}\right )=\alpha_{(U,\frakT,u)}(U_i)
$$
and the result follows. For the case $\underline{\calE}(X/\frakS),$ we have the following cartesian rectangles:
$$
\begin{tikzcd}
T_i \ar{r} \ar[swap]{d}{\Delta_{(U_i,\frakT_i,u_i)}} & T_V \ar{d}{\Delta_{(V,\frakT_V,u_V)}} \ar{rr} & & T \ar{d}{\Delta_{u}} \\
U_i' \ar{r} & V' \ar{rr} & & U'.
\end{tikzcd}
$$
We deduce an isomorphism
$$
T_i \xrightarrow{\sim} T_{U_i}.
$$
The rest of the proof is similar to the previous case.
\end{proof}

\begin{lemma}\label{Kokolemcocont}
Let $u:\calC \ra \calD$ be a functor between sites whose topologies are defined by pretopologies. Suppose that $u$ satisfies the following properties:
\begin{enumerate}
\item $u$ commutes with fiber products.
\item If $(X_i \ra X)_{i\in I}$ is a covering family for the pretopology of $\calC$ then $\left (u(X_i) \ra u(X) \right )_{i\in I}$ is a covering family in $\calD.$
\item If $\left (Y_i \xrightarrow{g_i} u(X) \right )_{i\in I}$ is a covering family in $\calD$ then there exists a covering family $(f_i:X_i \ra X)_{i\in I}$ in $\calC,$ such that $u(f_i)=g_i.$
\end{enumerate}
Then $u$ is continuous and cocontinuous and induces a morphism of topoi
$$
u:\widetilde{\calC} \ra \widetilde{\calD}
$$
such that $u^{-1}$ is the composition by $u.$
\end{lemma}

\begin{proof}
Since $u$ satisfies properties 1 and 2, its continuity results from (\cite{SGA43} III 1.6). The cocontinuity of $u$ results from property 3 and (\cite{SGA43} III 2.1 and II 1.4). The induced morphism of topoi results then from (\cite{SGA43} IV 4.7).
\end{proof}

\begin{proposition}
Let $(U,\frakT,u)$ be an object of $\calE(X/\frakS)$ (resp. $\underline{\calE}(X/S)$). The functor $\alpha_{(U,\frakT,u)}$ \eqref{eqKoko1861} induces a morphism of topoi
\begin{equation}\label{Kokoalphatopos}
\alpha_{(U,\frakT,u)}:U_{\text{ét}} \ra \widetilde{\calE}(X/\frakS)\ \left( \mathrm{resp. }\ \widetilde{\underline{\calE}}(X/\frakS)\right )
\end{equation}
such that the inverse image functor is composition by $\alpha_{(U,\frakT,u)}.$
\end{proposition}

\begin{proof}
It is sufficient to apply \ref{Kokolemcocont}. By \ref{era3paralphaKoko} (3), the functor $\alpha_{(U,\frakT,u)}$ sends covering families to covering families for the étale topology. Then, by \ref{propalphacomfiber} and \ref{propalphacocont}, the functor $\alpha_{(U,\frakT,u)}$ satisfies the properties of \ref{Kokolemcocont}. The result follows.
\end{proof}

\begin{parag}\label{era4logflattop}
We say that a family of morphisms $((U_i,\frakT_i,u_i) \ra (U,\frakT,u))_{i\in I}$ of $\calE(X/\frakS)$ (resp. $\underline{\calE}(X/\frakS)$) is a \emph{covering for the log flat topology} if $(U_i \ra U)_{i\in I}$ is an étale covering and, for all positive integers $n,$ $(\frakT_{i,n} \ra \frakT_n)_{i\in I}$ is a log flat covering (\ref{logflattop}). By \ref{era3fiberparag}, this defines a pretopology on $\calE(X/\frakS)$ (resp. $\underline{\calE}(X/\frakS)$). We call the corresponding topology the \emph{log flat topology} and we denote by $\widetilde{\calE}_{lf}(X/\frakS)$ (resp. $\widetilde{\underline{\calE}}_{lf}(X/\frakS)$) the corresponding topos.
\end{parag}

\begin{remark}
The log flat topology on $\calE(X/\frakS)$ (resp. $\underline{\calE}(X/\frakS)$) is finer then the étale topology (\ref{parettop}). Indeed, if $\left ((U_i,\frakT_i,u_i) \ra (U,\frakT,u) \right )_{i\in I}$ is an étale covering then $(\frakT_i \ra \frakT)_{i\in I}$ is an étale and strict covering.
\end{remark}

\subsection*{Crystals}

\begin{definition}
For all positive integers $n,$ we define a presheaf of rings $\Ox_{\calE(X/\frakS),n}$ (resp. $\Ox_{\underline{\calE}(X/\frakS),n}$) on $\calE(X/\frakS)$ (resp. $\underline{\calE}(X/\frakS)$) by
$$
\Ox_{\calE(X/\frakS),n}:(U,\frakT,u) \mapsto \Gamma(\frakT_n,\Ox_{\frakT_n})
$$
$$
(\op{resp.}\ \Ox_{\underline{\calE}(X/\frakS),n}: (U,\frakT,u)\mapsto \Gamma(\frakT_n,\Ox_{\frakT_n})).
$$
We denote by $\Ox_{\calE_{\mathrm{strict}}(X/\frakS),n}$ and $\Ox_{\underline{\calE}_{\mathrm{strict}}(X/\frakS),n}$ the restrictions of $\Ox_{\calE(X/\frakS),n}$ and $\Ox_{\underline{\calE}(X/\frakS),n}$ to $\calE_{\mathrm{strict}}(X/\frakS)$ and $\underline{\calE}_{\mathrm{strict}}(X/\frakS)$) respectively.
For $n=1,$ we denote $\Ox_{\calE(X/\frakS),1}$ (resp. $\Ox_{\underline{\calE}(X/\frakS),1},$ resp. $\Ox_{\calE_{\mathrm{strict}}(X/\frakS),1}$ resp. $\Ox_{\underline{\calE}_{\mathrm{strict}}(X/\frakS),1}$) simply by $\Ox_{\calE(X/\frakS)}$ (resp. $\Ox_{\underline{\calE}(X/\frakS)},$ resp. $\Ox_{\calE_{\mathrm{strict}}(X/\frakS)},$ resp. $\Ox_{\underline{\calE}_{\mathrm{strict}}(X/\frakS)}$).
\end{definition}

\begin{proposition}
For every positive integer $n,$ the presheaves $\Ox_{\calE(X/\frakS),n}$ and $\Ox_{\underline{\calE}(X/\frakS),n}$ are sheaves for the étale and log flat topologies.
\end{proposition}

\begin{proof}
Since the log flat topology if finer then the étale topology, it is sufficient to prove the result for the log flat topology. Let $n$ be a positive integer and 
$$((U_i,\frakT_i,u_i) \xrightarrow{(f_i,g_i)} (U,\frakT,u))_{i\in I}$$
a covering for the log flat topology. For any $i,j\in I,$ let $f_{ij}:\frakT_{ij}=\frakT_i\times_{\frakT}^{\op{log}}\frakT_j \ra \frakT$ be the canonical morphism. Since $(\frakT_{i,n} \xrightarrow{f_i}\frakT_n)_{i\in I}$ is a covering for the log flat topology (\ref{era4logflattop}), the sequence
$$0 \ra \Ox_{\frakT_n} \ra \prod_{i\in I}f_{i*}f_i^*\Ox_{\frakT_n} \rightrightarrows \prod_{i,j\in I}f_{ij*}f_{ij}^*\Ox_{\frakT_n}$$
is exact by \ref{eraflogflatdescent}. We deduce the exactness of the sequence
$$0 \ra \Gamma \left (\frakT_n,\Ox_{\frakT_n} \right ) \ra \prod_{i\in I}\Gamma \left (\frakT_{i,n},\Ox_{\frakT_{i,n}} \right ) \rightrightarrows \prod_{i,j\in I} \Gamma\left (\frakT_{ij,n},\Ox_{\frakT_{ij,n}}\right ).$$
\end{proof}

\begin{parag}\label{lindescentdata}
Let $n$ be a positive integer, $(U,\frakT,u)$ an object of $\calE(X/\frakS)$ (resp. $\underline{\calE}(X/\frakS)$) and $V$ an étale $U$-scheme. Let
$$(V,\frakT_V,u_V)=\alpha_{(U,\frakT,u)}(V).$$
Then
\begin{alignat*}{2}
\Gamma \left (V,\left (\Ox_{\calE(X/\frakS),n} \right )_{(U,\frakT,u)} \right ) = \Gamma \left (\frakT_{V,n},\Ox_{\frakT_{V,n}} \right ) 
= \Gamma \left (V, u_{V*}\Ox_{\frakT_{V,n}} \right ).
\end{alignat*}
It follows that
\begin{equation}\label{eqKoko18171}
\left (\Ox_{\calE(X/\frakS),n} \right )_{(U,\frakT,u)}=u_*\Ox_{\frakT_n}.
\end{equation}
The same is true for $\underline{\calE}(X/\frakS).$

Let $\F$ be a module of $\left (\widetilde{\calE}(X/\frakS),\Ox_{\calE(X/\frakS)}\right )$ or $\left (\widetilde{\underline{\calE}}(X/\frakS),\Ox_{\underline{\calE}(X/\frakS)}\right )$ and $\left (\F_{(U,\frakT,u)},\gamma_{\F,(f,g)} \right )$ the associated descent data (\ref{era3propdescentdata}). A morphism $(f,g):(U_1,\frakT_1,u_1) \ra (U_2,\frakT_2,u_2)$ induces a morphism of ringed topoi
\begin{equation}\label{Kokoftilda}
\widetilde{f}:\left (U_{1,\text{ét}},u_{1*}\Ox_{T_{1}} \right ) \ra \left ( U_{2,\text{ét}}, u_{2*}\Ox_{T_{2}} \right ).
\end{equation}
Note that this morphism $\widetilde{f}$ does not depend on $f:\frakT_1 \ra \frakT_2$ but rather on its special fiber.
The morphism
$$\gamma_{\F,(f,g)}:g^{-1}\F_{(U_2,\frakT_2,u_2)} \ra \F_{(U_1,\frakT_1,u_1)}$$
then induces a $u_{1*}\Ox_{T_1}$-linear morphism
\begin{equation}\label{cflindata}
c_{\F,(f,g)}: \widetilde{f}^*\F_{(U_2,\frakT_2,u_2)} \ra \F_{(U_1,\frakT_1,u_1)}.
\end{equation}
By \ref{era3propdescentdata}, if $(f,g)$ is cartesian, then
$$
\gamma_{\Ox_{\calE(X/\frakS)},(f,g)}:g^{-1}\Ox_{\calE(X/\frakS),(U_2,\frakT_2,u_2)} \ra \Ox_{\calE(X/\frakS),(U_1,\frakT_1,u_1)}
$$
is an isomorphism. This means that the morphism
$$
\gamma_{\Ox_{\calE(X/\frakS)},(f,g)}:g^{-1}u_{2*}\Ox_{U_2} \ra u_{1*}\Ox_{U_1}
$$
is an isomorphism. It follows that, if $(f,g)$ is cartesian,
$$
\widetilde{f}^{-1}=\widetilde{f}^*.
$$
By \ref{era3propdescentdata} and \ref{era4prop1723}, an $\Ox_{\calE(X/\frakS)}$-module (resp. $\Ox_{\underline{\calE}(X/\frakS)}$-module) $\F$ is equivalent to the data:
\begin{enumerate}
\item For every object $(U,\frakT,u)$ of $\calE(X/\frakS)$ (resp. $\underline{\calE}(X/\frakS)$), a module $\F_{(U,\frakT,u)}$ of $(U_{\text{ét}},u_*\Ox_T).$
\item For every morphism $(f,g):(U_1,\frakT_1,u_1) \ra (U_2,\frakT_2,u_2),$ a $u_{1*}\Ox_{T_1}$-linear morphism
$$c_{\F,f,g)}:\widetilde{f}^*\F_{(U_2,\frakT_2,u_2)} \ra \F_{(U_1,\frakT_1,u_1)},$$
\end{enumerate}
satisfying the following conditions:
\begin{enumerate}[(i)]
\item If $(f,g)$ is the identity then so is $c_{\F,(f,g)}.$
\item For all composable morphisms $(f_1,g_1):(U_1,\frakT_1,u_1) \ra (U_2,\frakT_2,u_2)$ and $(f_2,g_2):(U_2,\frakT_2,u_2) \ra (U_3,\frakT_3,u_3)$ in $\calE(X/\frakS)$ (resp. $\underline{\calE}(X/\frakS)$),
$$c_{\F,(f_1,g_1)} \circ \widetilde{f_1}^*c_{\F,(f_2,g_2)}=c_{\F,(f_2,g_2)\circ (f_1,g_1)}.$$
\item If $(f,g)$ is a cartesian morphism, then $c_{\F,(f,g)}$ is an isomorphism.
\end{enumerate}
This equivalence remains true if we replace $\calE(X/\frakS)$ and $\underline{\calE}(X/\frakS)$ with $\calE_{\mathrm{strict}}(X/\frakS)$ and $\underline{\calE}_{\mathrm{strict}}(X/\frakS)$ respectively.
We call $\left (\F_{(U,\frakT,u)},c_{(f,g)} \right )$ the \emph{linearized descent data associated with $\F$}.
\end{parag}

\begin{definition}\label{defcrys}
Let $\F$ be a module of $\widetilde{\calE}(X/\frakS),$ $\widetilde{\underline{\calE}}(X/\frakS),$
$\widetilde{\calE}_{lf}(X/\frakS),$ $\widetilde{\underline{\calE}}_{lf}(X/\frakS),$ $\widetilde{\calE}_{\mathrm{strict}}(X/\frakS)$ or $\widetilde{\underline{\calE}}_{\mathrm{strict}}(X/\frakS)$ and $\left (\F_{(U,\frakT,u)},c_{\F,(f,g)} \right )$ the associated linearized descent data (\ref{lindescentdata}).
\begin{enumerate}
\item We say that $\F$ is \emph{quasi-coherent} if $\F_{(U,\frakT,u)}$ is a quasi-coherent $u_*\Ox_{T}$-module for every object $(U,\frakT,u).$
\item We say that $\F$ is a crystal if $c_{\F,(f,g)}$ is an isomorphism for every morphism $(f,g).$
\end{enumerate}
We denote by $\calC(X/\frakS),$ $\underline{\calC}(X/\frakS),$ $\calC_{lf}(X/\frakS),$ $\underline{\calC}_{lf}(X/\frakS),$ $\calC_{\mathrm{strict}}(X/\frakS)$ and $\underline{\calC}_{\mathrm{strict}}(X/\frakS)$ the full subcategories of crystals of the categories of modules of $\widetilde{\calE}(X/\frakS),$ $\widetilde{\underline{\calE}}(X/\frakS),$
$\widetilde{\calE}_{lf}(X/\frakS),$ $\widetilde{\underline{\calE}}_{lf}(X/\frakS),$ $\widetilde{\calE}_{\mathrm{strict}}(X/\frakS)$ or $\widetilde{\underline{\calE}}_{\mathrm{strict}}(X/\frakS)$ respectively.
We also denote by $\calC^{\text{qcoh}}(X/\frakS),$ $\underline{\calC}^{\text{qcoh}}(X/\frakS),$ $\calC_{lf}^{\text{qcoh}}(X/\frakS),$ $\underline{\calC}_{lf}^{\text{qcoh}}(X/\frakS),$ $\calC_{\mathrm{strict}}^{\text{qcoh}}(X/\frakS)$ and $\underline{\calC}_{\mathrm{strict}}^{\text{qcoh}}(X/\frakS)$ the full subcategories of quasi-coherent crystals respectively.
\end{definition}

\begin{proposition}\label{equivcrysShiho}
The restriction functor
$$
R:\widetilde{\calE}(X/\frakS) \ra \widetilde{\calE}_{\mathrm{strict}}(X/\frakS) \ \left (\text{resp.}\ \underline{\widetilde{\calE}}(X/\frakS) \ra \underline{\widetilde{\calE}}_{\mathrm{strict}}(X/\frakS) \right )
$$
induces an equivalence of categories
$$
R_{\mathrm{crys}}:\begin{Bmatrix}\mathrm{Crystals\ of\ }\widetilde{\calE}(X/\frakS) \end{Bmatrix} \xrightarrow{\sim} \begin{Bmatrix}\mathrm{Crystals\ of\ }\widetilde{\calE}_{\mathrm{strict}}(X/\frakS) \end{Bmatrix}$$
$$\left ( \mathrm{resp.}\ \begin{Bmatrix}\mathrm{Crystals\ of\ }\widetilde{\underline{\calE}}(X/\frakS) \end{Bmatrix} \xrightarrow{\sim} \begin{Bmatrix}\mathrm{Crystals\ of\ }\widetilde{\underline{\calE}}_{\mathrm{strict}}(X/\frakS) \end{Bmatrix} \right ).
$$
\end{proposition}

\begin{proof}
We just prove the result for $\calE(X/\frakS)$ as the other case is similar. Let $f_1,f_2:\F \ra \calG$ be morphisms of crystals of $\widetilde{\calE}(X/\frakS)$ such that $R(f_1)=R(f_2).$ We want to prove that $f_1=f_2.$ Let $(U,\frakT,u)$ be an object of $\calE(X/\frakS).$ To prove that $f_{1,(U,\frakT,u)}=f_{2,(U,\frakT,u)},$ we can suppose that $U$ is affine. It follows that $T$ is also affine. Since $U\ra S$ is log smooth, there exists a log smooth logarithmic formal scheme $\frakU$ over $\frakS,$ fitting into the cartesian square
$$
\begin{tikzcd}
U \ar{r} \ar{d} & \frakU \ar{d} \\
S\ar{r} & \frakS.
\end{tikzcd}
$$
There exists then a morphism $\varphi:\frakT \ra \frakU$ fitting into the commutative diagram
$$
\begin{tikzcd}
U\ar{rr} & & \frakU \ar{d} \\
T\ar{r} \ar{u} & \frakT \ar{r} \ar{ur}{\varphi} & \frakS.
\end{tikzcd}
$$
We then get a morphism
$$
\left ( \op{Id}_U,\varphi \right ):(U,\frakT,u) \ra (U,\frakU,\op{Id}_U)
$$
in $\calE(X/\frakS)$ and we obtain a commutative diagram
$$
\begin{tikzcd}
\widetilde{\varphi}^*\F_{(U,\frakU,\op{Id}_U)} \ar{r} \ar[swap]{d}{\widetilde{\varphi}^*f_{i,(U,\frakU,\op{Id}_U)}} & \F_{(U,\frakT,u)} \ar{d}{f_{i,(U,\frakT,u)}} \\
\widetilde{\varphi}^*\calG_{(U,\frakU,\op{Id}_U)} \ar{r} & \calG_{(U,\frakT,u)},
\end{tikzcd}
$$
where $\widetilde{\varphi}$ is given in \eqref{Kokoftilda}. Since $\F$ and $\calG$ are crystals of $\widetilde{\calE}(X/\frakS),$ the horizontal morphisms are isomorphisms. Since $R(f_1)=R(f_2)$ and $(U,\frakU,\op{Id}_U)$ is an object of $\calE_{\mathrm{strict}}(X/\frakS),$ we have $f_{1,(U,\frakU,\op{Id}_U)}=f_{2,(U,\frakU,\op{Id}_U)}.$ It follows that
$$
f_{1,(U,\frakT,u)}=f_{2,(U,\frakT,u)}.
$$
We deduce that $f_1=f_2$ and $R_{\mathrm{crys}}$ is faithful.

Now let $\F$ and $\calG$ be crystals of $\widetilde{\calE}(X/\frakS)$ and $g:R(\F) \ra R(\calG)$ a morphism. We want to prove that there exists a morphism $f:\F \ra \calG$ in $\widetilde{\calE}(X/\frakS)$ such that $R(f)=g.$ Let $(U,\frakT,u)$ and $\varphi:\frakT \ra \frakU$ be as defined above. Just as above, we suppose that $U$ is affine. Note that the morphism of topoi
$$
\widetilde{\varphi}:\left (U_{\text{ét}},u_*\Ox_T \right ) \ra \left (U_{\text{ét}},\Ox_U \right )
$$
depends only on $u:T\ra U.$
We set
$$
g_{(U,\frakU,\op{Id}_U)}=f_{(U,\frakU,\op{Id}_U)},\ g_{(U,\frakT,u)}=\widetilde{\varphi}^*g_{(U,\frakU,\op{Id}_U)}.
$$
This is independant of the lifting $\varphi.$ We also prove that it is independant of the lifting $\frakU.$ Let $\frakU_1$ and $\frakU_2$ be log smooth logarithmic formal schemes over $\frakS$ fitting into cartesian squares
$$
\begin{tikzcd}
U \ar{r} \ar{d} & \frakU_i \ar{d} \\
S\ar{r} & \frakS.
\end{tikzcd}
$$
There exists morphisms $\varphi_1:\frakT\ra \frakU_1$ and $\varphi_{12}:\frakU_1 \ra \frakU_2$ fitting into the commutative diagram
$$
\begin{tikzcd}
 & & \frakU_2 \ar{d} \\
U \ar{r} \ar[bend right=-30]{urr} & \frakU_1 \ar{r}  \ar{ur}{\varphi_{12}} & \frakS \\
T \ar{u} \ar{r} & \frakT \ar{u}{\varphi_1} \ar{ur} & 
\end{tikzcd}
$$
Set $\varphi_2 =\varphi_{12} \circ \varphi_1.$
We hence obtain morphisms
$$
\left ( \op{Id}_U,\varphi_2 \right ):(U,\frakT,u) \xrightarrow{(\op{Id}_U,\varphi_1)} (U,\frakU_1,\op{Id}_U) \xrightarrow{(\op{Id}_U,\varphi_{12})} (U,\frakU_2,\op{Id}_U).
$$
Since $\F$ and $\calG$ are crystals of $\widetilde{\calE}(X/\frakS),$ we get a canonical isomorphism
$$
\widetilde{\varphi}_{12}^*g_{(U,\frakU_2,\op{Id}_U)} \xrightarrow{\sim} g_{(U,\frakU_1,\op{Id}_U)}
$$
and then, by applying $\widetilde{\varphi}_1^*,$ a canonical isomorphism
\begin{equation}\label{Shihoeq10}
\widetilde{\varphi}_2^*g_{(U,\frakU_2,\op{Id}_U)} \xrightarrow{\sim} \widetilde{\varphi}_1^*g_{(U,\frakU_1,\op{Id}_U)}.
\end{equation}
Now let $(U,\frakT,u)$ an object of $\calE(X/\frakS)$ where $U$ is not necessarily affine. Let $U=\bigcup_{i\in I}U_i$ an affine open cover of $U.$ By \ref{era3paralphaKoko} (1), we have an étale covering $\left (\alpha_{(U,\frakT,u)}\left (U_i \right ) \ra (U,\frakT,u) \right )_{i\in I}.$ The isomorphisms \eqref{Shihoeq10} satisfy the cocycle condition and hence prove that the morphisms $g_{\alpha_{(U,\frakT,u)}(U_i)}$ glue together into a morphism $g_{(U,\frakT,u)}.$ We then obtain the desired morphism $g:\F\ra \calG.$ The functor $R_{\mathrm{crys}}$ is hence full.

To prove the essential surjectivity, we consider a crystal $\calG$ of $\widetilde{\calE}_{\mathrm{strict}}(X/\frakS).$ For an object $(U,\frakT,u)$ of $\calE(X/\frakS)$ such that $U$ is affine, we consider $\varphi:\frakT \ra \frakU$ as defined above and set
$$
\F_{(U,\frakT,u)}=\widetilde{\varphi}^*\calG_{(U,\frakU,\op{Id}_U)}.
$$
By an argument similar to the one used for fullness, we define $\F_{(U,\frakT,u)}$ for any object $(U,\frakT,u)$ of $\calE(X/\frakS)$ and then obtain a crystal $\F$ of $\widetilde{\calE}(X/\frakS)$ extending $\calG.$
\end{proof}

\begin{parag}\label{equivRQDef}
Suppose that $X\ra S$ lifts to a log smooth morphism of framed logarithmic formal schemes $(\frakX,Q) \ra (\frakS,P).$ Consider the logarithmic formal schemes $R_{\frakX,1}$ and $Q_{\frakX}$ defined in \ref{parag86}. 
Let $R_1$ and $Q_1$ be the special fibers of $R_{\frakX,1}$ and $Q_{\frakX}$ respectively, $i_X:\underline{X} \ra X$ the canonical immersion, $q_1,q_2:R_{\frakX,1} \ra \frakX$ and $q_1',q_2':Q_{\frakX} \ra \frakX$ the canonical projections and $\lambda_R:R_1 \ra X$ and $\lambda_Q:\underline{Q_1} \ra X$ the morphisms given in \ref{erafprop13}.
By \ref{Qlogflat}, $(X,\frakX,\op{Id}_X)$ and $(X,R_{\frakX,1},\lambda_R)$ are objects of $\calE(X/\frakS)$ (resp. $(X,\frakX,i_X)$ and $(X,Q_{\frakX},\lambda_Q)$ are objects of $\underline{\calE}(X/\frakS)$).
The projection $q_i':Q_{\frakX} \ra \frakX$ induces by reducing modulo $p$ and functoriality a morphism $\underline{Q_1} \ra \underline{X}.$ It fits into the commutative diagram
$$
\begin{tikzcd}
 & & \underline{X} \ar[bend right=-30]{dddl}{i_X} \\
\underline{Q_1} \ar{r} \ar[bend right=-30]{urr} \ar{d}{\lambda_Q} & Q_1 \ar{d} & \\
X \ar{r} \ar[equal]{dr} & X\times_{S,[Q]}^{\op{log}}X \ar{d} & \\
& X &,
\end{tikzcd}
$$
where $X\ra X\times_{S,[Q]}^{\op{log}}X$ is the strict diagonal immersion and $X\times_{S,[Q]}^{\op{log}}X \ra X$ is the $i$th projection. We hence have a morphism
$$
(\op{Id}_X,q_i'):(X,Q_{\frakX},\lambda_Q) \ra (X,\frakX,i_X) \quad \text{in}\quad \underline{\calE}(X/\frakS).
$$
Let $\mathfrak{E}$ be a crystal of $\calC(X/\frakS)$ (resp. $\underline{\calC}(X/\frakS)$). We set $\calE=\mathfrak{E}_{(X,\frakX,\op{Id}_X)}$ (resp. $\calE=\mathfrak{E}_{(X,\frakX,i_X)}$). Since $\mathfrak{E}$ is a crystal, the morphisms $(\op{Id_X},q_i):(X,R_{\frakX,1},\lambda_R) \ra (X,\frakX,\op{Id}_X)$ (resp. $(\op{Id_X},q_i'):(X,Q_{\frakX},\lambda_Q) \ra (X,\frakX,i_X)$) induce isomorphisms of $\Ox_{R_1}$-modules (resp. $\Ox_{Q_1}$-modules)
$$c_{(\op{Id}_X,q_i)}:q_i^*\calE \xrightarrow{\sim} \mathfrak{E}_{(X,R_{\frakX,1},\lambda_R)}$$
$$
\left (\text{resp.}\ c_{(\op{Id}_X,q_i')}:q_i'^*\calE \xrightarrow{\sim} \mathfrak{E}_{(X,Q_{\frakX},\lambda_Q)}\right ).
$$
Let $\epsilon$ be the composition
$$\epsilon = \left (c_{\op{Id}_X,q_1} \right )^{-1} \circ c_{\op{Id}_X,q_2}:q_2^*\calE \xrightarrow{\sim} \mathfrak{E}_{(X,R_{\frakX,1},\lambda_R)} \xrightarrow{\sim} q_1^*\calE$$
$$
\left (\text{resp.}\ \epsilon = \left (c_{\op{Id}_X,q_1'} \right )^{-1} \circ c_{\op{Id}_X,q_2'}:q_2'^*\calE \xrightarrow{\sim} \mathfrak{E}_{(X,Q_{\frakX},\lambda_Q)} \xrightarrow{\sim} q_1'^*\calE \right ).
$$
Then $\epsilon$ is an $R_{\frakX,1}$-stratification (resp. $Q_{\frakX}$-stratification) on $\calE$ and so we obtain a functor
\begin{equation}\label{funct1} 
\begin{array}[t]{clc c}
\calC(X/\frakS) \left ( \text{resp.\ } \underline{\calC}(X/\frakS) \right )& \ra & \begin{Bmatrix}\Ox_{X}\text{-}\mathrm{modules\ with\ an}\\R_{\frakX,1}\text{-}\mathrm{stratification} \end{Bmatrix} & \left (\mathrm{resp.} \begin{Bmatrix}\Ox_{X}\text{-}\mathrm{modules\ with\ an}\\Q_{\frakX}\text{-}\mathrm{stratification} \end{Bmatrix} \right )\\
\mathfrak{E} & \mapsto & (\calE,\epsilon). &
\end{array}
\end{equation}
\end{parag}

\begin{proposition}\label{equivRQ}
Suppose that $X\ra S$ lifts to a log smooth morphism of framed logarithmic formal schemes $(\frakX,Q) \ra (\frakS,P).$ Consider the logarithmic formal schemes $R_{\frakX,1}$ and $Q_{\frakX}$ defined in \ref{parag86}. The functor \eqref{funct1} is an equivalence of categories.
\end{proposition}

\begin{proof}
This is similar to (\cite{DXU19} 8.10) so we just outline the proof for $\calC(X/\frakS)$ (the case $\underline{\calC}(X/\frakS)$ is similar).
Let $(\calE,\epsilon)$ be an $\Ox_X$-module equipped with an $R_{\frakX,1}$-stratification. Set
$$\calE = \mathfrak{E}_{(X,\frakX,\op{Id}_X)}.$$
Let $(U,\frakT,u)$ be an object of $\calE(X/\frakS)$ such that $U$ is affine. Since $u:T\ra U$ is by definition affine, so is $T.$ Since $\frakX \ra \frakS$ is log smooth and $T$ is affine, there exists a morphism $\varphi:\frakT \ra \frakX$ fitting into the commutative diagram
$$
\begin{tikzcd}
U \ar{rr} & & \frakX\ar{d} \\
T\ar{u}{u} \ar{r} & \frakT\ar{ur}{\varphi} \ar{r} & \frakS.
\end{tikzcd}
$$
We obtain a morphism $\varphi:(U,\frakT,u) \ra (X,\frakX,\op{Id}_X).$ We set
$$\mathfrak{E}_{(U,\frakT,u)}=\widetilde{\varphi}^*\mathfrak{E}_{(X,\frakX,\op{Id}_X)},$$
where $\widetilde{\varphi}$ is given in \eqref{Kokoftilda}.
This definition of $\mathfrak{E}_{(U,\frakT,u)}$ is independant of the choice of the lifting $\varphi$ up to a canonical isomorphism induced by the stratification $\epsilon.$

Let $(f,g):(U_1,\frakT_1,u_1) \ra (U_2,\frakT_2,u_2)$ be a morphism of $\calE(X/\frakS)$ such that $U_1$ and $U_2$ are affine. We have to define a $u_{1*}\Ox_{T_1}$-linear morphism
$$
\widetilde{f}^*\frakE_{(U_2,\frakT_2,u_2)} \ra \frakE_{(U_1,\frakT_1,u_1)}.
$$
By the argument above, there exists a morphism $\varphi_2$ fitting into the commutative diagram
$$
\begin{tikzcd}
U_2 \ar{rr} & & \frakX\ar{d} \\
T_2\ar{u}{u_2} \ar{r} & \frakT_2\ar{ur}{\varphi_2} \ar{r} & \frakS.
\end{tikzcd}
$$
Let $\varphi_1$ be the composition
$$\varphi_1:\frakT_1 \xrightarrow{f} \frakT_2 \xrightarrow{\varphi_2} \frakX.$$
We set
$$c_{\mathfrak{E},(f,g)}:\widetilde{f}^*\mathfrak{E}_{(U_2,\frakT_2,U_2)}=\widetilde{f}^*\widetilde{\varphi_2}^*\calE \xrightarrow{\sim}  \mathfrak{E}_{(U_1,\frakT_1,u_1)},$$
where the isomorphism is the canonical one. The isomorphism $c_{\mathfrak{E},(f,g)}$ is independant of the choices so 
we can glue together the morphisms $c_{\mathfrak{E},(f,g)}$ for $U_1$ and $U_2$ affine. The cocycle condition of $\epsilon$ implies that of the morphisms $c_{\mathfrak{E},(f,g)}.$ We thus obtain linearized descent datum $\left (\mathfrak{E}_{(U,\frakT,u)},c_{\mathfrak{E},(f,g)}\right )$ and then a crystal $\mathfrak{E}$ of $\widetilde{\calE}(X/\frakS).$
The correspondance $(\calE,\epsilon) \mapsto \mathfrak{E}$ is functorial and quasi-inverse to the one given in \ref{equivRQDef}.
\end{proof}

\begin{parag}
We have a functor
\begin{equation}\label{era3rhoShiho}
\rho:\begin{array}[t]{clc}
\underline{\calE}(X/\frakS) & \ra & \calE(X'/\frakS) \\
(U,\frakT,u) & \mapsto & \left ( U',\frakT,\Delta_{u} \right ),
\end{array}
\end{equation}
where $\Delta_u$ is defined in \ref{propShiho185}.
We will show in \ref{thmproof2} that, if $X\ra S$ is of Cartier type (or equivalently saturated, by \cite{Ogus2018} III 2.5.4), the functor $\rho$ is continuous and cocontinuous for the étale and log flat topologies. It follows, by (\cite{SGA43} IV 4.7), that $\rho$ defines two morphisms of topoi
\begin{alignat}{2}
C_{X/\frakS}&:\widetilde{\underline{\calE}}(X/\frakS) \ra \widetilde{\calE}(X'/\frakS), \label{era3C} \\
C_{X/\frakS,lf}&:\widetilde{\underline{\calE}}_{lf}(X/\frakS) \ra \widetilde{\calE}_{lf}(X'/\frakS). \label{era3Clf}
\end{alignat}
For any object $(U,\frakT,u)$ of $\underline{\calE}(X/\frakS),$ we have
\begin{alignat*}{2}
\left (C_{X/\frakS}^{-1}\Ox_{\calE(X'/\frakS)} \right )(U,\frakT,u) &= \Ox_{\calE(X'/\frakS)} \circ \rho (U,\frakT,u) \\
&= \Ox_{\calE(X'/\frakS)}(U',\frakT,\Delta_u) \\
&= \Gamma(T,\Ox_T) \\
&= \Ox_{\underline{\calE}(X/\frakS)}(U,\frakT,u).
\end{alignat*}
Then
$$
C_{X/\frakS}^{-1}\Ox_{\calE(X'/\frakS)}=\Ox_{\underline{\calE}(X/\frakS)}.
$$
Similarly, we have
$$
C_{X/\frakS,lf}^{-1}\Ox_{\calE(X'/\frakS)}=\Ox_{\underline{\calE}(X/\frakS)}.
$$
It follows that $C_{X/\frakS}$ and $C_{X/\frakS,lf}$ are morphisms of ringed topoi.
\end{parag}

\begin{lemma}\label{lemfiberprodcom}
Suppose that $X\ra S$ is of Cartier type. The functor $\rho:\underline{\calE}(X/\frakS) \ra \calE(X'/\frakS)$ defined in \eqref{era3rho} commutes with fiber products.
\end{lemma}

\begin{proof}
Let $(f_i,g_i):(U_i,\frakT_i,u_i) \ra (U,\frakT,u),$ $i=1,2,$ be two morphisms of $\underline{\calE}(X/\frakS).$ Since $X\ra S$ is of Cartier type,
$$
U_1'=U_1'',\ U_2'=U_2'',\ U'=U''.
$$
By definition,
$$(U_1,\frakT_1,u_1) \times_{(U,\frakT,u)}(U_2,\frakT_2,u_2)=(U_1\times_UU_2,\frakT_1\times_{\frakT}^{\op{log}}\frakT_2,v),$$
where $v$ is the composition
$$v:\underline{T_1\times_T^{\op{log}}T_2} \ra \underline{T_1}\times_{\underline{T}}^{\op{log}}\underline{T_2} \ra U_1\times_UU_2$$
of the canonical morphism $\underline{T_1\times_T^{\op{log}}T_2} \ra \underline{T_1}\times_{\underline{T}}^{\op{log}}\underline{T_2}$ (\ref{era3functoriality}) and the morphism induced by $u_1$ and $u_2.$ It follows that
$$\rho((U_1,\frakT_1,u_1) \times_{(U,\frakT,u)}(U_2,\frakT_2,u_2))=(U_1'\times_{U'}U_2',\frakT_1\times_{\frakT}^{\op{log}}\frakT_2,\Delta_v).$$
On the other hand,
$$\rho(U_1,\frakT_1,u_1) \times_{\rho(U,\frakT,u)}\rho(U_2,\frakT_2,u_2)=(U_1'\times_{U'}U_2',\frakT_1\times_{\frakT}^{\op{log}}\frakT_2,w),$$
where $w:T_1\times_T^{\op{log}}T_2 \ra U_1'\times_{U'}U_2'$ is induced by $\Delta_{u_1}:T_1 \ra U_1',$ $\Delta_{u_2}:T_2 \ra U_2'$ and $\Delta_{u}:T \ra U'.$
The composition
$$
T_1\times_T^{\op{log}}T_2 \xrightarrow{G_{T_1\times_T^{\op{log}}T_2}}  \underline{T_1\times_T^{\op{log}}T_2} \xrightarrow{v} U_1\times_UU_2,
$$
where $G_{T_1\times_T^{\op{log}}T_2}$ is defined in \ref{era4prop15},
is equal to the composition
$$
T_1\times_T^{\op{log}}T_2 \ra \underline{T_1}\times_{\underline{T}}^{\op{log}}\underline{T_2} \ra U_1\times_UU_2
$$
of the morphism induced by $G_{T_i}$ and $G_T$ and the morphism induced by $u_i$ and $u.$
The diagram
$$
\begin{tikzcd}
T_1\times_T^{\op{log}}T_2 \ar{r}{G_{T_1\times_T^{\op{log}}T_2}} \ar{dr}{w} \ar[bend right=30]{ddr} & \underline{T_1\times_T^{\op{log}}T_2} \ar{dr}{v} & \\
 & U_1'\times_{U'}U_2' \ar{r} \ar{d} & U_1\times_UU_2 \ar{d} \\
& S \ar{r}{F_S} & S
\end{tikzcd}
$$
is then commutative. It follows, by definition of $\Delta_v,$ that $\Delta_v=w.$
\end{proof}

\begin{lemma}\label{lemcoco}
Suppose that $X\ra S$ is of Cartier type.
Let $(U,\frakT,u)$ be an object of $\underline{\calE}(X/\frakS)$ and $g:(W,\frakZ,w) \ra \rho(U,\frakT,u)$ a morphism of $\calE(X'/\frakS).$ There exists a morphism $f:(V,\frakZ,v) \ra (U,\frakT,u)$ of $\underline{\calE}(X/\frakS)$ such that $\rho(f)=g.$ In addition, if $g$ is cartesian (resp. log flat), then so is $f.$
\end{lemma}

\begin{proof}
We have a commutative diagram
$$
\begin{tikzcd}
Z \ar{r} \ar[swap]{d}{w}  & T \ar{d}{\Delta_u} \\
W \ar{r} & U'.
\end{tikzcd}
$$
By \ref{piunivhomeo}, the morphism $\pi:X' \ra X$ \eqref{diag51} is a universal homeomorphism. It follows that the functor
$$
\text{ét}_{/X} \ra \text{ét}_{/X'},\ Y\mapsto Y\times_XX',
$$
induced by $\pi,$ is an equivalence of categories (\href{https://stacks.math.columbia.edu/tag/04DZ}{Theorem 04DZ}). Let $V$ be the unique étale $X$-scheme such that $V\times_XX'=W.$ We then have cartesian squares
$$
\begin{tikzcd}
W \ar{r} \ar{d} & V \ar{d} \\
X' \ar{r}{\pi} \ar{d} & X \ar{d} \\
S \ar{r}{F_S} & S.
\end{tikzcd}
$$ 
It follows that $W=V'.$ We have an $X'$-morphism $W=V' \ra U'$ and there exists a unique $X$-morphism $V \ra U$ whose base change by $\pi$ is $V' \ra U'.$
Consider the diagram
\begin{equation}\label{diagSH1}
\begin{tikzcd}
V' \ar{rrr} & & & V \ar{d} \\
Z \ar{u}{w} \ar{r}{G_Z} & \underline{Z} \ar[dashed]{urr}{v} \ar{r} & \underline{T} \ar{r}{u} & U,
\end{tikzcd}
\end{equation}
where $G_Z$ is defined in \ref{era4prop15}. The composition
$$
Z \xrightarrow{G_Z} \underline{Z} \ra \underline{T} \xrightarrow{u} U
$$
is equal to
$$
Z \ra T \xrightarrow{G_T} \underline{T} \xrightarrow{u} U.
$$
By definition of $\Delta_u$ \eqref{propShiho185} and \eqref{morphismShiho2}, this is equal to
$$
Z \ra T \xrightarrow{\Delta_u} U' \ra U.
$$
This is the same as
$$
Z \ra V' \ra U' \ra U.
$$
This is equal to
$$
Z \ra V' \ra V \ra U
$$
and the diagram \eqref{diagSH1} is commutative. Since $G_Z$ is inseparable, the dashed morphism $v:\underline{Z} \ra V$ exists and is unique by (\cite{Ogus2018} IV 3.3.7). In addition, $v$ is affine since $V' \ra V$ and $w$ are affine. Hence we have an object $(V,\frakZ,v)$ of $\underline{\calE}(X/\frakS)$ and a morphism $g:(V,\frakZ,v) \ra (U,\frakT,u).$ It remains to prove that
$$
\rho(V,\frakZ,v)=(W,\frakZ,w).
$$
This is equivalent to proving that $w=\Delta_v.$
Since $X\ra S$ is of Cartier type, $V'=V''.$ The morphism $\Delta_v$ is hence the unique morphism fitting into the diagram \eqref{morphismShiho2}, which we recall here:
$$
\begin{tikzcd}
Z \ar{r}{G_Z} \ar{dr}{\Delta_v} \ar[bend right=30]{ddr} & \underline{Z} \ar{dr}{v} & \\
 & V' \ar{r} \ar{d} & V \ar{d} \\
 & S \ar{r}{F_S} & S.
\end{tikzcd}
$$
The morphism $w:Z\ra V'$ fits in the same commutative diagram by \eqref{diagSH1} and the fact that $w$ is an $S$-morphism. It is clear that $g$ is cartesian (resp. log flat) if and only if $f$ is.
\end{proof}

\begin{proposition}\label{thmproof2}
Suppose that $X\ra S$ is of Cartier type. The functor $\rho:\underline{\calE}(X/\frakS) \ra \calE(X'/\frakS)$ defined in \eqref{era3rhoShiho} is continuous and cocontinuous for the étale topology \eqref{parettop} and the log flat topology \eqref{era4logflattop}.
\end{proposition}

\begin{proof}
A family $((U_i,\frakT_i,u_i) \ra (U,\frakT,u))_{i\in I}$ of $\underline{\calE}(X/\frakS)$ is a covering for the log flat topology (resp. étale topology) if and only if $(\rho(U_i,\frakT_i,u_i) \ra \rho(U,\frakT,u))_{i\in I}$ is a covering for the log flat topology (resp. étale topology). In addition, $\rho$ commutes with fiber products (\ref{lemfiberprodcom}). It follows, by (\cite{SGA43} III 1.6), that $\rho$ is continuous for both topologies. The cocontinuity follows from (\cite{SGA43} III 2.1) and \ref{lemcoco}.
\end{proof}

\begin{lemma}\label{lem1741}
For every object $(U,\frakT,u)$ of $\underline{\calE}(X/\frakS),$ the diagram
$$
\begin{tikzcd}
\text{ét}_{/U} \ar{rr}{\alpha_{(U,\frakT,u)}} \ar[swap, sloped]{d}{\sim} & & \underline{\calE}(X/\frakS) \ar{d}{\rho} \\
\text{ét}_{/U'} \ar[swap]{rr}{\alpha_{\rho(U,\frakT,u)}} & & \calE(X'/\frakS)
\end{tikzcd}
$$
is commutative, where the left vertical arrow is $V \mapsto V'.$
\end{lemma}

\begin{proof}
This is clear from the definitions of $\alpha_{(U,\frakT,u)}$ and $\alpha_{\rho(U,\frakT,u)}$ \eqref{eqKoko1861}.
\end{proof}

\begin{lemma}\label{Kokolem1833}
Suppose that $X\ra S$ is of Cartier type. Consider the morphisms of topoi
$$C_{X/\frakS}:\left ( \widetilde{\underline{\calE}}(X/\frakS), \Ox_{\underline{\calE}(X/\frakS)}\right ) \ra \left ( \widetilde{\calE}(X'/\frakS), \Ox_{\calE(X'/\frakS)}\right )$$
and
$$C_{X/\frakS,lf}:\left ( \widetilde{\underline{\calE}}_{lf}(X/\frakS), \Ox_{\underline{\calE}(X/\frakS)}\right ) \ra \left ( \widetilde{\calE}_{lf}(X'/\frakS), \Ox_{\calE(X'/\frakS)}\right ),$$
defined in \eqref{era3Clf}. Let $(U,\frakT,u)$ be an object of $\underline{\calE}(X/\frakS)$ and $\F$ a module of $\widetilde{\calE}(X'/\frakS)$ (resp. $\widetilde{\calE}_{lf}(X'/\frakS)$). We identify $\text{ét}_{/U}$ and $\text{ét}_{/U'}$ via the exact relative Frobenius $U\ra U'.$ We denote $C_{X/\frakS}$ (resp. $C_{X/\frakS,lf}$) by $C.$ Then
$$
(C^*\F)_{(U,\frakT,u)} = \F_{\rho(U,\frakT,u)}.
$$
In addition, if $(f,g)$ is a morphism in $\underline{\calE}(X/\frakS),$ then $\rho(f,g)=(f,g')$ and
$$c_{C^*\F,(f,g)}=c_{\F,(f,g')}.$$
\end{lemma}

\begin{proof}
This is a consequence of \ref{era3propdescentdata} and \ref{lem1741} as shown in the following computation:
\begin{alignat*}{2}
(C^*\F)_{(U,\frakT,u)} &= (C^*\F)\circ \alpha_{(U,\frakT,u)} \\
&= \F \circ \rho \circ \alpha_{(U,\frakT,u)} \\
&= \F \circ \alpha_{\rho(U,\frakT,u)} \\
&= \F_{\rho(U,\frakT,u)}.
\end{alignat*}
\end{proof}

\begin{proposition}\label{era3preserve}
Suppose that $X\ra S$ is of Cartier type. The inverse image functors of
$$C_{X/\frakS}:\left ( \widetilde{\underline{\calE}}(X/\frakS), \Ox_{\underline{\calE}(X/\frakS)}\right ) \ra \left ( \widetilde{\calE}(X'/\frakS), \Ox_{\calE(X'/\frakS)}\right )$$
and
$$C_{X/\frakS,lf}:\left ( \widetilde{\underline{\calE}}_{lf}(X/\frakS), \Ox_{\underline{\calE}(X/\frakS)}\right ) \ra \left ( \widetilde{\calE}_{lf}(X'/\frakS), \Ox_{\calE(X'/\frakS)}\right ),$$
defined in \eqref{era3Clf}, preserve crystals (resp. quasi-coherent modules).
\end{proposition}

\begin{proof}
In this proof, and for simplicity, we will denote $C_{X/\frakS}$ (resp. $C_{X/\frakS,lf}$) by $C.$ Let $(U,\frakT,u)$ be an object of $\underline{\calE}(X/\frakS)$ and $\F$ a module of $\widetilde{\calE}(X'/\frakS)$ (resp. $\widetilde{\calE}_{lf}(X'/\frakS)$). We canonically identify $\text{ét}_{/U}$ and $\text{ét}_{/U'}.$ By \ref{Kokolem1833}, we have
\begin{alignat*}{2}
(C^*\F)_{(U,\frakT,u)} &= \F_{\rho(U,\frakT,u)}.
\end{alignat*}
In addition, if $(f,g)$ is a morphism in $\underline{\calE}(X/\frakS),$ then $\rho(f,g)=(f,g')$ and
$$c_{C^*\F,(f,g)}=c_{\F,(f,g')}.$$
It follows that if $\F$ is a quasi-coherent module (resp. crystal), so is $C^*\F.$
\end{proof}

\begin{theorem}\label{thmKoko1926}
Suppose that $X\ra S$ is of Cartier type and that it lifts to a log smooth morphism of framed logarithmic formal schemes $(\frakX,Q) \ra (\frakS,P).$ We also suppose that the exact relative Frobenius $F_1:X\ra X'$ lifts to an $(\frakS,P)$-morphism of framed logarithmic formal schemes $(\frakX,Q) \ra (\frakX',Q'),$ such that $\frakX' \ra \frakS$ is log smooth. Consider the logarithmic formal schemes $R_{\frakX',1}$ and $Q_{\frakX}$ defined in \ref{parag86}, the morphism of formal groupoids $\nu:Q_{\frakX}\ra R_{\frakX',1}$ given in \ref{lem92}. Denote by $\nu_1:Q_1\ra R_1'$ the special fiber of $\nu$ and by $\mathcal{V}:\calR_1' \ra \calQ_1$ the corresponding morphism of Hopf algebras. Consider the diagram
\begin{equation}\label{diag185}
\begin{tikzcd}
\calC(X'/\frakS) \ar{r}{C_{X/\frakS}^*} \ar[swap, sloped]{d}{\sim} & \underline{\calC}(X/\frakS)  \ar[sloped]{d}{\sim} \\
\begin{Bmatrix}\Ox_{X'}\text{-modules\ with\ an}\\ \calR_1'\text{-stratification} \end{Bmatrix} \ar{r}{\Psi_0} & \begin{Bmatrix}\Ox_{X}\text{-modules\ with\ a}\\ \calQ_1\text{-stratification} \end{Bmatrix},
\end{tikzcd}
\end{equation}
where the vertical equivalences are given in \ref{equivRQ}, $C_{X/\frakS}$ is given in \eqref{era3C} and \ref{era3preserve} and $\Psi_0$ is defined by
$$\Psi_0:(\calE',\epsilon') \mapsto (F_1^*\calE',\mathcal{V}^*\epsilon').$$
Then the diagram \eqref{diag185} is commutative is commutative up to a canonical isomorphism.
\end{theorem}

\begin{proof}
If $\F$ is a module on $\calE(X'/\frakS)$ and $(f,g):(U_1,\frakT_1,u_1) \ra (U_2,\frakT_2,u_2)$ is a morphism of $\calE(X'/\frakS)$ then we have a $u_{1*}\Ox_{T_1}$-linear morphism \eqref{cflindata}
\begin{equation}\label{Oh181}
c_{\F,(f,g)}:\widetilde{f}^*\F_{(U_2,\frakT_2,u_2)} \ra \F_{(U_1,\frakT_1,u_1)},
\end{equation}
where 
$$
\widetilde{f}:\left (U_{1,\text{ét}},u_{1*}\Ox_{T_1}\right ) \ra \left (U_{2,\text{ét}},u_{2*}\Ox_{T_2}\right )
$$
is the morphism of ringed topoi \eqref{Kokoftilda}. The morphism \eqref{Oh181} is, by definition, equal to
\begin{equation}\label{Oh183}
c_{\F,(f,g)}:u_{1*}\Ox_{T_1} \otimes_{g^{-1}u_{2*}\Ox_{T_2}}g^{-1}\F_{(U_2,\frakT_2,u_2)} \ra \F_{(U_1,\frakT_1,u_1)},
\end{equation}
where the morphism $g^{-1}u_{2*}\Ox_{T_2} \ra u_{1*}\Ox_{T_1}$ is induced by the special fiber $f_1:T_1\ra T_2$ of $f.$
Let $Q_1$ and $R'_1$ be the special fibers of $Q_{\frakX}$ and $R_{\frakX',1}$ respectively, $q_1,q_2:Q_{\frakX} \ra \frakX$ and $r_1',r_2':R_{\frakX',1} \ra \frakX'$ the canonical projections and $\lambda_Q:\underline{Q_1} \ra X$ and $\lambda_{R'}:R'_1 \ra X'$ the canonical morphisms \eqref{erafprop13}.
Consider the object $(X,Q_{\frakX},\lambda_Q)$ of $\underline{\calE}(X/\frakS)$ and the objects $(X',R_{\frakX',1},\lambda_{R'})$ and $(X',\frakX',\op{Id}_{X'})$ of $\calE(X'/\frakS).$ As in \ref{equivRQDef}, $(X,\frakX,i_X)$ is an object of $\underline{\calE}(X/\frakS),$ where $i_X:\underline{X} \ra X$ is the canonical immersion, and the projections $q_1,q_2,r_1'$ and $r_2'$ define morphisms
\begin{alignat*}{2}
q_1,q_2:&(X,Q_{\frakX},\lambda_Q) \ra (X,\frakX,i_X), \\
r_1',r_2':&(X',R_{\frakX',1},\lambda_{R'}) \ra (X',\frakX',\op{Id}_{X'})
\end{alignat*}
of $\underline{\calE}(X/\frakS)$ and $\calE(X'/\frakS)$ respectively.
Let $\mathfrak{E}'$ be a crystal of $\calC(X'/\frakS).$ Denote by $\calR_1'$ the Hopf algebra corresponding to $R_1'.$ Then $\calR_1'=\lambda_{R'*}\Ox_{R_1'}.$ The morphisms \eqref{Oh183} corresponding to $r_1'$ and $r_2'$ are then equal to
$$
c_{\frakE',r_2'}:\calR_1'\otimes_{\Ox_{X'}}\frakE'_{(X',\frakX',\op{Id}_{X'})} \xrightarrow{\sim} \frakE'_{(X',R_{\frakX',1},\lambda_{R'})},
$$
$$
c_{\frakE',r_1'}:\frakE'_{(X',\frakX',\op{Id}_{X'})} \otimes_{\Ox_{X'}} \calR_1'\xrightarrow{\sim} \frakE'_{(X',R_{\frakX',1},\lambda_{R'})}.
$$
On one hand, by \ref{equivRQ}, the corresponding $\Ox_{X'}$-module with an $R_{\frakX',1}$-stratification $(\calE',\epsilon')$ is given as follows:
$$\calE'=\mathfrak{E}'_{(X',\frakX',\op{Id}_{X'})},$$
$$\epsilon':\calR_1'\otimes_{\Ox_{X'}}\calE' \xrightarrow{c_{\frakE',r_2'}} \mathfrak{E}'_{(X',R_{\frakX',1},\lambda_{R'})} \xrightarrow{c_{\frakE',r_1'}^{-1}}\calE' \otimes_{\Ox_{X'}} \calR_1'.$$
Its image by $\Psi_0$ is then the module
$$F_1^*\calE'=F_1^*\mathfrak{E}'_{(X',\frakX',\op{Id}_{X'})},$$
equipped with the stratification
$$\mathcal{V}^*\epsilon':\calQ_1\otimes_{\Ox_{X'}}\calE' \xrightarrow{\mathcal{V}^*c_{\frakE',r_2'}} \mathcal{V}^*\mathfrak{E}'_{(X',R_{\frakX',1},\lambda_{R'})}\xrightarrow{\mathcal{V}^*c_{\frakE',r_1'}^{-1}}\calE' \otimes_{\Ox_{X'}}\calQ_1.$$
This is equal to
$$
\mathcal{V}^*\epsilon':\calQ_1\otimes_{\Ox_{X}}(F_1^*\calE') \xrightarrow{\mathcal{V}^*c_{\frakE',r_2'}} \mathcal{V}^*\mathfrak{E}'_{(X',R_{\frakX',1},\lambda_{R'})}\xrightarrow{\mathcal{V}^*c_{\frakE',r_1'}^{-1}} (F_1^*\calE') \otimes_{\Ox_{X}}\calQ_1.
$$
On the other hand, by \ref{equivRQ}, the $\Ox_X$-module with a $Q_{\frakX}$-stratification $(\calE,\epsilon)$ corresponding to $\frakE:=C_{X/\frakS}^*\mathfrak{E}'$ is given as follows:
$$\calE=\left (C_{X/\frakS}^*\mathfrak{E}' \right )_{(X,\frakX,i_X)},$$
$$\epsilon:\calQ_1 \otimes_{\Ox_X} \calE\xrightarrow{c_{\frakE,q_2}} \left (C_{X/\frakS}^*\mathfrak{E}' \right )_{(X,Q_{\frakX},\lambda_Q)} \xrightarrow{c_{\frakE,q_1}^{-1}} \calE \otimes_{\Ox_X}\calQ_1.$$
By \ref{propShiho185}, we have
$$
\rho(X,\frakX,i_X)=(X',\frakX,F_1).
$$
By \ref{Kokolem1833}, we have
$$\calE=\left (C_{X/\frakS}^*\mathfrak{E}' \right )_{(X,\frakX,i_X)}=\mathfrak{E}'_{\rho(X,\frakX,i_X)}=\mathfrak{E}'_{(X',\frakX,F_1)}$$
and
$$
\left (C_{X/\frakS}^*\mathfrak{E}' \right )_{(X,Q_{\frakX},\lambda_Q)}=\frakE'_{\rho(X,Q_{\frakX},\lambda_Q)}=\frakE'_{(X',Q_{\frakX},\Delta_{\lambda_Q})}.
$$
The lifting $F:\frakX\ra \frakX'$ induces a morphism $(F,\op{Id}_{X'}):(X',\frakX,F_1) \ra (X',\frakX',\op{Id}_{X'})$ in $\calE(X'/\frakS).$ We then get a morphism of ringed topoi
$$
\widetilde{F}:(X'_{\text{ét}},F_{1*}\Ox_X) \ra (X'_{\text{ét}},\Ox_{X'}).
$$
Note that, since we identify $X_{\text{ét}}$ and $X'_{\text{ét}}$ via $F_1,$ we have $\widetilde{F}^*=F_1^*.$
We then get an isomorphism \eqref{cflindata}
\begin{equation}\label{zzcF}
c_{F_1}:F_1^*\calE' = F_1^*\mathfrak{E}'_{(X',\frakX',\op{Id}_{X'})} \xrightarrow{\sim} \mathfrak{E}'_{(X',\frakX,F_1)}=\calE.
\end{equation}
By \ref{lemKhaminei1}, $\nu:Q_{\frakX} \ra R_{\frakX',1}$ defines a morphism
\begin{equation}\label{zabb222}
\nu:\rho(X,Q_{\frakX},\lambda_Q)=(X',Q_{\frakX},\Delta_{\lambda_Q}) \ra (X',R_{\frakX',1},\lambda_{R'})
\end{equation}
in $\calE(X'/\frakS).$ We then have a morphism of ringed topoi
$$
\widetilde{\nu}: \left (X'_{\text{ét}}, \left (\Delta_{\lambda_Q} \right )_*\Ox_{Q_1} \right ) \ra \left (X'_{\text{ét}},\lambda_{R'*}\Ox_{R_1'} \right )=\left (X'_{\text{ét}},\calR_1' \right ).
$$
By \ref{lemKhaminei1}, we have
$$
\left (\Delta_{\lambda_Q} \right )_*\Ox_{Q_1}=F_{1*}q_{1*}\Ox_{Q_1}=F_{1*}\calQ_1.
$$
Since we identify $X_{\text{ét}}$ and $X'_{\text{ét}}$ via $F_1$ and $\frakE'$ is a crystal, we get an isomorphism
$$c_{\frakE',\nu}:\mathcal{V}^*\mathfrak{E}'_{(X',R_{\frakX',1},\lambda_{R'})} = \calQ_1 \otimes_{\calR_1'} \mathfrak{E}'_{(X',R_{\frakX',1},\lambda_{R'})} \xrightarrow{\sim} \mathfrak{E}'_{\rho(X,Q_{\frakX},\lambda_Q)},$$
fitting into the commutative diagram
$$
\begin{tikzcd}
\calQ_1 \otimes_{\Ox_X} \left (F_1^*\calE' \right ) \ar{rr}{\mathcal{V}^*c_{\frakE',r_2'}} \ar[swap]{d}{\op{Id} \otimes c_{F_1}}  & & \mathcal{V}^*\mathfrak{E}'_{(X',R_{\frakX',1},\lambda_{R'})}\ar{rr}{\mathcal{V}^*c_{\frakE,r_1'}^{-1}} \ar{d}{c_{\frakE',\nu}} & & \left (F_1^*\calE'\right ) \otimes_{\Ox_X}\calQ_1 \ar{d}{c_{F_1} \otimes \op{Id}} \\
\calQ_1 \otimes_{\Ox_X} \calE\ar{rr}{c_{\frakE,q_2}} & & \mathfrak{E}'_{\rho(X,Q_{\frakX},\lambda_Q)} \ar{rr}{c_{\frakE,q_1}^{-1}} & & \calE \otimes_{\Ox_X} \calQ_1.
\end{tikzcd}
$$
This proves that the isomorphism $c_{F_1}:F_1^*\calE' \xrightarrow{\sim} \calE$ \eqref{zzcF} underlies an isomorphism of stratified modules
$$
(F_1^*\calE',\nu^*\epsilon') \xrightarrow{\sim} (\calE,\epsilon).
$$
This finishes the proof.
\end{proof}

\begin{lemma}\label{lemKhaminei1}
Keep the hypothesis of \ref{thmKoko1926} and let $q_1,q_2:Q_1 \ra X$ be the canonical projections. Then $\Delta_{\lambda_Q}=F_1\circ q_1=F_1\circ q_2$ and
the diagram
\begin{equation}\label{Khaminei4}
\begin{tikzcd}
Q_1 \ar{r} \ar[swap]{d}{\Delta_{\lambda_Q}} & R_1' \ar{dl}{\lambda_{R'}} \\
X' & 
\end{tikzcd}
\end{equation}
is commutative.
\end{lemma}

\begin{proof}
Since $X\ra S$ is of Cartier type, the relative Frobenius $F_{X/S}$ is exact and hence $F_{X/S}=F_1.$
Then, by \eqref{morphismShiho2}, $\Delta_{\lambda_Q}$ is the unique morphism fitting into the commutative diagram
\begin{equation}\label{Khaminei1}
\begin{tikzcd}
Q_1\ar{r}{G_{Q_1}} \ar[swap]{dr}{\Delta_{\lambda_Q}} \ar[bend right=30]{rdd} & \underline{Q_1} \ar{dr}{\lambda_Q} & \\
 & X'\ar{d}\ar{r}{\pi} & X\ar{d} \\
 & S \ar{r}{F_S} & S.
\end{tikzcd}
\end{equation}
The composition $$Q_1 \xrightarrow{q_i} X \xrightarrow{F_1} X' \ra S,$$ where $X' \ra S$ is the canonical projection, and $Q_1 \ra S$ are clearly equal. It is hence sufficient to prove that the compositions
\begin{equation}\label{Khaminei2}
Q_1 \xrightarrow{q_i} X \xrightarrow{F_1} X'\xrightarrow{\pi} X
\end{equation}
and
\begin{equation}\label{Khaminei3}
Q_1 \xrightarrow{G_{Q_1}} \underline{Q_1} \xrightarrow{\lambda_{Q}} X
\end{equation}
are equal.
The composition \eqref{Khaminei2} is equal to
$$
Q_1\xrightarrow{q_i} X \xrightarrow{F_X} X.
$$
This is equal to
$$
Q_1 \xrightarrow{F_{Q_1}} Q_1 \xrightarrow{q_i} X
$$
and then to
$$
Q_1 \xrightarrow{G_{Q_1}} \underline{Q_1} \hookrightarrow Q_1 \xrightarrow{q_i} X.
$$
By definition of $\lambda_Q,$ this is equal to \eqref{Khaminei3}.
The commutativity of \eqref{Khaminei4} follows from the definition of $\nu$ and the fact that $\Delta_{\lambda_Q}=F_1\circ q_i.$
\end{proof}

\newpage

\section{Logarithmic Oyama topoi and the Cartier transform over a trivial log point}

\begin{parag}
In this section, we fix a perfect field $k$ of positive characteristic $p$ and denote its ring of Witt vectors by $W(k).$ We equip $\op{Spf}W(k)$ with the trivial logarithmic structure. Unless explicitly stated, all formal schemes considered are $p$-adic and all logarithmic structures are fs. If a gothic letter $\frakT$ denotes a $p$-adic logarithmic formal scheme, the corresponding roman letter $T$ will denote the logarithmic scheme obtained from $\frakT$ by reduction modulo $p$ \eqref{propkey}. We also denote by $\frakT_n,$ for any positive integer $n,$ the logarithmic scheme obtained from $\frakT$ by reduction modulo $p^n$ \eqref{propkey}.
Let $\frakS=\op{Spf}W(k)$ and $X\ra S$ a log smooth morphism of logarithmic schemes which is of Cartier type.
\end{parag}

\begin{parag}
By \ref{rem178} and \ref{propfT/S}, for an object $(U,\frakT,u)$ of $\underline{\calE}(X/\frakS),$ we have the following commutative diagram
$$
\begin{tikzcd}
U \ar[swap]{d}{F_{U/S}} & & \underline{T}\ar[swap]{ll}{u} \ar[swap]{d}{F_{\underline{T}/S}} \ar[hook]{rr} & & T\ar{d}{F_{T/S}} \ar{dll}{f_{T/S}} \\
U' & & \underline{T}' \ar{ll}{u'} \ar[hook]{rr} & & T',
\end{tikzcd}
$$
where the vertical arrows are the exact relative Frobenius morphisms (\ref{PFrob}), $U'=U\times_{S,F_S}S,$ $u':\underline{T}' \ra U'$ is the morphism induced by $u$ and $f_{T/S}$ is defined in \ref{propfT/S}. The commutativity of the right upper triangle follows from the commutativity of the right square, the lower right triangle and the fact that $\underline{T}' \ra T'$ is an immersion, hence a monomorphism.
In this case, the functor \eqref{era3rhoShiho} is given by
\begin{equation}\label{era3rho}
\rho:\begin{array}[t]{clc}
\underline{\calE}(X/\frakS) & \ra & \calE(X'/\frakS) \\
(U,\frakT,u) & \mapsto & (U',\frakT,u'\circ f_{T/S}).
\end{array}
\end{equation}
We will see in \ref{era3equivlfcrystals} that the direct image and inverse image functors of $C_{X/\frakS,lf}$ are equivalences of categories quasi-inverse to each other.  We will also see in \ref{lemdirectimage} that $C_{X/\frakS}^*$ induces a fully faithful functor between crystals. We start with some lemmas.
\end{parag}

\begin{lemma}[\cite{DXU19} 8.5]\label{era3lem1}
Let $u:T\ra U$ be an affine morphism of schemes. The functors
$$u_*:(T_{\text{ét}},\Ox_T) \ra (U_{\text{ét}},u_*\Ox_T),\ u^*:(U_{\text{ét}},u_*\Ox_T) \ra (T_{\text{ét}},\Ox_T)$$
induce equivalences of categories quasi-inverse to each other between the category of quasicoherent $\Ox_T$-modules of $T_{\text{ét}}$ \eqref{Not6} and the category of quasicoherent $u_*\Ox_T$-modules of $U_{\text{ét}}.$
\end{lemma}

\begin{proof}
Recall that, for every scheme $X,$ $\op{QCoh}(X_{\text{zar}})$ and $\op{QCoh}(X_{\text{ét}})$ are canonically equivalent (\href{https://stacks.math.columbia.edu/tag/03DX}{Proposition 03DX}).

We start by proving that $u_*$ is exact. Let $g:\calM \ra \calN$ be a surjective morphism of quasicoherent modules of $(T_{\text{ét}},\Ox_T).$ Let $V \ra U$ be an étale morphism such that $V$ is affine. The projection $T\times_UV \ra V$ is affine and $\calM$ and $\calN$ are quasi-coherent so the morphism $(u_*g)(V)=g(T\times_UV)$ is surjective. Since $u_*$ is also left exact, it is exact.

Now let $\calM$ be a quasicoherent module of $(T_{\text{ét}},\Ox_T).$ We will prove that the adjunction morphism $u^*u_*\calM \ra \calM$ is an isomorphism. For this, we may suppose that $U$ is affine. The scheme $T$ is then also affine and hence we have an exact sequence
$$\Ox_T^{\oplus I} \ra \Ox_T^{\oplus J} \ra \calM \ra 0.$$
Since $u_*$ is exact, we deduce the exactness of the sequence
$$u_*\left (\Ox_T^{\oplus I} \right ) \ra u_* \left (\Ox_T^{\oplus J} \right ) \ra u_*\calM \ra 0.$$
The quasi-coherent module $\Ox_T^{\oplus I}$ is associated to the $\Gamma(T,\Ox_T)$-module $\bigoplus_{i\in I}\Gamma(T,\Ox_T).$ Similarly, the quasi-coherent module $\left (u_*\Ox_T \right )^{\oplus I}$ is associated to the $\Gamma(U,\Ox_U)$-module $\bigoplus_{i\in I} \Gamma(U,u_*\Ox_T)=\bigoplus_{i\in I}\Gamma(T,\Ox_T).$ It follows that the canonical morphism
$$\left (u_*\Ox_T \right )^{\oplus I} \ra u_*\left ( \Ox_T ^{\oplus I} \right )$$
is clearly an isomorphism.
Hence the quasi-coherence of $u_*\calM.$
Consider the following commutative diagram of $\Ox_T$-modules
$$
\begin{tikzcd}
u^*\left (u_*\left (\Ox_T \right )^{\oplus I}\right )  \ar{r} \ar{d} & u^* \left (u_* \left (\Ox_T\right )^{\oplus J}\right )  \ar{r} \ar{d} & u^*u_*\calM \ar{d}\ar{r} & 0 \\
\Ox_T^{\oplus I} \ar{r} & \Ox_T^{\oplus J} \ar{r}& \calM \ar{r} & 0
\end{tikzcd}
$$
Since $u_*$ and $u^*$ are both right exact, the rows are exact. The middle and left vertical arrows are clearly isomorphisms so the right vertical arrow is also an isomorphism.

Now let $\calN$ be a quasi-coherent $u_*\Ox_T$-module of $U_{\text{ét}}.$ We will prove that the adjunction morphism $\calN \ra u_*u^*\calN$ is an isomorphism. For this, we may work étale locally on $U$ and suppose that we have an exact sequence
$$\left (u_*\Ox_T\right )^{\oplus I} \ra \left (u_*\Ox_T \right )^{\oplus J} \ra \calN \ra 0.$$
Consider the following commutative diagram of $u_*\Ox_T$-modules
$$
\begin{tikzcd}
\left (u_*\Ox_T \right )^{\oplus I} \ar{r} \ar{d} & \left (u_*\Ox_T\right )^{\oplus J} \ar{r} \ar{d} & \calN \ar{r} \ar{d} & 0 \\
u_*u^*\left (\left (u_*\Ox_T\right )^{\oplus I}\right ) \ar{r} & u_*u^*\left (\left (u_*\Ox_T\right )^{\oplus J} \right ) \ar{r} & u_*u^*\calN \ar{r} & 0
\end{tikzcd}
$$
The rows are exact and the left and middle vertical arrows are clearly isomorphisms. The result follows.
\end{proof}

\begin{lemma}\label{era3lem2}
Any crystal $\F$ of $\calC^{\text{qcoh}}(X/\frakS)$ (resp. $\underline{\calC}^{\text{qcoh}}(X/\frakS)$) is a sheaf for the log flat topology of $\calE(X/\frakS)$ (resp. $\underline{\calE}(X/\frakS)$). 
\end{lemma}

\begin{proof}
We just prove the result for a crystal $\F$ of $\calC^{\text{qcoh}}(X/\frakS)$ as the proof for $\underline{\calC}^{\text{qcoh}}(X/\frakS)$ is similar. Let $(\F_{(U,\frakT,u)},c_{(f,g)})$ be the linearized descent data associated with $\F$ and let $((U_i,\frakT_i,u_i) \xrightarrow{(f_i,g_i)} (U,\frakT,u))_{i\in I}$ be a log flat covering, where $f_i:\frakT_i \ra \frakT$ and $g_i:U_i \ra U.$ For any $i,j\in I,$ let
$$(U_{ij},\frakT_{ij},u_{ij})=(U_i,\frakT_i,u_i)\times_{(U,\frakT,u)}(U_j,\frakT_j,u_j)$$
and
$$(f_{ij},g_{ij}):(U_{ij},\frakT_{ij},u_{ij}) \ra (U,\frakT,u)$$
the canonical morphism.
We will abusively denote by $f_i:T_i\ra T$ and $f_{ij}:T_{ij} \ra T$ the special fiber of $f_i$ and $f_{ij}$ respectively.
We have to prove that the sequence
\begin{equation}\label{era3exseq1}
0 \ra \F(U,\frakT,u) \ra \prod_{i\in I}\F(U_i,\frakT_i,u_i) \rightrightarrows \prod_{i,j\in I}\F(U_{ij},\frakT_{ij},u_{ij})
\end{equation}
is exact.
By \ref{era3lem1}, we can consider $\F_{(U_i,\frakT_i,u_i)},$ $\F_{(U,\frakT,u)}$ and $\F_{(U_{ij},\frakT_{ij},u_{ij})}$ as modules of $(T_{i,\text{ét}},\Ox_{T_i}),$ $(T_{\text{ét}},\Ox_T)$ and $(T_{ij,\text{ét}},\Ox_{T_{ij}})$ respectively. The sequence \eqref{era3exseq1} becomes equal to
\begin{equation}\label{era3exseq1c}
0 \ra \F_{(U,\frakT,u)}(T) \ra \prod_{i\in I}\F_{(U_i,\frakT_i,u_i)}(T_i) \rightrightarrows \prod_{i,j\in I}\F_{(U_{ij},\frakT_{ij},u_{ij})}(T_{ij}).
\end{equation}
Since $\F$ is a crystal, the morphisms
$$c_{(f_i,g_i)}:f_i^*\F_{(U,\frakT,u)} \xrightarrow{\sim} \F_{(U_i,\frakT_i,u_i)}$$
$$c_{(f_{ij},g_{ij})}:f_{ij}^*\F_{(U,\frakT,u)} \xrightarrow{\sim} \F_{(U_{ij},\frakT_{ij},u_{ij})}$$
are isomorphisms. Let $\calG=\F_{(U,\frakT,u)}.$ The sequence \eqref{era3exseq1c} becomes equal to
\begin{equation}
0 \ra \calG(T) \ra \prod_{i\in I}(f_i^*\calG)(T_i) \rightrightarrows \prod_{i,j\in I} (f_{ij}^*\calG)(T_{ij}).
\end{equation}
It is thus sufficient to prove the exactness of the sequence
\begin{equation}\label{era3exseq2}
0 \ra \calG \ra \prod_{i\in I}f_{i*}f_i^*\calG \rightrightarrows \prod_{i,j\in I} f_{ij*}f_{ij}^*\calG.
\end{equation}
This follows from \ref{eraflogflatdescent}.
\end{proof}

\begin{proposition}\label{era3propequivlfetale}
The canonical functors
$$\widetilde{\calE}_{lf}(X/\frakS) \ra \widetilde{\calE}(X/\frakS),\ \widetilde{\underline{\calE}}_{lf}(X/\frakS) \ra \widetilde{\underline{\calE}}(X/\frakS)$$
induce equivalences of categories
$$\calC_{lf}^{\text{qcoh}}(X/\frakS) \xrightarrow{\sim} \calC^{\text{qcoh}}(X/\frakS),\ \underline{\calC}_{lf}^{\text{qcoh}}(X/\frakS) \xrightarrow{\sim} \underline{\calC}^{\text{qcoh}}(X/\frakS).$$
\end{proposition}

\begin{proof}
We only prove the result for $\calE(X/\frakS)$ as the proof for $\underline{\calE}(X/\frakS)$ is similar. It is clear that $\widetilde{\calE}_{lf}(X/\frakS) \ra \widetilde{\calE}(X/\frakS)$ induces a functor
\begin{equation}\label{eraf17291}
\calC_{lf}^{\text{qcoh}}(X/\frakS) \ra \calC^{\text{qcoh}}(X/\frakS).
\end{equation}
The functor $\widetilde{\calE}_{lf}(X/\frakS) \ra \widetilde{\calE}(X/\frakS)$ is clearly fully faithful and hence so is \eqref{eraf17291}. It is thus sufficient to prove that a crystal of $\calC^{\text{qcoh}}(X/\frakS)$ is a sheaf for the log flat topology. This is proved in \ref{era3lem2}.
\end{proof}

\begin{proposition}[\cite{Oyama} 4.2.1]\label{Oyamathm}
Let $\calC$ and $\calD$ be two sites, such that the topology of $\calD$ is defined by a pretopology, and $u:\calC\ra \calD$ a functor. If $u$ satisfies the following conditions:
\begin{enumerate}
\item $u$ is fully faithful.
\item $u$ is continuous and cocontinuous.
\item For any object $Y$ of $\calD,$ there exists a covering of the form $(u(X_i)\ra Y)_{i\in I}.$
\end{enumerate}
then $u$ induces an equivalence of topoi
$$u=(u^*,u_*):\widetilde{\calC} \ra \widetilde{\calD},$$
such that $u^*$ is composition with $u$ and $u_*$ is a right adjoint of $u^*$ for presheaves.
\end{proposition}

\begin{proposition}\label{thmproof1}
The functor $\rho:\underline{\calE}(X/\frakS) \ra \calE(X'/\frakS)$ defined in \eqref{era3rho} is fully faithful.
\end{proposition}

\begin{proof}
Let $(f_1,g_1),(f_2,g_2):(U_1,\frakT_1,u_1) \ra (U_2,\frakT_2,u_2)$ be two morphisms in $\underline{\calE}(X/\frakS)$ such that $\rho(f_1,g_1)=\rho(f_2,g_2).$ Since
$$\rho(f_i,g_i)=(f_i,g_i')$$
for $i=1,2,$ we get $g_1'=g_2'$ and so $g_1=g_2.$ The functor $\rho$ is then faithful. We now prove its fullness. Let $(U_1,\frakT_1,u_1)$ and $(U_2,\frakT_2,u_2)$ be two objects of $\underline{\calE}(X/\frakS)$ and $(f,h):\rho(U_1,\frakT_1,u_1) \ra \rho(U_2,\frakT_2,u_2)$ a morphism in $\calE(X'/\frakS).$ There exists a unique morphism $g:U_1 \ra U_2$ such that $h=g'.$ In the diagram
$$
\begin{tikzcd}
T_1 \ar{r}{f_1} \ar[swap]{d}{f_{T_1/S}}  & T_2\ar{d}{f_{T_2/S}} \\
\underline{T_1'} \ar{r}{\underline{f_1}'} \ar[swap]{d}{u_1'} & \underline{T_2'}\ar{d}{u_2'} \\
U_1' \ar{r}{g'} & U_2',
\end{tikzcd}
$$
the rectangle and the upper square are commutative.
Since $f_{T_1/S}$ is an epimorphism (\ref{propfT/S}), the lower square is commutative. Since the functor $T \mapsto T'$ is faithful, we get
$$u_2\circ \underline{f_1} = g\circ u_1.$$
\end{proof}

\begin{lemma}\label{thmproof3}
For any object $(W,\frakT,v)$ of $\calE(X'/\frakS),$ there exists a covering of $(W,\frakT,v)$ for the log flat topology \eqref{era4logflattop} of the form
$$\left ( \rho(U_i,\frakT_i,u_i) \ra (W,\frakT,v)\right )_{i\in I}.$$
\end{lemma}

\begin{proof}
Let $F_1:X\ra X'$ be the relative Frobenius morphism of $X$ over $S$ (which is automatically exact since $X\ra S$ is saturated by \cite{Ogus2018} III 2.5.4). 
There exists a unique strict étale morphism $V \ra X$ such that $W=V\times_XX'=V'.$ The assertion of the proposition being étale local on $W,$ we may suppose that there exists a $p$-adic logarithmic formal scheme $\frakV$ and $\frakV'$ log smooth over $\frakS,$ fitting into the cartesian squares
$$
\begin{tikzcd}
V \ar{r} \ar{d} & \frakV \ar{d} & & V' \ar{r} \ar{d} & \frakV' \ar{d} \\
S \ar{r}& \frakS & & S \ar{r} & \frakS
\end{tikzcd}
$$
We may also suppose that the restriction $F_1:V \ra V'$ and $v:T \ra V'$ lift, respectively, to morphisms $F:\frakV \ra \frakV'$ and $\frakT \ra \frakV'$ of $p$-adic logarithmic formal schemes such that the diagrams
$$
\begin{tikzcd}
V' \ar{rr} & & \frakV' \ar{d} & & V' \ar{rr} & & \frakV' \ar{d} \\
V \ar{r} \ar{u}{F_1} & \frakV \ar{r} \ar{ur}{F} & \frakS & & T\ar{u}{v} \ar{r} & \frakT \ar{ur} \ar{r} & \frakS
\end{tikzcd}
$$
are commutative.
Set $\frakZ=\frakT \times_{\frakV'}^{\op{log}} \frakV$ and $u:\underline{Z} \ra Z \ra V$ the composition of the canonical immersion $\underline{Z} \ra Z$ and the morphism $Z \ra V$ induced by the projection $\frakZ \ra \frakV.$ Since $v:T \ra V'=W$ and the canonical morphism $Z \ra T\times_{V'}V$ are affine (\cite{Ogus2018} III 2.1.6), the morphism $u$ is also affine. We obtain an object $(V,\frakZ,u)$ of $\underline{\calE}(X/\frakS)$ and a morphism
$$\rho(V,\frakZ,u)=(V',\frakZ,u'\circ f_{Z/S}) \ra (V',\frakT,v),$$
given by the projection $\frakZ \ra \frakT$ and $\op{Id}_{V'}.$ The result then follows from \ref{logflatfiber} applied to $F$ and the fact that the projection $Z \ra T$ is the base change of $F_1:V \ra V',$ which is a log flat covering (\ref{Flogflatcover}).
\end{proof}

\begin{theorem}\label{era3equivlfcrystals}
The functor $\rho:\underline{\calE}(X/\frakS) \ra \calE(X'/\frakS),$ defined in \eqref{era3rho}, induces an equivalence of topoi
$$
C_{X/S,lf}:\widetilde{\underline{\calE}}_{lf}(X/\frakS) \ra \widetilde{\calE}_{lf}(X'/\frakS).
$$
\end{theorem}

\begin{proof}
This is a consequence of \ref{Oyamathm}, \ref{thmproof1}, \ref{thmproof2} and $\ref{thmproof3}.$
\end{proof}

\begin{theorem}\label{lemdirectimage}
The functor $C_{X/\frakS,lf}^*$ induces a fully faithful functor
$$\calC^{\text{qcoh}}(X'/\frakS) \ra \underline{\calC}^{\text{qcoh}}(X/\frakS).$$
\end{theorem}

\begin{proof}
For a ringed topos $(T,\Ox),$ we denote by $\op{Mod}(T,\Ox)$ the category of $\Ox$-modules of $T.$
By \ref{era3preserve} and \ref{era3propequivlfetale}, the commutative diagram
$$
\begin{tikzcd}
\op{Mod}\left ( \widetilde{\calE}_{lf}(X'/\frakS), \Ox_{\calE(X'/\frakS)}\right ) \ar{r}{C_{X/\frakS,lf}^*} \ar[hook]{d} & \op{Mod}\left ( \widetilde{\underline{\calE}}_{lf}(X/\frakS), \Ox_{\underline{\calE}(X/\frakS)}\right ) \ar[hook]{d} \\
\op{Mod}\left ( \widetilde{\calE}(X'/\frakS), \Ox_{\calE(X'/\frakS)}\right ) \ar{r}{C_{X/\frakS}^*} & \op{Mod}\left ( \widetilde{\underline{\calE}}(X/\frakS), \Ox_{\underline{\calE}(X/\frakS)}\right )
\end{tikzcd}
$$
induces
$$
\begin{tikzcd}
\calC^{\text{qcoh}}_{lf}(X'/\frakS) \ar{r}{C_{X/\frakS,lf}^*} \ar[swap,sloped]{d}{\sim} & \underline{\calC}^{\text{qcoh}}_{lf}(X/\frakS) \ar[sloped]{d}{\sim} \\
\calC^{\text{qcoh}}(X'/\frakS)\ar{r}{C_{X/\frakS}^*} & \underline{\calC}^{\text{qcoh}}(X/\frakS).
\end{tikzcd}
$$
By \ref{era3equivlfcrystals}, $C_{X/\frakS,lf}^*$ and the functor it induces between crystals is fully faithful, hence the result.
\end{proof}

\newpage

\section{Indexed logarithmic Oyama topoi}

\begin{parag}\label{intro19}
In this section, we keep the hypothesis and notation of \ref{intro18}. For an étale $X$-scheme $U,$ we equip $U$ with the logarithmic structure $\left (\calM_{X}\right )_{|U},$ restriction of $\calM_X$ on $U,$ so that $U\ra X$ is strict.
For every logarithmic scheme $T,$ we denote by $\upmu_T$ the sheaf $\ov{\calM}_T^{gp}.$
\end{parag}

\begin{parag}\label{paragu}
The functors
$$
u_X':\begin{array}[t]{clc}
\calE(X'/\frakS) & \ra & \text{ét}_{/X'} \\
(U,\frakT,u) & \mapsto & U
\end{array},\ 
\underline{u}_X:\begin{array}[t]{clc}
\underline{\calE}(X/\frakS) & \ra & \text{ét}_{/X} \\
(U,\frakT,u) & \mapsto & U
\end{array}
$$
clearly satisfy properties 1 and 2 of \ref{Kokolemcocont}. If $(U,\frakT,u)$ is an object of $\calE(X'/\frakS)$ (resp. $\underline{\calE}(X/\frakS)$) and $(V_i \ra U)_{i\in I}$ is an étale covering, then, by \ref{era3paralphaKoko} (1),
$$\left (\alpha_{(U,\frakT,u)}(V_i) \ra (U,\frakT,u) \right )_{i\in I}$$
is an étale covering. It follows that the functors $u_X'$ and $\underline{u}_X$ satisfy property 3 of \ref{Kokolemcocont}. We conclude that
$u_X'$ and $\underline{u}_X$ are continuous and cocontinuous for the étale topologies and they induce morphisms of topoi, that we abusively denote by
\begin{equation}\label{Kokou}
u_X':\widetilde{\calE}(X'/\frakS) \ra X'_{\text{ét}},\ \underline{u}_X:\widetilde{\underline{\calE}}(X/\frakS) \ra X_{\text{ét}},
\end{equation}
such that the inverse image functors are composition with $u_X'$ and $\underline{u}_X$ respectively. 

We denote by $\mathbbm{T}'$ (resp. $\underline{\mathbbm{T}}$) the split fibered $\widetilde{\calE}(X'/\frakS)$-category (resp. $\widetilde{\underline{\calE}}(X/\frakS)$-category) defined, for every sheaf $\F$ of $\widetilde{\calE}(X'/\frakS)$ (resp. $\widetilde{\underline{\calE}}(X/\frakS)$), by the localized category $\widetilde{\calE}(X'/\frakS)_{/\F}$ (resp. $\widetilde{\underline{\calE}}(X/\frakS)_{/\F}$) and, for every morphism $\varphi:\F \ra \calG,$ by the inverse image functor of the morphism of topoi
\begin{equation}\label{Kokotopoi1}
\widetilde{\calE}(X'/\frakS)_{/\F} \ra \widetilde{\calE}(X'/\frakS)_{/\calG} \left( \text{resp.}\ \widetilde{\underline{\calE}}(X/\frakS)_{/\F} \ra \widetilde{\underline{\calE}}(X/\frakS)_{/\calG}\right )
\end{equation}
induced by $\varphi.$ By (\cite{Giraud} II 3.4.4), $\mathbbm{T}'$ (resp. $\underline{\mathbbm{T}}$) is a stack. We consider the fibered category $\textgoth{T}' \ra X'_{\text{ét}}$ (resp. $\underline{\textgoth{T}} \ra X_{\text{ét}}$) defined by base change of $\mathbbm{T}' \ra \widetilde{\calE}(X'/\frakS)$ (resp. $\underline{\mathbbm{T}} \ra \widetilde{\underline{\calE}}(X/\frakS)$) by the functor $u_X'^{-1}$ (resp. $\underline{u}_X^{-1}$):
\begin{equation}\label{stack1}
\begin{tikzcd}
\textgoth{T}' \ar{r} \ar{d} & \mathbbm{T}' \ar{d} & \underline{\textgoth{T}} \ar{r} \ar{d} & \underline{\mathbbm{T}} \ar{d} \\
X'_{\text{ét}} \ar{r}{u_X'^{-1}} & \widetilde{\calE}(X'/\frakS) & X_{\text{ét}} \ar{r}{\underline{u}_X^{-1}} & \widetilde{\underline{\calE}}(X/\frakS). 
\end{tikzcd}
\end{equation}
By (\cite{Giraud} II 3.1.1), they are stacks.
Let $U$ be an étale $X$-scheme and denote by $U^a$ the associated sheaf of $X_{\text{ét}}.$ By definition, for any object $(V,\frakZ,w)$ of $\widetilde{\underline{\calE}}(X/\frakS),$ we have
$$
\underline{u}_X^{-1}\left (U^a\right )(V,\frakZ,w) = U^a(V)=\op{Hom}_{X}(V,U).
$$
A morphism $(V,\frakZ,w) \ra \underline{u}_X^{-1}(U^a)$ is then equivalent to an $X$-morphism $V\ra U.$ It follows that
$$
\begin{array}[t]{clc}
\underline{\calE}(X/\frakS)_{/\underline{u}_X^{-1}(U^a)} & \ra & \underline{\calE}(U/\frakS) \\
\left ((V,\frakZ,w) \ra \underline{u}_X^{-1}(U^a) \right ) & \mapsto & (V,\frakZ,w)
\end{array}
$$
is an equivalence of categories. This functor clearly satisfies the conditions of \ref{Kokolemcocont} for the étale topology on the target and the topology on the source induced by that of $\underline{\calE}(X/\frakS).$ It is hence continuous and cocontinuous and induces an equivalence of topoi
\begin{equation}\label{stackequiv1}
\widetilde{\underline{\calE}}(X/\frakS)_{/\underline{u}_X^{-1}(U^a)} \xrightarrow{\sim} \widetilde{\underline{\calE}}(U/\frakS).
\end{equation}
Similarly, if $U$ is an étale $X'$-scheme, then we have an equivalence of topoi
\begin{equation}\label{stackequiv2}
\widetilde{\calE}(X'/\frakS)_{/u_X'^{-1}(U^a)} \xrightarrow{\sim} \widetilde{\calE}(U/\frakS).
\end{equation}
For every morphism of étale $X$-schemes $g:V\ra U,$ the functor
$$
\underline{\rho_g}:\underline{\calE}(V/\frakS) \ra \underline{\calE}(U/\frakS),
$$
defined by composition with $g,$ is clearly continuous and cocontinuous. It then induces a morphism of topoi, abusively denoted
\begin{equation}\label{eqrhoEST1}
\underline{\rho_g}:\widetilde{\underline{\calE}}(V/\frakS) \ra \widetilde{\underline{\calE}}(U/\frakS),
\end{equation}
whose inverse image functor is the composition with the functor $\underline{\rho_g}.$
The morphism of topoi $\underline{\rho_g}$ corresponds to the localization morphism of topoi $\underline{\mathbbm{T}}_{V^a} \ra \underline{\mathbbm{T}}_{U^a}$ \eqref{Kokotopoi1}.
Similarly, if $g:V\ra U$ is a morphism of étale $X'$-schemes, then we have a continuous and cocontinuous functor
$$
\rho_{g}:\calE(V/\frakS) \ra \calE(U/\frakS),
$$
that induces a morphism of topoi, abusively denoted by
\begin{equation}\label{eqrhoEST2}
\rho_{g}:\widetilde{\calE}(V/\frakS) \ra \widetilde{\calE}(U/\frakS),
\end{equation}
and that corresponds to the localization morphism $\mathbbm{T}'_{V^a} \ra \mathbbm{T}'_{U^a}$ \eqref{Kokotopoi1}.
If no ambiguity arises, we will denote $\rho_g$ and $\underline{\rho_g}$ by $\rho_{U,V}$ and $\underline{\rho_{U,V}}$ respectively.
\end{parag}

\begin{parag}\label{paragMuzan194}
Let $\calI_X$ be a sheaf of $X_{\text{ét}}.$ We identify the small étale sites of $X$ and $X'$ via the exact relative Frobenius $F_1:X\ra X'$ and hence consider $\calI_X$ also as a sheaf of $X'_{\text{ét}}.$ We set
$$\calI'=u_X'^{-1}\calI_X=\calI_X\circ u_X',\ \underline{\calI}=\underline{u}_X^{-1}\calI_X=\calI_X\circ \underline{u}_X.$$
Let $U$ be an étale $X$-scheme and $s\in \Gamma(U,\calI_X).$
For every object $(V,\frakZ,v)$ of $\calE(X'/\frakS)$ (resp. $\underline{\calE}(X/\frakS)$), we have
$$s_{|V} \in \Gamma(V,\calI_X)=\Gamma\left ( (V,\frakZ,v),\calI'\right )\ \left (\text{resp.}\ \Gamma \left ((V,\frakZ,v),\underline{\calI}\right )\right ).$$
So we have functors
$$
\alpha_s':\begin{array}[t]{clc}
\calE(U'/\frakS) & \ra & \calE(X'/\frakS)_{/\calI'} \\
(V,\frakZ,v) & \mapsto & \left ( \rho_{X',U'}(V,\frakZ,v),s_{|V}\right ),
\end{array}
$$
$$
\underline{\alpha}_s:\begin{array}[t]{clc}
\underline{\calE}(U/\frakS) & \ra & \underline{\calE}(X/\frakS)_{/\underline{\calI}} \\
(V,\frakZ,v) & \mapsto & \left ( \underline{\rho_{X,U}}(V,\frakZ,v),s_{|V}\right ).
\end{array}
$$
Again, by \ref{Kokolemcocont}, $\alpha_s'$ and $\underline{\alpha}_s$ are continuous and cocontinuous so they induce morphisms of topoi
\begin{equation}\label{alphakraz}
\alpha_s':\widetilde{\calE}(U'/\frakS)  \ra  \widetilde{\calE}(X'/\frakS)_{/\calI'},\ 
\underline{\alpha}_s:\widetilde{\underline{\calE}}(U/\frakS) \ra  \widetilde{\underline{\calE}}(X/\frakS)_{/\underline{\calI}},
\end{equation}
such that the inverse image functors are composition by $\alpha_s'$ and $\underline{\alpha}_s$ respectively.

The morphisms of topoi $\alpha'_s$ and $\underline{\alpha}_s$ can be defined equivalently as follows: let $U$ be an étale $X$-scheme (resp. étale $X'$-scheme) and $s\in \Gamma(U,\calI_X).$ The local section $s$ corresponds to a morphism $s:U^a\ra \calI_X$ and so we have a morphism
$$
\underline{u}_X^{-1}(s):\underline{u}_X^{-1}(U^a) \ra \underline{u}_X^{-1}(\calI_X)=\underline{\calI}\ \left (\text{resp.}\ u_X'^{-1}(s):u_X'^{-1}(U^a) \ra u_X'^{-1}(\calI_X)=\calI'\right ).
$$
By \eqref{stackequiv1} (resp. \eqref{stackequiv2}), it induces, by localization, a morphism of topoi:
$$
\underline{\alpha}_s:\widetilde{\underline{\calE}}(U/\frakS) \ra \widetilde{\underline{\calE}}(X/\frakS)_{/\underline{\calI}}\ \left (\text{resp.}\ \alpha'_s:\widetilde{\calE}(U/\frakS)^{-1}(U^a) \ra \widetilde{\calE}(X'/\frakS)_{/\calI'} \right ).
$$
Finally, if $V \ra U$ is an étale morphism of $X$-schemes then the diagrams
$$
\begin{tikzcd}
\calE(V'/\frakS) \ar{r}{\rho_{U',V'}} \ar[swap]{dr}{\alpha_{s_{|V}}'} & \calE(U'/\frakS), \ar{d}{\alpha_s'} & \underline{\calE}(V/\frakS) \ar{r}{\underline{\rho_{U,V}}} \ar[swap]{dr}{\underline{\alpha}_{s_{|V}}} & \underline{\calE}(U/\frakS) \ar{d}{\underline{\alpha}_s} \\
 & \calE(X'/\frakS)_{/\calI'} & & \underline{\calE}(X/\frakS)_{\underline{\calI}}
\end{tikzcd}
$$
are commutative.
\end{parag}

\begin{remark}\label{remShiho194}
If $\calI_X$ is an étale sheaf on $X,$ we can restrict the sheaves $\calI'$ and $\underline{\calI},$ defined in \ref{paragMuzan194}, to the subsites $\calE_{\mathrm{strict}}(X'/\frakS)$ and $\underline{\calE}_{\mathrm{strict}}(X/\frakS)$ respectively. We denote by $\calI'_{\mathrm{strict}}$ and $\underline{\calI}_{\mathrm{strict}}$ these restrictions. We then have localized topoi
$$
\calE_{\mathrm{strict}}(X'/\frakS)_{/\calI'_{\mathrm{strict}}},\ \underline{\calE}_{\mathrm{strict}}(X/\frakS)_{\underline{\calI}_{\mathrm{strict}}}.
$$
For a local section $s$ of $\calI_X$ over an étale $X$-scheme $U,$ we have functors
\begin{equation}\label{alphakrazShiho}
\underline{\alpha}^{\mathrm{strict}}_s:\widetilde{\underline{\calE}}_{\mathrm{strict}}(U/\frakS) \ra \widetilde{\underline{\calE}}_{\mathrm{strict}}(X/\frakS)_{/\underline{\calI}_{\mathrm{strict}}}\ \left (\text{resp.}\ \alpha'^{\mathrm{strict}}_s:\widetilde{\calE}_{\mathrm{strict}}(U'/\frakS) \ra \widetilde{\calE}_{\mathrm{strict}}(X'/\frakS)_{/\calI'_{\mathrm{strict}}} \right )
\end{equation}
and similar properties to those mentioned in \ref{paragMuzan194}.
\end{remark}

\begin{proposition}\label{indexeddescentdata}
Let $\calI_X$ be a sheaf of $X_{\text{ét}}$ and consider the sheaves $\calI'$ and $\underline{\calI}$ defined in \ref{paragMuzan194}. A presheaf of sets $\F$ on $\calE(X'/\frakS)_{/\calI'}$ (resp. $\underline{\calE}(X/\frakS)_{/\underline{\calI}}$) is equivalent to the data, for every section $s\in \Gamma(U,\calI_X)$ over an étale $X$-scheme $U,$ of a presheaf of sets $\F_s$ on $\calE(U'/\frakS)$ (resp. $\underline{\calE}(U/\frakS)$), such that, for every morphism $g:V\ra U$ of étale $X$-schemes and every section $s\in \Gamma(U,\calI_X),$ we have isomorphisms
$$\rho_{g'}^{-1}\F_s \xrightarrow{\sim} \F_{s_{|V}}\ \left (\text{resp.}\ \underline{\rho_g}^{-1}\F_s \xrightarrow{\sim} \F_{s_{|V}}\right ),$$
satisfying cocycle conditions.
\end{proposition}

\begin{proof}
If $\F$ is a presheaf of sets on $\calE(X'/\frakS)_{/\calI'},$ $U\ra X$ is an étale morphism and $s\in \Gamma(U,\calI_X)$ then we set $\F_s=\F\circ \alpha_s'.$
Conversly, if the presheaves $\F_s$ are given, we define $\F$ by
$$
\F:\begin{array}[t]{clc} 
\calE(X'/\frakS)_{/\calI'} & \ra & \op{Set} \\
(U,\frakT,u,s) & \mapsto & \F_s(\rho_{X',U}(U,\frakT,u)).
\end{array}
$$
If $(U_1,\frakT_1,u_1,s_1) \ra (U_2,\frakT_2,u_2,s_2)$ is a morphism in $\calE(X'/\frakS)_{/\calI'},$ then $s_{2|U_1}=s_1$ and $\rho_{X',U_2}(U_1,\frakT_1,u_1) \ra \rho_{X',U_2}(U_2,\frakT_2,u_2)$ is a morphism in $\calE(U_2/\frakS).$ We thus have a map
$$\F_{s_2}\circ \rho_{X',U_2}(U_2,\frakT_2,u_2) \ra \F_{s_2}\circ \rho_{X',U_2}(U_1,\frakT_1,u_1).$$
Since
$$
\rho_{X',U_2}=\rho_{U_1,U_2}\circ \rho_{X',U_1},
$$
we have
\begin{alignat*}{2}
\F_{s_2} \circ \rho_{X',U_2}=\F_{s_2}\circ \rho_{U_1,U_2}\circ \rho_{X',U_1} &= \left (\rho_{U_1,U_2}^{-1}\F_{s_2}\right )\circ \rho_{X',U_1}\\
&\xrightarrow{\sim}\F_{s_{2|U_1}}\circ \rho_{X',U_1}\\
&=\F_{s_1}\circ \rho_{X',U_1},
\end{alignat*}
we deduce a functorial map
$$\F(U_2,\frakT_2,u_2,s_2) \ra \F(U_1,\frakT_1,u_1,s_1).$$
These two constructions are clearly inverse to each other.
The proof for $\underline{\calE}(X/\frakS)_{/\underline{\calI}}$ is similar.
\end{proof}

\begin{remark}\label{remShiho196}
Proposition \ref{indexeddescentdata} remains true if we replace $\calE(X'/\frakS)_{/\calI'}$ and $\underline{\calE}(X/\frakS)_{\underline{\calI}}$ with $\calE_{\mathrm{strict}}(X'/\frakS)_{/\calI'_{\mathrm{strict}}}$ and $\underline{\calE}_{\mathrm{strict}}(X/\frakS)_{/\underline{\calI}_{\mathrm{strict}}}$ respectively, $\calE(U'/\frakS)$ and $\underline{\calE}(U/\frakS)$ with $\calE_{\mathrm{strict}}(U'/\frakS)$ and $\underline{\calE}_{\mathrm{strict}}(U/\frakS)$ respectively and $\underline{\alpha}_s$ and $\alpha'_s$ with $\underline{\alpha}^{\mathrm{strict}}_s$ and $\alpha'^{\mathrm{strict}}_s$ respectively.
\end{remark}

\begin{definition}\label{inddescentdata}
Let $\calI_X$ be a sheaf of $X_{\text{ét}}$ and consider the sheaves $\calI'$ and $\underline{\calI}$ defined in \ref{paragMuzan194} and the sheaves $\calI'_{\mathrm{strict}}$ and $\underline{\calI}_{\mathrm{strict}}$ defined in \ref{remShiho194}. For a presheaf $\F$ on $\calE(X'/\frakS)_{/\calI'},$ $\underline{\calE}(X/\frakS)_{/\underline{\calI}},$ $\calE_{\mathrm{strict}}(X'/\frakS)_{/\calI'_{\mathrm{strict}}}$ or $\underline{\calE}_{\mathrm{strict}}(X/\frakS)_{/\underline{\calI}_{\mathrm{strict}}},$ we call \emph{descent data associated to $\F$} the presheaves $\F_s$ given in \ref{indexeddescentdata} and \ref{remShiho196}.
\end{definition}

\begin{proposition}\label{iranlem145}
Let $\calI_X$ be a sheaf of $X_{\text{ét}}$ and consider the sheaves $\calI'$ and $\underline{\calI}$ defined in \ref{paragMuzan194}. Let $\F$ be a presheaf on $\calE(X'/\frakS)_{/\calI'}$ (resp. $\underline{\calE}(X/\frakS)_{/\underline{\calI}}$) and $(\F_s)$ the associated descent data \eqref{inddescentdata}.  Then the presheaf $\F$ is a sheaf if and only if $\F_s$ is a sheaf for all $s.$
\end{proposition}

\begin{proof}
This follows from \eqref{stackequiv1}, \eqref{stackequiv2} and the fact that the fibered categories $\mathbbm{T}'$ and $\underline{\mathbbm{T}}$ \eqref{stack1} are stacks.
\end{proof}

\begin{remark}
Proposition \ref{iranlem145} remains true if we replace $\calE(X'/\frakS)_{/\calI'}$ and $\underline{\calE}(X/\frakS)_{\underline{\calI}}$ with $\calE_{\mathrm{strict}}(X'/\frakS)_{/\calI'}$ and $\underline{\calE}_{\mathrm{strict}}(X/\frakS)_{/\underline{\calI}}$ respectively.
\end{remark}

\begin{parag}
Let $\calI_X$ be a sheaf of $X_{\text{ét}}$ and consider the sheaves $\calI'$ and $\underline{\calI}$ defined in \ref{paragMuzan194}.
Let $(U,\frakT,u)$ be an object of $\calE(X'/\frakS)$ (resp. $\underline{\calE}(X/\frakS)$). Consider the morphism of topoi $\alpha_{(U,\frakT,u)}$ \eqref{Kokoalphatopos} and denote by $\calI_U$ the restriction of $\calI_X$ on $U.$ We have
\begin{alignat*}{2}
\alpha_{(U,\frakT,u)}^{-1}\calI' &= \calI' \circ \alpha_{(U,\frakT,u)}=\calI_U\\
(\text{resp.}\ \alpha_{(U,\frakT,u)}^{-1}\underline{\calI} &= \underline{\calI} \circ \alpha_{(U,\frakT,u)}=\calI_U).
\end{alignat*}
The morphism $\alpha_{(U,\frakT,u)}$ then induces a morphism of topoi
\begin{equation}\label{deltatop}
\delta_{(U,\frakT,u),\calI_X}:U_{\text{ét}/\calI_U} \ra \widetilde{\calE}(X'/\frakS)_{/\calI'} \left (\text{resp.\ }\widetilde{\underline{\calE}}(X/\frakS)_{/\underline{\calI}} \right ).
\end{equation}
Its inverse image functor is composition by the following functor, that we abusively denote by $\delta_{(U,\frakT,u),\calI_X}:$
$$
\delta_{(U,\frakT,u),\calI_X}:
\begin{array}[t]{clc}
\text{ét}_{/U\times \calI_X} & \ra & \calE(X'/\frakS)_{/\calI'} \left (\text{resp.\ }\underline{\calE}(X/\frakS)_{/\underline{\calI}}\right )\\
(V,s) & \mapsto & \left (\alpha_{(U,\frakT,u)}(V),s\right ),
\end{array}
$$
where $V$ is an étale $U$-scheme and $s\in \Gamma(V,\calI_X)$ (note that we consider $\calI_X$ as sheaf of $X'_{\text{ét}}$ via $F_1:X\ra X'$) and $\alpha_{(U,\frakT,u)}$ is given in \eqref{eqKoko1861}. If $(U_1,\frakT_1,u_1) \ra (U_2,\frakT_2,u_2)$ is a morphism of $\calE(X'/\frakS)$ (resp. $\underline{\calE}(X/\frakS)$), we consider the functor
$$
j_{g,\calI_X}:\begin{array}[t]{cll}
\text{ét}_{/U_1\times \calI_X} & \ra & \text{ét}_{/U_2\times \calI_X} \\
(V \ra U_1\times \calI_X) & \mapsto & \left (V \ra U_1\times \calI_X \xrightarrow{g\times \op{Id}} U_2\times \calI_X \right ).
\end{array}
$$
If $V$ is an étale $U_1$-scheme then we have a morphism
$$
\beta_{(f,g)}(V):\alpha_{(U_1,\frakT_1,u_1)}(V) \ra \alpha_{(U_2,\frakT_2,u_2)}\circ j_g(V),
$$
where $\beta_{(f,g)}$ is given in \eqref{defbeta115} and $j_g$ in \eqref{jg}.
We thus get a morphism
\begin{equation}\label{defbeta115ind}
\beta_{(f,g),\calI_X}:\delta_{(U_1,\frakT_1,u_1)} \ra \delta_{(U_2,\frakT_2,u_2)} \circ j_{g,\calI_X}.
\end{equation}
For étale morphisms $V\ra U \ra X'$ (resp. $V\ra U \ra X$), an object $(V,\frakZ,v)$ of $\calE(U/\frakS)$ (resp. $\underline{\calE}(U/\frakS)$) and a section $s\in \Gamma(U,\calI_X),$ we have the following commutative diagram:
\begin{equation}\label{diag188}
\begin{tikzcd}
\text{ét}_{/V} \ar{rr}{\alpha_{(V,\frakZ,v)}} \ar[swap]{dd}{W\mapsto (W,s_{|W})} & & \calE(U/\frakS) \left (\text{resp.\ }\underline{\calE}(U/\frakS)\right )\ar{dd}{\alpha_s'\left (\text{resp.\ }\underline{\alpha}_s\right )} \\
 & & \\
\text{ét}_{/V\times \calI_X} \ar[swap]{rr}{\delta_{(V,\frakZ,v),\calI_X}} & & \calE(X'/\frakS)_{/\calI'} \left (\text{resp.\ }\underline{\calE}(X/\frakS)_{/\underline{\calI}}\right ).
\end{tikzcd}
\end{equation}
For a positive integer $n$ and $\calI_X=\upmu_X^n,$ we will denote $\delta_{(U,\frakT,u),\calI_X}$ simply by $\delta_{(U,\frakT,u),n}.$ If, furthermore $n=1,$ we denote $\delta_{(U,\frakT,u),n}$ simply by $\delta_{(U,\frakT,u)}.$
\end{parag}

\begin{remark}
Let $\calI_X$ be an étale sheaf on $X$ and consider the sheaves $\calI'$ and $\underline{\calI}$ defined in \ref{paragMuzan194} and the sheaves $\calI'_{\mathrm{strict}}$ and $\underline{\calI}_{\mathrm{strict}}$ defined in \ref{remShiho194}. If $(U,\frakT,u)$ is an object of $\calE_{\mathrm{strict}}(X'/\frakS)$ (resp. $\underline{\calE}_{\mathrm{strict}}(X/\frakS)$), then, by \ref{remShiho1}, the functor \eqref{deltatop} factors through $\widetilde{\calE}_{\mathrm{strict}}(X'/\frakS)_{/\calI'_{\mathrm{strict}}}$ (resp. $\widetilde{\underline{\calE}}_{\mathrm{strict}}(X/\frakS)_{/\underline{\calI}_{\mathrm{strict}}})$, yielding a morphism of topoi, abusively denoted
\begin{equation}
\delta_{(U,\frakT,u),\calI_X}^{\mathrm{strict}}:U_{\text{ét}/\calI_U} \ra \widetilde{\calE}_{\mathrm{strict}}(X'/\frakS)_{/\calI'_{\mathrm{strict}}} \left (\text{resp.\ }\widetilde{\underline{\calE}}_{\mathrm{strict}}(X/\frakS)_{/\underline{\calI}_{\mathrm{strict}}} \right ).
\end{equation}
\end{remark}

\begin{proposition}\label{inddescentdata2}
Let $\calI_X$ be an étale sheaf on $X$ an consider it as an étale sheaf on $X'$ via the exact relative Frobenius $X\ra X'.$
For an étale $X$-scheme or an étale $X'$-scheme $U,$ we denote by $\calI_U$ the restriction of $\calI_X$ on $U.$ A presheaf $\F$ on $\calE(X'/\frakS)_{/\calI'}$ (resp. $\underline{\calE}(X/\frakS)_{/\underline{\calI}}$) is equivalent to the following data:
\begin{enumerate}
\item For every object $(U,\frakT,u)$ of $\calE(X'/\frakS)$ (resp. $\underline{\calE}(X/\frakS)$), a presheaf $\F_{(U,\frakT,u)}$ on $\text{ét}_{/U\times \calI_X}.$
\item For every morphism $(f,g):(U_1,\frakT_1,u_1) \ra (U_2,\frakT_2,u_2)$ in $\calE(X'/\frakS)$ (resp. $\underline{\calE}(X/\frakS)$), a morphism of presheaves on $\text{ét}_{/U_1\times \calI_X},$
$$\gamma_{\F,(f,g)}:g^{-1}\F_{(U_2,\frakT_2,u_2)} \ra \F_{(U_1,\frakT_1,u_1)},$$
\end{enumerate}
where $g$ denotes abusively the morphism of topoi $U_{1,\text{ét}/\calI_{U_1}} \ra U_{2,\text{ét}/\calI_{U_2}}$
induced by $g:U_1\ra U_2,$
such that
\begin{enumerate}[(i)]
\item For every object $(U,\frakT,u)$ of $\calE(X/\frakS)$ (resp. $\underline{\calE}(X/\frakS)$),
$$\gamma_{\F,\op{Id}_{(U,\frakT,u)}}=\op{Id}_{\F_{(U,\frakT,u)}}.$$
\item For all composable morphisms $(f_1,g_1):(U_1,\frakT_1,u_1) \ra (U_2,\frakT_2,u_2)$ and $(f_2,g_2):(U_2,\frakT_2,u_2) \ra (U_3,\frakT_3,u_3)$ in $\calE(X'/\frakS)$ (resp. $\underline{\calE}(X/\frakS)$),
$$\gamma_{\F,(f_1,g_1)} \circ g_1^{-1}\gamma_{\F,(f_2,g_2)}=\gamma_{\F,(f_2,g_2)\circ (f_1,g_1)}.$$
\item If $(f,g)$ is a cartesian morphism then $\gamma_{\F,(f,g)}$ is an isomorphism.
\end{enumerate}
This equivalence is given as follows: for a presheaf $\F$ on $\calE(X'/\frakS)_{/\calI'}$ (resp. $\underline{\calE}(X/\frakS)_{/\underline{\calI}}$),
$$\F_{(U,\frakT,u)}=\F \circ \delta_{(U,\frakT,u),\calI_X}$$
and $\gamma_{\F,(f,g)}$ is induced by $\beta_{(f,g),\calI_X}$ \eqref{defbeta115ind}. Conversely, given a data $(\F_{(U,\frakT,u)},\gamma_{\F,(f,g)})$ as above, the presheaf is defined by
$$\F(U,\frakT,u,s)=\F_{(U,\frakT,u)}(U,s)$$
and, for a morphism $(f,g):(U_1,\frakT_1,u_1,s_1) \ra (U_2,\frakT_2,u_2,s_2),$ the morphism
$$\F(U_2,\frakT_2,u_2,s_2) \ra \F(U_1,\frakT_1,u_1,s_1)$$
is equal to $\gamma_{\F,(f,g)}(U_1,s_1).$
\end{proposition}

\begin{proof}
Similar to \ref{era3propdescentdata}.
\end{proof}

\begin{remark}
Proposition \ref{inddescentdata2} remains true if we replace $\calE(X'/\frakS)_{/\calI'}$ and $\underline{\calE}(X/\frakS)_{\underline{\calI}}$ with $\calE_{\mathrm{strict}}(X'/\frakS)_{/\calI'_{\mathrm{strict}}}$ and $\underline{\calE}_{\mathrm{strict}}(X/\frakS)_{/\underline{\calI}_{\mathrm{strict}}}$ respectively and if we replace $\calE(X'/\frakS)$ and $\underline{\calE}(X/\frakS)$ with $\calE_{\mathrm{strict}}(X'/\frakS)$ and $\underline{\calE}_{\mathrm{strict}}(X/\frakS)$ respectively.
\end{remark}

\begin{definition}\label{definddescentdata}
For a presheaf of sets $\F$ on $\calE(X'/\frakS)_{/\calI'}$ (resp. $\underline{\calE}(X/\frakS)_{/\underline{\calI}},$ $\calE_{\mathrm{strict}}(X'/\frakS)_{/\calI'_{\mathrm{strict}}},$ $\underline{\calE}_{\mathrm{strict}}(X/\frakS)_{/\underline{\calI}_{\mathrm{strict}}}$), we call \emph{indexed descent data associated with $\F$} the data $(\F_{(U,\frakT,u)},\gamma_{\F,(f,g)})$ given in \ref{inddescentdata2}.
\end{definition}

\begin{proposition}\label{prop1409}
Let $\F$ be a presheaf of sets on $\calE(X'/\frakS)_{/\calI'}$ (resp. $\underline{\calE}(X/\frakS)_{/\underline{\calI}},$ $\calE_{\mathrm{strict}}(X'/\frakS)_{/\calI'},$ $\underline{\calE}_{\mathrm{strict}}(X/\frakS)_{/\underline{\calI}}$) and $(\F_{(U,\frakT,u)},\gamma_{\F,(f,g)})$ the indexed descent data associated with $\F$ \eqref{definddescentdata}. Then $\F$ is a sheaf for the étale topology if and only if $\F_{(U,\frakT,u)}$ is a sheaf of $U_{\text{ét}/\calI_U}$ for every $(U,\frakT,u).$
\end{proposition}

\begin{proof}
The proof is similar to \ref{era4prop1723}.
\end{proof}

\begin{parag}\label{kraz1911}
Consider the localization morphisms of topoi
$$j_{\calI'}:\widetilde{\calE}(X'/\frakS)_{/\calI'} \ra \widetilde{\calE}(X'/\frakS),\ j_{\underline{\calI}}:\widetilde{\underline{\calE}}(X/\frakS)_{/\underline{\calI}} \ra \widetilde{\underline{\calE}}(X/\frakS),$$
$$
j_{\calI'_{\mathrm{strict}}}:\widetilde{\calE}_{\mathrm{strict}}(X'/\frakS)_{/\calI'_{\mathrm{strict}}} \ra \widetilde{\calE}_{\mathrm{strict}}(X'/\frakS),\ j_{\underline{\calI}_{\mathrm{strict}}}:\widetilde{\underline{\calE}}_{\mathrm{strict}}(X/\frakS)_{/\underline{\calI}_{\mathrm{strict}}} \ra \widetilde{\underline{\calE}}_{\mathrm{strict}}(X/\frakS)
$$
given in \ref{indeq2}.
We equip the topos $\widetilde{\calE}(X'/\frakS)_{/\calI'}$ (resp. $\widetilde{\underline{\calE}}(X/\frakS)_{/\underline{\calI}}$) with the ring $\Ox_{\calI'}=j_{\calI'}^{-1}\Ox_{\calE(X'/\frakS)}$ (resp. $\Ox_{\underline{\calI}}=j_{\underline{\calI}}^{-1}\Ox_{\underline{\calE}(X/\frakS)}$). Notice that, if $s\in \Gamma(U,\calI_X)$ is a local section over an étale $X'$-scheme (resp. $X$-scheme) $U,$ then
$$\alpha_s'^{-1}\Ox_{\calI'}=\rho_{X',U}^{-1}\Ox_{\calE(X'/\frakS)}=\Ox_{\calE(U/\frakS)},$$
$$\left (\text{resp.}\ \underline{\alpha}_s^{-1}\Ox_{\underline{\calI}}=\underline{\rho_{X,U}^{-1}}\Ox_{\underline{\calE}(X/\frakS)}=\Ox_{\underline{\calE}(U/\frakS)}\right ).$$
If $\F$ is a module of $\widetilde{\calE}(X'/\frakS)_{/\calI'}$ (resp. $\widetilde{\underline{\calE}}(X/\frakS)_{/\underline{\calI}}$) then $\alpha_s'^{-1}\F$ (resp. $\underline{\alpha}_s^{-1}\F$) is a module on $\calE(U/\frakS)$ (resp. $\underline{\calE}(U/\frakS)$).

We similarly define the rings $\Ox_{\calI'_{\mathrm{strict}}}$ and $\Ox_{\underline{\calI}_{\mathrm{strict}}}.$ We then have similar properties to those given above.
\end{parag}

\begin{lemma}
Let $(U,\frakT,u)$ be an object of $\calE(X'/\frakS)$ (resp. $\underline{\calE}(X/\frakS)$). Consider the morphism of topoi $\delta_{(U,\frakT,u),\calI_X}$ \eqref{deltatop}. Set $\Ox=\Ox_{\calE(X'/\frakS)}$ (resp. $\Ox=\Ox_{\underline{\calE}(X/\frakS)}$) and $j=j_{\calI'}$ (resp $j=j_{\underline{\calI}}$). Then
$$
\delta_{(U,\frakT,u),\calI_X}^{-1}j^{-1}\Ox=\left (u_*\Ox_T\right )_{/\calI_X}.
$$
\end{lemma}

\begin{proof}
For $\calE(X'/\frakS),$ this follows from \eqref{eqKoko18171} and the commutative diagram
$$
\begin{tikzcd}
\text{ét}_{/U\times \calI_X} \ar{rr}{\delta_{(U,\frakT,u),\calI_X}} \ar[swap]{d}{j_{\calI_X}} & & \calE(X'/\frakS)_{/\calI'} \ar{d}{j_{\calI'}} \\
\text{ét}_{/U} \ar[swap]{rr}{\alpha_{(U,\frakT,u)}} & & \calE(X'/\frakS).
\end{tikzcd}
$$
The proof for $\underline{\calE}(X/\frakS)$ is similar.
\end{proof}

\begin{parag}
Let $(f,g):(U_1,\frakT_1,u_1) \ra (U_2,\frakT_2,u_2)$ be a morphism in $\calE(X'/\frakS)$ (resp. $\underline{\calE}(X/\frakS)$). We have a morphism of ringed topoi \eqref{Kokoftilda}
\begin{equation}
\widetilde{f}:\left (U_{1,\text{ét}},u_{1*}\Ox_{T_{1}} \right ) \ra \left ( U_{2,\text{ét}}, u_{2*}\Ox_{T_{2}} \right ).
\end{equation}
We denote by $\calI_{U_i}$ the restriction of $\calI_X$ to $U_i.$ We have the following commutative diagram of topoi
$$
\begin{tikzcd}
U_{1,\text{ét}/\calI_{U_1}} \ar{r}{g_{\calI}} \ar[swap]{d}{j_{\calI_{U_1}}} & U_{2,\text{ét}/\calI_{U_2}} \ar{d}{j_{\calI_{U_2}}} \\
U_{1,\text{ét}} \ar{r}{g} & U_{2,\text{ét}}, 
\end{tikzcd}
$$
where the upper arrow is induced by $g.$ We then have
$$
\left (g^{-1}u_{2*}\Ox_{T_2} \right )_{\calI_{U_1}} =g_{\calI}^{-1}\left ( \left (u_{2*}\Ox_{T_2} \right )_{\calI_{U_2}} \right ).
$$
It follows that the canonical morphism
$$
g^{-1}u_{2*}\Ox_{T_2}=\widetilde{f}^{-1}u_{2*}\Ox_{T_2} \ra u_{1*} \Ox_{T_1}
$$
induces
$$
g_{\calI}^{-1}\left ( \left (u_{2*}\Ox_{T_2} \right )_{\calI_{U_2}} \right ) \ra \left ( u_{1*}\Ox_{T_1} \right )_{\calI_{U_1}}.
$$
We then get a morphism of ringed topoi
\begin{equation}\label{gmu}
\widetilde{f}_{\calI}: \left (U_{1,\text{ét}/\calI_{U_1}},\left ( u_{1*}\Ox_{T_1} \right )_{\calI_{U_1}} \right ) \ra \left (U_{2,\text{ét}/\calI_{U_2}},  \left (u_{2*}\Ox_{T_2} \right )_{\calI_{U_2}}\right ).
\end{equation}
\end{parag}

\begin{proposition}\label{indlindescentdata2}
Let $\calI_X$ be an étale sheaf on $X$ an consider it as an étale sheaf on $X'$ via the exact relative Frobenius $X\ra X'.$ For an étale $X$-scheme $U,$ we denote by $\calI_U$ the restriction of $\calI_X$ on $U.$ A module $\F$ on $\calE(X'/\frakS)_{/\calI'}$ (resp. $\underline{\calE}(X/\frakS)_{/\underline{\calI}}$) is equivalent to the following data:
\begin{enumerate}
\item For every object $(U,\frakT,u)$ of $\calE(X'/\frakS)$ (resp. $\underline{\calE}(X/\frakS)$), a module $\F_{(U,\frakT,u)}$ of
$$\left (U_{\text{ét}/ \calI_U},\left (u_*\Ox_T \right )_{\calI_U}\right ).$$
\item For every morphism $(f,g):(U_1,\frakT_1,u_1) \ra (U_2,\frakT_2,u_2)$ in $\calE(X'/\frakS)$ (resp. $\underline{\calE}(X/\frakS)$), a morphism of modules of $U_{1,\text{ét}/ \calI_{U_1}},$
\begin{equation}\label{cind}
c_{\F,(f,g)}:\widetilde{f}_{\calI}^*\F_{(U_2,\frakT_2,u_2)} \ra \F_{(U_1,\frakT_1,u_1)},
\end{equation}
\end{enumerate}
where $\widetilde{f}_{\calI}$ is given in \eqref{gmu}.
such that
\begin{enumerate}[(i)]
\item For every object $(U,\frakT,u)$ of $\calE(X/\frakS)$ (resp. $\underline{\calE}(X/\frakS)$),
$$c_{\F,\op{Id}_{(U,\frakT,u)}}=\op{Id}_{\F_{(U,\frakT,u)}}.$$
\item For all composable morphisms $(f_1,g_1):(U_1,\frakT_1,u_1) \ra (U_2,\frakT_2,u_2)$ and $(f_2,g_2):(U_2,\frakT_2,u_2) \ra (U_3,\frakT_3,u_3)$ in $\calE(X'/\frakS)$ (resp. $\underline{\calE}(X/\frakS)$),
$$c_{\F,(f_1,g_1)} \circ \widetilde{f}_{1,\calI}^{*}c_{\F,(f_2,g_2)}=c_{\F,(f_2,g_2)\circ (f_1,g_1)}.$$
\item If $(f,g)$ is a cartesian morphism then $c_{\F,(f,g)}$ is an isomorphism.
\end{enumerate}
This equivalence is given as follows: for a module $\F$ on $\calE(X'/\frakS)_{/\calI'}$ (resp. $\underline{\calE}(X/\frakS)_{/\underline{\calI}}$),
$$\F_{(U,\frakT,u)}=\F \circ \delta_{(U,\frakT,u),\calI_X}$$
and $c_{\F,(f,g)}$ is induced by $\gamma_{\F,(f,g)}$ \eqref{inddescentdata2}. Conversely, given a data $(\F_{(U,\frakT,u)},c_{\F,(f,g)})$ as above, the module is defined by
$$\F(U,\frakT,u,s)=\F_{(U,\frakT,u)}(U,s)$$
and, for a morphism $(f,g):(U_1,\frakT_1,u_1,s_1) \ra (U_2,\frakT_2,u_2,s_2),$ the morphism
$$\F(U_2,\frakT_2,u_2,s_2) \ra \F(U_1,\frakT_1,u_1,s_1)$$
is equal to $c_{\F,(f,g)}(U_1,s_1).$
\end{proposition}

\begin{proof}
The proof is similar to \ref{lindescentdata}.
\end{proof}

\begin{remark}\label{remShiho1920}
Proposition \ref{indlindescentdata2} remains true if we replace $\calE(X'/\frakS),$ $\underline{\calE}(X/\frakS),$ $\calI',$ $\underline{\calI}$ and $\delta_{(U,\frakT,u),\calI_X}$ with $\calE_{\mathrm{strict}}(X'/\frakS),$ $\underline{\calE}_{\mathrm{strict}}(X/\frakS),$ $\calI'_{\mathrm{strict}},$ $\underline{\calI}_{\mathrm{strict}}$ and $\delta_{(U,\frakT,u),\calI_X}^{\mathrm{strict}}$ respectively.
\end{remark}

\begin{definition}\label{defindlindescentdata}
For a module $\F$ on $\calE(X'/\frakS)_{/\calI'}$ (resp. $\underline{\calE}(X/\frakS)_{/\underline{\calI}}$), we call \emph{indexed linearized descent data associated with $\F$} the data $(\F_{(U,\frakT,u)},c_{\F,(U,\frakT,u)})$ defined in \ref{indlindescentdata2}.
\end{definition}

\begin{parag}
Let $p_1,p_2:\calI_X^2 \ra \calI_X$ be the canonical projections. These projections induce morphisms of ringed topoi, that we abusively denote by
$$
p_1,p_2:\widetilde{\calE}(X'/\frakS)_{/\calI'^2} \ra \widetilde{\calE}(X'/\frakS)_{/\calI'} \left (\text{resp.}\ \widetilde{\underline{\calE}}(X/\frakS)_{/\underline{\calI}^2} \ra \widetilde{\underline{\calE}}(X/\frakS)_{/\underline{\calI}}\right ).
$$
For modules $\F$ and $\calG$ of $\widetilde{\calE}(X'/\frakS)_{/\calI'}$ (resp. $\widetilde{\underline{\calE}}(X/\frakS)_{/\underline{\calI}}$), we denote by $\F\boxtimes \calG$ the exterior tensor product
$$
\F\boxtimes \calG=p_1^*\F \otimes_{\Ox_{\calI'^2}} p_2^*\calG \left (\text{resp.}\ \F\boxtimes \calG=p_1^*\F \otimes_{\Ox_{\underline{\calI}^2}} p_2^*\calG\right ).
$$
\end{parag}

\begin{proposition}\label{prop1916Ko}
Let $\calG$ and $\F$ be two $\Ox_{\calI'}$-modules of $\widetilde{\calE}(X'/\frakS)_{/\calI'}$ (resp. two $\Ox_{\underline{\calI}}$-modules of $\widetilde{\underline{\calE}}(X/\frakS)_{/\underline{\calI}}$) and $(\calG_s),$ $(\F_s)$ and $\left (\calG\boxtimes\F \right )_{t}$ (resp. $\left (\calG\boxtimes\F \right )_{t}$) the corresponding descent data \eqref{inddescentdata}. Then, for any local sections $s,t\in \Gamma(U,\calI_X)$ over an étale $X'$-scheme (resp. $X$-scheme) $U,$
$$\left (\calG\boxtimes\F \right )_{(s,t)}=\calG_s\otimes_{\Ox_{\calE(U/\frakS)}}\F_t$$
$$\left (\text{resp.}\ \left (\calG\boxtimes\F \right )_{(s,t)}=\calG_s\otimes_{\Ox_{\underline{\calE}(U/\frakS)}}\F_t\right ).$$
\end{proposition}

\begin{proof}
Let $p_1,p_2:\calI'^2 \ra \calI'$ be the canonical projections. The proposition results from the commutativity of the following diagram, for any étale $X'$-scheme $U$ and $s_1,s_2\in \Gamma(U,\calI_X):$
$$
\begin{tikzcd}
\widetilde{\calE}(U/\frakS) \ar{rr}{\alpha_{(s_1,s_2)}'} \ar[swap]{drr}{\alpha_{s_i}'} & & \widetilde{\calE}(X'/\frakS)_{/\calI'^2} \ar{d}{p_i}  \\
 & & \widetilde{\calE}(X'/\frakS)_{/\calI'} 
\end{tikzcd}
$$
and from a similar diagram for $\underline{\calE}(X/\frakS).$
\end{proof}

\begin{parag}
Suppose that $\calI_X$ is a sheaf of monoids. Let $p_1,p_2:\calI'^2\ra \calI'$ and $p_{1,X},p_{2,X}:\calI_X^2 \ra \calI_X$ be the canonical projections and $\sigma:\calI'^2\ra \calI'$ and $\sigma_X:\calI_X^2\ra \calI_X$ the addition maps. Then
$$
u_X'^{-1}p_i=p_{i,X},\ u_X'^{-1}\sigma=\sigma_X,
$$
where $u_X'$ is given in \eqref{Kokou}.
\end{parag}

\begin{proposition}\label{prop1415}
Let $\F$ be a presheaf of sets on $\calE(X'/\frakS)_{/\calI'}$ (resp. $\underline{\calE}(X/\frakS)_{/\underline{\calI}}$) and $$(\F_{(U,\frakT,u)},\gamma_{\F,(f,g)})$$ the indexed descent data associated with $\F$ \eqref{definddescentdata}. Then, if $\F$ has a $\calI'$-indexed (resp. $\underline{\calI}$-indexed) algebra structure, the sheaves $\F_{(U,\frakT,u)}$ have $\calI_U$-indexed algebra structures for any object $(U,\frakT,u)$ of $\calE(X'/\frakS)_{/\calI'}$ (resp. $\underline{\calE}(X/\frakS)_{/\underline{\calI}}$). The same is true for indexed modules.
\end{proposition}

\begin{proof}
Let $p_1,p_2:\calI'^2\ra \calI'$ and $p_{1,X},p_{2,X}:\calI_X^2 \ra \calI_X$ be the canonical projections and $\sigma:\calI'^2\ra \calI'$ and $\sigma_X:\calI_X^2\ra \calI_X$ the addition maps.
For any object $(U,\frakT,u)$ of $\calE(X'/\frakS),$ we have the following commutative diagram
$$
\begin{tikzcd}
U_{\text{ét}/\calI_U^2} \ar{rr}{\delta_{(U,\frakT,u),\calI_X^2}} \ar[swap]{d}{p_{i,U}\ (\text{resp.}\ \sigma_U)} & & \widetilde{\calE}(X'/\frakS)_{/\calI'^2} \ar{d}{p_{i}\ (\text{resp.}\ \sigma)} \\
U_{\text{ét}/\calI_U} \ar[swap]{rr}{\delta_{(U,\frakT,u),\calI_X}} & & \widetilde{\calE}(X'/\frakS)_{/\calI'}
\end{tikzcd}
$$
This implies that, if $\F$ has a $\calI'$-indexed algebra structure, then $\F_{(U,\frakT,u)}=\delta_{(U,\frakT,u),\calI_X}^*\F$ has a $\calI_U$-indexed algebra structure. The proof for $\underline{\calE}(X/\frakS)_{/\underline{\calI}}$ and for indexed modules is similar.
\end{proof}

\begin{proposition}\label{propKoko1915}
Let $\calA$ be a $\calI'$-indexed algebra on $\calE(X'/\frakS)_{/\calI'}$ and $\calE$ and $\F$ two $\calI'$-indexed $\calA$-modules \eqref{inddef}. For any object $(U,\frakT,u)$ of $\calE(X'/\frakS),$ there exists a canonical isomorphism
$$\calE_{(U,\frakT,u)}\circledast_{\calA_{(U,\frakT,u)}}\F_{(U,\frakT,u)} \xrightarrow{\sim} \left (\calE \circledast_{\calA}\F \right )_{(U,\frakT,u)}.$$
\end{proposition}

\begin{proof}
Let $p_1,p_2:\calI'^2\ra \calI'$ and $p_{1,X},p_{2,X}:\calI_X^2 \ra \calI_X$ be the canonical projections and $\sigma:\calI'^2 \ra \calI'$ and $\sigma_X:\calI_X^2\ra \calI_X$ the addition maps.
Applying \ref{lemtop15} to
$$
\begin{tikzcd}
U_{\text{ét}/\calI_U^2} \ar{rr}{\delta_{(U,\frakT,u),\calI_X^2}} \ar[swap]{d}{\sigma_U} & & \widetilde{\calE}(X'/\frakS)_{/\calI'^2} \ar{d}{\sigma} \\
U_{\text{ét}/\calI_U} \ar[swap]{rr}{\delta_{(U,\frakT,u),\calI_X}} & & \widetilde{\calE}(X'/\frakS)_{/\calI'},
\end{tikzcd}
$$
we get a canonical isomorphism
$$\sigma_{U!}\delta_{(U,\frakT,u),\calI_X^2}^*\left (\calE\boxtimes \F \right ) \xrightarrow{\sim} \delta_{(U,\frakT,u),\calI_X}^*\sigma_!\left (\calE\boxtimes \F \right ).$$
Then, by the commutative diagram
$$
\begin{tikzcd}
U_{\text{ét}/\calI_U^2} \ar{rr}{\delta_{(U,\frakT,u),\calI_X^2}} \ar[swap]{d}{p_{i,X}} & & \widetilde{\calE}(X'/\frakS)_{/\calI'^2} \ar{d}{p_i} \\
U_{\text{ét}/\calI_U} \ar[swap]{rr}{\delta_{(U,\frakT,u),\calI_X}} & & \widetilde{\calE}(X'/\frakS)_{/\calI'},
\end{tikzcd}
$$
we get a canonical isomorphism
\begin{equation}\label{feq012}
\sigma_{U!}\left (\calE_{(U,\frakT,u)}\boxtimes \F_{(U,\frakT,u)} \right ) \xrightarrow{\sim} \left (\sigma_!\left (\calE\boxtimes \F \right )\right )_{(U,\frakT,u)}.
\end{equation}
Let $\tau:\calI^3 \ra \calI$ and $\tau_X:\calI_X^3\ra \calI_X$ be the addition maps.
Recall that, by \eqref{exseqindtens}, we have an exact sequence of $\Ox_{\calI}$-modules
$$\tau_!\left (\calA \boxtimes \calE\boxtimes \F \right )\ra \sigma_!\left (\calE \boxtimes \F\right ) \ra \calE \circledast_{\calA}\F \ra 0.$$
Since $\delta_{(U,\frakT,u)}^*$ is right exact, we obtain the exact sequence
$$\left (\tau_!\left (\calA \boxtimes \calE\boxtimes \F \right )\right )_{(U,\frakT,u)}\ra \left (\sigma_!\left (\calE \boxtimes \F\right )\right )_{(U,\frakT,u)} \ra \left (\calE \circledast_{\calA}\F\right )_{(U,\frakT,u)} \ra 0.$$
For a sheaf $\calG$ of $\widetilde{\calE}(X'/\frakS)_{/\calI'},$ we denote $\calG_{(U,\frakT,u)}$ simply by $\calG_U.$
By \eqref{feq012} and a similar isomorphism for $\tau,$ we get a commutative diagram
$$
\begin{tikzcd}
\tau_{U!}\left (\calA_U \boxtimes \calE_U\boxtimes \F_U \right )\ar{r} \ar[sloped]{d}{\sim} & \sigma_{U!}\left (\calE_U \boxtimes \F_U\right ) \ar{r} \ar[sloped]{d}{\sim} & \calE_U  \circledast_{\calA_U }\F_U \ar{r} \ar[dashed]{d} & 0 \\
\left (\tau_!\left (\calA \boxtimes \calE\boxtimes \F \right )\right )_{(U,\frakT,u)}\ar{r} & \left (\sigma_!\left (\calE \boxtimes \F\right )\right )_{(U,\frakT,u)} \ar{r} & \left (\calE \circledast_{\calA}\F\right )_{(U,\frakT,u)} \ar{r} & 0
\end{tikzcd}
$$
with exact lines.
Hence the existence of the dashed arrow and the fact that it is an isomorphism.
\end{proof}

\subsection*{Crystals and the Cartier transform}

\begin{definition}
Let $\F$ be an $\Ox_{\calI'}$-module of $\widetilde{\calE}(X'/\frakS)_{/\calI'}$ (resp. an $\Ox_{\underline{\calI}}$-module of $\widetilde{\underline{\calE}}(X/\frakS)_{/\underline{\calI}},$ an $\Ox_{\calI'_{\mathrm{strict}}}$-module of $\widetilde{\calE}_{\mathrm{strict}}(X'/\frakS)_{/\calI'_{\mathrm{strict}}},$ an $\Ox_{\underline{\calI}_{\mathrm{strict}}}$-module of $\widetilde{\underline{\calE}}_{\mathrm{strict}}(X/\frakS)_{/\underline{\calI}_{\mathrm{strict}}}$) \eqref{kraz1911} and $(\F_s)$ the corresponding descent data \eqref{inddescentdata}. Recall that for every local section $s\in \Gamma(U,\calI_X)$ over an étale $X'$-scheme (resp. $X$-scheme) $U,$ the sheaf $\F_s$ is a module on $\calE(U/\frakS)$ (resp. $\underline{\calE}(U/\frakS),$ $\calE_{\mathrm{strict}}(U'/\frakS),$ $\underline{\calE}_{\mathrm{strict}}(U/\frakS)$). We say that $\F$ is \emph{quasi-coherent} (resp. a \emph{crystal}) if $\F_s$ is quasi-coherent (resp. a crystal) \eqref{defcrys} for every local section $s$ of $\calI_X.$
\end{definition}

\begin{proposition}\label{equivShihoind}
The restriction functor
$$
\widetilde{\calE}(X'/\frakS) \ra \widetilde{\calE}_{\mathrm{strict}}(X'/\frakS)
$$
$$
\left ( \text{resp.}\ \widetilde{\underline{\calE}}(X/\frakS) \ra \widetilde{\underline{\calE}}_{\mathrm{strict}}(X/\frakS) \right )
$$
induces an equivalence between the subcategories of crystals
$$
\begin{Bmatrix}\text{Crystals\ of\ }\widetilde{\calE}(X'/\frakS)_{/\calI'} \end{Bmatrix} \xrightarrow{\sim} \begin{Bmatrix}\text{Crystals\ of\ }\widetilde{\calE}_{\mathrm{strict}}(X'/\frakS)_{/\calI'_{\mathrm{strict}}} \end{Bmatrix}
$$
$$
\left (\text{resp.}\ \begin{Bmatrix}\text{Crystals\ of\ }\widetilde{\underline{\calE}}(X/\frakS)_{/\underline{\calI}} \end{Bmatrix} \xrightarrow{\sim} \begin{Bmatrix}\text{Crystals\ of\ }\widetilde{\underline{\calE}}_{\mathrm{strict}}(X/\frakS)_{/\underline{\calI}_{\mathrm{strict}}} \end{Bmatrix}
\right )
$$
\end{proposition}

\begin{proof}
Note that, for any local section $s\in \Gamma \left (U,\calI_X \right )$ over an étale $X$-scheme $U,$ we clearly have the commutative diagram
$$
\begin{tikzcd}
\widetilde{\calE}(X'/\frakS)_{/\calI'} \ar[swap]{d}{\alpha_s'^{-1}} \ar{r}{R} & \widetilde{\calE}_{\mathrm{strict}}(X'/\frakS)_{/\calI'_{\mathrm{strict}}} \ar{d}{\left (\alpha_s'^{\mathrm{strict}} \right )^{-1}} \\
\widetilde{\calE}(U' /\frakS) \ar{r}{R_U}  & \widetilde{\calE}_{\mathrm{strict}}(U'/\frakS),
\end{tikzcd}
$$
where the horizontal arrows are the restriction functors and the vertical ones are the inverse image functors of the morphisms of topoi \eqref{alphakraz} and \eqref{alphakrazShiho}. It follows that, if $\F$ is a crystal of $\widetilde{\calE}(X'/\frakS)_{/\calI'},$ then
$$
R(\F)_s=R_U(\F_s).
$$
The module $R(\F)_s$ is then a crystal for all local sections $s$ and hence $R(\F)$ is a crystal. The functor $R$ then induces
$$
R:\begin{Bmatrix}\text{Crystals\ of\ }\widetilde{\calE}(X'/\frakS)_{/\calI'} \end{Bmatrix} \ra \begin{Bmatrix}\text{Crystals\ of\ }\widetilde{\calE}_{\mathrm{strict}}(X'/\frakS)_{/\calI'_{\mathrm{strict}}} \end{Bmatrix}.
$$
The fact that it is an equivalence follows from \ref{equivcrysShiho}. The proof for $\underline{\calE}(X/\frakS)$ is similar.
\end{proof}

\begin{parag}\label{AcrysShiho}
We take the sheaf $\upmu_X=\ov{\calM}_X^{gp}$ and define the sheaf $\underline{\upmu}$ on $\underline{\calE}(X/\frakS)$ from $\upmu_X$ as in \ref{paragMuzan194} i.e. $\underline{\upmu}=\underline{u}_X^{-1}\upmu_X$ where $\underline{u}_X$ is given in \eqref{Kokou}. We also consider the restriction $\underline{\upmu}_{\mathrm{strict}}$ of $\underline{\upmu}$ to the subsite $\underline{\calE}_{\mathrm{strict}}(X/\frakS).$

For an object $(U,\frakT,u)$ of $\underline{\calE}_{\mathrm{strict}}(X/\frakS),$ we identify the étale sites on $T$ and $\underline{T}$ via the canonical nilpotent immersion $\underline{T} \ra T.$ The strict morphism $u:\underline{T} \ra U$ and the canonical immersion $\underline{T} \ra T,$ which is strict, induce an isomorphism
$$
u^{-1}\upmu_U=u^{-1}\ov{\calM}_U^{gp} \xrightarrow{\sim} \ov{\calM}_{\underline{T}}^{gp} \xrightarrow{\sim} \ov{\calM}_T^{gp}=\upmu_T.
$$
We then have a morphism of ringed topoi, abusively denoted
$$
u:\left (T_{\text{ét}/\upmu_T},\Ox_{T,\upmu_T} \right ) \xrightarrow{\sim} \left (\underline{T}_{\text{ét}/\upmu_{T}},\Ox_{T,\upmu_T} \right ) \ra \left (U_{\text{ét}/\upmu_U},\left (u_*\Ox_T\right )_{\upmu_U} \right ).
$$
We set
\begin{equation}\label{Pandaeq2021}
\frakA_{X,(U,\frakT,u)}=u_*\calA_{T}.
\end{equation}
Let
$$
(f,g):(U_1,\frakT_1,u_1) \ra (U_2,\frakT_2,u_2)
$$
be a morphism in $\underline{\calE}_{\mathrm{strict}}(X/\frakS).$ The strict morphism $f:\frakT_1 \ra \frakT_2$ induces a strict morphism $f_1:T_1 \ra T_2,$ which yields, by (\cite{Lor2000} I 3.1), an isomorphism
$$
f_1^*\calA_{T_2} \xrightarrow{\sim} \calA_{T_1}.
$$
This is equivalent to an isomorphism
\begin{equation}\label{Pandaeq2022}
c_{\frakA_X,(f,g)}:\widetilde{f}_{\upmu}^*\frakA_{X,(U_2,\frakT_2,u_2)} \xrightarrow{\sim} \frakA_{X,(U_1,\frakT_1,u_1)},
\end{equation}
where $\widetilde{f}_{\upmu}$ is given in \eqref{gmu}. The data \eqref{Pandaeq2021} and \eqref{Pandaeq2022} satisfy the conditions of \ref{remShiho1920} and hence are the descent data of a module $\frakA_X$ of $\widetilde{\underline{\calE}}_{\mathrm{strict}}(X/\frakS)_{/\underline{\upmu}_{\mathrm{strict}}}.$ Since \eqref{Pandaeq2022} is an isomorphism, $\frakA_X$ is a crystal. By \ref{equivShihoind}, it extends to a crystal $\frakA_X$ of $\widetilde{\underline{\calE}}(X/\frakS)_{/\underline{\upmu}}.$
\end{parag}

\begin{proposition}\label{propcrystens}
Let $\F$ and $\calG$ be $\Ox_{\calI'}$-modules of $\widetilde{\calE}(X'/\frakS)_{/\calI'}$ (resp. an $\Ox_{\underline{\calI}}$-modules of $\widetilde{\underline{\calE}}(X/\frakS)_{/\underline{\calI}}$) \eqref{kraz1911}.
\begin{enumerate}
\item If $\F$ and $\calG$ are both quasi-coherent modules or both crystals, then so is $\F \otimes_{\Ox_{\calI'}}\calG$ (resp. $\F\otimes_{\Ox_{\underline{\calI}}}\calG$).
\item Let $h_X:\calJ_X \ra \calI_X$ be a morphism of sheaves of $X_{\text{ét}},$ $h':\calJ':=u_X'^{-1}\calJ_X \ra \calI'$ (resp. $\underline{h}:\underline{\calJ}:=\underline{u}_X^{-1}\calJ_X \ra \underline{\calI}$) the morphism induced by $h_X$ and
$$
h':\widetilde{\calE}(X'/\frakS)_{/\calJ'} \ra \widetilde{\calE}(X'/\frakS)_{/\calI'}\ \left (\widetilde{\underline{\calE}}(X/\frakS)_{/\underline{\calJ}} \ra \widetilde{\underline{\calE}}(X/\frakS)_{/\underline{\calI}} \right )
$$
the morphism of ringed topoi induced by $h'$ (resp. $\underline{h}$). If $\F$ is a quasi-coherent module or a crystal of $\widetilde{\calE}(X'/\frakS)_{/\calI'}$ (resp. $\widetilde{\underline{\calE}}(X/\frakS)_{/\underline{\calI}}$), then so is $h'^*\F$ (resp. $\underline{h}^*\F$).
\end{enumerate}
\end{proposition}

\begin{proof}
We only prove the case of $\widetilde{\calE}(X'/\frakS)_{/\calI'}$ as the other is similar. Let $(\F_s)$ and $(\calG_s)$ be descent data corresponding to $\F$ and $\calG$ respectively \eqref{inddescentdata}. Then, for every local section $s\in \Gamma(U,\calI_X)$ over an étale $X'$-scheme $U,$ we have a canonical isomorphism
$$
\left (\F \otimes_{\Ox_{\calI'}} \calG \right )_s=\alpha_s'^*\left (\F \otimes_{\Ox_{\calI'}}\calG \right )\xrightarrow{\sim} \left (\alpha_s'^*\F \right ) \otimes_{\Ox_{\calE(U/\frakS)}} \left (\alpha_s'^*\calG \right )=\F_s \otimes_{\Ox_{\calE(U/\frakS)}}  \calG_s.
$$
By the commutativity of the diagram
$$
\begin{tikzcd}
\widetilde{\calE}(U/\frakS) \ar{r}{\alpha_s'} \ar[swap]{dr}{\alpha_{h_X(s)}'} & \widetilde{\calE}(X'/\frakS)_{/\calJ'} \ar{d}{h'} \\
 & \widetilde{\calE}(X'/\frakS)_{\calI'}, 
\end{tikzcd}
$$
we have
$$
(h'^*\F)_s=\F_{h_X(s)}.
$$
The result follows.
\end{proof}

\begin{proposition}\label{propcrystensstar}
Let $\calA$ be a $\calI'$-indexed (resp. $\underline{\calI}$-indexed) algebra of $\widetilde{\calE}(X'/\frakS)_{/\calI'}$ (resp. of $\widetilde{\underline{\calE}}(X/\frakS)_{/\underline{\calI}}$) and $\F$ and $\calG$ two $\calI'$-indexed (resp. $\underline{\calI}$-indexed) $\calA$-modules \eqref{inddef}. If $\calA,$ $\F$ and $\calG$ are crystals then so is $\F \circledast_{\calA} \calG.$
\end{proposition}

\begin{proof}
We only prove the case of $\widetilde{\calE}(X'/\frakS)_{/\calI'}$ as the other is similar. Let $s\in \Gamma(U,\calI_X)$ be a section over an étale $X$-scheme $U.$ We have to prove that $\left (\F \circledast_{\calA} \calG\right )_s=\alpha_s'^*\left (\F \circledast_{\calA} \calG\right )$ is a crystal of $\widetilde{\calE}(U'/\frakS).$ Let $(f,g):(U_1,\frakT_1,u_1) \ra (U_2,\frakT_2,u_2)$ be a morphism of $\calE(U'/\frakS).$ We have to prove that the morphism
\begin{equation}\label{eqkraz1923}
c_{\alpha_s'^* \left (\F \circledast_{\calA} \calG \right )}:\widetilde{f}^* \left ( \left (\alpha_s'^* \left (\F \circledast_{\calA} \calG\right )\right )_{(U_2,\frakT_2,u_2)}\right ) \ra \left (\alpha_s'^*\left (\F \circledast_{\calA} \calG\right ) \right ) _{(U_1,\frakT_1,u_1)} 
\end{equation}
is an isomorphism, where
$$
\widetilde{f}:\left (U_{1,\text{ét}},u_{1*}\Ox_{T_{1}} \right ) \ra \left ( U_{2,\text{ét}}, u_{2*}\Ox_{T_{2}} \right )
$$
is given in \eqref{Kokoftilda}. Recall the morphism of topoi
$$
\widetilde{f}_{\calI}: \left (U_{1,\text{ét}/\calI_{U_1}},\left ( u_{1*}\Ox_{T_1} \right )_{\calI_{U_1}} \right ) \ra \left (U_{2,\text{ét}/\calI_{U_2}},  \left (u_{2*}\Ox_{T_2} \right )_{\calI_{U_2}}\right )
$$
is given in \eqref{gmu}.
Consider the morphism of ringed topoi
$$
\Lambda_i:U_{i,\text{ét}} \ra U_{i,\text{ét}/\calI_{U_i}},
$$
induced by the continuous and cocontinuous functor
$$
\text{ét}_{/U_i} \ra \text{ét}_{/U_{i,\text{ét}}\times \calI_{U_i}},\ V \mapsto (V,s_{|V}).
$$
By \eqref{diag188}, the morphism \eqref{eqkraz1923} becomes
$$
\widetilde{f}^*\left (\Lambda_2^*\delta_{(U_2,\frakT_2,u_2),\calI_X}^*\left (\F \circledast_{\calA} \calG \right ) \right ) \ra \Lambda_1^* \delta_{(U_1,\frakT_1,u_1),\calI_X}^* \left (\F \circledast_{\calA} \calG\right ).
$$
By \ref{propKoko1915}, this is isomorphic to
$$
\widetilde{f}^*\left (\Lambda_2^*\left (\F_{(U_2,\frakT_2,u_2)} \circledast_{\calA_{(U_2,\frakT_2,u_2)}} \calG_{(U_2,\frakT_2,u_2)} \right ) \right ) \ra \Lambda_1^* \left (\F_{(U_1,\frakT_1,u_1)} \circledast_{\calA_{(U_1,\frakT_1,u_1)}} \calG_{(U_1,\frakT_1,u_1)}\right ).
$$
By the commutativity of the diagram
$$
\begin{tikzcd}
U_{1,\text{ét}} \ar{rr}{\widetilde{f}} \ar[swap]{d}{\Lambda_1}  & & U_{2,\text{ét}} \ar{d}{\Lambda_2} \\
U_{1,\text{ét}/\calI_{U_1}} \ar{rr}{\widetilde{f}_{\calI}} & & U_{2,\text{ét}/\calI_{U_2}},
\end{tikzcd}
$$
it is isomorphic to
$$
\Lambda_1^*\widetilde{f}_{\calI}^* \left (\F_{(U_2,\frakT_2,u_2)} \circledast_{\calA_{(U_2,\frakT_2,u_2)}} \calG_{(U_2,\frakT_2,u_2)} \right )  \ra \Lambda_1^* \left (\F_{(U_1,\frakT_1,u_1)} \circledast_{\calA_{(U_1,\frakT_1,u_1)}} \calG_{(U_1,\frakT_1,u_1)}\right ).
$$
By \ref{Koko829}, it is isomorphic to
$$
\Lambda_1^*\widetilde{f}_{\calI}^*\F_{(U_2,\frakT_2,u_2)} \circledast_{\Lambda_1^*\widetilde{f}_{\calI}^*\calA_{(U_2,\frakT_2,u_2)}} \Lambda_1^*\widetilde{f}_{\calI}^*\calG_{(U_2,\frakT_2,u_2)}   \ra \Lambda_1^* \F_{(U_1,\frakT_1,u_1)} \circledast_{\Lambda_1^*\calA_{(U_1,\frakT_1,u_1)}} \Lambda_1^*\calG_{(U_1,\frakT_1,u_1)}.
$$
Since $\calA,$ $\F$ and $\calG$ are crystals, we have isomorphisms
\begin{alignat*}{2}
\Lambda_1^*\widetilde{f}_{\calI}^*\calA_{(U_2,\frakT_2,u_2)} & \xrightarrow{\sim} \Lambda_1^* \calA_{(U_1,\frakT_1,u_1)},\\
\Lambda_1^*\widetilde{f}_{\calI}^*\F_{(U_2,\frakT_2,u_2)} & \xrightarrow{\sim} \Lambda_1^* \F_{(U_1,\frakT_1,u_1)},\\
\Lambda_1^*\widetilde{f}_{\calI}^*\calG_{(U_2,\frakT_2,u_2)} & \xrightarrow{\sim} \Lambda_1^* \calG_{(U_1,\frakT_1,u_1)}.
\end{alignat*}
The result follows.
\end{proof}

\begin{parag}\label{paragEST1}
In the rest of this section, we take $\calI_X=\upmu_X=\ov{\calM}_X^{gp}.$ Let $V$ be an étale $X$-scheme.
We identify the small étale sites of $V$ and $V'$ via the exact relative Frobenius $V\ra V'$ and hence consider $\upmu_V=\ov{\calM}_V^{gp}$ as a sheaf of $V'_{\text{ét}}.$ We set
$$\upmu'_V=u_V'^{-1}\upmu_V=\upmu_V\circ u_V',\ \underline{\upmu}_V=\underline{u}_V^{-1}\upmu_V=\upmu_V\circ \underline{u}_V.$$
For $V=X,$ we denote $\upmu'_V$ (resp. $\underline{\upmu}_V$) simply by $\upmu'$ (resp. $\underline{\upmu}$).
Let $n$ be a positive integer, $U$ an étale $X$-scheme. Since $u_X'^{-1}$ and $\underline{u}_X^{-1}$ commute with finite products, we have
$$
\upmu'^n=u_X'^{-1}\left (\upmu_X^n\right )=\upmu_X^n\circ u_X',\ \underline{\upmu}^n=\underline{u}_X^{-1}\left (\upmu_X^n\right )=\upmu_X^n \circ \underline{u}_X.
$$
We easily prove that the morphism of topoi \eqref{eqrhoEST1} (resp. \eqref{eqrhoEST2})
$$
\underline{\rho_{X,V}}:\widetilde{\underline{\calE}}(V/\frakS) \ra \widetilde{\underline{\calE}}(X/\frakS)\ \left (\text{resp.}\ \rho_{X',V'}:\widetilde{\calE}(V'/\frakS) \ra \widetilde{\calE}(X'/\frakS)\right )
$$
induces a morphism of topoi
\begin{equation}\label{eqREST}
\underline{R_{X,V}}:\widetilde{\underline{\calE}}(V/\frakS)_{/\underline{\upmu}_V} \ra \widetilde{\underline{\calE}}(X/\frakS)_{/\underline{\upmu}}\ \left (\text{resp.}\ R_{X',V'}:\widetilde{\calE}(V'/\frakS)_{/\upmu'_V} \ra \widetilde{\calE}(X'/\frakS)_{\upmu'}\right ).
\end{equation}
For a local section $s\in \Gamma \left ( W,\upmu_X \right )$ over an étale $V$-scheme $W,$ we clearly have the commutative diagram
$$
\begin{tikzcd}
\widetilde{\underline{\calE}}(W/\frakS) \ar{drr} \ar{d} & & \\
\widetilde{\underline{\calE}}(V/\frakS)_{\underline{\upmu}_V} \ar[swap]{rr}{\underline{R_{X,V}}} & & \widetilde{\underline{\calE}}(X/\frakS)_{/\underline{\upmu}},
\end{tikzcd}
$$
where the vertical and oblique arrows are given in \eqref{alphakraz}. By a similar diagram for $R_{X',V'},$ we deduce that the inverse image functors $\underline{R_{X,V}}^{-1}$ and $R_{X',V'}^{-1}$ preserve crystals.
By \eqref{stackequiv1} and \eqref{stackequiv2}, the morphisms $\underline{R_{X,V}}$ and $R_{X',V'}$ identify with the localization morphisms
\begin{equation}\label{eqRlocEST}
\underline{R_{X,V}}:\widetilde{\underline{\calE}}(X/\frakS)_{/\underline{u}_X^{-1}\left (V^a\right ) \times \underline{\upmu}} \ra \widetilde{\underline{\calE}}(X/\frakS)_{/\underline{\upmu}},\ R_{X',V'}:\widetilde{\calE}(V'/\frakS)_{/u_X'^{-1}\left (V'^a\right ) \times \upmu'} \ra \widetilde{\calE}(X'/\frakS)_{\upmu'}.
\end{equation}
It follows, by \ref{Locprop1}, that the inverse image functors $\underline{R_{X,V}}^{-1}$ and $R_{X',V'}^{-1}$ preserve structures of indexed algebras and indexed modules.
\end{parag}

\begin{parag}\label{paragEST2}
Consider the split fibered $\widetilde{\calE}(X'/\frakS)_{/\upmu'}$-category $\mathbbm{E}'$ (resp. $\widetilde{\underline{\calE}}(X/\frakS)_{/\underline{\upmu}}$-category $\underline{\mathbbm{E}}$) defined, for every sheaf $\F$ of $\widetilde{\calE}(X'/\frakS)_{/\upmu'}$ (resp. $\widetilde{\underline{\calE}}(X/\frakS)_{/\underline{\upmu}}$), by the localized category $\left (\widetilde{\calE}(X'/\frakS)_{/\upmu'} \right )_{/\F}$ (resp. $\left (\widetilde{\underline{\calE}}(X/\frakS)_{/\underline{\upmu}}\right )_{/\F}$) and, for every morphism $\varphi:\F \ra \calG,$ by the inverse image functor of the morphism of topoi
\begin{equation}
\left (\widetilde{\calE}(X'/\frakS)_{/\upmu'} \right )_{/\F} \ra \left (\widetilde{\calE}(X'/\frakS)_{/\upmu'} \right )_{/\calG} \left( \text{resp.}\ \left (\widetilde{\underline{\calE}}(X/\frakS)_{\underline{\upmu}}\right )_{/\F} \ra \left (\widetilde{\underline{\calE}}(X/\frakS)_{/\underline{\upmu}}\right )_{/\calG}\right )
\end{equation}
induced by $\varphi.$ By (\cite{Giraud} II 3.4.4), $\mathbbm{E}'$ (resp. $\underline{\mathbbm{E}}$) is a stack.
\end{parag}

\begin{parag}
Recall the functor $\rho:\underline{\calE}(X/\frakS) \ra \calE(X'/\frakS)$ \eqref{era3rhoShiho} and the functors $u_X'$ and $\underline{u}_X$ defined in \ref{paragu}. By definition, these functors fit into the commutative diagram
$$
\begin{tikzcd}
\underline{\calE}(X/\frakS) \ar[swap]{d}{\underline{u}_X} \ar{r}{\rho} & \calE(X'/\frakS) \ar{d}{u_X'} \\
\text{ét}_{/X} \ar{r}{\sim} & \text{ét}_{/X'},
\end{tikzcd}
$$
where the lower equivalence of sites is induced by the exact relative Frobenius $F_1:X\ra X'.$
Since these functors are continuous and cocontinuous, we deduce that the diagram
\begin{equation}\label{diag14151}
\begin{tikzcd}
\widetilde{\underline{\calE}}(X/\frakS) \ar{rr}{C_{X/\frakS}} \ar[swap]{d}{\underline{u}_X} & & \widetilde{\calE}(X'/\frakS) \ar{d}{u_X'} \\
X_{\text{ét}} \ar{rr}{\sim} & & X'_{\text{ét}},
\end{tikzcd}
\end{equation}
where $C_{X/\frakS}$ is given in \eqref{era3C} and $u_X'$ and $\underline{u}_X$ are defined in \eqref{Kokou}, is commutative up to a canonical isomorphism.
Recall that $\upmu'=u_X'^{-1}F_{1*}\upmu_X.$ Then
$$C_{X/\frakS}^{-1}\upmu'=\underline{u}_X^{-1}F_1^{-1}F_{1*}\upmu_X=\underline{u}_X^{-1}\upmu_X=\underline{\upmu}.$$
Therefore \eqref{diag14151} induces a commutative diagram
\begin{equation}\label{diagkraz2}
\begin{tikzcd}
\widetilde{\underline{\calE}}(X/\frakS)_{/\underline{\upmu}} \ar{rr}{C_{X/\frakS,\upmu}} \ar{d}{\underline{u}_{X,\underline{\upmu}}} & & \widetilde{\calE}(X'/\frakS)_{/\upmu'} \ar{d}{u'_{X,\upmu'}} \\
X_{\text{ét}/\upmu_X} \ar{rr}{\sim} & & X'_{\text{ét}/F_{1*}\upmu_X}.
\end{tikzcd}
\end{equation}
The continuous and cocontinuous functor $\rho$ also induces a continuous and cocontinuous functor
$$
\rho_{\upmu}:\underline{\calE}(X/\frakS)_{/\underline{\upmu}} \ra \calE(X'/\frakS)_{/\upmu'}.
$$
The inverse image functor $C_{X/\frakS,\upmu}^{-1}$ is then given by composition with $\rho_{\upmu}.$
For every local section $s$ of $\upmu_X$ over an étale $X$-scheme $U,$ we have a bigger commutative diagram:
\begin{equation}\label{diagkraz19253}
\begin{tikzcd}
\widetilde{\underline{\calE}}(U/\frakS) \ar{dd}{\underline{u}_X} \ar{rr}{C_{U/\frakS}} \ar{dr}{\underline{\alpha}_s} & & \widetilde{\calE}(U'/\frakS) \ar[dashed]{dd}[pos=0.20]{u_X'} \ar{dr}{\alpha_s'} &  \\
 & \widetilde{\underline{\calE}}(X/\frakS)_{/\underline{\upmu}} \ar{rr}[pos=0.20]{C_{X/\frakS,\upmu}} \ar{dd}[pos=0.20]{\underline{u}_{X,\underline{\upmu}}} & & \widetilde{\calE}(X'/\frakS)_{/\upmu'} \ar{dd}{u'_{X,\upmu'}} \\
U_{\text{ét}} \ar{dr} \ar[dashed]{rr}[pos=0.20]{\sim} & & U'_{\text{ét}} \ar[dashed]{dr} &  \\
 & X_{\text{ét}/\upmu_X} \ar{rr}{\sim} & & X'_{\text{ét}/\upmu_X},
\end{tikzcd}
\end{equation}
where $U_{\text{ét}} \ra X_{\text{ét}/\upmu_X}$ and $U'_{\text{ét}} \ra X'_{\text{ét}/\upmu_X}$ are the localization morphisms corresponding to
$$
s\in \Gamma(U,\upmu_X)=\Gamma(U',F_{1*}\upmu_X).
$$
\end{parag}

\begin{proposition}\label{commCXU}
Let $U$ be an étale $X$-scheme. Then, the diagram
$$
\begin{tikzcd}
\widetilde{\underline{\calE}}(U/\frakS)_{/\underline{\upmu}_U} \ar{d} \ar{rr}{C_{U/\frakS,\upmu}} & & \widetilde{\calE}(U'/\frakS)_{/\upmu'_U} \ar{d} \\
\widetilde{\underline{\calE}}(X/\frakS)_{/\underline{\upmu}} \ar[swap]{rr}{C_{X/\frakS,\upmu}} & & \widetilde{\calE}(X'/\frakS)_{/\upmu'},
\end{tikzcd}
$$
where $\upmu'_U$ and $\underline{\upmu}_U$ are defined in \ref{paragEST1}, the horizontal arrows are defined in \eqref{diagkraz2} and the vertical arrows are \eqref{eqREST}, is commutative.
\end{proposition}

\begin{proof}
This follows from the commutativity of
$$
\begin{tikzcd}
\widetilde{\underline{\calE}}(U/\frakS) \ar[swap]{d}{\underline{\rho_{X,U}}} \ar{rr}{C_{U/\frakS}} & & \widetilde{\calE}(U'/\frakS) \ar{d}{\rho_{X',U'}} \\
\widetilde{\underline{\calE}}(X/\frakS) \ar[swap]{rr}{C_{X/\frakS}} & & \widetilde{\calE}(X'/\frakS).
\end{tikzcd}
$$
\end{proof}

\begin{proposition}\label{ommek4}
Let $\frakE$ be a module of $\widetilde{\calE}(X'/\frakS)_{/\upmu'}.$ For any object $(U,\frakT,u)$ of $\underline{\calE}(X/\frakS),$ there exists a canonical isomorphism
$$
\left (C_{X/\frakS,\upmu}^{-1} \frakE \right )_{(U,\frakT,u)}\xrightarrow{\sim} \frakE_{\rho(U,\frakT,u)}.
$$
\end{proposition}

\begin{proof}
This follows from the commutativity of
$$
\begin{tikzcd}
\text{ét}_{/U\times \upmu_X} \ar{rr}{\delta_{(U,\frakT,u)}} \ar[swap]{d}{(V,s) \mapsto (V',s)} & & \underline{\calE}(X/\frakS)_{/\underline{\upmu}} \ar{d}{\rho_{\upmu}} \\
\text{ét}_{/U'\times \upmu_X} \ar[swap]{rr}{\delta_{\rho(U,\frakT,u)}} & & \calE(X'/\frakS)_{/\upmu'}.
\end{tikzcd}
$$
\end{proof}

\newpage

\section{Indexed logarithmic Cartier transform}

\begin{parag}
In this section, just as in \ref{intro18}, we fix a perfect field $k$ of positive characteristic $p$ and denote by $W(k)$ its ring of Witt vectors. We equip $\op{Spf}W(k)$ with the trivial logarithmic structure. Unless explicitly stated, all formal schemes considered are $p$-adic and all logarithmic structures are fs. If a gothic letter $\frakT$ denotes a $p$-adic logarithmic formal scheme, the corresponding roman letter $T$ will denote the logarithmic scheme obtained from $\frakT$ by reduction modulo $p$ \eqref{propkey}. We also denote by $\frakT_n,$ for any positive integer $n,$ the logarithmic scheme obtained from $\frakT$ by reduction modulo $p^n$ \eqref{propkey}.

Let $\frakS$ be a logarithmic formal scheme, log flat and locally of finite type over $\op{Spf}W(k),$ and $X\ra S$ a log smooth morphism of logarithmic schemes, which is of Cartier type.

For an étale $X$-scheme $U,$ we equip $U$ with the logarithmic structure $\left (\calM_{X}\right )_{|U},$ restriction of $\calM_X$ on $U,$ so that $U\ra X$ is strict.
For every logarithmic scheme $T,$ we denote by $\upmu_T$ the sheaf $\ov{\calM}_T^{gp}.$ We consider the monoids $\upmu_X,$ $\upmu'$ and $\underline{\upmu}$ defined in \ref{paragEST1}.

We will sometimes need to assume one of following two hypothesis: 
\begin{itemize}
\item[$(H_1)$] The morphism $X\ra S$ lifts to a log smooth morphism $f:(\frakX,Q) \ra (\frakS,P)$ of framed logarithmic formal schemes and $X'$ lifts to a framed logarithmic formal scheme $(\frakX',Q')$ over $(\frakS,P),$ such that $\frakX' \ra \frakS$ is log smooth.
\item[$(H_2)$] The hypothesis $(H_1)$ is satisfied and the exact relative Frobenius $F_1:X\ra X'$ lifts to an $(\frakS,P)$-morphism of framed logarithmic formal schemes $F:(\frakX,Q) \ra (\frakX',Q').$
\end{itemize}
When we suppose that $(H_1)$ is satisfied, we consider the logarithmic formal schemes $R_{\frakX',1}$ and $Q_{\frakX}$ defined in \ref{parag86} and we let $Q_1$ and $R_1'$ be the special fibers of $Q_{\frakX}$ and $R_{\frakX',1}$ respectively and $\calQ_1$ and $\calR_1'$ the corresponding Hopf algebras and $\lambda_Q:\underline{Q_1} \ra X$ and $\lambda_{R'}:R_1' \ra X'$ the morphisms given in \ref{erafprop13}. When we suppose that $(H_2)$ is satisfied, we consider the morphism of Hopf algebras $\calV:\calR_1' \ra \calQ_1$ induced by $\nu:Q_{\frakX}\ra R_{\frakX',1}$ \eqref{lem92}.
We start by studying the crystals in the indexed logarithmic Oyama topoi, introduced in the previous section, and then construct an indexed version of the Cartier transform.
\end{parag}

\begin{proposition}\label{equivRQind}
Assume $(H_1)$ is satisfied. There exists a canonical equivalence of categories between the category of crystals of $\widetilde{\underline{\calE}}(X/\frakS)_{/\underline{\upmu}}$ (resp. $\widetilde{\calE}(X'/\frakS)_{/\upmu'}$) and the category of $\Ox_{X,\upmu_X}$-modules equipped with a $\calQ_{1,\upmu_X}$-stratification (resp. $\Ox_{X',\upmu_X}$-modules equipped with an $\calR_{1,\upmu_X}'$-stratification). In addition, these equivalences preserve structures of $\upmu_X$-indexed algebras and $\upmu_X$-indexed modules.
\end{proposition}

\begin{proof}
A crystal $\mathfrak{E}$ of $\widetilde{\underline{\calE}}(X/\frakS)_{/\underline{\upmu}}$ is equivalent to its descent data $(\mathfrak{E}_s)$ (\ref{inddescentdata}). By definition, for any étale $X$-scheme $U$ and $s\in \Gamma(U,\upmu_X),$ the $\Ox_{\underline{\calE}(U/\frakS)}$-module $\mathfrak{E}_s$ is a crystal of $\widetilde{\underline{\calE}}(U/\frakS).$ Let $\frakU$ be the unique étale formal $\frakX$-scheme such that $U=\frakU\times_{\frakX}X.$ Consider the logarithmic formal scheme $Q_{\frakU}$ obtained from $\frakU \ra \frakS$ just as $Q_{\frakX}$ is obtained from $\frakX \ra \frakS.$ By \ref{equivRQ}, $\mathfrak{E}_s$ corresponds to an $\Ox_U$-module $\F_s$ equipped with a $Q_{\frakU}$-stratification $\epsilon_s.$ Since $\calQ_{\frakU}=\calQ_{\frakX|\frakU}$ \eqref{QU}, by \ref{indremark73}, the objects $(\F_s,\epsilon_s)$ correspond to an $\Ox_{X,\upmu_X}$-module $\F$ equipped with a $\calQ_1$-stratification $\epsilon.$ This proves the equivalence.

Suppose now that $\mathfrak{E}$ has a structure of a $\upmu_X$-indexed algebra given by a morphism
$$\pi:\mathfrak{E}\boxtimes \mathfrak{E} \ra \sigma^*\mathfrak{E},$$
where $\sigma:\upmu^2 \ra \upmu$ is the addition map. By \ref{propcrystens}, $\mathfrak{E}\boxtimes \mathfrak{E}$ and $\sigma^*\mathfrak{E}$ are crystals. For every local sections $s,t$ of $\upmu_X$ over an étale $X$-scheme $U,$ by applying $\underline{\alpha}_{(s,t)}^{-1}$ and \ref{prop1916Ko}, we obtain a morphism of crystals
$$\pi_{(s,t)}:\mathfrak{E}_s\otimes_{\Ox_{\underline{\calE}(U/\frakS)}}\mathfrak{E}_t \ra \mathfrak{E}_{s+t}.$$
This corresponds, by \ref{equivRQ}, to a morphism
$$\F_s \otimes_{\Ox_U}\F_t \ra \F_{s+t}$$
of $\Ox_U$-modules with $Q_{\frakU}$-stratification. These morphisms correspond, by \ref{indremark73}, to a morphism
$$\F\boxtimes \F\ra \sigma_X^*\F,$$
where $\sigma_X:\upmu_X^2 \ra \upmu_X$ is the addition map. This defines a $\upmu_X$-indexed algebra structure on $\F.$
The proof for $\widetilde{\calE}(X'/\frakS)_{/\upmu'}$ and for indexed modules is similar.
\end{proof}

\begin{parag}\label{parag1419}
Assume $(H_1)$ is satisfied. As in \ref{equivRQDef}, we have canonical projections $q_1,q_2:(X,Q_{\frakX},\lambda_Q) \ra (X,\frakX,i_X)$ of $\underline{\calE}(X/\frakS),$ where $i_X:\underline{X} \ra X$ is the canonical immersion, and $q_1',q_2':(X',R_{\frakX',1},\lambda_{R'}) \ra (X',\frakX',\op{Id}_{X'})$ of $\calE'(X/\frakS).$
Let $\F$ be a module of $\widetilde{\calE}(X'/\frakS)_{/\upmu'}.$ The morphism $c_{\F,q_i'}$ \eqref{cind} is equal to
$$
c_{\F,q_i'}:\left (\lambda_{R'*}\Ox_{R_1'} \right )_{\upmu_X} \otimes_{\Ox_{X',\upmu_X}} \F_{(X',\frakX',\op{Id}_{X'})} \ra \F_{(X',R_{\frakX',1},\lambda_{R'})},
$$
where $\Ox_{X'} \ra \lambda_{R'*}\Ox_{R_1'}$ is induced by $q_i'.$ Similarly, if $\F$ is a module of $\widetilde{\underline{\calE}}(X/\frakS)_{/\underline{\upmu}}$ the morphism $c_{\F,q_i}$ is equal to
$$
c_{\F,q_i}:\left (\lambda_{Q*}\Ox_{Q_1} \right )_{\upmu_X} \otimes_{\Ox_{X,\upmu_X}} \F_{(X,\frakX,i_{X})} \ra \F_{(X,Q_{\frakX},\lambda_{Q})}.
$$
Recall that $\calR_1'=\lambda_{R'*}\Ox_{R_1'}$ and $\calQ_1=\lambda_{Q*}\Ox_{Q_1}.$
By commutativity of \eqref{diag188}, the equivalence given in \ref{equivRQind} can be explicitly given as follows:
To a crystal $\mathfrak{E}$ of $\widetilde{\underline{\calE}}(X/\frakS)_{/\underline{\upmu}}$ (resp. $\widetilde{\calE}(X'/\frakS)_{/\upmu'}$), corresponds the $\Ox_{X,\upmu_X}$-module $\mathfrak{E}_{(X,\frakX,i_X)}$ (resp. the $\Ox_{X',\upmu_X}$-module $\mathfrak{E}_{(X',\frakX',\op{Id}_{X'})}$) equipped with the $\calQ_{1,\upmu_X}$-stratification
$$\epsilon:\calQ_{1,\upmu_X} \otimes_{\Ox_{X,\upmu_X}} \mathfrak{E}_{(X,\frakX,i_X)} \xrightarrow{c_{\mathfrak{E},q_2}} \mathfrak{E}_{(X,Q_{\frakX},\lambda_Q)} \xrightarrow{c_{\mathfrak{E},q_1}^{-1}} \mathfrak{E}_{(X,\frakX,i_X)} \otimes_{\Ox_{X,\upmu_X}}\calQ_{1,\upmu_X}.$$
(resp. the $\calR_{1,\upmu_X}'$-stratification
$$\epsilon':\calR'_{1,\upmu_X} \otimes_{\Ox_{X',\upmu_X}} \mathfrak{E}_{(X',\frakX',\op{Id}_{X'})} \xrightarrow{c_{\mathfrak{E},q_2'}} \mathfrak{E}_{(X',R_{\frakX',1},\lambda_{R'})} \xrightarrow{c_{\mathfrak{E},q_1'}^{-1}} \mathfrak{E}_{(X',\frakX',\op{Id}_{X'})})\otimes_{\Ox_{X',\upmu_X}}\calR'_{1,\upmu_X})$$
\end{parag}

\begin{proposition}\label{propEST1831}
Assume $(H_1)$ is satisfied. Let $U$ be an étale $X$-scheme. The diagrams
$$
\begin{tikzcd}
\begin{Bmatrix}
\mathrm{Crystals\ of\ }\widetilde{\underline{\calE}}(X/\frakS)_{/\underline{\upmu}}
\end{Bmatrix}
\ar{r}{\sim} \ar{d} &
\begin{Bmatrix}
\Ox_{X,\upmu_X}\text{-}\mathrm{modules\ equipped} \\
\mathrm{with\ a\ }\calQ_{1,\upmu_X}\text{-}\mathrm{stratification}
\end{Bmatrix}
 \ar{d} \\
\begin{Bmatrix}
\mathrm{Crystals\ of\ }\widetilde{\underline{\calE}}(U/\frakS)_{/\underline{\upmu}_U}
\end{Bmatrix} 
\ar{r}{\sim} &
\begin{Bmatrix}
\Ox_{U,\upmu_U}\text{-}\mathrm{modules\ equipped} \\
\mathrm{with\ a\ }\calQ_{1|U,\upmu_U}\text{-}\mathrm{stratification}
\end{Bmatrix}
\end{tikzcd}
$$
and
$$
\begin{tikzcd}
\begin{Bmatrix}
\mathrm{Crystals\ of\ }\widetilde{\calE}(X'/\frakS)_{/\upmu'}
\end{Bmatrix}
\ar{r}{\sim} \ar{d} &
\begin{Bmatrix}
\Ox_{X',\upmu_X}\text{-}\mathrm{modules\ equipped} \\
\mathrm{with\ a\ }\calR'_{1,\upmu_X}\text{-}\mathrm{stratification}
\end{Bmatrix}
 \ar{d} \\
\begin{Bmatrix}
\mathrm{Crystals\ of\ }\widetilde{\calE}(U'/\frakS)_{/\upmu'_U}
\end{Bmatrix}
\ar{r}{\sim} &
\begin{Bmatrix}
\Ox_{U',\upmu_U}\text{-}\mathrm{modules\ equipped} \\
\mathrm{with\ a\ }\calR'_{1|U',\upmu_U}\text{-}\mathrm{stratification}
\end{Bmatrix},
\end{tikzcd}
$$
where the horizontal arrows are the equivalences given in \ref{equivRQind}, the right vertical arrows are the restriction functors and the left arrows are induced by the inverse image functors of \eqref{eqREST},
are commutative.
\end{proposition}

\begin{proof}
We just check the commutativity of the first diagram. The second is similar. For any local section $s\in \Gamma\left (V,\upmu_X \right )$ over an étale $X$-scheme $V,$ we have the diagram
$$
\begin{tikzcd}
\begin{Bmatrix}
\mathrm{Crystals\ of\ }\widetilde{\underline{\calE}}(X/\frakS)_{/\underline{\upmu}}
\end{Bmatrix}
\ar{r}{\sim} \ar{d} &
\begin{Bmatrix}
\Ox_{X,\upmu_X}\text{-}\mathrm{modules\ equipped} \\
\mathrm{with\ a\ }\calQ_{1,\upmu_X}\text{-}\mathrm{stratification}
\end{Bmatrix}
 \ar{d} \\
\begin{Bmatrix}
\mathrm{Crystals\ of\ }\widetilde{\underline{\calE}}(U/\frakS)_{/\underline{\upmu}_U}
\end{Bmatrix} 
\ar{r}{\sim} \ar[swap]{d}{\alpha_s^{-1}} &
\begin{Bmatrix}
\Ox_{U,\upmu_U}\text{-}\mathrm{modules\ equipped} \\
\mathrm{with\ a\ }\calQ_{1|U,\upmu_U}\text{-}\mathrm{stratification}
\end{Bmatrix} \ar{d}{(\F,\epsilon ) \mapsto (\F_s,\epsilon_s)} \\
\begin{Bmatrix}
\mathrm{Crystals\ of\ }\widetilde{\underline{\calE}}(V/\frakS)
\end{Bmatrix} 
\ar{r}{\sim} &
\begin{Bmatrix}
\Ox_{V}\text{-}\mathrm{modules\ equipped} \\
\mathrm{with\ a\ }\calQ_{1|V}\text{-}\mathrm{stratification}
\end{Bmatrix}.
\end{tikzcd}
$$
The lower and outer squares are commutative by definition of the equivalences \ref{equivRQind}. This implies the commutativity of the upper square.
\end{proof}

\begin{parag}\label{defQstrat}
Here, we give an alternate description of the crystal $\frakA_X$ defined in \ref{AcrysShiho} under the hypothesis $(H_1).$ Assume $(H_1)$ is satisfied. Let $q_1,q_2:Q_{\frakX} \ra \frakX$ and $q_1',q_2':R_{\frakX',1} \ra \frakX'$ be the canonical projections. By definition of $Q_{\frakX},$ the projections $q_1$ and $q_2$ are strict and so they induce, by (\cite{Lor2000} I 3.1 page 273), isomorphisms of $\upmu_X$-indexed algebras
\begin{equation}\label{Khaminei2151}
c_i:q_i^*\calA_X \xrightarrow{\sim} \calA_{Q_1}.
\end{equation}
We get a $\calQ_{1,\upmu_X}$-stratification
\begin{equation}\label{eqeAQ}
\epsilon_{\calA_X,Q}:q_2^*\calA_X \xrightarrow{c_2} \calA_{Q_1} \xrightarrow{c_1^{-1}} q_1^*\calA_X.
\end{equation}
Note that $\epsilon_{\calA_X,Q}$ is a morphism of $\upmu_X$-indexed algebras.
We denote by $\frakA$ the crystal of $\widetilde{\underline{\calE}}(X/\frakS)_{/\underline{\upmu}}$ corresponding to $\left (\calA_X,\epsilon_{\calA_X,Q} \right )$ \eqref{equivRQind}. By \ref{equivRQind}, $\frakA$ is a $\upmu$-indexed algebra. By construction,
\begin{equation}\label{eq14201}
\frakA_{(X,\frakX,i_X)}=\calA_X,\ \frakA_{(X,Q_{\frakX},\lambda_Q)}=\calA_{Q_1}.
\end{equation}
The crystal $\frakA$ is hence equal to the crystal $\frakA_X$ defined in \ref{AcrysShiho}.
\end{parag}

\begin{remark}\label{remkraz1}
By definition, the stratification $\epsilon_{\calA_X,Q}$ is an isomorphism of $\upmu_X$-indexed algebras.
\end{remark}

\begin{parag}\label{paragBglue}
There exist étale coverings $\left ( \frakS_i \ra \frakS \right )_{i\in I}$ and $\left (X_{ij} \ra X_i:=X\times_SS_i \right )_{j\in I_i}$ satisfying:
\begin{itemize}
\item $\frakS_i$ has a chart $\frakS_i \ra B\langle P_i \rangle.$
\item $X_{ij} \ra S_i$ has a chart $A_1[Q_{ij}] \ra A_1[P_i]$ such that the morphism $X_{ij} \ra S_i\times_{A_1[P_i]}A_1[Q_{ij}],$ induced by $X_{ij} \ra S_i$ and the chart $X_{ij} \ra A_1[Q_{ij}],$ is étale and strict and $P_i^{gp} \ra Q_{ij}^{gp}$ is injective and the torsion subgroup of its cokernel has a finite order coprime with $p.$
\end{itemize}
Let $Q_{ij}'$ and $Q_{ij}''$ be the monoids defined from $P_i \ra Q_{ij}$ as in \ref{PFrob} and $X_{ij}''$ the fiber product, in the category of fine logarithmic schemes, of $X_{ij}\ra S_i$ and the Frobenius $F_{S_i}:S_i \ra S_i.$ By \ref{PFrob}, we have a chart $X_{ij}'' \ra A_1[Q_{ij}''],$ $X_{ij}'=X_{ij}''\times_{A_1[Q_{ij}'']}A_1[Q_{ij}']$ and the canonical morphism $P_i \ra Q_{ij}'$ is a chart of $X_{ij}' \ra S_i.$
For $i\in I$ and $j\in I_i,$ we have seen, in \ref{paragu}, that the functor
\begin{equation}\label{eqShiho1}
\begin{array}[t]{clc}
\calE(X_{ij}'/\frakS) & \ra & \text{ét}_{/X_{ij}'} \\
(U,\frakT,u) & \mapsto & U
\end{array}
\end{equation}
is continuous and cocontinuous and hence induces a morphism of topoi
$$
u_{ij}:\widetilde{\calE}(X_{ij}'/\frakS) \ra X'_{{ij},\text{ét}},
$$
such that $u_{ij}^{-1}$ is given by composition with \eqref{eqShiho1}. We set
$$
\upmu_{ij}'=u_{ij}^{-1} \left (\upmu_{X_{ij}} \right ) = u_{ij}^{-1} \left (\ov{\calM}_{X_{ij}}^{gp} \right ).
$$
Similarly, for $i,j\in I,$ $k\in I_i$ and $l\in I_j,$ we have a morphism of topoi
$$
u_{ikjl}:\widetilde{\calE}\left ( X_{ik}'\times_{X'}X_{jl}'/\frakS \right ) \ra \left (X_{ik}'\times_{X'}X_{jl}' \right )_{\text{ét}}
$$
and we set
$$
\upmu_{ikjl}'=u_{ikjl}^{-1} \left (\upmu_{X_{ik}\times_XX_{jl}} \right ).
$$
By the lemma \ref{lemShiho208} below, the morphism $X_{ij}' \ra S_i\times_{A_1[P_i]}A_1[Q_{ij}']$ is étale and strict and so there exists a unique étale and strict morphism of logarithmic formal schemes $\frakX'_{ij} \ra \frakS_i\times_{B\langle P_i \rangle}B\langle Q_{ij}' \rangle,$ fitting into a cartesian square
$$
\begin{tikzcd}
X_{ij}' \ar{r} \ar{d} & \frakX_{ij}' \ar{d} \\
S_i \times_{A_1[P_i]}A_1[Q_{ij}'] \ar{r} & \frakS_i\times_{B\langle P_i \rangle }B\langle Q_{ij}' \rangle .
\end{tikzcd}
$$
We have the dilatation $R_{\frakX_{ij}',1},$ defined from $\left ( \frakX_{ij}',Q_{ij}' \right ) \ra \left ( \frakS_i, P_i \right ),$ as in \ref{parag86}. Let $R_{ij}'$ be the logarithmic scheme obtained by reduction of $R_{\frakX_{ij}',1}$ modulo $p$ and $\calR_{ij}'$ its structural ring. Now consider the $\Ox_{X_{ij}',\upmu_{X_{ij}}}$-module $\calB_{X_{ij}/S}$ \eqref{Bgroupe2} equipped with the stratification
$$
\calR'_{ij,\upmu_{X_{ij}}} \otimes_{\Ox_{X_{ij}',\upmu_{X_{ij}}}} \calB_{X_{ij}/S} \xrightarrow{\sim} \calB_{X_{ij}/S} \otimes_{\Ox_{X_{ij}',\upmu_{X_{ij}}}} \calR'_{{ij},\upmu_{X_{ij}}}
$$
that exchanges factors (note that the two projections $R_{ij}' \ra X'_{ij}$ are equal by \ref{erafprop13}). By \ref{equivRQind}, it corresponds to a crystal $\frakB_{ij}$ of $\widetilde{\calE}\left ( X_{ij}' / \frakS \right )_{/\upmu_{ij}'}.$ By \ref{equivRQind}, $\frakB_{ij}$ is a $\upmu_{ij}$-indexed algebra.
We prove in \ref{propBglue} that the crystals $\frakB_{ij}$ glue together into a crystal $\frakB$ of $\widetilde{\calE}\left (X'/\frakS \right )_{/\upmu'}.$
\end{parag}

\begin{lemma}\label{lemShiho208}
Keep the notation of \ref{paragBglue}. For $i\in I$ and $j\in I_i,$ the morphism $X_{ij}' \ra S_i\times_{A_1[P_i]}A_1[Q_{ij}'],$ induced by $X_{ij}' \ra S_i$ and the chart $X_{ij}' \ra A_1[Q_{ij}'],$ is étale and strict.
\end{lemma}

\begin{proof}
Consider the cartesian squares in the category of fine logarithmic schemes:
$$
\begin{tikzcd}
T_{ij} \ar{r} \ar{d} & S_i \times_{A_1[P_i]}A_1[Q_{ij}] \ar{r} \ar{d} & A_1[Q_{ij}] \ar{d} \\
S_i \ar[swap]{r}{F_{S_i}} & S_i \ar{r} & A_1[P_i].
\end{tikzcd}
$$
The composition
$$
S_i \xrightarrow{F_{S_i}} S_i \rightarrow A_1[P_i]
$$
is equal to
$$
S_i \ra A_1[P_i] \xrightarrow{F_{A_1[P_i]}} A_1[P_i].
$$
It follows that we also have the following cartesian squares in the category of fine logarithmic schemes:
$$
\begin{tikzcd}
T_{ij} \ar{r} \ar{d} & A_1[Q_{ij}''] \ar{r} \ar{d} & A_1[Q_{ij}] \ar{d} \\
S_i \ar{r} & A[P_i] \ar[swap]{r}{F_{A_1[P_i]}} & A_1[P_i].
\end{tikzcd}
$$
The morphism $X_{ij}\ra S_i \times_{A_1[P_i]}A_1[Q_{ij}]$ is étale and strict hence so is its base change by $F_{S_i}$ in the category of fine logarithmic schemes, which is $X_{ij}''\ra T_{ij}=S_i \times_{A_1[P_i]}A_1[Q_{ij}''].$ The result then follows by base change by the canonical $A_1[Q_{ij}'] \ra A_1[Q_{ij}''].$
\end{proof}

\begin{proposition}\label{propBglue}
Keep the notation of \ref{paragBglue}. For $i\in I,$ let
$$
Res_{ij}:\widetilde{\calE}\left (X_{ij}' / \frakS \right )_{/\upmu'_{ij}}  \ra  \widetilde{\calE} \left (X' / \frakS \right )_{/\upmu'}
$$
be the morphism of topoi induced by $\rho_{X',X_{ij}'}$ \eqref{eqrhoEST2}. There exists a unique crystal $\frakB_X$ of $\widetilde{\calE}(X'/\frakS)_{/\upmu'}$ such that, for all $i\in I$ and $j\in I_i,$ there exists a canonical isomorphism
$$
Res_{ij}^{-1}\frakB_X \xrightarrow{\sim} \frakB_{ij}.
$$
In addition, $\frakB_X$ is naturally equipped with a structure of a $\upmu'$-indexed algebra.
\end{proposition}

\begin{proof}
For $i,j\in I$ and $k\in I_i,$ $l\in I_j,$ the morphism of topoi $\rho_{X_{ik}',X_{ik}'\times_{X'}X_{jl}'}$ \eqref{eqrhoEST2} induces a morphism of topoi
$$
Res_{ikjl}: \widetilde{\calE}\left (X_{ik}' \times_{X'} X_{jl}' / \frakS \right )_{/\upmu'_{ikjl}}  \ra  \widetilde{\calE} \left (X_{ik}' / \frakS \right )_{/\upmu'_{ik}}
$$
such that $Res_{ikjl}^{-1}$ is given by composition with
\begin{equation}\label{eqShiho2}
\begin{array}[t]{clc}
\calE\left (X_{ik}' \times_{X'} X_{jl}' / \frakS \right )_{/\upmu_{ikjl}'} & \ra & \calE \left (X_{ik}' / \frakS \right )_{/\upmu_{ik}'} \\
(U,\frakT,u,s) & \mapsto & \left (\rho_{X_{ik}',X_{ik}'\times_{X'}X_{jl}'}\left (U,\frakT,u \right ),s \right ),
\end{array}
\end{equation}
where $(U,\frakT,u)$ is an object of $\calE\left (X_{ik}' \times_{X'} X_{jl}' / \frakS \right )$ and $s\in \Gamma(U,\upmu_X).$
We construct an isomorphism
$$
Res_{ikjl}^{-1} \frakB_{ik} \xrightarrow{\sim} Res_{jlik}^{-1}\frakB_{jl}.
$$
Let $(U,\frakT,u)$ be an object of $\calE \left (X_i' \times_{X'} X_j' / \frakS \right ).$ We abusively denote $\rho_{X_{ik}',X_{ik}'\times_{X'}X_{jl}'}(U,\frakT,u)$ by $(U,\frakT,u).$ By lemma \ref{lemShiho209} below,
$$
\left (Res_{ikjl}^{-1}\frakB_{ik} \right )_{(U,\frakT,u)}=\frakB_{ik,(U,\frakT,u)}. 
$$
Suppose that $U$ is affine, $\frakT$ is then affine and there exists a morphism $\frakT \ra \frakX_{ik}'$ fitting into the commutative diagram
$$
\begin{tikzcd}
U \ar{r} & X_{ik}'\times_{X'}X_{jl}' \ar{r} & \frakX_{ik}' \ar{d} \\
T \ar{u} \ar{r} & \frakT \ar{ur} \ar{r} & \frakS.
\end{tikzcd}
$$
We obtain a morphism
$$
(U,\frakT,u) \ra (X_{ik}',\frakX_{ik}',\op{Id}_{X_{ik}'})
$$
of $\calE\left ( X_{ik}'/\frakS \right ).$ The corresponding morphism of ringed topoi
$$
\varphi_{ik}: \left ( U_{\text{ét}/\upmu_U},\left (u_*\Ox_T \right )_{\upmu_U} \right ) \ra \left ( X_{{ik},\text{ét}/\upmu_{X_{ik}}}',\Ox_{X_{ik}',\upmu_{X_{ik}}}\right )
$$
does not depend on the choice of the lifting $\frakT \ra \frakX_{ik}'.$ Since $\frakB_{ik}$ is a crystal, we have an isomorphism
\begin{equation}\label{eqShiho3}
\varphi_{ik}^*\calB_{X_{ik}/S}=\varphi_{ik}^*\frakB_{{ik},(X_{ik}',\frakX_{ik}',\op{Id}_{X_{ik}'})} \xrightarrow{\sim} \frakB_{{ik},(U,\frakT,u)}.
\end{equation}
Similarly, we have an isomorphism
\begin{equation}\label{eqShiho4}
\varphi_{jl}^*\calB_{X_{jl}/S} \xrightarrow{\sim} \frakB_{{jl},(U,\frakT,u)}.
\end{equation}
Since $U \ra X_{ik}'$ and $U\ra X_{jl}'$ factor through $U \ra X_{ik}'\times_{X'}X_{jl}',$ the morphisms of ringed topoi $\varphi_{ik}$ and $\varphi_{jl}$ factor through
$$
\varphi_{ikjl}:\left ( U_{\text{ét}/\upmu_U},\left (u_*\Ox_T \right )_{\upmu_U} \right ) \ra \left ( \left ( X_{ik}'\times_{X'}X_{jl}' \right )_{\text{ét}/\upmu_{X_{ik}\times_XX_{jl}}},\Ox_{X_{ik}'\times_{X'}X_{jl}',\upmu_{X_{ik}\times_XX_{jl}}}\right ).
$$
By \ref{lemBres}, the isomorphisms \eqref{eqShiho3} and \eqref{eqShiho4} induce then isomorphisms
$$
\varphi_{ikjl}^*\calB_{X_{ik}\times_XX_{jl}/S} \xrightarrow{\sim} \frakB_{ik,(U,\frakT,u)},
$$
$$
\varphi_{ikjl}^*\calB_{X_{ik}\times_XX_{jl}/S} \xrightarrow{\sim} \frakB_{jl,(U,\frakT,u)}.
$$
We obtain an isomorphism
$$
\frakB_{ik,(U,\frakT,u)} \xrightarrow{\sim} \varphi_{ikjl}^*\calB_{X_{ik}\times_XX_{jl}/S} \xrightarrow{\sim} \frakB_{jl,(U,\frakT,u)}.
$$
These isomorphisms satisfy the cocyclicity condition, hence the result. If $\sigma_{ik}':\upmu_{ik}'^2 \ra \upmu_{ik}'$ and $\sigma':\upmu'^2\ra \upmu'$ are the addition maps, then the morphisms
$$
\frakB_{ik} \boxtimes \frakB_{ik} \ra \sigma_{ik}'^*\frakB_{ik},
$$
defining the $\upmu_{ik}'$-indexed algebra structure on $\frakB_{ik},$ induce
$$
\frakB_X\boxtimes \frakB_X \ra \sigma'^*\frakB_X.
$$
This morphism defines the $\upmu'$-indexed algebra structure on $\frakB_X.$
\end{proof}

\begin{lemma}\label{lemShiho209}
Keep the notation of \ref{propBglue}. For any sheaf $\F$ of $\widetilde{\calE}(X'/\frakS)_{/\upmu'}$ and any object $(U,\frakT,u)$ of $\calE(X'_{ik}/\frakS),$ we have
$$
\left (Res_{ik}^{-1}\F \right )_{(U,\frakT,u)}=\F_{\rho_{X',X_{ik}'}\left (U,\frakT,u \right )}.
$$
\end{lemma}

\begin{proof}
This follows from the obvious commutativity of
$$
\begin{tikzcd}
\text{ét}_{/U\times \upmu_U} \ar{rr}{\delta_{(U,\frakT,u)}} \ar[swap]{drr}{\delta_{\rho_{X',X_{ik}'}(U,\frakT,u)}} & & \widetilde{\calE}(X_{ik}'/\frakS)_{/\upmu_{ik}'} \ar{d}{Res_{ik}} \\
 & & \widetilde{\calE}(X'/\frakS)_{/\upmu'}.
\end{tikzcd}
$$
\end{proof}

\begin{proposition}\label{AXUESTiso}
Consider an étale $X$-scheme $U$ and the morphisms $\underline{R_{X,U}}$ and $R_{X',U'}$ \eqref{eqREST}.
Then we have canonical isomorphisms
$$
\underline{R_{X,U}}^{-1}\frakA_X \xrightarrow{\sim} \frakA_U,\ R_{X',U'}^{-1}\frakB_X \xrightarrow{\sim}\frakB_U.
$$
\end{proposition}

\begin{proof}
By \ref{propEST1831}, we have the commutative diagram
$$
\begin{tikzcd}
\begin{Bmatrix}
\mathrm{Crystals\ of\ }\widetilde{\underline{\calE}}(X/\frakS)_{/\underline{\upmu}}
\end{Bmatrix}
\ar{r}{\sim} \ar[swap]{d}{\underline{R_{X,U}}^{-1}} &
\begin{Bmatrix}
\Ox_{X,\upmu_X}\text{-}\mathrm{modules\ equipped} \\
\mathrm{with\ a\ }\calQ_{1,\upmu_X}\text{-}\mathrm{stratification}
\end{Bmatrix}
 \ar{d} \\
\begin{Bmatrix}
\mathrm{Crystals\ of\ }\widetilde{\underline{\calE}}(U/\frakS)_{/\underline{\upmu}_U}
\end{Bmatrix} 
\ar{r}{\sim} &
\begin{Bmatrix}
\Ox_{U,\upmu_U}\text{-}\mathrm{modules\ equipped} \\
\mathrm{with\ a\ }\calQ_{1|U,\upmu_U}\text{-}\mathrm{stratification}
\end{Bmatrix}.
\end{tikzcd}
$$
The crystal $\frakA_X$ is an object of the upper left category. It corresponds to the stratified $\Ox_{X,\upmu_X}$-module $\left (\calA_X,\epsilon_{\calA_X,Q}\right ).$ This is sent by the right vertical arrow to $\left (\calA_U,\epsilon_{\calA_U,Q}\right ).$ This corresponds by the lower equivalence to $\frakA_U.$ It follows that we have a canonical isomorphism
$$
\underline{R_{X,U}}^{-1}\frakA_X \xrightarrow{\sim} \frakA_U.
$$
The isomorphism for $\frakB_X$ is proved similarly.
\end{proof}

\begin{lemma}\label{lemkraz}
Let $\Delta:X\ra Y$ be the exact diagonal immersion, $\calI$ its ideal, $p_1,p_2:Y \ra X$ the canonical projections, $\ov{x} \ra X$ a geometric point, and $a_1,\hdots a_d$ generators of $\calI_{\ov{y}}$ that reduce to a basis of the $\Ox_{X,\ov{x}}$-module $\calI_{\ov{x}}/\calI_{\ov{x}}^2.$ Equip $\calI$ with a structure of $\Ox_X$-module via the first projection $p_1.$ For $I\in \N^d,$ set 
$$
a^I=\prod_{i=1}^da_i^{I_i}.
$$
Then, for any positive integer $n$ and $b\in \calI_{\ov{x}},$ there exist $c_I\in \Ox_{X,\ov{x}}$ such that
$$
b-\sum_{0<|I|\le n}c_Ia^I\in \calI_{\ov{x}}^{n+1}.
$$
\end{lemma}

\begin{proof}
The canonical morphism $Y\ra X\times_S^{\op{log}}X$ is log étale and the projections $X\times^{\op{log}}_SX\ra X$ are log smooth. It follows that the projections $Y \ra X$ are log smooth. Since they are strict, they are smooth and $\Delta$ is a regular immersion. It follows that we have a canonical isomorphism
\begin{equation}\label{SymgrKoko}
S^{\bullet}\left (\calI/\calI^2\right )\xrightarrow{\sim} \bigoplus_{k\ge 0} \calI^k/\calI^{k+1}.
\end{equation}
There exist $c_1,\hdots,c_d\in \Ox_{X,\ov{x}}$ such that
$$
b-\sum_{i=1}^dc_ia_i \in \calI_{\ov{x}}^2.
$$
By \eqref{SymgrKoko}, $\left (a^I\right )_{I\in \N^d, |I|=2}$ is a basis of the $\Ox_{X,\ov{x}}$-module $\calI_{\ov{x}}^2/\calI_{\ov{x}}^3.$ Then there exist $c_I\in \Ox_{X,\ov{x}}$ such that
$$
b-\sum_{i=1}^dc_ia_i-\sum_{I\in \N^d, |I|=2}c_Ia^I\in \calI_{\ov{x}}^3.
$$
The result is then proved by induction.
\end{proof}

\begin{proposition}\label{propkraz1932}
Assume $(H_1)$ is satisfied. Let $s\in \Gamma(U,\calM^{gp}_X)$ be a section over an étale $X$-scheme $U$ and denote by $\ov{s}$ its image in $\Gamma(U,\upmu_X).$ Let $\Delta:X\ra Y$ be the exact diagonal immersion, $\calI$ its ideal, $p_1,p_2:Y \ra X$ the canonical projections, $g:Q_1 \ra Y$ and $\iota:X\ra Q_1$ the canonical morphisms and $q_1=p_1\circ g,$ $q_2=p_2\circ g.$ Consider a geometric point $\ov{x} \ra U.$
Then, for every $m\in \calM_{X,\ov{x}},$ there exists $a \in \calI_{\Delta(\ov{x})}$ such that
\begin{equation}\label{Kokooooh}
p_{2,\Delta(\ov{x})}^{\flat}m-p_{1,\Delta(\ov{x})}^{\flat}m=\alpha_{Y,\Delta(\ov{x})}^{-1}(1+a).
\end{equation}
In addition, if $e_s$ is the basis \eqref{eqAbase} of the $\Ox_U$-module $\calA_{X,\ov{s}}$ and $b$ is a local section of $\Ox_U,$ then
\begin{equation}\label{eqepsilonAQ2}
\epsilon_{\calA_X,Q,\ov{s}}(1\otimes be_s) = \left (e_s\otimes g^{\#}\left (p_2^{\#}(b)(1+a)\right )\right ).
\end{equation}
Moreover, we suppose that the hypothesis \ref{loccoord} is satisfied and set $\eta^{I}:=\prod_{i=1}^d\eta_i^{I_i}$ and $\eta_{(q)}^J:=\prod_{i=1}^d\eta_{i(q)}^{J_i}$ for $I\in \llbracket 0,p-1 \rrbracket^d$ and $J\in \N^d.$ We have a basis $\left (\eta^I\eta_{(q)}^J \right )_{\substack{I\in \llbracket 0,p-1 \rrbracket^d \\ J\in \N^d}}$ of the $\Ox_X$-module $\calQ_1.$ Let $\left (\varphi_{I,J} \right )_{\substack{I\in \llbracket 0,p-1 \rrbracket^d \\ J\in \N^d}}$ be its dual basis. Then, for any local section $x$ of $\calA_{X,\ov{s}},$ we have
\begin{equation}\label{eqepsilonAQ}
\epsilon_{\calA_X,Q,\ov{s}}(1\otimes x) = \sum_{I\in \llbracket 0,p-1 \rrbracket^d} (\varphi_{I,0}\cdot x) \otimes \eta^I,
\end{equation}
where $\varphi_{I,0}\cdot x$ is the action of $\varphi_{I,0}$ on $x$ via the $\calQ_1$-stratification $\epsilon_{\calA_X,Q,\ov{s}}.$
\end{proposition}

\begin{proof}
Let $\upmu_X=\ov{\calM}_X^{gp}$ and $\upmu_{Q_1}=\ov{\calM}_{Q_1}^{gp}.$ The projection $q_i:Q_1 \ra X$ induces a morphism of $Q_{1,\text{ét}/q_i^{-1}\upmu_X}$
$$
q_i^{-1}\calA_X \ra \calA_{Q_1},
$$
and then a morphism of $X_{\text{ét}/\upmu_X}$
$$
\calA_X\ra \iota^{-1}\calA_{Q_1}.
$$
Taking the fibers over $\ov{s}\in \Gamma(U,\upmu_X)$ on both sides and then the stalks over a geometric point $\ov{x} \ra U,$ we get
$$
c_i:\calA_{X,\ov{s},\ov{x}} \ra \left (\iota^{-1}\calA_{Q_1}\right )_{\ov{s},\ov{x}}
$$
Consider $\ov{s}$ as an element of $\upmu_{X,\ov{x}}\xrightarrow{\sim} \upmu_{Q_1,\iota(\ov{x})}$ (the isomorphism comes from the fact that $\iota$ is strict). Recall the description $\calA_{X,\ov{s},\ov{x}}=\Ox_{X,\ov{x}}\times \calM_{X,\ov{x},\ov{s}}^{gp}/\sim$ (\eqref{eq8422b} and \eqref{eq8422c}) where $\calM_{X,\ov{x},\ov{s}}^{gp}$ is the fiber of $\calM_{X,\ov{x}}^{gp}$ over $\ov{s}.$ Similarly
$$
\left (\iota^{-1}\calA_{Q_1}\right )_{\ov{s},\ov{x}}=\left (\Ox_{Q_1,\iota(\ov{x})} \times \calM_{Q_1,\iota(\ov{x}),\ov{s}}^{gp} \right )/\sim.
$$
We have, for all $m\in \calM_{X,\ov{x},\ov{s}}^{gp},$
$$
c_i(1,m)=(1,q_{i,\iota(\ov{x})}^{\flat}m).
$$
To lighten notation, for a morphism $f$ of sheaves of $X_{\text{ét}}$ (resp. $Y_{\text{ét}},$ $Q_{1,\text{ét}}$), we denote $f_{\ov{x}}$ (resp. $f_{\Delta(\ov{x})},$ $f_{\iota(\ov{x})}$) simply by $f.$ For every $m\in \calM_{X,\ov{x}} \subset \calM_{X,\ov{x}}^{gp},$ by the exact sequence \eqref{era2exactseq}, there exists a local section $a$ of $\Delta^{-1}\calI$ such that, in $\calM_{Y,\Delta(\ov{x})}^{gp},$
\begin{equation}
p_2^{\flat}m-p_1^{\flat}m=\alpha_Y^{-1}(1+a).
\end{equation}
It follows that, in $\calM^{gp}_{Q_1,\iota(\ov{x})},$ we have
$$
q_2^{\flat}m-q_1^{\flat}m=g^{\flat}(p_2^{\flat}m-p_1^{\flat}m)=g^{\flat}\alpha_Y^{-1}(1+a)=\alpha_{Q_1}^{-1}g^{\#}(1+a).
$$
If $m$ is over $\ov{s},$ then, in $\left (\iota^{-1}\calA_{Q_1}\right )_{\ov{s},\ov{x}},$
$$
(1,q_2^{\flat}m)=(g^{\#}(1+a),q_1^{\flat}m).
$$
We finally get, in $\left (\calA_{X,\ov{s}}\otimes_{\Ox_U}\calQ_1\right )_{\ov{x}},$
$$
\epsilon_{\calA_X,Q,\ov{s}}(1\otimes(1,m)) = ((1,m)\otimes g^{\#}(1+a)).
$$
If $b$ is a local section of $\Ox_U,$ then
\begin{equation}
\begin{alignedat}{2}
\epsilon_{\calA_X,Q,\ov{s}}(1\otimes(b,s)) &= \epsilon_{\calA_X,Q,\ov{s}}(q_2^{\#}(b)\otimes(1,s)) \\
&= q_2^{\#}(b) \epsilon_{\calA_X,Q,\ov{s}}(1\otimes(1,s)) \\
&= q_2^{\#}(b) \left ((1,s)\otimes g^{\#}(1+a)\right ) \\
&= \left ((1,s)\otimes g^{\#}\left (p_2^{\#}(b)(1+a)\right )\right ).
\end{alignedat}
\end{equation}
Suppose that the hypothesis of \ref{loccoord} is satisfied. For $1\le i\le d,$ we have, in $\calQ_1,$
$$
\eta_i^p=p\eta_{i(q)}=0.
$$
Then $\calI^{pd} \Ox_{Q_1}=0.$
Therefore, by \ref{lemkraz}, $g^{\#}(a)$ is an $\Ox_X$-linear combination of $\eta^{I}:=\prod_{i=1}^d\eta_i^{I_i}$ for $I\in \llbracket 0,p-1 \rrbracket^d$ (recall that the Hopf algebra $\calQ_1$ is considered as an $\Ox_X$-module via the first projection). It follows that, for any local section $x$ of $\calA_{X,\ov{s}},$ there exist local sections $x_I \in \calA_{X,\ov{s}}$ such that
\begin{equation*}
\epsilon_{\calA_X,Q,\ov{s}}:\begin{array}[t]{clc}
\calQ_1 \otimes_{\Ox_U} \calA_{X,\ov{s}} & \ra & \calA_{X,\ov{s}} \otimes_{\Ox_U} \calQ_1 \\
1\otimes x & \mapsto & \sum_{I\in \llbracket 0,p-1 \rrbracket^d} x_I \otimes \eta^I.
\end{array}
\end{equation*}
By definition, the action $\varphi_{I,J}\cdot x$ is the image of $1\otimes x$ by the composition
$$
\calQ_1 \otimes_{\Ox_U}\calA_{X,\ov{s}} \xrightarrow{\epsilon_{\calA_X,Q,\ov{s}}} \calA_{X,Q,\ov{s}} \otimes_{\Ox_U} \calQ_1 \xrightarrow{\op{Id}\otimes \varphi_{I,J}} \calA_{X,Q,\ov{s}}.
$$
It is then clear that $x_I=\varphi_{I,0}\cdot x.$
\end{proof}

\begin{lemma}\label{lem1924Koko}
Assume $(H_1)$ and the hypothesis \ref{loccoord} are satisfied and consider the basis
$$\left (\eta^I\eta_{(q)}^J\right )_{\substack{I\in \llbracket 0,p-1\rrbracket^d\\ J\in \N^d}}$$
of the $\Ox_X$-module $\calQ_1,$ where
$$
\eta^I=\prod_{i=1}^d\eta_i^{I_i},\ \eta_{(q)}^J=\prod_{i=1}^d\eta_{(q)}^{J_i}.
$$
Let $\left (\varphi_{I,J} \right )$ be its dual basis and extend the notation $\varphi_{I,J}$ by setting $\varphi_{I,J}=0$ if $I\in \Z^d$ and $I\notin \llbracket 0,p-1\rrbracket^d.$
Consider the ring structure on $\calQ_1^{\vee}$ induced by the Hopf algebra structure on $\calQ_1$ \eqref{HopfQ}. Let $\epsilon_i$ be the multi-index of $\N^d$ whose all coefficients are zero except for the $i^{th}$ which is equal to $1.$ Then
$$
\varphi_{I,0}\cdot \varphi_{\epsilon_i,0}=I_i\varphi_{I,0}+(I_i+1)\varphi_{I+\epsilon_i,0}-\varphi_{I-(p-1)\epsilon_i,\epsilon_i}.
$$
In particular, if $J\in \llbracket 0,p-1\rrbracket^d$ and $1\le i\le d$ such that $J_i\ge 1,$ then
\begin{equation}\label{eq1925Koko}
\varphi_{J-\epsilon_i,0}\cdot \varphi_{\epsilon_i,0}=(J_i-1)\varphi_{J-\epsilon_i,0}+J_i\varphi_{J,0}.
\end{equation}
\end{lemma}

\begin{proof}
By definition, the product $\varphi_{I,0}\cdot \varphi_{\epsilon_i,0}$ is equal to the composition
$$
\calQ_1 \xrightarrow{\delta} \calQ_1 \otimes_{\Ox_X}\calQ_1 \xrightarrow{\op{Id}\otimes \varphi_{\epsilon_i,0}} \calQ_1 \xrightarrow{\varphi_{I,0}} \Ox_X,
$$
where $\delta$ is given in \ref{HopfQ}. Let $\alpha \in \llbracket 0,p-1 \rrbracket ^d$ and $\beta\in \N^d.$ Let $\calK \subset \calQ_1$ be the ideal generated by $\eta_1,\hdots,\eta_d,\eta_{1(q)},\hdots,\eta_{d(q)}.$ The linear form $\varphi_{\epsilon_i,0}:\calQ_1 \ra \Ox_X$ factors clearly through the canonical projection $\pi:\calQ_1 \ra \calQ_1/\calK^2.$ By \ref{HopfQ}, we have
\begin{alignat*}{2}
\delta \left (\eta^{\alpha}\eta_{(q)}^{\beta} \right ) &= \prod_{j=1}^d\left (1\otimes \eta_j+\eta_j\otimes 1+\eta_j\otimes \eta_j\right )^{\alpha_j} \\
& \times \prod_{j=1}^d\left (1\otimes \eta_{j(q)}+\sum_{0<b+c<p} \frac{(-1)^{b+c}}{b+c}\begin{pmatrix}b+c \\ b \end{pmatrix} \eta_j^{b+c}\otimes \eta_j^{p-b} +\eta_{j(q)}\otimes 1\right )^{\beta_j}.
\end{alignat*}
Then, in $\calQ_1 \otimes_{\Ox_X} \left ( \calQ_1/\calK^2\right ),$ we have
\begin{alignat*}{2}
\left (\pi \otimes \op{Id}_{\calQ_1} \right )\circ \delta \left (\eta^{\alpha}\eta_{(q)}^{\beta} \right ) &= \prod_{j=1}^d\left (1\otimes \eta_j+\eta_j\otimes 1+\eta_j\otimes \eta_j\right )^{\alpha_j} \\
& \times \prod_{j=1}^d\left (1\otimes \eta_{j(q)}- \eta_j^{p-1}\otimes \eta_j +\eta_{j(q)}\otimes 1\right )^{\beta_j}.
\end{alignat*}
By the binomial formula, we have
\begin{alignat*}{2}
\left (1\otimes \eta_j+\eta_j\otimes 1+\eta_j\otimes \eta_j\right )^{\alpha_j} &= \sum_{\substack{r,s,t\in \N \\ r+s+t=\alpha_j}} \frac{\alpha_j!}{r!s!t!} \eta_j^{s+t}\otimes \eta_j^{r+t}, \\
\left (1\otimes \eta_{j(q)}- \eta_j^{p-1}\otimes \eta_j +\eta_{j(q)}\otimes 1\right )^{\beta_j} &= \sum_{\substack{u,v,w\in \N \\ u+v+w=\beta_j}}\frac{\beta_j!}{u!v!w!}(-1)^v \times \left (\eta_{j(q)}^{w}\eta_j^{(p-1)v} \right )\otimes\left ( \eta_{j(q)}^{u}\eta_j^{v} \right ).
\end{alignat*}
We have the same equalities with multi-indices:
\begin{alignat*}{2}
\prod_{j=1}^d\left (1\otimes \eta_j+\eta_j\otimes 1+\eta_j\otimes \eta_j\right )^{\alpha_j}&=\sum_{\substack{r,s,t\in\N^d \\ r+s+t=\alpha}} \frac{\alpha!}{r!s!t!} \eta^{s+t}\otimes \eta^{r+t}, \\
\prod_{j=1}^d\left (1\otimes \eta_{j(q)}- \eta_j^{p-1}\otimes \eta_j +\eta_{j(q)}\otimes 1\right )^{\beta_j} &= \sum_{\substack{u,v,w\in \N^d \\ u+v+w=\beta}}\frac{\beta!}{u!v!w!}(-1)^{|v|} \times \left (\eta_{(q)}^{w}\eta^{(p-1)v} \right )\otimes\left ( \eta_{(q)}^{u}\eta^{v} \right ).
\end{alignat*}
Suppose that $\alpha_i \ge 1$ and $\beta_i\ge 1.$
Then, the only terms whose image doesn't vanish by $\op{Id}\otimes \varphi_{\epsilon_i,0}$ are those that either correspond to
$$
\begin{cases}r+t=\epsilon_i\\
u=v=0, \end{cases}\ \text{or}\ 
\begin{cases}
r+t=0\\
u=0\\
v=\epsilon_i. 
\end{cases}
$$
These cases are equivalent to
$$
\begin{cases}r=\epsilon_i, t=0, s=\alpha-\epsilon_i\\
u=v=0, w=\beta,\end{cases}, \begin{cases}t=\epsilon_i, r=0, s=\alpha-\epsilon_i\\
u=v=0, w=\beta,\end{cases}\ \text{or}\ 
\begin{cases}
r=t=0, s=\alpha\\
u=0\\
v=\epsilon_i\\
w=\beta-\epsilon_i.
\end{cases}
$$
They yield
$$
\alpha_i\left (\eta^{\alpha-\epsilon_i}\otimes \eta_i+\eta^{\alpha}\otimes \eta_i \right )\left (\eta_{(q)}^{\beta}\otimes 1\right )-\beta_i(\eta^{\alpha}\otimes 1)\left (\left (\eta^{(p-1)\epsilon_i}\eta_{(q)}^{\beta-\epsilon_i} \right )\otimes \eta_i \right ).
$$
If $\alpha_i=0$ and $\beta_i\ge 1,$ the only term that doesn't vanish by $\op{Id}\otimes \varphi_{\epsilon_i,0}$ is
$$
-\beta_i(\eta^{\alpha}\otimes 1)\left (\left (\eta^{(p-1)\epsilon_i}\eta_{(q)}^{\beta-\epsilon_i} \right )\otimes \eta_i \right ).
$$
If $\alpha_i \ge 1$ and $\beta_i=0,$ the only term that doesn't vanish by $\op{Id}\otimes \varphi_{\epsilon_i,0}$ is
$$
\alpha_i\left (\eta^{\alpha-\epsilon_i}\otimes \eta_i+\eta^{\alpha}\otimes \eta_i \right )\left (\eta_{(q)}^{\beta}\otimes 1\right ).
$$
If $\alpha_i=\beta_i=0,$ all terms vanish by $\op{Id}\otimes \varphi_{\epsilon_i,0}.$
It follows that
$$
\left (\varphi_{I,0}\cdot \varphi_{\epsilon_i,0} \right )\left (\eta^{\alpha}\eta_{(q)}^{\beta} \right )=\begin{cases}I_i+1\ \text{if}\ \alpha=I+\epsilon_i\ \text{and}\ \beta=0,\\
I_i\ \text{if}\ \alpha=I\ \text{and}\ \beta=0,\\
-1\ \text{if}\ \alpha=I-(p-1)\epsilon_i\ \text{and}\ \beta=\epsilon_i.
\end{cases}
$$
The result follows.
\end{proof}

\begin{proposition}\label{Kokoprop1925}
Assume $(H_1)$ is satisfied. The stratification $\epsilon_{\calA_X,Q}$ \eqref{eqeAQ} is $\calB_{X/S}$-linear, where $\calB_{X/S}$ is defined in \ref{Bgroupe2}.
\end{proposition}

\begin{proof}
Let $s\in \Gamma(U,\calM_X^{gp})$ be a local section of $\calM_X^{gp}$ over an étale $X$-scheme $U$ and $\ov{s}\in \Gamma(U,\upmu_X)$ its image. By definition, $\epsilon_{\calA_X,Q}$ is a morphism of $\upmu_X$-indexed algebras \eqref{remkraz1}. It is thus sufficient to prove that, for any local section $b$ of $\calB_{X/S,\ov{s}},$
$$
\epsilon_{\calA_X,Q,\ov{s}}(1\otimes b)=b\otimes 1.
$$
The stratification $\epsilon_{\calA_X,Q,\ov{s}}$ defines a $\calQ_1^{\vee}$-module structure on $\calA_{X,\ov{s}},$ given explicitly by
$$
\begin{array}[t]{clc}
\calQ_1^{\vee} \otimes_{\Ox_U} \calA_{X,\ov{s}} & \ra & \calA_{X,\ov{s}} \\
\varphi\otimes a & \mapsto & \left (\op{Id} \otimes \varphi \right )\circ \epsilon_{\calA_X,Q,\ov{s}}(1\otimes a).
\end{array}
$$
We may work étale locally on $X$ and suppose that the hypothesis \ref{loccoord} is satisfied. We take again the notation of \ref{lem1924Koko} and denote by $\varphi_{I,J}\cdot b$ the action of $\varphi_{I,J}$ on $b.$
By \eqref{eqepsilonAQ}, we have
$$
\epsilon_{\calA_X,Q,\ov{s}}(1\otimes b)=\sum_{I\in \llbracket 0,p-1\rrbracket^d}\left (\varphi_{I,0}\cdot b\right )\otimes \eta^I.
$$
We have to prove that $\varphi_{I,0}\cdot b=0$ for all $I\in\llbracket 0,p-1\rrbracket^d$ such that $|I|>0.$ By \eqref{eq1925Koko} and by induction on $|I|,$ it is sufficient to prove that $\varphi_{\epsilon_i,0}\cdot b=0$ for $1\le i\le d.$ 
Recall that $e_s=(1,s)$ \eqref{eqAbase} is a basis for the $\Ox_U$-module $\calA_{X,\ov{s}}.$ It follows that there exists a local section $t$ of $\Ox_X$ such that
$$
b=te_s.
$$
By \eqref{eq8431Koko}, we have
$$
dt+t\op{dlog}s=0.
$$
Let $X\ra Y$ be the exact diagonal immersion \eqref{prop45}, $\calI$ its ideal, $\iota:X\ra Q_1$ and $g:Q_1\ra Y$ the canonical morphisms and $p_1,p_2:Y\ra X$ and $q_1,q_2:Q_1 \ra X$ the canonical projections. By \eqref{eqepsilonAQ2} and \eqref{Kokooooh}, there exists a local section $a$ of $\Delta^{-1}\calI$ such that
$$
\begin{cases}
\epsilon_{\calA_X,Q,\ov{s}}\left (1\otimes te_s\right )=\left (e_s\otimes g^{\#}(p_2^{\#}(t)(1+a)) \right ),\\
p_2^{\flat}(s)=\alpha_Y^{-1}(1+a)+p_1^{\flat}(s).
\end{cases}
$$
We have a canonical isomorphism of $\Ox_X$-modules
$$
\Delta^{-1}\left (\calI/\calI^2\right )\xrightarrow{\sim}\omega^1_{X/S}.
$$
Let $\ov{x} \ra X$ be a geometric point.
By \ref{parag46}, the differential $\op{dlog}s \in \omega^1_{X/S,\ov{x}}$ corresponds to the class of $a\in \calI_{\Delta(\ov{x})}$ modulo $\calI_{\Delta(\ov{x})}^2.$ The local section $dt$ corresponds to the class of $p_{2,\Delta(\ov{x})}^{\#}(t)-p_{1,\Delta(\ov{x})}^{\#}(t).$ Since $dt+t\op{dlog}s=0,$
$$
p_{2,\Delta(\ov{x})}^{\#}(t)-p_{1,\Delta(\ov{x})}^{\#}(t)+p_{2,\Delta(\ov{x})}^{\#}(t)a=p_{2,\Delta(\ov{x})}^{\#}(t)(1+a)-p_{1,\Delta(\ov{x})}^{\#}(t)
$$
belongs to $\calI_{\Delta(\ov{x})}^2.$ Since
$$
\varphi_{\epsilon_i,0}(q_{1,\iota(\ov{x})}^{\#}(t))=t\varphi_{\epsilon_i,0}(1)=0,
$$
we deduce that, for $1\le i\le d,$
$$
\varphi_{\epsilon_i,0}\left (g_{\Delta(\ov{x})}^{\#}(p_{2,\Delta(\ov{x})}^{\#}(t)(1+a))\right )=0.
$$
By \ref{eqepsilonAQ2}, we have
\begin{equation*}
\epsilon_{\calA_X,Q,\ov{s}}(1\otimes b) =\epsilon_{\calA_X,Q,\ov{s}}(1\otimes te_s) = e_s\otimes g^{\#}\left (p_2^{\#}(t)(1+a)\right ).
\end{equation*}
Then,
$$
\varphi_{\epsilon_i,0}\cdot b=e_s\otimes \varphi_{\epsilon_i,0}\left (g^{\#}\left (p_2^{\#}(t)(1+a)\right )\right )=0.
$$
This finishes the proof.
\end{proof}

\begin{remark}\label{actiondepsQ}
Assume $(H_1)$ is satisfied and consider the special fiber $h:P_0 \ra Q_1$ of the morphism \eqref{wiwkraz}. The pull-back
$$
h^*\epsilon_{\calA_X,Q}:(q_2\circ h)^*\calA_X \ra (q_1\circ h)^*\calA_X
$$
is, by definition, equal to the $P_0$-stratification $\epsilon_{\calA_X}$ defined by Montagnon in \cite{Mon} 4.4.1 page 46.
\end{remark}

\begin{proposition}\label{propkraz1934}
Assume $(H_2)$ is satisfied. Consider the morphism of Hopf algebras $\calV:\calR_1' \ra \calQ_1$ induced by $\nu:Q_{\frakX}\ra R_{\frakX',1}$ \eqref{lem92}. Then, the diagram
\begin{equation}\label{diag186kraz}
\begin{tikzcd}
\begin{Bmatrix}\text{Crystals\ of\ }\widetilde{\calE}(X'/\frakS)_{/\upmu'} \end{Bmatrix} \ar{r}{C^{-1}_{X/\frakS,\upmu}} \ar[swap, sloped]{d}{\sim} & \begin{Bmatrix}\text{Crystals\ of\ }\widetilde{\underline{\calE}}(X/\frakS)_{/\underline{\upmu}} \end{Bmatrix}  \ar[sloped]{d}{\sim} \\
\begin{Bmatrix}\Ox_{X',\upmu_X}\text{-modules\ with\ an}\\ \calR_{1,\upmu_X}'\text{-stratification} \end{Bmatrix} \ar{r}{\Psi_{\upmu}} & \begin{Bmatrix}\Ox_{X,\upmu_X}\text{-modules\ with\ a}\\ \calQ_{1,\upmu_X}\text{-stratification} \end{Bmatrix},
\end{tikzcd}
\end{equation}
where the vertical equivalences are given in \ref{equivRQind}, $C_{X/\frakS,\upmu}$ is given in \eqref{diagkraz2} and $\Psi_{\upmu}(\calE',\epsilon') = \left (F_1^*\calE',\calV^*\epsilon'\right )$ for an $\Ox_{X',\upmu_X}$-module $\calE'$ equipped with an $\calR_{1,\upmu_X}'$-stratification $\epsilon',$ is 2-commutative.
\end{proposition}

\begin{proof}
Let $U$ be an étale $X$-scheme and $s\in \Gamma(U,\upmu_X)=\Gamma(U',F_{1*}\upmu_X).$ Consider the diagram
$$
\begin{tikzcd}
\begin{Bmatrix}\text{Crystals\ of}\\ \widetilde{\calE}(X'/\frakS)_{/\upmu'} \end{Bmatrix} \ar{rr}{C^{-1}_{X/\frakS,\upmu}} \ar[swap, sloped]{dd}{\sim} \ar{dr}{\alpha_s'^*} & &  \begin{Bmatrix}\text{Crystals\ of} \\ \widetilde{\underline{\calE}}(X/\frakS)_{/\underline{\upmu}} \end{Bmatrix}  \ar[pos=0.20,dashed,sloped]{dd}{\sim} \ar{dr}{\underline{\alpha}_s^*} & \\
 & \calC(U'/\frakS) \ar[pos=0.20]{rr}{C_{U/\frakS}^{-1}} \ar[pos=0.20,swap,sloped]{dd}{\sim} & & \underline{\calC}(U/\frakS) \ar[sloped] {dd}{\sim} \\
\begin{Bmatrix}\Ox_{X',\upmu_X}\text{-mod}\\ \text{with\ an}\\ \calR_{1,\upmu_X}'\text{-strat.} \end{Bmatrix} \ar[pos=0.20,dashed]{rr}{\Psi_{\upmu}} \ar{dr} & & \begin{Bmatrix}\Ox_{X,\upmu_X}\text{-mod}\\ \text{with\ a}\\ \calQ_{1,\upmu_X}\text{-strat.} \end{Bmatrix} \ar[dashed]{dr} & \\
& \begin{Bmatrix}\Ox_{U'}\text{-mod}\\ \text{with\ an}\\ \calR_{1|U}'\text{-strat.} \end{Bmatrix} \ar{rr}{\Psi_{0}} & & \begin{Bmatrix}\Ox_{U}\text{-mod}\\ \text{with\ a}\\ \calQ_{1|U}\text{-strat.} \end{Bmatrix},
\end{tikzcd}
$$
where $\alpha_s'$ and $\underline{\alpha}_s$ are given in \eqref{alphakraz}, $\calC(U'/\frakS)$ and $\underline{\calC}(U/\frakS)$ are defined in \ref{defcrys}, $\Psi_0$ is given in \ref{thmKoko1926} and the lower oblique lines take a module with a stratification to its fiber over $s.$
The front square is commutative by \ref{thmKoko1926}. The upper square is commutative by definition of $C^{-1}_{X/\frakS,\upmu}$ \eqref{diagkraz19253}. The two side squares are commutative by construction of the equivalences of \ref{equivRQind} and the lower square is clearly commutative.
This being true for all local sections $s$ of $\upmu_X,$ the commutativity of the back square then follows from the fact that two $\Ox_{X,\upmu_X}$-modules with $\calQ_{1,\upmu_X}$-stratifications $(\calE,\epsilon)$ and $(\calE',\epsilon')$ are equal if and only if, for any local section $s$ of $\upmu_X,$ $(\calE_s,\epsilon_s)=(\calE'_s,\epsilon'_s)$ \eqref{indremark73}.
\end{proof}

\begin{proposition}
There exists a canonical morphism
$$C_{X/\frakS,\upmu}^{-1}\frakB_X \ra \frakA_X,$$
where $\frakB_X$ and $\frakA_X$ are defined in \ref{defQstrat} and \ref{propBglue} and $C_{X/\frakS,\upmu}$ is defined in \eqref{diagkraz2}.
\end{proposition}

\begin{proof}
Assume $(H_2)$ is satisfied. Consider the morphism of Hopf algebras $\calV:\calR_1' \ra \calQ_1$ induced by the morphism $\nu:Q_{\frakX}\ra R_{\frakX',1}$ given in \ref{lem92}. By \ref{propkraz1934}, we have a 2-commutative diagram
\begin{equation}\label{diag186}
\begin{tikzcd}
\begin{Bmatrix}\text{Crystals\ of\ }\widetilde{\calE}(X'/\frakS)_{/\upmu'} \end{Bmatrix} \ar{r}{C^{-1}_{X/\frakS,\upmu}} \ar[swap, sloped]{d}{\sim} & \begin{Bmatrix}\text{Crystals\ of\ }\widetilde{\underline{\calE}}(X/\frakS)_{/\underline{\upmu}} \end{Bmatrix}  \ar[sloped]{d}{\sim} \\
\begin{Bmatrix}\Ox_{X',\upmu_X}\text{-modules\ with\ an}\\ \calR_{1,\upmu_X}'\text{-stratification} \end{Bmatrix} \ar{r}{\Psi_{\upmu}} & \begin{Bmatrix}\Ox_{X,\upmu_X}\text{-modules\ with\ a}\\ \calQ_{1,\upmu_X}\text{-stratification} \end{Bmatrix}
\end{tikzcd}
\end{equation}
The lower arrow $\Psi_{\upmu}$ is defined by
$$\Psi_{\upmu}(\calE',\epsilon')=\left (F_1^*\calE',\calV^*\epsilon'\right )=\left (\Ox_{X,\upmu_X}\otimes_{\Ox_{X',\upmu_X}}\calE',\calV^*\epsilon'\right ).$$
Consider the morphism
$$
\sigma_R:\begin{array}[t]{clc}
\calR_{1,\upmu_X}' \otimes_{\Ox_{X',\upmu_X}} \calB_{X/S} & \ra & \calB_{X/S} \otimes_{\Ox_{X',\upmu_X}} \calR_{1,\upmu_X}'\\
a\otimes b & \mapsto & b\otimes a.
\end{array}
$$
Clearly, $\sigma_R$ is an $\calR_{1,\upmu_X}'$-stratification on the $\Ox_{X',\upmu_X}$-module $\calB_{X/S}.$
By \ref{qF}, we have a canonical isomorphism
$$
\calQ_{1,\upmu_X} \otimes_{\Ox_{X',\upmu_X}}\calB_{X/S} 
\xrightarrow{\sim} \calB_{X/S}\otimes_{\Ox_{X',\upmu_X}}\calQ_{1,\upmu_X}.
$$
Consider the composition
\begin{equation}\label{eqsigmaQkraz}
\begin{alignedat}{2}
\sigma_Q:
\calQ_{1,\upmu_X} \otimes_{\Ox_{X,\upmu_X}} (F_1^*\calB_{X/S}) &\xrightarrow{\sim} \calQ_{1,\upmu_X} \otimes_{\Ox_{X',\upmu_X}}\calB_{X/S} \\
&\xrightarrow{\sim} \calB_{X/S}\otimes_{\Ox_{X',\upmu_X}}\calQ_{1,\upmu_X} \\
&\xrightarrow{\sim} (F_1^*\calB_{X/S}) \otimes_{\Ox_{X,\upmu_X}} \calQ_{1,\upmu_X}.
\end{alignedat}
\end{equation}
Clearly, $\sigma_Q$ is a $\calQ_{1,\upmu_X}$-stratification on the $\Ox_{X,\upmu_X}$-module $(F_1^*\calB_{X/S})$ and we have
$$
\Psi_{\upmu}\left (\calB_{X/S},\sigma_R\right )=\left (F_1^*\calB_{X/S},\sigma_Q\right ).
$$
Note that $\Psi_{\upmu}\left (\calB_{X/S},\sigma_R\right )$ is independant of the choice of the lifting $F.$
By \ref{Kokoprop1925}, the diagram
$$
\begin{tikzcd}
\calQ_{1,\upmu_X} \otimes_{\Ox_{X',\upmu_X}} \calB_{X/S} \ar{d} \ar{r}{\sigma_Q} & \calB_{X/S} \otimes_{\Ox_{X',\upmu_X}} \calQ_{1,\upmu_X} \ar{d} \\
\calQ_{1,\upmu_X}\otimes_{\Ox_{X,\upmu_X}} \calA_X \ar{r}{\epsilon_{\calA_X,Q}} & \calA_X \otimes_{\Ox_{X,\upmu_X}}\calQ_{1,\upmu_X},
\end{tikzcd}
$$
where the vertical arrows are the canonical ones, is commutative.
It follows that the canonical morphism
$$\Ox_{X,\upmu_X}\otimes_{\Ox_{X',\upmu_X}}\calB_{X/S} \ra \calA_X$$
is compatible with the $\calQ_{1,\upmu_X}$-stratifications. We deduce a canonical morphism, independant of the lifting $F,$ of $\upmu$-indexed algebras
$$C_{X/\frakS,\upmu}^{-1}\frakB_X \ra \frakA_X.$$
\end{proof}

\begin{parag}
The functor \eqref{diagkraz2}
$$
C^{-1}_{X/\frakS,\upmu}:\begin{Bmatrix}\text{Crystals\ of\ }\widetilde{\calE}(X'/\frakS)_{/\upmu'} \end{Bmatrix} \ra  \begin{Bmatrix}\text{Crystals\ of\ }\widetilde{\underline{\calE}}(X/\frakS)_{/\underline{\upmu}} \end{Bmatrix}
$$
induces a functor
\begin{equation}\label{CESTeq}
C^{-1}_{X/\frakS,\upmu}:\begin{Bmatrix}\text{Crystals\ of\ }\upmu'\text{-indexed}\\ \frakB_X\text{-modules\ of\ }\widetilde{\calE}(X'/\frakS)_{/\upmu'} \end{Bmatrix} \ra \begin{Bmatrix}\text{Crystals\ of\ }\underline{\upmu}\text{-indexed}\\ C_{X/\frakS,\upmu}^{-1}\frakB_X\text{-modules\ of\ }\widetilde{\underline{\calE}}(X/\frakS)_{/\underline{\upmu}} \end{Bmatrix}.
\end{equation}
By \ref{propcrystensstar}, we have a well-defined linearization functor
\begin{equation}\label{ffL}
L_{\mathrm{crys}}:\begin{array}[t]{clc}
\begin{Bmatrix}\text{Crystals\ of\ }\underline{\upmu}\text{-indexed}\\ C_{X/\frakS,\upmu}^{-1}\frakB_X\text{-modules\ of\ }\widetilde{\underline{\calE}}(X/\frakS)_{/\underline{\upmu}} \end{Bmatrix} & \ra & \begin{Bmatrix}\text{Crystals\ of\ }\underline{\upmu}\text{-indexed}\\ \frakA_X\text{-modules\ of\ }\widetilde{\underline{\calE}}(X/\frakS)_{/\underline{\upmu}} \end{Bmatrix} \\
\mathfrak{E} & \mapsto & \frakA_X \circledast_{C_{X/\frakS,\upmu}^{-1}\frakB_X}\mathfrak{E}.
\end{array}
\end{equation}
Set
\begin{equation}\label{ffC}
C^{-1}=L_{\mathrm{crys}}\circ C_{X/\frakS,\upmu}^{-1}:\begin{array}[t]{clc}
\begin{Bmatrix}\text{Crystals\ of\ }\upmu'\text{-indexed}\\ \frakB_X\text{-modules\ of\ }\widetilde{\calE}(X'/\frakS)_{/\upmu'} \end{Bmatrix} & \ra & \begin{Bmatrix}\text{Crystals\ of\ }\underline{\upmu}\text{-indexed}\\ \frakA_X\text{-modules\ of\ }\widetilde{\underline{\calE}}(X/\frakS)_{/\underline{\upmu}} \end{Bmatrix} \\
\mathfrak{E} & \mapsto & \frakA_X \circledast_{C_{X/\frakS,\upmu}^{-1}\frakB_X}C_{X/\frakS,\upmu}^{-1}\mathfrak{E}.
\end{array}
\end{equation}
\end{parag}

\begin{definition}~  
\begin{enumerate}
\item A $\upmu_X$-indexed $\calB_{X/S}$-module $\calE'$ equipped with an $\calR_{1,\upmu_X}'$-stratification
$$\epsilon':\calR_{1,\upmu_X}'\otimes_{\Ox_{X',\upmu_X}}\calE' \ra \calE'\otimes_{\Ox_{X',\upmu_X}}\calR_{1,\upmu_X}'$$
is said to be \emph{admissible} if $\epsilon'$ is $\calB_{X/S}$-linear.
\item Let $\frakA$ be a $\upmu_X$-indexed algebra equipped with a $\calQ_{1,\upmu_X}$-stratification $\epsilon_{\frakA}$ and $\calE$ a $\upmu_X$-indexed $\frakA$-module equipped with a $\calQ_{1,\upmu_X}$-stratification
$$\epsilon:\calQ_{1,\upmu_X}\otimes_{\Ox_{X,\upmu_X}}\calE \ra \calE\otimes_{\Ox_{X,\upmu_X}}\calQ_{1,\upmu_X}.$$
We have a morphism
\begin{alignat*}{2}
&\epsilon_{\frakA}\boxtimes \op{Id}_{\calE}=p_1^*\epsilon_{\frakA}\otimes p_2^*\op{Id}_{\calE}:\\ 
&\calQ_{1,\upmu_X^2}\otimes_{\Ox_{X,\upmu_X^2}} p_1^*\frakA \otimes_{\Ox_{X,\upmu_X^2}}p_2^*\calE &\ra p_1^*\frakA \otimes_{\Ox_{X,\upmu_X^2}}\calQ_{1,\upmu_X^2}\otimes_{\Ox_{X,\upmu_X^2}}p_2^*\calE,
\end{alignat*}
which identifies with
$$
\epsilon_{\frakA}\boxtimes \op{Id}_{\calE}:\calQ_{1,\upmu_X^2}\otimes_{\Ox_{X,\upmu_X^2}}\left (\frakA \boxtimes  \calE\right ) \ra \frakA \boxtimes \left (\calQ_{1,\upmu_X}\otimes_{\Ox_{X,\upmu_X}}\calE\right ).
$$
Similarly, we have a morphism
$$
\op{Id}_{\frakA} \boxtimes \epsilon:\frakA \boxtimes (\calQ_{1,\upmu_X}\otimes_{\Ox_{X,\upmu_X}} \calE) \ra (\frakA \boxtimes \calE)\otimes_{\Ox_{X,\upmu_X^2}}\calQ_{1,\upmu_X^2}.
$$
The stratification $\epsilon$ is said to be \emph{$\frakA$-admissible} if it is compatible with $\epsilon_{\frakA}$ i.e. if the diagram
\begin{equation}\label{admindQstratdiag}
\begin{tikzcd}
\calQ_{1,\upmu_X^2}\otimes_{\Ox_{X,\upmu_X^2}}(\frakA\boxtimes \calE) \ar{d}{\epsilon_{\frakA}\boxtimes \op{Id}_{\calE}} \ar{rr}{\op{Id}_{\calQ_{1,\upmu_X^2}}\otimes \pi} & & \calQ_{1,\upmu_X^2}\otimes_{\Ox_{X,\upmu_X^2}} \sigma_X^*\calE \ar{dd}{\sigma_X^*\epsilon} \\
\frakA \boxtimes (\calQ_{1,\upmu_X}\otimes_{\Ox_{X,\upmu_X}} \calE) \ar{d}{\op{Id}_{\frakA}\boxtimes \epsilon} & & \\
(\frakA \boxtimes \calE)\otimes_{\Ox_{X,\upmu_X^2}}\calQ_{1,\upmu_X^2}\ar{rr}{\pi\otimes \op{Id}_{\calQ_{1,\upmu_X^2}}} & & \sigma_X^*\calE\otimes_{\Ox_{X,\upmu_X^2}}\calQ_{1,\upmu_X^2}
\end{tikzcd}
\end{equation}
is commutative, where $\sigma_X:\upmu_X^2 \ra \upmu_X$ is the addition map and $\pi:\frakA \boxtimes \calE \ra \sigma_X^*\calE$ is the morphism defining the $\upmu_X$-indexed $\frakA$-module structure on $\calE.$
If $\frakA=\calA_X$ and $\epsilon_{\frakA}=\epsilon_{\calA_X,Q}$ \eqref{defQstrat}, then we say that $\epsilon$ is \emph{admissible}.
\end{enumerate}
\end{definition}

\begin{proposition}\label{equivRQAB}
Assume $(H_1)$ is satisfied. The canonical equivalence of categories between the category of crystals of $\widetilde{\underline{\calE}}(X/\frakS)_{/\underline{\upmu}}$ (resp. $\widetilde{\calE}(X'/\frakS)_{/\upmu'}$) and the category of $\Ox_{X,\upmu_X}$-modules equipped with a $\calQ_{1,\upmu_X}$-stratification (resp. $\Ox_{X',\upmu_X}$-modules equipped with an $\calR_{1,\upmu_X}'$-stratification), given in \ref{equivRQind}, induces a canonical equivalence of categories between the full subcategory of crystals of $\upmu$-indexed $\frakA_X$-modules of $\widetilde{\underline{\calE}}(X/\frakS)_{/\underline{\upmu}}$ (resp. crystals of $\upmu'$-indexed $\frakB_X$-modules of $\widetilde{\calE}(X'/\frakS)_{/\upmu'}$) and the full subcategory of $\upmu_X$-indexed $\calA_X$-modules equipped with an admissible $\calQ_{1,\upmu_X}$-stratification (resp. $\upmu_X$-indexed $\calB_{X/S}$-modules equipped with an admissible $\calR_{1,\upmu_X}'$-stratification).
\end{proposition}

\begin{proof}
Let $\frakE$ be a crystal of $\widetilde{\underline{\calE}}(X/\frakS)_{/\underline{\upmu}}$ and $(\F,\epsilon)$ the corresponding $\Ox_{X,\upmu_X}$-module with a $\calQ_{1,\upmu_X}$-stratification. We take again the notation of \ref{parag1419}. To prove the proposition, it is sufficient to prove that $\frakE$ is equipped with a structure of $\upmu$-indexed $\frakA_X$-module if and only if $\epsilon$ is admissible. Suppose that $\frakE$ is a $\upmu$-indexed $\frakA_X$-module. By \ref{parag1419}, $\F=\frakE_{(X,\frakX,i_X)},$ where $i_X:\underline{X} \ra X$ is the canonical immersion, and $\epsilon$ is the composition
$$\epsilon:\calQ_{1,\upmu_X} \otimes_{\Ox_{X,\upmu_X}} \mathfrak{E}_{(X,\frakX,i_X)} \xrightarrow{c_{\mathfrak{E},q_2}} \mathfrak{E}_{(X,Q_{\frakX},\lambda_Q)} \xrightarrow{c_{\mathfrak{E},q_1}^{-1}} \mathfrak{E}_{(X,\frakX,i_X)}\otimes_{\Ox_{X,\upmu_X}}\calQ_{1,\upmu_X}.$$
By \ref{prop1415} and \eqref{eq14201}, $\F$ is a $\upmu_X$-indexed $\calA_X$-module. Let $\sigma:\upmu^2\ra \upmu$ and $\sigma_X:\upmu_X^2 \ra \upmu_X$ be the addition maps and $\frakA_X \boxtimes \frakE \ra \sigma^*\frakE$ the morphism defining the $\upmu$-indexed $\frakA_X$-algebra structure on $\frakE.$ This morphism and the projections $q_1$ and $q_2$ yield the commutative diagram
$$
\begin{tikzcd}
\calQ_{1,\upmu_X^2} \otimes_{\Ox_{X,\upmu_X^2}}\left (\frakA_X \boxtimes \frakE \right )_{(X,\frakX,i_X)} \ar[swap]{d}{c_{\frakA_X\boxtimes \frakE,q_2}} \ar{r} & \calQ_{1,\upmu_X^2} \otimes_{\Ox_{X,\upmu_X^2}}\left ( \left (\sigma^*\frakE\right )_{(X,\frakX,i_X)}\right ) \ar{d}{c_{\sigma^*\frakE,q_2}} \\
\left (\frakA_X\boxtimes\frakE \right )_{(X,Q_{\frakX},\lambda_Q)} \ar{r} \ar[swap]{d}{c_{\frakA_X\boxtimes \frakE,q_1}^{-1}} & \left (\sigma^*\frakE \right )_{(X,Q_{\frakX},\lambda_Q)} \ar{d}{c_{\sigma^*\frakE,q_1}^{-1}} \\
\left (\frakA_X \boxtimes \frakE \right )_{(X,\frakX,i_X)} \otimes_{\Ox_{X,\upmu_X^2}}\calQ_{1,\upmu_X^2}  \ar{r} & \left ( \left (\sigma^*\frakE\right )_{(X,\frakX,i_X)}\right )\otimes_{\Ox_{X,\upmu_X^2}}\calQ_{1,\upmu_X^2}.
\end{tikzcd}
$$
This proves that $\epsilon$ is admissible.

Conversly, suppose that $\F$ is a $\upmu_X$-indexed $\calA_X$-module and that $\epsilon$ is admissible. The stratifications $\epsilon$ and $\epsilon_{\calA_X,Q}$ \eqref{eqeAQ} canonically induce stratifications on $\calA_X\boxtimes \F$ and $\sigma_X^*\F.$ The crystals corresponding to $\calA_X \boxtimes \F$ and $\sigma_X^*\F$ are $\frakA_X \boxtimes \frakE$ and $\sigma^*\frakE$ respectively. The admissibility of $\epsilon$ proves that the morphism $\calA_X\boxtimes \F \ra \sigma_X^*\F,$ defining the $\upmu_X$-indexed $\calA_X$-algebra structure on $\F,$ is compatible with the stratifications of $\calA_X\boxtimes \F$ and $\F$ and hence induces, by \ref{equivRQind}, a morphism of crystals
$$
\frakA_X \boxtimes \frakE \ra \sigma^*\frakE.
$$
This yields the desired $\upmu$-indexed $\frakA_X$-module structure on $\frakE.$ These two constructions are quasi-inverse to each other. The proof for crystals of $\widetilde{\calE}(X'/\frakS)_{/\upmu'}$ is similar.
\end{proof}

\begin{lemma}\label{ommek1}
Let $\calP$ be a presheaf of $\Ox_X$-modules, $\calE$ an $\Ox_X$-module and $\F$ the presheaf
$$
\F:U\mapsto \calE(U) \otimes_{\Ox_X(U)}\calP(U).
$$
Then, there exists a canonical isomorphism
$$
a\F \xrightarrow{\sim} \calE \otimes_{\Ox_X} a\calP,
$$
where $a$ is the associated sheaf functor.
\end{lemma}

\begin{proof}
The canonical morphism of presheaves $\calP \ra a\calP$ induces a canonical morphism of presheaves
$$
 \F \ra \calE \otimes_{\Ox_X}a\calP,
$$
and then a canonical morphism of $\Ox_X$-modules
\begin{equation}\label{omek2}
a\F \ra \calE \otimes_{\Ox_X} a\calP.
\end{equation}
Conversly, since $a$ is left exact (\cite{SGA43} II 3.4), the canonical $\Ox_X$-bilinear morphism of presheaves
$$
\calE \times  \calP \ra \F \ra a\F
$$
factors through $\calE \times a\calP.$ It then induces a canonical $\Ox_X$-linear morphism
\begin{equation}\label{omek3}
\calE\otimes_{\Ox_X}a\calP \ra a\F.
\end{equation}
By checking the stalks, the morphisms \eqref{omek2} and \eqref{omek3} are inverse to each other.
\end{proof}

\begin{proposition}\label{propdraft1}
Assume $(H_1)$ is satisfied. Let $(\calE_1,\epsilon_1)$ and $(\calE_2,\epsilon_2)$ be two $\Ox_{X,\upmu_X}$-modules equipped with $\calQ_{1,\upmu_X}$-stratifications.
Then $\epsilon=p_1^*\epsilon_1 \otimes_{\calQ_{1,\upmu_X^2}} p_2^*\epsilon_2$ (note that we tensor over $\calQ_{1,\upmu_X^2}$) is a $\calQ_{1,\upmu_X^2}$-stratification on $\calE=\calE_1 \boxtimes \calE_2.$
\end{proposition}

\begin{proof}
First, we have a commutative diagram
$$
\begin{tikzcd}
\left (\calQ_{1,\upmu_X^2} \otimes_{\Ox_{X,\upmu_X^2}} p_1^*\calE_1 \right ) \otimes_{\calQ_{1,\upmu_X^2}} \left (\calQ_{1,\upmu_X^2} \otimes_{\Ox_{X,\upmu_X^2}} p_2^*\calE_2 \right ) \ar{d}{\epsilon}  \ar{rr}{\sim} & & \calQ_{1,\upmu_X^2} \otimes_{\Ox_{X,\upmu_X^2}} \left (\calE_1 \boxtimes \calE_2 \right ) \ar{d}\\
\left (p_1^*\calE_1 \otimes_{\Ox_{X,\upmu_X^2}} \calQ_{1,\upmu_X^2} \right ) \otimes_{\calQ_{1,\upmu_X^2}} \left (p_2^*\calE_2 \otimes_{\Ox_{X,\upmu_X^2}} \calQ_{1,\upmu_X^2} \right ) \ar{rr}{\sim} & & \left (\calE_1 \boxtimes \calE_2 \right )\otimes_{\Ox_{X,\upmu_X^2}} \calQ_{1,\upmu_X^2},
\end{tikzcd}
$$
where the horizontal isomorphisms are the canonical ones. This way, we consider $\epsilon$ as an isomorphism
$$
\epsilon:\calQ_{1,\upmu_X^2} \otimes_{\Ox_{X,\upmu_X^2}} \left (\calE_1 \boxtimes \calE_2 \right ) \ra \left (\calE_1 \boxtimes \calE_2 \right )\otimes_{\Ox_{X,\upmu_X^2}} \calQ_{1,\upmu_X^2}.
$$
To prove that it is a stratification, the only difficult part is to check the cocyclicity condition. The cocyclicity condition satisfied by $\epsilon_i$ implies the commutativity of the diagram
$$
\begin{tikzcd}
\calQ_{1,\upmu_X^2} \otimes \calQ_{1,\upmu_X^2} \otimes p_i^* \calE_i \ar{rr}{\delta_{\upmu_X^2}^*p_i^*\epsilon_i} \ar[swap]{dr}{\op{Id}\otimes p_i^*\epsilon_i} & & p_i^*\calE_i \otimes \calQ_{1,\upmu_X^2}\otimes \calQ_{1,\upmu_X^2}, \\
 & \calQ_{1,\upmu_X^2} \otimes p_i^*\calE_i \otimes \calQ_{1,\upmu_X^2} \ar[swap]{ur}{p_i^*\epsilon_i \otimes \op{Id}} &
\end{tikzcd}
$$
where the tensor products are over $\Ox_{X,\upmu_X^2}$ and $\delta$ is defined in \ref{HopfQ}. Tensoring both diagrams for $i=1,2$ over $\calQ_{1,\upmu_X^2},$ we get the cocycle condition for $\epsilon.$
\end{proof}

\begin{proposition}\label{propdraft2}
Assume $(H_1)$ is satisfied. Let $(\calE,\epsilon)$ be an $\Ox_{X,\upmu_X^2}$-module equipped with a $\calQ_{1,\upmu_X^2}$-stratification. Then $\sigma_!\epsilon$ is a $\calQ_{1,\upmu_X}$-stratification on $\sigma_!\calE.$ 
\end{proposition}

\begin{proof}
By (\cite{SGA43} 12.11 b), we have canonical isomorphisms
$$
\sigma_! \left (\calQ_{1,\upmu_X} \otimes  \calE \right ) \xrightarrow{\sim} \calQ_{1,\upmu_X} \otimes  \sigma_!\calE,\ \sigma_! \left (\calE \otimes \calQ_{1,\upmu_X} \right ) \xrightarrow{\sim} \sigma_!\calE \otimes \calQ_{1,\upmu_X}.
$$
Via these isomorphisms, $\sigma_!\epsilon$ identifies with an isomorphism
$$
\sigma_!\epsilon:\calQ_{1,\upmu_X} \otimes  \sigma_!\calE \ra \sigma_!\calE \otimes \calQ_{1,\upmu_X}.
$$
We start by checking the cocycle condition: applying $\sigma_!$ to the cocyclicity diagram of $\epsilon,$ we get a commutative diagram
$$
\begin{tikzcd}
\sigma_! \left (\calQ_{1,\upmu_X^2} \otimes \calQ_{1,\upmu_X^2} \otimes  \calE \right )\ar{rr}{\sigma_!\delta_{\upmu_X^2}^*\epsilon} \ar[swap]{dr}{\sigma_! \left (\op{Id}\otimes \epsilon\right )} & & \sigma_! \left (\calE \otimes \calQ_{1,\upmu_X^2}\otimes \calQ_{1,\upmu_X^2}\right ), \\
 & \sigma_! \left (\calQ_{1,\upmu_X^2} \otimes \calE \otimes \calQ_{1,\upmu_X^2}\right ) \ar[swap]{ur}{\sigma_! \left (\epsilon \otimes \op{Id}\right )} &
\end{tikzcd}
$$
where the tensor products are over $\Ox_{X,\upmu_X^2}$ and $\delta$ is defined in \ref{HopfQ}. By (\cite{SGA43} 12.11 b), we have a canonical isomorphism
$$
\sigma_! \left (\calQ_{1,\upmu_X^2} \otimes \calQ_{1,\upmu_X^2} \otimes  \calE \right ) \xrightarrow{\sim} \calQ_{1,\upmu_X} \otimes \calQ_{1,\upmu_X} \otimes  \sigma_!\calE,
$$
where the tensor products in the right-hand side are over $\Ox_{X,\upmu_X}.$ By similar isomorphisms for the other terms of the diagram, we deduce the commutativity of the diagram
$$
\begin{tikzcd}
\calQ_{1,\upmu_X} \otimes \calQ_{1,\upmu_X} \otimes  \sigma_!\calE \ar{rr}{\sigma_!\delta_{\upmu_X^2}^*\epsilon} \ar[swap]{dr}{\op{Id}\otimes \sigma_!\epsilon} & & \sigma_!\calE \otimes \calQ_{1,\upmu_X}\otimes \calQ_{1,\upmu_X}. \\
 & \calQ_{1,\upmu_X} \otimes \sigma_!\calE \otimes \calQ_{1,\upmu_X} \ar[swap]{ur}{\sigma_!\epsilon \otimes \op{Id}} &
\end{tikzcd}
$$
It remains to prove that $\sigma_!\delta_{\upmu_X^2}^*\epsilon$ identifies with $\delta_{\upmu_X}^*\sigma_!\epsilon.$
We have a commutative diagram of ringed topoi
$$
\begin{tikzcd}
\left (X_{\text{ét}/\upmu_X^2}, \calQ_{1,\upmu_X^2}\otimes \calQ_{1,\upmu_X^2} \right ) \ar[swap,bend right=80]{dd}{j_{\upmu_X^2}} \ar{rr}{\delta_{\upmu_X^2}} \ar[swap]{d}{\sigma} & & \left (X_{\text{ét}/\upmu_X^2}, \calQ_{1,\upmu_X^2} \right ) \ar{d}{\sigma} \ar[bend right=-80]{dd}{j_{\upmu_X^2}} \\
\left (X_{\text{ét}/\upmu_X}, \calQ_{1,\upmu_X} \otimes \calQ_{1,\upmu_X} \right ) \ar[swap]{d}{j_{\upmu_X}} \ar{rr}{\delta_{\upmu_X}} & & \left (X_{\text{ét}/\upmu_X}, \calQ_{1,\upmu_X} \right ) \ar{d}{j_{\upmu_X}} \\
\left (X_{\text{ét}}, \calQ_{1}\otimes \calQ_{1} \right ) \ar{rr}{\delta} & & \left (X_{\text{ét}}, \calQ_{1} \right ).
\end{tikzcd}
$$
By \ref{lemtop15} (2), $\sigma_!\delta_{\upmu_X^2}^*\epsilon$ identifies with $\delta_{\upmu_X}^*\sigma_!\epsilon.$ The cocyclicity condition of $\sigma_!\epsilon$ is hence proved. To prove that $\sigma_!\epsilon$ is a stratification, it remains to check that
$$
\pi_{\upmu_X}^*\sigma_!\epsilon=\op{Id}_{\sigma_!\calE},
$$
where $\pi:\Ox_X\ra \calQ_1$ is the counit of the Hopf algebra $\calQ_1.$ For that, it is sufficient to prove that $\pi_{\upmu_X}^*\sigma_!\epsilon$ identifies with $\sigma_!\pi_{\upmu_X^2}^*\epsilon.$ This is again proved using \ref{lemtop15} (2) and the commutative diagram
$$
\begin{tikzcd}
\left (X_{\text{ét}/\upmu_X^2},  \calQ_{1,\upmu_X^2} \right ) \ar{rr}{\pi_{\upmu_X^2}} \ar[swap]{d}{\sigma} & & \left (X_{\text{ét}/\upmu_X^2}, \Ox_{X,\upmu_X^2} \right ) \ar{d}{\sigma} \\
\left (X_{\text{ét}/\upmu_X}, \calQ_{1,\upmu_X} \right ) \ar{rr}{\pi_{\upmu_X}} & & \left (X_{\text{ét}/\upmu_X}, \Ox_{X,\upmu_X} \right ).
\end{tikzcd}
$$
\end{proof}

\begin{corollaire}\label{cordraft1}
Assume $(H_1)$ is satisfied. Let $(\calE_1,\epsilon_1)$ and $(\calE_2,\epsilon_2)$ be two $\Ox_{X,\upmu_X}$-modules equipped with $\calQ_{1,\upmu_X}$-stratifications and $\epsilon=p_1^*\epsilon_1 \otimes_{\calQ_{1,\upmu_X^2}} p_2^*\epsilon_2$ the stratification on $\calE=\calE_1 \boxtimes \calE_2$ as in \ref{propdraft1}. Then $\sigma_!\epsilon$ is a $\calQ_{1,\upmu_X}$-stratification on $\sigma_!\calE.$
\end{corollaire}

\begin{proof}
This is a consequence of \ref{propdraft1} and \ref{propdraft2}.
\end{proof}

\begin{proposition}\label{propdraft3}
Assume $(H_1)$ is satisfied. Let $\frakA$ a $\upmu_X$-indexed algebra and $(\calE_1,\epsilon_1)$ and $(\calE_2,\epsilon_2)$ two $\upmu_X$-indexed $\frakA$-modules equipped with $\calQ_{1,\upmu_X}$-stratifications and $\epsilon=p_1^*\epsilon_1 \otimes_{\calQ_{1,\upmu_X^2}} p_2^*\epsilon_2$ the stratification on $\calE=\calE_1 \boxtimes \calE_2$ as in \ref{propdraft1}. If $\epsilon_1$ and $\epsilon_2$ are $\frakA$-linear, then the stratification $\sigma_!\epsilon$ \eqref{cordraft1} induces a $\calQ_{1,\upmu_X}$-stratification
$$
\epsilon_1 \circledast_{\calQ_{1,\upmu_X}\otimes_{\Ox_{X,\upmu_X}}\frakA}\epsilon_2:
\calQ_{1,\upmu_X} \otimes_{\Ox_{X,\upmu_X}}\left (\calE_1\circledast_{\frakA}\calE_2 \right ) \ra \left (\calE_1\circledast_{\frakA}\calE_2 \right )\otimes_{\Ox_{X,\upmu_X}}\calQ_{1,\upmu_X}
$$
on $\calE_1\circledast_{\frakA}\calE_2.$
\end{proposition}

\begin{proof}
Let $s\in \Gamma(U,\upmu_X)$ be a local section over an étale $X$-scheme $U.$ By \ref{indpropshriek} (3), the fiber $\sigma_!\left (\calE_1 \boxtimes \calE_2 \right )_s$ is the sheaf associated to the presheaf $\calP$ on $\text{ét}_{/U}$ defined by
$$
\calP:V \mapsto \bigoplus_{\substack{\alpha,\beta \in \Gamma(V,\upmu_X) \\ \alpha+\beta=s_{|V}}} \calE_{1,\alpha}(V)\otimes_{\Ox_{X}(V)}\calE_{2,\beta}(V).
$$
By \ref{ommek1}, the fiber $\left (\calQ_{1,\upmu_X} \otimes_{\Ox_{X,\upmu_X}} \sigma_!\left (\calE_1\boxtimes \calE_2 \right ) \right )_s=\calQ_{1|U} \otimes_{\Ox_U} \sigma_!\left (\calE_1 \boxtimes \calE_2\right )_s$ is then the sheaf associated to the presheaf on $\text{ét}_{/U}$ defined by
$$
V \mapsto \calQ_1(V) \otimes_{\Ox_X(V)}\bigoplus_{\substack{\alpha,\beta \in \Gamma(V,\upmu_X) \\ \alpha+\beta=s_{|V}}} \calE_{1,\alpha}(V)\otimes_{\Ox_{X}(V)}\calE_{2,\beta}(V).
$$
The morphism $\left (\sigma_!\epsilon \right )_s$ is then associated to the morphism of presheaves
$$
\varphi_s:\calQ_{1|U} \otimes_{\Ox_U}\calP \ra \calP\otimes_{\Ox_U} \calQ_{1|U}
$$
defined, for every étale $U$-scheme $V,$ by
$$
\begin{tikzcd}
\calQ_1(V) \otimes_{\Ox_X(V)} \left (\displaystyle\bigoplus_{\substack{\alpha,\beta \in \Gamma(V,\upmu_X) \\ \alpha+\beta=s_{|V}}} \calE_{1,\alpha}(V)\otimes_{\Ox_{X}(V)}\calE_{2,\beta}(V) \right ) \ar{d} & 1 \otimes x \otimes y \ar[mapsto]{d} \\
\left (\displaystyle\bigoplus_{\substack{\alpha,\beta \in \Gamma(V,\upmu_X) \\ \alpha+\beta=s_{|V}}} \calE_{1,\alpha}(V)\otimes_{\Ox_{X}(V)}\calE_{2,\beta}(V) \right ) \otimes_{\Ox_X(V)} \calQ_1(V) & \epsilon_{1,\alpha}(1\otimes x)\otimes \epsilon_{2,\beta}(1\otimes y),
\end{tikzcd}
$$
where $x$ and $y$ are local sections of $\calE_{1,\alpha}$ and $\calE_{2,\beta}$ over $V$ and
$$
\epsilon_{1,\alpha}(1\otimes x)\otimes \epsilon_{2,\beta}(1\otimes y)
$$
is a local section of
$$
\left (\calE_{1,\alpha}\otimes_{\Ox_U} \calQ_{1|U}\right ) \otimes_{\calQ_{1|U}} \left (\calE_{2,\beta} \otimes_{\Ox_U} \calQ_{1|U} \right ) \xrightarrow{\sim} \left (\calE_{1,\alpha}\otimes_{\Ox_U} \calE_{2,\beta} \right )\otimes_{\Ox_U} \calQ_{1|U}.
$$
For every $\alpha,\beta,\gamma \in \Gamma(V,\upmu_X)$ such that $\alpha+\beta+\gamma=s_{|V}$ and $x\in \calE_{1,\alpha}(V),$ $y\in \calE_{2,\beta}(V),$ $a\in \frakA_{\gamma}(V),$ we have
\begin{alignat*}{2}
\varphi_s(1\otimes x\otimes (by)) &= \epsilon_{1,\alpha}(1\otimes x)\otimes \epsilon_{2,\beta+\gamma}(1\otimes ay) \\
&=\epsilon_{1,\alpha}(1\otimes x)\otimes a\epsilon_{2,\beta}(1\otimes y), \\
\varphi_s(1\otimes (bx)\otimes y) &= \epsilon_{1,\alpha+\gamma}(1\otimes ax)\otimes \epsilon_{2,\beta}(1\otimes y) \\
&= \left (a\epsilon_{1,\alpha}(1\otimes x)\right )\otimes \epsilon_{2,\beta}(1\otimes y).
\end{alignat*}
We deduce that, when composing $\varphi_s$ with the canonical morphism
$$\calP\otimes_{\Ox_U} \calQ_{1|U} \ra \left (\calE_1 \circledast_{\frakA}\calE_2 \right )_s \otimes \calQ_{1|U},$$
the images of $\varphi_s(1\otimes x\otimes (by))$ and $\varphi_s(1\otimes (bx)\otimes y)$ are equal.
It follows, by \ref{indpropshriek} (3), that $\sigma_!\epsilon$ induces a $\calQ_{1,\upmu_X}$-linear isomorphism
\begin{equation}
\epsilon_1 \circledast_{\calQ_{1,\upmu_X}\otimes_{\Ox_{X,\upmu_X}}\frakA}\epsilon_2:
\calQ_{1,\upmu_X} \otimes_{\Ox_{X,\upmu_X}}\left (\calE_1\circledast_{\frakA}\calE_2 \right ) \ra \left (\calE_1\circledast_{\frakA}\calE_2 \right )\otimes_{\Ox_{X,\upmu_X}}\calQ_{1,\upmu_X}.
\end{equation}
The fact that it is a stratification follows from the fact that $\sigma_!\epsilon$ is.
\end{proof}

\begin{proposition}\label{propdraft4}
Assume $(H_1)$ is satisfied. Let $\frakB \ra \frakA$ a morphism of $\upmu_X$-indexed algebras and $(\calE,\epsilon)$ an $\upmu_X$-indexed $\frakB$-module equipped with $\frakB$-linear $\calQ_{1,\upmu_X}$-stratification. Suppose that $\frakA$ is equipped with a $\frakB$-linear $\calQ_{1,\upmu_X}$-stratification $\epsilon_{\frakA}$ satisfying
$$
\epsilon_{\frakA}(1\otimes (ab))=\epsilon_{\frakA}(1\otimes a)\epsilon_{\frakA}(1\otimes b)
$$
for any local sections $a$ and $b$ of $\frakA.$ Then the stratification
$$
\epsilon_1=\epsilon_{\frakA} \circledast_{\frakB \otimes_{\Ox_{X,\upmu_X}} \calQ_{1,\upmu_X}} \epsilon: \calQ_{1,\upmu_X} \otimes_{\Ox_{X,\upmu_X}} \left (\frakA \circledast_{\frakB}\calE \right ) \xrightarrow{\sim} \left (\frakA \circledast_{\frakB}\calE \right )  \otimes_{\Ox_{X,\upmu_X}} \calQ_{1,\upmu_X},
$$
defined in \ref{propdraft3}, is $\frakA$-admissible.
\end{proposition}

\begin{proof}
Let $\epsilon_2=p_1^*\epsilon_{\frakA} \otimes_{\calQ_{1,\upmu_X^2}} p_2^*\epsilon$ the stratification on $\calE_2=\frakA \boxtimes \calE$ as in \ref{propdraft1}. Since $\epsilon_{\frakA}$ and $\epsilon$ are $\frakB$-linear, then the stratification $\sigma_!\epsilon_2$ \eqref{cordraft1} induces the $\calQ_{1,\upmu_X}$-stratification $\epsilon_1$ on $\frakA\circledast_{\frakB}\calE.$
We have to prove the commutativity of the diagram \eqref{admindQstratdiag}
\begin{equation*}
\begin{tikzcd}
\calQ_{1,\upmu_X^2}\otimes_{\Ox_{X,\upmu_X^2}}(\frakA\boxtimes \left (\frakA \circledast_{\frakB}\calE \right )) \ar{d}{\epsilon_{\frakA}\boxtimes \op{Id}_{\left (\frakA \circledast_{\frakB}\calE \right )}} \ar{rr}{\op{Id}_{\calQ_{1,\upmu_X^2}}\otimes \pi} & & \calQ_{1,\upmu_X^2}\otimes_{\Ox_{X,\upmu_X^2}} \sigma^*\left (\frakA \circledast_{\frakB}\calE \right ) \ar{dd}{\sigma^*\epsilon_1} \\
\frakA \boxtimes (\calQ_{1,\upmu_X}\otimes_{\Ox_{X,\upmu_X}} \left (\frakA \circledast_{\frakB}\calE \right )) \ar{d}{\op{Id}_{\frakA}\boxtimes \epsilon_1} & & \\
(\frakA\boxtimes \left (\frakA \circledast_{\frakB}\calE \right ))\otimes_{\Ox_{X,\upmu_X^2}}\calQ_{1,\upmu_X^2}\ar{rr}{\pi\otimes \op{Id}_{\calQ_{1,\upmu_X^2}}} & & \sigma^*\left (\frakA \circledast_{\frakB}\calE \right )\otimes_{\Ox_{X,\upmu_X^2}}\calQ_{1,\upmu_X^2},
\end{tikzcd}
\end{equation*}
where $\pi$ is the morphism defining the $\upmu_X$-indexed $\frakA$-algebra structure on $\frakA \circledast_{\frakB} \calE.$
Since $\sigma^*$ is right exact, it is sufficient to prove the commutativity of the diagram
$$
\begin{tikzcd}
\calQ_{1,\upmu_X^2} \otimes_{\Ox_{X,\upmu_X^2}} \left (\frakA \boxtimes\sigma_!\left (\frakA\boxtimes \calE \right )\right ) \ar[swap]{d}{\epsilon_{\frakA}\boxtimes \op{Id}} \ar{r} & \calQ_{1,\upmu_X^2} \otimes_{\Ox_{X,\upmu_X^2}}  \sigma^*\sigma_{!}\left (\frakA \boxtimes \calE\right ) \ar{dd}{\sigma^*\sigma_{!}\epsilon_{2}} 
\\
 \frakA \boxtimes \left (\calQ_{1,\upmu_X} \otimes_{\Ox_{X,\upmu_X}} \sigma_!\left (\frakA\boxtimes\calE\right )\right ) \ar[swap]{d}{\op{Id}\boxtimes \sigma_!\epsilon_{2}} & 
\\
 \left (\frakA \boxtimes \left (\sigma_!\left (\frakA\boxtimes\calE\right )\right )\right )  \otimes_{\Ox_{X,\upmu_X^2}}  \calQ_{1,\upmu_X^2} \ar{r} &  \sigma^*\sigma_{!}\left (\frakA \boxtimes \calE\right ) \otimes_{\Ox_{X,\upmu_X^2}} \calQ_{1,\upmu_X^2},
\end{tikzcd}
$$
where the upper and lower arrows are induced by the $\upmu_X$-indexed $\frakA$-algebra structure on $\sigma_!\left (\frakA\boxtimes \calE\right )$ \eqref{krazzzz}.
For that, it is sufficient to prove the commutativity of the fiber of the previous diagram over all local sections of $\upmu_X^2.$ Let $s,t\in \Gamma(U,\upmu_X)$ be local sections over an étale $X$-scheme $U.$ It is then sufficient to prove the commutativity of
\begin{equation}\label{diagdraft2}
\begin{tikzcd}
\calQ_{1|U} \otimes \frakA_s \otimes \left (\sigma_!\left (\frakA\boxtimes \calE \right )\right )_t \ar[swap]{d}{\epsilon_{\frakA,s}\otimes \op{Id}} \ar{r} & \calQ_{1|U} \otimes \left (\sigma_{!}\left (\frakA \boxtimes \calE\right )\right )_{s+t} \ar{dd}{\left (\sigma_{!}\epsilon_{2}\right )_{s+t}} 
\\
 \frakA_s \otimes \calQ_{1|U} \otimes \left (\sigma_!\left (\frakA\boxtimes\calE\right )\right )_t \ar[swap]{d}{\op{Id}\otimes \left (\sigma_!\epsilon_{2}\right )_t} & 
\\
 \frakA_s \otimes \left (\sigma_!\left (\frakA\boxtimes\calE\right )\right )_t  \otimes \calQ_{1|U} \ar{r} &  \left (\sigma_{!}\left (\frakA \boxtimes \calE\right )\right )_{s+t} \otimes \calQ_{1|U},
\end{tikzcd}
\end{equation}
where the tensor products are over $\Ox_U.$
By \ref{indpropshriek} (3), the fiber $\sigma_!\left (\frakA \boxtimes \calE \right )_s$ is the sheaf associated to the presheaf $\calP_s$ on $\text{ét}_{/U}$ defined by
$$
\calP_s:V \mapsto \bigoplus_{\substack{\alpha,\beta \in \Gamma(V,\upmu_X) \\ \alpha+\beta=s_{|V}}} \frakA_{\alpha}(V)\otimes_{\Ox_{X}(V)}\calE_{\beta}(V).
$$
If $V$ is an étale $U$-scheme, $\alpha,\beta \in \Gamma(V,\upmu_X)$ and $b$ and $y$ are local sections of $\frakA_{\alpha}$ and $\calE_{\beta}$ respectively, then
\begin{alignat*}{2}
\epsilon_{2,(\alpha,\beta)}(1\otimes b\otimes y) = \epsilon_{\frakA,\alpha}(1\otimes b) \otimes \epsilon_{\beta}(1\otimes y) &\in \left ( \frakA_{\alpha} \otimes_{\Ox_V} \calQ_{1|V} \right ) \otimes_{\calQ_{1|V}} \left ( \calE_{\beta} \otimes_{\Ox_V} \calQ_{1|V} \right )\\ &\xrightarrow{\sim} \left (\frakA_{\alpha} \otimes_{\Ox_V} \calE_{\beta} \right )\otimes_{\Ox_V}\calQ_{1|V}.
\end{alignat*}
Then
$$
\left (\sigma_!\epsilon_2 \right )_{\alpha+\beta}(1\otimes b \otimes y)=\epsilon_{\frakA,\alpha}(1\otimes b) \otimes \epsilon_{\beta}(1\otimes y).
$$
Let $\ov{x} \ra U$ be a geometric point. We have then a canonical isomorphism of stalks
$$
\left (\sigma_!\left (\frakA \boxtimes \calE \right )\right )_{t,\ov{x}} \xrightarrow{\sim} \calP_{t,\ov{x}}.
$$
For simplicity, when tensoring over $\Ox_{X,\ov{x}},$ we use the notation $\otimes$ instead of $\otimes_{\Ox_{X,\ov{x}}}.$
To prove the commutativity of \eqref{diagdraft2}, it is sufficient to prove the commutativity of the diagram
\begin{equation}\label{diagdraft1}
\begin{tikzcd}
\calQ_{1,\ov{x}} \otimes \frakA_{s,\ov{x}} \otimes \calP_{t,\ov{x}} \ar[swap]{d}{\epsilon_{\frakA,s,\ov{x}}\otimes \op{Id}} \ar{r} & \calQ_{1,\ov{x}} \otimes \calP_{s+t,\ov{x}} \ar{dd}{\left (\sigma_{!}\epsilon_{2}\right )_{s+t,\ov{x}}} 
\\
 \frakA_{s,\ov{x}} \otimes \calQ_{1,\ov{x}} \otimes \calP_{t,\ov{x}} \ar[swap]{d}{\op{Id}\otimes \left (\sigma_!\epsilon_{2}\right )_{t,\ov{x}}} & 
\\
 \frakA_{s,\ov{x}} \otimes \calP_{t,\ov{x}}  \otimes \calQ_{1,\ov{x}} \ar{r} &  \calP_{s+t,\ov{x}} \otimes \calQ_{1,\ov{x}}.
\end{tikzcd}
\end{equation}
If $V$ is an étale $U$-scheme, $\alpha,\beta \in \Gamma(V,\upmu_X)$ such that $\alpha+\beta=t_{|V},$ and if $a,$ $b$ and $y$ are local sections of $\frakA_{s},$ $\frakA_{\alpha}$ and $\calE_{\beta}$ respectively, then the images of the element $1\otimes a \otimes b\otimes y$ of $\calQ_{1,\ov{x}} \otimes \frakA_{s,\ov{x}} \otimes \calP_{t,\ov{x}}$ by the arrows of the previous diagram are as follows:
$$
\begin{tikzcd}
1\otimes a \otimes b \otimes y \ar[mapsto]{r} \ar[mapsto]{d} & 1\otimes (ab) \otimes y \ar[mapsto]{dd} \\
\epsilon_{\frakA,s,\ov{x}}(1\otimes a) \otimes b \otimes y \ar[mapsto]{d} & \\
\begin{matrix}\epsilon_{\frakA,s,\ov{x}} (1\otimes a) \otimes \epsilon_{2,(\alpha,\beta),\ov{x}}(1\otimes b \otimes y) \\
= \epsilon_{\frakA,s,\ov{x}}(1\otimes a) \otimes \epsilon_{\frakA,\alpha,\ov{x}}(1\otimes b) \otimes \epsilon_{\beta,\ov{x}}(1\otimes y). \end{matrix} \ar[mapsto]{r} & \begin{matrix} \left (\epsilon_{\frakA,s,\ov{x}}(1\otimes a) \cdot \epsilon_{\frakA,\alpha,\ov{x}}(1\otimes b) \right ) \otimes \epsilon_{\beta,\ov{x}}(1\otimes y) \\
= \epsilon_{\frakA,s+\alpha,\ov{x}}(1\otimes ab) \otimes \epsilon_{\beta,\ov{x}}(1\otimes y), \end{matrix}
\end{tikzcd}
$$
where
\begin{alignat*}{2}
\epsilon_{\frakA,s,\ov{x}} (1\otimes a) \otimes \epsilon_{2,(\alpha,\beta),\ov{x}}(1\otimes b \otimes y) &\in \left (\frakA_{s,\ov{x}} \otimes \calQ_{1,\ov{x}} \right ) \otimes_{\calQ_{1,\ov{x}}} \left (\frakA_{\alpha,\ov{x}} \otimes \calE_{\beta,\ov{x}} \otimes \calQ_{1,\ov{x}} \right ) \\
& \xrightarrow{\sim} \left (\frakA_{s,\ov{x}}\otimes \frakA_{\alpha,\ov{x}} \otimes \calE_{\beta,\ov{x}} \right ) \otimes \calQ_{1,\ov{x}}.
\end{alignat*}
This proves the commutativity of \eqref{diagdraft1} and hence the $\frakA$-admissibility of $\epsilon.$
\end{proof}

\begin{parag}\label{ffLcrys}
Assume $(H_1)$ is satisfied. Let $\calE'$ be a $\upmu_X$-indexed $F_1^*\calB_{X/S}$-module equipped with a $\calB_{X/S}$-linear $\calQ_{1,\upmu_X}$-stratification
$$
\epsilon':\calQ_{1,\upmu_X}\otimes_{\Ox_{X,\upmu_X}}\calE' \ra \calE' \otimes_{\Ox_{X,\upmu_X}} \calQ_{1,\upmu_X}.
$$
Recall the stratification \eqref{eqeAQ}
$$
\epsilon_{\calA_X,Q}:\calQ_{1,\upmu_X}\otimes_{\Ox_{X,\upmu_X}}\calA_X \ra \calA_X \otimes_{\Ox_{X,\upmu_X}} \calQ_{1,\upmu_X}.
$$
By \ref{propdraft3} we have a stratification
$$
\epsilon=\epsilon_{\calA_X,Q} \circledast_{\calQ_{1,\upmu_X}\otimes_{\Ox_{X,\upmu_X}}F_1^*\calB_{X/S}}\epsilon'.
$$
Its source is
$$
\left (\calQ_{1,\upmu_X} \otimes_{\Ox_{X,\upmu_X}} \calA_X \right ) \circledast_{\calQ_{1,\upmu_X}\otimes_{\Ox_{X,\upmu_X}}F_1^*\calB_{X/S}} \left (\calQ_{1,\upmu_X}\otimes_{\Ox_{X,\upmu_X}} \calE' \right ).
$$
By \ref{kraz827} and \eqref{indiso1}, this source is canonically isomorphic to
$$
\calQ_{1,\upmu_X} \otimes_{\Ox_{X,\upmu_X}} \left (\calA_X \circledast_{F_1^*\calB_{X/S}} \calE'\right ) = \left (\calQ_{1,\upmu_X} \otimes_{\Ox_{X,\upmu_X}} \calA_X \right )\circledast_{F_1^*\calB_{X/S}} \calE'.
$$
Similarly, the target is isomorphic to
$$
\calA_X \circledast_{F_1^*\calB_{X/S}}\calE'\otimes_{\Ox_{X,\upmu_X}}\calQ_{1,\upmu_X},
$$
considered as a $\upmu_X$-indexed module over $\calQ_{1,\upmu_X}\otimes_{\Ox_{X,\upmu_X}}F_1^*\calB_{X/S}$ via the morphism 
$$
\sigma_Q:\calQ_{1,\upmu_X}\otimes_{\Ox_{X,\upmu_X}}F_1^*\calB_{X/S} \ra F_1^*\calB_{X/S} \otimes_{\Ox_{X,\upmu_X}}\calQ_{1,\upmu_X}$$
given in \eqref{eqsigmaQkraz}.
By \ref{propdraft4}, $\epsilon$ is an admissible $\calQ_{1,\upmu_X}$-stratification.
We obtain a linearization functor
\begin{equation}\label{eqLstrat1942f}
\begin{array}[t]{clc}
L_{\mathrm{strat}}:\begin{Bmatrix}\underline{\upmu}\text{-indexed\ }F_1^*\calB_{X/S}\text{-modules\ with\ a}\\ \calB_{X/S}\text{-linear}\\Q_{1,\upmu_X}\text{-stratification} \end{Bmatrix} & \ra & \begin{Bmatrix}\underline{\upmu}\text{-indexed\ }\calA_X\text{-modules\ with}\\ \text{an\ admissible}\\Q_{1,\upmu_X}\text{-stratification} \end{Bmatrix} \\
 & & \\
(\calE',\epsilon') & \mapsto & \left (\calA_X\circledast_{F_1^*\calB_{X/S}}\calE',\epsilon \right ).
\end{array}
\end{equation}
\end{parag}

\begin{theorem}\label{thmKoko20125}
Assume $(H_2)$ is satisfied. Keep the notation of \ref{parag1419}. Consider the diagram
\begin{equation}\label{diag187}
\begin{tikzcd}
\begin{Bmatrix}\mathrm{Crystals\ of\ }\upmu'\text{-}\mathrm{indexed}\\ \frakB_X\text{-}\mathrm{modules\ of\ }\widetilde{\calE}(X'/\frakS)_{/\upmu'} \end{Bmatrix} \ar{r}{C^{-1}} \ar[swap,sloped]{d}{\sim} & \begin{Bmatrix}\mathrm{Crystals\ of\ }\underline{\upmu}\text{-}\mathrm{indexed}\\ \frakA_X\text{-}\mathrm{modules\ of\ }\widetilde{\underline{\calE}}(X/\frakS)_{/\underline{\upmu}} \end{Bmatrix} \ar[sloped]{d}{\sim}\\
\begin{Bmatrix}\upmu_X\text{-}\mathrm{indexed\ }\calB_{X/S}\text{-}\mathrm{modules\ with}\\ \mathrm{an\ admissible}\\ \calR_{1,\upmu_X}'\text{-}\mathrm{stratification} \end{Bmatrix} \ar{r}{L_{\mathrm{strat}} \circ \Psi_{\upmu}} & \begin{Bmatrix}\upmu_X\text{-}\mathrm{indexed\ }\calA_X\text{-}\mathrm{modules\ with}\\ \mathrm{an\ admissible}\\Q_{1,\upmu_X}\text{-}\mathrm{stratification} \end{Bmatrix},
\end{tikzcd}
\end{equation}
where the vertical equivalences are given in \ref{equivRQAB} and $C^{-1},$ $L_{\mathrm{strat}}$ and $\Psi_{\upmu}$ are given respectively in \eqref{ffC}, \ref{ffLcrys} and \eqref{diag186}.
The diagram \eqref{diag187} is 2-commutative.
\end{theorem}

\begin{proof}
By definition of $C^{-1}$ \eqref{ffC}, the diagram \eqref{diag187} factors as follows:
$$
\begin{tikzcd}
\begin{Bmatrix}\mathrm{Crystals\ of\ }\\ \upmu'\text{-indexed}\\ \frakB_X\text{-modules\ of\ }\\ \widetilde{\calE}(X'/\frakS)_{/\upmu'} \end{Bmatrix} \ar{r}{C_{X/\frakS,\upmu}^{-1}} \ar[bend right=-30]{rr}{C^{-1}} \ar[swap,sloped]{d}{\sim} & \begin{Bmatrix}\mathrm{Crystals\ of\ }\\ \underline{\upmu}\text{-indexed}\\ C_{X/\frakS,\upmu}^{-1}\frakB_X\text{-modules\ of\ }\\ \widetilde{\underline{\calE}}(X/\frakS)_{/\underline{\upmu}} \end{Bmatrix} \ar{r}{L_{\mathrm{crys}}} \ar[swap,sloped]{d}{\sim} \ar{d}{\omega} & \begin{Bmatrix}\mathrm{Crystals\ of\ } \\ \underline{\upmu}\text{-indexed}\\ \frakA_X\text{-modules\ of\ } \\ \widetilde{\underline{\calE}}(X/\frakS)_{/\underline{\upmu}} \end{Bmatrix} \ar[swap, sloped]{d}{\sim} \ar{d}{\omega'} \\
\begin{Bmatrix}\upmu_X\text{-indexed\ }\\ \calB_{X/S}\text{-modules\ with}\\ \mathrm{an\ admissible}\\ \calR_{1,\upmu_X}'\text{-stratification} \end{Bmatrix} \ar{r}{\Psi_{\upmu}} \ar[swap,bend right=30]{rr}{L_{\mathrm{strat}} \circ \Psi_{\upmu}} & \begin{Bmatrix}\upmu_X\text{-indexed\ } \\ F_1^*\calB_{X/S}\text{-modules\ with\ a}\\ \calB_{X/S}\text{-linear}\\ \calQ_{1,\upmu_X}\text{-stratification} \end{Bmatrix} \ar{r}{L_{\mathrm{strat}}} & \begin{Bmatrix}\upmu_X\text{-indexed\ } \\ \calA_X\text{-modules\ with}\\ \mathrm{an\ admissible}\\Q_{1,\upmu_X}\text{-stratification} \end{Bmatrix}.
\end{tikzcd}
$$
The 2-commutativity of the left square was proved in \ref{propkraz1934}. It is hence sufficient to prove the 2-commutativity of the right square. For simplicity, when we tensor over $\Ox_{X,\upmu_X},$ we use the notation $\otimes$ instead of $\otimes_{\Ox_{X,\upmu_X}}.$
Let $\frakE$ be an object of the upper middle category. Keep the notation of \ref{parag1419}. By definition, we have
$$
L_{\mathrm{crys}}(\frakE)=\frakA_X \circledast_{C_{X/\frakS,\upmu}^{-1}\frakB_X}\mathfrak{E},
$$
and
$$
\omega'\left (\frakA_X \circledast_{C_{X/\frakS,\upmu}^{-1}\frakB_X}\mathfrak{E} \right )= \left ( \calE',\epsilon'\right ),
$$
where $\calE'=\left (\frakA_X \circledast_{C_{X/\frakS,\upmu}^{-1}\frakB_X}\mathfrak{E}\right )_{(X,\frakX,i_X)}$ and $\epsilon'$ is the composition
$$
\epsilon':\calQ_{1,\upmu_X} \otimes \calE' \xrightarrow{c_{\frakA_X \circledast_{C_{X/\frakS,\upmu}^{-1}\frakB_X}\frakE,q_2}} \left (\frakA_X \circledast_{C_{X/\frakS,\upmu}^{-1}\frakB_X}\mathfrak{E}\right )_{(X,Q_{\frakX},\lambda_Q)} \xrightarrow{c_{\frakA_X \circledast_{C_{X/\frakS,\upmu}^{-1}\frakB_X}\mathfrak{E},q_1}^{-1}}\calE'  \otimes \calQ_{1,\upmu_X}.
$$
On the other hand, we have
$$
\omega(\frakE)=\left (\F,\epsilon \right ),
$$
where $\F=\frakE_{(X,\frakX,i_X)}$ and $\epsilon$ is the composition
$$
\epsilon:\calQ_{1,\upmu_X} \otimes \F \xrightarrow{c_{\mathfrak{E},q_2}} \mathfrak{E}_{(X,Q_{\frakX},\lambda_Q)} \xrightarrow{c_{\mathfrak{E},q_1}^{-1}} \F \otimes \calQ_{1,\upmu_X}.
$$
By \eqref{eqLstrat1942f}, we have
$$
L_{\mathrm{strat}}\left (\F,\epsilon \right )=\left (\calA_X\circledast_{F_1^*\calB_{X/S}}\F,\epsilon_{\calA_X,Q}\circledast_{\calQ_{1,\upmu_X}\otimes F_1^*\calB_{X/S}}\epsilon\right ),
$$
where the target of
$$
\epsilon_{\calA_X,Q}\circledast_{\calQ_{1,\upmu_X}\otimes F_1^*\calB_{X/S}}\epsilon
$$
is considered as a $\upmu_X$-indexed module over $\calQ_{1,\upmu_X}\otimes_{\Ox_{X,\upmu_X}}F_1^*\calB_{X/S}$ via the morphism 
$$
\sigma_Q:\calQ_{1,\upmu_X}\otimes F_1^*\calB_{X/S} \ra F_1^*\calB_{X/S} \otimes \calQ_{1,\upmu_X}$$
given in \eqref{eqsigmaQkraz}.
We have to prove the existence of an isomorphism
\begin{equation}\label{eqtakriz185555}
\left (\calE',\epsilon' \right ) \xrightarrow{\sim} \left (\calA_X\circledast_{F_1^*\calB_{X/S}}\F,\epsilon_{\calA_X,Q}\circledast_{\calQ_{1,\upmu_X}\otimes F_1^*\calB_{X/S}}\epsilon\right ).
\end{equation}
By \ref{propKoko1915}, there exists canonical isomorphisms
$$
\calE'=\left (\frakA_X \circledast_{C_{X/\frakS,\upmu}^{-1}\frakB_X}\mathfrak{E}\right )_{(X,\frakX,i_X)} \xrightarrow{\sim} \calA_X\circledast_{ \left ( C_{X/\frakS,\upmu}^{-1}\frakB_X\right )_{(X,\frakX,i_X)}} \F,
$$
$$
\left (\frakA_X \circledast_{C_{X/\frakS,\upmu}^{-1}\frakB_X}\mathfrak{E}\right )_{(X,Q_{\frakX},\lambda_Q)} \xrightarrow{\sim} \frakA_{X,(X,Q_{\frakX},\lambda_Q)} \circledast_{\left ( C_{X/\frakS,\upmu}^{-1}\frakB_X \right )_{(X,Q_{\frakX},\lambda_Q)}}\frakE_{(X,Q_{\frakX},\lambda_Q)}.
$$
By \ref{ommek4}, we have canonical isomorphisms
$$
\left ( C_{X/\frakS,\upmu}^{-1}\frakB_X\right )_{(X,\frakX,i_X)} \xrightarrow{\sim} \frakB_{X,\rho(X,\frakX,i_X)},
$$
$$
\left ( C_{X/\frakS,\upmu}^{-1}\frakB_X \right )_{(X,Q_{\frakX},\lambda_Q)} \xrightarrow{\sim} \frakB_{X,\rho(X,Q_{\frakX},\lambda_Q)}.
$$
The lifting $F:\frakX \ra \frakX'$ and the morphism $\nu:Q_{\frakX} \ra R_{\frakX',1}$ induce morphisms \eqref{zabb222}, abusively denoted $F$ and $\nu,$
$$
F:\rho(X,\frakX,i_X) = (X',\frakX,F_1) \ra (X',\frakX',\op{Id}_{X'}),
$$
$$
\nu:\rho(X,Q_{\frakX},\lambda_Q) \ra (X',R_{\frakX',1},\lambda_{R'}),
$$
such that the diagram
$$
\begin{tikzcd}
\rho(X,Q_{\frakX},\lambda_Q) \ar{r}{\nu} \ar[swap]{d}{\rho(q_i)} & (X',R_{\frakX',1},\lambda_{R'}) \ar{d}{r_i} \\
\rho(X,\frakX,i_X) \ar[swap]{r}{F} & (X',\frakX',\op{Id}_{X'})
\end{tikzcd}
$$
is commutative, where $r_1$ and $r_2$ are the canonical projections.
Since $\frakB_X$ is a crystal, we have isomorphisms
$$
c_{\frakB_X,F}:F_1^*\calB_{X/S}=F_1^*\frakB_{X,(X',\frakX',\op{Id}_{X'})} \xrightarrow{\sim} \frakB_{X,\rho(X,\frakX,i_X)},
$$
$$
c_{\frakB_X,\nu}:\calQ_{1,\upmu_X} \otimes_{\calR'_{1,\upmu_X}} \frakB_{X,(X',R_{\frakX',1},\lambda_{R'})} \xrightarrow{\sim} \frakB_{X,\rho(X,Q_{\frakX},\lambda_Q)},
$$
$$
c_{\frakB_X,r_2}:\calR_{1,\upmu_X}' \otimes_{\Ox_{X',\upmu_X}} \calB_{X/S} \xrightarrow{\sim} \frakB_{X,(X',R_{\frakX',1},\lambda_{R'})},
$$
$$
c_{\frakB_X,r_1}:  \calB_{X/S}\otimes_{\Ox_{X',\upmu_X}} \calR_{1,\upmu_X}'\xrightarrow{\sim} \frakB_{X,(X',R_{\frakX',1},\lambda_{R'})}.
$$
We deduce canonical isomorphisms
$$
\calE'=\left (\frakA_X \circledast_{C_{X/\frakS,\upmu}^{-1}\frakB_X}\mathfrak{E}\right )_{(X,\frakX,i_X)} \xrightarrow{\sim} \calA_X\circledast_{F_1^*\calB_{X/S}} \F,
$$
\begin{alignat*}{2}
\left (\frakA_X \circledast_{C_{X/\frakS,\upmu}^{-1}\frakB_X}\mathfrak{E}\right )_{(X,Q_{\frakX},\lambda_Q)} &\xrightarrow{\sim} \frakA_{X,(X,Q_{\frakX},\lambda_Q)} \circledast_{\calQ_{1,\upmu_X} \otimes_{\calR'_{1,\upmu_X}} \frakB_{X,(X',R_{\frakX',1},\lambda_{R'})}}\frakE_{(X,Q_{\frakX},\lambda_Q)} \\
&\xrightarrow{\sim}  \frakA_{X,(X,Q_{\frakX},\lambda_Q)} \circledast_{\calQ_{1,\upmu_X} \otimes_{\Ox_{X',\upmu_X}} \calB_{X/S}}\frakE_{(X,Q_{\frakX},\lambda_Q)},
\end{alignat*}
fitting into the commutative diagram
$$
\begin{tikzcd}
\calQ_{1,\upmu_X} \otimes \calE' \ar{r}{\sim} \ar[swap]{dd}{c_{\frakA_X \circledast_{C_{X/\frakS,\upmu}^{-1}\frakB_X}\frakE,q_2}} &
\calQ_{1,\upmu_X} \otimes \left (\calA_X\circledast_{F_1^*\calB_{X/S}} \F\right ) \ar[sloped]{d}{\sim} \\
 &
\left (\calQ_{1,\upmu_X} \otimes \calA_X \right )\circledast_{\calQ_{1,\upmu_X} \otimes F_1^*\calB_{X/S}}  \left ( \calQ_{1,\upmu_X} \otimes \F\right ) \ar{d}{c_{\frakA_X,q_2}\circledast c_{\frakE,q_2}} \\
\left (\frakA_X \circledast_{C_{X/\frakS,\upmu}^{-1}\frakB_X}\mathfrak{E}\right )_{(X,Q_{\frakX},\lambda_Q)} \ar{r}{\sim} \ar[swap]{dd}{c_{\frakA_X \circledast_{C_{X/\frakS,\upmu}^{-1}\frakB_X}\mathfrak{E},q_1}^{-1}}  &
 \frakA_{X,(X,Q_{\frakX},\lambda_Q)} \circledast_{\calQ_{1,\upmu_X} \otimes F_1^*\calB_{X/S}}\frakE_{(X,Q_{\frakX},\lambda_Q)}  \ar{d}{c_{\frakA_X,q_1}^{-1}\circledast c_{\frakE,q_1}^{-1}} \\
&
\left ( \calA_X \otimes \calQ_{1,\upmu_X} \right )\circledast_{ F_1^*\calB_{X/S} \otimes \calQ_{1,\upmu_X} }  \left (  \F \otimes \calQ_{1,\upmu_X} \right ) \ar[sloped]{d}{\sim} \\
\calE' \otimes \calQ_{1,\upmu_X} \ar{r}{\sim} & \left (\calA_X\circledast_{F_1^*\calB_{X/S}} \F \right )\otimes \calQ_{1,\upmu_X}.
\end{tikzcd}
$$
We deduce the desired isomorphism of stratified $\upmu_X$-indexed $\calA_X$-modules \eqref{eqtakriz185555}.
\end{proof}

\begin{parag}\label{kraz1935}
In this paragraph, we assume $(H_1)$ is satisfied and we define a $\upmu_X$-indexed algebra structure on $\calQ_{1,\upmu_X}^{\vee} \otimes_{\Ox_{X,\upmu_X}} \calA_X.$
Since $\calA_X$ is an invertible $\Ox_{X,\upmu_X}$-module, we have a canonical isomorphism
\begin{equation}\label{isoQHom}
\calQ_{1,\upmu_X}^{\vee} \otimes_{\Ox_{X,\upmu_X}} \calA_X \xrightarrow{\sim} \mathscr{Hom}_{\Ox_{X,\upmu_X}}\left ( \calQ_{1,\upmu_X},\calA_X \right ).
\end{equation}
Note that, by (\cite{SGA43} IV 12.3), we have a canonical isomorphism
$$
\calQ_{1,\upmu_X}^{\vee} \xrightarrow{\sim} \left (\calQ_1^{\vee} \right )_{\upmu_X}.
$$
We define a multiplication on $\mathscr{Hom}_{\Ox_{X,\upmu_X}}\left ( \calQ_{1,\upmu_X},\calA_X \right )$ as follows: for local sections $s,t\in \Gamma(U,\upmu_X)$ over an étale $X$-scheme $U,$ and global sections $f$ and $g$ of
$$\mathscr{Hom}_{\Ox_{X,\upmu_X}}\left ( \calQ_{1,\upmu_X},\calA_X \right )_s\xrightarrow{\sim} \mathscr{Hom}_{\Ox_{U}}\left ( \calQ_{1|U},\calA_{X,s} \right )$$
and
$$\mathscr{Hom}_{\Ox_{X,\upmu_X}}\left ( \calQ_{1,\upmu_X},\calA_X \right )_t\xrightarrow{\sim}\mathscr{Hom}_{\Ox_{U}}\left ( \calQ_{1|U},\calA_{X,t} \right )$$
respectively (where the isomorphisms follow from (\cite{SGA43} IV 12.3)). We abusively denote $\calQ_{1|U}$ simply by $\calQ_1.$ Define $f\circ g$ to be the composition
$$
\calQ_1 \xrightarrow{\delta} \calQ_1 \otimes_{\Ox_U} \calQ_1 \xrightarrow{\op{Id}\otimes g} \calQ_1 \otimes_{\Ox_U} \calA_{X,t} \xrightarrow{\epsilon_{\calA_X,Q,t}} \calA_{X,t}\otimes_{\Ox_U} \calQ_1  \xrightarrow{\op{Id}\otimes f} \calA_{X,t} \otimes_{\Ox_U} \calA_{X,s} \ra \calA_{X,s+t}.
$$
The only difficulty we encounter when proving that this multiplication defines an indexed algebra structure on $\mathscr{Hom}_{\Ox_{X,\upmu_X}}\left ( \calQ_{1,\upmu_X},\calA_X \right )$ and hence on $\calQ_{1,\upmu_X}^{\vee} \otimes_{\Ox_{X,\upmu_X}} \calA_X,$ is the associativity. This is proved exactly the same way as in \ref{era2stralgD}.

Suppose that the hypothesis \ref{loccoord} is satisfied. We have a basis $\left (\eta^K\eta_{(q)}^L \right )_{\substack{K\in \llbracket 0,p-1 \rrbracket^d\\ L\in \N^d}}$ for the $\Ox_X$-module $\calQ_1.$ We denote by $\left (\varphi_{K,L} \right )_{\substack{K\in \llbracket 0,p-1 \rrbracket^d \\ L\in \N^d}}$ its dual basis. Let $a$ be a local section of $\calA_{X,s}$ and $(I,J) \in \llbracket 0,p-1 \rrbracket^d \times \N^d.$ Then $\varphi_{I,J}\otimes 1$ is a section of $\calQ_{1}^{\vee} \otimes_{\Ox_{U}} \calA_{X,0}.$ Let us compute the product
$$
(\varphi_{I,J} \otimes 1)\cdot (1 \otimes a).
$$
By definition, $1\otimes a$ corresponds, by \eqref{isoQHom}, to the morphism
$$
\theta_a:\begin{array}[t]{clc}
\calQ_1 & \ra & \calA_{X,s} \\
\eta^{\alpha} \eta_{(q)}^{\beta} & \mapsto & \begin{cases} a\ \text{if }\alpha=\beta=0,\\ 0\ \text{else.} \end{cases} 
\end{array}
$$
By \eqref{isoQHom}, the product $(\varphi_{I,J} \otimes 1)\cdot (1 \otimes a)$ corresponds then to the composition
$$
\calQ_1 \xrightarrow{\delta} \calQ_1 \otimes_{\Ox_U} \calQ_1 \xrightarrow{\op{Id}\otimes \theta_a} \calQ_1 \otimes_{\Ox_U} \calA_{X,s} \xrightarrow{\epsilon_{\calA_X,Q,s}} \calA_{X,s}\otimes_{\Ox_U} \calQ_1  \xrightarrow{\op{Id}\otimes \varphi_{I,J}} \calA_{X,s}.
$$
By \ref{HopfQ}, $\eta^{\alpha}\eta_{(q)}^{\beta}$ is sent, by $\delta,$ to
\begin{alignat*}{2}
&\prod_{i=1}^d\left (1\otimes \eta_i+\eta_i\otimes 1 + \eta_i\otimes \eta_i \right )^{\alpha_i} \\
\times & \prod_{i=1}^d \left (1\otimes \eta_{i(q)}+\sum_{0<b+c<p} \frac{(-1)^{b+c}}{b+c}\begin{pmatrix}b+c \\ b \end{pmatrix} \eta_i^{b+c}\otimes \eta_i^{p-b} +\eta_{i(q)}\otimes 1 \right )^{\beta_i}.
\end{alignat*}
By developing the right hand side, this is sent by $\op{Id}\otimes \theta_a$ to
$$
\eta^{\alpha}\eta_{(q)}^{\beta} \otimes a.
$$
By \eqref{eqepsilonAQ}, this is sent, by $\epsilon_{A_X,Q,s},$ to
$$
\sum_{K\in \llbracket 0,p-1 \rrbracket^d} \left (\varphi_{K,0}\cdot a\right )\otimes \eta^{K+\alpha}\eta_{(q)}^{\beta}.
$$
This is then sent, by $\op{Id}\otimes \varphi_{I,J},$ to
$$
\left (\varphi_{I-\alpha,0}\cdot a\right )\delta_{\beta,J}.
$$
Recall that, for $1\le i\le d,$ we have $\eta_i^p=p\eta_{i(q)}=0.$
We conclude that
\begin{equation}\label{eq20311}
(\varphi_{I,J} \otimes 1)\cdot (1 \otimes a)=\sum_{K+L=I}\left (\varphi_{K,0}\cdot a \right )\varphi_{L,J}.
\end{equation}
\end{parag}

\begin{proposition}\label{prop1932Ko}
Assume $(H_1)$ is satisfied. Let $(\calE,\epsilon)$ be a $\upmu_X$-indexed $\calA_X$-module equipped with a $\calQ_{1,\upmu_X}$-stratification and consider the corresponding $\calQ_{1,\upmu_X}^{\vee}$-module structure on $\calE.$ This structure extends to a structure of $\left (\calQ_{1,\upmu_X}^{\vee} \otimes_{\Ox_{X,\upmu_X}} \calA_X\right )$-module if and only if $\epsilon$ is admissible.
\end{proposition}

\begin{proof}
Denote by $\sigma:\upmu_X^2 \ra \upmu_X$ the addition map and by $\pi:\calA_X\boxtimes \calE \ra \sigma^*\calE$ the morphism defining the $\upmu_X$-indexed $\calA_X$-module structure on $\calE.$ Keep the notation of \ref{kraz1935}. Let $s,t\in \Gamma(U,\upmu_X)$ be local sections over an étale $X$-scheme $U$ and $a$ and $x$ local sections of $\calA_{X,s}$ and $\calE_t$ respectively. We have
$$
\epsilon_{s+t}(1\otimes (ax)) = \sum_{\substack{I \in \llbracket 0,p-1 \rrbracket^d \\ J\in \N^d}} \left (\varphi_{I,J}\cdot (ax) \right )\otimes \eta^I\eta_{(q)}^J.
$$
By \eqref{eqepsilonAQ}, we have
$$
\epsilon_{\calA_X,Q,s}(1\otimes a)=\sum_{K\in \llbracket 0,p-1 \rrbracket^d} \left (\varphi_{K,0}\cdot a\right )\otimes \eta^K.
$$
The image of $1\otimes a\otimes x$ by the composition
$$
\calQ_1 \otimes_{\Ox_U} \calA_{X,s} \otimes_{\Ox_U} \calE_t \xrightarrow{\epsilon_{\calA_X,Q,s}\otimes \op{Id}} \calA_{X,s}\otimes_{\Ox_U}\calQ_1\otimes_{\Ox_U}\calE_t \xrightarrow{\op{Id}\otimes \epsilon_t} \calA_{X,s}\otimes_{\Ox_U}\calE_t\otimes_{\Ox_U}\calQ_1
$$
is then equal to
$$
\sum_{\substack{I,K \in \llbracket 0,p-1 \rrbracket^d \\ J\in \N^d}} \left (\varphi_{K,0}\cdot a \right )\otimes \left (\varphi_{I,J}\cdot x\right ) \otimes \eta^{I+K}\eta_{(q)}^J.
$$
This is sent, by the morphism $\pi_{(s,t)}\otimes \op{Id}_{\calQ_1}:\calA_{X,s}\otimes_{\Ox_U}\calE_t\otimes_{\Ox_U}\calQ_1 \ra \calE_{s+t}\otimes_{\Ox_U}\calQ_1,$ to
$$
\sum_{\substack{I,K \in \llbracket 0,p-1 \rrbracket^d \\ J\in \N^d}} \left (\varphi_{K,0}\cdot a \right ) \left (\varphi_{I,J}\cdot x\right ) \otimes \eta^{I+K}\eta_{(q)}^J.
$$
Recall that, for $1\le i\le d,$ we have $\eta_i^p=p\eta_{i(q)}=0.$
It follows that $\epsilon$ is admissible if and only if, for all $I\in\llbracket 0,p-1 \rrbracket^d$ and $J\in \N^d,$
$$
\varphi_{I,J}\cdot (ax)=\sum_{K+L=I}\left (\varphi_{K,0}\cdot a \right ) \left (\varphi_{L,J}\cdot x\right ).
$$
We conclude by \eqref{eq20311}.
\end{proof}

\begin{definition}
Let
$$
S=\left (S^{\bullet}\calT_{X'/S} \right )_{\upmu_X},\ \widehat{\Gamma}=\left (\widehat{\Gamma}^{\bullet}\calT_{X'/S} \right )_{\upmu_X}.
$$
Recall the $\upmu_X$-indexed algebra $\widetilde{\calD}_{X/S}=\calA_X \otimes_{\Ox_{X,\upmu_X}} \calD_{X/S,\upmu_X}$ \eqref{Dtilda1713} and the $\Ox_X$-algebra $\calD_{X/S}^{\gamma}=\calD_{X/S} \otimes_{S^{\bullet}\calT_{X'/S}} \w{\Gamma}^{\bullet} \calT_{X'/S}$ \eqref{eq12331bis}.
We denote
$$
\widetilde{\calD}_{X/S}^{\gamma}=\widetilde{\calD}_{X/S}\otimes_S\widehat{\Gamma}=\calA_X\otimes_{\Ox_{X,\upmu_X}}\calD_{X/S,\upmu_X}^{\gamma}.
$$
\end{definition}

\begin{parag}\label{ommek5}
Assume $(H_1)$ is satisfied. By \ref{indremark73}, a locally PD-nilpotent $\calD_{X/S,\upmu_X}^{\gamma}$-module $\calE$ \eqref{deflocnilind} is equivalent to the data, for every local section $s$ of $\upmu_X$ over an étale $X$-scheme $U,$ of a locally PD-nilpotent $\calD_{U/S}^{\gamma}$-module $\calE_s$ satisfying the conditions given in \ref{indremark73}. By \ref{thm1237}, the data $\calE_s$ is equivalent to the data of an $\Ox_U$-module $\calE_s$ equipped with a $\calQ_{1|U}$-stratification $\epsilon_s.$ This is in turn equivalent to an $\Ox_{X,\upmu_X}$-module $\calE$ equipped with a $\calQ_{1,\upmu_X}$-stratification. We obtain an equivalence of categories
\begin{equation}\label{starrr1}
\begin{Bmatrix}\mathrm{Locally\ PD}\text{-}\mathrm{nilpotent}\\\calD_{X/S,\upmu_X}^{\gamma}\text{-}\mathrm{modules} \end{Bmatrix} \xrightarrow{\sim} \begin{Bmatrix}\Ox_{X,\upmu_X}\text{-}\mathrm{modules\ equipped\ with\ a\ }\calQ_{1,\upmu_X}\text{-}\mathrm{stratification} \end{Bmatrix}.
\end{equation}
Similarly, we have an equivalence of categories
\begin{equation}\label{starrr2}
\begin{Bmatrix}\mathrm{Locally\ PD}\text{-}\mathrm{nilpotent}\\ \left (\w{\Gamma}^{\bullet}\calT_{X'/S}\right )_{\upmu_X}\text{-}\mathrm{modules} \end{Bmatrix} \xrightarrow{\sim} \begin{Bmatrix}\Ox_{X',\upmu_X}\text{-}\mathrm{modules\ equipped\ with\ an\ }\calR'_{1,\upmu_X}\text{-}\mathrm{stratification} \end{Bmatrix}.
\end{equation}
\end{parag}

\begin{proposition}\label{propKoko19129}
Assume $(H_1)$ is satisfied. The equivalences \eqref{starrr1} and \eqref{starrr2} induce equivalences of categories
$$
\begin{Bmatrix}\mathrm{Locally\ PD}\text{-}\mathrm{nilpotent}\\\widetilde{\calD}_{X/S}^{\gamma}\text{-}\mathrm{modules} \end{Bmatrix} \xrightarrow{\sim} \begin{Bmatrix}\upmu_X\text{-}\mathrm{indexed\ }\calA_X\text{-}\mathrm{modules\ equipped}\\ \mathrm{with\ an\ admissible\ }\calQ_{1,\upmu_X}\text{-}\mathrm{stratification} \end{Bmatrix},
$$
$$
\begin{Bmatrix}\mathrm{Locally\ PD}\text{-}\mathrm{nilpotent}\\\left (\left (\w{\Gamma}^{\bullet}\calT_{X'/S}\right )_{\upmu_X}\otimes_{\Ox_{X',\upmu_X}}\calB_{X/S} \right )\text{-}\mathrm{modules} \end{Bmatrix} \xrightarrow{\sim} \begin{Bmatrix}\upmu_X\text{-}\mathrm{indexed\ }\calB_{X/S}\text{-}\mathrm{modules\ equipped}\\ \mathrm{with\ an\ admissible\ }\calR'_{1,\upmu_X}\text{-}\mathrm{stratification} \end{Bmatrix}.
$$
\end{proposition}

\begin{proof}
The first equivalence results from \ref{prop1932Ko} and \eqref{eq12331}. The second one follows from the isomorphism \eqref{eq1182} and the $\calB_{X/S}$-linearity of admissible $\calR'_{1,\upmu_X}$-stratifications.
\end{proof}

\begin{theorem}\label{THMXXaa}
Assume $(H_2)$ is satisfied. The diagram
\begin{equation}
\begin{tikzcd}
\begin{Bmatrix}\upmu_X\text{-}\mathrm{indexed\ }\calB_{X/S}\text{-}\mathrm{modules\ with}\\ \mathrm{an\ admissible}\\ \calR_{1,\upmu_X}'\text{-}\mathrm{stratification} \end{Bmatrix} \ar[swap,sloped]{d}{\sim} \ar{d}{\omega'} \ar{r}{L_{\mathrm{strat}} \circ \Psi_{\upmu}} & \begin{Bmatrix}\upmu_X\text{-}\mathrm{indexed\ }\calA_X\text{-}\mathrm{modules\ with}\\ \mathrm{an\ admissible}\\Q_{1,\upmu_X}\text{-}\mathrm{stratification} \end{Bmatrix} \ar[sloped]{d}{\sim} \ar[swap]{d}{\omega} \\
\begin{Bmatrix}\mathrm{Locally\ PD}\text{-}\mathrm{nilpotent\ } \upmu_X\text{-}\mathrm{indexed}\\\left (\left (\w{\Gamma}^{\bullet}\calT_{X'/S}\right )_{\upmu_X}\otimes_{\Ox_{X',\upmu_X}}\calB_{X/S} \right )\text{-}\mathrm{modules} \end{Bmatrix} \ar{r}{\Phi_{\upmu}} & \begin{Bmatrix}\mathrm{Locally\ PD}\text{-}\mathrm{nilpotent\ }\upmu_X\text{-}\mathrm{indexed}\\\widetilde{\calD}_{X/S}^{\gamma}\text{-}\mathrm{modules} \end{Bmatrix},
\end{tikzcd}
\end{equation}
where the vertical arrows are given in \ref{propKoko19129} and $L_{\mathrm{strat}},$ $\Psi_{\upmu}$ and $\Phi_{\upmu}$ are given respectively in \eqref{eqLstrat1942f}, \eqref{diag186} and \ref{propKoko1717}, is 2-commutative.
\end{theorem}

\begin{proof}
We may suppose that the hypothesis \ref{loccoord} is satisfied. Let $(\partial_1',\hdots,\partial_d') \in \calT_{X'/S}$ be the dual basis of $(\op{dlog}m_1',\hdots,\op{dlog}m_d')$ and denote by $\partial_I$ the differential operators of $\calD_{X/S},$ obtained from $\op{dlog}m_1,\hdots,\op{dlog}m_d$ \eqref{P2}. For $I\in \N^d,$ set
\begin{alignat*}{1}
\partial'^{[I]}=\prod_{i=1}^d\partial_i'^{[I_i]} \in \Gamma^{\bullet}\calT_{X'/S},\ \eta_{(r)}'^I=\prod_{i=1}^d\eta_{i(r)}'^{I_i}\in \calR_1', \\
\eta^I=\prod_{i=1}^d\eta_i^{I_i}\in \calQ_1,\ \eta_{(q)}^I=\prod_{i=1}^d\eta_{i(q)}^{I_i}\in \calQ_1.
\end{alignat*}
The family $\left (\eta_{(r)}'^{I} \right )_{I\in \N^d}$ is a basis for the $\Ox_X$-module $\calR_1'.$
Let $\left (\varphi_I' \right )_{I\in \N^d}$ be its dual basis. We also have a basis $\left (\eta^K\eta_{(q)}^L \right )_{\substack{K\in \llbracket 0,p-1 \rrbracket^d\\ L\in \N^d}}$ for the $\Ox_X$-module $\calQ_1.$ We denote by $\left (\varphi_{K,L} \right )_{\substack{K\in \llbracket 0,p-1 \rrbracket^d \\ L\in \N^d}}$ its dual basis. Let $\nu:Q_1 \ra R_1'$ be the morphism induced by \eqref{nu} and $\calV:\calR_1' \ra \calQ_1$ the corresponding morphism of Hopf algebras.
Before starting the proof, let us recall isomorphisms that we will need:
\begin{itemize}
\item By (\cite{SGA43} IV 12.3), there exists canonical isomorphisms
$$
\left (\calR'^{\vee}_1 \right )_{\upmu_X} \xrightarrow{\sim} \left (\calR'_{1,\upmu_X} \right )^{\vee},\ \left (\calQ^{\vee}_1 \right )_{\upmu_X} \xrightarrow{\sim} \left (\calQ_{1,\upmu_X} \right )^{\vee}.
$$
\item We have an isomorphism \eqref{eq12331}
\begin{equation}\label{eqmuzanf1}
\calD_{X/S,\upmu_X}^{\gamma} = \calD_{X/S,\upmu_X} \otimes_{\left (S^{\bullet}\calT_{X'/S} \right )_{\upmu_X}} \left (\w{\Gamma}^{\bullet}\calT_{X'/S}\right )_{\upmu_X}\xrightarrow{\sim} \calQ_{1,\upmu_X}^{\vee},
\end{equation}
by which, $(-1)^{|J|}I! \partial_I \otimes \partial'^J$ corresponds to $\varphi_{I,J}.$ 
\item We have an isomorphism \eqref{eq1182}
\begin{equation}\label{eqMuzanf2}
\left (\w{\Gamma}^{\bullet}\calT_{X'/S}\right )_{\upmu_X} \xrightarrow{\sim} \calR_{1,\upmu_X}'^{\vee},
\end{equation}
by which, the local section $\partial'^{[I]}$ corresponds to $\varphi'_{I}.$
\end{itemize}
Let $(\calE',\epsilon')$ be an object of the upper left category i.e. $\calE'$ is a $\upmu_X$-indexed $\calB_{X/S}$-module equipped with an admissible $\calR_{1,\upmu_X}'$-stratification
$$
\epsilon':\calR_{1,\upmu_X}' \otimes_{\Ox_{X',\upmu_X}} \calE' \ra \calE' \otimes_{\Ox_{X',\upmu_X}} \calR_{1,\upmu_X}'.
$$
We have to prove the existence of an isomorphism
$$
\omega \circ L_{\mathrm{strat}} \circ \Psi_{\upmu} \left ( \calE',\epsilon' \right ) \xrightarrow{\sim} \Phi_{\upmu} \circ \omega' \left ( \calE' , \epsilon' \right ).
$$
We start by expliciting the action of $\widetilde{\calD}_{X/S}^{\gamma}$ on both of the these modules.
Let $s,t,u\in \Gamma(U,\upmu_X)$ be sections of $\upmu_X$ over an étale $X$-scheme $U$ and $a,$ $b$ and $x'$ local sections of $\calA_{X,s},$ $\calB_{X/S,u}$ and $\calE'_t$ respectively.
The image of $(\calE',\epsilon')$ by $\omega'$ is the $\upmu_X$-indexed $\left (\left (\w{\Gamma}^{\bullet}\calT_{X'/S}\right )_{\upmu_X}\otimes_{\Ox_{X',\upmu_X}}\calB_{X/S} \right )$-module $\calE'.$ Let us explicit the actions of $\left (\w{\Gamma}^{\bullet}\calT_{X'/S}\right )_{\upmu_X}$ and $\calB_{X/S}$ on $\calE':$
\begin{itemize}
\item First, the action of $\calB_{X/S}$ on $\calE'$ comes from the structure of $\upmu_X$-indexed $\calB_{X/S}$-module on $\calE'.$
\item Second, the action of $\left (\w{\Gamma}^{\bullet}\calT_{X'/S}\right )_{\upmu_X}$ on $\calE'$ is given by
$$
\left (\w{\Gamma}^{\bullet}\calT_{X'/S}\right )_{\upmu_X} \times \calE' \xrightarrow{\sim} \calR_{1,\upmu_X}'^{\vee}\times \calE' \ra \calE',
$$
where the first isomorphism is induced by \eqref{eqMuzanf2} and the second morphism is
$$
\begin{array}[t]{clc}
\calR_{1,\upmu_X}'^{\vee} \times \calE' & \ra & \calE' \\
(\varphi,x') & \mapsto & \left (\op{Id}_{\calE'} \otimes \varphi \right ) \circ \epsilon'(1\otimes x').
\end{array}
$$
By the $\calB_{X/S}$-linearity of $\epsilon',$ we have
\begin{alignat*}{2}
\varphi \cdot (bx') &= \left (\op{Id}_{\calE'}\otimes \varphi \right )\circ \epsilon'(1\otimes (bx')) \\
&=\left (\op{Id}_{\calE'} \otimes \varphi\right )(b \epsilon'(1\otimes x')) \\
&= b\left (\op{Id}_{\calE'} \otimes \varphi\right )( \epsilon'(1\otimes x')) = b\left (\varphi'\cdot x' \right ).
\end{alignat*}
By the isomorphism \eqref{eqMuzanf2}, $\partial'^{[I]}$ corresponds to $\varphi'_I.$
The action of
$$
\partial'^{[I]}\otimes b\in \left (\left (\w{\Gamma}^{\bullet}\calT_{X'/S}\right )_{\upmu_X}\otimes_{\Ox_{X',\upmu_X}}\calB_{X/S} \right )_u=\w{\Gamma}^{\bullet}\calT_{U'/S}\otimes_{\Ox_{U'}}\calB_{X/S,u}
$$
on $x'\in \calE'_t$ is then given by
\begin{equation}\label{eq20z}
\left (\partial'^{[I]}\otimes b\right )\cdot x'=b\left (\varphi_I'\cdot x' \right ) \in \calE'_{u+t}.
\end{equation}
\end{itemize}
Let $\calE=\Phi_{\upmu}(\calE').$ Then $\calE=\calA_X \circledast_{\calB_{X/S}}\calE'$ (\eqref{ren6} and \ref{kraz827}). Let us explicit the action of $\widetilde{\calD}^{\gamma}_{X/S}$ on $\calE.$ For that, we explicit the actions of $\calD_{X/S,\upmu_X}^{\gamma}$ and $\calA_X$ on $\calE:$
\begin{itemize}
\item The action of $\calA_X$ comes from the canonical structure of $\upmu_X$-indexed $\calA_X$-module on $\calE.$
\item By \ref{indproptensprod}, we can consider $a\otimes x'$ as a local section of $\calE_{s+t}.$ Let $P_0$ be the special fiber of $P_{\frakX/\frakS,0}$ and $p_1,p_2:P_0 \ra X$ the canonical projections. By \eqref{actionpartialax}, the local sections $c_1,\hdots,c_d$ of $\Ox_X,$ appearing in \eqref{eqKoko1323}, such that the action of $\partial_{\epsilon_i}\in \calD_{X/S}$ on $a\otimes x'$ is given by
\begin{equation}\label{eq20zz}
\partial_{\epsilon_i}\cdot (a\otimes x')=(\partial_{\epsilon_i}\cdot a) \otimes x' +a\otimes (\partial'_{i}\cdot x')+\sum_{j=1}^d\partial_{\epsilon_i}(p_2^{\#}c_j-p_1^{\#}c_j) a \otimes (\partial'_j\cdot x'),
\end{equation}
where $\partial_{\epsilon_i}\cdot a$ denotes the action of $\partial_{\epsilon_i} \in \calD_{X/S}$ on $a\in \calA_{X,s}$ given by the canonical connection $d_{\calA_X}.$
By \eqref{actionpartialax}, the action $\partial'^{[L]} \cdot (a\otimes x')$ of $\partial'^{[L]} \in \Gamma^{\bullet} \calT_{X'/S}$ is given, for $L\in \N^d,$ by
\begin{equation}\label{eq20zztt}
\partial'^{[L]} \cdot (a\otimes x')=a \otimes \left ( \partial'^{[L]} \cdot x' \right ).
\end{equation}
\end{itemize}
To conclude, the module structure on $\Phi_{\upmu} \circ \omega' \left (\calE',\epsilon' \right )$ is given by
\begin{equation}\label{Muzanf3}
\begin{alignedat}{2}
\partial_{\epsilon_i}\cdot (a\otimes x') &= (\partial_{\epsilon_i}\cdot a) \otimes x' +a\otimes (\partial'_{i}\cdot x')+\sum_{j=1}^d\partial_{\epsilon_i}(p_2^{\#}c_j-p_1^{\#}c_j) a \otimes (\partial'_j\cdot x'), \\
\partial'^{[L]} \cdot (a\otimes x') &= a \otimes \left ( \partial'^{[L]} \cdot x' \right ).
\end{alignedat}
\end{equation}
Let us now determine $\omega \circ L_{\mathrm{strat}}\circ \Psi_{\upmu}\left (\calE',\epsilon'\right ):$
We consider $\calQ_{1,\upmu_X} \otimes_{\Ox_{X,\upmu_X}} \calA_X$ as a $\upmu_X$-indexed $\left (\calR_{1,\upmu_X}' \otimes_{\Ox_{X',\upmu_X}} \calB_{X/S} \right )$-module via $\nu:Q_1 \ra R_1'$ and the canonical embedding $\calB_{X/S} \hookrightarrow \calA_X.$
The image of $(\calE',\epsilon')$ by $\Psi_{\upmu}$ \eqref{diag186} is 
$$
\Psi_{\upmu}\left (\calE',\epsilon'\right )=\left (F_1^*\calE',\calV^*\epsilon' \right ).
$$
This is sent, by $L_{\text{strat}}$ \eqref{eqLstrat1942f}, to
$$
\left (\calA_X \circledast_{F_1^*\calB_{X/S}} F_1^*\calE',\epsilon_{\calA_X,Q}\circledast_{\calQ_{1,\upmu_X} \otimes_{\Ox_{X',\upmu_X}}\calB_{X/S} } \calV^*\epsilon' \right ).
$$
By \ref{kraz827}, we have canonical isomorphisms
\begin{alignat*}{2}
\calA_X \circledast_{F_1^*\calB_{X/S}} F_1^*\calE' \xrightarrow{\sim} & \calA_X \circledast_{\calB_{X/S}} \calE', \\
\epsilon_{\calA_X,Q}\circledast_{\calQ_{1,\upmu_X} \otimes_{\Ox_{X',\upmu_X}}\calB_{X/S} } \calV^*\epsilon'
 \xrightarrow{\sim} & \epsilon_{\calA_X,Q}\circledast_{\calV^*\left (\calR_{1,\upmu_X}' \otimes_{\Ox_{X',\upmu_X}}\calB_{X/S}\right ) } \calV^*\epsilon'\\
\xrightarrow{\sim} &  \epsilon_{\calA_X,Q} \circledast_{\calR_{1,\upmu_X}' \otimes_{\Ox_{X',\upmu_X}} \calB_{X/S} } \epsilon'.
\end{alignat*}
We conclude that the image of $(\calE',\epsilon')$ by $L_{\mathrm{strat}} \circ \Psi_{\upmu}$ is $(\calE,\epsilon)$ where
$$
\begin{cases}
\calE=\calA_X \circledast_{\calB_{X/S}}\calE' \\
\epsilon = \epsilon_{\calA_X,Q} \circledast_{\calR_{1,\upmu_X}' \otimes_{\Ox_{X',\upmu_X}} \calB_{X/S} } \epsilon'.
\end{cases}
$$
By \eqref{eqepsilonAQ}, we have
\begin{alignat*}{2}
\epsilon_{\calA_X,Q,s}(1\otimes a) &= \sum_{K \in \llbracket 0,p-1 \rrbracket^d} (\varphi_{K,0}\cdot a)\otimes \eta^K, \\
\epsilon_t'(1\otimes x') &= \sum_{I\in \N^d}(\varphi'_I\cdot x')\otimes \eta_{(r)}'^I.
\end{alignat*}
It follows that
\begin{equation}\label{eq2030}
\epsilon_{s+t} (1\otimes a\otimes x')=\sum_{\substack{K \in \llbracket 0,p-1 \rrbracket^d \\ I\in \N^d}} (\varphi_{K,0}\cdot a)\otimes (\varphi_I'\cdot x') \otimes \eta^K \calV \left (\eta_{(r)}^I \right ).
\end{equation}
By \eqref{eqKoko1323}, we have
\begin{alignat*}{2}
\calV \left (\eta_{i(r)} \right ) &= \eta_{i(q)} + \sum_{k=1}^{p-1}\frac{(-1)^{k+1}}{k}\eta_i^k + q_2^{\#}c_i-q_1^{\#}c_i \\
&= \eta_{i(q)} + \sum_{k=1}^{p-1}\frac{(-1)^{k+1}}{k}\eta_i^k + \sum_{\substack{\alpha \in \llbracket 0,p-1 \rrbracket^d \\ \beta\in \N^d}} \varphi_{\alpha,\beta}(q_2^{\#}c_i-q_1^{\#}c_i)\eta^{\alpha}\eta_{(q)}^{\beta},
\end{alignat*}
where $q_1,q_2:Q_1 \ra X$ are the canonical projections. Note that the local sections $c_i$ are the same in \eqref{eq20zz}.
Let $X\ra Y$ be the exact diagonal immersion. The projection $q_i$ factors, by definition, through $Y.$ It follows, by \ref{lemkraz}, that $\varphi_{\alpha,\beta}(q_2^{\#}c_i-q_1^{\#}c_i)=0$ whenever $\beta \neq 0$ and if $(\alpha,\beta)=(0,0).$ Then
\begin{equation}\label{eqfinalMuzan33}
\calV \left (\eta_{i(r)} \right ) = \eta_{i(q)} + \sum_{k=1}^{p-1}\frac{(-1)^{k+1}}{k}\eta_i^k + \sum_{\substack{\alpha \in \llbracket 0,p-1 \rrbracket^d \\ \alpha\neq 0}} \varphi_{\alpha,0}(q_2^{\#}c_i-q_1^{\#}c_i)\eta^{\alpha}.
\end{equation}
Denote by $\calK_0 \subset \calQ_1$ and $\calK_1 \subset \calQ_1$ the ideals generated by $\eta_1,\hdots,\eta_d$ and by $\eta_1,\hdots,\eta_d,\eta_{1(q)}, \hdots, \eta_{d(q)}$ respectively. Also denote by $\pi_0:\calQ_1 \ra \calQ_1/\calK_0,$ $\pi_1:\calQ_1 \ra \calQ_1/\calK_1^2$ the canonical projections. We have, by \eqref{eqfinalMuzan33},
$$
\pi_0 \circ \calV\left (\eta_{i(r)} \right ) = \eta_{i(q)}.
$$
It follows that, for $L\in \N^d,$
$$
\varphi_{0,L}\left ( \eta^K \calV \left ( \eta_{(r)}^I \right ) \right ) = \begin{cases}
0\ \text{if\ }(K,I) \neq (0,L) \\
1\ \text{if\ }(K,I)=(0,L).
\end{cases}
$$
Then, by \eqref{eq2030},
\begin{equation}\label{eqkraz25}
\left ( \op{Id} \otimes \varphi_{0,L} \right )\circ \epsilon_{s+t}(1 \otimes a \otimes x')=a \otimes \left ( \varphi_L' \cdot x' \right ).
\end{equation}
We also have, by \eqref{eqfinalMuzan33},
$$
\pi_1\circ \calV\left (\eta_{i(r)} \right ) = \eta_{i(q)} + \eta_i + \varphi_{\epsilon_1,0}(q_2^{\#}c_1-q_1^{\#}c_1)\eta_1+ \hdots +\varphi_{\epsilon_d,0}(q_2^{\#}c_d-q_1^{\#}c_d)\eta_d.
$$ 
By \eqref{eq2030}, we have
\begin{equation}\label{eq2030kraz}
\begin{alignedat}{2}
\left (\op{Id}\otimes \pi_1 \right )\circ \epsilon_{s+t} (1\otimes a\otimes x') =& \sum_{1\le i\le d } (\varphi_{\epsilon_i,0}\cdot a)\otimes x' \otimes \eta_i + \sum_{1\le i \le d}a\otimes \left (\varphi_{\epsilon_i}'\cdot x' \right )\otimes \pi_1 \circ \calV\left (\eta_{i(r)} \right ) \\
=& \sum_{1\le i\le d} (\varphi_{\epsilon_i,0}\cdot a)\otimes x' \otimes \eta_i + a\otimes \left (\varphi_{\epsilon_i}'\cdot x' \right )\otimes \left (\eta_{i(q)}+\eta_i \right ) \\
&+ \sum_{1\le i,j\le d} a\otimes \left (\varphi_{\epsilon_i}'\cdot x' \right )\otimes \left (\varphi_{\epsilon_j,0}\left (q_2^{\#}c_j-q_1^{\#}c_j\right )\eta_j \right ).
\end{alignedat}
\end{equation}
By \eqref{eqkraz25} and \eqref{eq2030kraz}, we get
\begin{equation}\label{eq20zzz}
\begin{alignedat}{2}
\varphi_{\epsilon_i,0}\cdot (a\otimes x') &= (\varphi_{\epsilon_i,0}\cdot a)\otimes x' + a\otimes (\varphi_{\epsilon_i}'\cdot x')+\sum_{j=1}^d \varphi_{\epsilon_j,0}(q_2^{\#}c_j-q_1^{\#}c_j) a\otimes (\varphi'_{\epsilon_j}\cdot x'),\\
\varphi_{0,L} \cdot (a\otimes x') &= a \otimes (\varphi'_{L} \cdot x').
\end{alignedat}
\end{equation}
Since $\varphi_{I,J}$ corresponds to $(-1)^{|J|}I! \partial_I \otimes \partial'^{[J]}$ by \eqref{eqmuzanf1} and $\varphi_I'$ corresponds to $\partial'^{[I]}$ by \eqref{eqMuzanf2}, we get
\begin{equation}\label{eq20zzzv2}
\begin{alignedat}{2}
\partial_{\epsilon_i}\cdot (a\otimes x') &= (\partial_{\epsilon_i}\cdot a)\otimes x' + a\otimes (\partial_i'\cdot x')+\sum_{j=1}^d \varphi_{\epsilon_j,0}(q_2^{\#}c_j-q_1^{\#}c_j) a\otimes (\partial_j'\cdot x'),\\
\partial'^{[L]} \cdot (a\otimes x') &= (-1)^{|L|}a \otimes (\partial'^{[L]} \cdot x').
\end{alignedat}
\end{equation}
We have thus determined $\omega \circ L_{\mathrm{strat}} \circ \Psi_{\upmu} \left ( \calE' , \epsilon' \right ).$

As we see, there are two differences between \eqref{Muzanf3} and \eqref{eq20zzzv2}: a difference of $(-1)^{|L|}$ in the second equality and, in the first equality, we have $\partial_{\epsilon_i}(p_2^{\#}c-p_1^{\#}c)$ on one hand and $\varphi_{\epsilon_j,0}(q_2^{\#}c-q_1^{\#}c)$ in the other. To fix the first difference, we introduce the antipode isomorphism
$$
\sigma:\begin{array}[t]{clc}
\left (\w{\Gamma}^{\bullet}\calT_{X'/S}\right )_{\upmu_X} & \xrightarrow{\sim} & \left (\w{\Gamma}^{\bullet}\calT_{X'/S}\right )_{\upmu_X} \\
\partial'_i & \mapsto & -\partial'_i,
\end{array}
$$
and the functor $T$ defined by base change by $\sigma$
$$
T: \calM \mapsto \calM \otimes_{\left (\w{\Gamma}^{\bullet}\calT_{X'/S}\right )_{\upmu_X},\sigma} \left (\w{\Gamma}^{\bullet}\calT_{X'/S}\right )_{\upmu_X}.
$$
Since $\sigma$ is an isomorphism, we have an isomorphism
$$
\Phi_{\upmu} \circ \omega' \left ( \calE' , \epsilon' \right ) \xrightarrow{\sim} T \circ \Phi_{\upmu} \circ \omega' \left ( \calE' , \epsilon' \right ).
$$
In addition, $\sigma \left ( \partial'^{[L]} \right ) = (-1)^{|L|}\partial'^{[L]}.$ It follows that the actions of $\partial_{\epsilon_i}$ and $\partial'^{[L]}$ on the local section $a\otimes x'$ of $T \circ \Phi_{\upmu} \circ \omega' \left ( \calE' , \epsilon' \right )$ are given by
\begin{equation}\label{eq20zzzv3}
\begin{alignedat}{2}
\partial_{\epsilon_i}\cdot (a\otimes x') &= (\partial_{\epsilon_i}\cdot a)\otimes x' + a\otimes (\partial_i'\cdot x')+\sum_{j=1}^d \varphi_{\epsilon_j,0}(q_2^{\#}c_j-q_1^{\#}c_j) a\otimes (\partial_j'\cdot x'),\\
\partial'^{[L]} \cdot (a\otimes x') &= a \otimes (\partial'^{[L]} \cdot x').
\end{alignedat}
\end{equation}
The last thing to check is that, for a local section $c$ of $\Ox_X,$
$$
\partial_{\epsilon_i}(p_2^{\#}c-p_1^{\#}c) = \varphi_{\epsilon_j,0}(q_2^{\#}c-q_1^{\#}c).
$$
Recall the morphism of Hopf algebras $u:\calQ_1 \ra \calP_0$ \eqref{Muzanu1}. By definition of $u,$ we have $u \circ q_i^{\#}=p_i^{\#}.$ Thus
$$
\partial_{\epsilon_i}(p_2^{\#}c-p_1^{\#}c)=\partial_{\epsilon_i} \circ u \left ( q_2^{\#}c-q_1^{\#}c \right ).
$$
It is hence sufficient to prove that, as linear forms on $\calQ_1,$ we have
$$
\partial_{\epsilon_i} \circ u = \varphi_{\epsilon_i,0}.
$$
This is immediate since $u(\eta_j)=\eta_j$ for all $1\le j\le d$ \eqref{ueta}.

By \eqref{Muzanf3} and \eqref{eq20zzzv3}, we conclude that we have an isomorphism
$$
T \circ \Phi_{\upmu} \circ \omega' \left ( \calE' , \epsilon' \right ) \xrightarrow{\sim} \omega \circ L_{\mathrm{strat}} \circ \Psi_{\upmu} \left ( \calE' , \epsilon' \right ).
$$
\end{proof}

\begin{corollaire}\label{THMXX}
Assume $(H_2)$ is satisfied. The diagram
\begin{equation}
\begin{tikzcd}
\begin{Bmatrix}\mathrm{Crystals\ of\ }\upmu'\text{-}\mathrm{indexed}\\ \frakB_X\text{-}\mathrm{modules\ of\ }\widetilde{\calE}(X'/\frakS)_{/\upmu'} \end{Bmatrix} \ar{r}{C^{-1}} \ar[swap,sloped]{d}{\sim} & \begin{Bmatrix}\mathrm{Crystals\ of\ }\underline{\upmu}\text{-}\mathrm{indexed}\\ \frakA_X\text{-}\mathrm{modules\ of\ }\widetilde{\underline{\calE}}(X/\frakS)_{/\underline{\upmu}} \end{Bmatrix} \ar[sloped]{d}{\sim}\\
\begin{Bmatrix}\upmu_X\text{-}\mathrm{indexed\ }\calB_{X/S}\text{-}\mathrm{modules\ with}\\ \mathrm{an\ admissible}\\ \calR_{1,\upmu_X}'\text{-}\mathrm{stratification} \end{Bmatrix} \ar[swap,sloped]{d}{\sim} \ar{d}{\omega'} \ar{r}{L_{\mathrm{strat}} \circ \Psi_{\upmu}} & \begin{Bmatrix}\upmu_X\text{-}\mathrm{indexed\ }\calA_X\text{-}\mathrm{modules\ with}\\ \mathrm{an\ admissible}\\Q_{1,\upmu_X}\text{-}\mathrm{stratification} \end{Bmatrix} \ar[sloped]{d}{\sim} \ar[swap]{d}{\omega} \\
\begin{Bmatrix}\mathrm{Locally\ PD}\text{-}\mathrm{nilpotent\ } \upmu_X\text{-}\mathrm{indexed}\\\left (\left (\w{\Gamma}^{\bullet}\calT_{X'/S}\right )_{\upmu_X}\otimes_{\Ox_{X',\upmu_X}}\calB_{X/S} \right )\text{-}\mathrm{modules} \end{Bmatrix} \ar{r}{\Phi_{\upmu}} & \begin{Bmatrix}\mathrm{Locally\ PD}\text{-}\mathrm{nilpotent\ }\upmu_X\text{-}\mathrm{indexed}\\\widetilde{\calD}_{X/S}^{\gamma}\text{-}\mathrm{modules} \end{Bmatrix},
\end{tikzcd}
\end{equation}
where the upper vertical equivalences are given in \ref{equivRQAB}, the lower vertical arrows are given in \ref{propKoko19129} and $C^{-1},$ $L_{\mathrm{strat}},$ $\Psi_{\upmu}$ and $\Phi_{\upmu}$ are given respectively in \eqref{ffC}, \eqref{eqLstrat1942f}, \eqref{diag186} and \ref{propKoko1717}, is 2-commutative. In addition, the functor $C^{-1}$ is an equivalence of categories.
\end{corollaire}

\begin{proof}
The commutativity of the upper square was proved in \ref{thmKoko20125}. The lower square's commutativity was proved in \ref{THMXXaa}. Since $\Phi_{\upmu}$ is an equivalence of categories \eqref{propKoko1717}, so is $C^{-1}.$
\end{proof}

\begin{lemma}\label{lemEST3}
Assume $(H_1)$ is satisfied. Let $U$ be an étale $X$-scheme. Consider the diagram
$$
\begin{tikzcd}
\begin{Bmatrix}\text{Crystals\ of\ }\underline{\upmu}\text{-indexed}\\ C_{X/\frakS,\upmu}^{-1}\frakB_X\text{-modules\ of\ }\widetilde{\underline{\calE}}(X/\frakS)_{/\underline{\upmu}} \end{Bmatrix} \ar[swap]{dd}{\mathfrak{E} \mapsto \frakA_X \circledast_{C_{X/\frakS,\upmu}^{-1}\frakB_X}\mathfrak{E}} \ar{r}{\underline{R_{X,U}}^{-1}} & \begin{Bmatrix}\text{Crystals\ of\ }\underline{\upmu}_U\text{-indexed}\\ C_{U/\frakS,\upmu}^{-1}\frakB_U\text{-modules\ of\ }\widetilde{\underline{\calE}}(U/\frakS)_{/\underline{\upmu}_U} \end{Bmatrix} \ar{dd}{\mathfrak{E}  \mapsto  \frakA_U \circledast_{C_{U/\frakS,\upmu}^{-1}\frakB_U}\mathfrak{E}} \\ 
 & \\
\begin{Bmatrix}\text{Crystals\ of\ }\underline{\upmu}\text{-indexed}\\ \frakA_X\text{-modules\ of\ }\widetilde{\underline{\calE}}(X/\frakS)_{/\underline{\upmu}} \end{Bmatrix} \ar[swap]{r}{\underline{R_{X,U}}^{-1}} & \begin{Bmatrix}\text{Crystals\ of\ }\underline{\upmu}_U\text{-indexed}\\ \frakA_U\text{-modules\ of\ }\widetilde{\underline{\calE}}(U/\frakS)_{/\underline{\upmu}_U} \end{Bmatrix},
\end{tikzcd}
$$
where the horizontal arrows are given in \eqref{eqREST} and the linearization functors are given in \eqref{ffL}. Then, this diagram is commutative.
\end{lemma}

\begin{proof}
For simplicity, we denote $R_{X',U'}$ and $\underline{R_{X,U}}$ by $R'$ and $\underline{R}$ respectively. Let $\frakE$ be a crystal of $C_{X/\frakS,\upmu}^{-1}\frakB_X$-modules. We have to prove the existence of a canonical isomorphism
\begin{equation}\label{almostthere1}
\underline{R}^{-1} \left ( \frakA_X \circledast_{C_{X/\frakS,\upmu}^{-1}\frakB_X}\mathfrak{E} \right ) \xrightarrow{\sim} \frakA_U \circledast_{C_{U/\frakS,\upmu}^{-1}\frakB_U}\underline{R}^{-1}\mathfrak{E}.
\end{equation}
In \ref{paragEST1}, we have seen that $\underline{R}$ is a localization morphism. By \ref{ren3}, we deduce a canonical isomorphism
$$
\underline{R}^{-1} \left ( \frakA_X \circledast_{C_{X/\frakS,\upmu}^{-1}\frakB_X}\mathfrak{E} \right ) \xrightarrow{\sim} \left ( \underline{R}^{-1}  \frakA_X \right ) \circledast_{\underline{R}^{-1}C_{X/\frakS,\upmu}^{-1}\frakB_X} \left (  \underline{R}^{-1} \mathfrak{E} \right ).
$$
By \ref{AXUESTiso}, we have canonical isomorphisms
$$
\underline{R}^{-1}\frakA_X \xrightarrow{\sim} \frakA_U,
$$
$$
\underline{R}^{-1}C_{X/\frakS,\upmu}^{-1}\frakB_X \xrightarrow{\sim} C_{U/\frakS,\upmu}^{-1}R'^{-1}\frakB_X \xrightarrow{\sim} C_{U/\frakS,\upmu}^{-1}\frakB_U.
$$
The isomorphism \eqref{almostthere1} follows.
\end{proof}

\begin{lemma}\label{lemEST2}
Assume $(H_1)$ is satisfied. Let $U$ be an étale $X$-scheme. Consider the diagram
$$
\begin{tikzcd}
\begin{Bmatrix}\text{Crystals\ of\ }\upmu'\text{-indexed}\\ \frakB_X\text{-modules\ of\ }\widetilde{\calE}(X'/\frakS)_{/\upmu'} \end{Bmatrix} \ar{r}{R_{X',U'}^{-1}} \ar[swap]{d}{C_{X/\frakS,\upmu}^{-1}} & \begin{Bmatrix}\text{Crystals\ of\ }\upmu'_U\text{-indexed}\\ \frakB_U\text{-modules\ of\ }\widetilde{\calE}(U'/\frakS)_{/\upmu'_U} \end{Bmatrix} \ar{d}{C_{U/\frakS,\upmu}^{-1}} \\
\begin{Bmatrix}\text{Crystals\ of\ }\underline{\upmu}\text{-indexed}\\ C_{X/\frakS,\upmu}^{-1}\frakB_X\text{-modules\ of\ }\widetilde{\underline{\calE}}(X/\frakS)_{/\underline{\upmu}} \end{Bmatrix} \ar[swap]{dd}{\mathfrak{E} \mapsto \frakA_X \circledast_{C_{X/\frakS,\upmu}^{-1}\frakB_X}\mathfrak{E}} \ar[swap]{r}{\underline{R_{X,U}}^{-1}} & \begin{Bmatrix}\text{Crystals\ of\ }\underline{\upmu}_U\text{-indexed}\\ C_{U/\frakS,\upmu}^{-1}\frakB_U\text{-modules\ of\ }\widetilde{\underline{\calE}}(U/\frakS)_{/\underline{\upmu}_U} \end{Bmatrix} \ar{dd}{\mathfrak{E}  \mapsto  \frakA_U \circledast_{C_{U/\frakS,\upmu}^{-1}\frakB_U}\mathfrak{E}} \\ 
 & \\
\begin{Bmatrix}\text{Crystals\ of\ }\underline{\upmu}\text{-indexed}\\ \frakA_X\text{-modules\ of\ }\widetilde{\underline{\calE}}(X/\frakS)_{/\underline{\upmu}} \end{Bmatrix} \ar[swap]{r}{\underline{R_{X,U}}^{-1}} & \begin{Bmatrix}\text{Crystals\ of\ }\underline{\upmu}_U\text{-indexed}\\ \frakA_U\text{-modules\ of\ }\widetilde{\underline{\calE}}(U/\frakS)_{/\underline{\upmu}_U} \end{Bmatrix},
\end{tikzcd}
$$
where the horizontal arrows are given in \eqref{eqREST}, $C_{X/\frakS,\upmu}^{-1}$ and $C_{U/\frakS,\upmu}^{-1}$ are given in \eqref{CESTeq} and the linearization functors are given in \eqref{ffL}. Then, this diagram is commutative.
\end{lemma}

\begin{proof}
The commutativity of the upper square follows from \ref{commCXU} and the commutativity of the lower square follows from \ref{lemEST3}.
\end{proof}

\begin{theorem}\label{THMXXYY}
The functor \eqref{ffC} 
$$
C^{-1}:
\begin{Bmatrix}\text{Crystals\ of\ }\upmu'\text{-indexed}\\ \frakB_X\text{-modules\ of\ }\widetilde{\calE}(X'/\frakS)_{/\upmu'} \end{Bmatrix}  \ra  \begin{Bmatrix}\text{Crystals\ of\ }\underline{\upmu}\text{-indexed}\\ \frakA_X\text{-modules\ of\ }\widetilde{\underline{\calE}}(X/\frakS)_{/\underline{\upmu}} \end{Bmatrix}
$$
is an equivalence of categories.
\end{theorem}

\begin{proof}
We can work étale locally on $\frakS$ and hence suppose that $\frakS$ has a chart $\frakS \ra B\langle P\rangle.$
Let $\left ( U_i \ra X \right )_{i\in I}$ be an étale covering such that, for every $i\in I,$ there exists a log smooth logarithmic formal scheme $\frakU_i$ over $\frakS_i$ fitting into a cartesian diagram
$$
\begin{tikzcd}
U_i \ar{r} \ar{d} & \frakU_i \ar{d} \\
S \ar{r} & \frakS.
\end{tikzcd}
$$
We can suppose that $\frakU_i$ is affine and $\frakU_i \ra \frakS$ has a chart $P \ra Q_i$ such that $P^{gp} \ra Q_i^{gp}$ is injective, its cokernel has a finite order coprime woth $p$ and the morphism $\frakU_i \ra \frakS\times_{B\langle P \rangle}B\langle Q_i \rangle,$ induced by $\frakU_i \ra \frakS$ and the chart $\frakU_i \ra B\langle Q_i \rangle,$ is strict and étale. By \ref{PFrob}, we have a chart $P \ra Q_i'$ of $U_i' \ra S$ and, by \ref{lemShiho208}, the morphism $U_i' \ra S\times_{A_1[P]}A_1[Q_i'],$ induced by $U_i' \ra S$ and the chart $U_i' \ra A_1[Q_i'],$ is strict and étale.
Then, there exists, for every $i\in I,$ a unique strict étale logarithmic formal scheme $\frakU'_i$ over $\frakS \times_{B\langle P\rangle} B\langle Q_i' \rangle,$ fitting into a cartesian square
$$
\begin{tikzcd}
U_i' \ar{r} \ar{d} & \frakU_i' \ar{d} \\
\frakS \times_{A_1[ P]} A_1[ Q_i'] \ar{r} & \frakS \times_{B\langle P\rangle} B\langle Q_i' \rangle.
\end{tikzcd}
$$
For every $i\in I,$ there exists a lift $\frakU_i \ra \frakU_i'$ fitting into the commutative diagram
$$
\begin{tikzcd}
U_i' \ar{rrr} & & & \frakU_i' \ar{d} \\
U_i \ar{u} \ar{r} & \frakU_i \ar{urr} \ar{r} & \frakS \times_{B\langle P\rangle} B\langle Q_i \rangle \ar{r} & \frakS \times_{B\langle P\rangle} B\langle Q_i' \rangle.
\end{tikzcd}
$$
By \ref{lemEST2}, we have, for every $i\in I,$ a commutative diagram
$$
\begin{tikzcd}
\begin{Bmatrix}\text{Crystals\ of\ }\upmu'\text{-indexed}\\ \frakB_X\text{-modules\ of\ }\widetilde{\calE}(X'/\frakS)_{/\upmu'} \end{Bmatrix} \ar{r}{R_i'^{-1}} \ar[swap]{dd}{C^{-1}} & \begin{Bmatrix}\text{Crystals\ of\ }\upmu'_{U_i}\text{-indexed}\\ \frakB_{U_i}\text{-modules\ of\ }\widetilde{\calE}(U_i'/\frakS)_{/\upmu'_{U_i}} \end{Bmatrix} \ar{dd}{C_i^{-1}} \\
 & \\
\begin{Bmatrix}\text{Crystals\ of\ }\underline{\upmu}\text{-indexed}\\ \frakA_X\text{-modules\ of\ }\widetilde{\underline{\calE}}(X/\frakS)_{/\underline{\upmu}} \end{Bmatrix} \ar[swap]{r}{\underline{R}_i^{-1}} & \begin{Bmatrix}\text{Crystals\ of\ }\underline{\upmu}_{U_i}\text{-indexed}\\ \frakA_{U_i}\text{-modules\ of\ }\widetilde{\underline{\calE}}(U_i/\frakS)_{/\underline{\upmu}_{U_i}} \end{Bmatrix},
\end{tikzcd}
$$
where $R_i'=R_{X',U_i'}$ and $\underline{R}_i=\underline{R_{X,U_i}}.$
By \ref{THMXX}, the functor $C_i^{-1}$ is an equivalence of categories.
For $i,j\in I,$ we have a commutative diagram
\begin{equation}\label{eqdiagESTn}
\begin{tikzcd}
\begin{Bmatrix}\text{Crystals\ of\ }\upmu'_{U_i}\text{-indexed}\\ \frakB_{U_i}\text{-modules\ of\ }\widetilde{\calE}(U_i'/\frakS)_{/\upmu'_{U_i}} \end{Bmatrix} \ar{r}{R_{ij}'^{-1}} \ar[swap]{dd}{C_i^{-1}} & \begin{Bmatrix}\text{Crystals\ of\ }\upmu'_{U_i\times_XU_j}\text{-indexed}\\ \frakB_{U_i\times_XU_j}\text{-modules\ of\ }\widetilde{\calE}(U_i'\times_{X'}U_j'/\frakS)_{/\upmu'_{U_i\times_XU_j}} \end{Bmatrix} \ar{dd}{C_{ij}^{-1}=C_{ji}^{-1}} \\
 & \\
\begin{Bmatrix}\text{Crystals\ of\ }\underline{\upmu}\text{-indexed}\\ \frakA_{U_i}\text{-modules\ of\ }\widetilde{\underline{\calE}}(U_i/\frakS)_{/\underline{\upmu}_{U_i}} \end{Bmatrix} \ar[swap]{r}{\underline{R}_{ij}^{-1}} & \begin{Bmatrix}\text{Crystals\ of\ }\underline{\upmu}_{U_i\times_XU_j}\text{-indexed}\\ \frakA_{U_i\times_XU_j}\text{-modules\ of\ }\widetilde{\underline{\calE}}(U_i\times_XU_j/\frakS)_{/\underline{\upmu}_{U_i\times_XU_j}} \end{Bmatrix},
\end{tikzcd}
\end{equation}
where $R_{ij}'=R_{U_i',U_i'\times_{X'}U_j'}$ and $\underline{R}_{ij}=\underline{R_{U_i,U_i\times_XU_j}}.$

We start by proving that $C^{-1}$ is essentially surjective. Let $\frakE$ be a crystal of $\underline{\upmu}$-indexed $\frakA_X$-modules. Since $C_i^{-1}$ is essentially surjective, there exists a crystal $\frakF_i$ of $\upmu'_{U_i}$-indexed $\frakB_{U_i}$-modules such that
$$
C_i^{-1}\frakF_i=\underline{R}_i^{-1}\frakE.
$$
There exists isomorphisms
$$
\underline{R}_{ij}^{-1}C_i^{-1}\frakF_i=\underline{R}_{ij}^{-1} \underline{R}_i^{-1} \frakE \xrightarrow{\sim} \underline{R}_{ji}^{-1}\underline{R}_j^{-1}\frakE=\underline{R}_{ji}^{-1}C_j^{-1}\frakF_j.
$$
satisfying the cocycle conditions.
We deduce, by \eqref{eqdiagESTn}, isomorphisms
$$
C_{ij}^{-1}R_{ij}'^{-1}\frakF_i \xrightarrow{\sim} C_{ij}^{-1}R_{ji}'^{-1}\frakF_j
$$
satisfying the cocycle conditions. Since $C_{ij}^{-1}$ is an equivalence of categories, we obtain isomorphisms
$$
R_{ij}'^{-1}\frakF_i \xrightarrow{\sim} R_{ji}'^{-1}\frakF_j
$$
satisfying the cocyle conditions. By the stack $\mathbbm{E}'$ \eqref{paragEST2}, the crystals $\frakF_i$ glue together into a crystal $\frakF$ of $\upmu'$-indexed $\frakB_X$-modules such that $R_i'^{-1}\frakF$ is canonically isomorphic to $\frakF_i.$
We then have
$$
\underline{R}_i^{-1}C^{-1}\frakF=C_i^{-1}R_i'^{-1}\frakF\xrightarrow{\sim} C_i^{-1}\F_i=\underline{R}_i^{-1}\frakE.
$$
We deduce an isomorphism
$$
C^{-1}\frakF \xrightarrow{\sim} \frakE,
$$
hence the essential surjectivity of $C^{-1}.$ The fullness is proved in a similar way by descending morphisms.

We now prove that $C^{-1}$ is faithful. Let $f,g:\frakF_1 \ra \frakF_2$ be a morphism of crystals of $\upmu'$-indexed $\frakB_X$-modules such that $C^{-1}(f)=C^{-1}(g).$ We get
$$
\underline{R}_i^{-1}C^{-1}(f)=\underline{R}_i^{-1}C^{-1}g.
$$
This is equivalent to
$$
C_i^{-1}R_i'^{-1}(f)=C_i^{-1}R'^{-1}_i(g).
$$
Since $C_i^{-1}$ is faithful, we get
$$
R_i'^{-1}(f)=R_i'^{-1}(g).
$$
This is true for all $i\in I.$ By the stack $\mathbbm{E}'$ \eqref{paragEST2}, we deduce that $f=g.$
\end{proof}

\newpage

\appendix

\section{Appendix: \texorpdfstring{$p^n$}{p} %
     -connections and stratifications}

The goal of this appendix is to prove proposition \ref{equivstrat}. We keep the notations of section 11 and we start by proving some lemmas.

\begin{lemma}\label{lem116f}
Let $k\ge 1$ and $n\ge 0$ be integers and suppose \ref{loccoord} is satisfied. For any $1\le i\le d,$ let $\xi_i$ be the image of $\frac{\widetilde{\eta}_i}{p^n} \in \calP_{\frakX / \frakS , n}$ in $\calP_{n,k}:=\calP_{\frakX/\frakS,n}/(p^k)$ and set, for any $I=(I_1,\hdots,I_d)\in \N^d,$
$$\xi^{[I]}=\prod_{i=1}^d\xi_i^{[I_i]}.$$
Let $\delta:\calP_{n,k} \ra \calP_{n,k}\otimes_{\Ox_{\frakX_k}}\calP_{n,k}$ the comultiplication map of the Hopf algebra $\calP_{n,k}.$ Then
$$\delta\left (\xi^{[I]}\right )=\sum_{\substack{a,b,c\in \N^d \\ a+b+c=I}} p^{n|c|} c!\begin{pmatrix}a+c \\ c \end{pmatrix} \begin{pmatrix}b+c \\ c \end{pmatrix}\xi^{[b+c]}\otimes \xi^{[a+c]}.$$
\end{lemma}

\begin{proof}
Let $\widetilde{\delta} : \calP_{\frakX / \frakS,n} \ra \calP_{\frakX / \frakS,n} \otimes_{\Ox_{\frakX}} \calP_{\frakX / \frakS,n}$ be the comultiplication map. By \ref{HopffrakP}, we have
\begin{alignat*}{2}
p^n\widetilde{\delta}\left ( \frac{\widetilde{\eta}_i}{p^n} \right ) = \widetilde{\delta} \left ( \widetilde{\eta}_i \right ) &= 1 \otimes \widetilde{\eta}_i + \widetilde{\eta}_i \otimes 1 + \widetilde{\eta}_i \otimes \widetilde{\eta}_i \\
&= p^n \left ( 1 \otimes \frac{\widetilde{\eta}_i}{p^n} + \frac{\widetilde{\eta}_i}{p^n} \otimes 1 + p^n \frac{\widetilde{\eta}_i}{p^n} \otimes \frac{\widetilde{\eta}_i}{p^n} \right ).
\end{alignat*}
By the flatness of $P_{\frakX / \frakS,n}$ over $\op{Spf} \Z_p$ \eqref{propflat}, we deduce that
$$
\widetilde{\delta}\left ( \frac{\widetilde{\eta}_i}{p^n} \right ) = 1 \otimes \frac{\widetilde{\eta}_i}{p^n} + \frac{\widetilde{\eta}_i}{p^n} \otimes 1 + p^n \frac{\widetilde{\eta}_i}{p^n} \otimes \frac{\widetilde{\eta}_i}{p^n}.
$$
It follows that
$$\delta(\xi_i)=1\otimes \xi_i+\xi_i\otimes 1+p^n\xi_i\otimes \xi_i.$$
Then, for $n\in \N,$
\begin{alignat*}{2}
\delta(\xi_i^{[n]}) &=  \left (1\otimes \xi_i+\xi_i\otimes 1+p^n\xi_i\otimes \xi_i\right )^{[n]} \\
&= \sum_{\substack{a,b,c\in \N \\a+b+c=n}}(1\otimes \xi_i)^{[a]}(\xi_i \otimes 1)^{[b]}(p^n\xi_i\otimes \xi_i)^{[c]} \\
&= \sum_{\substack{a,b,c\in \N \\a+b+c=n}}p^{nc}c!(\xi_i^{[b]}\xi_i^{[c]})\otimes (\xi_i^{[a]}\xi_i^{[c]}) \\
&= \sum_{\substack{a,b,c\in \N \\a+b+c=n}}p^{nc}c!\begin{pmatrix}a+c\\c \end{pmatrix}\begin{pmatrix}b+c \\ c \end{pmatrix} \xi_i^{[b+c]}\otimes \xi_i^{[a+c]}.
\end{alignat*}
The same holds for multi-indices hence the result.
\end{proof}

\begin{proof}[Proof of \ref{equivstrat}]
The proof is similar to \ref{prop39} so we just prove how (3) implies (2) and give a sketch of the other equivalences.
For simplicity, we suppose $k=1$ and we drop the subscript $k$ from our notation.

Suppose we are given the data (1). The data (2) is then obtained by setting
$$
\theta_l:\begin{array}[t]{clclc}
\calE & \ra & \calP_n^{\{l\}}\otimes_{\Ox_X} \calE & \xrightarrow{\epsilon_l} & \calE \otimes_{\Ox_X}\calP_n^{\{l\}} \\
x & \mapsto & 1\otimes x & \mapsto & \epsilon_l(1\otimes x).
\end{array}
$$
Conversly, if (2) is given then $\epsilon_l$ is the $\calP_n^{\{l\}}$-linearization of $\theta_l.$

Suppose we are given the data (2). The data (3) is obtained as follows: let
$$
\nabla:\begin{array}[t]{clc}
\calE & \ra & \calE\otimes_{\Ox_X}\calP_n^{\{1\}} \\
x & \mapsto & \theta_1( x)-x\otimes 1.
\end{array}
$$
Since $\theta_0=\op{Id}_{\calE},$ we get
$$\nabla(\calE) \subset \calE \otimes_{\Ox_X}\ov{\calI}_n^{\{1\}}.$$
In addition, if we denote by $p_1,p_2:P_{n} \ra X$ the canonical projections, then, for any local sections $a$ and $x$ of $\Ox_X$ and $\calE$ respectively, we have
\begin{alignat*}{2}
\nabla(ax) &= \theta_1(ax)-(ax)\otimes 1 = p_2^{\#}(a)\theta_1(x)-p_1^{\#}(a) (x\otimes 1) \\
&= p_2^{\#}(a) \left ( \theta_1(x)-x\otimes 1 \right )+ x\otimes \left ( p_2^{\#}(a) - p_1^{\#}(a) \right ) = a\nabla(x)+x\otimes d'(a).
\end{alignat*}
By \ref{propfinalMuzan1}, $\nabla$ corresponds to a $p^n$-connection on $\calE.$ The integrability follows from the commutativity of \eqref{diagstrconn1}.
Suppose now that an integrable $p^n$-connection $\calE \ra \calE \otimes_{\Ox_X}\omega^1_{X/S}$ is given. Let
$$\nabla:\calE \ra \calE \otimes_{\Ox_X}\ov{\calI}_{n}^{\{1\}}$$
be the corresponding connection given by \ref{propfinalMuzan1}. We construct the morphisms $\theta_l$ of (2) étale locally on $X.$ Suppose that the hypothesis \ref{loccoord} is satisfied and consider the notation introduced in \ref{lem116f}. We have a canonical isomorphism
$$\calP_{n}\xrightarrow{\sim} \Ox_X\langle \xi_1,\hdots,\xi_d \rangle.$$
Denote by $\left ( \partial_I \right )_{I\in \N^d}$ the dual basis of $\left ( \xi^{[I]} \right )_{I\in \N^d}.$ Let $\epsilon_1$ be the $\calP_n$-linear morphism defined by
$$
\epsilon_1:\begin{array}[t]{clc}
\calP_n^{\{1\}} \otimes_{\Ox_X}\calE & \ra & \calE\otimes_{\Ox_X}\calP_n^{\{1\}} \\
1\otimes x & \mapsto & \nabla(x)+x\otimes 1.
\end{array}
$$
We define an $\Ox_X$-linear morphisms
$$
\nabla_l:\mathscr{Hom}_{\Ox_X} \left ( \calP_n^{ \{l\} },\Ox_X \right ) \ra \mathscr{Hom}_{\Ox_X} \left ( \calP_n^{ \{l\} }\otimes_{\Ox_X} \calE,\calE \right ).
$$
For that, it is sufficient to define $\nabla_l$ on the basis $\left (\partial_I \right )_{|I| \le l}.$
For $l,l'\in \N$ and $I\in \N^d$ such that $l\ge l' \ge |I|,$ we consider a linear morphism of
$\mathscr{Hom}_{\Ox_X} \left ( \calP_n^{ \{l'\} }\otimes_{\Ox_X} \calE,\calE \right )$
as a linear morphism of
$\mathscr{Hom}_{\Ox_X} \left ( \calP_n^{ \{l\} }\otimes_{\Ox_X} \calE,\calE \right )$
via the canonical projection $\calP_n^{ \{ l \} } \ra \calP_n^{ \{ l' \} }.$ 
For any $\Ox_X$-linear morphism $f:\calP_n^{\{1\}}\ra \Ox_X,$ let $\nabla_1(f)=(\op{Id}_{\calE}\otimes f)\circ \epsilon_1.$
The integrability of $\nabla$ implies that the differential operators $\nabla_1(\partial_{\epsilon_i})$ pairwise commute and we can define $\nabla_{l}(\partial_N),$ for any multi-index $N=(n_1,\hdots,n_d)\in\mathbb{N}^d$ of length $|N|\le l,$ by
$$
\nabla_l(\partial_N)=\prod_{i=1}^d\prod_{j=0}^{n_i-1}(\nabla_1(\partial_{\epsilon_i})-p^nj),
$$
where $p^nj$ denotes $p^nj \op{Id}_{\calE},$ considered as a differential operator of order 1.
Then
$$
\nabla_l\left (\partial_I \right ) = \nabla_{l'}\left ( \partial_I \right ).
$$
We can hence drop the subscript $l$ in the rest of this proof. We can also show that, for $I,J\in \N^d,$
$
\nabla \left ( \partial_I \circ \partial_J \right ) = \nabla \left ( \partial_I \right ) \circ \nabla \left ( \partial_J \right ).
$
Set
$
\theta_l:\begin{array}[t]{clc}
\calE & \ra & \calE \otimes_{\Ox_X}\calP_n^{\{l\}}. \\
x & \mapsto & \sum_{|I| \le l}\nabla(\partial_I)(1\otimes x) \otimes \xi^{[I]}
\end{array}
$
We check the commutativity of the diagram \eqref{diagstrconn1}. On one hand, by \ref{lem116f}, we have
\begin{equation}\label{eqfinalMuzan5}
\begin{alignedat}{2}
(\op{Id}_{\calE}\otimes \delta^{l,l'})\circ \theta_{l+l'}(x) &= (\op{Id}_{\calE}\otimes \delta^{l,l'})\left (\sum_{|I|\le l+l'} \nabla(\partial_I)(1\otimes x)\otimes \xi^{[I]} \right ) \\
&= \sum_{|I|\le l+l'}\sum_{\substack{a+b+c=I \\ |b+c|\le l \\ |a+c|\le l'}} p^{n|c|} c!\begin{pmatrix}a+c \\ c \end{pmatrix} \begin{pmatrix}b+c \\ c \end{pmatrix}\nabla(\partial_I)(1\otimes x) \otimes \xi^{[b+c]}\otimes \xi^{[a+c]}.
\end{alignedat}
\end{equation}
On the other hand, we have
\begin{equation}\label{eqfinalMuzan6}
\begin{alignedat}{2}
\theta_{l'}(x) &= \sum_{|I|\le l'} \nabla(\partial_I)(1\otimes x)\otimes \xi^{[I]}. \\
(\theta_l \otimes \op{Id})\circ \theta_{l'}(x) &= \sum_{\substack{|I|\le l' \\ |J|\le l}} \nabla(\partial_J)(1\otimes (\nabla(\partial_I)(1\otimes x)))\otimes \xi^{[J]}\otimes \xi^{[I]} \\
&= \sum_{\substack{|I|\le l' \\ |J|\le l}} \nabla(\partial_J) \circ \nabla(\partial_I)(1\otimes x)\otimes \xi^{[J]}\otimes \xi^{[I]} \\
&= \sum_{\substack{|I|\le l' \\ |J|\le l}} \nabla(\partial_J \circ \partial_I)(1\otimes x)\otimes \xi^{[J]}\otimes \xi^{[I]}.
\end{alignedat}
\end{equation}
We have to prove that \eqref{eqfinalMuzan5} and \eqref{eqfinalMuzan6} are equal. For that, we prove that, for $I,J \in \N^d$ such that $|J| \le l$ and $|I| \le l',$
$$
\partial_J \circ \partial_I = \sum_{\substack{a+c=I \\ b+c=J}} p^{n|c|}c!\begin{pmatrix}a+c \\ c \end{pmatrix} \begin{pmatrix}b+c \\ c \end{pmatrix}\partial_{a+b+c}.
$$
We proceed by induction on $|J|.$ For $|J|=0,$ the equality is clear.
Fix an integer $m$ and suppose the result is true for all $J$ and $I$ such that $|J|=m.$ Let $1\le i\le d.$
Then
\begin{alignat*}{2}
\partial_{J+\epsilon_i}\circ \partial_I &= (\partial_{\epsilon_i}\circ \partial_J -p^nJ_i\partial_J)\circ \partial_I \\
&= \sum_{\substack{a+c=I\\ b+c=J}} p^{n|c|}c! \begin{pmatrix}a+c \\ c \end{pmatrix} \begin{pmatrix}b+c \\ c \end{pmatrix} \left (\partial_{\epsilon_i} \circ \partial_{a+b+c}-p^n J_i \partial_{a+b+c} \right ) \\
&=  \sum_{\substack{a+c=I\\ b+c=J}} p^{n|c|}c! \begin{pmatrix}a+c \\ c \end{pmatrix} \begin{pmatrix}b+c \\ c \end{pmatrix} \left (\partial_{a+b+c+\epsilon_i} +p^n(a_i+b_i+c_i) \partial_{a+b+c}-p^n J_i \partial_{a+b+c} \right ) \\
&= \sum_{\substack{a+c=I\\ b+c=J}} p^{n|c|}c! \begin{pmatrix}a+c \\ c \end{pmatrix} \begin{pmatrix}b+c \\ c \end{pmatrix} \left (\partial_{a+b+c+\epsilon_i}+p^n a_i \partial_{a+b+c} \right ).
\end{alignat*}
Applying this to $\xi^{[K]}$ for $K\in \N^d,$ and considering the convention $\begin{pmatrix} \alpha \\ \beta \end{pmatrix}=0$ if $\beta<0,$ we get
\begin{align}\label{yeaa1}
\begin{split}
p^{n|I+J-K|+n}(I+J-K+\epsilon_i)! \begin{pmatrix}I\\K-J-\epsilon_i \end{pmatrix} \begin{pmatrix}J\\ K-I-\epsilon_i \end{pmatrix}\\+p^{n|I+J-K|+n}(I+J-K)!\begin{pmatrix}I\\K-J\end{pmatrix} \begin{pmatrix}J\\ K-I \end{pmatrix}(K_i-J_i).
\end{split}
\end{align}
Applying
$$
\sum_{\substack{a+c=I \\ b+c=J+\epsilon_i}} p^{n|c|}c!\begin{pmatrix}a+c \\ c \end{pmatrix} \begin{pmatrix}b+c \\ c \end{pmatrix}\partial_{a+b+c}
$$
to $\xi^{[K]},$ we get
\begin{equation}\label{yeaa2}
p^{n|I+J-K|+n}(I+J-K+\epsilon_i)!\begin{pmatrix}I \\ K-J-\epsilon_i \end{pmatrix} \begin{pmatrix}J+\epsilon_i \\ K-I \end{pmatrix}.
\end{equation}
To show that \eqref{yeaa1} and \eqref{yeaa2} are equal, recalling the notation \ref{Not9}, we have to prove the equality
\begin{alignat*}{2}
&(I_i+J_i-K_i+1) \begin{pmatrix}I_i\\K_i-J_i-1 \end{pmatrix} \begin{pmatrix}J_i\\ K_i-I_i-1 \end{pmatrix}+\begin{pmatrix}I_i\\K_i-J_i\end{pmatrix} \begin{pmatrix}J_i\\ K_i-I_i \end{pmatrix}(K_i-J_i) \\
=& (I_i+J_i-K_i+1)\begin{pmatrix}I_i \\ K_i-J_i-1 \end{pmatrix} \begin{pmatrix}J_i+1 \\ K_i-I_i \end{pmatrix}.
\end{alignat*}
This is equivalent to
\begin{alignat*}{2}
&\begin{pmatrix}I_i\\K_i-J_i\end{pmatrix} \begin{pmatrix}J_i\\ K_i-I_i \end{pmatrix}(K_i-J_i) \\
=& (I_i+J_i-K_i+1)\begin{pmatrix}I_i \\ K_i-J_i-1 \end{pmatrix} \left ( \begin{pmatrix}J_i+1 \\ K_i-I_i \end{pmatrix} - \begin{pmatrix}J_i\\ K_i-I_i-1 \end{pmatrix} \right ).
\end{alignat*}
But
$$
\begin{pmatrix}J_i+1 \\ K_i-I_i \end{pmatrix} - \begin{pmatrix}J_i\\ K_i-I_i-1 \end{pmatrix}= \begin{pmatrix}J_i \\ K_i-I_i \end{pmatrix}.
$$
So we have to prove that
\begin{equation*}
\begin{pmatrix}I_i\\K_i-J_i\end{pmatrix} (K_i-J_i) = (I_i+J_i-K_i+1)\begin{pmatrix}I_i \\ K_i-J_i-1 \end{pmatrix}.
\end{equation*}
This is immediate.
\end{proof}


\begin{thebibliography}{99}

\bibitem{Ahmed2010}{\sc A. Abbes}, Éléments de Géométrie Rigide: Volume I. Construction et Étude Géométrique des Espaces Rigides. Suisse: Birkhauser, 2010.

\bibitem{SGA43}{\sc M. Artin, A. Grothendieck, J. L. Verdier}, Théorie des topos et cohomologie étale des schémas. Tome 1, Séminaire de Géométrie Algébrique du Bois-Marie 1963–1964 (SGA 4). Lecture Notes in Mathematics 305, Springer 1973.

\bibitem{Ber74}{\sc P. Berthelot}, Cohomologie cristalline des schémas de caractéristique $p>0,$ Springer-Verlag, LNM 407, (1974).

\bibitem{Ogus78}{\sc P. Berthelot, A. Ogus}, Notes on Crystalline Cohomology. (MN-21). Princeton University Press, 1978.

\bibitem{Bourbaki}{\sc N. Bourbaki}, Algèbre commutative, Chapitres 1 à 4. Springer Berlin, Heidelberg, 2007.

\bibitem{Giraud}{\sc J. Giraud}, Cohomologie non abélienne. Springer Berlin, Heidelberg, 1971.

\bibitem{EGA1}{\sc A. Grothendieck, J. Dieudonné}, Éléments de Géométrie algébrique I. Publications Mathématiques de l'IHÉS, Tome 4 (1966), pp. 5-214.

\bibitem{EGA43}{\sc A. Grothendieck}, Éléments de géométrie algébrique : IV. Étude locale des schémas et des morphismes de schémas, Troisième partie. Publications Mathématiques de l'IHÉS, Tome 28 (1966), pp. 5-255.

\bibitem{Raynaud71}{\sc A. Grothendieck, M. Raynaud}, Revêtements étales et groupe fondamental (SGA 1). Lecture notes in mathematics. Volume 224. Springer.

\bibitem{INT}{\sc L. Illusie, C. Nakayama, T. Tsuji}, On log flat descent. Proc. Japan Acad. Ser. A Math. Sci. 89 (1) 1 - 5, January 2013.

\bibitem{Kat89}{\sc K. Kato}, Logarithmic structures of Fontaine-Illusie. Algebraic Analysis, Geometry and Number Theory (Baltimore, MD, 1988), 1989 pp. 191-224.

\bibitem{Kat19}{\sc K. Kato}, Logarithmic structures of Fontaine-Illusie. II. Logarithmic flat topology. ArXiv:1905.10678.

\bibitem{Saito04}{\sc K. Kato, T. Saito}, On the conductor formula of Bloch. Publications Mathématiques de l'IHÉS, Volume 100 (2004), pp. 5-151.

\bibitem{Katz}{\sc N.M. Katz}, Nilpotent connections and the monodromy theorem: Applications of a result of turrittin. Publications Mathématiques de L’Institut des Hautes Scientifiques 39, 175–232 (1970).

\bibitem{Lor2000}{\sc P. Lorenzon}, Indexed algebras associated to a log structure and a theorem of p-descent on log schemes. manuscripta math. 101, 271–299 (2000).

\bibitem{Mon}{\sc C. Montagnon}, Généralisation de la théorie arithmétique des D-modules à la géométrie logarithmique, Ph.D. thesis, L'université de Rennes I (2002).

\bibitem{Ogus94}{\sc A. Ogus}, F-crystals, Griffiths Transversality, and the Hodge Decomposition. Volume 221 de Asterisque / Société mathématique de France, 1994.

\bibitem{Ogus2018}{\sc A. Ogus}, Lectures on Logarithmic Geometry. Cambridge University Press, 2018.

\bibitem{OgusVol}{\sc A. Ogus, V. Vologodsky}, Nonabelian Hodge theory in characteristic $p$. Publications Mathématiques de l'IHÉS, Volume 106 (2007), pp. 1-138.

\bibitem{Ohkawa}{\sc S. Ohkawa}, On logarithmic nonabelian Hodge theory of higher level in characteristic $p$. Rend. Sem. Mat. Univ. Padova 134 (2015), pp. 47–91

\bibitem{Oyama}{\sc H. Oyama}, PD Higgs crystals and Higgs cohomology in characteristic $p.$ J. Algebraic Geom. 26 (2017), 735-802 

\bibitem{Schepler}{\sc D. Schepler}, Logarithmic nonabelian Hodge theory in characteristic $p,$ Ph.D. thesis (2005).

\bibitem{Shiho}{\sc A. Shiho}, Notes on generalizations of local Ogus-Vologodsky correspondence. J. Math. Sci. Univ. Tokyo 22 (2015),
793–875.

\bibitem{Tsuji}{\sc T. Tsuji}, Saturated morphisms of logarithmic schemes. Tunisian Journal of Mathematics. Vol. 1, No 2, 2019.

\bibitem{DXU19}{\sc D. Xu}, Lifting the Cartier transform of Ogus-Vologodsky modulo $p^n.$ Mémoires de la Société Mathématique de France 163 (2019).

\bibitem{SP}{\sc The stacks project authors}, The stacks project.


\end{thebibliography}
\end{document}